\documentclass{amsart}
\usepackage{enumitem}
\usepackage[figuresleft]{rotating}
\usepackage{wrapfig}
\usepackage[letterpaper,margin=1in,centering,marginparwidth=.9in,marginparsep=.01in]{geometry}
\usepackage{amsthm,amsmath,amssymb,amsfonts,mathtools,microtype,cancel}
\usepackage[utf8]{inputenc} 
\usepackage[T1]{fontenc}    
\usepackage[hidelinks,pagebackref]{hyperref}       
\usepackage{booktabs}       
\usepackage{nicefrac}       
\usepackage[widespace]{fourier} 
\usepackage[dvipsnames]{xcolor}
\usepackage{pgf}

\allowdisplaybreaks[1]

\newcommand\hidesubsection[2][]{{\def\addcontentsline##1##2##3{}\expandafter\ifx\expandafter\relax\detokenize{#1}\relax\subsection*{#2}\else\subsection*{#2\textnormal{\ (#1)}}\fi\strut}}
\DeclareMathOperator{\sinc}{sinc}
\DeclareSymbolFontAlphabet{\mathcal}   {symbols}
\DeclareMathAlphabet{\mathcal}{OMS}{Zplm}{m}{n}
\SetMathAlphabet{\mathcal}{bold}{OMS}{zplm}{b}{n}

\usepackage{graphicx}
\graphicspath{ {figs/} }
\usepackage{doi}
\usepackage{placeins}
\usepackage{caption}
\usepackage{subcaption}
\usepackage{tikz}
\usepackage{tkz-euclide}
\usepackage{yfonts}
\usepackage[capitalize,nameinlink,noabbrev]{cleveref}
\usetikzlibrary{spy}
\usetikzlibrary {angles,quotes}
\def\pcref#1{\cref{#1}, page \pageref{#1}}
\newtheorem{theorem}{Theorem}[section]
\newtheorem{ltheorem}{Theorem} 
\newtheorem{corollary}[theorem]{Corollary}
\newtheorem{proposition}[theorem]{Proposition}
\newtheorem{lemma}[theorem]{Lemma}
\theoremstyle{definition}
\newtheorem{definition}[theorem]{Definition}
\newtheorem{notation}[theorem]{Notation}
\crefname{notation}{Notation}{Notations}
\crefname{ltheorem}{Theorem}{Theorem}
\newtheorem{remark}[theorem]{Remark}
\numberwithin{equation}{section}
\usetikzlibrary{intersections, calc, angles}
\def\bigstrut{\rule[-0.4\baselineskip]{0pt}{1.2\baselineskip}}
\def\Bigstrut{\rule[-0.5\baselineskip]{0pt}{1.3\baselineskip}}
\newif\ifcomments
\newcounter{commentlabel}
  \newlength{\commentlift}
\makeatletter\@mparswitchfalse\def\commentfl@g{%
    \vbox to\z@{%
      \setlength{\fboxrule}{0.75pt}%
      \setlength{\fboxsep}{0.75pt}%
      \setlength{\commentlift}{1ex}%
      \addtolength{\commentlift}{\fboxrule}%
      \addtolength{\commentlift}{\fboxsep}%
\vss\color{blue}\rlap{\rlap{\vrule\@height\commentlift\@width\fboxrule}\raise \commentlift%
      \hbox{\fcolorbox{blue}{blue!3}{\normalfont\tiny\bfseries\thecommentlabel}}}}}
  \DeclareRobustCommand{\COMMENT}[1]{\commentstrue\stepcounter{commentlabel}\unskip%
\hypertarget{\the\value{section}-\the\value{commentlabel}}{}%
\edef\WRITECOM##1{\noexpand\write\COM{\thecommentlabel, S.##1}}\WRITECOM{\the\value{section}, \noexpand\hyperlink{\the\value{section}-\the\value{commentlabel}}{p.\the\value{page}}: {#1}}\@bsphack%
    \commentfl@g%
    \marginpar{\fontfamily{EBGaramond-TLF}\selectfont\noindent\lineskip=1pt\lineskiplimit=\maxdimen\raggedright%
    \tiny\def\Hrule{\strut\hrule height 1pt\strut}\textbf{\color{blue}\llap{\textbullet}\thecommentlabel:}{\color{blue}\thinspace#1\endgraf}}%
  \@esphack}%
  \DeclareRobustCommand{\COMMENTINLINE}[1]{\commentstrue\stepcounter{commentlabel}%
    \write\COM{\thecommentlabel, S.\the\value{section}, p.\the\value{page}: {#1}}\@bsphack%
      \noindent\setlength{\fboxsep}{0pt}\newline\fcolorbox{blue}{lightgray!20}{\begin{minipage}{\textwidth}\tiny\textbf{\color{red}\thecommentlabel}:\thinspace\Tiny\bfseries#1\par\end{minipage}}\allowbreak%
  \@esphack}\makeatother
\newwrite\COM
\immediate\openout\COM=comments
\AtEndDocument{\ifcomments\def\red{\color{red}}\section*{Marginal Comments}\label{commentstex}{\let\tiny\relax\let\Tiny\relax\let\Hrule\newline\def\pullout#1, S.{\strut\llap{\textbf{#1: }}\S\ }\obeylines\parindent0pt\everypar{\normalsize\normalfont\pullout}\immediate\closeout\COM\def\MT_extended_eqref:n {\eqref}\def\math@bb{\mathbb }\input comments }\fi}
\makeatother
\let\NEWPAGE\newpage
\def\COMMENT#1{}\def\COMMENTINLINE#1{}\let\NEWPAGE\relax
\def\red{\color{red}}
\def\green{\color{green!30!black}}
\let\implies\Rightarrow
\newcommand{\arc}[1]{\wideparen{#1}}
\newcommand{\Int}{\operatorname{int}}
\newcommand{\Fc}{\mathcal F}

\let\Hrule\relax

\newcommand{\p}{\ensuremath{\text{\upshape p}}}
\newcommand{\phistar}{\zeta}
\newcommand{\Phistar}{\varsigma}
\newcommand{\Bcal}{\ensuremath{\mathcal{B}}}
\newcommand{\Fcal}{\ensuremath{\mathcal{F}}}
\newcommand{\Scal}{\ensuremath{\mathcal{S}}}
\newcommand{\constant}{\ensuremath{\Psi}}
\newcommand{\const}{\ensuremath{\text{const.}}}
\newcommand{\dfn}{\ensuremath{\mathbin{{:}{=}}}}
\newcommand{\nfd}{\ensuremath{\mathbin{{=}{:}}}}
\newcommand{\aeq}{\ensuremath{\mathbin{\overset{\text{\upshape\ae}}{=}}}}
\newcommand{\Mout}{{\color{red}M}^{\text{{\color{red}\upshape out}}}}
\newcommand{\Min}{{\color{blue}M}^{\text{{\color{blue}\upshape in}}}}
\newcommand{\Mrout}{\ensuremath{{\color{red}M^{\text{\upshape out}}_r}}}
\newcommand{\Mrin}{\ensuremath{{\color{blue}M^{\text{\upshape in}}_r}}}
\newcommand{\MRout}{\ensuremath{{\color{red}M^{\text{\upshape out}}_R}}}
\newcommand{\MRin}{\ensuremath{{\color{blue}M^{\text{\upshape in}}_R}}}
\newcommand{\Nin}{\ensuremath{{\color{blue}N^{\text{\upshape in}}}}}
\newcommand{\Nout}{\ensuremath{{\color{red}N^{\text{\upshape out}}}}}
\def\PET{\mathbb{Pet}}
\def\LEM{\mathbb{Lem}}
\newcommand{\boris}{\strut\textcolor{red}{\emph{\bfseries--Boris says: }}}
\newcommand{\wentao}{\strut\textcolor{green!70!black}{\emph{\bfseries--Wentao says: }}}
\title{Uniformly hyperbolic strictly convex billiards}
\author{Wentao Fan}
\address{School of Mathematical Sciences (Tianjiabing Education College)\\
Hebei Normal University\\
No. 20 South 2nd Ring East Road\\ Yuhua District, Shijiazhuang, Hebei, 050024, China}
\email{fanwentao@hebtu.edu.cn\\
wentao.fan@tufts.edu}
\author{Boris Hasselblatt}
\address{Department of Mathematics\\
        Tufts University\\
        Medford, MA 02155}
\email{boris.hasselblatt@tufts.edu}

\hypersetup{
pdftitle={Uniform Hyperbolicity and ergodicity of asymmetric lemon Billiard},
pdfkeywords={Expansion, Uniform Hyperbolic, Wojtkowski Billiards, Cone, Wave Fronts Curvature, singularity curve, local ergodic theorem, global ergodic},
}
\begin{document}
\begin{abstract}
Can a strictly convex billiard be ergodic? For the billiard in a suitable intersection of two circles, we establish \emph{uniform} expansion of vectors in an invariant cone family under the return to a well-chosen section. This is the first and most onerous step towards a forthcoming proof of ergodicity (the Chen--Mohr--Zhang--Zhang conjecture) and statistical properties. Together, this appears to be the first time that a strictly convex billiard is shown to be ergodic.
\end{abstract}
\maketitle
\setcounter{tocdepth}{2}
\section{Main result}
\label{sec:IntroMainResult}
\label{sec:imr}
\strut\COMMENT{\small\bfseries\green For arXiv posting and submission define \textbackslash COMMENT and \textbackslash NEWPAGE to do nothing---around line 106 in the TeX file}Billiard tables \cite{WhatIsABilliard} with hyperbolic time-evolution are those for which one can establish good statistical properties, beginning with ergodicity, a form of measure-theoretic indecomposability \cite{WhatIsErgodic}. Hyperbolicity, however, relies on dispersing boundary pieces, or on restrictions on focusing boundary arcs that are incompatible with strict convexity. A review article on ``Some new surprises in chaos'' \cite{BunimovichSurprise} presented \emph{asymmetric lemon billiards} as shown in \cref{fig:BilliardTable} as ``a new class of chaotic focusing billiards\dots that clearly violates the main condition considered to be necessary for chaos in focusing
billiards.'' These billiards are strictly convex, but even that article dared not address their ergodicity. 

So, can \emph{strictly} convex billiards be ergodic? 
\cref{RMASUH} affirmatively answers a related open question: is there a strictly convex billiard that admits a global section with \emph{uniformly} hyperbolic return map? 

For over a decade it has been conjectured that asymmetric lemon billiards are ergodic. This article is devoted to the most onerous step towards a resolution of this conjecture---we establish \emph{uniform} hyperbolicity  (not only nonzero Lyapunov exponents everywhere) of the return map to a suitable section in \cref{RMASUH}. Later work will control singularity curves (which here look different than elsewhere in the literature---near vertical with sign change of slope), and then use Hopf chains to establish ergodicity and stronger statistical properties.
\begin{ltheorem}\label{RMASUH}
For a feasible $\phistar\in \big(0,\tan^{-1}{({1}/{3})}\big)$ as in \cref{def:ApettrajectoryExpansion}, there exists an $R_{\text{\upshape HF}}(r,\phistar)$
such that if $R>R_{\text{\upshape HF}}(r,\phistar)$ in \cref{fig:BilliardTable}, 
then:
\begin{enumerate}
    \item The ``half-quadrant'' cone family $C_x\dfn\big\{(d\phi,d\theta)\bigm|\frac{d\theta}{d\phi}\in [0,1]\big\}$ is strictly invariant under the return map $\hat F$ to the global section $\hat{M}$ (from \eqref{def:Mhat} and \cref{fig:hatUplusUminus}) of next to last or sole collisions with $C_r$.
    \item This return map is \emph{uniformly} hyperbolic: there exist $c>0$ and $\Lambda>1$ such that if $x\in\hat{M}$ is nonsingular, $dx\in C_x$ and $n\in\mathbb{N}$, then $\|D_x\hat{F}^n(dx)\|>c\Lambda^n\|dx\|$ (\cref{Corollary:EuclideanExpansionMhat,corollary:MroutreturnUniformExpansion,Col:EuclideanExpansionMhat1}).
\end{enumerate}
\end{ltheorem}
    \begin{figure}[ht]
    \begin{center}
    \begin{tikzpicture}[scale=0.7]
    \tkzDefPoint(-6,0){OR};
    \tkzDefPoint(0,0){Or};
    \tkzDefPoint(2,0){C};
    \pgfmathsetmacro{\XValueArc}{2*sin(60)}
    \tkzDefPoint(\XValueArc,-1){A};
    \tkzDefPoint(\XValueArc,1){B};
    \clip(-14,-2.1) rectangle (4, 2.2);
    \node[] at ([shift={(0.8,1.0)}]B) {\textcolor{blue}{\small circle $C_R$}};
    \node[] at ([shift={(-2.4,-1.6)}]Or) {\small circle $C_r$};
    \tkzDrawArc[very thick](OR,A)(B);
    \tkzDrawArc[very thick](Or,B)(A);
    \tkzDrawArc[thick,dashed, blue](OR,B)(A);
    \tkzDrawArc[thick,dashed](Or,A)(B);
    \tkzDrawPoints(OR,Or,A,B);
    \tkzLabelPoint[right](A){$A$};
    \tkzLabelPoint[right](B){$B$};
    \tkzLabelPoint[below](OR){$O_R$};
    \tkzLabelPoint[below](Or){$O_r$};
    \draw [very thin, dashed] (OR) --(5,0);
    \coordinate (C) at (5,0);
    \tkzDrawSegments[very thin,dashed](OR,B);
    \tkzDrawSegments[very thin,dashed](Or,B);
    \tkzDrawSegments[very thin,dashed](Or,A);
    \tkzLabelSegment[above left=0.5pt and 3pt](OR,B){$R$};
    \tkzLabelSegment[above left=-10pt and -5pt](Or,A){$r$};
    \draw (0.6,0) arc (0:30:0.6);
    \node[] at (15:1)  {$\phistar$};
    \node[] at (-3.2,-0.3){$b$};
    \tkzMarkAngle[arc=lll,ultra thin,size=1.9](C,OR,B);
    \tkzLabelPoint[right]([shift={(2.4,0.15)}]OR){\small $\Phistar$};
    \end{tikzpicture}
    \end{center}
    \caption{Lemon billiard Table} 
    \label{fig:BilliardTable}
    \end{figure}
Mathematical billiard systems describe a point particle moving freely in a domain with piecewise smooth boundary and reflecting off the boundary the way a photon does in a mirror, with equal angles before and after collision (specular reflection) \cite{WhatIsABilliard,Tabachnikov}. These are dynamical systems which Birkhoff introduced as a simplified limit of a range of mechanical systems with two degrees of freedom in which ``the formal side, usually so formidable in dynamics, almost completely disappears, and only the interesting qualitative questions need to be considered''  \cite[p.~288f]{HF}. For instance, under contraction of the $z$-coordinate (say), an ellipsoid will eventually be squashed to an elliptic billiard in the $xy$-plane. In terms of paradigms for dynamical features and methods of study, billiards can be grouped into convex billiards (elliptic), polygonal billiards (parabolic), and dispersing billiards (hyperbolic) \cite{Playground}. The latter are the subject of the present article. The core concern and criterion is the existence of a cone field in the tangent bundle which is invariant under the (differential of the) dynamics and (at least eventually) strictly so. One highlight of the study of these kinds of billiards are criteria for nonuniform hyperbolicity (i.e., for the existence of an eventually strictly invariant cone field) in the presence of focusing boundary circle arcs (the circle to which such an arc belongs must lie inside the billiard table \cite{MR848647,ICLE1985}; see also \cite[p.~144]{Tabachnikov}) by  Wojtkowski and Bunimovich; notably, these criteria preclude strict convexity (\cref{WCriteriaNotConvex}). Another is the work by Chernov and Markarian \cite{cb}, which carefully presents the needed technology and is a central reference for the present undertaking. In more recent years, significant efforts have been directed at the study of billiards which are expected to be hyperbolic but do not satisfy those criteria for hyperbolicity. Our main result was motivated by recent work which established the existence (and eventual strict invariance) of such cone fields for particular kinds of such billiards---the \emph{lemon billiards} shown in \cref{fig:BilliardTable} for suitable parameters \cite{oaedhalb16,hoalb}. 

\cref{RMASUH} takes a substantial and essential next step by establishing \emph{uniform} hyperbolicity when the asymmetry is large enough (for $R$ sufficiently larger than $300\,000r$). This is not a small tweak up but crucial leap for achieving a core aim associated with hyperbolic dynamics: to show that the natural invariant measure is ergodic, and more, i.e., that statistical methods are suitable for describing the complexity of the time evolution. That aim will be realized in subsequent work, which thereby proves the Chen--Mohr--Zhang--Zhang conjecture \cite[Observations, p.\ 3]{eceifcb13}. Specifically, uniform expansion provides the basis for carrying out the Sinai--Chernov Ansatz, the standard method for establishing ergodicity from hyperbolicity.\COMMENT{This is the one place where one could mention Young towers if it is indeed the case that uniform hyperbolicity is a basis for them.}

It is important to clarify what exactly is proved to be uniformly hyperbolic in \cref{RMASUH}. Our introduction of billiards suggested a flow, and indeed, the billiard flow is the unit-speed (say) geodesic flow on the region in question with respect to the Euclidean metric, plus collisions: in the tangent space at a boundary point, vectors are identified that are related by reflection over the tangent line to the boundary at that point. Instead, we  study the billiard \emph{map}, which associates with an inward-pointing (post-collision) vector at a boundary point the vector after the next collision with the boundary. Thus, the phase space (or state space) consists of the product of the boundary (a topological circle) with the interval $[0,\pi]$ of possible angles of inward-pointing vectors with the tangent line. In fact, wanting a smooth structure, our phase space is the union of two rectangles whose bases are the two circle arcs. 

It is, however, not this billiard map we prove to be uniformly hyperbolic---and indeed, it is not. Diameters of $C_r$ are nonhyperbolic periodic points. We instead choose a subset  $\hat{M}$ (in \eqref{def:Mhat}, shown in red in \cref{fig:hatUplusUminus}) for which we prove that the first-return map is uniformly hyperbolic and that almost every point of phase space visits this set (therefore called a global section) in past and future, so images of  $\hat{M}$ under the billiard map cover almost all of the phase space (but not, for instance, the aforementioned nonhyperbolic periodic points). Here, the first-return map associates with any point of $\hat{M}$ the first of its successive images under the billiard map that lies in $\hat{M}$. ``Almost all'' is sufficient for our eventual purpose: in the sense of ergodic theory, $\hat{M}$ is a global section, and ergodicity of the return map implies ergodicity of the natural Liouville volume in the phase space for the full billiard map, and also for the billiard flow.

We hasten to note that the billiard map is not well-defined at corners because the tangent line and hence the billiard map has a jump discontinuity there. Usually, corners are eliminated from the domain of the billiard map. However, since its one-sided limits exist, it is natural to consider the continuous extension to the closure of each of the two components of the phase space. This turns out to be useful, for instance, in \cref{SBS1petalbilliardorbits}.

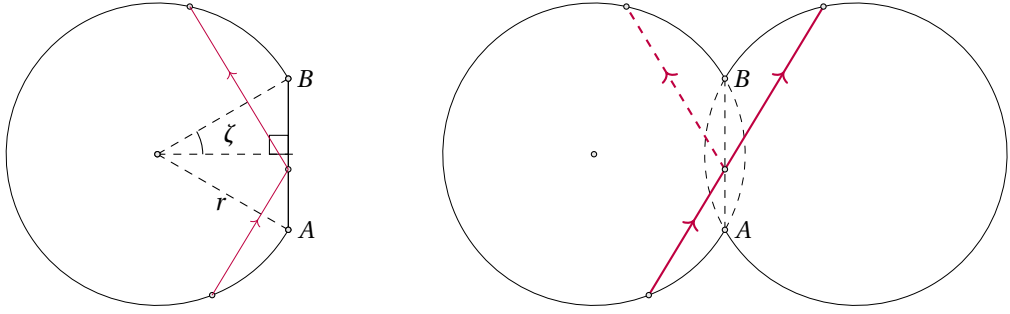
\begin{figure}[hb]
\strut\hfill
    \begin{tikzpicture}
    \tkzDefPoint(0,0){Or};
    \tkzDefPoint(2,0){C};
    \pgfmathsetmacro{\XValueArc}{2*sin(60)}
    \tkzDefPoint(\XValueArc,-1){A};
    \tkzDefPoint(\XValueArc,1){B};
    \tkzDrawArc[name path=Cr](Or,B)(A);
    \tkzDrawPoints(Or,A,B);
    \draw[name path=lineAB,thin](A)--(B);
    \tkzLabelPoint[right](A){$A$};
    \tkzLabelPoint[right](B){$B$};
    \draw [name path=xaxis, thin, dashed] (Or) --(1.83,0);
    \path[name intersections={of=xaxis and lineAB,by={C}}];
    \coordinate (Or_mirror) at ($(Or)!2!(C)$);
    \draw[thin](A)--(B);
    \tkzMarkRightAngle(B,C,Or);
    \tkzDrawSegments[thin,dashed](Or,B);
    \tkzDrawSegments[thin,dashed](Or,A);
    \tkzLabelSegment[above left=-10pt and -5pt](Or,A){$r$};
    \draw (0.6,0) arc (0:30:0.6);
    \node[] at (15:1)  {$\phistar$};
    
    \coordinate (P) at (\XValueArc,-0.2);
    \tkzDefPoint(2*cos(-1.2),2*sin(-1.2)){T};
    \tkzDefPointBy[projection=onto B--A](T)\tkzGetPoint{H}; 
    \coordinate (TReflect) at ($(T)!2.0!(H)$);
    \coordinate (QQ) at ($(TReflect)!2.4!(P)$);
    \path[name path=linePQQ](P)--(QQ);
    \path[name intersections={of=Cr and linePQQ,by={Q}}];
    \tkzDrawPoints(Or,A,B,P,T,Q);
    \begin{scope}[very thin,decoration={
    markings,
    mark=at position 0.6 with {\arrow{>}}}
    ] 
    \draw[purple,postaction={decorate}](T)--(P);
    \draw[purple,postaction={decorate}](P)--(Q);
    \end{scope}

    \end{tikzpicture}\hfill
    \begin{tikzpicture}
    \tkzDefPoint(0,0){Or};
    \tkzDefPoint(2,0){C};
    \pgfmathsetmacro{\XValueArc}{2*sin(60)}
    \tkzDefPoint(\XValueArc,-1){A};
    \tkzDefPoint(\XValueArc,1){B};
    \tkzDrawArc[name path=Cr](Or,B)(A);
    \draw[name path=lineAB,thin,dashed](A)--(B);
    \tkzLabelPoint[right](A){$A$};
    \tkzLabelPoint[right](B){$B$};
    \path[name path=xaxis, thin, dashed] (Or) --(1.83,0);
    \path[name intersections={of=xaxis and lineAB,by={C}}];
    \coordinate (Or_mirror) at ($(Or)!2!(C)$);
    \tkzDrawArc[name path=mirrorCr](Or_mirror,A)(B);
    \tkzDrawArc[dashed](Or,A)(B);
    \tkzDrawArc[dashed](Or_mirror,B)(A);
        \coordinate (P) at (\XValueArc,-0.2);
    \tkzDefPoint(2*cos(-1.2),2*sin(-1.2)){T};
    \tkzDefPointBy[projection=onto B--A](T)\tkzGetPoint{H};

    \coordinate (TReflect) at ($(T)!2.0!(H)$);
    \coordinate (QQ) at ($(TReflect)!2.4!(P)$);
    \path[name path=linePQQ](P)--(QQ);
    \path[name intersections={of=Cr and linePQQ,by={Q}}];
    \tkzDefPointBy[projection=onto B--A](Q)\tkzGetPoint{G};
    \coordinate (QReflect) at ($(Q)!2.0!(G)$);
    \begin{scope}[thick,decoration={
    markings,
    mark=at position 0.6 with {\arrow{>}}}
    ] 
    \draw[purple,postaction={decorate}](T)--(P);
    \draw[purple,postaction={decorate}](P)--(QReflect);
    \draw[purple,postaction={decorate},dashed](P)--(Q);
    \end{scope}
    \tkzDrawPoints(Or,A,B,P,T,Q,QReflect);
    \end{tikzpicture}
    \hfill\strut
    \caption{The 1-petal billiard and its unfolding to the 2-petal billiard}
\label{fig:WojtkowskiBilliardTable}
\end{figure}

Lemon billiards may have first appeared in \cite{Tomsovic} in the form of the intersection of two circles of equal radius. \cite{oaedhalb16} proved their hyperbolicity. We here consider \emph{asymmetric lemon billiards} (see \cref{fig:BilliardTable}), which consist of the intersection of two disks of which one contains the centers of both and were introduced by Chen--Mohr--Zhang--Zhang, who conjectured that they are ergodic when the radii differ greatly \cite[Observations, p.\ 3]{eceifcb13} (this was stated as a conjecture in the submitted manuscript, but the journal insisted on renaming it to ``observation''). The extreme cases of this are on one hand a circular billiard, which is not hyperbolic (indeed, integrable), and on the other hand a disk with a small part of the boundary flattened to a straight line (\cref{fig:WojtkowskiBilliardTable}). The latter is significantly more tractable than lemon billiards---it satisfies the Wojtkowski criteria---and we call it the 1-petal billiard \cite{MR357736,MR530154}. When unfolded along the line segment, it becomes the two-petal billiard shown in \cref{fig:WojtkowskiBilliardTable}, and any focusing circle arc of the billiard table is such that the disk in whose boundary it lies is contained in the billiard table. This produces a defocusing (or overfocusing) mechanism which traces back to Bunimovich \cite{MR357736}, Wojtkowski  \cite{MR848647},  Markarian \cite{MR954676}, Donnay \cite{MR1133266}, and again Bunimovich \cite{MR1179172}. By design, lemon billiards do not satisfy these criteria, so hyperbolicity must be established by arguments closely tied to the specifics of the billiard table. Consequently, the present work adapts techniques from the theory of dispersing billiards without being able to invoke ready-made theorems (other than the abundant general toolkit provided by \cite{cb}).

From this perspective, asymmetric lemon billiards can be viewed as a perturbation in which the straight boundary segment of the 1-petal billiard in \cref{fig:WojtkowskiBilliardTable} is replaced by    
an arc of a circle with large radius; these systems do not satisfy the Wojtkowski criteria and cannot be unfolded to do so in the way shown in \cref{fig:WojtkowskiBilliardTable} for the 1-petal billiard. Nonetheless, it was shown relatively recently, that they possess an invariant cone family \cite{oaedhalb16,hoalb}, which implies (nonuniform) hyperbolicity. This can be expressed in terms of invariant cone families but is defined as the property that the natural measure has nonzero Lyapunov exponents, i.e., pointwise eventual exponential expansion.

\cref{RMASUH} constitutes the next step, which is to establish \emph{uniform} expansion 
of vectors in cones of a suitable family under application of the return map to a suitable section. This is \emph{uniform} hyperbolicity: uniform exponential lower bounds for expansion a.e. 
We note that the tractability of the 1-petal billiard is being used in handling \emph{some} of  our arguments. At the same time, this is not (and likely cannot be) the overall strategy. 
 
Future work will show that the singularity curves are properly aligned, that the billiard map is ergodic (proving the Chen--Mohr--Zhang--Zhang conjecture \cite[Observations, p.\ 3]{eceifcb13}) and indeed Bernoulli, and we expect exponential mixing. While there are other sections for which one can prove uniform hyperbolicity of the return map, it is worth emphasizing that the choice of section provides a novel way of handling configurations of singularity curves that made them seem intractable before. 
\begin{remark}
Some comments about \cref{RMASUH}.

The ``next to last'' choice may seem odd but is an important device for handling the challenging configuration of some singularity curves (and hence invariant cones), which occurs near the left and right edges in \cref{fig:hatUplusUminus}. This choice keeps the section away from those. This issue is of particular importance in follow-on work to establish ergodicity; see below. 

Note also that it is essential here to study the return map to a section because the billiard map on the whole is not uniformly hyperbolic. Indeed, there are parabolic periodic orbits (those colliding only with $\Gamma_r$), though these only occur for small periods.

In the present context (planar billiards), strict invariance of a cone family implies expansion (and this is more generally true in symplectic contexts). However, ergodicity arguments require \emph{uniform} expansion, and this does not follow from eventual strict invariance of cones but only from some form of \emph{uniformly} strict invariance of cone families (that is, (1) does not imply (2) in \cref{RMASUH}); this would require a novel definition but be no easier to prove than uniform expansion itself. (The proofs would likely be equivalent to ours but more cumbersome.)
\end{remark}

One price we pay for the upgrade to uniform hyperbolicity is that $R_{\text{\upshape HF}}$ in \cref{RMASUH} has to be quite large. 
The main result of \cite{hoalb} is that for $\phistar\in(0,\frac{\pi}{2})$ and
\begin{equation}\label{eqJZHypCond0}
R>R_{\text{JZ}}(r,\phistar)\dfn r\cdot\max\big\{\frac{14.6}{\min{\{\phistar,\frac{\pi}{2}-\phistar\}}},\frac{147}{\sin^2{\phistar}},1773.7\big\},
\end{equation} the return map to the global section \(\hat M\) from \cref{def:SectionSetsDefs} preserves the cones defined by the positive quadrant in \cref{corollary:InvQuadrant}
\cite[Equations (3.1) and (3.2)]{hoalb}.\COMMENT{Refer to here when "1700" arises.\newline\wentao{Agree. Work In Progress}}

The proof of \cref{RMASUH} establishes that for a.e.\  $\phistar\in\big(0,\tan^{-1}(1/3)\big)$, there exists an
\begin{equation}\label{eqJZHypCond}
    R_{\text{\upshape HF}}(r,\phistar)
    \ge\frac{30000r}{\sin^2\phistar}\ge R_{\text{JZ}}(r,\phistar)\text{ such that }R>R_{\text{\upshape HF}}(r,\phistar)\Rightarrow\text{ hyperbolicity,}
\end{equation} 
that is, the return maps on certain global sections are uniformly hyperbolic.\COMMENT{It will be good to keep thinking about what references to \cref{eqJZHypCond} we need.\newline\wentao{Agree. Work in Progress}} Note that $\phistar\in\big(0,\tan^{-1}(1/3)\big)$ implies $\sin^2\phistar<1/10$, so our assumptions imply in particular that $R>300\,000r$.

Unlike in \cite{hoalb}, we have no closed form for $R_{\text{\upshape HF}}(r,\phistar)$ because its dependence on $\phistar$ is significantly more subtle and related to the need to exclude some values of $\phistar$. These are not numerous and can be described as follows: The chord $\overline{AB}$ is part of a periodic orbit in the 1-petal billiard bounded by $C_r$ that does not hit (the interior of) $\arc{AB}$, so it subtends an angle $2\pi/n$ of $C_r$ for some integer $n$. 
Note that the feasible $\phistar$ set is defined by the 1-petal billiard in \cref{def:ApettrajectoryExpansion} so it is independent of $R$. We provide an algorithm in \cref{def:AcomputationForRHF} to determine $R_{\text{\upshape HF}}$. This value is not meant to be optimal, nor does it seem urgent to find the optimal value.

We add a few remarks about invariant cone fields, beginning with a look at circular billiards, from which, after all, lemon billiards are built.

A billiard trajectory in a circle consists of secant lines, all of which have the same angle with the boundary; equivalently, the segments of an orbit are all tangent to the same concentric circle. (Either the radius of the circle or the constant collision angle provide a constant of motion, which makes this a completely integrable system in the sense of classical mechanics.) The billiard map accordingly preserves collision angles. The action on the phase space is easy to describe. The phase space is the product of a circle (the boundary) and an interval (of post-collision angles). In this cylinder, the billiard map rotates circles corresponding to a given angle; up to reparametrization, it does so in a linear way.\begin{figure}[ht]
    \centering
    \includegraphics[width=0.5\linewidth]{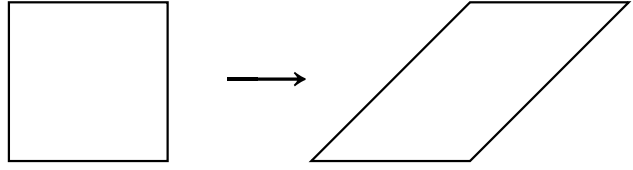}
    \caption{Circle biliard map}
    \label{fig:billiardincircle}
\end{figure}
Thus, the billiard map is a linear twist, and so is its differential. The positive quadrant is an invariant cone field; it is mapped into a proper subset of the positive quadrant by the billiard map. However, this cone field is not \emph{strictly} invariant in the needed sense: the image of a cone does not lie in the \emph{interior} of the first quadrant, no matter how many iterates of the billiard map are applied. This is related to the growth of positive vectors: under iteration, these become arbitrarily long but the growth is merely linear in time, whereas hyperbolic dynamics involves exponential growth.

\begin{figure}[hbt]
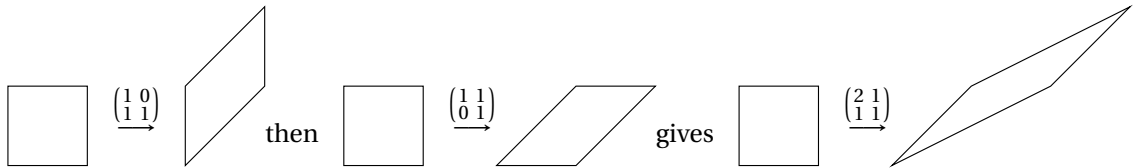

    \centering
\resizebox{.9\textwidth}{!}{\tikz{\draw(0,0)--(0,1)--(1,1)--(1,0)--cycle;}\quad\raise2.5ex\hbox{\(\overset{\left(\begin{smallmatrix}1&0\\1&1\end{smallmatrix}\right)}\longrightarrow\)}\quad\tikz{\draw(0,0)--(0,1)--(1,2)--(1,1)--cycle;}\raise2ex\hbox{then}\quad\tikz{\draw(0,0)--(0,1)--(1,1)--(1,0)--cycle;}\quad\raise2.5ex\hbox{\(\overset{\left(\begin{smallmatrix}1&1\\0&1\end{smallmatrix}\right)}\longrightarrow\)}\tikz{\draw(0,0)--(1,1)--(2,1)--(1,0)--cycle;}\raise2ex\hbox{gives}\quad\tikz{\draw(0,0)--(0,1)--(1,1)--(1,0)--cycle;}\quad\raise2.5ex\hbox{\(\overset{\left(\begin{smallmatrix}2&1\\1&1\end{smallmatrix}\right)}\longrightarrow\)}\tikz{\draw(0,0)--(1,1)--(3,2)--(2,1)--cycle;}}    \caption{Linking twist maps}
    \label{fig:linkedtwist}
\end{figure}

There is a way of seeing why a combination of circle arcs might produce exponential complexity: through the linking of twist maps. Preceding the action in \cref{fig:billiardincircle} by a vertical twist will move the bottom edge of the positive quadrant strictly inside it. This is illustrated in \cref{fig:linkedtwist}.
Notably, iteration of this linked twist map produces exponential growth of vectors in time.

At the most basic level, this describes what happens as a trajectory of a lemon billiard transitions between the two circle arcs. While the work here is nothing like as simple as this suggests, this consideration makes it natural to build a global section from transitions between the two circle arcs, and the reader is now prepared to find places where these transitions prove essential. It may be well to mention again that \cite{hoalb} produced (with great effort) a cone field for lemon billiards that is invariant and eventually strictly invariant: for almost every point, some iterate will send the positive quadrant into the interior of the positive quadrant at the image point. Without those authors blazing a trail, the present work would not have been undertaken. Its object is to strengthen this by obtaining immediate and uniform strict invariance---when attending to the return to a well-chosen subset (the section $\hat{M}$ in \eqref{def:Mhat}).
\begin{proof}[Proof of \cref{RMASUH}]
As is common in the literature, we prove expansion of the so-called p-metric (\cref{def:SectionSetsDefs}), and expansion with respect to the Euclidean metric then follows (\cref{Col:EuclideanExpansionMhat,Col:EuclideanExpansionMhat1,corollary:ExpansionBetweenMhatminusMhatplus,corollary:MroutreturnUniformExpansion}).

For the return map $\hat{F}$ to the section $\hat{M}$ in \eqref{def:Mhat} and $x\in\hat{M}$, the return orbit segment $x, \Fc(x),\cdots,\Fc^{\sigma(x)}(x)\in\hat{M}$ in \cref{def:MhatReturnOrbitSegment} with $x_1=(\Phi_1,\theta_1)\in M_R$ in \eqref{eqPtsFromeqMhatOrbSeg} and for all $dx\in C_x$, there are the following cases.
    
In case (a0) of \eqref{eqMainCases}, $\sin\theta_1>\sqrt{4r/R}$, which corresponds to ``large'' collision angles. Here, $\frac{\|D_x\hat{F}(dx)\|_{\p}}{\|dx\|_{\p}}>1+\lambda_c$ with constant $\lambda_c>0$ from \cref{TEcase_a0}. Thus we get expansion in this situation.

Not so in cases (a1)(b)(c) of \eqref{eqMainCases}, where we can only control how much contraction can occur. \cref{TEcase_a1,TEcase_b,TEcase_c} give  $\frac{\|D_x\hat{F}(dx)\|_{\p}}{\|dx\|_{\p}}>0.05$, but  \cref{def:ApettrajectoryExpansion} gives a uniform $\overline{l}>0$ (determined by the 1-petal billiard in \cref{lemma:PETExpExpansion}, thus independent of $R$) such that 
\[\frac{\|D_x\hat{F}^{\overline{l}+1}(dx)\|_{\p}}{\|dx\|_{\p}}>0.05\times22.5=1.125.\]

Let $\overline{\Lambda}=\min\{1+\lambda_c,1.125\}$, where $\lambda_c$ is from \cref{TEcase_a0}. \cref{def:ApettrajectoryExpansion,def:AcomputationForRHF}, and \eqref{eq:xorbitA0Condition}, ensure that among every  $\overline{l}+1$ consecutive returns to $\hat{M}$ in \eqref{eqMhatOrbSeg}, there is at most one that does not fall into case (a0) of \eqref{eqMainCases}. Therefore,\[\frac{\|D_x\hat{F}^{n}(dx)\|_{\p}}{\|dx\|_{\p}}>(\overline{\Lambda})^{\left\lfloor\frac{n}{\overline{l}+1}\right\rfloor}(0.05)\] for all nonsingular $x$ and $n\ge1$. 

Thus, $\frac{\|D_x\hat{F}^{n}(dx)\|_{\p}}{\|dx\|_{\p}}>(\overline{\Lambda})^{\frac{n}{\overline{l}+1}-1}(0.05)=c\Lambda^n$, where $c=(0.05)/\overline{\Lambda}$, $\Lambda=\overline{\Lambda}^{\frac{1}{1+\overline{l}}}$.
\end{proof}
Let us justify a claim made at the very beginning:
\begin{proposition}\label{WCriteriaNotConvex}
A piecewise $C^4$ planar billiard that satisfies the Wojtkowski criteria for hyperbolicity \cite[p.~144]{Tabachnikov} is not strictly convex.
\end{proposition}
\begin{proof}
Consider a strictly convex piecewise $C^4$ planar billiard. Each boundary arc is focusing by strict convexity. One part of the criterion is that at a common end-point of two such arcs the interior angle must exceed $\pi$, which is incompatible with convexity. Thus, there is only one arc, i.e., the billiard is $C^4$. Another part of the criterion is that for any focusing arc (hence in this context for the entire boundary) the radius of curvature as a function of arc length is strictly concave, i.e., has strictly monotone first derivative. This is incompatible with periodicity of the first derivative.
\end{proof}
\subsection*{Acknowledgments}
We owe much gratitude to Pengfei Zhang and Xin Jin. This work would not have begun without their seminal article \cite{hoalb}. Their support and encouragement were invaluable, notably conversations with Pengfei Zhang at various stages of this project. Moreover, Xin Jin served on the candidacy examination committee for Wentao Fan. Conversations with Marco Lenci, Domokos Sz\'asz, and Maciej Wojtkowski were also quite helpful in directing us to pertinent works which proved useful \cite{MR1636897,OKPPHB}.\COMMENT{Insert any other pertinent references here.\wentao{all set.}} Finally, the support and perspective of Hongkun Zhang and Alexey Korepanov provided a significant and invaluable boost as this article neared completion.

\subsection*{Structure of this paper}
The proof of \cref{RMASUH} above invoked results that will be established in subsequent parts of this work, all the way through \cref{sec:SACTWFELBWB}. These and subsequent arguments are organized as follows.

\cref{sec:ECE} states lower bounds for the differential of the return map in a case-by-case fashion (\cref{eqMainCases,TEcase_a0,TEcase_a1,TEcase_b,TEcase_c}) and establishes that those cases in which these lower bounds do not produce expansion are confined to a small portion of the section (\cref{contraction_region}). This is sufficient information to show strict invariance of the positive quadrant (\cref{corollary:InvQuadrant}; see also \cref{corollary:halfquadrantcone}) and strict expansion in it upon iteration (\cref{def:ApettrajectoryExpansion}). The proofs of these items in \cref{sec:ECE} are carried out in later sections as follows. 

\cref{subsec:Monotonicity} prepares for delicate estimates across various parameter ranges by establishing monotonicity results which reduce the needed estimates to considerations at the endpoints of the various ranges. This runs all the way to page~\pageref{sec:ProofTECase_a0_a1}. 


\cref{sec:ProofTECase_a0_a1} then establishes expansion in cases (a0) and (a1). These may be the easiest cases, but this takes up pages \pageref{subsec:AFELB}--\pageref{sec:ProofTECase_b}. 

\cref{sec:ProofTECase_b} covers Case (b) (\cref{TEcase_b,subsec:ProofTECase_b}). 

\cref{sec:ProofTECase_c} finally tackles the remaining case by proving \cref{TEcase_c}. 

One important deferred item remains for \cref{sec:SACTWFELBWB}, which is to control the dynamics in the small ``bad region’’ by approximating those trajectories by their counterpart for the 1-petal billiard (\cref{fig:WojtkowskiBilliardTable}). This necessitates the assumption that the larger of the two radii is so much larger than the other.

We repeat that establishing uniform hyperbolicity is the most laborious step towards ergodicity. Next steps (in forthcoming work) will follow the blueprint in \cite{cb}.
\tableofcontents\section{Preliminaries}\label{subsec:NotationsdefinitionsAndpropositions}
We begin with basic notations and observations.
\begin{notation}[Configuration parameters, phase space, singularities, and sections]\label{def:BasicNotations}\label{eqMroutMrin}\label{def:LemonTableConfiguration}\label{def:BasicNotations0}\label{def:StandardCoordinateTable}
     The (asymmetric) lemon billiard table is the intersection of two disks with radii $r$ and $R$, respectively,\footnote{The notations would be slightly lighter if we were to choose $r=1$ throughout, but we retain both radii to conform with the notation in other articles on asymmetric lemon billiards.} and centers $O_r$ and $O_R$ at distance $b$ as in \Cref{fig:BilliardTable}.
     \begin{itemize}
\item When we use cartesian coordinates, they will be chosen such that $O_r=(0,0)$ and $O_R=(0,b)$ (see, for example, \cref{fig:CaseAMonotone,fig:R1_R2_d_0,fig:lengthfraction,fig:twoTangentLines,fig:QonUpperHalfCircle}).
\item The pertinent notations are as in \Cref{fig:BilliardTable}: $A$, $B$ are the intersection points, $2\phistar=\measuredangle A O_r B$ and $2\Phistar=\measuredangle A O_R B$, so $R\sin\Phistar=r\sin\phistar$.
\item We denote by $\Gamma_r$ and $\Gamma_R$  the two smooth pieces of the table boundary with radii $r$, $R$, respectively, and by $C_r$ and $C_R$ the circles to which they belong.
\item  The phase space $M\dfn M_r\cup M_R$ consists of inward-pointing vectors on \(\Gamma_r\cup\Gamma_R\), where 
    $M_r=[\phistar,2\pi-\phistar]\times[0,\pi]$ and $M_R=[-\Phistar,\Phistar]\times[0,\pi]$ for the collisions on $\Gamma_r$, $\Gamma_R$, respectively. Here, the position coordinates are the angle parameters $\phi\in[\phistar,2\pi-\phistar]$ on $C_r$ and  $\Phi\in [-\Phistar,\Phistar]$ on $C_R$. The direction coordinates $\theta\in [0,\pi]$ are the collision angle with the tangent direction. We can also write $\theta=\pi/2-\varphi$, where $\varphi$ is the collision angle with the inward normal vector. Then $d\phi=\nicefrac{ds}{r}$, $d\Phi=\nicefrac{ds}{R}$, where $ds$ is the arc-length differential on $\Gamma_r$ and $\Gamma_R$, respectively. The tangent direction on $\Gamma_r$ and $\Gamma_R$ is given by the counter-clockwise orientation on those boundary.
    \item A piece of a billiard table boundary is a focusing or dispersive arc if it has negative or positive curvatures given the counterclockwise orientation of the arc and a flat arc if it has 0 curvature \cite[Section 2.1, Equation (2.4)]{cb}. The lemon billiard has two focusing boundary pieces.
    \item For $x\in M$, we denote by $p(x)$ its position on the lemon billiard table.
    \item Denote by $\Fc$ the billiard map which sends a vector in \(M\) to the outward vector after the next collision---provided that this next collision is not at a corner \cite[Equation (2.18)]{cb}. (This is the return map of the billiard flow to the boundary \cite[Sections 2.9--2.12]{cb}.
    \item The boundary of the phase space is $\partial M\dfn \Scal_0\dfn\partial M_r\cup\partial M_R$

    $\Scal_1\dfn\big\{x\text{ is an interior point of }M\text{ with }\Fc(x)\in \Scal_0\big\}\cup \Scal_0$, 
    
    $\Scal_{-1}\dfn\big\{x\text{ is an interior point of } M \text{ with }\Fc^{-1}(x)\in \Scal_0\big\}\cup \Scal_0$. 
    \item Recursively define $\Scal_k\dfn \Scal_{k-1}\cup\big\{x\bigm|x\in M\setminus \Scal_{k-1},\Fc(x)\in \Scal_{k-1}\big\}$, $\Scal_{-k}\dfn \Scal_{-k+1}\cup\big\{x\bigm|x\in M\setminus \Scal_{-k+1},\Fc(x)\in \Scal_{-k+1}\big\}$ for $k\ge 1$ \cite[Equation (2.27)]{cb}. 
    \item Singular points. The forward and backward singularity sets are $\Scal_{\infty}\dfn \bigcup_{k=0}^{\infty}\Scal_k$ and $\Scal_{-\infty}\dfn \bigcup_{k=0}^{\infty}\Scal_{-k}$ respectively \cite[Equation (2.28)]{cb}. Those singular points are either ``grazing'' vectors on $\Gamma_r\cup\Gamma_R$  (with $\theta=0\text{ or }\pi$) tangent to the billiard boundary (this does happen in the lemon billiard) or those vectors coming from a corner or going to a corner. The singularity sets have Lebesgue measure 0 \cite[Section 2.11 and equation (2.28)]{cb}, i.e., a.e. $x\in M$ is nonsingular.
    \item The billiard map $\Fc$ is a diffeomorphism from $M\setminus \Scal_{1}$ to $M\setminus \Scal_{-1}$ \cite[Theorem 2.33 and equation (2.26)]{cb}. 
    \item We consider the sections \begin{align*}
        \Mrout&\dfn\Int M_r\cap \Fc^{-1}(\Int M_R)\text{ and }\Mrin\dfn\Int M_r\cap \Fc(\Int M_R), \\\MRout&\dfn\Int M_R\cap \Fc^{-1}(\Int M_r)\text{ and }\MRin\dfn M_R\cap \Fc(\Int M_r)
    \end{align*}
shown in \cref{fig:Mr,fig:MR}. Here  ``interior'' merely omits corners. These sections consist of the last and first collisions with $\Gamma_r$ and $\Gamma_R$, respectively. $\Mrin,\Mrout$ are global sections of $M$. That is, $\bigcup_{n\ge0}\Fc^n(\Mrin)\aeq M\aeq\bigcup_{n\ge0}\Fc^n(\Mrout)$ since every $x=(\phi,\theta)\in M_r$ (except for finitely many $\theta$ such that $\theta/\pi$ is rational) has an $\Fc$-orbit across $\Mrout$ and $\Mrin$.

\end{itemize}
\end{notation}
\begin{remark}
We noted that the billiard map $\Fc$ is not always well defined on the boundary of the phase space, i.e., for the corners (and in other cases, grazing collisions) in the table. Since the angle at corners of the lemon billiard does not divide $\pi$, there are two ways to continuously extend the definition of the billiard map to the phase-space boundary by considering the collision as one with one or the other of the adjacent boundary pieces (see \cite[Section 2.8]{cb}, \cref{SBS1petalbilliardorbits}).
\end{remark}
\begin{figure}[ht]
\begin{minipage}{.6\textwidth}
\begin{center}
    \begin{tikzpicture}[xscale=.6,yscale=.6]
        \pgfmathsetmacro{\PHISTAR}{0.8}
        \pgfmathsetmacro{\WIDTH}{7}
        \pgfmathsetmacro{\LENGTH}{7}
        \tkzDefPoint(-\LENGTH,0.5*\WIDTH){LLM};
        \tkzDefPoint(-\LENGTH,0){LLL};
        \tkzDefPoint(-\LENGTH+\PHISTAR,0){LL};
        \tkzDefPoint(-\LENGTH+\PHISTAR,\WIDTH){LU};
        \tkzDefPoint(-\LENGTH,\WIDTH){LLU};
        \tkzDefPoint(\LENGTH-\PHISTAR,0){RL};
        \tkzDefPoint(\LENGTH-\PHISTAR,\WIDTH){RU};
        \tkzDefPoint(\LENGTH,0){RRL};
        \tkzDefPoint(\LENGTH,\WIDTH){RRU};
        \tkzDefPoint(-\LENGTH+\PHISTAR,\PHISTAR){X};
        \tkzDefPoint(-\LENGTH+\PHISTAR,\WIDTH-\PHISTAR){IX};
        \tkzDefPoint(\LENGTH-\PHISTAR,\PHISTAR){IY};
        \tkzDefPoint(\LENGTH-\PHISTAR,\WIDTH-\PHISTAR){Y};
        \tkzLabelPoint[below](LL){$\phistar$};
        \tkzLabelPoint[below](RL){$2\pi-\phistar\qquad$};
        \tkzLabelPoint[left](LLL){0};
        \tkzLabelPoint[left](LLM){$\frac{\pi}{2}$};
        \tkzLabelPoint[left](LLU){$\pi$};
        \tkzLabelPoint[below](LLL){0};
        \tkzLabelPoint[below](RRL){\ $2\pi$};
        \tkzLabelPoint[right](IY){$x_*$};
        \tkzLabelPoint[left](IX){$y_*$};
        \tkzLabelPoint[right](Y){$I x_*$};
        \tkzLabelPoint[left](X){$I y_*$};

        \draw [thick] (LL) --(LU);
        \draw [thick] (LU) --(RU);
        \draw [thick] (RL) --(RU);
        \draw [thick] (LL) --(RL);
        \draw [thin,dashed](LLL)--(LL);
        \draw [thin,dashed](LLU)--(LU);
        \draw [thin,dashed](LLU)--(LLL);
        \draw [thin,dashed](RRL)--(RL);
        \draw [thin,dashed](RRU)--(RU);
        \draw [thin,dashed](RRU)--(RRL);
        \fill [blue, opacity=10/30](LL) -- (Y) -- (RU) -- (X) -- cycle;
        \fill [red, opacity=10/30](LU) -- (IY) -- (RL) -- (IX) -- cycle;
        \begin{scope}
            \clip (LL) -- (X) -- (RU) -- cycle;
            \draw[fill=blue] circle[at=(X),radius=0.35];
        \end{scope}
        \begin{scope}
            \clip (RU) -- (Y) -- (LL) -- cycle;
            \draw[fill=blue] circle[at=(Y),radius=0.35];
        \end{scope}
        \begin{scope}
            \clip (RL) -- (IX) -- (LU) -- cycle;
            \draw[fill=red] circle[at=(IX),radius=0.35];
        \end{scope}
        \begin{scope}
            \clip (LU) -- (IY) -- (RL) -- cycle;
        \draw[fill=red] circle[at=(IY),radius=0.35];
        \end{scope}
    \end{tikzpicture}
    \end{center}
    \caption{The region $M_r$. The \textcolor{blue}{blue} parallelogram region is \Mrin, and the \textcolor{red}{red} one is \Mrout. The small sectors are `bad' regions (\cref{REMBadRegions}).} 
    \label{fig:Mr}    
\end{minipage}\hfill
\begin{minipage}{.35\textwidth}
\begin{center}
    \begin{tikzpicture}[scale=.28]
        \pgfmathsetmacro{\PHISTAR}{2.5};
        \pgfmathsetmacro{\WIDTH}{15};
        \pgfmathsetmacro{\LENGTH}{7.5};
        \tkzDefPoint(-\LENGTH,0.5*\WIDTH){LLM};
        \tkzDefPoint(-\LENGTH,0){LLL};
        \tkzDefPoint(-\PHISTAR,0){LL};
        \tkzDefPoint(-\PHISTAR,\WIDTH){LU};
        \tkzDefPoint(-\LENGTH,\WIDTH){LLU};
        \tkzDefPoint(\PHISTAR,0){RL};
        \tkzDefPoint(\PHISTAR,\WIDTH){RU};
        \tkzDefPoint(\LENGTH,0){RRL};
        \tkzDefPoint(\LENGTH,\WIDTH){RRU};
        \tkzDefPoint(-\PHISTAR,\PHISTAR){X};
        \tkzDefPoint(-\PHISTAR,\WIDTH-\PHISTAR){IX};
        \tkzDefPoint(\PHISTAR,\WIDTH-\PHISTAR){Y};
        \tkzDefPoint(\PHISTAR,\PHISTAR){IY};
        \draw [thick] (LL) --(LU);
        \draw [thick] (LU) --(RU);
        \draw [thick] (RL) --(RU);
        \draw [thick] (LL) --(RL);
        \draw [thin,dashed](LLL)--(LL);
        \draw [thin,dashed](LLU)--(LU);
        \draw [thin,dashed](LLU)--(LLL);
        \draw [thin,dashed](RRL)--(RL);
        \draw [thin,dashed](RRU)--(RU);
        \draw [thin,dashed](RRU)--(RRL);
         \fill [blue, opacity=10/30](LL) -- (IY) -- (RU) -- (IX) -- cycle;
        \fill [red, opacity=10/30](LU) -- (X) -- (RL) -- (Y) -- cycle;
        \tkzLabelPoint[below](LL){$-\Phistar$};
        \tkzLabelPoint[left](LLL){$0$};
        \tkzLabelPoint[below](RL){$+\Phistar$};
        \tkzLabelPoint[left](LLU){$\pi$};
        \tkzLabelPoint[left](LLM){${\pi}/{2}$};
    \end{tikzpicture}
\end{center}
  \captionsetup{width=\linewidth}
    \caption{The region $M_R$. The \textcolor{red}{red} region is \MRout, the \textcolor{blue}{blue} region is \MRin.}
    \label{fig:MR}    
\end{minipage}
\end{figure}

\begin{remark}[Symmetries]\label{DefSymmetries}
    \cref{fig:Mr,fig:MR,fig:MrN} exhibit some symmetry, and there are indeed two symmetries in the lemon billiard system that give rise to this. One of these is \emph{reversibility} of the billiard \emph{flow} (which holds for all billiard systems): changing the sign of the velocity vector in any billiard reverses the motion and the collision sequence. Equivalently, the flip map $I$ from \cref{def:SectionSetsDefs} conjugates the billiard \emph{map} to its inverse: $I\Fc=\Fc^{-1}I$.
    \begin{equation}\label{eq:InversionSymmetryDef}
        \begin{aligned}
        I\colon M_r\sqcup  M_R&\rightarrow M_r\sqcup M_R,\\ 
        (\phi,\theta)&\mapsto(\phi,\pi-\theta)\in M_r,\text{ if }(\phi,\theta)\in M_r.\\
        (\Phi,\theta)&\mapsto(\Phi,\pi-\theta)\in M_R,\text{ if }(\Phi,\theta)\in M_R.
        \end{aligned}
    \end{equation}
(Throughout, ``$\sqcup$'' denotes a union of sets known to be disjoint.)    In the phase space, this symmetry is manifested as an up-down symmetry, which interchanges blue and red in \cref{fig:Mr,fig:MR}.
    
The other symmetry is the left-right symmetry of the phase space, which projects to the mirror 
symmetry of the lemon billiard table with respect to the line through the centers of the defining circles. It is defined as follows.
    \begin{equation}\label{eq:FoldSymmetryDef}
    \begin{aligned}
        J\colon M_r\sqcup M_R&\rightarrow M_r\sqcup M_R,\\ 
        (\phi,\theta)&\mapsto(2\pi-\phi,\theta)\in M_r,\text{ if }(\phi,\theta)\in M_r.\\
        (\Phi,\theta)&\mapsto(-\Phi,\theta)\in M_R,\text{ if }(\Phi,\theta)\in M_R.
    \end{aligned}
    \end{equation}
This also interchanges blue and red in \cref{fig:Mr,fig:MR} and conjugates the billiard map to its inverse: $J\Fc=\Fc^{-1}J$. $I$, $J$ are symmetries between $\Mrin$ and $\Mrout$, as well as between $\MRin$ and $\MRout$. Moreover, $I \circ J=J\circ I$ is the lift of the mirror symmetry of the billiard table; it reflects the position and the collision direction. When restricted to $M_r$, it is the symmetry with respect to its center $(\pi,\pi/2)$.  When restricted to $M_R$, it is the symmetry with respect to its center $(0,\pi/2)$. Also, $I \circ J\circ\Fc=\Fc\circ I\circ J$, i.e., the lift of the mirror symmetry of the billiard (unsurprisingly) commutes with $\Fc$.
\end{remark}

\begin{notation}[p-metric, grading of sections, ``bad'' regions]\label{def:SectionSetsDefs}\hfill
\begin{itemize}
\item Define the \emph{p-metric} by $\|dx\|_\p=\sin{\theta}|ds|$ for $dx\in T_x(M)$, where $ds$ is the differential of the arc-length parameter \cite[page 58]{cb}.\footnote{Here ``p'' stands for ``\textbf perpendicular,'' per Leonid Bunimovich---and this is a \textbf pseudometric.}
\item For $(\phi,\theta)=x\in M$, define the flip (or reversal) operation $I x=(\phi,\pi-\theta)$. This interchanges pre- and post-collision directions. In particular, for the corner points $y_*=(\phistar,\pi-\phistar)$ and $x_*=(2\pi-\phistar,\phistar)$ in \cref{fig:MrN}, we have $I y_*=(\phistar,\phistar)$ and $I x_*=(2\pi-\phistar,\pi-\phistar)$. (See \cref{DefSymmetries} below.) $y_*$ and $x_*$ in \cref{fig:Mr,fig:MrN} correspond to the two trajectories $\overrightarrow{AB}$ and $\overrightarrow{BA}$ between the corners $A$, $B$ of the billiard table in \cref{fig:BilliardTable}.
\item$B(Ix_*,\delta_*)$, $B(Iy_*,\delta_*)$, $B(x_*,\delta_*)$, $B(y_*,\delta_*)$ are the neighborhoods of radius $\delta_*\dfn4.5\sqrt{r/R}$

\COMMENT{\color{red}Is $\delta_*$ ideal notation? Is 17 a reasonable choice? In the proof of \cref{Prop:Udisjoint}, we could replace it with 89.\wentao{We may discuss. This number has a theoretic smallest asymptotic value.}}
\item ``Bad'' regions as in \cref{fig:MrN}:
\[\Nin\dfn\big(B(Ix_*,\delta_*)\cup B(Iy_*,\delta_*)\big)\cap \Mrin\text{
and }\Nout\dfn\big(B(x_*,\delta_*)\cup B(y_*,\delta_*)\big)\cap \Mrout.\]
\item For $x\in \Mrin$, let 
\[
m(x)\dfn\max\big\{m\ge 0\mid p(\Fc^i(x))\in\operatorname{int}\Gamma_r\text{ for }0\le i\le m\big\} 
\]
be the number of subsequent collisions with $\Gamma_r$.
\item For $x\in \MRin$, let 
\[
m(x)\dfn\max\big\{m\ge 0\mid p(\Fc^i(x))\in\operatorname{int}\Gamma_R\text{ for }0\le i\le m\big\} 
\]
be the number of subsequent collisions with $\Gamma_R$.
\item Let $\Min_{r,n}\dfn\big\{x\in \Mrin\big| m(x)=n\big\}$ and $\Mout_{r,n}\dfn I(\Min_{r,n})$ (see \Cref{fig:MrN}).\label{def-Mri}
We note that
\begin{equation}\label{EQMrRinoutgrading}\begin{aligned}\Mrin&=(\bigsqcup_{n=0}^{\infty}\Min_{r,n})\cup\{\text{singular points in }\Mrin\}\aeq\bigsqcup_{n=0}^{\infty}\Min_{r,n},\\ 
\Mrout&=(\bigsqcup_{n=0}^{\infty}\Mout_{r,n})\cup\{\text{singular points in }\Mrout\}\aeq\bigsqcup_{n=0}^{\infty}\Mout_{r,n},\\ 
\MRin&=(\bigsqcup_{n=0}^{\infty}\Min_{R,n})\cup\{\text{singular points in }\MRin\}\aeq\bigsqcup_{n=0}^{\infty}\Min_{R,n},\\
\MRout&=(\bigsqcup_{n=0}^{\infty}\Mout_{R,n})\cup\{\text{singular points in }\Mrout\}\aeq\bigsqcup_{n=0}^{\infty}\Mout_{R,n}.\end{aligned}\end{equation}
\item The domain \begin{equation}\label{def:Mhat}\hat{M}_{-}\dfn\hat M\dfn\big(\bigsqcup_{n\ge1}\Fc^{n-1}\Min_{r,n}\big)\sqcup\Min_{r,0}\aeq\Fc^{-1}\big(\Mrout\smallsetminus{\Mrin}\big)\sqcup\big(\Mrout\cap\Mrin\big)\end{equation} of next to last or sole collisions with $\Gamma_r$ is shown in \cref{fig:MrN,fig:hatUplusUminus},\COMMENT{Might it be possible to have a picture of $\hat M$ in the \emph{\Large\color{red} configuration space?\wentao{still think... seems not needed. any reference is needed to this pic?}}} and $\hat{F}$ is the return map to $\hat{M}$.\label{def:Mhatdefinition} This is a global section, i.e., the orbit of almost every point meets it in both past and future.
\item 
As shown in \cref{fig:hatUplusUminus},
\begin{equation}\label{def:symetricsectionMhatPlusMinus}
I\hat M=J\hat M=\hat{M}_{+}\dfn\hat{M}_1\dfn\big(\bigsqcup_{n\ge1}\Fc\Min_{r,n}\big)\sqcup\Min_{r,0}\aeq(\Mrin\cap\Mrout)\sqcup\Fc(\Mrin\smallsetminus\Mrout).
\end{equation}
\end{itemize}
\end{notation}
\begin{remark}\label{remark:M1hatMhattwoglobalsection}
    $\hat{M}_{\pm}$ are global sections, i.e., $M\aeq\bigcup_{k\in\mathbb{Z}}\Fc^{k}(\hat{M}_{\pm})$.
\end{remark}
\begin{figure}[ht]
\begin{center}
    \begin{tikzpicture}[xscale=1.0,yscale=1.0]
        \pgfmathsetmacro{\PHISTAR}{0.68}
        \pgfmathsetmacro{\WIDTH}{7}
        \pgfmathsetmacro{\LENGTH}{7}
        \tkzDefPoint(-\LENGTH,0.5*\WIDTH){LLM};
        \tkzDefPoint(-\LENGTH,0){LLL};
        \tkzDefPoint(-\LENGTH+\PHISTAR,0){LL};
        \tkzDefPoint(-\LENGTH+\PHISTAR,\WIDTH){LU};
        \tkzDefPoint(-\LENGTH,\WIDTH){LLU};
        \tkzDefPoint(\LENGTH-\PHISTAR,0){RL};
        \tkzDefPoint(\LENGTH-\PHISTAR,\WIDTH){RU};
        \tkzDefPoint(\LENGTH,0){RRL};
        \tkzDefPoint(\LENGTH,\WIDTH){RRU};
        \tkzDefPoint(-\LENGTH+\PHISTAR,\PHISTAR){X};
        \tkzDefPoint(-\LENGTH+\PHISTAR,\WIDTH-\PHISTAR){IX};
        \tkzDefPoint(\LENGTH-\PHISTAR,\PHISTAR){IY};
        \tkzDefPoint(\LENGTH-\PHISTAR,\WIDTH-\PHISTAR){Y};
        \tkzDefPoint(0,0.5*\WIDTH+0.5*\PHISTAR){UMrin01};
        \tkzDefPoint(\PHISTAR,0.5*\WIDTH){RightMrin01};
        \tkzDefPoint(0,0.5*\WIDTH-0.5*\PHISTAR){LMrin01};
        \tkzDefPoint(-1.0*\PHISTAR,0.5*\WIDTH){LeftMrin01};
        \fill [red, opacity=10/30](UMrin01) -- (RightMrin01) -- (LMrin01) -- (LeftMrin01) -- cycle;
        \fill [blue, opacity=10/30](UMrin01) -- (RightMrin01) -- (LMrin01) -- (LeftMrin01) -- cycle;
        \tkzDrawPoints(UMrin01,RightMrin01,LMrin01,LeftMrin01);
        \tkzLabelPoint[left](LeftMrin01){K1};
        \tkzLabelPoint[below](LMrin01){K2};
        \tkzLabelPoint[right](RightMrin01){K3};
        \tkzLabelPoint[above](UMrin01){K4};
 
        \tkzDefPoint(\LENGTH-\PHISTAR,0.5*\WIDTH){RUhatUplus};
        \tkzDefPoint(\LENGTH-\PHISTAR,0.5*\WIDTH-0.5*\PHISTAR){RLhatUplus};
        \tkzDefPoint(-1.0*\LENGTH+3.0*\PHISTAR,\PHISTAR){LUhatUplus};
        \fill [blue, opacity=10/30](RUhatUplus) -- (RLhatUplus) -- (LL) -- (LUhatUplus) -- cycle;
        \begin{scope}
            \clip (LL) -- (LUhatUplus) -- (RUhatUplus) -- cycle;
            \draw[fill=blue] circle[at=(LUhatUplus),radius=0.15];
        \end{scope}
        \node[] at ([shift={(0.24,-0.33)}]LUhatUplus) {\small $\Fc(N^{\text{\upshape in}})$};
        \node[] at ([shift={(-0.03,0.2)}]LUhatUplus) {\tiny $LF=(3\phistar,\phistar)$};

        \tkzDefPoint(-\LENGTH+\PHISTAR,0.5*\WIDTH){LLhatUminus};
        \tkzDefPoint(-\LENGTH+\PHISTAR,0.5*\WIDTH+0.5*\PHISTAR){LUhatUminus};
        \tkzDefPoint(1.0*\LENGTH-3.0*\PHISTAR,\WIDTH-\PHISTAR){RLhatUminus};
        \fill [blue, opacity=10/30](LLhatUminus) -- (LUhatUminus) -- (RU) -- (RLhatUminus) -- cycle;
        \begin{scope}
            \clip (RU) -- (RLhatUminus) -- (LLhatUminus) -- cycle;
            \draw[fill=blue] circle[at=(RLhatUminus),radius=0.15];
        \end{scope}
         \begin{scope}
            \clip (RU) -- (RLhatUminus) -- (LLhatUminus) -- cycle;
            \draw[fill=blue] circle[at=(RLhatUminus),radius=0.15];
        \end{scope}
         \node[] at ([shift={(0.1,-0.13)}]RLhatUminus) {\small $\Fc(N^{\text{\upshape in}})$};

        \tkzLabelPoint[below](LL){$LL=(\phistar,0)$};
        \tkzLabelPoint[below]([shift={(-0.5,0)}]RL){$2\pi-\phistar$};
        \tkzLabelPoint[left](LLL){0};
        \tkzLabelPoint[left](LLM){${\pi}/{2}$};
        \tkzLabelPoint[left](LLU){$\pi$};
        \tkzLabelPoint[below](RRL){$2\pi$};
        \tkzLabelPoint[right](IY){$x_*$};
        \tkzLabelPoint[left](IX){$y_*$};
        \tkzLabelPoint[right](Y){$I x_*$};
        \tkzLabelPoint[left](X){$IY=I y_*$};
        
        \draw [thick] (LL) --(LU);
        \draw [thick] (LU) --(RU);
        \draw [thick] (RL) --(RU);
        \draw [thick] (LL) --(RL);
        \draw [thin,dashed](LLL)--(LL);
        \draw [thin,dashed](LLU)--(LU);
        \draw [thin,dashed](LLU)--(LLL);
        \draw [thin,dashed](RRL)--(RL);
        \draw [thin,dashed](RRU)--(RU);
        \draw [thin,dashed](RRU)--(RRL);
        \draw [thin,dashed](LL)--(Y);
        \draw [thin,dashed](RU)--(X);
        \draw [thin,dashed](LU) -- (IY);
        \draw [thin,dashed](RL) -- (IX);
        \tkzDefPoint(-\LENGTH+\PHISTAR,\WIDTH/2-\PHISTAR/4-\PHISTAR/4){M1LL};
        \tkzDefPoint(-\LENGTH+\PHISTAR,\WIDTH/2-\PHISTAR/4+\PHISTAR/4){M1LU};
        \tkzDefLine[parallel=through M1LU](M1LL,RL)\tkzGetPoint{M11RU};
        \path[name path=lineM1U] (M11RU) -- (M1LU);
        \path[name path=lineM0L] (IX) -- (RL);
        \path[name intersections={of=lineM1U and lineM0L,by={M11RU}}];
        \fill [red, opacity=10/30](RL) -- (M1LL) -- (M1LU) -- (M11RU) -- cycle;
        \begin{scope}
            \clip (RL) -- (M11RU) -- (M1LU) -- cycle;
            \draw[fill=red] circle[at=(M11RU),radius=0.15];
        \end{scope}
        \node[] at ([shift={(-0.14,0.18)}]M11RU) {\small $\Fc^{-1}(N^{\text{\upshape out}})$};
        \path[name path=RBoundry](RL) -- (RU);
        \tkzDefLine[parallel=through LU](M1LL,RL)\tkzGetPoint{UM11RU};
        \path[name path=lineUM1U](LU)--(UM11RU);
        \path[name intersections={of=lineUM1U and RBoundry,by={UM1RU}}];
        \tkzDefShiftPoint[UM1RU](0,-\PHISTAR/2){UM1RL};
        \tkzDefLine[parallel=through UM1RL](UM1RU,LU)\tkzGetPoint{UM11LL};
        \path[name path=lineLUIY](LU) -- (IY);
        \path[name path=UM11LLUM1RL](UM11LL)--(UM1RL);
        \path[name intersections={of=lineLUIY and UM11LLUM1RL,by={UM1LL}}];
        \fill [red, opacity=10/30](LU) -- (UM1RU) -- (UM1RL) -- (UM1LL) -- cycle;
        \begin{scope}
            \clip (LU) -- (UM1LL) -- (UM1RL) -- cycle;
            \draw[fill=red] circle[at=(UM1LL),radius=0.15];
        \end{scope}
        \node[] at ([shift={(-0.12,-0.18)}]UM1LL) {\small $\Fc^{-1}(N^{\text{\upshape out}})$};
        \tkzDefPoint(-\LENGTH+\PHISTAR,\WIDTH/3){L_LUM2};
        \tkzDefShiftPoint[L_LUM2](0,-\PHISTAR/3){L_LLM2};
        \tkzDefLine[parallel=through L_LUM2](L_LLM2,RL)\tkzGetPoint{Far_L_RUM2};
        \path[name path=line_L_LUM2_Far_L_RUM2] (L_LUM2)--(Far_L_RUM2);
        \path[name path=line_RL_M1LL] (RL)--(M1LL);
        \path[name intersections={of=line_L_LUM2_Far_L_RUM2 and line_RL_M1LL,by={L_RUM2}}];
        \tkzDefPoint(\LENGTH-\PHISTAR,2*\WIDTH/3){R_RLM2};
        \tkzDefShiftPoint[R_RLM2](0,\PHISTAR/3){R_RUM2};
        \tkzDefLine[parallel=through R_RLM2](R_RUM2,LU)\tkzGetPoint{Far_R_LLM2};
        \path[name path=line_R_RLM2_Far_R_LLM2] (R_RLM2)--(Far_R_LLM2);
        \path[name path=line_LU_UM1RU] (LU)--(UM1RU);
        \path[name intersections={of=line_R_RLM2_Far_R_LLM2 and line_LU_UM1RU,by={R_LLM2}}];

        \tkzDefPoint(\LENGTH-\PHISTAR,1*\WIDTH/3){MR_RLM2};
        \tkzDefShiftPoint[MR_RLM2](0,\PHISTAR/3){MR_RUM2};
        \path[name path=line_IY_LU] (IY)--(LU);
        \tkzDefLine[parallel=through MR_RLM2](R_RLM2,R_LLM2)\tkzGetPoint{far_MR_LLM2};
        \tkzDefLine[parallel=through MR_RUM2](R_RLM2,R_LLM2)\tkzGetPoint{far_MR_LUM2};
        \path[name path=line_MR_RLM2_far_MR_LLM2] (far_MR_LLM2)--(MR_RLM2);
        \path[name path=line_MR_RUM2_far_MR_LUM2] (far_MR_LUM2)--(MR_RUM2);
        \path[name path=line_LU_IY] (LU)--(IY);
        \path[name intersections={of=line_LU_IY and line_MR_RLM2_far_MR_LLM2,by={MR_LLM2}}];
        \path[name intersections={of=line_LU_IY and line_MR_RUM2_far_MR_LUM2,by={MR_LUM2}}];

        \tkzDefPoint(-\LENGTH+\PHISTAR,2*\WIDTH/3){ML_LUM2};
        \tkzDefShiftPoint[ML_LUM2](0,-\PHISTAR/3){ML_LLM2};
        \tkzDefLine[parallel=through ML_LUM2](R_LLM2,R_RLM2)\tkzGetPoint{far_ML_RUM2};
        \tkzDefLine[parallel=through ML_LLM2](R_LLM2,R_RLM2)\tkzGetPoint{far_ML_RLM2};
        \path[name path=line_ML_LUM2_far_ML_RUM2] (ML_LUM2)--(far_ML_RUM2);
        \path[name path=line_ML_LLM2_far_ML_RLM2] (ML_LLM2)--(far_ML_RLM2);
        \path[name path=line_RL_IX] (RL)--(IX);
        \path[name intersections={of=line_RL_IX and line_ML_LUM2_far_ML_RUM2,by={ML_RUM2}}];
        \path[name intersections={of=line_RL_IX and line_ML_LLM2_far_ML_RLM2,by={ML_RLM2}}];
    \begin{scope}
      \clip (LL) -- (X) -- (RU) -- cycle;
      \draw[fill=green] circle[at=(X),radius=0.2];
    \end{scope}
    \begin{scope}
      \clip (RU) -- (Y) -- (LL) -- cycle;
      \draw[fill=green] circle[at=(Y),radius=0.2];
    \end{scope}
    \begin{scope}
      \clip (RL) -- (IX) -- (LU) -- cycle;
      \draw[fill=yellow] circle[at=(IX),radius=0.2];
    \end{scope}
    \begin{scope}
      \clip (LU) -- (IY) -- (RL) -- cycle;
      \draw[fill=yellow] circle[at=(IY),radius=0.2];
    \end{scope}
    
    \node[] at (0,0.5*\WIDTH)  {$\color{white}M^{\text{\upshape in}}_{r,0}$};
    \node[] at (0,0.25*\WIDTH)  {\tiny $\Fc(M^{\text{\upshape in}}_{r,2})$};
    \node[] at (0,0.75*\WIDTH)  {\tiny $\Fc(M^{\text{\upshape in}}_{r,2})$};
    
    \path[name path=line_RU_X] (RU)--(X);
    \path[name intersections={of=line_RU_X and lineM1U,by={Mr1Label}}];

    \path[name intersections={of=line_RU_X and line_L_LUM2_Far_L_RUM2,by={Mr2Label}}];
    \node[] at ([shift={(0.2,-0.05)}]X) {\tiny $N^{\text{\upshape in}}$};
    \node[] at ([shift={(-0.14,0.13)}]Y) {\tiny $N^{\text{\upshape in}}$};
    \node[] at ([shift={(0.25,0.1)}]IX) {\tiny $N^{\text{\upshape out}}$};
    \node[] at ([shift={(-0.08,-0.1)}]IY) {\tiny $N^{\text{\upshape out}}$};
    \tkzDrawPoints(LL,X,LUhatUplus);
    
    \end{tikzpicture}
    \caption{$\hat{M}=\hat{M}_{-}=(\Mrout\cap\Mrin)\sqcup\Fc^{-1}(\Mrout\smallsetminus\Mrin)$ is the union of the center diamond $M^{\text{\upshape in}}_{r,0}=\Mrout\cap\Mrin$ and the strips with negative slope.
    $\hat{M}_1=\hat{M}_{+}=(\Mrout\cap\Mrin)\sqcup\Fc(\Mrin\smallsetminus\Mrout)$ is the mirror image. Also shown: $\Fc(\Nin)$ and $\Fc^{-1}(\Nout)$.}\label{fig:hatUplusUminus}.
\end{center}
\end{figure}

\cref{EQMrRinoutgrading} implies that a nonsingular $x\in\hat{M}$ has a return orbit segment as follows, with $x_0\in\Mrout$, $x_1\in\MRin$, and $x_2\in\Mrin$.
\begin{figure}[ht]
\begin{center}
    \begin{tikzpicture}[xscale=1.0,yscale=1.0]
        \pgfmathsetmacro{\PHISTAR}{0.68}
        \pgfmathsetmacro{\WIDTH}{7}
        \pgfmathsetmacro{\LENGTH}{7}
        \tkzDefPoint(-\LENGTH,0.5*\WIDTH){LLM};
        \tkzDefPoint(-\LENGTH,0){LLL};
        \tkzDefPoint(-\LENGTH+\PHISTAR,0){LL};
        \tkzDefPoint(-\LENGTH+\PHISTAR,\WIDTH){LU};
        \tkzDefPoint(-\LENGTH,\WIDTH){LLU};
        \tkzDefPoint(\LENGTH-\PHISTAR,0){RL};
        \tkzDefPoint(\LENGTH-\PHISTAR,\WIDTH){RU};
        \tkzDefPoint(\LENGTH,0){RRL};
        \tkzDefPoint(\LENGTH,\WIDTH){RRU};
        \tkzDefPoint(-\LENGTH+\PHISTAR,\PHISTAR){X};
        \tkzDefPoint(-\LENGTH+\PHISTAR,\WIDTH-\PHISTAR){IX};
        \tkzDefPoint(\LENGTH-\PHISTAR,\PHISTAR){IY};
        \tkzDefPoint(\LENGTH-\PHISTAR,\WIDTH-\PHISTAR){Y};
        \tkzLabelPoint[below](LL){$\phistar$};
        \tkzLabelPoint[below]([shift={(-0.5,0)}]RL){$2\pi-\phistar$};
        \tkzLabelPoint[left](LLL){0};
        \tkzLabelPoint[left](LLM){$\frac{\pi}{2}$};
        \tkzLabelPoint[left](LLU){$\pi$};
        \tkzLabelPoint[below](LLL){0};
        \tkzLabelPoint[below](RRL){$2\pi$};
        \tkzLabelPoint[right](IY){$x_*$};
        \tkzLabelPoint[left](IX){$y_*$};
        \tkzLabelPoint[right](Y){$I x_*$};
        \tkzLabelPoint[left](X){$I y_*$};
        
        \draw [thick] (LL) --(LU);
        \draw [thick] (LU) --(RU);
        \draw [thick] (RL) --(RU);
        \draw [thick] (LL) --(RL);
        \draw [thin,dashed](LLL)--(LL);
        \draw [thin,dashed](LLU)--(LU);
        \draw [thin,dashed](LLU)--(LLL);
        \draw [thin,dashed](RRL)--(RL);
        \draw [thin,dashed](RRU)--(RU);
        \draw [thin,dashed](RRU)--(RRL);
        \fill [blue, opacity=17/30](LL) -- (Y) -- (RU) -- (X) -- cycle;
        \fill [red, opacity=17/30](LU) -- (IY) -- (RL) -- (IX) -- cycle;
        \tkzDefPoint(-\LENGTH+\PHISTAR,\WIDTH/2-\PHISTAR/4-\PHISTAR/4){M1LL};
        \tkzDefPoint(-\LENGTH+\PHISTAR,\WIDTH/2-\PHISTAR/4+\PHISTAR/4){M1LU};
        \tkzDefLine[parallel=through M1LU](M1LL,RL)\tkzGetPoint{M11RU};
        \path[name path=lineM1U] (M11RU) -- (M1LU);
        \path[name path=lineM0L] (IX) -- (RL);
        \path[name intersections={of=lineM1U and lineM0L,by={M11RU}}];
        \fill [red, opacity=10/30](RL) -- (M1LL) -- (M1LU) -- (M11RU) -- cycle;
        \path[name path=RBoundry](RL) -- (RU);
        \tkzDefLine[parallel=through LU](M1LL,RL)\tkzGetPoint{UM11RU};
        \path[name path=lineUM1U](LU)--(UM11RU);
        \path[name intersections={of=lineUM1U and RBoundry,by={UM1RU}}];
        \tkzDefShiftPoint[UM1RU](0,-\PHISTAR/2){UM1RL};
        \tkzDefLine[parallel=through UM1RL](UM1RU,LU)\tkzGetPoint{UM11LL};
        \path[name path=lineLUIY](LU) -- (IY);
        \path[name path=UM11LLUM1RL](UM11LL)--(UM1RL);
        \path[name intersections={of=lineLUIY and UM11LLUM1RL,by={UM1LL}}];
        \fill [red, opacity=10/30](LU) -- (UM1RU) -- (UM1RL) -- (UM1LL) -- cycle;
        \tkzDefPoint(-\LENGTH+\PHISTAR,\WIDTH/3){L_LUM2};
        \tkzDefShiftPoint[L_LUM2](0,-\PHISTAR/3){L_LLM2};
        \tkzDefLine[parallel=through L_LUM2](L_LLM2,RL)\tkzGetPoint{Far_L_RUM2};
        \path[name path=line_L_LUM2_Far_L_RUM2] (L_LUM2)--(Far_L_RUM2);
        \path[name path=line_RL_M1LL] (RL)--(M1LL);
        \path[name intersections={of=line_L_LUM2_Far_L_RUM2 and line_RL_M1LL,by={L_RUM2}}];
        \fill [green, opacity=15/30](RL)--(L_LLM2)--(L_LUM2)--(L_RUM2) -- cycle;
        \tkzDefPoint(\LENGTH-\PHISTAR,2*\WIDTH/3){R_RLM2};
        \tkzDefShiftPoint[R_RLM2](0,\PHISTAR/3){R_RUM2};
        \tkzDefLine[parallel=through R_RLM2](R_RUM2,LU)\tkzGetPoint{Far_R_LLM2};
        \path[name path=line_R_RLM2_Far_R_LLM2] (R_RLM2)--(Far_R_LLM2);
        \path[name path=line_LU_UM1RU] (LU)--(UM1RU);
        \path[name intersections={of=line_R_RLM2_Far_R_LLM2 and line_LU_UM1RU,by={R_LLM2}}];
        \fill [green, opacity=15/30](LU)--(R_RUM2)--(R_RLM2)--(R_LLM2)-- cycle;

        \tkzDefPoint(\LENGTH-\PHISTAR,1*\WIDTH/3){MR_RLM2};
        \tkzDefShiftPoint[MR_RLM2](0,\PHISTAR/3){MR_RUM2};
        \path[name path=line_IY_LU] (IY)--(LU);
        \tkzDefLine[parallel=through MR_RLM2](R_RLM2,R_LLM2)\tkzGetPoint{far_MR_LLM2};
        \tkzDefLine[parallel=through MR_RUM2](R_RLM2,R_LLM2)\tkzGetPoint{far_MR_LUM2};
        \path[name path=line_MR_RLM2_far_MR_LLM2] (far_MR_LLM2)--(MR_RLM2);
        \path[name path=line_MR_RUM2_far_MR_LUM2] (far_MR_LUM2)--(MR_RUM2);
        \path[name path=line_LU_IY] (LU)--(IY);
        \path[name intersections={of=line_LU_IY and line_MR_RLM2_far_MR_LLM2,by={MR_LLM2}}];
        \path[name intersections={of=line_LU_IY and line_MR_RUM2_far_MR_LUM2,by={MR_LUM2}}];
        \fill [green, opacity=15/30](MR_RUM2)--(MR_RLM2)--(MR_LLM2)--(MR_LUM2)-- cycle;

        \tkzDefPoint(-\LENGTH+\PHISTAR,2*\WIDTH/3){ML_LUM2};
        \tkzDefShiftPoint[ML_LUM2](0,-\PHISTAR/3){ML_LLM2};
        \tkzDefLine[parallel=through ML_LUM2](R_LLM2,R_RLM2)\tkzGetPoint{far_ML_RUM2};
        \tkzDefLine[parallel=through ML_LLM2](R_LLM2,R_RLM2)\tkzGetPoint{far_ML_RLM2};
        \path[name path=line_ML_LUM2_far_ML_RUM2] (ML_LUM2)--(far_ML_RUM2);
        \path[name path=line_ML_LLM2_far_ML_RLM2] (ML_LLM2)--(far_ML_RLM2);
        \path[name path=line_RL_IX] (RL)--(IX);
        \path[name intersections={of=line_RL_IX and line_ML_LUM2_far_ML_RUM2,by={ML_RUM2}}];
        \path[name intersections={of=line_RL_IX and line_ML_LLM2_far_ML_RLM2,by={ML_RLM2}}];
        \fill [green, opacity=15/30](ML_RUM2)--(ML_RLM2)--(ML_LLM2)--(ML_LUM2)-- cycle;
    \begin{scope}
      \clip (LL) -- (X) -- (RU) -- cycle;
      \draw[fill=blue] circle[at=(X),radius=0.35];
    \end{scope}
    \begin{scope}
      \clip (RU) -- (Y) -- (LL) -- cycle;
      \draw[fill=blue] circle[at=(Y),radius=0.35];
    \end{scope}
    \begin{scope}
      \clip (RL) -- (IX) -- (LU) -- cycle;
      \draw[fill=red] circle[at=(IX),radius=0.35];
    \end{scope}
    \begin{scope}
      \clip (LU) -- (IY) -- (RL) -- cycle;
      \draw[fill=red] circle[at=(IY),radius=0.35];
    \end{scope}
    
    \node[] at (0,0.5*\WIDTH)  {$\color{white}M^{\text{\upshape in}}_{r,0}$};
    
    \path[name path=line_RU_X] (RU)--(X);
    \path[name intersections={of=line_RU_X and lineM1U,by={Mr1Label}}];
    \node[] at ([shift={(0.07,-0.18)}]Mr1Label) {\small $M^{\text{\upshape in}}_{r,1}$};
    \node[] at ([shift={(2.1*\LENGTH/3,\WIDTH/3.27)}]Mr1Label) {\small $M^{\text{\upshape in}}_{r,1}$};

    \path[name intersections={of=line_RU_X and line_L_LUM2_Far_L_RUM2,by={Mr2Label}}];
    \node[] at ([shift={(0.14,-0.16)}]Mr2Label) {\small $M^{\text{\upshape in}}_{r,2}$};
    \node[] at ([shift={(3.2*\LENGTH/3,\WIDTH*0.48)}]Mr2Label) {\small $M^{\text{\upshape in}}_{r,2}$};
    \node[] at ([shift={(0.2,-0.05)}]X) {\tiny $\color{white}N^{\text{\upshape in}}$};
    \node[] at ([shift={(-0.14,0.13)}]Y) {\tiny\color{white} $N^{\text{\upshape in}}$};
    \node[] at ([shift={(0.25,0.1)}]IX) {\tiny\color{white} $N^{\text{\upshape out}}$};
    \node[] at ([shift={(-0.08,-0.1)}]IY) {\tiny\color{white} $\!\!\!\!\!\!\!\!\!\!\!\!N^{\text{\upshape out}}$};
    \end{tikzpicture}
    \caption{\cref{fig:Mr} with $\phistar<\pi/6$ and more detail. The two pink strips are $\Fc^{-1}(\Mrout\setminus\Mrin)$ and have slope $-1/4$. 
    The green strips are $\Fc^{-2}(\Mrout\setminus(\Mout_{r,0}\cup\Mout_{r,1}\cup\Mout_{r,2}))$ and have slope $-1/6$.\newline The intersections $\Min_{r,0}$, $\textcolor{blue}{M^{\text{\upshape in}}_{r,1}}$ and $\textcolor{blue}{M^{\text{\upshape in}}_{r,2}}$ are marked.  \newline The leftmost point of the closure of \textcolor{blue}{$M^{\text{\upshape in}}_{r,0}$} is $(\pi-\phistar,\frac{\pi}{2})$. The leftmost point of the closure of \textcolor{blue}{$M^{\text{\upshape in}}_{r,1}$} is $(\frac{2\pi}{3}-\phistar,\frac{\pi}{3})$. The leftmost point of the closure of \textcolor{blue}{$M^{\text{\upshape in}}_{r,2}$} is $(\frac{\pi}{2}-\phistar,\frac{\pi}{4})$\newline$\hat{M}$ is the union of the two pink strips and the center rhombus $\Mout_{r,0}=\Min_{r,0}=\Mrin\cap\Mrout$.\newline\Nin\ and \Nout\ are neighborhoods of $x_*,y_*,Ix_*,Iy_*$ respectively.}\label{fig:MrN}.
\end{center}
\end{figure}
\begin{notation}[Returns to $\hat M$]\label{def:MhatReturnOrbitSegment}
To study $\hat F$, we write the orbit segment of $x\in\hat{M}$ to its first return to $\hat M$ as 
\begin{equation}\label{eqMhatOrbSeg}
        \underbracket{x}_{\in \hat{M}},\underbracket{\Fc(x),\cdots}_{\notin\hat M},\underbracket{\Fc^{\sigma(x)}(x)=\hat{F}(x)}_{\in \hat{M}},\text{ where }\sigma(x)\dfn\inf\{k>0\mid \Fc^{k}(x)\in\hat{M}\}.  
\end{equation}
By \eqref{def:Mhat}, either $x\in\Fc^{-1}(\Mrout\setminus\Mrin)$, in which case we set $x_0\dfn\Fc(x)\in\Mrout\setminus\Mrin$, or otherwise $x_0\dfn x\in\Mrout\cap\Mrin$. Hence, this orbit segment includes an $x_0\in\Mrout$ among its first two points. Then $x_1\dfn \Fc(x_0)\in \MRin$ and, taking $n_1\dfn m(x_1)$ from \cref{def:SectionSetsDefs}, $x_2\dfn\Fc^{n_1+1}(x_1)\in \Mrin$ (because  $\hat M\subset M_r$).

Then $m(x_2)$ determines whether $\hat{F}(x)=x_2$ or $\hat{F}(x)=\Fcal^{m(x_2)-1}(x_2)$ as follows.
    
Case i): $m(x_2)\le1$, that is $x_2\in(\textcolor{blue}{M^{\text{\upshape{in}}}_{r,0}}\cup\textcolor{blue}{M^{\text{\upshape in}}_{r,1}})\subset\hat{M}$. Then $x_2=\hat{F}(x)$.
    
Case ii): $m(x_2)\ge2$, that is $x_2\in{\Min_{r,n}}$ for some $n=m(x_2)\ge2$. Then $\Fc^{m(x_2)}(x_2)\in\Mrout\setminus\Mrin$ and $\Fc^{m(x_2)-1}(x_2)=\hat{F}(x)\in\Fc^{-1}(\Mrout\setminus\Mrin)\subset\hat{M}$ since $\Fc^i(x_2)\notin\Fc^{-1}(\Mrout\setminus\Mrin)$ for $i=0,\cdots,m(x_2)-2$.
\end{notation}

\begin{remark}[``Bad regions'']\label{REMBadRegions}
    One inconspicuous feature of \cref{fig:Mr,fig:MrN} is worth attention. It turns out that the small regions \Nin\ and \Nout\ in \cref{def:SectionSetsDefs,fig:MrN} contain those points which are problematic with respect to the expansion we wish to establish (\cref{contraction_region}). Therefore we now prove, in effect, that they are as small as they look in \cref{fig:MrN} in the sense that they do not meet the other regions shown. Studying the dynamics on these regions entails significant effort and constrains the parameters for which we can prove uniform hyperbolicity. In effect, that work entails treating the lemon billiard as a perturbation of a 1-petal billiard in \cref{sec:SACTWFELBWB}. However, it is only there that this helps; a perturbative overall strategy does not seem promising.
\end{remark}\begin{proposition}[Small bad regions]\label{Prop:Udisjoint}
    If $\phistar\in\big(0,\tan^{-1}(1/3)\big)$ and $R$ satisfies \eqref{eqJZHypCond}, then $\Nin\cap\big(\textcolor{blue}{M^{\text{\upshape in}}_{r,0}}\cup{\textcolor{blue}{M^{\text{\upshape in}}_{r,1}}}\cup{\textcolor{blue}{M^{\text{\upshape in}}_{r,2}}}\big)=\emptyset$ and $\Nout\cap\big(\textcolor{red}{M^{\text{\upshape in}}_{r,0}}\cup\textcolor{red}{M^{\text{\upshape in}}_{r,1}}\cup\textcolor{red}{M^{\text{\upshape in}}_{r,2}}\big)=\emptyset$.
\end{proposition}
\begin{proof}
$\Nin$ has two components shown in \cref{fig:MrN}, a left component near $(\phistar,\phistar)$, and a right component near $(2\pi-\phistar,\pi-\phistar)$. They are symmetric with respect to the center of $M_r$. ($I \circ J$ in \cref{DefSymmetries}).
    
The green strips in \cref{fig:MrN} have slope $-1/6$ and the blue strip  \Mrin\  has slope $1/2$. \textcolor{blue}{$M^{\text{\upshape in}}_{r,2}$} has two components. The left component is the interior of the parallelogram with vertices $(\pi/2-\phistar,\pi/4)$, $(\pi/2-\phistar+\phistar/2,\pi/4+\phistar/4)$, $(\pi/2+\phistar,\pi/4)$ and $(\pi/2+\phistar-\phistar/2,\pi/4-\phistar/4)$ among which $(\pi/2-\phistar,\pi/4)$ is the leftmost. The left and right components are symmetric with respect to the center of $M_r$ (i.e., under $I \circ J$, see \cref{DefSymmetries}). It suffices to show that the left components of \textcolor{blue}{$M^{\text{\upshape in}}_{r,2}$} and $\Nin$ are disjoint.

Note that in \cref{fig:MrN} by elementary geometry, the Euclidean distance between $(\pi/2-\phistar,\pi/4)$ and $(\phistar,\phistar)$ is the smallest Euclidean distance from any point in the closure of the left component of $\textcolor{blue}{M^{\text{\upshape in}}_{r,2}}$ to $(\phistar,\phistar)=Iy_*$, and this smallest Euclidean distance is $\sqrt{(\nicefrac\pi2-2\phistar)^2+(\nicefrac\pi4-\phistar)^2}=\sqrt{5}(\nicefrac{\pi}{4}-\phistar)$ since $0<\phistar<\tan^{-1}(1/3)<\pi/4$. Thus,\COMMENT{It does not matter, but 500 is a bizarre choice; why not 7500, which is 30000/4?? Then $7500\cdot5\cdot(.46)^2=7935>89^2$, so  17 can be replaced by 89. I don't suggest doing that, but this leaves room for a nicer number than 17, such as 20. Alternatively, we need much less than 30\,000 to make \emph{this} argument work.\wentao{OK, now a lot of things are changed, The $\delta_*$ is also changed. Noting is bizarre now.}} 
\[
\frac Rr\overbracket>^{\eqref{eqJZHypCond0}}1700\overbracket{>}^{\mathclap{\substack{0<\phistar<\tan^{-1}(1/3)\Rightarrow \pi/4-\tan^{-1}1/3=\tan^{-1}1/2>0.46\Rightarrow \frac{4.5^2}{5(\pi/4-\phistar)^2}<\frac{21}{5(\pi/4-\tan^{-1}(1/3))^2}<1700\\|}}}\frac{4.5^2}{5(\pi/4-\phistar)^2},
\]

    which implies $\sqrt{5}(\nicefrac{\pi}{4}-\phistar)>\delta_*=4.5\sqrt{r/R}$. Hence, the left components of $\textcolor{blue}{M^{\text{\upshape in}}_{r,2}}$ and \Nin are disjoint. Thus $\Nin\cap\textcolor{blue}{M^{\text{\upshape in}}_{r,2}}=\emptyset$. 
    
Similar arguments about the distance between the leftmost points of the closure of $\textcolor{blue}{M^{\text{\upshape in}}_{r,0}}$ and $I y_*$ and the distance between the leftmost points of the closure of $\textcolor{blue}{M^{\text{\upshape in}}_{r,1}}$ and $I y_*$ yield $\Nin\cap\textcolor{blue}{M^{\text{\upshape in}}_{r,0}}=\emptyset$ and $\Nin\cap\textcolor{blue}{M^{\text{\upshape in}}_{r,1}}=\emptyset$. Thus, we have $\Nin\cap\big(\textcolor{blue}{M^{\text{\upshape in}}_{r,0}}\cup\textcolor{blue}{M^{\text{\upshape in}}_{r,1}}\cup\textcolor{blue}{M^{\text{\upshape in}}_{r,2}}\big)=\emptyset$, and, by symmetry $I \circ J$ in \cref{DefSymmetries},  $\Nout\cap\big(\textcolor{red}{M^{\text{\upshape out}}_{r,0}}\cup\textcolor{red}{M^{\text{\upshape out}}_{r,1}}\cup\textcolor{red}{M^{\text{\upshape out}}_{r,2}}\big)=\emptyset$.
\end{proof}
\begin{remark}\label{remark:MoreUdisjoint}
    In fact, using the same proof we see that if $0<\phistar<\tan^{-1}\big(1/3\big)$ with $R>\frac{30000r}{\sin^2\phistar}$, then $\Nin\cap{\textcolor{blue}{M^{\text{\upshape in}}_{r,3}}}$, $\Nin\cap{\textcolor{blue}{M^{\text{\upshape in}}_{r,4}}}$, $\Nout\cap\textcolor{red}{M^{\text{\upshape out}}_{r,3}}$, $\Nout\cap\textcolor{red}{M^{\text{\upshape out}}_{r,4}}$ are also empty.
\end{remark}


Unless stated otherwise, we henceforth assume $0<\phistar<\tan^{-1}{(1/3)}$ and \eqref{eqJZHypCond}.

\begin{remark}[Uniform hyperbolicity of return map on other sections]\label{rmk:othersectionUH}
\cref{thm:UniformExpansionOnMhat1,Col:EuclideanExpansionMhat1} give the uniform hyperbolicity of the return map $\hat{F}_1$ on $\hat{M}_1$. \cref{thm:MroutreturnUniformExpansion,corollary:MroutreturnUniformExpansion} give the uniform hyperbolicity of the return map $\tilde{F}$ on $\Mrout$.
Although the sections $\hat{M}$ and $\hat{M}_1$ in in \eqref{def:Mhat} and \eqref{def:symetricsectionMhatPlusMinus} are not connected (\cref{fig:MrN,fig:hatUplusUminus}), in \cref{thm:UniformExpansionOnMhat,thm:UniformExpansionOnMhat1} the unstable invariant cone on those sections is inside the positive quadrant cone with associated quadratic norm being nonnegative. This will be convenient for our local ergodicity proof.

Choosing reference sections $\Mrout$ also works fine for local ergodicity proof and we also prove \cref{thm:MroutreturnUniformExpansion,corollary:MroutreturnUniformExpansion} the uniform hyperbolicity for the return map on $\Mrout$. In this case, with the stable cone on $\Nout$ is larger than the quadrant cone,\COMMENT{\Large At the end of the paper it says "Using symmetries $I$ and $J$, we conclude that the stable cone family"---so here there should not be large stable cones.\wentao{thanks, I mean here in the problematic $\Nin$ unstable cone is large, in the problematic $\Nout$ stable cone is large. In the end, since the unstable cone in $\Nout$ is in quadrants, by symmetry the stable cone in $\Nin$ is in quadrants.}} we have to exclude those ``bad'' regions inside $\Nout$ in the verification for L4, L3 local ergodicity conditions in \cite{MR3082535}.\COMMENT{Is this sentence illuminating in any way, given the previous comment?} 

Because they come with nicer cones, we still use the more fragmented sections. This choice will not cause difficulties in the proof of global ergodicity so long as the set of double-order singularities is countable \cite[Proposition 6.19] {cb} (Also See \cite[Theorem 6.20]{cb}). 
\end{remark}

\section[Expansion estimates]{Expansion estimates}\label{sec:ECE}
\label{sec:FPCFICBMD}
This section gives lower bounds for the differential of the billiard map in a case-by-case fashion (\cref{eqMainCases,TEcase_a0,TEcase_a1,TEcase_b,TEcase_c}) and establishes that those cases in which these lower bounds do not produce expansion are confined to a small portion of the section (\cref{contraction_region}). This is sufficient information to show strict invariance of the families of cones in \cref{thm:UniformExpansionOnMhat,Corollary:EuclideanExpansionMhat,thm:UniformExpansionOnMhat1,thm:MroutreturnUniformExpansion} (\cref{corollary:InvQuadrant}; see also \cref{corollary:halfquadrantcone}) and strict expansion in them upon iteration (\cref{def:ApettrajectoryExpansion}). The proofs are carried out in later sections.

Toward establishing uniform expansion in invariant cones, a core step is to estimate the expansion obtained in the transitions from one of the boundary arcs to the other. This turns out to require a case-by-case approach, which already occurred when the existence of an invariant cone family was established in \cite{hoalb}. We begin with a presentation of the main cases and the notation and prior results about invariant cones.
\subsection{Invariant cones}
\label{subsec:pclf}
For $x\in\hat{M}$ consider the orbit segment \eqref{eqMhatOrbSeg} (notations as in \cref{def:BasicNotations,def:SectionSetsDefs,def:MhatReturnOrbitSegment}) and the points in it found right after \eqref{eqMhatOrbSeg}:
\begin{equation}\label{eqPtsFromeqMhatOrbSeg}
    x_{0}=(\phi_0,\theta_0)\in\Mrout,\ x_1=(\Phi_1,\theta_1)\dfn \Fc(x_0)\in \MRin
    \text{ and }x_{2}=(\phi_2,\theta_2)=\Fc^{n_1+1}(x_1)\in\Mrin\text{ with }n_1\ge0.
\end{equation}
\begin{figure}[h]
\begin{center}
\begin{sideways}
\begin{tikzpicture}[xscale=0.53, yscale=0.53]
    \tkzDefPoint(0,0){Or};
    \tkzDefPoint(0,-4.8){Y};
    \pgfmathsetmacro{\rradius}{4.5};
    \pgfmathsetmacro{\bdist}{8.7};
    \tkzDefPoint(0,\bdist){OR};
    \pgfmathsetmacro{\phistardeg}{40};
    \pgfmathsetmacro{\XValueArc}{\rradius*sin(\phistardeg)};
    \pgfmathsetmacro{\YValueArc}{\rradius*cos(\phistardeg)};
    \clip(-4.9,-13) rectangle (6.6, 11.5);
    \tkzDefPoint(\XValueArc,-1.0*\YValueArc){B};
    \tkzDefPoint(-1.0*\XValueArc,-1.0*\YValueArc){A};
    \pgfmathsetmacro{\Rradius}{veclen(\XValueArc,\YValueArc+\bdist)};
    \draw[name path=Cr,dashed,ultra thin] (Or) circle (\rradius);
    \tkzDrawArc[thick](Or,B)(A);
    \tkzDrawArc[name path=CR, thick](OR,A)(B);
    \tkzDrawArc[blue, dashed, ultra thin](OR,B)(A);
   
    \pgfmathsetmacro{\PXValue}{\Rradius*sin(7)};
    \pgfmathsetmacro{\PYValue}{\bdist-\Rradius*cos(7)};
    \tkzDefPoint(\PXValue,\PYValue){P};
    \draw[name path=PQ_Far,blue,ultra thin](P)--+(70:7)coordinate(Q_Far);
    \draw[name path=PT_Far,blue,ultra thin](P)--+(124:9);
    \draw[name path=PT0_Far,blue,dashed,ultra thin](P)--+(-56:7);
    \draw[name path=PQ0_Far,blue,dashed,ultra thin](P)--+(-110:9);
    \draw[red,dashed,->](P)--+(7:4)coordinate(PposX);
    \draw[name path=PPnegX,red,dashed](P)--+(-173:6)coordinate(PnegX);
    \draw[name path=linePOR,red,dashed,->](P)--(OR);
     \path[name path= linePOr] (P)--(Or);
    \tkzDefPoint(0,-5.5){Y};
    \tkzDefPoint(0,\bdist-\Rradius){C};
    \draw[red,dashed](OR)--(C);
    \tkzMarkAngle[arc=lll,ultra thin,size=1.7](Or,OR,P);
    \path[name intersections={of=Cr and PQ_Far,by={Q}}];
    \draw[,blue,ultra thin,dashed](Q)--+(70:7);
    \path[name intersections={of=Cr and PQ0_Far,by={Q0}}];
    \path[name intersections={of=Cr and PT0_Far,by={T0}}];
    \path[name intersections={of=Cr and PT_Far,by={T}}];
    \tkzMarkAngle[arc=lll,very thin,size=3.2](PposX,P,Q);
    \tkzMarkRightAngle(OR,P,PnegX);
    \node[rotate=-90] at ($(OR)+(-85:2.3)$)  {$\Phi_1$};
    \node[rotate=-90] at ($(P)+(30:3.5)$)  {\Large {$\theta_1$}};

    \node[rotate=-90] at ([shift={(0.4,2.2)}]PposX) {\textcolor{blue}{\small circle $C_R$}};
    
    \tkzLabelPoint[left,rotate=-90](A){$A$};
    \tkzLabelPoint[left,rotate=-90](B){$B$};
    \tkzLabelPoint[left,rotate=-90](OR){$O_R$};
    \tkzLabelPoint[below,rotate=-90](Or){$O_r$};
    \tkzLabelPoint[above,rotate=-90]([shift={(0,0.35)}]P){$P$};
    \tkzLabelPoint[above,rotate=-90](Q){$Q$};
    \tkzLabelPoint[left,rotate=-90](T){$T$};
    \tkzLabelPoint[below right,rotate=-90](Q0){$Q_0$};
    \tkzLabelPoint[right,rotate=-90](T0){$T_0$};
    \tkzLabelPoint[right,rotate=-90](PposX){$TP$};
    \tkzDefLine[perpendicular=through Or](P,PnegX)\tkzGetPoint{M_far};
    \path[name path= lineOrMfar] (Or)--(M_far);
    \path[name path= PPnegX] (P)--(PnegX);
    \path[name path= PPposX] (P)--(PposX);
    \path[name intersections={of=lineOrMfar and PPnegX,by={M}}];
    \path[name intersections={of=Cr and PPnegX,by={A0}}];
    \path[name intersections={of=Cr and PPposX,by={B0}}];
    \coordinate (Or_mirror) at ($(Or)!2!(M)$);
    \tkzDrawArc[name path=mirrorCr, ultra thin, dashed](Or_mirror,A0)(B0);
    \tkzLabelPoint[below,rotate=-90](Or_mirror){$O'_r$};

\path[name intersections={of=mirrorCr and PQ0_Far,by={Q1}}];
    \path[name intersections={of=mirrorCr and PT0_Far,by={T1}}];
    \tkzLabelPoint[above,rotate=-90](Q1){$Q_1$};
    \tkzLabelPoint[above,rotate=-90](T1){$T_1$};
    \tkzDrawPoints(A,Or,OR,P,B,T,Q,Or_mirror,A0,T0,Q0,T1,Q1);
    \begin{scope}[ultra thin,decoration={markings, mark=at position 0.75 with {\arrow{>}}}
    ] 
    \draw[blue,postaction={decorate}](P)--(Q);
    \draw[blue,postaction={decorate}](T)--(P);
    \end{scope}
    \tkzLabelPoint[above,rotate=-90]([shift={(0.15*\rradius,1.1*\rradius)}]Or){\Large $\Gamma_r$};
    \tkzLabelPoint[above,rotate=-90]([shift={(0.35*\rradius,-1.05*\rradius)}]Or_mirror){\Large $\Gamma_r'$};
\end{tikzpicture}
\end{sideways}
\captionsetup{type=figure}
\caption{$n_1=0$ case: the trajectory hitting $P=p(x_1)=p(\Fcal^{-1}(x_2))$ with $T=p(\Fcal^{-1}(x_1))$ and $Q=p(x_2)$. $\tau_0=|TP|$, $\tau_1=|PQ|$. With $O'_r$ is symmetric to $O_r$ and arc $\Gamma'_r$ is symmetric to arc $\Gamma_r$ with respect to the tangent line to $\Gamma_R$ at $P$.\newline 
Extend line segment $TP$ to meet the circle $C_r$ at $T_0$ and the arc $\Gamma'_r$ at $T_1$. 
\newline Extend line segment $QP$ to meet the circle $C_r$ at $Q_0$ and the arc $\Gamma'_r$ at $Q_1$.
}\label{fig:lengthfunctonFigure}
\end{center}
\end{figure}
\begin{notation}\label{def:fpclf} 
For $x_0\in\Mrout$, $x_1\in\MRin$ and $x_2\in\Mrin$, the free-path and chord-length functions (see \cref{fig:lengthfunctonFigure}) are:
\begin{itemize}
    \item $\tau_0$ is the distance between $p(x_0)\in\Gamma_r$ and $p(\Fc(x_0))\in\Gamma_R$.
    \item $\tau_1$ is the distance between $p(x_2)\in\Gamma_r$ and $p(\Fc^{-1}(x_2))\in\Gamma_R$, or $L_1\dfn\tau_1$ to denote length.
    \item $d_0\dfn r\sin{\theta_0}$ is the half-length of the chord that connects the two intercepts of $C_r$ with the line that goes through $p(x_0)$ and $p(\Fc(x_0))$, and in \cref{fig:lengthfunctonFigure} $|TT_0|$ is such a chord.
    \item $L_0\dfn2d_0-\tau_0$ and in \cref{fig:lengthfunctonFigure}, $|PT_0|=L_0$.
    \item $d_1\dfn R\sin{\theta_1}$ is the half-length of the chord that connects the two intercepts of $C_R$ with the line that goes through $p(x_1)$ to $p(\Fc(x_1))$. In \cref{fig:lengthfunctonFigure}, the line $PQ$ has two intersections with $C_R$. $P=p(x_1)$ is one of the two intersections.
    \item $d_2\dfn r\sin{\theta_2}$ is the half-length of the chord that connects the two intercepts of $C_r$ with the line that goes through $p(x_2)$ and $p(\Fc^{-1}(x_2))$, and in \cref{fig:lengthfunctonFigure} $|QQ_0|$ is such a chord.
\end{itemize}
\end{notation}
\begin{remark}\label{remark:deffpclf}
$d_0,d_1,d_2,\tau_0,\tau_1$ are continuous functions of nonsingular $x_1\in\MRin\subset\Int M_R$. In \cref{subsec:continuousextension,subsec:A1contractionCase,sec:ProofTECase_c}, they are continuously extended to the boundary of $M_R$. Because any free path ( i.e., a segment between collisions) is contained in the disk surrounded by $C_r$ and also inside the disk surrounded by $C_R$, we have the following.
\[2d_2>\tau_1=L_1,\quad 2d_0-\tau_0=L_0>0,\quad2d_1>\tau_0, \quad2d_1>\tau_1.
\]
\end{remark}
\begin{remark}[\cref{fig:lengthfunctonFigure}]\label{remark:symmetrywiththetangentline}
If $n_1=0$, let $P=p(x_1)=p(\Fcal^{-1}(x_2))$ as in \cref{fig:lengthfunctonFigure}, then with $|TP|=\tau_0$, $|PQ|=\tau_1$, $|TT_0|=2d_0$ since $|PQ|$ is symmetric to $|PT_1|$ with respect to the tangent line to $\Gamma_R$ at $P$, we see $\tau_0+\tau_1-2d_0=|TP|+|PQ|-|TT_0|=|TP|+|PT_1|-|TT_0|=|TT_1|-|T_0T_1|=|T_0T_1|>0$. Similarly, by symmetry $|Q_1P|=|TP|=\tau_0$ and $|QQ_0|=2d_2$, we have $\tau_0+\tau_1-2d_2=|TP|+|PQ|-|QQ_0|=|Q_1P|+|PQ|-|QQ_0|=|QQ_1|-|Q_0Q_1|=|Q_0Q_1|>0$. Therefore, we have the following inequalities.
\[
    \tau_0+\tau_1-2d_0>0,\quad \tau_0+\tau_1-2d_2>0,\quad \tau_0+\tau_1-d_0-d_2>0.
\]
\end{remark}

Our case-by-case approach is related to different features of the derivative matrix of the return map---the good scenarios are those when all entries are positive or all are negative; in both cases, the positive quadrant is preserved. 
We introduce these cases now. (It may be good to remember here that $R>300\,000r$, hence $\sqrt{{r}/{R}}>547r/R$.)
\begin{equation}\label{eqMainCases}
    \begin{cases}\text{Case (a) in \eqref{eqPtsFromeqMhatOrbSeg}:}&d_1=R\sin{\theta_1}\ge 2r,\quad\begin{cases}
        \text{Case (a0):}&\sin{(\theta_1)}\ge2\sqrt{{r}/{R}},\\\text{Case (a1):}&\sin{(\theta_1)}\in [{2r}/{R},2\sqrt{{r}/{R}}).
    \end{cases}\\
    \text{Case (b) in \eqref{eqPtsFromeqMhatOrbSeg}:}&d_1=R\sin{\theta_1}< 2r\text{ and }n_1=0\text{, i.e., the orbit segment has exactly one collision on }\Gamma_R,\\
    \text{Case (c) in \eqref{eqPtsFromeqMhatOrbSeg}:}&d_1=R\sin{\theta_1}< 2r\text{ and }n_1\ge1\text{, i.e., the orbit segment has $n_1+1\ge 2$ collisions on }\Gamma_R.
    \end{cases}
\end{equation}
\cref{JZProp3.6a,JZProp3.6b,JZProp3.6c} below (contained in \cite[Proposition 3.6]{hoalb}) provide invariance of the first-quadrant cone family. They assume that $R$ satisfies \eqref{eqJZHypCond} and apply to the points in \eqref{eqPtsFromeqMhatOrbSeg} from \eqref{eqMhatOrbSeg}.
\begin{proposition}[{\cite[Proposition 3.6]{hoalb}}]\label{JZProp3.6a} In Case {(a)}\label{caseA}, $x_2=\Fc(x_1)\in\Mrin$, $D_{x_0}\Fc^2$ is a 2 by 2 matrix with entries depending on $\tau_0,\tau_1,d_0,d_1,d_2$, and all 4 entries are negative \cite[equations (4.2) through (4.5)]{hoalb}.
\end{proposition}
\begin{proposition}[{\cite[Proposition 3.6]{hoalb}}]\label{JZProp3.6b} In Case {(b)}\label{caseB}, $x_2=\Fc(x_1)\in\Mrin$, $x_{-1}:=\Fc^{-1}(x_0)\in M_r$, $x_3:=\Fc(x_2)\in M_r$, $D_{x_{-1}}\Fc^4$ is a 2 by 2 matrix with entries depending on $\tau_0,\tau_1,d_0,d_1,d_2$, and all 4 entries are negative \cite[Equations (4.10) through (4.19)]{hoalb}.
\end{proposition}
\begin{proposition}[{\cite[Proposition 3.6]{hoalb}}]\label{JZProp3.6c} In Case {(c)}\label{caseC},  $x_3\dfn \Fc(x_2)\in M_r$, $x_{-1}\dfn \Fc^{-1}(x_0)\in M_r$, $D_{x_{-1}}\Fc^{n_1+4}$ is a 2 by 2 matrix with entries depending on $\tau_0,\tau_1,d_0,d_1,d_2,n_1$, and  all 4 entries are positive in this case \cite[Section 5]{hoalb}.
\end{proposition}
\begin{remark}
The calculations of the entries of the derivative matrices in these three cases can be found in \cite[Sections 3--5]{hoalb} except that in subcase (a0) \cref{TEcase_a0}, we need to extend the definition of $D_{x_0}\Fc^2$ to the boundary of $M_R$ (\cref{subsec:continuousextension}).
\end{remark}
Unlike the invariance of cones, expansion of vectors in cones (the main novelty of this paper) is quite difficult outside Case~(a0) in \eqref{eqMainCases} (see \cref{TEcase_a0,TEcase_a1,TEcase_b,TEcase_c} below). Let us note here that those harder cases are confined to a small part of the phase space: cases {(a1) (\Cref{TEcase_a1}), (b) (\Cref{TEcase_b}), and  (c) (\Cref{TEcase_c}}) in \eqref{eqMainCases} arise only when $x_0\in\Nout$ and $x_2\in\Nin$ (see \cref{def:SectionSetsDefs}). This is \cref{contraction_region} below (page \pageref{contraction_region}).

\subsection{Control of \texorpdfstring{$D\hat F$}{Control of return map differential}}
We now proceed to the central object of this article, control of $D\hat F$---which is clearly crucial for \cref{RMASUH}. We will eventually establish expansion of the p-metric, beginning with \cref{TEcase_a0,TEcase_a1,TEcase_b,TEcase_c}, one proposition for each of the four cases in \eqref{eqMainCases}. The three main cases ultimately take up one section each (\cref{sec:ProofTECase_a0_a1,sec:ProofTECase_b,sec:ProofTECase_c}).

{\bfseries Standing assumptions for \cref{TEcase_a0,TEcase_a1,TEcase_b,TEcase_c}.} Fix $r,R$ satisfying the hyperbolicity condition \eqref{eqJZHypCond}, fix $\phistar\in(0,\tan^{-1}{(1/3)})$, $\hat F$ from \cref{def:SectionSetsDefs}, and the cone family $C_x$ from \cref{RMASUH}.

\begin{proposition}[Case (a0)]\label{TEcase_a0} 
There exists a constant $\lambda_c>0$ such that 
if $\sin{\theta_1}\ge\sqrt{{4r}/{R}}$ for the points in \eqref{eqPtsFromeqMhatOrbSeg} with nonsingular $x$, then
\begin{enumerate}
\item\label{itemTEcase_a0-1} the cone family $C_x$ is strictly invariant under $D\hat{F}$, i.e., $D_{x}\hat{F}(C_x)\subset\text{\big\{interior of }C_{\hat{F}(x)}\text{\big\}}$,
\item\label{itemTEcase_a0-2} $\frac{\|D_{x}\hat{F}(dx)\|_\p}{\|dx\|_\p}>1+\lambda_c$ for all $dx\in C_x$ and $\lambda_c>0$ is a constant determined by $r$, $R$, $\phistar$ configuration.
\end{enumerate}
\end{proposition}
\begin{proof}
These are the conclusions of \cref{theorem:UniformExpansionA0} below (see page \pageref{theorem:UniformExpansionA0} and the proof in \cref{sec:ProofTECase_a0_a1}).
\end{proof}
It is the kind of statement one would wish to prove in all subsequent cases as well, but instead, the corresponding propositions have lower bounds less than 1, necessitating subsequent work. Here they are. 
\begin{proposition}[Case (a1)]\label{TEcase_a1} 
If $\sin{\theta_1}\in [{2r}/{R},\sqrt{{4r}/{R}})$ for the points in \eqref{eqPtsFromeqMhatOrbSeg} with nonsingular $x$, then
\begin{enumerate}
    \item\label{itemTEcase_a1-1} the cone family $C_x$ is strictly invariant under $D\hat{F}$, i.e., $D_{x}\hat{F}(C_x)\subset\text{\big\{interior of }C_{\hat{F}(x)}\text{\big\}}$,
    \item\label{itemTEcase_a1-2} $\frac{\|D_{x}\hat{F}(dx)\|_\p}{\|dx\|_\p}>0.26$, for all $dx\in C_x$.
\end{enumerate}
\end{proposition}
\begin{proof}
This is \cref{theorem:contractioncontrolA1} below (see page \pageref{theorem:contractioncontrolA1} and its proof in \cref{sec:ProofTECase_a0_a1}).
\end{proof} Note that this does not establish expansion; it controls how far from expansion the return map is in this case. Subsequent iterates are needed to convert this control into actual expansion. The same goes for the subsequent cases; the next case is the most daunting because the lower bound is quite small. For ease of exposition, we will take the liberty to refer to these lower bounds as lower bounds for expansion even when they are less than 1.
\begin{proposition}[Case (b)]\label{TEcase_b}
If 
$\sin{\theta_1}<{2r}/{R}$ and $n_1=0$ for the points in \eqref{eqPtsFromeqMhatOrbSeg}, then
\begin{enumerate}
    \item\label{itemTEcase_b1} the cone family $C_x$ is strictly invariant under $D\hat{F}$, i.e., $D_{x}\hat{F}(C_x)\subset\text{\big\{interior of }C_{\hat{F}(x)}\text{\big\}}$,
    \item\label{itemTEcase_b2} $\frac{\|D_{x}\hat{F}(dx)\|_\p}{\|dx\|_\p}>0.05$ for all $dx\in C_x$.
\end{enumerate}
\end{proposition}
\begin{proof}
\eqref{itemTEcase_b1} is  \cref{corollary:halfquadrantcone}. \eqref{itemTEcase_b2} is \cref{theorem:proofofTEcase_b}.
\end{proof}
\begin{proposition}[Case (c)]\label{TEcase_c}
If $\sin{\theta_1}<{2r}/{R}$ and $n_1\ge1$ for the points in \eqref{eqPtsFromeqMhatOrbSeg}, then
\begin{enumerate}
    \item\label{itemTEcase_c1} the cone family $C_x$ is strictly invariant under $D\hat{F}$, i.e., $D_{x}\hat{F}(C_x)\subset\text{\big\{interior of }C_{\hat{F}(x)}\text{\big\}}$.
    \item\label{itemTEcase_c2} $\frac{\|D_{x}\hat{F}(dx)\|_\p}{\|dx\|_\p}>1-1743\big(r/R\big)-15450\big(r/R\big)^2>0.9$, for all $dx\in C_x$.
    \item\label{itemTEcase_c3} when $n_1\ge4$, $\frac{\|D_x\hat{F}(dx)\|_\p}{\|dx\|_\p}>\frac{0.96}{1.01}n_1^2+5.655n_1-8.628>29.1$, for all $dx\in C_x$ (\cref{subsec:expansionLarge}).
\end{enumerate}
\end{proposition}
\begin{proof}
\eqref{itemTEcase_c1} is \cref{corollary:halfquadrantcone}. \eqref{itemTEcase_c2} follows from \cref{lemma:AformularForLowerBoundOfExpansionForn1largerthan0,proposition:casecexpansionn1eq3,proposition:casecexpansionn1eq2,proposition:casecexpansionn1eq1}. \eqref{itemTEcase_c3} is obtained from  \cref{lemma:AformularForLowerBoundOfExpansionForn1largerthan0,proposition:casecexpansionn1ge4}.
\end{proof}
\subsection{Expansion for an iterate}
\Cref{TEcase_a0,TEcase_a1,TEcase_b,TEcase_c} in the previous section showed invariance of cones and:\newline
(a0): Expansion $>1+\lambda_c(r,R,\phistar)>1$.\newline
(a1): ``Expansion'' $>0.26$.\newline
(b): ``Expansion'' $>0.05$.\newline
(c): ``Expansion'' $>-.9$, and $>\frac{0.96}{1.01}n_1^2+5.6n_1-8.6>29$ if $n_1\ge4$.\newline
We next obtain from this a strictly invariant cone family on $\hat M$ (\cref{corollary:InvQuadrant}) and then state the uniform expansion implied by  \Cref{TEcase_a0,TEcase_a1,TEcase_b,TEcase_c,contraction_region} and three approximation theorems in  \Cref{sec:SACTWFELBWB}.

As noted above, collisions between transitions are often needed to obtain expansion. We now present this mechanism. The first ingredient is \cref{contraction_region}, which narrows down the possibilities for the problematic cases, i.e., those other than Case~(a0).
From \cref{contraction_region,Prop:Udisjoint,caseA,caseB,caseC} we get a strictly invariant cone family on $\hat M$ as follows.
\begin{proposition}\label{corollary:InvQuadrant}
If $0<\phistar<\tan^{-1}{\big(1/3\big)}$ and $R$ satisfies \eqref{eqJZHypCond}, then the quadrant cone family $\mathcal{Q}_x(\mathrm{I,III})\dfn\big\{\frac{d\theta}{d\phi}\in[0,+\infty]\big\}$ is strictly invariant under the differential $D\hat{F}$ of the return map from \cref{def:SectionSetsDefs}.
\end{proposition}
\begin{proof}
As noted before \cref{def:MhatReturnOrbitSegment}, nonsingular $x\in\hat{M}$ have a return orbit segment containing $x_0\in\Mrout$, $x_1=(\Phi_1,\theta_1)\in\MRin$, $x_2\in\Mrin$. By \cref{def:MhatReturnOrbitSegment} we have:
\[
x_0\in\Mrout\xrightarrow{\Fc}x_1\in\MRin\xrightarrow{\Fc}\cdots \xrightarrow{\Fc}\Fc^{n_1}(x_1)\in\MRout\xrightarrow{\Fc}x_2\in\Mrin.
\]
In cases (a1), (b) and (c) of \eqref{eqMainCases}, \cref{contraction_region} gives $x_0\in\Nout$ and $x_2\in\Nin$, then $x\notin\Mrin\cap\Mrout$ because otherwise $x=x_0\in \Nout\cap \textcolor{red}{M^{\text{\upshape{out}}}_{r,0}}=\emptyset$ by \cref{Prop:Udisjoint}. Thus, by \cref{def:SectionSetsDefs}, we have $x\in \Fc^{-1}(\Mrout\cap\Mrin)$ and $x_0=\Fc(x)\in\Mrout\setminus\Mrin$ .

On the other hand, $x_2\in\Nin$ and \cref{Prop:Udisjoint} imply $x_2\notin\big(\textcolor{blue}{M_{r,0}^{\text{\upshape in}}}\cup\textcolor{blue}{M_{r,1}^{\text{\upshape in}}}\cup\textcolor{blue}{M_{r,2}^{\text{\upshape in}}}\big)$, so $m(x_2)\ge3$, hence $x_2\notin\Mrout$ and $x_3\dfn\Fc(x_2)\notin\Mrout$, and $\Fc^2(x_2)=\Fc(x_3)\notin\Mrout$. Hence, the \cref{def:MhatReturnOrbitSegment} return orbit segment in cases (a1), (b) and (c) is 
\begin{equation}\label{eq:a1bccasesorbitsegment}
\underbracket{x\in\hat{M}\xrightarrow{\Fc}x_0\in\Mrout}_{\text{collisions on }\Gamma_r}\xrightarrow{\Fc}\underbracket{x_1\in\MRin\xrightarrow{\Fc}\cdots \xrightarrow{\Fc}\Fc^{n_1}(x_1)\in\MRout}_{\text{collisions on }\Gamma_R}\xrightarrow{\Fc}\underbracket{x_2\in\Mrin\xrightarrow{\Fc}x_3\cdots\xrightarrow{\substack{m(x_2)-2\\\text{ steps}}}\Fc^{m(x_2)-2}(x_3)=\hat{F}(x)}_{\text{collisions on }\Gamma_r},
\end{equation} with 
$D_{x_3}\Fc^{m(x_2)-2}=\begin{pmatrix}
1 & 2(m(x_2)-2) \\
0 & 1 
\end{pmatrix}$ in $\phi\theta$-coordinates. Hence, $m(x_2)\ge3$ indicates that
\begin{equation}\label{eq:mx2ge3halfquadrantinvariant}
D_{x_3}\Fc^{m(x_2)-2}\big(\Int\mathcal{Q}_{x_3}(\mathrm{I,III})\big)\subset\Int\big\{(d\phi,d\theta)\in \mathcal{Q}_{\hat{F}(x)}(\mathrm{I,III})\bigm|\frac{d\theta}{d\phi}\in[0,1]\big\}.
\end{equation}\newline In case (a1) of \eqref{eqMainCases}, \cref{caseA} says that $D_{x_0}\Fc^2$ is a negative matrix with $n_1=0$. Since $x$, $x_0$ are the collisions on the same arc $\Gamma_r$, $D_x\Fc=\begin{pmatrix}
1 & 2 \\
0 & 1 
\end{pmatrix}$. Since $x_2$, $x_3$ are the collisions on the same arc $\Gamma_r$, $D_{x_2}\Fc=\begin{pmatrix}
1 & 2 \\
0 & 1 
\end{pmatrix}$. Therefore, $D_x\Fc^{n_1+4}=D_x\Fc^4\overbracket{=}^{\mathrlap{\text{chain rule}}}D_{x_2}\Fc\cdot D_{x_0}\Fc^2\cdot D_x\Fc$ is a negative matrix. So $D_x\Fc^{n_1+4}\big(\mathcal{Q}_x(\mathrm{I,III})\big)=D_x\Fc^4\big(\mathcal{Q}_x(\mathrm{I,III})\big)\subset \Int\mathcal{Q}_{x_3}(\mathrm{I,III})$.
\newline In case (b) of \eqref{eqMainCases}, \cref{caseB} says that $D_x\Fc^4$ is a negative matrix with $n_1=0$, so $D_x\Fc^{n_1+4}\big(\mathcal{Q}_x(\mathrm{I,III})\big)=D_x\Fc^4\big(\mathcal{Q}_x(\mathrm{I,III})\big)\subset \Int\mathcal{Q}_{x_3}(\mathrm{I,III})$.\newline In case (c) of \eqref{eqMainCases}, \cref{caseC} says that $D_x\Fc^{4+n_1}$ is a positive matrix with $n_1\ge1$, so $D_x\Fc^{4+n_1}\big(\mathcal{Q}_x(\mathrm{I,III})\big)\subset \Int\mathcal{Q}_{x_3}(\mathrm{I,III})$.

Therefore, in cases (a1), (b) and (c) of \eqref{eqMainCases}, by the chain rule, the quadrant is strictly invariant. In fact, since $m(x_2)\ge3$, we have the image of the quadrant under $D_x\hat{F}$ is 
\begin{equation}\label{eq:strictlyinvcone_casea1bc}
    \begin{aligned}
        D_x\hat{F}(\mathcal{Q}_x(\mathrm{I,III}))&= D_x\Fc^{n_1+4+m(x_2)-2}(\mathcal{Q}_x(\mathrm{I,III}))\overbracket{\subset}^{\mathclap{\substack{D_x\Fc^{4+n_1}\big(\mathcal{Q}_x(\mathrm{I,III})\big)\subset \Int\mathcal{Q}_{x_3}(\mathrm{I,III})\text{ and chain rule}}}} D_{x_3}\Fc^{m(x_2)-2}(\Int\mathcal{Q}_{x_3}(\mathrm{I,III}))\\&\overbracket{\subset}^{\mathclap{\text{\eqref{eq:mx2ge3halfquadrantinvariant}}}}\Big\{\text{interior of }\big\{(d\phi,d\theta)\in \mathcal{Q}_{\hat{F}(x)}(\mathrm{I,III})\bigm|\frac{d\theta}{d\phi}\in[0,1]\big\}\Big\}\nfd\mathcal{HQ}_{\hat{F}(x)}(\mathrm{I,III})
    \end{aligned}
\end{equation}

Now we check case (a0) of \eqref{eqMainCases}. There are 4 subcases.

Subcase 1) $x\in\Fc^{-1}(\Mrout\setminus\Mrin)$ and $x_2\in\Mrout\cap\Mrin$, that is, $m(x_2)=0$.

In this subcase, $x_2=\hat{F}(x)$. Since $x_0=\Fc(x)$, and $p(x)$, $p(\Fc(x))$ are on the same boundary $\Gamma_r$, we have $D_{x}\Fc=\begin{pmatrix}
1 & 2 \\
0 & 1 
\end{pmatrix}$ in $\phi\theta$-coordinates. Thus, the quadrant is invariant under $D\Fc$. i.e. $D_x\Fc(\mathcal{Q}_x(\mathrm{I,III}))\subset\mathcal{Q}_{x_0}(\mathrm{I,III})$. By \cref{caseA}, $D_{x_0}\Fc^2$ is a negative matrix, $D_{x_0}\Fc^2(\mathcal{Q}_{x_0}(\mathrm{I,III}))\subset\Int\mathcal{Q}_{x_2}(\mathrm{I,III})=\Int\mathcal{Q}_{\hat{F}(x)}(\mathrm{I,III})$. Hence, by the chain rule, $D_x\Fc^3(\mathcal{Q}_x(\mathrm{I,III}))\subset D_{x_0}\Fc^2(\mathcal{Q}_{x_0}(\mathrm{I,III}))\subset\Int\mathcal{Q}_{x_2}(\mathrm{I,III})=\Int\mathcal{Q}_{\hat{F}(x)}(\mathrm{I,III})$.

Subcase 2) $x\in\Fc^{-1}(\Mrout\setminus\Mrin)$ and $x_2\notin\Mrout\cap\Mrin$, that is, $m(x_2)\ge1$. 

From $x_2$ to $\hat{F}(x)$ are the $m(x_2)$ collisions on the same arc, the differential $D_{x_2}\Fc^{m(x_2)-1}=\begin{pmatrix}
1 & 2(m(x_2)-2) \\
0 & 1 
\end{pmatrix}$ in the coordinates
$(\phi,\theta)$. Then for the same reasons as in subcase 1) and by \cref{caseA}. The quadrant is strictly invariant.

Subcase 3) $x\in\Mrout\cap\Mrin$ and $x_2\in\Mrout\cap\Mrin$, that is, $m(x_2)=0$. The \cref{caseA} immediately indicates that the quadrant is strictly invariant.

Subcase 4) $x\in\Mrout\cap\Mrin$ and $x_2\notin\Mrout\cap\Mrin$, that is, $m(x_2)\ge1$. The same reasons as in subcase 2) and by \cref{caseA} give that the quadrant is strictly invariant.
\end{proof}
\begin{remark}
    We will see that even the half-quadrant family $\big\{\frac{d\theta}{d\phi}\in[0,1]\big\}$ is strictly invariant (\cref{corollary:halfquadrantcone}).
\end{remark}

\subsection{Region of contraction and boundaries for (a0), (a1) (b) (c) cases regions in phase space}\label{subsec:regionofcontraction}\label{subsec:ProofContractionRegion}\hfill

We first give more detailed \cref{fig:x0a0a1bcboundary,fig:x1A0A1BCboundary,fig:Fn1(x1)A0A1BCboundary,fig:x2a0a1bcboundary} for the orbit segment $x_0\in\Mrout\xrightarrow{\Fc}x_1\longrightarrow\cdots\Fc^{n_1}(x_1)\xrightarrow{\Fc} x_2$ in the (a0)(a1)(b)(c) cases of \eqref{eqMainCases} for the $x_0$, $x_1$, $x_2$ defined in the return orbit segment \eqref{eqMhatOrbSeg}. Then we conclude this section with the previously advertised control of the region of contraction.

\begin{figure}[ht]
\begin{minipage}{.6\textwidth}
\begin{center}
\begin{tikzpicture}[xscale=0.75, yscale=0.75]
    \tkzDefPoint(0,0){Or};
    \tkzDefPoint(0,-4.8){Y};
    \pgfmathsetmacro{\rradius}{4.5};
    \pgfmathsetmacro{\bdist}{6.3};
    \tkzDefPoint(0,\bdist){OR};
    \pgfmathsetmacro{\phistardeg}{40};
    \pgfmathsetmacro{\XValueArc}{\rradius*sin(\phistardeg)};
    \pgfmathsetmacro{\YValueArc}{\rradius*cos(\phistardeg)};
    \tkzDefPoint(\XValueArc,-1.0*\YValueArc){B};    \tkzDefPoint(-1.0*\XValueArc,-1.0*\YValueArc){A};
    \pgfmathsetmacro{\Rradius}{veclen(\XValueArc,\YValueArc+\bdist)};
    \tkzDrawArc[name path=Cr,thick](Or,B)(A);
    \tkzDrawArc[name path=CR,thick](OR,A)(B);
    \tkzDrawArc[name path=Cr2,dashed ,very thin](Or,A)(B);
    \pgfmathsetmacro{\PXValue}{\Rradius*sin(7)};
    \pgfmathsetmacro{\PYValue}{\bdist-\Rradius*cos(7)};
    \tkzDefPoint(\PXValue,\PYValue){P};
    \tkzDefPoint(0,\bdist-\Rradius){C};
    \draw[ultra thin,dashed](OR)--(P);
    \draw[ultra thin,dashed](OR)--(A);
    \draw[ultra thin,dashed](OR)--(C);
     \draw[ultra thin,dashed](OR)--(B);
    \tkzDefPointBy[projection=onto OR--P](A)\tkzGetPoint{M};
    \coordinate (Q_far) at ($(A)!2!(M)$);
    \path[name path=linePQ] (P)--(Q_far);
    \path[name intersections={of=Cr and linePQ,by={Q}}];
    \tkzDefPointBy[rotation=center P angle 90](OR) \tkzGetPoint{P2}
    \draw[red,ultra thin,dashed](P)--($(P)!0.4!(P2)$);
    \draw[red,ultra thin,dashed](P)--($(P2)!1.35!(P)$);
    \coordinate (P1) at ($(P2)!1.35!(P)$);
    \tkzMarkRightAngle(OR,P,P1);
    \begin{scope}[ultra thin,decoration={
    markings,
    mark=at position 0.5 with {\arrow{>}}}
    ] 
    \draw[blue,postaction={decorate}] (P)--(Q);
    \end{scope}
    \tkzLabelPoint[above left](P){$P=p(y)$};
    \tkzLabelPoint[above right](A){$A$};
    \tkzLabelPoint[above right](Q){$Q=p(\Fc(y))$};
    \tkzLabelPoint[above left](B){$B$};
    \tkzLabelPoint[left](Or){$O_r$};
    \tkzLabelPoint[above](OR){$O_R$};
    \tkzLabelPoint[below](P1){$P_1$};
    \tkzMarkAngle[arc=lll,ultra thin,size=2.7,blue](P,OR,B);
     \node[blue] at ($(OR)+(0.9,-3.4)$)  {\Large $\Phistar-\Phi$};
     
    \node[] at ($(OR)+(3.2,0.4)$)  {\small $\measuredangle O_rO_RP=\Phi$, $\measuredangle AO_RB=2\Phistar$.};
    \node[] at ($(OR)+(3.2,-0.2)$) {\small $\measuredangle{PO_RB}=\Phistar-\Phi$};
    \node[] at ($(OR)+(3.5,-0.8)$)  {\small $\measuredangle P_1PB=0.5\measuredangle{PO_RB}=0.5(\Phistar-\Phi)$};
    \node[] at ($(OR)+(3.5,-1.4)$)  {\small $\measuredangle P_1PQ=\theta_1$};
    \path[name path=linePQ] (P)--($(Q)!2!(P)$);
    \path[name intersections={of=Cr2 and linePQ,by={Q0}}];
    \tkzLabelPoint[left](Q0){$Q_0$};
    \draw[blue,ultra thin,dashed](P)--(Q0);
    \draw[blue,ultra thin,dashed](P)--(B);
    \draw[ultra thin,dashed](Q0)--(B);
    \draw[blue, ultra thin,dashed](B)--(Q);
    \draw[ultra thin,dashed](B)--(Or);
    \draw[ultra thin,dashed](Q)--(Or);
    \tkzDrawPoints(A,B,P,Or,OR,Q,P1,Q0);
    
\end{tikzpicture}
\caption{$y=(\Phi,\theta_1)\in\MRout$ with $\theta_1<\sin^{-1}(\sqrt{4r/R})$ in cases (a1),(b),(c) of \cref{fig:Fn1(x1)A0A1BCboundary}.}\label{fig:caseA1BCContractionRegion}
\end{center}
\end{minipage}\hfill
\begin{minipage}{.4\textwidth}
    \begin{tikzpicture}[xscale=0.7, yscale=0.7]
        \pgfmathsetmacro{\PHISTAR}{1.5};
        \pgfmathsetmacro{\Dlt}{0.6};
        \pgfmathsetmacro{\Dltplus}{\Dlt+0.6}
        \pgfmathsetmacro{\WIDTH}{12};
        \tkzDefPoint(-\PHISTAR,0){LL};
        \tkzDefPoint(-\PHISTAR,\WIDTH){LU};
        \tkzDefPoint(\PHISTAR,0){RL};
        \tkzDefPoint(\PHISTAR,\WIDTH){RU};
        \tkzDefPoint(-\PHISTAR,\PHISTAR){X};
        \tkzDefPoint(\PHISTAR,\WIDTH-\PHISTAR){Y};
        \tkzDefPoint(-\PHISTAR,\WIDTH-\PHISTAR){IX};
        \tkzDefPoint(\PHISTAR,\PHISTAR){IY};
        \tkzDefPoint(-\PHISTAR,\PHISTAR+\Dlt){LLD};
        \tkzDefPoint(\PHISTAR,\PHISTAR+\Dlt){RLD};
        \tkzDefPoint(-\PHISTAR,\PHISTAR+\Dlt){LLD};
        \tkzDefPoint(\PHISTAR,\PHISTAR+\Dlt){RLD};
        \tkzDefPoint(-\PHISTAR,\WIDTH-\PHISTAR-\Dlt){LUD};
        \tkzDefPoint(\PHISTAR,\WIDTH-\PHISTAR-\Dlt){RUD};
        \tkzDefPoint(-\PHISTAR,\PHISTAR+\Dltplus){LLDp};
        \tkzDefPoint(\PHISTAR,\PHISTAR+\Dltplus){RLDp};
        \tkzDefPoint(-\PHISTAR,\WIDTH-\PHISTAR-\Dltplus){LUDp};
        \tkzDefPoint(\PHISTAR,\WIDTH-\PHISTAR-\Dltplus){RUDp};
        \tkzDefPoint(-\PHISTAR,\PHISTAR+\Dltplus+0.8){LLDdp};
        \tkzDefPoint(\PHISTAR,\PHISTAR+\Dltplus+0.8){RLDdp};
        \tkzDefPoint(-\PHISTAR,\WIDTH-\PHISTAR-\Dltplus-0.8){LUDdp};
        \tkzDefPoint(\PHISTAR,\WIDTH-\PHISTAR-\Dltplus-0.8){RUDdp};
        \tkzDefPoint (0.3*\PHISTAR,\PHISTAR+\Dltplus){montonestart};
        \tkzDefPoint (0.3*\PHISTAR,0.333*\WIDTH){montoneend};
        \draw [thin] (LL) --(LU);
        \draw [thin] (LU) --(RU);
        \draw [thin] (RL) --(RU);
        \draw [thin] (LL) --(RL);
        \draw [blue,ultra thin,dashed] (LL) --(IY);
        \draw [red,ultra thin,dashed]  (X)--(RL);
        \draw [blue,ultra thin,dashed] (RU) --(IX);
        \draw [red,ultra thin,dashed]  (LU)--(Y);
        \draw [purple,thick]  (LLD)--(RLD);
        \draw [purple,thick]  (LUD)--(RUD);
        \draw [green,ultra thin]  (LLDp)--(RLDp);
        \draw [green,ultra thin]  (LUDp)--(RUDp);
        \tkzLabelPoint[right](RU){$\pi$};
        \tkzLabelPoint[right](RL){$0$};
        \tkzLabelPoint[right](RLD){\small $\sin{\theta_1}=2r/R$, $d_1=2r$};
        \tkzLabelPoint[right](RUD){\small $\sin{\theta_1}=2r/R$, $d_1=2r$};
        \tkzLabelPoint[right](RLDp){\small $\sin{\theta_1}=\sqrt{4r/R}$};
        \tkzLabelPoint[right](RUDp){\small $\sin{\theta_1}=\sqrt{4r/R}$};
        \coordinate (RM) at ($(RL)!0.5!(RU)$);
        \tkzLabelPoint[right](RM){$\frac{\pi}{2}$};
        \coordinate (UpperM) at ($(LU)!0.5!(RU)$);
        \coordinate (LowerM) at ($(LL)!0.5!(RL)$);
        \coordinate (LM) at ($(LU)!0.5!(LL)$);
        \tkzLabelPoint[below](LL){$-\Phistar$};
        \tkzLabelPoint[below](RL){$+\Phistar$};
        \tkzLabelPoint[below](LowerM){$0$};
        
        \tkzDefPoint(0,0.5*\PHISTAR){LMend};
        \tkzDefPoint(0,\WIDTH-0.5*\PHISTAR){HMend};
        \tkzDefPoint(-\PHISTAR,0.5*\PHISTAR){thetamin};
        \tkzDefPoint(\PHISTAR,\WIDTH-0.5*\PHISTAR){thetamax};
        \node[] at ($(LMend)!0.7!(thetamin)$) {\small (c)};
        \node[] at ($(HMend)!0.7!(thetamax)$) {\small (c)};
        \fill [red, opacity=10/30] (RU) -- (HMend)--(Y)-- cycle;

        \fill [red, opacity=10/30] (LL) -- (LMend)--(X)-- cycle;
        \fill [green, opacity=10/30](LLD) -- (RLD) -- (IY) -- (LMend) -- (X) -- cycle;
        \fill [green, opacity=10/30](LUD) -- (RUD) -- (Y) -- (HMend) -- (IX) -- cycle;       
        
        \fill [blue, opacity=10/30](LLDp) -- (RLDp) -- (RLD) -- (LLD) -- cycle;

        \fill [blue, opacity=10/30](LUDp) -- (RUDp) -- (RUD) -- (LUD) -- cycle;

        \fill [orange, opacity=10/30](LUDp) -- (RUDp) -- (RLDp) -- (LLDp) -- cycle;
        
        \node[] at ($(RUDp)!0.5!(LUD)$) {\small (a1)};
        \node[] at ($(RLDp)!0.5!(LLD)$) {\small (a1)};

        \node[] at ($(LM)!0.5!(RM)$) {\small (a0)};
        
        \tkzDrawPoints(X,Y,IX,IY,LLD,LUD,RLD,RUD,LMend,HMend,thetamin,thetamax);
        \node[] at ([shift={(0,0.5*\PHISTAR)}]LMend) {\small (b)};
        \node[] at ([shift={(0,-0.5*\PHISTAR)}]HMend) {\small (b)};
    \end{tikzpicture}
    \caption{$\MRin$ divided for $x_1$ in cases (a0)(a1)(b)(c).}\label{fig:X1BirthdayCake}
\end{minipage}\hfill
\end{figure}

\begin{remark}\label{rmk:CaseA1Bcoordinates}[For \cref{fig:caseA1BCContractionRegion}]
In the coordinate system defined in \cref{def:StandardCoordinateTable} with $O_r=(0,0)$ and $O_R=(0,b)$, for $x_1=(\Phi_1,\theta_1)\in\MRout\cap\Mrin$ with $\theta_1<\sin^{-1}(\sqrt{4r/R})<\pi/2$ (for the lower (a1)(b) component in phase space \cref{fig:Fn1(x1)A0A1BCboundary}), we have the following coordinate in billiard table \cref{fig:caseA1BCContractionRegion}.
\begin{equation}\label{eq:newXqXt}
\begin{aligned}
Q\dfn p(x_2)& =(x_Q,y_Q)=(r\sin{\phi_2},-r\cos{\phi_2}),\\
P\dfn p(x_1) & =(x_P,y_P)=(R\sin{\Phi_1},b-R\cos{\Phi_1}).
\end{aligned}
\end{equation}    
\end{remark}

\begin{figure}[ht]
\begin{minipage}{0.75\textwidth}
\begin{center}    
\includegraphics[scale=2.1]{figs/X0InCaseA0A1BCBoundary.pdf}\caption{$x_0\in\Mrout$ near $(2\pi-\phistar,\phistar)$ in \eqref{eqMhatOrbSeg} for the cases (a0)(a1)(b)(c) in \eqref{eqMainCases}.}\label{fig:x0a0a1bcboundary}
\end{center}
\end{minipage}\hfill
\begin{minipage}{0.25\textwidth}
\begin{center}
    \begin{tikzpicture}[scale=0.42]
        \pgfmathsetmacro{\PHISTAR}{1.5};
        \pgfmathsetmacro{\Dlt}{0.6};
        \pgfmathsetmacro{\Dltplus}{\Dlt+0.6}
        \pgfmathsetmacro{\WIDTH}{12};
        \tkzDefPoint(-\PHISTAR,0){LL};
        \tkzDefPoint(-\PHISTAR,\WIDTH){LU};
        \tkzDefPoint(\PHISTAR,0){RL};
        \tkzDefPoint(\PHISTAR,\WIDTH){RU};
        \tkzDefPoint(-\PHISTAR,\PHISTAR){X};
        \tkzDefPoint(\PHISTAR,\WIDTH-\PHISTAR){Y};
        \tkzDefPoint(-\PHISTAR,\WIDTH-\PHISTAR){IX};
        \tkzDefPoint(\PHISTAR,\PHISTAR){IY};
        \tkzDefPoint(-\PHISTAR,\PHISTAR+\Dlt){LLD};
        \tkzDefPoint(\PHISTAR,\PHISTAR+\Dlt){RLD};
        \tkzDefPoint(-\PHISTAR,\PHISTAR+\Dlt){LLD};
        \tkzDefPoint(\PHISTAR,\PHISTAR+\Dlt){RLD};
        \tkzDefPoint(-\PHISTAR,\WIDTH-\PHISTAR-\Dlt){LUD};
        \tkzDefPoint(\PHISTAR,\WIDTH-\PHISTAR-\Dlt){RUD};
        \tkzDefPoint(-\PHISTAR,\PHISTAR+\Dltplus){LLDp};
        \tkzDefPoint(\PHISTAR,\PHISTAR+\Dltplus){RLDp};
        \tkzDefPoint(-\PHISTAR,\WIDTH-\PHISTAR-\Dltplus){LUDp};
        \tkzDefPoint(\PHISTAR,\WIDTH-\PHISTAR-\Dltplus){RUDp};
        \tkzDefPoint(-\PHISTAR,\PHISTAR+\Dltplus+0.8){LLDdp};
        \tkzDefPoint(\PHISTAR,\PHISTAR+\Dltplus+0.8){RLDdp};
        \tkzDefPoint(-\PHISTAR,\WIDTH-\PHISTAR-\Dltplus-0.8){LUDdp};
        \tkzDefPoint(\PHISTAR,\WIDTH-\PHISTAR-\Dltplus-0.8){RUDdp};
        \tkzDefPoint (0.3*\PHISTAR,\PHISTAR+\Dltplus){montonestart};
        \tkzDefPoint (0.3*\PHISTAR,0.333*\WIDTH){montoneend};
        \draw [ultra thin, dashed] (X) --(LU);
        \draw [ultra thin, dashed] (LU) --(RU);
        \draw [ultra thin, dashed] (RL) --(Y);
        \draw [ultra thin, dashed] (LL) --(RL);
        \draw [very thick] (LL) --(IY);
        \draw [very thick] (LL) --(IX);
        \draw [green,ultra thin,dashed]  (X)--(RL);
        \draw [very thick] (RU) --(IX);
        \draw [very thick] (RU) --(IY);
        \draw [green,ultra thin,dashed]  (LU)--(Y);
        \draw [red,thick]  (LLD)--(RLD);
        \draw [red,thick]  (LUD)--(RUD);
        \draw [blue,thick]  (LLDp)--(RLDp);
        \draw [blue,thick]  (LUDp)--(RUDp);
        \tkzLabelPoint[right](RU){$\pi$};
        \tkzLabelPoint[right](RL){$0$};
        \tkzLabelPoint[right](RLD){\small $\sin{\theta_1}=\frac{2r}{R}$};
        \tkzLabelPoint[right](RUD){\small $\sin{\theta_1}=\frac{2r}{R}$};
        \tkzLabelPoint[above right](RLDp){\small $\sin{\theta_1}=\sqrt{\frac{4r}{R}}$};
        \tkzLabelPoint[below right](RUDp){\small $\sin{\theta_1}=\sqrt{\frac{4r}{R}}$};
        \coordinate (RM) at ($(RL)!0.5!(RU)$);
        \tkzLabelPoint[right](RM){\Large $\frac{\pi}{2}\xrightarrow{\Fc^{n_1}}$};
        \coordinate (LMid) at ($(LL)!0.5!(LU)$);
        \tkzLabelPoint[left](LMid){\small$\xrightarrow{\Fc}$};
        \coordinate (UpperM) at ($(LU)!0.5!(RU)$);
        \coordinate (LowerM) at ($(LL)!0.5!(RL)$);
        \coordinate (LM) at ($(LU)!0.5!(LL)$);
        \tkzLabelPoint[below](LL){$-\Phistar$};
        \tkzLabelPoint[below](RL){$+\Phistar$};
        \tkzLabelPoint[below](LowerM){$0$};
        
        \tkzDefPoint(0,0.5*\PHISTAR){LMend};
        \tkzDefPoint(0,\WIDTH-0.5*\PHISTAR){HMend};
        \tkzDefPoint(-\PHISTAR,0.5*\PHISTAR){thetamin};
        \tkzDefPoint(\PHISTAR,\WIDTH-0.5*\PHISTAR){thetamax};
        \node[] at ($(LMend)!0.7!(thetamin)$) {\small (c)};
        \node[] at ($(HMend)!0.7!(thetamax)$) {\small (c)};

        \draw [green,thick]  (HMend)--(Y);
        \draw [green,thick]  (LMend)--(X);

       
        


        
        \node[] at ($(RUDp)!0.5!(LUD)$) {\small (a1)};
        \node[] at ($(RLDp)!0.5!(LLD)$) {\small (a1)};

        \node[] at ($(LM)!0.5!(RM)$) {\small $\mathfrak{U}$: (a0)};
        
        \tkzDrawPoints(X,Y,IX,IY,LLD,LUD,RLD,RUD,LMend,HMend,RLDp,RUDp,LLDp,LUDp);
        \node[] at ([shift={(0,0.5*\PHISTAR)}]LMend) {\small (b)};
        \node[] at ([shift={(0,-0.5*\PHISTAR)}]HMend) {\small (b)};
    \end{tikzpicture}\caption{\small $x_1\in\MRin$ for the cases (a0)(a1)(b)(c) in \eqref{eqMainCases}.}\label{fig:x1A0A1BCboundary}
\end{center}
\end{minipage}\hfill
\end{figure}
\begin{figure}[ht]
\begin{minipage}{0.25\textwidth}
\begin{center}
    \begin{tikzpicture}[scale=0.37]
        \pgfmathsetmacro{\PHISTAR}{1.5};
        \pgfmathsetmacro{\Dlt}{0.6};
        \pgfmathsetmacro{\Dltplus}{\Dlt+0.6}
        \pgfmathsetmacro{\WIDTH}{12};
        \tkzDefPoint(-\PHISTAR,0){LL};
        \tkzDefPoint(-\PHISTAR,\WIDTH){LU};
        \tkzDefPoint(\PHISTAR,0){RL};
        \tkzDefPoint(\PHISTAR,\WIDTH){RU};
        \tkzDefPoint(-\PHISTAR,\PHISTAR){X};
        \tkzDefPoint(\PHISTAR,\WIDTH-\PHISTAR){Y};
        \tkzDefPoint(-\PHISTAR,\WIDTH-\PHISTAR){IX};
        \tkzDefPoint(\PHISTAR,\PHISTAR){IY};
        \tkzDefPoint(-\PHISTAR,\PHISTAR+\Dlt){LLD};
        \tkzDefPoint(\PHISTAR,\PHISTAR+\Dlt){RLD};
        \tkzDefPoint(-\PHISTAR,\PHISTAR+\Dlt){LLD};
        \tkzDefPoint(\PHISTAR,\PHISTAR+\Dlt){RLD};
        \tkzDefPoint(-\PHISTAR,\WIDTH-\PHISTAR-\Dlt){LUD};
        \tkzDefPoint(\PHISTAR,\WIDTH-\PHISTAR-\Dlt){RUD};
        \tkzDefPoint(-\PHISTAR,\PHISTAR+\Dltplus){LLDp};
        \tkzDefPoint(\PHISTAR,\PHISTAR+\Dltplus){RLDp};
        \tkzDefPoint(-\PHISTAR,\WIDTH-\PHISTAR-\Dltplus){LUDp};
        \tkzDefPoint(\PHISTAR,\WIDTH-\PHISTAR-\Dltplus){RUDp};
        \tkzDefPoint(-\PHISTAR,\PHISTAR+\Dltplus+0.8){LLDdp};
        \tkzDefPoint(\PHISTAR,\PHISTAR+\Dltplus+0.8){RLDdp};
        \tkzDefPoint(-\PHISTAR,\WIDTH-\PHISTAR-\Dltplus-0.8){LUDdp};
        \tkzDefPoint(\PHISTAR,\WIDTH-\PHISTAR-\Dltplus-0.8){RUDdp};
        \tkzDefPoint (0.3*\PHISTAR,\PHISTAR+\Dltplus){montonestart};
        \tkzDefPoint (0.3*\PHISTAR,0.333*\WIDTH){montoneend};
        \draw [thin] (LL) --(LU);
        \draw [thin] (LU) --(RU);
        \draw [thin] (RL) --(RU);
        \draw [thin] (LL) --(RL);
        \draw [thick]  (X)--(RL);
        \draw [thick]  (LU)--(Y);
        \draw [red,thick]  (LLD)--(RLD);
        \draw [red,thick]  (LUD)--(RUD);
        \draw [blue,thick]  (LLDp)--(RLDp);
        \draw [blue,thick]  (LUDp)--(RUDp);
        \tkzLabelPoint[right](RU){$\pi$};
        \tkzLabelPoint[right](RL){$0$};
        \tkzLabelPoint[left](LLD){\small $\sin{\theta_1}=2r/R$};
        \tkzLabelPoint[left](LUD){\small $\sin{\theta_1}=2r/R$};
        \tkzLabelPoint[above left](LLDp){\small $\sin{\theta_1}=\sqrt{4r/R}$};
        \tkzLabelPoint[below left](LUDp){\small $\sin{\theta_1}=\sqrt{4r/R}$};
        \coordinate (RM) at ($(RL)!0.5!(RU)$);
        \tkzLabelPoint[right](RM){\Large $\frac{\pi}{2}$ $\xrightarrow{\Fc}$};
        \coordinate (UpperM) at ($(LU)!0.5!(RU)$);
        \coordinate (LowerM) at ($(LL)!0.5!(RL)$);
        \coordinate (LM) at ($(LU)!0.5!(LL)$);
        \tkzLabelPoint[below](LL){$-\Phistar$};
        \tkzLabelPoint[below](RL){$+\Phistar$};
        \tkzLabelPoint[below](LowerM){$0$};
        
        \tkzDefPoint(0,0.5*\PHISTAR){LMend};
        \tkzDefPoint(0,\WIDTH-0.5*\PHISTAR){HMend};
        \tkzDefPoint(\PHISTAR,0.5*\PHISTAR){thetamin};
        \tkzDefPoint(-\PHISTAR,\WIDTH-0.5*\PHISTAR){thetamax};
        \node[] at ($(LMend)!0.7!(thetamin)$) {\small (c)};
        \node[] at ($(HMend)!0.7!(thetamax)$) {\small (c)};

        \node[] at ($(RUDp)!0.5!(LUD)$) {\small $(a1)$};
        \node[] at ($(RLDp)!0.5!(LLD)$) {\small (a1)};

        \node[] at ($(LM)!0.5!(RM)$) {\small $\mathfrak{U}$: (a0)};

        \tkzLabelPoint[left](LM){\small $\xrightarrow{\Fc^{n_1}=\begin{pmatrix}1 & 2n_1\\ 0 & 1\end{pmatrix}}$};
        
        \draw [orange,thick]  (HMend)--(IX);
        \draw [orange,thick]  (LMend)--(IY);

       
        


        

        \fill [red, opacity=2/30](LU) -- (Y) -- (RL) -- (X) -- (LU) -- cycle;
        
        \node[] at ([shift={(0,0.5*\PHISTAR)}]LMend) {\small (b)};
        \node[] at ([shift={(0,-0.5*\PHISTAR)}]HMend) {\small (b)};
    \end{tikzpicture}
    \caption{\small $\Fc^{n_1}(x_1)\in\MRout$. Note in the cases (a0)(a1)(b) of \eqref{eqMainCases} $n_1=0$.}\label{fig:Fn1(x1)A0A1BCboundary}
\end{center}
\end{minipage}\hfill
\begin{minipage}{0.75\textwidth}
\begin{center}    
\includegraphics[scale=1.71]{figs/X2InCaseA0A1BCBoundary.pdf}\caption{$x_2\in\Mrin$ near $(\phistar,\phistar)$ in \eqref{eqMhatOrbSeg} for the cases (a0)(a1)(b)(c) in \eqref{eqMainCases}.}\label{fig:x2a0a1bcboundary}\hfill
\end{center}
\end{minipage}
\end{figure}

\begin{proposition}[\Nin size]\label{proposition:NinsizeInA1BCcases}
    Suppose that our lemon billiard satisfies \eqref{eqJZHypCond}. Let $y=(\Phi,\theta_1)\in\MRout$ and with $\sin{\theta_1}<\sqrt{4r/R}$, that is, $y$ is in the (a1)(b)(c) regions in \cref{fig:Fn1(x1)A0A1BCboundary}, then we have the following.
    \begin{enumerate}
        \item\label{item01NinsizeInA1BCcases} If $\theta_1<\pi/2$, that is, $y$ is in the lower component of (a1)(b)(c) $\MRout$ regions in \cref{fig:Fn1(x1)A0A1BCboundary}, then $(\phi_2,\theta_2)\nfd \Fc(y)\in \Mrin$ has $\mid\phi_2-\phistar\mid<2.001\sqrt{4r/R}$ and $\mid\theta_2-\phistar\mid<1.002\sqrt{4r/R}$.
        \item\label{item02NinsizeInA1BCcases} If $\theta_1>\pi/2$, that is, $y$ is in the upper component of (a1)(b)(c) $\MRout$ regions in \cref{fig:Fn1(x1)A0A1BCboundary}, then $(\phi_2,\theta_2)\nfd \Fc(y)\in \Mrin$ has $\mid\phi_2-(2\pi-\phistar)\mid<2.001\sqrt{4r/R}$ and $\mid\theta_2-(\pi-\phistar)\mid<1.002\sqrt{4r/R}$.
    \end{enumerate}
\end{proposition}
\begin{proof}
We note that if $y$ is in the upper component of the (a1)(b)(c) $\MRout$ regions in \cref{fig:Fn1(x1)A0A1BCboundary}), then by \cref{DefSymmetries} $I\circ J(y)$ is in the lower component of the (a1)(b)(c) $\MRout$ regions in \cref{fig:Fn1(x1)A0A1BCboundary}) and $\Fc(I\circ J(y))=I\circ J (\Fc(y))$. Therefore, \eqref{item02NinsizeInA1BCcases} follows \eqref{item01NinsizeInA1BCcases}. We only need to prove \eqref{item01NinsizeInA1BCcases}. 

In the billiard table \cref{fig:caseA1BCContractionRegion} using the coordinate system of \cref{def:StandardCoordinateTable} with $O_r=(0,0)$ and $O_R=(0,b)$, we have $P\dfn p(y)\in\Gamma_R$, $Q\dfn p(\Fc(y))=(r\sin\phi_2,-r\cos\phi_2)\in\Gamma_r$. Given counterclockwise orientation on $\Gamma_r$, by the definition of collision angle $\theta_2$ we have $\theta_2=\measuredangle{(T_Q\Gamma_r,\overrightarrow{PQ})}$. Note that in the billiard table, $\overrightarrow{PQ}$ direction angle is $\Phi+\theta_1$ and $T_Q\Gamma_r$ the tangent vector of $\Gamma_r$ at $Q$ has the direction angle being $\phi_2$. Hence $\theta_2=\phi_2-(\Phi+\theta_1)$.

Note that in \cref{fig:caseA1BCContractionRegion} we suppose that $PP_1$ is the tangential line at $Q$ of $\Gamma_R$ and $Q_0$ is the intersection of line $PQ$ and circle $C_r$. Then, using elementary geometry, we can get the following: \[\measuredangle P_1PB=0.5\measuredangle{PO_RB}=0.5(\Phistar-\Phi)\text{, and } \measuredangle P_1PQ=\theta_1.\] Therefore, in \cref{fig:caseA1BCContractionRegion} we get $\measuredangle BPQ=\measuredangle P_1PQ-\measuredangle P_1PB=\theta_1-0.5(\Phistar-\Phi)$. Using elementary geometry again, we get \begin{equation}\label{eq:observation1(a1)(b)(c)region}
    0<\phi_2-\phistar=\measuredangle{BO_rQ}=2\measuredangle{BQ_0Q}=2(\measuredangle{BPQ}-\overbracket{\measuredangle{Q_0BP}}^{>0})<2\measuredangle{BPQ}=2(\theta_1-0.5\Phistar+0.5\Phi).
\end{equation}

Also, we see that for $(\Phi,\theta_1)=y\in\MRout$ in \cref{fig:Fn1(x1)A0A1BCboundary}, $\theta_1-(0.5\Phistar-0.5\Phi)>0$  always holds since $\big\{(\Phi,\theta_1)\mid\theta_1=0.5\Phistar-0.5\Phi\big\}$ is the line segment as the lower bound of $\MRout$. Together with $\theta_2=\phi_2-(\Phi+\theta_1)$, we get the following inequalities. \begin{align*}
    -1.002\sqrt{4r/R}&\underbracket{<}_{\mathrlap{R>1700r:\, r/R+\sqrt{4r/R}<1.00125\sqrt{4r/R}}}-1.0005(r/R+\sqrt{4r/R})\overbracket{<}^{\mathclap{\text{since \eqref{eqJZHypCond0}: }R>1700r\text{ and }\arcsin \text{ is convex on }[0,1]}}-\big[\sin^{-1}(r\sin{\phistar}/R)+\sin^{-1}(\sqrt{4r/R})\big]\\
    &\underbracket{<}_{\mathllap{\Phi\in(-\Phistar,\Phistar),\, R\sin\Phistar=r\sin\phistar}}-(\Phi+\theta_1)\underbracket{<}_{\mathrlap{\phi_2-\phistar>0}}\theta_2-\phistar=(\phi_2-\phistar)-(\Phi+\theta_1)\\
    &\underbracket{<}_{\mathclap{\eqref{eq:observation1(a1)(b)(c)region}}}2(\theta_1-0.5\Phistar+0.5\Phi)-(\Phi+\theta_1)=\theta_1-\Phistar<\theta_1<\sin^{-1}(\sqrt{4r/R})\underbracket{<}_{\mathclap{\text{since \eqref{eqJZHypCond0}: }R>1700r\text{ and }\arcsin \text{ is convex on }[0,1]}}1.0005\sqrt{4r/R}.
\end{align*}

Then we further notice that for $(\Phi,\theta_1)=y\in\MRout$ with $\theta_1<\sin^{-1}(\sqrt{4r/R})$ in (a1)(b)(c) lower component regions in \cref{fig:Fn1(x1)A0A1BCboundary}, $\theta_1+0.5\Phi<\sin^{-1}(\sqrt{4r/R})+0.5\Phistar$. Therefore \eqref{eq:observation1(a1)(b)(c)region} gives \begin{equation}\label{eq:observation2(a1)(b)(c)region}
0<\phi_2-\phistar<2\sin^{-1}{(\sqrt{4r/R})}<2\cdot 1.0005\sqrt{4r/R}.\end{equation}
Therefore, we have $\mid\phi_2-\phistar\mid<2.001\sqrt{4r/R}$ and $\mid\theta_2-\phistar\mid<1.002\sqrt{4r/R}$.
\end{proof}

\begin{proposition}[Set size for (b)(c) cases inside \Nin ]\label{proposition:NinsizeInBCcases}
    Suppose that our lemon billiard satisfies \eqref{eqJZHypCond}. Let $y=(\Phi,\theta_1)\in\MRout$ and with $\sin{\theta_1}<2r/R$, that is, $y$ is in the (b)(c) regions in \cref{fig:Fn1(x1)A0A1BCboundary}, then we have the following.
    \begin{enumerate}
        \item\label{item01NinsizeInBCcases} If $\theta_1<\pi/2$, that is, $y$ is in the lower component of (b)(c) $\MRout$ regions in \cref{fig:Fn1(x1)A0A1BCboundary}, then $(\phi_2,\theta_2)\nfd \Fc(y)\in \Mrin$ has $\mid\phi_2-\phistar\mid<2.001\cdot2r/R$ and $\mid\theta_2-\phistar\mid<1.002\cdot2r/R$.
        \item\label{item02NinsizeInBCcases} If $\theta_1>\pi/2$, that is, $y$ is in the upper component of (b)(c) $\MRout$ regions in \cref{fig:Fn1(x1)A0A1BCboundary}, then $(\phi_2,\theta_2)\nfd \Fc(y)\in \Mrin$ has $\mid\phi_2-(2\pi-\phistar)\mid<2.001\cdot2r/R$ and $\mid\theta_2-(\pi-\phistar)\mid<1.002\cdot2r/R$.
    \end{enumerate}
\end{proposition}
\begin{proof}
    The proof is the same as \cref{proposition:NinsizeInA1BCcases}.
\end{proof}

\begin{proposition}[\Nout size]\label{proposition:NoutsizeInA1BCcases}
    Suppose that our lemon billiard satisfies \eqref{eqJZHypCond}. Let $y=(\Phi,\theta_1)\in\MRin$ and with $\sin{\theta_1}<\sqrt{4r/R}$, that is, $y$ is in the (a1)(b)(c) regions in \cref{fig:x1A0A1BCboundary}, then we have the following.
    \begin{enumerate}
        \item\label{item01NoutsizeInA1BCcases} If $\theta_1<\pi/2$, that is, $y$ is in the lower component of (a1)(b)(c) $\MRin$ regions in \cref{fig:x1A0A1BCboundary}, then $(\phi_0,\theta_0)\nfd \Fc^{-1}(y)\in \Mrout$ has $\mid\phi_0-(2\pi-\phistar)\mid<2.001\sqrt{4r/R}$ and $\mid\theta_0-\phistar\mid<1.002\sqrt{4r/R}$.
         \item\label{item02NoutsizeInA1BCcases} If $\theta_1>\pi/2$, that is, $y$ is in the upper component of (a1)(b)(c) $\MRin$ regions in \cref{fig:x1A0A1BCboundary}, then $(\phi_0,\theta_0)\nfd \Fc^{-1}(y)\in \Mrout$ has $\mid\phi_0-\phistar\mid<2.001\sqrt{4r/R}$ and $\mid\theta_0-(\pi-\phistar)\mid<1.002\sqrt{4r/R}$.
    \end{enumerate}
\end{proposition}
\begin{proof}
If $\theta_1<\pi/2$, that is, $y$ is in the lower component of (a1)(b)(c) $\MRin$ regions in \cref{fig:x1A0A1BCboundary}), then $J(y)$ is in the lower component of (a1)(b)(c) $\MRout$ according to the definition of symmetry $J$ in \cref{DefSymmetries}. Again by \cref{DefSymmetries}, we know $J^{-1}(\Fc(J(y)))=\Fc^{-1}(y)$. Therefore, \cref{proposition:NinsizeInA1BCcases}\eqref{item01NinsizeInA1BCcases} implies that \[\Fc^{-1}(y)=J^{-1}\Fc(Jy)\overbracket{\in}^{\mathclap{\text{\cref{proposition:NinsizeInA1BCcases}\eqref{item01NinsizeInA1BCcases}}}}J^{-1}\big(\big\{(\phi,\theta)\bigm|\mid\phi-\phistar\mid<2.001\sqrt{4r/R}\text{ and }\mid\theta-\phistar\mid<1.002\sqrt{4r/R}\big\}\big).\] 

By $J^{-1}=J$ according to the definition of symmetry $J$ in \cref{DefSymmetries} we get $(\phi_0,\theta_0)\nfd \Fc^{-1}(y)\in \Mrout$ has $\mid\phi_0-(2\pi-\phistar)\mid<2.001\sqrt{4r/R}$ and $\mid\theta_0-\phistar\mid<1.002\sqrt{4r/R}$. This proves \eqref{item01NoutsizeInA1BCcases}.

If $\theta_1>\pi/2$, that is, $y$ is in the upper component of (a1)(b)(c) $\MRin$ regions in \cref{fig:x1A0A1BCboundary}, then $I(y)$ is in the lower component of (a1)(b)(c) $\MRout$ according to the definition of symmetry $I$ in \cref{DefSymmetries}. Again by \cref{DefSymmetries}, we know $I^{-1}(\Fc(I(y)))=\Fc^{-1}(y)$. Therefore, \cref{proposition:NinsizeInA1BCcases}\eqref{item01NinsizeInA1BCcases} implies that \[\Fc^{-1}(y)=I^{-1}\Fc(Iy)\overbracket{\in}^{\mathclap{\text{\cref{proposition:NinsizeInA1BCcases}\eqref{item01NinsizeInA1BCcases}}}}I^{-1}\big(\big\{(\phi,\theta)\bigm|\mid\phi-\phistar\mid<2.001\sqrt{4r/R}\text{ and }\mid\theta-\phistar\mid<1.002\sqrt{4r/R}\big\}\big).\]

By $I^{-1}=I$ according to the definition of symmetry $I$ in \cref{DefSymmetries} we get $(\phi_0,\theta_0)\nfd \Fc^{-1}(y)\in \Mrout$ has $\mid\phi_0-\phistar\mid<2.001\sqrt{4r/R}$ and $\mid\theta_0-(\pi-\phistar)\mid<1.002\sqrt{4r/R}$. This proves \eqref{item02NoutsizeInA1BCcases}.
\end{proof}

\begin{proposition}[Set size for (b)(c) cases inside \Nout ]\label{proposition:NoutsizeInBCcases}
    Suppose that our lemon billiard satisfies \eqref{eqJZHypCond}. Let $y=(\Phi,\theta_1)\in\MRin$ and with $\sin{\theta_1}<\sqrt{4r/R}$, that is, $y$ is in the (b)(c) regions in \cref{fig:x1A0A1BCboundary}, then we have the following.
    \begin{enumerate}
        \item\label{item01NoutsizeInBCcases} If $\theta_1<\pi/2$, that is, $y$ is in the lower component of (b)(c) $\MRin$ regions in \cref{fig:x1A0A1BCboundary}, then $(\phi_0,\theta_0)\nfd \Fc^{-1}(y)\in \Mrout$ has $\mid\phi_0-(2\pi-\phistar)\mid<2.001\cdot2r/R$ and $\mid\theta_0-\phistar\mid<1.002\cdot 2r/R$.
         \item\label{item02NoutsizeInBCcases} If $\theta_1>\pi/2$, that is, $y$ is in the upper component of (b)(c) $\MRin$ regions in \cref{fig:x1A0A1BCboundary}, then $(\phi_0,\theta_0)\nfd \Fc^{-1}(y)\in \Mrout$ has $\mid\phi_0-\phistar\mid<2.001\cdot 2r/R$ and $\mid\theta_0-(\pi-\phistar)\mid<1.002\cdot 2r/R$.
    \end{enumerate}
\end{proposition}
\begin{proof}
    The proof reasoning is the  same as \cref{proposition:NoutsizeInA1BCcases} due to the $I(y)$ $J(y)$ being the $y$ in \cref{proposition:NinsizeInBCcases}.
\end{proof}

Now from \cref{proposition:NinsizeInA1BCcases,proposition:NoutsizeInA1BCcases}, we get the following corollary.

\begin{corollary}[Contraction region]\label{contraction_region}
Cases {(a1) (\Cref{TEcase_a1}), (b) (\Cref{TEcase_b}), and  (c) (\Cref{TEcase_c}}) in \eqref{eqMainCases} arise only when $x_0\in\Nout$ and $x_2\in\Nin$ (see \eqref{eqPtsFromeqMhatOrbSeg} and \cref{def:SectionSetsDefs}).
\end{corollary}
\begin{proof}
    We have \(x_0=(\phi_0,\theta_0)\) , \(x_1=(\Phi,\theta_1)\) and \(x_2=(\phi_2,\theta_2)\) from \eqref{eqPtsFromeqMhatOrbSeg}. With $x_*=(2\pi-\phistar,\phistar),\,y_*=(\phistar,\pi-\phistar),\,Ix_*=(2\pi-\phistar,\pi-\phistar),\,Iy_*=(\phistar,\phistar)$ given in \cref{def:SectionSetsDefs,fig:Mr}, by conclusions of \cref{proposition:NoutsizeInA1BCcases,proposition:NinsizeInA1BCcases} we have get following.
    
    In cases (a1)(b)(c), if $\theta_1<\sin^{-1}(\sqrt{4r/R})$, then 
    \begin{equation}\label{eq:caseA1BCNinNout1}
        \begin{aligned}
    \|x_2-Iy_*\|=&\sqrt{(\phi_2-\phistar)^2+(\theta_2-\phistar)^2}\underbracket{<}_{\mathclap{\text{\cref{proposition:NinsizeInA1BCcases}\eqref{item01NinsizeInA1BCcases}}}}\sqrt{4.002^2+2.004^2}\sqrt{r/R}<4.5\sqrt{r/R}=\delta_*,\\
    \text{ and }\|x_0-x_*\|=&\sqrt{\big(\phi_0-(2\pi-\phistar)\big)^2+(\theta_0-\phistar)^2}\underbracket{<}_{\mathllap{\text{\cref{proposition:NoutsizeInA1BCcases}\eqref{item01NoutsizeInA1BCcases}}}}\sqrt{4.002^2+2.004^2}\sqrt{r/R}<4.5\sqrt{r/R}=\delta_*.
    \end{aligned}
    \end{equation}

    In cases (a1)(b)(c), if $\theta_1>\pi-\sin^{-1}
    (\sqrt{4r/R})$, then 
    \begin{equation}\label{eq:caseA1BCNinNout2}
        \begin{aligned}
    \|x_2-Ix_*\|=&\sqrt{(\phi_2-(2\pi-\phistar))^2+(\theta_2-(\pi-\phistar))^2}\underbracket{<}_{\mathclap{\text{\cref{proposition:NinsizeInA1BCcases}\eqref{item02NinsizeInA1BCcases}}}}\sqrt{4.002^2+2.004^2}\sqrt{r/R}<4.5\sqrt{r/R}=\delta_*,\\
    \text{ and }\|x_0-y_*\|=&\sqrt{\big(\phi_0-\phistar\big)^2+(\theta_0-(\pi-\phistar))^2}\underbracket{<}_{\mathllap{\text{\cref{proposition:NoutsizeInA1BCcases}\eqref{item02NoutsizeInA1BCcases}}}}\sqrt{4.002^2+2.004^2}\sqrt{r/R}<4.5\sqrt{r/R}=\delta_*.
    \end{aligned}
    \end{equation}

We get \(x_0\in\Nout=\big(B(x_*,\delta_*)\cup B(y_*,\delta_*)\big)\cap \Mrout\)
and \(x_2\in\Nin=\big(B(Ix_*,\delta_*)\cup B(Iy_*,\delta_*)\big)\cap \Mrin\), 
where
$B(\cdot,\delta_*)$ are the neighborhoods of radius $\delta_*={4.5}\sqrt{r/R}$.
\end{proof}

\begin{remark}[Boundaries for (a0)(a1)(b)(c) in \cref{fig:x0a0a1bcboundary,fig:x1A0A1BCboundary,fig:Fn1(x1)A0A1BCboundary,fig:x2a0a1bcboundary} and optimal $\delta_*$]\label{rmk:RegionContractionA0A1Boundary}\hfill

Note that our choice of $\delta_*=4.5\sqrt{r/R}$ is not the smallest radius to suffice the $\delta_*$-neighborhoods of $x_*$, $y_*$, $Ix_*$, $Iy_*$ to include the $x_0$, $x_2$ in (a1)(b)(c) cases of \eqref{eqMainCases} orbit segments. 

The optimal $\delta_*$ is $\sin^{-1}(2\sqrt{r/R})+\sin^{-1}(r\sin\phistar/R)$ which can be seen in the following \cref{proposition:RegionContractionA0A1Boundary}. In this article we will not give a rigorous proof for \cref{proposition:RegionContractionA0A1Boundary} and will not use those facts of \cref{proposition:RegionContractionA0A1Boundary}, \eqref{item01ContractionRegion}, \eqref{item02ContractionRegion}, \eqref{item03ContractionRegion}, \eqref{item04ContractionRegion}, \eqref{item05ContractionRegion}, \eqref{item06ContractionRegion}. Readers will find that its proof comes from the tools/ingredients in \cref{subsec:Monotonicity}. 

In \cref{fig:x0a0a1bcboundary,fig:x2a0a1bcboundary}, the the arch singularity curve (in \cref{proposition:RegionContractionA0A1Boundary}\eqref{item05ContractionRegion}) as the marked boundary between case (b) and (c) for the $x_2$ region near $(\phistar,\phistar)$, and the arch singularity curve (in \cref{proposition:RegionContractionA0A1Boundary}\eqref{item06ContractionRegion}) as the marked boundary between case (b) and (c) for the $x_0$ region near $(2\pi-\phistar,\phistar)$ are to be discussed in next article. 

\cref{proposition:RegionContractionA0A1Boundary}\eqref{item03ContractionRegion} improves \cite[equation (3.14)]{hoalb}: $\|x_0-x_*\|<14r/(R\sin\phistar)$ to $<\sin^{-1}(2r/R)+\sin^{-1}(r\sin\phistar/R)$. \cref{proposition:RegionContractionA0A1Boundary}\eqref{item05ContractionRegion} improves \cite[equation (3.21)] {hoalb}: $\|x_0-x_*\|<5.84(r\sin\phistar)/R$ to $<2\sin^{-1}(r\sin\phistar/R)$.

Since these curves in \cref{fig:x0a0a1bcboundary,fig:x2a0a1bcboundary,fig:x1A0A1BCboundary,fig:Fn1(x1)A0A1BCboundary} marked as boundaries between case (b) (c) regions and boundaries between case (a1) (b) regions are well discussed in \cref{proposition:RegionContractionA0A1Boundary}, in possible future work readers can use these facts to improve our this paper estimates of $R_{\text{\upshape HF}}$ and $R_{\text{\upshape JZ}}$ in \cite{hoalb}. And we conjecture and believe our computations/conclusions can be generalized to apply to the lemon billiards with $\tan^{-1}(1/3)\le\phistar<\pi/2$.

The python code to generate \cref{fig:x0a0a1bcboundary,fig:x2a0a1bcboundary} is in \cite{LemonBilliardsCode}.
\end{remark}

\begin{proposition}\label{proposition:RegionContractionA0A1Boundary}Suppose that our lemon billiards satisfies \eqref{eqJZHypCond} and the straight line segment boundary of (a0) and (a1): $\mathfrak{L}_1\dfn\big\{(\Phi,\theta_1)\in M_R\bigm|\theta_1=\sin^{-1}(\sqrt{4r/R}),\Phi\in(-\Phistar,+\Phistar)\big\}$ and boundary of (a1) and (b): $\mathfrak{L}_2\dfn\big\{(\Phi,\theta_1)\in M_R\bigm|\theta_1=\sin^{-1}(2r/R),\Phi\in(-\Phistar,+\Phistar)\big\}$ (see \cref{fig:x1A0A1BCboundary,fig:Fn1(x1)A0A1BCboundary}), then we have the following.\hfill
    \begin{enumerate}
        \item\label{item01ContractionRegion} The image $\Fc(\mathfrak{L}_1)$ is a curve segment near $(\phistar,\phistar)\in \partial\Mrin$, and the curve segment $\Fc(\mathfrak{L}_1)$ has one endpoint $(\phistar,\phistar-\sin^{-1}(\sqrt{4r/R})-\Phistar)\in\partial\Mrin$ and the tangent slope of every point on the curve segment is in $\big(1,\infty\big)$ in the $\phi\theta$-coordinates (see \cref{fig:x2a0a1bcboundary}).
        \item\label{item02ContractionRegion} The preimage $\Fc^{-1}(\mathfrak{L}_1)$ is a curve segment near $(2\pi-\phistar,\phistar)\in \partial\Mrout$, and the curve segment $\Fc^{-1}(\mathfrak{L}_1)$ has one endpoint $(2\pi-\phistar,\phistar-\sin^{-1}(\sqrt{4r/R})-\Phistar)\in\partial\Mrin$ and the tangent slope of every point on the curve segment is in $\big(-\infty,-1\big)$ in the $\phi\theta$-coordinates (see \cref{fig:x0a0a1bcboundary}).
        \item\label{item03ContractionRegion} The image $\Fc(\mathfrak{L}_2)$ is a curve segment near $(\phistar,\phistar)\in \partial\Mrin$, and the curve segment $\Fc(\mathfrak{L}_1)$ has one endpoint $(\phistar,\phistar-\sin^{-1}(2r/R)-\Phistar)\in\partial\Mrin$ and the tangent slope of every point on the curve segment is in $\big(1,\infty\big)$ in the  $\phi\theta$-coordinates (see \cref{fig:x2a0a1bcboundary}).
        \item\label{item04ContractionRegion} The preimage $\Fc^{-1}(\mathfrak{L}_2)$ is a curve segment near $(2\pi-\phistar,\phistar)\in \partial\Mrout$, and the curve segment $\Fc^{-1}(\mathfrak{L}_1)$ has one endpoint $(2\pi-\phistar,\phistar-\sin^{-1}(2r/R)-\Phistar)\in\partial\Mrin$ and the tangent slope of every point on the curve segment is in $\big(-\infty,-1\big)$ in the $\phi\theta$-coordinates (see \cref{fig:x0a0a1bcboundary}).
        \item\label{item05ContractionRegion}
        In \cref{fig:x0a0a1bcboundary}, the boundary between case (b) and (c) for the $x_2$ region near $(\phistar,\phistar)$ is a curve connecting $(\phistar,\phistar-2\Phistar)$ and $(\phistar,\phistar-\frac{1}{2}\Phistar)$ with tangent slope varying in $(-\infty,-1/2+O(r/R))\cup(1,+\infty]$ with $O(r/R)$ as the Landau symbol.
        \item\label{item06ContractionRegion}
        In \cref{fig:x2a0a1bcboundary}, the boundary between case (b) and (c) for the $x_2$ region near $(2\pi-\phistar,\phistar)$ is a curve connecting $(2\pi-\phistar,\phistar-2\Phistar)$ and $(\phistar,\phistar-\frac{1}{2}\Phistar)$ with tangent slope varying in $(1/2-O(r/R), \infty)\cup [-\infty,-1)$ with $O(r/R)$ as the Landau symbol.
    \end{enumerate}
\end{proposition}

\section{Monotonicity}\label{sec:FPCLA}\label{subsec:Monotonicity}\hfill

The proofs of the results in  \Cref{sec:ECE} (\Cref{RMASUH,TEcase_a0,TEcase_a1,TEcase_b,TEcase_c,contraction_region}) and \cref{def:ApettrajectoryExpansion,RMASUHwithProof,Corollary:EuclideanExpansionMhat} below require delicate estimates across various parameter ranges. In order to make those feasible, we now establish monotonicity of various of the quantities involved, so the inequalities will later only have to be established at end-points of parameter ranges. Readers will probably be inclined to first skip this section and refer back to it as needed later.
For ease of reference, with the length functions $L_1=\tau_1$, $L_0$, $\tau_0$, $d_0$ given in \cref{def:fpclf} here is the list of results proved in this section.
\begin{itemize}
    \item$\theta\mapsto L_1$ increases (\pcref{monotone_p1})
\item$\theta\mapsto L_0$ decreases (\pcref{monotone_p2})
\item$\Phi\mapsto\tau_1$, $\Phi\mapsto L_0$, $\Phi\mapsto L_0/\tau_1$  decrease (\pcref{monotone_p3}, \pcref{monotone_p5})
\item$\Phi\mapsto\tau_0$, $\Phi\mapsto d_0$ increase (\pcref{monotone_p3}, \pcref{monotone_p4})\COMMENT{\red\bfseries Check that the list is correct and complete.\wentao{Yes, with a bit editing.}}
\end{itemize}\begin{figure}[ht]
\begin{center}
    \begin{tikzpicture}[xscale=0.7,yscale=0.7]
        \pgfmathsetmacro{\Dlt}{0.5}
        \pgfmathsetmacro{\PHISTAR}{3};
        \pgfmathsetmacro{\WIDTH}{12};
        \pgfmathsetmacro{\LENGTH}{5};
        \tkzDefPoint(-\PHISTAR,0){LL};
        \tkzDefPoint(-\PHISTAR,\WIDTH){LU};
        \tkzDefPoint(-\PHISTAR,0.5*\WIDTH){LM};
        \tkzDefPoint(\PHISTAR,0){RL};
        \tkzDefPoint(\PHISTAR,\WIDTH){RU};
        \tkzDefPoint(\PHISTAR,0.5*\WIDTH){RM};
        \tkzDefPoint(-\PHISTAR,\PHISTAR){X};
        \tkzDefPoint(-\PHISTAR,\WIDTH-\PHISTAR){IX};
        \tkzDefPoint(\PHISTAR,\WIDTH-\PHISTAR){Y};
        \tkzDefPoint(\PHISTAR,\PHISTAR){IY};
        \tkzDefPoint(-\PHISTAR,\PHISTAR+\Dlt){LLD};
        \tkzDefPoint(-\PHISTAR,\WIDTH-\PHISTAR-\Dlt){LUD};
        \tkzDefPoint(\PHISTAR,\PHISTAR+\Dlt){RLD};
        \tkzDefPoint(\PHISTAR,\WIDTH-\PHISTAR-\Dlt){RUD};
        \tkzDefPoint(0,0.5*\WIDTH){CP};
        \tkzDefPoint(0,\PHISTAR+\Dlt){DLM};
        \tkzDefPoint(0,\WIDTH-\PHISTAR-\Dlt){DUM};
        \draw [thick] (DLM) --(DUM);
        \draw [thin] (LL) --(LU);
        \draw [thin] (LU) --(RU);
        \draw [thin] (RL) --(RU);
        \draw [thin] (LL) --(RL);
        \draw [blue,thin] (LL) --(IX);
        \draw [blue,thin,dashed] (LL) --(IY);
        \draw [blue,thick,dashed] (IX) --(RU);
        \draw [blue,thick] (IY) --(RU);
        \draw [red,thin,dashed]  (X)--(RL);
        \draw [red,thin,dashed]  (Y)--(LU);
        \draw [purple,thick] (LUD)--(RUD);
        \draw [purple,thick] (LLD)--(RLD);
        \draw [purple,thick] (LM)--(RM);
        \coordinate (x1W) at ($(X)!0.7!(RL)$);
        \coordinate (x1E) at ($(IY)!0.7!(RL)$);
        \draw [green,thick] (x1W)--(x1E);
        \tkzLabelPoint[right](x1E){\small $ (\Phi,\theta=\constant_3<{\pi}/{2})\in\MRout$}
        \coordinate (x2W) at ($(X)!0.3!(LL)$);
        \coordinate (x2E) at ($(IY)!0.3!(LL)$);
        \draw [green,thick] (x2W)--(x2E);
        \tkzLabelPoint[left](x2W){\small $ (\Phi,\theta=\constant_4)\in\MRin$}
        \node[] at (0.5*\PHISTAR, 0.5*\PHISTAR+0.5*\Dlt+0.25*\WIDTH){\Large $R_1$};
        \node[] at (0.5*\PHISTAR, \WIDTH-0.5*\PHISTAR-0.5*\Dlt-0.25*\WIDTH){\Large $R_4$};
        \node[] at (-0.5*\PHISTAR,         0.5*\PHISTAR+0.5*\Dlt+0.25*\WIDTH){\Large $R_2$};
        \node[] at (-0.5*\PHISTAR,         \WIDTH-0.5*\PHISTAR-0.5*\Dlt-0.25*\WIDTH){\Large $R_3$};
        \fill [orange, opacity=18/30](DLM) -- (RLD) -- (RM) -- (CP) -- cycle;
        \tkzDefPoint(0.7*\PHISTAR,\PHISTAR+\Dlt){X1L};
        \tkzDefPoint(0.7*\PHISTAR,0.47*\WIDTH){X1U};
        \tkzDefPoint(\PHISTAR,0.47*\WIDTH){X1UR};
        \draw [green,thick] (X1L)--(X1U);
        \draw [red,ultra thin,dashed] (X1U)--(X1UR);
        \tkzLabelPoint[right](X1UR){$\theta=\frac{\pi}{2}-(\sin^{-1}{(\frac{R\sin{\Phi}}{|O_rP|})}-\Phi)=\frac{\pi}{2}-\measuredangle{O_rPO_R}>\frac{\pi}{2}-\phistar$};
        \node[] at (0.5*\PHISTAR, 0.52*\WIDTH){ \small$ (\Phi=\constant_1,\theta)$};
        \tkzLabelPoint[below](LL){$-\Phistar$};
        \tkzLabelPoint[left](LL){$0$};
        \tkzLabelPoint[below](RL){$+\Phistar$};
        \tkzLabelPoint[left](LU){$\pi$};
        \tkzLabelPoint[left](LM){${\pi}/{2}$};
        \tkzLabelPoint[right](RLD){$\theta=\sin^{-1}{({2r}/{R})}$};
        \tkzLabelPoint[right](RUD){$\theta=\pi-\sin^{-1}{({2r}/{R})}$};
        \tkzDefPoint(0,0.35*\WIDTH){HL};
        \tkzDefPoint(\PHISTAR,0.35*\WIDTH){HR};
        \draw[green,thick](HL)--(HR);
        \tkzLabelPoint[right](HR){\small $(\Phi\ge0,\theta=\constant_2<{\pi}/{2})\in R_1\subset\MRout\cap\MRin$}
    \end{tikzpicture}
\end{center}
    \caption{$\mathfrak D^\circ$ is the union of $R_1\dfn\big\{(\Phi,\theta)\bigm|0\le\Phi<\Phistar,\,\sin^{-1}{({2r}/{R})}\le\theta\le{\pi}/{2}\big\}\subset\MRin\cap\Mrout$,\newline
    $R_2\dfn\big\{(\Phi,\theta)\bigm|-\Phistar<\Phi<0,\,\sin^{-1}{({2r}/{R})}\le\theta\le{\pi}/{2}\big\}$,\newline
    $R_3\dfn\big\{(\Phi,\theta)\bigm|-\Phistar<\Phi\le0,\,\pi-\sin^{-1}{({2r}/{R})}\ge\theta\ge{\pi}/{2}\big\}$,\newline
    $R_4\dfn\big\{(\Phi,\theta)\bigm|\Phistar>\Phi>0,\,\pi-\sin^{-1}{({2r}/{R})}\ge\theta\ge{\pi}/{2}\big\}$.}\label{fig:R1}
\end{figure}

For $x_0\in \Mrout$, $\Fc(x_0)=x_1=(\Phi,\theta)\in \MRin$
 (see \cref{fig:MR}). If $\sin^{-1}{(2r/R)}\le\theta\le\frac{\pi}{2}$, then $x_1\in R_1\cup R_2$.
 The following propositions establish monotonicity when we shift a line with a fixed angle crossing the flatter boundary or rotate a line around a fixed point on the flatter boundary. i.e., We move $x_1\in\MRin$ horizontally or vertically in \cref{fig:R1}. Note that by symmetry, the following conclusions on $R_1$ and $R_2$ of the phase space shown in \Cref{fig:R1} will apply to regions $R_3$ and $R_4$. The conclusion of the case on $\theta \in[0,\pi/2]$ will also imply the same conclusion of the case on $\theta \in[\frac{\pi}{2},\pi]$. The subsequent discussions are limited to $\theta\le\pi/2$.
\begin{figure}[ht]
\begin{center}
\begin{tikzpicture}
    \tkzDefPoint(0,0){Or};
    \pgfmathsetmacro{\rradius}{2.3};
    \pgfmathsetmacro{\bdist}{6};
    \tkzDefPoint(0,\bdist){OR};
    \pgfmathsetmacro{\phistardeg}{40};
    \pgfmathsetmacro{\XValueArc}{\rradius*sin(\phistardeg)};
    \pgfmathsetmacro{\YValueArc}{\rradius*cos(\phistardeg)};
    \tkzDefPoint(\XValueArc,-1.0*\YValueArc){B};
    \tkzDefPoint(-1.0*\XValueArc,-1.0*\YValueArc){A};
    \pgfmathsetmacro{\Rradius}{veclen(\XValueArc,\YValueArc+\bdist)};
\clip (-\Rradius,-3.2) rectangle (\Rradius,6.5);
    \draw[name path=Cr,thin] (Or) circle (\rradius);
    \draw[name path=CR,ultra thin,blue!50] (OR) circle (\Rradius);
    \pgfmathsetmacro{\PXValue}{\Rradius*sin(5)};
    \pgfmathsetmacro{\PYValue}{\bdist-\Rradius*cos(5)};
    \tkzDefPoint(\PXValue,\PYValue){P};
    \path[name path=lineQ] (P) -- +(30:7);
    \draw [blue,ultra thin](P) -- +(30:7);
    \path[name path=lineT] (P) -- +(-20:3);
    \draw [blue,ultra thin](P) -- +(-20:3);
    \path[name path=lineT0] (P) -- +(160:7);
    \draw [blue,ultra thin](P) -- +(160:7);
    \path[name intersections={of=Cr and lineQ,by={Q}}];
    \path[name intersections={of=Cr and lineT,by={T}}];
    \path[name intersections={of=Cr and lineT0,by={T0}}];
    \path[name intersections={of=CR and lineT0,by={T1}}];
    
    \tkzDefLine[orthogonal=through P](P,OR)\tkzGetPoint{LF};
    \tkzDefLine[perpendicular=through P](OR,P)\tkzGetPoint{RF};
    \draw[-latex,draw=red,ultra thin,dashed] (OR)--(P);
    \draw[-latex,draw=red,ultra thin,dashed] (P)--(LF);
    \draw[-latex,draw=red,ultra thin,dashed] (P)--(RF);
    
    \tkzLabelPoint[above left,xshift=-3,yshift=3](P){\large $P$};
    \tkzLabelPoint[right,xshift=1](Q){$Q$};
    \tkzLabelPoint[below,xshift=2,yshift=-2](T){\large $T$};
    \tkzLabelPoint[right,above,xshift=8.5](T0){\large $T_0$};
    \tkzLabelPoint[left,xshift=-3.5](Or){$O_r$};
    \tkzLabelPoint[above](OR){$O_R$};
    \tkzDefPointBy[projection=onto P--Q](OR)\tkzGetPoint{M1};
    \tkzDefPointBy[projection=onto P--T0](OR)\tkzGetPoint{M2};
    \tkzDefPointBy[projection=onto P--Q](Or)\tkzGetPoint{M3};
    \tkzDefPointBy[projection=onto P--T0](Or)\tkzGetPoint{M4};
    \tkzDefPointBy[projection=onto OR--M1](Or)\tkzGetPoint{M5};
     \tkzDefPointBy[projection=onto OR--M2](Or)\tkzGetPoint{M6};
    \draw[red,dashed,ultra thin] (OR)--(M1);
    \draw[red,dashed,ultra thin] (OR)--(M2);
    \draw[red,dashed,ultra thin] (Or)--(M3);
    \draw[red,dashed,ultra thin] (Or)--(M4);
    \draw[red,dashed,ultra thin] (Or)--(M5);
    \draw[red,dashed,ultra thin] (Or)--(M6);
    \tkzDrawPoints(OR,Or,P,Q,T,T0,M1,M2,M3,M4,M5,M6);
    \draw[purple](Or)--(OR);
    \tkzLabelPoint[below,right,xshift=4](M1){\large E};
    \tkzLabelPoint[below,left,xshift=-3,yshift=-3](M2){\large F};
    \tkzLabelPoint[above,right,xshift=0.5,yshift=12](M3){\large G};
    \tkzLabelPoint[below](M4){\large H};
    \tkzLabelPoint[above,right,xshift=0.9,yshift=1.8](M5){\large I};
    \tkzLabelPoint[left,xshift=-8.5](M6){\large J};
    \tkzLabelPoint[right,below,xshift=50](LF){Tangent line at P};
    \tkzLabelSegment[left,xshift=-4](OR,Or){$b$}
    \tkzMarkRightAngle(OR,M1,P);
    \tkzMarkRightAngle(OR,M2,T1)
    \tkzMarkRightAngle(Or,M3,Q);
    \tkzMarkRightAngle(Or,M4,T1);
    \tkzMarkRightAngle(Or,M5,M1);
    \tkzMarkRightAngle(Or,M6,M2);
    \tkzMarkAngle[arc=lll,ultra thin,size=1](P,OR,M1);
    \tkzMarkAngle[arc=lll,blue,ultra thin,size=2.5](Or,OR,P);
    \tkzMarkAngle[arc=lll,ultra thin,size=1.5](M6,OR,P);
    
    \node at ($(OR) + (-75:1.3)$)  {\Large $\theta$};
    \node at ($(OR) + (-88:3.2)$)  {
    \Large \textcolor{blue}{$\Phi$}};
    \node at ($(OR) + (-100:2.1)$)  {\Large $\theta$};   
    \coordinate (Y) at (0,-2.8);
    \draw[red,ultra thin,dashed](Or)--(Y);
    \tkzLabelPoint[below](Y){$Y$}
    \draw[very thin,dashed](Or)--(P);
\end{tikzpicture}
    \caption{
    $T_0$, $P$, $Q$ are the positions $p(x_0)$, $p(x_1)$, $p(x_2)$ in billiard table.\\
    $E$ is the foot of $PQ$'s perpendicular through $O_R$. $F$ is the foot of $PT_0$'s perpendicular through $O_R$.\\
    $G$ is the foot of $PQ$'s perpendicular through $O_r$. $H$ is the foot of $PT_0$'s perpendicular through $O_r$.\\
    $I $ is the foot of $O_RE$'s perpendicular through $O_r$.$J$ is the foot of $O_RF$'s perpendicular through $O_r$.}\label{fig:CaseAMonotone}
\end{center}
\end{figure}

\begin{proposition}[$\theta\mapsto L_1$ increases]\label{monotone_p1}
With $\constant_1\in[0,\Phistar)$ fixed, for $x_1$ on the vertical green line in  \cref{fig:R1}, $P=p(x_1)$ is independent of $\theta$ and $L_1=\tau_1$ is a strictly increasing function of $\theta\in[\sin^{-1}{({2r}/{R})},{\pi}/{2}]$ (notations as in \cref{def:fpclf}).
\end{proposition}
\begin{proof}
We show that $\frac{d L_1}{d\theta}>0$. Use cartesian coordinates with $O_r=(0,0)$ and $O_R=(0,b)$. Let $Q=p(\Fc(x_1))$ then $\tau_1=L_1=|PQ|$  (\cref{def:fpclf}). Since $x_1=(\Phi,\theta)\in\MRin$,  $P=p(x_1)$ has cartesian coordinates  
\begin{equation} \label{xP}
P=(x_P,y_P)=(R\sin{\Phi},b-R\cos(\Phi)),
\end{equation}
and likewise $Q$ has coordinates
\begin{equation}\label{xQ}
Q=(x_Q,y_Q)=(x_P+L_1\cos{(\Phi+\theta)},y_P+L_1\sin{(\Phi+\theta)}.
\end{equation}
Thus,
\begin{equation}\label{eqLsquare}
r^2=x^2_Q+y^2_Q=L^2_1+[2x_P\cos{(\Phi+\theta)}+2y_P\sin{(\Phi+\theta)}]L_1+x^2_P+y^2_P.
\end{equation}
Since $\Phi$ is held constant, and hence so are $x_P$, and $y_P$, we have $\frac{dx_P}{d\theta}=0$ and $\frac{dy_P}{d\theta}=0$. Differentiating this last equation with respect to $\theta$ and using the chain rule therefore yields $$2L_1\frac{d L_1}{d\theta}+L_1[-2x_P\sin{(\Phi+\theta)}+2y_P\sin{(\Phi+\theta)}]+\frac{d L_1}{d\theta}[2x_P\cos{(\Phi+\theta)}+2y_P\sin{(\Phi+\theta)}]=0,$$
or
\begin{equation} \label{dL1Dtheta}
\frac{d L_1}{d\theta}=\frac{x_P\sin{(\Phi+\theta)}-y_P\cos{(\Phi+\theta)}}{L_1+x_P\cos{(\Phi+\theta)}+y_P\sin{(\Phi+\theta)}}L_1\stackrel{\eqref{xP}}{=}\frac{\overbracket{R\cos{(\theta)}-b\cos{(\Phi+\theta)}}^{\mathclap{=\text{distance}(O_r,\text{line }PQ)\text{ (proof below)}}}}{\underbracket{L_1+b\sin{(\Phi+\theta)}-R\sin{\theta}}_{=d_2>0\text{ (proof below)}}}L_1>0,
\end{equation}
so  if $x_1=(\Phi=\text{constant},\theta)$, then $\tau_1=L_1$ is an increasing function of $\theta\in[\sin^{-1}{(2r/R)},\pi/2]$.

It remains to compute the last numerator and denominator in \eqref{dL1Dtheta}.

To prove that the denominator satisfies $L_1+b\sin{(\Phi+\theta)}-R\sin{\theta}=d_2$, let $x_2=(\Phi_2,\theta_2)$, and observe from \cref{fig:CaseAMonotone}  that $p(x_2)=Q$, 
$|PQ|=L_1=\tau_1$, $|PE|=R\sin\theta=d_1$, $|GQ|=d_2=r\sin{\theta_2}$, $|GE|=b\sin{(\theta+\Phi)}$.

There are six possible configurations of the points $P,G,Q,E$ on the line $PQ$, depending on the lengths $\tau_1,d_1,d_2$.

Case 1) $\tau_1\le d_2\le d_1$: the ordering is $\overrightarrow{GPQE}$, and $|PQ|+|GE|-|PE|=|GQ|$.

Case 2) $\tau_1\le d_1\le d_2$: the ordering is $\overrightarrow{GPQE}$, and $|PQ|+|GE|-|PE|=|GQ|$. This case does not occur when $d_1\ge2r$.

Case 3) $d_1\le\tau_1\le d_2$: the ordering is $\overrightarrow{GPEQ}$, and $|PQ|+|GE|-|PE|=|GQ|$. This case does not occur when $d_1\ge2r$.

Case 4) $d_2\le\tau_1\le d_1$: the ordering is $\overrightarrow{PGQE}$, and $|PQ|+|GE|-|PE|=|GQ|$.

Case 5) $d_2\le d_1<\tau_1$: the ordering is $\overrightarrow{PGEQ}$ or $\overrightarrow{PEGQ}$, but since $\tau_1\le2d_1,\tau_1\le2d_2,\tau_1\le d_1+d_2$, the ordering   $\overrightarrow{PEGQ}$ cannot happen. So, for the ordering $\overrightarrow{PGEQ}$ observe that $|PQ|+|GE|-|PE|=|GE|+|EQ|=|GQ|$, even though this case also does not occur when $d_1\ge2r$.

Case 6) $d_1\le d_2<\tau_1$: reasoning as in case 5) we see that the ordering is $\overrightarrow{PGEQ}$, hence $|PQ|+|GE|-|PE|=|GE|+|EQ|=|GQ|$, even though this case does not occur when $d_1\ge2r$.

Hence in all six cases, $L_1+b\sin{(\theta+\Phi)}-R\sin{\theta}=|PQ|+|GE|-|PE|=|GQ|=d_2=r\sin{\theta_2}$.



Lastly, we prove that the numerator $R\cos{(\theta)}-b\cos{(\Phi+\theta)}=\text{distance}(O_r,\text{line }PQ)$. In \cref{fig:CaseAMonotone} suppose that a line $l_1$ passing through $O_r$ is parallel to the line $PQ$. For $\Phi\in[0,\Phistar)$ fixed, we observe that: if $\pi/2-\theta<\Phi$, then $|O_RI|=b\cos{(\theta+\Phi)}$, $|O_RE|=R\cos\theta$, distance$(O_r,\text{line }PQ)=|IE|=|O_RE|-|O_RI|=R\cos\theta-b\cos{(\theta+\Phi)}$. And if $\pi/2-\theta\ge\Phi$, then $|O_RI|=-b\cos{(\theta+\Phi)}$, $|O_RE|=R\cos\theta$, distance$(O_r,\text{line }PQ)=|IE|=|O_RE|+|O_RI|=R\cos\theta-b\cos{(\theta+\Phi)}$. Therefore, \(R\cos\theta-b\cos{(\theta+\Phi)}=\text{distance}(O_r,\text{line }PQ)>0\).
\end{proof}
\begin{proposition}[$\theta\mapsto L_0$ decreases]\label{monotone_p2}
With $\constant_1\in[0,\Phistar)$ fixed, for $x_1=(\Phi_1=\constant_1,\theta)$ on the vertical green line in the region $R_1=\big\{(\Phi,\theta)\in \MRin \bigm| 0\le\Phi<\Phistar,\;\sin^{-1}{(2r/R)}\le\theta\le\pi/2\big\}$ in \Cref{fig:R1} (notation as in \Cref{def:fpclf}), $P=p(x_1)$ is held a constant, thus independent of $\theta$, and 
$L_0=2d_0-\tau_0$ is a strictly decreasing function of $\theta\in\big[\sin^{-1}{(2r/R)},\ \pi/2-(\sin^{-1}{(R\sin{\Phi_1}/\rho)}-\Phi_1))\big]$, where $\rho\dfn|O_rP|$ (in \cref{fig:CaseAMonotone}; this is independent of $\theta$ since $P=p(x_1)$ is).
\end{proposition}
\begin{remark}\label{rmk:monotone_p2}\hfill
\begin{enumerate}
    \item\label{item02rmkmonotone_p2} For $\Phi_1\in[0,\Phistar)$, in \cref{fig:CaseAMonotone} we see $\Phi=\Phi_1=\measuredangle{O_rO_RP}$, $\measuredangle{PO_rY}=\sin^{-1}{(\frac{R\sin\Phi_1}{\rho})}$. Therefore, $\measuredangle O_rPO_R=\measuredangle{PO_rY}-\measuredangle{O_rO_RP}=\sin^{-1}{(\frac{R\sin\Phi_1}{\rho})}-\Phi_1$. Since $0\le\measuredangle{O_rPO_R}\le\measuredangle{YO_rP}<\phistar$, $0\le\sin^{-1}{(\frac{R\sin\Phi_1}{\rho})}-\Phi_1<\phistar$,
    \[
        \pi/2-(\sin^{-1}{(R\sin{\Phi_1}/\rho)}-\Phi_1))>\pi/2-\phistar.
    \] This fact is referenced in the proof of \cref{prop:casea0smallthetaR1}.\COMMENT{If it is this specific fact, then maybe the equation itself should be labeled and referred to, instead of the full remark.}
    \item\label{item03rmkmonotone_p2} In \cref{fig:CaseAMonotone} suppose that a line $l$ passing through $O_r$ is parallel to the line $PT_0$. Therefore, distance($l,O_R$)$=|O_RJ|=b\cos(\theta-\Phi)$, distance($\text{line } PT_0,O_R$)$=|O_RF|=R\cos\theta$. For $\Phi\in[0,\Phistar)$ fixed, we observe that in \cref{fig:CaseAMonotone}, the order of $O_R,F,J$ depends on the comparing of $\measuredangle{T_0PO_R}$ with $\measuredangle{O_rO_RP}$. If $\measuredangle{T_0PO_R}=\measuredangle{O_rPO_R}$, then $PT_0$ passes through $O_r$, so $J$ and $F$ are the same point, $|O_RJ|=|O_RF|$.
    
    if $\pi/2-\theta=\measuredangle{T_0PO_R}>\measuredangle{O_rPO_R}=\measuredangle{YO_rP}-\measuredangle{O_rO_RP}=\sin^{-1}{(\frac{R\sin\Phi}{\rho})}-\Phi$, then $|O_RJ|<|O_RF|$.
    
    And if $\pi/2-\theta=\measuredangle{T_0PO_R}\le\measuredangle{O_rPO_R}=\measuredangle{YO_rP}-\measuredangle{O_rO_RP}=\sin^{-1}{(\frac{R\sin\Phi}{\rho})}-\Phi$, then $|O_RJ|\ge |O_RF|$.

\end{enumerate}
\end{remark}
\begin{proof}
With fixed $\Phi=\constant_1$, in cartesian coordinates (\cref{def:StandardCoordinateTable}) for  \cref{fig:CaseAMonotone}, 
\begin{equation} \label{xPP}
        P=(x_P,y_P)=(R\sin{\Phi},b-R\cos(\Phi)),
\end{equation}
is independent of $\theta$. 

By \Cref{def:fpclf}, we have $L_0=2d_0-
\tau_0=|TP|$ in \cref{fig:CaseAMonotone}.
And by elementary geometry, we can observe from \cref{fig:CaseAMonotone} that if $x_1=(\Phi,\theta)\in \MRin$, then $\theta>\Phi$ and\COMMENT{\textcolor{red}{Add {\bfseries correct} punctuation in ALL displays!}\wentao{Ok.}}
\begin{equation}\label{xT}
        r^2=x^2_T+y^2_T\overbracket=^{\mathrlap{T=(x_T,y_T)=(x_P+L_0\cos{(\theta-\Phi)},y_P-L_0\sin{(\theta-\Phi)}}\quad}L^2_0+[2x_P\cos{(\theta-\Phi)}-2y_P\sin{(\theta-\Phi)}]L_0+x^2_P+y^2_P.
\end{equation}
Since $\Phi$, $x_P$, $y_P$ are constants, $\frac{dx_P}{d\theta}=0$, $\frac{dy_P}{d\theta}=0$. Differentiating \eqref{xT} with respect to $\theta$ and using the chain rule gives
\[
2L_0\frac{dL_0}{d\theta}-[2x_{P}\sin{(\theta-\Phi)}+2y_{P}\cos(\theta-\Phi)]L_0+[2x_{P}\cos{(\theta-\Phi)-2y_{P}\sin{(\theta-\Phi)}}]\frac{dL_0}{d\theta}=0,
\]
so
\begin{equation} \label{dL2Dtheta}
\begin{aligned}
\frac{d L_0}{d\theta}  &=\frac{x_{P}\sin{(\theta-\Phi)}+y_{P}\cos{(\theta-\Phi)}}{L_0+x_{P}\cos{(\theta-\Phi)}-y_{P}\sin{(\theta-\Phi)}}L_0\\
&\overbracket{=}^{\mathllap{\eqref{xPP}}}\frac{R\sin{(\Phi)}\sin{(\theta-\Phi)}-R\cos{\Phi}\cos{(\theta-\Phi)}+b\cos{(\theta-\Phi)}}{L_0+R\sin{\Phi}\cos{(\theta-\Phi)}+R\cos{\Phi}\sin{(\theta-\Phi)}-b\sin{(\theta-\Phi)}}L_0=\frac{\overbracket{b\cos(\theta-\Phi)-R\cos{\theta}}^{\mathclap{<0\text{ if } 0<\theta<\frac{\pi}{2}-\sin^{-1}{(\frac{R\sin\Phi}{\rho})}\text{\cref{rmk:monotone_p2}\eqref{item03rmkmonotone_p2}}}}}{\underbracket{L_0-b\sin{(\theta-\Phi)+R\sin{\theta}}}_{=d_0>0\text{ (proof below)}}}L_0<0.
 \end{aligned}
\end{equation}

Now, it remains to compute the numerator and denominator of the term on the right-hand side of the last equality in \eqref{dL2Dtheta}.

For the numerator, \cref{rmk:monotone_p2}\eqref{item03rmkmonotone_p2}: if $0<\theta<\frac{\pi}{2}-\sin^{-1}{(\frac{R\sin{\Phi}}{\rho})}+\Phi$, then $b\cos{(\Phi-\theta)}-R\cos\theta=|O_RJ|-|O_RF|<0$.

For the denominator, 
\cref{fig:CaseAMonotone} gives $L_0=|TP|$,  $b\sin{(\theta-\Phi)}=|HF|$, $|PF|=R\sin\theta=d_1$, and $|TH|=d_0$.
The ordering of $T,P,H,F$  on the line $TP$ is determined by the lengths $L_0=2\tau_0-d_0,d_1,d_0$. Since $d_1\ge2r$, there are two cases according to the sign of $d_0-\tau_0$:
\newline$d_1\ge2r> d_0\ge2d_0-\tau_0$: the ordering is $\overrightarrow{TPHF}$, and $L_0-b\sin{(\theta-\Phi)}+R\sin{\theta}=|TP|-|HF|+|PF|=|TH|=d_0$.
\newline$d_1\ge2r>2d_0-\tau_0>d_0$: the ordering is $\overrightarrow{THPF}$, and $|TP|-|HF|+|PF|=|TH|=d_0$ (this case will not occur when $x_1\in R_1$).

Hence, $\frac{dL_0}{d\theta}=\frac{|O_RJ|-|O_RF|}{d_0}L_0<0$, so $L_0=2d_0-\tau_0$ is a decreasing function of $\theta$.
\end{proof}

\begin{proposition}[$\Phi\mapsto\tau_1$, $\Phi\mapsto L_0$, $\Phi\mapsto-\tau_0$ decrease]\label{monotone_p3}
Referring to \cref{fig:R1,def:StandardCoordinateTable}:
\begin{enumerate}
     \item\label{itemmonotone_p3-1}Suppose $x_1=(\Phi,\theta)$ is on a horizontal line segment in $\MRout$ with $\theta<\pi/2$ constant, and $P=p(x_1)=(x_P,y_P)$, $Q=p(\Fc(x_1))=(x_Q,y_Q)$. (Two such horizontal line segments are shown in green: $\big\{(\Phi,\theta=\constant_4<\frac{\pi}{2})\in \MRout,\Phi\in [-\Phistar,\Phistar]\big\}$ and $\big\{(\Phi,\theta=\constant_3<\Phistar)\in \MRout\big\}$.) Then $\tau_1$ (\cref{def:fpclf}) is a function of $\Phi$ with
\begin{equation}\label{eq1monotone_p3}
    \frac{d\tau_1}{d\Phi}=-\frac{bx_Q}{d_2}.
\end{equation}
If $\Fc(x_1)=(\phi_2,\theta_2)\in\Mrin$, then $d_2=r\sin{\theta_2}$, $x_Q=r\sin{\phi_2}$ (\cref{def:fpclf}), and if $x_Q>0$, then $\Phi\mapsto\tau_1$ decreases.

\item\label{itemmonotone_p3-2} Suppose $x_1$ is on a horizontal line segment $\big\{(\Phi,\theta=\constant_2<\frac{\pi}{2})\subseteq R_1,\ \Phi\in [0,\Phistar]\big\}\subset\MRin$ with $\sin^{-1}{(\frac{2r}{R})}\le \constant_2\le\frac{\pi}{2}$. The intersection $T$ of $\Gamma_R$ with the line through $T_0=p(\Fc^{-1}(x_1))$ and $P$ has cartesian coordinates $(x_T,y_T)$. Then $L_0=2d_0-\tau_0$ is a function of $\Phi$ with 
\begin{equation}\label{eq2monotone_p3}
    \frac{dL_0}{d\Phi}=-\frac{bx_T}{d_0}.
\end{equation} If $x_T>0$, then $\Phi\mapsto L_0$ decreases.
\item\label{itemmonotone_p3-3}Suppose $x_1=(\Phi,\theta)$ is on a horizontal line segment in $\MRin$ with $\theta<\pi/2$ constant. Let $P=p(x_1)=(x_P,y_P)$, $T_0=p(\Fc^{-1}(x_1))=(x_{T_0},y_{T_0})$. (Two examples of such horizontal line segments are shown in green: $\big\{(\Phi,\theta=\constant_4<\frac{\pi}{2})\in \MRin,\Phi\in [-\Phistar,\Phistar]\big\}$ and $\big\{(\Phi,\theta=\constant_4<\Phistar)\in \MRin\big\}$.) Then $\tau_0$ (\cref{def:fpclf}) is a function of $\Phi$ with
\begin{equation}\label{eq3monotone_p3}
    \frac{d\tau_0}{d\Phi}=-\frac{bx_{T_0}}{d_0}.
\end{equation}
If $\Fc^{-1}(x_1)=(\phi_0,\theta_0)\in\Mrout$, then $d_0=r\sin{\theta_0}$, $x_{T_0}=r\sin{\phi_0}$ (\cref{def:fpclf}), so if $x_{T_0}<0$, then $\Phi\mapsto\tau_0$ increases.
\end{enumerate}
\end{proposition}
\begin{proof}
\textbf{Proof of \eqref{itemmonotone_p3-1}}: In \cref{fig:CaseAMonotone}, $L_1=\tau_1=|PQ|,L_0=2d_0-\tau_0=|PT|$. The coordinates of $P,Q,T$ are the same as in \cref{monotone_p1,monotone_p2}. The equations in\eqref{xP}\eqref{eqLsquare} are still true.

Differentiating \eqref{xP} with respect to $\Phi$ gives

\begin{equation}
\begin{aligned}
    & \frac{dx_P}{d\Phi}=R\cos{\Phi},\\
    & \frac{dy_P}{d\Phi}=R\sin{\Phi}.
\end{aligned}
\end{equation}

Differentiating \eqref{eqLsquare} with respect to $\Phi$ and using the chain rule yields
\begin{equation} \label{eq:dphiL1}
\begin{multlined}[.85\textwidth]
 2L_1\frac{dL_1}{d\Phi}+2\frac{dL_1}{d\Phi}\big[x_P\cos{(\theta+\Phi)}+y_P\sin{(\theta+\Phi)}\big]\\+2L_1\big[R\cos{\theta}-x_P\sin{(\theta+\Phi)+y_P\cos{(\theta+\Phi)}}\big]
 +2x_PR\cos{\Phi}+2y_PR\sin{\Phi}=0,
\end{multlined}
\end{equation}
that is,
\begin{equation}\label{eq:dphiL1_1}
\begin{aligned}
\big(L_1+b\sin{(\theta+\Phi)}-R\sin{\theta}\big)\frac{dL_1}{d\Phi} &=
{-L_1(R\cos{(\theta)}-\overbracket{x_P}^{\mathclap{=R\sin{\Phi}}}\sin{(\theta+\Phi)}+\overbracket{y_P}^{\mathclap{=b-R\cos{\Phi}}}\cos{(\theta+\Phi)})-bR\sin{\Phi}}
\\
& =
{-L(R\cos{\theta}-R\cos{\theta}+b\cos{(\theta+\Phi)})-bR\sin{\Phi}}
\\
& =
{-bL_1\cos{(\theta+\Phi)}-bR\sin{\Phi}}
.\end{aligned}
\end{equation}
Note that if $x_1\in\MRout$, then $\theta+\Phi>0$. In \cref{fig:CaseAMonotone}, $|PQ|=L_1=\tau_1$, $|PE|=d_1$, $|GQ|=d_2$, and $|GE|=b\sin{(\theta+\Phi)}$.

Using the same proof as in \cref{monotone_p1}, note that in all six cases of that proposition, we have$$L_1+b\sin{(\theta+\Phi)}-R\sin{\theta}=d_2.$$
Then by \eqref{xP} and from \cref{fig:CaseAMonotone}, $$L_1\cos{(\theta+\Phi)}+R\sin{\Phi}=(x_Q-x_P)+x_P=x_Q.$$
Therefore, we have proved item \eqref{itemmonotone_p3-1}, i.e.,
$$\frac{d\tau_1}{d\Phi}=\frac{dL_1}{d\Phi}=-b\frac{x_Q}{d_2}.$$
\newline \textbf{Proof of \eqref{itemmonotone_p3-2}}:
Differentiating \eqref{xT} with respect to $\Phi$ yields $\frac{dL_0}{d\Phi}=-b\frac{R\sin\Phi+L_0\cos{(\theta-\Phi)}}{L_0+R\sin\theta-b\sin{(\theta-\Phi)}}$.

From \cref{fig:CaseAMonotone} and equation \eqref{xP}, we get $R\sin\Phi+L_0\cos{(\theta-\Phi)}=x_P+(x_T-x_P)=x_T$, and from \eqref{dL2Dtheta} $d_0=L_0+R\sin\theta-b\sin{(\theta-\Phi)}$,
therefore $\frac{d(2d_0-\tau_0)}{d\Phi}=\frac{dL_0}{d\Phi}=-b\frac{x_T}{d_0}.$ This proves \eqref{itemmonotone_p3-2}.\COMMENT{Immediately after proof of (1) explain how it implies (3). That is currently missing. Refer to remark on "reversal symmetry" as needed\red\small\wentao{So far (3) is used once later.}}
\newline \textbf{Proof of \eqref{itemmonotone_p3-3}}: The mirror symmetry $J:(\Phi,\theta=\Psi_4)\mapsto (-\Phi,\theta=\Psi_4)$ in \eqref{eq:FoldSymmetryDef} diffeomorphically maps the line segment in $\mathfrak{l}\subset\MRin$ with $\theta$ being a constant $\Psi_4$ to a line segment $J(\mathfrak{l})\subset\MRout$ with $\theta$ being the same constant $\Psi_4$. Since concerning $\tau_0$ as a function of $x_1\in\MRin$ and $\tau_1$ as a function on $\MRout$, $\tau_0$ and $\tau_1$ satisfy $\tau_0\restriction_{\MRin}=\tau_1\restriction_{\MRout}\circ J$. 
When $\tau_0$ is restricted to the line segment $\mathfrak{l}$ with $\Phi$ as parameter, $\frac{dJ}{d\Phi}=-1$ and $\frac{d\tau_0}{d\Phi}$ satisfies the following.
\begin{align*}
    \frac{d\tau_0}{d\Phi}\Bigm|_{x_1=(\Phi,\theta=\Psi_4)}=\frac{d(\tau_1\circ J)}{d\Phi}\overbracket{=}^{\text{chain rule}}\frac{d\tau_1}{d\Phi}\Bigm|_{J(x_1)}\overbracket{\frac{d J}{d\Phi}}^{=-1}=-\frac{d\tau_1}{d\Phi}\Bigm|_{J(x_1)}.
\end{align*}
Let $\Fc J(x_1)\nfd x_{2,J}=(\phi_{2,J},\theta_{2,J})\in\Mrin$, then, when restricted to the line segment $J(\mathfrak{l})$, \eqref{itemmonotone_p3-1} gives $$\frac{d\tau_1}{d\Phi}\Bigm|_{J(x_1)}=-b\frac{r\sin{\phi_{2,J}}}{r\sin{\theta_{2,J}}}.$$ Hence $\frac{d\tau_0}{d\Phi}\Bigm|_{x_1=(\Phi,\theta=\Psi_4)}=b\frac{r\sin{\phi_{2,J}}}{r\sin{\theta_{2,J}}}$. Since $\Fc^{-1}(x_1)=x_0=(\phi_0,\theta_0)$ and $J\Fc^{-1}=\Fc J$, it gives $$x_{2,J}=\Fc J(x_1)=J\Fc^{-1}(x_1)=J(x_0),$$, by \eqref{eq:FoldSymmetryDef} which means that $\phi_0=2\pi-\phi_{2,J}$ and $\theta_0=\theta_{2,J}$. Thus, \[\frac{d\tau_0}{d\Phi}\Bigm|_{x_1=(\Phi,\theta=\Psi_4)}=b\frac{r\sin{\phi_{2,J}}}{r\sin{\theta_{2,J}}}=-b\frac{r\sin{\phi_0}}{r\sin{\theta_0}}=-b\frac{x_{T_0}}{d_0}.\qedhere\]
\end{proof}

\begin{proposition}[$\Phi\mapsto d_0$ increases]\label{monotone_p4}
Given $0<\phistar<\tan^{-1}(1/3)$ and $R$ which meets the hyperbolicity condition \eqref{eqJZHypCond}, $x\in\hat{M}$ has an orbit segment defined in \cref{def:fpclf} with $x_1\in \MRin$. If $x_1$ moves along a horizontal line segment: $\{\,x_1=(\Phi,\,\theta=\constant_4\,)\,\}\subset\MRin$ in \cref{fig:R1} for some $\constant_4<\sin^{-1}{(\sqrt{4r/R})}$,\COMMENT{It would be good to explain the origin of this upper bound when it first appears.\wentao{It's used in \eqref{eq:new74}.}} then $d_0$ is an increasing function of $\Phi$.
\end{proposition}
\begin{proof}
In the \cref{fig:R1_R2_d_0} and in the \cref{def:StandardCoordinateTable} coordinate system, 
$P=p(x_1)=(x_P,y_P)=(R\sin(\Phi),b-R\cos{\Phi})$ for $x_1=(\Phi,\theta)\in\MRin$. And $x_0=\Fc^{-1}(x_1)=(\phi_0,\theta_0)\in\Mrout, Q=p(x_0)=(x_Q,y_Q)=(r\sin\phi_0,-r\cos\phi_0)$. $T$ is the another intercept of line $PQ$ with circle $C_r$. Thus, the collision angle \(\theta\) is the angle between $\overrightarrow{TP}$ and the tangent direction of $\Gamma_R$ at $P$ given the counterclockwise orientation on $\Gamma_R$. 


By \cref{proposition:NoutsizeInA1BCcases}\eqref{item01NoutsizeInA1BCcases} we get the following.
\[
\big|\phi_0-(2\pi-\phistar)\big |\overbracket{<}^{\mathclap{\text{\cref{proposition:NoutsizeInA1BCcases}\eqref{item01NoutsizeInA1BCcases}}}}2.001\sqrt{4r/R}\overbracket{<}^{\mathclap{\text{\eqref{eqJZHypCond0}: }R>1700r}}0.097.
\]


Hence, $\phi_0-2\pi+\phistar>-0.097$, so $ \phi_0>2\pi-\phistar-0.097>2\pi-0.0971-\tan^{-1}(1/3)>5.86>3\pi/2$. This means $Q$ is on the lower left half of the circle $C_r$ in \cref{fig:R1_R2_d_0}. Hence $x_Q<0$, $y_Q<0$.


The range of $\Phi$ is determined by the two constraints $P\in\Gamma_R$ and $T\notin\Gamma_R$ in \Cref{fig:R1_R2_d_0}.\COMMENT{\red\small Better to just omit the previous broken sentence?\wentao{OK}}
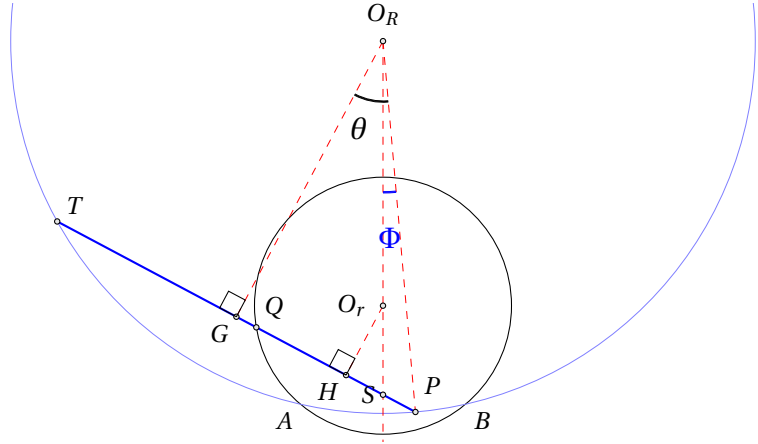
\begin{figure}[ht]
\begin{center}
\begin{tikzpicture}
    \tkzDefPoint(0,0){Or};
    \pgfmathsetmacro{\rradius}{1.7};
    \pgfmathsetmacro{\bdist}{3.5};
    \tkzDefPoint(0,\bdist){OR};
    \pgfmathsetmacro{\phistardeg}{40};
    \pgfmathsetmacro{\XValueArc}{\rradius*sin(\phistardeg)};
    \pgfmathsetmacro{\YValueArc}{\rradius*cos(\phistardeg)};
    \tkzDefPoint(\XValueArc,-1.0*\YValueArc){B};
    \tkzDefPoint(-1.0*\XValueArc,-1.0*\YValueArc){A};
    \pgfmathsetmacro{\Rradius}{veclen(\XValueArc,\YValueArc+\bdist)};
\clip(-\Rradius,-1.8) rectangle (\Rradius,4);
    \draw[name path=Cr,thin] (Or) circle (\rradius);
    \draw[name path=CR,ultra thin,blue!50] (OR) circle (\Rradius);
    \pgfmathsetmacro{\PXValue}{\Rradius*sin(5)};
    \pgfmathsetmacro{\PYValue}{\bdist-\Rradius*cos(5)};
    \pgfmathsetmacro{\LPXValue}{-\Rradius*sin(70)};
    \pgfmathsetmacro{\LPYValue}{\bdist-\Rradius*cos(70)};
    \tkzDefPoint(\PXValue,\PYValue){P};
    \tkzDefPoint(0,-2){Y};
    \draw[red,ultra thin,dashed](OR)--(Y);
    \path[name path=lineQ] (P) -- +(152:10);
    \path[name path=lineY] (Or) -- (Y);
    \path[name intersections={of=CR and lineQ,by={P0}}];
    \path[name intersections={of=Cr and lineQ,by={Q}}];
    \path[name intersections={of=lineQ and lineY,by={S}}];
    \draw[blue,thick](P)--(P0);
    \tkzDefPointBy[projection=onto P--Q](OR)\tkzGetPoint{G};
    \tkzDefPointBy[projection=onto P--Q](Or)\tkzGetPoint{H};
    \draw[red,ultra thin,dashed](Or)--(H);
    \draw[red,ultra thin,dashed](OR)--(G);
    \draw[red,ultra thin,dashed](OR)--(P);
    
    \tkzMarkAngle[arc=lll,ultra thin,size=0.8](G,OR,P);
    \tkzMarkAngle[blue,arc=lll,ultra thin,size=2](Y,OR,P);
    \node at ($(OR) + (-105:1.2)$)  {\Large $\theta$};
    \node at ($(OR) + (-88:2.6)$)  {
    \Large \textcolor{blue}{$\Phi$}};
    \tkzLabelPoint[above,yshift=3](OR){$O_R$};
    \tkzLabelPoint[left,xshift=-3](Or){$O_r$};
    \tkzLabelPoint[above,right,yshift=10](P){$P$};
    \tkzLabelPoint[left](S){$S$};
    \tkzLabelPoint[below left](G){$G$};
    \tkzLabelPoint[below left](H){$H$};
    \tkzLabelPoint[above right](P0){$T$};
    \tkzLabelPoint[above right](Q){$Q$};
    \tkzLabelPoint[below left](A){$A$};
    \tkzLabelPoint[below right](B){$B$};
    \tkzDrawPoints(OR,Or,P,Q,S,G,H,P0);
    \tkzMarkRightAngle(Or,H,P0);
    \tkzMarkRightAngle(OR,G,P0);
\end{tikzpicture}
\end{center}
    \caption{$PT$ moves along $\Gamma_R$ with constant collision angle $\theta$ small. $G$ is the perpendicular foot of $O_R$ on $PT$. $H$ is the perpendicular foot of $O_r$ on $PT$.}\label{fig:R1_R2_d_0}
\end{figure}

Let $S=(0,y)$ be the intersection of line $PT$ and the $y$-axis (in \cref{fig:R1_R2_d_0}). From \eqref{xP} and by checking the slopes of $\overrightarrow{SP}$, we see that $\frac{y-b+R\cos{\Phi}}{-R\sin{\Phi}}=-\tan{(\theta-\Phi)}$. This gives $y=b-R\frac{\cos{\theta}}{\cos{(\theta-\Phi)}}$.  
Now,
$x_1=(\Phi,\,\theta=\constant_4\,)\in\MRin$ in \cref{fig:R1}, so
\[
0<\theta-\Phi<\theta+\Phistar\overbracket{=}^{\strut\mathclap{\text{\cref{def:BasicNotations0}}}}\theta+\sin^{-1}(r\sin\phistar/R)\overbracket{<}^{\mathclap{\theta=\constant_4<\sin^{-1}\sqrt{4r/R}}}\sin^{-1}{\sqrt{4r/R}}+\sin^{-1}(r\sin\phistar/R)\overbracket{<}^{\mathclap{\strut\text{\eqref{eqJZHypCond}: }R>1700r}}\pi/2.
\]
Also, $R$, $b$, $\theta$ are constants, so $y$ increases with $\Phi$, that is, if $\Phi$ increases, then $y$ increases. Then we will show that $y$ is always negative.

The slopes of lines $PQ$ and $PS$ are the same: $\frac{y-y_Q}{0-x_Q}=\frac{y_P-y_Q}{x_P-x_Q}$, which gives $y=y_Q-x_Q\frac{y_P-y_Q}{x_P-x_Q}$. The positions of $P,Q$ in the table imply $x_P-x_Q>0$, $y_P-y_Q<0$ (so the fraction is negative). Since $Q$ is in the lower left quarter circle,  $x_Q<0>y_Q$, hence 
$y=y_Q-x_Q\frac{y_P-y_Q}{x_P-x_Q}<y_Q<0$, that is, $y$ is always negative as $x_1$ moves along the horizontal line: $\{x_1=(\Phi,\,\theta=\constant_4\,)\,\}\subset\MRin$ in \cref{fig:R1} with some $\constant_4<\sin^{-1}{(\sqrt{4r/R})}$. Hence, if $\Phi$ increases, then $|y|$ decreases.

Note that we also have $\frac{|O_rH|}{|O_RG|}=\frac{|O_rS|}{|O_RS|}=\frac{|y|}{|y|+b}=1-\frac{b}{b+|y|}\overbracket{>}^{\mathclap{|y|>0}}0$, and $|O_RG|=R\cos\theta$ remains constant since $\theta$ is held as constant. If $|y|$ decreases, then $|O_rH|$ decreases.

As $P$ moves on $\Gamma_R$ with $P=(R\sin\Phi,b-R\cos\Phi)\in\Gamma_R$ and $T\in\Gamma_r$ such that $\overrightarrow{TP}$ maintains a constant $\theta=\constant_4$ as the angle with tangent direction of $\Gamma_R$ at $P$, if $\Phi$ increases, then $|y|$ decreases and $0<|O_rH|$ decreases, so $d_0=\sqrt{r^2-|O_rH|^2}$ increases. Hence $d_0$ is an increasing function of $\Phi$.
\end{proof}
\begin{proposition}[$\Phi\mapsto L_0/\tau_1$ decreases]\label{monotone_p5}
Given $\phistar\in(0,\tan^{-1}{(1/3)})$ and $R$ satisfying \eqref{eqJZHypCond}, if $x_1=(\Phi, \theta)$ moves on the horizontal line segment $\big\{(\Phi, \theta=\sin^{-1}{\sqrt{4r/R}}) \bigm| \Phi\in[0,\Phistar] \big\}$ in the region $R_1$ of \cref{fig:R1}, then the length ratio $\frac{2d_0-\tau_0}{\tau_1}$ is a decreasing function of $\Phi$.
\end{proposition}
The proof is strenuous; it ends on page \pageref{eq:34}.\COMMENT{Shortcuts would be nice. Probably no easy ones, though.\wentao{Maybe AI could help..}}
\begin{proof}
Consider \cref{fig:lengthfraction} with cartesian coordinates such that $O_r=(0,0)$, $O_R=(0,b)$. Let $P=p(x_1)$ with $x_1=(\Phi,\theta=\sin^{-1}{\sqrt{4r/R}})\in R_1\subset\MRin$. As shown in \Cref{fig:lengthfraction}, $T_0$ is the position $p(\Fc^{-1}(x_1))$, and  $Q$ is the position $p(\Fc(x_1))$. The extension of $T_0P$ intersects with $\Gamma_r$ at $T$. $Q_1$ is the point of the horizontal line through $Q$ intersecting $\Gamma_r$. $T_1$ is the point of the horizontal line through $T$ intersecting $\Gamma_r$.
Suppose 
\begin{equation}\label{eq:15}
\begin{split}
    & Q=(x_Q,y_Q)=(r\sin{\alpha},-\cos{\alpha}),\\
    & T=(x_T,y_T)=(r\sin{\beta},-\cos{\beta}).
\end{split}
\end{equation}

It is clear that $\pi/2\ge\alpha>\phistar>\beta>0$ from \cref{fig:lengthfraction}. Because $0\le\Phi\le\Phistar=\sin^{-1}{(r\sin\phistar/R)}<\sin^{-1}{\sqrt{4r/R}}=\theta$, the slope of the line $T_0P$ is $\tan{(\Phi-\theta)}<0$, the slope of the line $Q_0P$ is $\tan{(\Phi+\theta)}>0$ which is larger than $\tan{\Big(\Phi+(\Phistar-\Phi)/2\Big)}>0$, and both slopes are bounded. If $\alpha>\pi$, then $\theta\ge\pi/2$. If $3\pi/2>\alpha>\pi/2$, then by \cref{lemma:returnorbittheta1Large} $\theta\in(\frac{\pi}{6},\frac{5\pi}{6})$. With \eqref{eqJZHypCond}: $R>1700r$, it is impossible for $\theta=\sin^{-1}{\sqrt{4r/R}}$ to satisfy $\theta\ge \pi/2$ or $\theta\in(\frac{\pi}{6},\frac{5\pi}{6})$.

Next, as shown in \Cref{fig:lengthfraction}, we pick $P_1$ on segment $QQ_1$ and $P_2$ on segment $TT_1$ such that $\measuredangle{PP_1Q_1}=\measuredangle{YO_rQ}=\alpha,\measuredangle{PP_2T}=\measuredangle{YO_rT}=\beta$. By elementary geometry,  $\measuredangle{T_1T_0T}=\measuredangle{PP_2T}=\beta$ and $\measuredangle{Q_1Q_0Q}=\measuredangle{PP_1Q}=\pi-\alpha$, so $\Delta TP_2P$ is similar to $\Delta TT_0T_1$ and $\Delta QP_1P$ is similar to $\Delta QQ_0Q_1$. Therefore, \begin{equation}\label{eq:16}
\begin{split}
& \frac{1}{|P_1Q|}=\frac{|Q_1Q|}{|Q_0Q|\times|PQ|}=\frac{x_Q}{d_2|PQ|},\\
& \frac{1}{|P_2T|}=\frac{|T_1T|}{|T_0T|\times|PT|}=\frac{x_T}{d_0|PT|}.
\end{split}
\end{equation}
Since $|PQ|=\tau_1=L_1$ and $|PT|=2d_0-\tau_0=L_0$, the quotient rule gives
\begin{equation}
\label{eq:18}
\frac{d}{d\Phi}\big(\frac{2d_0-\tau_0}{\tau_1}\big)=\frac{1}{L^2_1}\big(\overbracket{\underbracket{\frac{dL_0}{d\Phi}}_{\qquad\mathllap{\overset{\eqref{eq3monotone_p3}}=-b{x_T}/{d_0}}}L_1-\underbracket{\frac{d
L_1}{d\Phi}}_{\mathrlap{\overset{\eqref{eq1monotone_p3}}=-b{x_Q}/{d_2}}\qquad} L_0}^{\mathclap{=b\big(-\frac{x_T}{d_0}L_1+\frac{x_Q}{d_2}L_0\big)
=bL_1L_0\big(-\frac{x_T}{d_0L_0}+\frac{x_Q}{d_2L_1}\big)\qquad\qquad}}\big)
=\frac{bL_0}{L_1}\big(\overbracket{-\underbracket{\frac{x_T}{d_0|PT|}}_{\mathclap{\overset{\eqref{eq:16}}=1/|P_2T|\quad}}+\underbracket{\frac{x_Q}{d_2|PQ|}}_{\mathclap{\quad\overset{\eqref{eq:16}}=1/|P_1Q|}}}^{\mathclap{\qquad\qquad\overset{\eqref{eq:16}}=\frac{|P_2T|-|P_1Q|}{|P_1Q||P_2T|}\overset{\eqref{eq:16}}=\frac{x_Qx_T}{d_0d_2L_0L_1}\big(|P_2T|-|P_1Q|\big)}}\big)=\frac{bx_Qx_T}{d_0d_2L^2_1}\big(|P_2T|-|P_1Q|\big).
\end{equation}
\begin{figure}[ht]
\begin{tikzpicture}
    \tkzDefPoint(0,0){Or};
    \tkzDefPoint(0,-4.5){Y};
    \pgfmathsetmacro{\rradius}{4.5};
    \pgfmathsetmacro{\bdist}{19};
    \tkzDefPoint(0,\bdist){OR};
    \pgfmathsetmacro{\phistardeg}{40};
    \pgfmathsetmacro{\XValueArc}{\rradius*sin(\phistardeg)};
    \pgfmathsetmacro{\YValueArc}{\rradius*cos(\phistardeg)};
    \tkzDefPoint(\XValueArc,-1.0*\YValueArc){B};
    \tkzDefPoint(-1.0*\XValueArc,-1.0*\YValueArc){A};
    \pgfmathsetmacro{\Rradius}{veclen(\XValueArc,\YValueArc+\bdist)};
    \clip(-\rradius,-\rradius-0.3) rectangle (\rradius,0.3);
    \draw[name path=Cr,ultra thin,dashed] (Or) circle (\rradius);
    \tkzDrawArc[name path=CR](OR,A)(B);
    \tkzDrawArc[](Or,B)(A);
    \pgfmathsetmacro{\PXValue}{\Rradius*sin(3)};
    \pgfmathsetmacro{\PYValue}{\bdist-\Rradius*cos(3)};
    \tkzDefPoint(\PXValue,\PYValue){P};
    \path[name path=lineQ] (P) -- +(18:3.5);
    \path[name path=lineQL] (P) -- +(198:4);
    \path[name path=lineT] (P) -- +(-12:4);
    \path[name path=lineTL] (P) -- +(-192:5);
    \path[name intersections={of=Cr and lineQ,by={Q}}];
    \path[name intersections={of=Cr and lineQL,by={QL}}];
    \path[name intersections={of=Cr and lineT,by={T}}];
    \path[name intersections={of=Cr and lineTL,by={TL}}];
    \path[name path=lineQLH] (Q) -- +(-8,0);
    \path[name path=lineTLH] (T) -- +(-8,0);
    \path[name intersections={of=Cr and lineQLH,by={QLH}}];
    \path[name intersections={of=Cr and lineTLH,by={TLH}}];
    \draw[red,ultra thin](Q)--(QL);
    \draw[blue,ultra thin](T)--(TL);
    \draw[name path=QQLH,red,ultra thin,dashed](Q)--(QLH);
    \draw[name path=TTLH,blue,ultra thin,dashed](T)--(TLH);
    \draw[purple,ultra thin,dashed](Or)--(Y);
    \coordinate (PQ) at ($(P)-(Q)$);
    \coordinate (QQL)at ($(Q)-(QL)$);
    \coordinate (QQLH) at ($(Q)-(QLH)$);

    \coordinate (PT) at ($(P)-(T)$);
    \coordinate (TTL)at ($(T)-(TL)$);
    \coordinate (TTLH) at ($(T)-(TLH)$);
    

    \draw[name path=lineP1,red,dashed] let \p1=($(PQ)$),
    \p2=($(QQL)$),
    \p3=($(QQLH)$),
    \n1={veclen(\x1,\y1)},
    \n2={veclen(\x2,\y2)},
    \n3={veclen(\x3,\y3)}
    in (P)--([shift={(-\n1*\n2/\n3,0)}]Q);
    
    \path[name intersections={of=QQLH and lineP1,by={P1}}];

    \draw[name path=lineP2,blue,dashed] let \p1=($(PT)$),
    \p2=($(TTL)$),
    \p3=($(TTLH)$),
    \n1={veclen(\x1,\y1)},
    \n2={veclen(\x2,\y2)},
    \n3={veclen(\x3,\y3)}
    in (P)--([shift={(-\n1*\n2/\n3,0)}]T);

     \path[name intersections={of=TTLH and lineP2,by={P2}}];

    \tkzDrawPoints(Or,P,Q,QL,QLH,A,B,P1,P2,T,TL,TLH,Y);
    \draw[red,ultra thin,dashed](Or)--(T);
    \draw[red,ultra thin,dashed](Or)--(Q);
    
    
    \tkzLabelPoint[left,xshift=-3](Or){$O_r$};
    \tkzLabelPoint[above left,yshift=0.1,xshift=-0.3](P){$P$};
    \tkzLabelPoint[right,xshift=1](Q){$Q$};
    \tkzLabelPoint[right,xshift=1](T){$T$};
    \tkzLabelPoint[above,yshift=1](TL){$T_0$};
    \tkzLabelPoint[left,xshift=-1](QL){$Q_0$};
    \tkzLabelPoint[left,xshift=-1](QLH){$Q_1$};
    \tkzLabelPoint[left,xshift=-1](TLH){$T_1$};
    \tkzLabelPoint[above right](P1){$P_1$};
    \tkzLabelPoint[below](P2){$P_2$};
    \tkzLabelPoint[left](Y){$Y$};
    \tkzLabelPoint[left](A){$A$};
    \tkzLabelPoint[right](B){$B$};
    \tkzMarkAngle[arc=lll,ultra thin,size=1.7](Y,Or,Q);
    \tkzMarkAngle[blue,arc=lll,ultra thin,size=0.7](Y,Or,T);
    \tkzMarkAngle[arc=lll,ultra thin,size=0.3](QLH,P1,P);
    \tkzMarkAngle[blue,arc=lll,ultra thin,size=0.3](T,P2,P);
    \node at ($(Or) + (-65:1.9)$)  {\Large $\alpha$};
    \node at ($(Or) + (-78:1.0)$)  {
    \Large \textcolor{blue}{$\beta$}};
    \node at ($(P1) + (-160:0.5)$)  {$\alpha$};
    \node at ($(P2) + (-20:0.45)$)  {
    \textcolor{blue}{$ \beta$}};
    \draw[name path=lineE,blue,dashed] (T)-- +(1.5,0);
    \tkzDefLine[perpendicular=through Q](T,TLH)\tkzGetPoint{QY};
    \path[name path=LineQY](Q)--(QY);
    \path[name intersections={of=LineQY and lineE,by={E}}];
    \draw[red,dashed,ultra thin](Q) -- (E);    
    \path[name path=LinePY](P)--+(90:1.5);
    \path[name path=LineP2Y](P2)--+(90:1.5);
    \path[name intersections={of=LinePY and QQLH,by={F}}];
    \path[name intersections={of=LineP2Y and QQLH,by={G}}];
    \draw[red,dashed,ultra thin](P)--(F);
    \draw[red,dashed,ultra thin](P2)--(G);
    \tkzMarkRightAngle(T,E,Q);
    \tkzMarkRightAngle(P,F,QLH);
    \tkzMarkRightAngle(P2,G,QLH);
    \tkzDrawPoints(E,F,G);
    \tkzLabelPoint[below](E){$E$};
    \tkzLabelPoint[above](F){$H$};
    \tkzLabelPoint[above](G){$G$};
    \path[name path=LineOrY](Or)--(Y);
    \path[name intersections={of=LineOrY and CR,by={M}}];
    \tkzLabelPoint[below](M){$M$};
    \tkzDrawPoints(M);
\end{tikzpicture}
\caption{$\tau_1$ and $2d_0-\tau_0$ ratio analysis with $\beta<\alpha$}\label{fig:lengthfraction}
\end{figure}
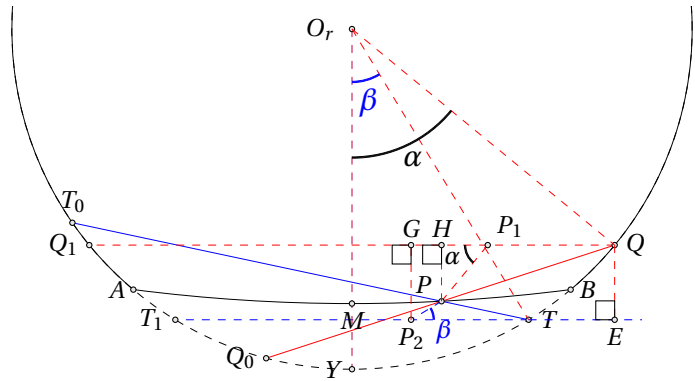
We now compute this last expression (see \eqref{eq:23}) and then show that it is negative, as required.

In \cref{fig:lengthfraction}, $QQ_1,TT_1$ are the horizontal lines, and $QE\perp TE$, $PH\perp HQ$, $P_2G\perp GQ$.
The $Q,T$ coordinates satisfy:
\begin{equation}\label{eq:hthqyqytyp}
    \begin{aligned}
    -r\cos{\alpha}&=y_Q=y_P+h_Q,\\
    -r\cos{\beta}&=y_T=y_P-h_T.
    \end{aligned}
\end{equation}
When $(\Phi,\theta)=x_1\in\Mrin$, $P=p(x_1)$\COMMENT{$(\Phi,\theta)=x_1\in\Mrin,P=p(x_1)$ must be replaced by $(\Phi,\theta)=x_1\in\Mrin$, $P=p(x_1)$, and likewise everywhere else throughout the text, and also with periods. Commas and periods must be followed by space or line break.\newline\red\small\wentao{Agree, Work in progress}}, then $\theta>\Phi$ from either \cref{fig:MR} or \cref{fig:R1}. Hence for $P=p(x_1)$, $\overrightarrow{PT}$ always has negative slope, and $h_Q>0$, $h_T>0$. Then \cref{fig:lengthfraction} further shows that
\begin{equation}\label{eq:anglelengthobservationforalphabeta}
\begin{aligned}
(i)&\quad |P_2G|=|QE|=h_T+h_Q,\\
(ii)&\quad |PH|=h_Q,|HP_1|=h_Q\cot{\alpha},\\
(iii)&\quad |GH|=h_T\cot{\beta},\\
(iv)&\quad \measuredangle{QTE}=\pi-\measuredangle{T_1TO_r}-\measuredangle{O_rTQ}=\pi-(\frac{\pi}{2}-\beta)-(\frac{\pi}{2}-\frac{\alpha-\beta}{2})=\frac{\alpha+\beta}{2},\\
(v)&\quad |TE|\overbracket{=}^{\mathllap{\text{ \cref{fig:lengthfraction}: }TE\perp QE}}|QE|\cot\angle{QTE}\overset{(iv)}=|QE|\cot{(\frac{\alpha+\beta}{2})}\overset{(i)}{=}(h_T+h_Q)\cot{(\frac{\alpha+\beta}{2})},\\
(vi)&\quad |GH|+|HP_1|+|P_1Q|=|GQ|=|P_2E|=|P_2T|+|TE|.
\end{aligned}
\end{equation}
From equations \eqref{eq:anglelengthobservationforalphabeta} and\eqref{eq:hthqyqytyp} and some trigonometry, we further obtain:
\begin{equation}\label{eq:21}
\begin{aligned}
|P_2T|-|P_1Q| &\overbracket{=}^{\mathclap{\text{\eqref{eq:anglelengthobservationforalphabeta} (vi)}}} |GH|+|HP_1|-|TE|\overbracket{=}^{\mathclap{\text{\eqref{eq:anglelengthobservationforalphabeta} (ii)(iii)(v)}}}h_T\cot{\beta}+h_Q\cot{\alpha}-(h_T+h_Q)\cot{(\frac{\alpha+\beta}{2})}\\
&\overbracket{=}^{\mathclap{\text{\eqref{eq:hthqyqytyp}}}} (y_P+r\cos{\beta})\cot{\beta}-(r\cos{\alpha}+y_P)\cot{\alpha}-(r\cos{\beta}-r\cos{\alpha})\cot{(\frac{\alpha+\beta}{2})}\\
&=(y_P+r\cos{\beta})\cot{\beta}-(r\cos{\alpha}+y_P)\cot{\alpha}+2r\sin{(\frac{\beta-\alpha}{2})}\sin{(\frac{\beta+\alpha}{2})}\cot{(\frac{\alpha+\beta}{2})}\\
&=(y_P+r\cos{\beta})\cot{\beta}-(r\cos{\alpha}+y_P)\cot{\alpha}+2r\sin{(\frac{\beta-\alpha}{2})}\cos{(\frac{\alpha+\beta}{2})}\\
&=(y_P+r\cos{\beta})\cot{\beta}-(r\cos{\alpha}+y_P)\cot{\alpha}-r\sin\alpha+r\sin\beta\\
&=y_P\cot{\beta}-y_P\cot{\alpha}+r\frac{\cos^2\beta}{\sin\beta}-r\frac{\cos^2\alpha}{\sin\alpha}-r\sin\alpha+r\sin\beta\\
&=y_P\cot{\beta}-y_P\cot{\alpha}+r\frac{\cos^2\beta+\sin^2\beta}{\sin\beta}-r\frac{\cos^2\alpha+\sin^2\alpha}{\sin\alpha}\\
&= y_P\cot{\beta}-y_P\cot{\alpha}-r\frac{1}{\sin{\alpha}}+r\frac{1}{\sin{\beta}}=y_P\frac{\cos\beta}{\sin\beta}-y_P\frac{\cos\alpha}{\sin\alpha}+r\frac{\sin{\alpha}-\sin{\beta}}{\sin{\alpha}\sin{\beta}}\\
&= y_P\frac{\sin{\alpha}\cos{\beta}-\cos{\alpha}\sin{\beta}}{\sin{\alpha}\sin{\beta}}+r\frac{\sin{\alpha}-\sin{\beta}}{\sin{\alpha}\sin{\beta}}= y_P\frac{\sin{(\alpha-\beta)}}{\sin{\alpha}\sin{\beta}}+r\frac{\sin{\alpha}-\sin{\beta}}{\sin{\alpha}\sin{\beta}}\\
&=y_P\frac{\sin{(\alpha-\beta)}}{\sin{\alpha}\sin{\beta}}+r\frac{2\sin{(\frac{\alpha-\beta}{2})}\cos{(\frac{\alpha+\beta}{2})}}{\sin{\alpha}\sin{\beta}}=y_P\frac{\sin{(\alpha-\beta)}}{\sin{\alpha}\sin{\beta}}+r\frac{2\sin{(\frac{\alpha-\beta}{2})}\cos{(\frac{\alpha+\beta}{2})\cos{(\frac{\alpha-\beta}{2})}}}{\sin{\alpha}\sin{\beta}\cos{(\frac{\alpha-\beta}{2})}}\\
&=y_P\frac{\sin{(\alpha-\beta)}}{\sin{\alpha}\sin{\beta}}+r\frac{\sin{(\alpha-\beta)}\cos{(\frac{\alpha+\beta}{2})}}{\sin{\alpha}\sin{\beta}\cos{(\frac{\alpha-\beta}{2})}}=\frac{\sin{(\alpha-\beta)}}{\sin{\alpha}\sin{\beta}}\big(y_P+r\frac{\cos{(\frac{\alpha+\beta}{2})}}{\cos{(\frac{\alpha-\beta}{2})}}\big).
\end{aligned}
\end{equation}
In \cref{fig:twoTangentLines} below,  the tangent lines at $T,Q$ intersect at $S$. Then by $|O_rS|=\frac{|O_rQ|}{\cos{\measuredangle{QO_rS}}}=\frac{r}{\cos{\frac{\alpha-\beta}{2}}}$ and $\measuredangle{MO_rS}=\frac{\alpha+\beta}{2}$
\begin{equation}\label{eq:xsys}
S=(x_S,y_S)=(|O_rS|\sin\angle{MO_rS},-|O_rS|\cos\angle{MO_rS})=(r\frac{\sin{(\frac{\alpha+\beta}{2})}}{\cos{(\frac{\alpha-\beta}{2})}},-r\frac{\cos{(\frac{\alpha+\beta}{2})}}{\cos{(\frac{\alpha-\beta}{2})}}),
\end{equation}
hence $|P_2T|-|P_1Q|\overset{\eqref{eq:21}}=\frac{\sin{(\alpha-\beta)}}{\sin{\alpha}\sin{\beta}}(y_P-y_S)$. Combine equations ~\eqref{eq:15}, \eqref{eq:18}, \eqref{eq:21}, and \eqref{eq:xsys} to obtain
\begin{equation}\label{eq:23}
\frac{d}{d\Phi}(\frac{2d_0-\tau_0}{\tau_1})=\frac{bx_Qx_T}{d_0d_2\tau^2_1}\frac{\sin{(\alpha-\beta)}}{\sin{\alpha}\sin{\beta}}(y_P-y_S)\overset{\eqref{eq:15}}=\frac{br^2\sin{(\alpha-\beta)}}{d_0d_2\tau^2_1}(y_P-y_S).
\end{equation}
\Cref{lemma:fractionLemma1,lemma:fractionLemma2} below and the Intermediate-Value Theorem imply the following. If $p(x_1)=P$ moves from the midpoint $M$ of $\arc{AB}$ to the corner $B$ with constant small crossing angle with the arc $\arc{AB}$, then $y_P\le y_S$ with equality  only when $P$ is at the mid-point $M$ or the corner $B$. Since $\alpha>\beta$, the right-hand side of \eqref{eq:23} is negative (except for two points), so $\frac{2d_0-\tau_0}{\tau_1}$ is a decreasing function of $\Phi\in[0,\Phistar]$.
\end{proof}

\begin{figure}[h]
\begin{tikzpicture}
    \tkzDefPoint(0,0){Or};
    \tkzDefPoint(0,-4.3){Y};
    \pgfmathsetmacro{\rradius}{4.5};
    \pgfmathsetmacro{\bdist}{19};
    \tkzDefPoint(0,\bdist){OR};
    \pgfmathsetmacro{\phistardeg}{40};
    \pgfmathsetmacro{\XValueArc}{\rradius*sin(\phistardeg)};
    \pgfmathsetmacro{\YValueArc}{\rradius*cos(\phistardeg)};
    \tkzDefPoint(\XValueArc,-1.0*\YValueArc){B};
    \tkzDefPoint(-1.0*\XValueArc,-1.0*\YValueArc){A};
    \pgfmathsetmacro{\Rradius}{veclen(\XValueArc,\YValueArc+\bdist)};
     \clip(-\rradius,-\rradius-0.3) rectangle (\rradius,0.3);
    \draw[name path=Cr,ultra thin,dashed] (Or) circle (\rradius);
    \tkzDrawArc[name path=CR](OR,A)(B);
    \tkzDrawArc[](Or,B)(A);
    \pgfmathsetmacro{\PXValue}{\Rradius*sin(3)};
    \pgfmathsetmacro{\PYValue}{\bdist-\Rradius*cos(3)};
    \tkzDefPoint(\PXValue,\PYValue){P};
    \path[name path=lineQ] (P) -- +(18:3.5);
    \path[name path=lineQL] (P) -- +(198:4);
    \path[name path=lineT] (P) -- +(-12:4);
    \path[name path=lineTL] (P) -- +(-192:5);
    \path[name intersections={of=Cr and lineQ,by={Q}}];
    \path[name intersections={of=Cr and lineQL,by={QL}}];
    \path[name intersections={of=Cr and lineT,by={T}}];
    \path[name intersections={of=Cr and lineTL,by={TL}}];
    \draw[red,ultra thin,dashed] (P) --(Q);
    \draw[blue,ultra thin,dashed] (P) --(T);
    \draw[red,ultra thin,dashed] (Or) --(Q);
    \draw[blue,ultra thin,dashed] (Or) --(T);
    \tkzDefLine[orthogonal=through Q](Q,Or)\tkzGetPoint{QI};
    \tkzDefLine[orthogonal=through T](Or,T)\tkzGetPoint{TI};
    \path[name path=ST] (T)--(TI);
    \path[name path=SQ] (Q)--(QI);
    \path[name intersections={of=ST and SQ,by={S}}];
    \draw[blue,ultra thin,dashed] (S) --(T);
    \draw[red,ultra thin,dashed] (S) --(Q);
    \draw[red,ultra thin,dashed] (Or)--(Y);
    \draw[blue,ultra thin,dashed](TL)--(P);
    \tkzMarkRightAngle(Or,T,S);
    \tkzMarkRightAngle(Or,Q,S);
    \tkzLabelPoint[below](P){$P$};
    \tkzLabelPoint[below](T){$T$};
    \tkzLabelPoint[right](Q){$Q$};
    \tkzLabelPoint[right](S){$S$};
    \tkzLabelPoint[left](Or){$O_r$};
    \tkzLabelPoint[below](Y){$Y$};
    \tkzLabelPoint[right](TL){$T_0$};
    \path[name path=LineOrY](Or)--(Y);
    \path[name intersections={of=LineOrY and CR,by={M}}];
    \path[name path=LineHS](S)--+(-5,0);
    \path[name intersections={of=LineHS and LineOrY,by={M0}}];
    \draw[red,ultra thin,dashed](M0)--(S);
    \tkzLabelPoint[below](M){$M$};
    \tkzLabelPoint[below left](A){$A$};
    \tkzLabelPoint[below](B){$B$};
    \tkzLabelPoint[above left](M0){$M_0$};
    \tkzDrawPoints(Or,P,Q,T,S,M,A,B,TL,M0);
    \tkzMarkAngle[arc=lll,ultra thin,size=1.7](Y,Or,Q);
    \tkzMarkAngle[blue,arc=lll,ultra thin,size=0.7](Y,Or,T);
    \node at ($(Or) + (-65:1.9)$)  {\Large $\alpha$};
    \node at ($(Or) + (-78:1.0)$)  {
    \Large \textcolor{blue}{$\beta$}};
\end{tikzpicture}
\caption{Tangent lines at $T$ and $Q$ intersecting at $S$}\label{fig:twoTangentLines}
\end{figure}
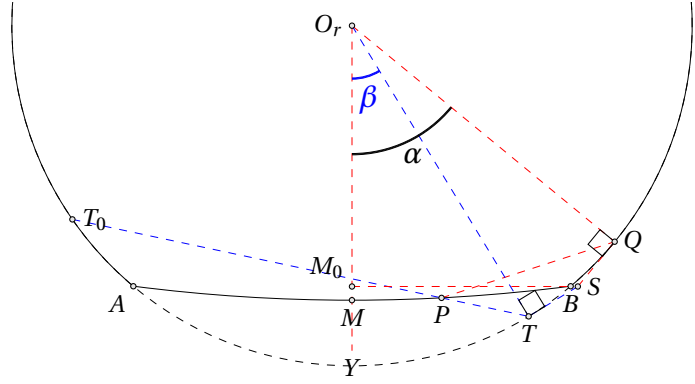

\begin{lemma}[$y_P\ne y_S$ in \cref{fig:twoTangentLines}]\label{lemma:fractionLemma1}With the notations in \cref{def:LemonTableConfiguration}, suppose that our lemon billiard satisfies the condition \eqref{eqJZHypCond} with the coordinate system in \cref{def:StandardCoordinateTable} such that $O_r=(0,0)$ and $O_R=(0,b)$. In \cref{fig:twoTangentLines}, let $P$, $Q$, $T_0$ be the positions in the table of $x_1=(\Phi,\theta)\in \Mrin$, $x_2=\Fc(x_1)\in \Mrin$, $x_0=\Fc^{-1}(x_1)\in \Mrout$ and $M$ be the midpoint of $\Gamma_R$. Suppose $T$ is the intersection of the extension of the line $\overrightarrow{T_0P}$ with the circle $C_r$. Let $S=(x_S,y_S)$ be the intersection of the tangent line of circle $C_r$ at $T$ with the tangent line of circle $C_r$ at $Q$. If $\overrightarrow{PQ}$ moves on $\arc{MB}$ with a constant crossing angle $\theta\in[\sin^{-1}{(2r/R)},\sin^{-1}{(\sqrt{4r/R})}]$ with the tangent direction of $\Gamma_R$ at $P=(x_P,y_P)$ and $P$ is the interior point of the arc segment $\arc{MB}$, then $y_P\ne y_S$.
\end{lemma}
\begin{proof}
In \cref{fig:twoTangentLines}, suppose that the horizontal line through $S$ intersects the $y$-axis at $M_0$. From the coordinates of $Q$,$T$,$S$ in \cref{fig:twoTangentLines}, the equations~\eqref{eq:15} and \eqref{eq:xsys}, we have $M_0=(0,y_S)=(0,-r\frac{\cos{(\frac{\alpha+\beta}{2})}}{\cos{(\frac{\alpha-\beta}{2})}})$.

Note that since $\measuredangle{O_rM_0S}=\measuredangle{O_rTS}=\pi/2$, $O_r,M_0,T,S$ are all on the circle with $O_rS$ as its diameter. Hence $\measuredangle{TM_0S}=\measuredangle{TO_rS}=\frac{\alpha-\beta}{2}$, because the angles at the circumference subtended by the same arc are equal and $S$ is on the bisector of $\angle{TO_rQ}$.

Similarly, since $\measuredangle{O_rM_0S}=\measuredangle{O_rQS}=\pi/2$, $O_r,M_0,S,Q$ are all on the circle with $O_rS$ as its diameter. Hence $\measuredangle{QM_0S}=\measuredangle{QO_rS}=\frac{\alpha-\beta}{2}$. Therefore, 
\begin{equation}\label{eq:AngeEquationForcontradiction1}
\measuredangle{QM_0S}=\measuredangle{TM_0S}=\frac{\alpha-\beta}{2}
\end{equation} in \cref{fig:twoTangentLines} and \cref{fig:M0PS}.

To prove that $y_P\ne y_S$ in \Cref{fig:twoTangentLines} we now assume $y_S=y_P$ in \cref{fig:M0PS} for purposes of contradiction.\COMMENT{Look for further shortcuts\wentao{It's quite concise now and still thinking.}} This means that $P$ is on the line segment $M_0S$.
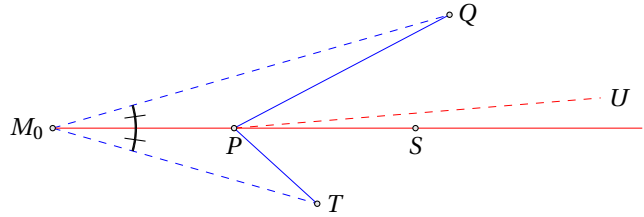
\begin{figure}[ht]
\begin{center}
\begin{tikzpicture}
\tkzDefPoint(-3,0){M0};
\tkzDefPoint(-0.6,0){P};
\tkzDefPoint(1.8,0){S};
\tkzDefPoint(2.25,1.5){Q};
\tkzDefPoint(0.5,-1){T};
\tkzDefPoint(4.25,0.4){U};
\draw[blue,ultra thin,dashed](M0)--(T);
\draw[blue,ultra thin,dashed](M0)--(Q);
\draw[red,ultra thin](M0)--(S)--+(3.0,0);
\draw[red,ultra thin,dashed](P)--(U);
\draw[blue,ultra thin](P)--(Q);
\draw[blue,ultra thin](P)--(T);
\tkzDrawPoints(M0,P,S,Q,T);
\tkzLabelPoint[left](M0){$M_0$};
\tkzLabelPoint[below](P){$P$};
\tkzLabelPoint[below](S){$S$};
\tkzLabelPoint[right](Q){$Q$};
\tkzLabelPoint[right](T){$T$};
\tkzLabelPoint[right](U){$U$};
\tkzMarkAngle[arc=lll,ultra thin,size=1.1,mark=|](P,M0,Q);
\tkzMarkAngle[arc=lll,ultra thin,size=1.1,mark=|](T,M0,P);
\end{tikzpicture}
\end{center}
    \caption{Assume $P$ were on the line  $M_0S$ which is a horizontal line}\label{fig:M0PS}
\end{figure}
In \cref{fig:M0PS} denote by $PU$ the tangent line at P of the billiard boundary $\Gamma_R=\arc{AB}$. $PU$ has positive slope and bisects $\angle{QPT}$.

Since $\theta>\Phi$, $PT$ has a negative slope, while $PQ$ has a positive slope. And since $\sin^{-1}(2r/R)\le\theta\le\sin^{-1}{(\sqrt{4r/R})}$, $0\le\Phi\le\Phistar$ are both small, it is clear that $\measuredangle{TPM_0}=\pi-\theta+\Phi$, $\measuredangle{QPM_0}=\pi-\theta-\Phi$ are always obtuse angles when $P$ is on $\arc{MB}$ ($\Phi$ varies in $[0,\Phistar]$). Then by the sine law, we get the followng.
\begin{equation}\label{eq:AngeEquationForcontradiction2}
\frac{\sin{\measuredangle{PTM_0}}}{|M_0P|}=\frac{\sin{\measuredangle{PM_0T}}}{|PT|}\text{ and } \frac{\sin{\measuredangle{PQM_0}}}{|M_0P|}=\frac{\sin{\measuredangle{PM_0Q}}}{|PQ|}.
\end{equation}

By assumption $M_0$, $P$, $S$ are on the same line,  \eqref{eq:AngeEquationForcontradiction1} gives $\measuredangle{PM_0Q}=\measuredangle{PM_0T}$. By \cref{remark:symmetrywiththetangentline}: $|PT|=2d_0-\tau_0<\tau_1=|PQ|$, \eqref{eq:AngeEquationForcontradiction2} givies $\sin{\measuredangle{PTM_0}}=|M_0P|\frac{\sin{\measuredangle{PM_0T}}}{|PT|}>|M_0P|\frac{\sin{\measuredangle{PM_0Q}}}{|PQ|}=\sin{\measuredangle{PQM_0}}$. Since $\angle{PTM_0}$ and $\angle{PQM_0}$ are acute, we get $\measuredangle{PTM_0}>\measuredangle{PQM_0}$. 

Therefore, by \eqref{eq:AngeEquationForcontradiction1}: $\measuredangle{TM_0S}=\measuredangle{QM_0S}$, we also get
\begin{equation}\label{eq:27}
\measuredangle{TPS}=\measuredangle{PTM_0}+\measuredangle{TM_0S}>\measuredangle{PQM_0}+\measuredangle{QM_0S}=\measuredangle{QPS}.
\end{equation}
But $PU$ is the bisector of $\angle{QPT}$ with positive slope, i.e., $\measuredangle{SPU}>0$, hence
\begin{equation}\label{eq:28}
\begin{aligned}
\measuredangle{QPS}&=\measuredangle{QPU}+\measuredangle{SPU}>\measuredangle{QPU}=\measuredangle{TPU}>\measuredangle{TPU}-\measuredangle{SPU}=\measuredangle{TPS},
\end{aligned}
\end{equation}
is contrary to \eqref{eq:27}, so $P$, $S$ cannot be on the same horizontal line ($P$ cannot be on the line $M_0P$), hence $y_P\ne y_S$.
\end{proof}
\begin{lemma}[$y_S-y_P$ is increasing at $\Phi=0$ (for \cref{monotone_p5})]\label{lemma:fractionLemma2}
In the context of \cref{lemma:fractionLemma1}, that is, $S=(x_S,y_S)$ as the intersection of tangent lines of $C_r$ at $Q$ and $T$, suppose that $\overrightarrow{PQ}$ moves in $\arc{MB}$ with a constant crossing angle $\theta\in[\sin^{-1}{(2r/R)},\sin^{-1}{(\sqrt{4r/R})}]$ with the tangent direction of $\Gamma_R$ in $P$. Then $P=(x_P,y_P)=(R\sin{\Phi},b-R\sin{\Phi})$ with $\Phi\in[0,\Phistar)$, and $y_P,y_S$ from \eqref{eq:xsys} are continuous differentiable functions of $\Phi$ such that $\frac{dy_P}{d\Phi}\bigr|_{
\Phi=0}<\frac{dy_S}{d\Phi}\bigr|_{\Phi=0}$, that is, $y_S-y_P$ is increasing at $M$.
\end{lemma}
\begin{proof}
When $\Phi=0$, i.e., $P=M$ is at the midpoint of the arc $\Gamma_R$, then $x_P=0$, so \[\frac{dy_P}{d\Phi}\bigr|_{\Phi=0}=\frac{d(b-R\cos{\Phi})}{d\Phi}=R\sin{\Phi}=x_P=0.\] 
It remains to show that $\frac{dy_S}{d\Phi}\bigr|_{\Phi=0}>0$, and this is a brute-force computation.
First, differentiate \eqref{eq:xsys} with respect to $\Phi$, using the chain rule and trigonometry:
\begin{equation}\label{eq:30}
\begin{aligned}
\frac{dy_S}{d\Phi}&=\frac{d}{d\Phi}\Big[\frac{-r\cos{(\frac{\alpha+\beta}{2})}}{\cos{(\frac{\alpha-\beta}{2})}}\Big]\\&=\frac{-r}{\cos^{2}{(\frac{\alpha-\beta}{2})}}\Big[-\cos{(\frac{\alpha-\beta}{2})}\sin{(\frac{\alpha+\beta}{2})}\frac{d}{d\Phi}(\frac{\alpha+\beta}{2})+\cos{(\frac{\alpha+\beta}{2})}\sin{(\frac{\alpha-\beta}{2})}\frac{d}{d\Phi}(\frac{\alpha-\beta}{2})\Big]\\&=\frac{r}{2\cos^{2}{(\frac{\alpha-\beta}{2})}}\Big[\sin{(\frac{\alpha+\beta}{2})}\cos{(\frac{\alpha-\beta}{2})}-\sin{(\frac{\alpha-\beta}{2})}\cos{(\frac{\alpha+\beta}{2})}\Big]\frac{d\alpha}{d\Phi}\\
& \qquad+\frac{r}{2\cos^{2}{(\frac{\alpha-\beta}{2})}}\Big[\sin{(\frac{\alpha+\beta}{2})}\cos{(\frac{\alpha-\beta}{2})}+\sin{(\frac{\alpha-\beta}{2})}\cos{(\frac{\alpha+\beta}{2})}\Big]\frac{d\beta}{d\Phi}\\
& =\frac{r}{\cos^{2}{(\frac{\alpha-\beta}{2})}}\Big(\frac{1}{2}\sin{\beta}\frac{d\alpha}{d\Phi}+\frac{1}{2}\sin{\alpha}\frac{d\beta}{d\Phi}\Big)
\end{aligned}
\end{equation}
We next find the derivatives of \(\alpha\) and \(\beta\), with $L_1=\tau_1,L_0=2d_0-\tau_0$ from \cref{def:fpclf}.
In \cref{fig:twoTangentLines},
\begin{equation}\label{eq:29}
\begin{split}
& Q=(x_Q,y_Q) =(r\sin{\alpha},-r\cos{\alpha}),\\
& T=(x_T,y_T) =(r\sin{\beta},-r\cos{\beta}),
\end{split}
\end{equation}
where $\pi/2\ge\alpha>\beta>0$ are smooth functions of $\Phi$.
Also in \cref{fig:twoTangentLines},  $\alpha-\theta-\Phi$ is the angle between $\overrightarrow{PQ}$ and the tangent $\overrightarrow{SQ}$, so 
\begin{equation}\label{eqd2}
    d_2=r\sin{(\alpha-\theta-\Phi)}.
\end{equation} 
Observing in \cref{fig:CaseAMonotone,fig:twoTangentLines}, we have the following coordinate equations for the points $P,Q$. 

\begin{gather}\label{eq:gatheredcoordinatesEquationsPQT_1}
r\sin{\alpha}=x_Q=x_P+L_1\cos{(\Phi+\theta)}=R\sin\Phi+L_1\cos{(\Phi+\theta)}\\\label{eq:gatheredcoordinatesEquationsPQT_2}
\overbracket{b-y_Q}^{\mathllap{y\text{-coordinate difference between } O_R\text{ and }Q.}}=\overbracket{-L_1\sin{(\Phi+\theta)}}^{\mathrlap{y\text{-coordinate difference between } P\text{ and }Q.}}+\underbracket{R\cos\Phi}_{\mathclap{y\text{-coordinate difference between } O_R\text{ and }P.}}
\end{gather}
Since $\theta$ is held as a constant, differentiating \eqref{eq:gatheredcoordinatesEquationsPQT_1} with respect to $\Phi$ in leftmost and rightmost sides, using \Cref{monotone_p3}\eqref{itemmonotone_p3-1}, \eqref{eq:29}, chain rule and trigonometry give
\begin{align*}
r\cos{\alpha}\frac{d\alpha}{d\Phi}&\overbracket{=}^{\mathclap{\substack{\text{Differentiate \eqref{eq:gatheredcoordinatesEquationsPQT_1}}\\|}}}\frac{d}{d\Phi}(R\sin{\Phi})-L_1\sin{(\Phi+\theta)}+\frac{dL_1}{d\Phi}\cos{(\Phi+\theta)}\\
&\overbracket{=}^{\mathclap{\text{\Cref{monotone_p3}\eqref{itemmonotone_p3-1}}}}R\cos{\Phi}-L_1\sin{(\Phi+\theta)}-\frac{bx_Q}{d_2}\cos{(\Phi+\theta)}\overbracket{=}^{\mathclap{\eqref{eq:gatheredcoordinatesEquationsPQT_2}}}-y_Q+b-\frac{b\sin\alpha}{\sin{(\alpha-\theta-\Phi)}}\cos{(\Phi+\theta)}\\
&\overbracket{=}^{\mathclap{\text{\eqref{eq:29}}}}r\cos{\alpha}-b\frac{\sin\alpha\cos{(\Phi+\theta)}-\sin{(\alpha-\theta-\Phi)}}{\sin{(\alpha-\theta-\Phi)}}=r\cos{\alpha}-b\frac{\sin{(\Phi+\theta)}\cos{\alpha}}{\underbracket{\sin{(\alpha-\theta-\Phi)}}_{\overset{\eqref{eqd2}}{=}d_2/r}},
\end{align*}
so 
\begin{equation}\label{eq:32}
\frac{d\alpha}{d\Phi}=
1-\frac{b\sin{(\Phi+\theta)}}{d_2}.
\end{equation}
In \cref{fig:twoTangentLines}, $\beta+\theta-\Phi$ is the angle between  $\overrightarrow{TS}$, and $\overrightarrow{T_0T}$, so
\begin{equation}\label{eqd0}
    d_0=r\sin{(\beta+\theta-\Phi)}.
\end{equation}
\cref{fig:CaseAMonotone,fig:twoTangentLines} give the following coordinate equations for the points $P,Q$. 
\begin{gather}\label{eq:gatheredcoordinatesEquationsPQT_3}
r\sin{\beta}=x_T=x_P+L_0\cos{(\theta-\Phi)}=R\sin\Phi+L_0\cos{(\theta-\Phi)}\\\label{eq:gatheredcoordinatesEquationsPQT_4}
\underbracket{b-y_T}_{\mathllap{y\text{-coordinate difference between }O_R\text{ and }T}}=\overbracket{R\cos\Phi}^{\mathclap{y\text{-coordinate difference between }O_R\text{ and }P}}+\underbracket{L_0\sin{(\theta-\Phi)}}_{\mathclap{y\text{-coordinate difference between }P\text{ and }T}}
\end{gather}
Since $\theta$ is held constant, differentiating \eqref{eq:gatheredcoordinatesEquationsPQT_3} with respect to $\Phi$ in the left and right sides, using \Cref{monotone_p3}\eqref{itemmonotone_p3-2}, \eqref{eq:29}, chain rule and trigonometry will further give
\[
\begin{aligned}
    r\cos\beta\frac{d\beta}{d\Phi}&\overbracket{=}^{\mathclap{\substack{\text{Differentiate \eqref{eq:gatheredcoordinatesEquationsPQT_3}}\\|}}}\frac{d}{d\Phi}(R\sin\Phi)+L_0\sin{(\theta-\Phi)}+\cos{(\theta-\Phi)}\frac{dL_0}{d\Phi}\\
    &\overbracket{=}^{\mathclap{\text{\Cref{monotone_p3}\eqref{itemmonotone_p3-2}}}}R\cos\Phi+L_0\sin{(\theta-\Phi)}-\frac{bx_T}{d_0}\cos{(\theta-\Phi)}\overbracket{=}^{\mathclap{\eqref{eq:gatheredcoordinatesEquationsPQT_4},\eqref{eqd0}}}b-y_T-\frac{bx_T\cos{(\theta-\Phi)}}{r\sin{(\beta+\theta-\Phi)}}\\
    &\overbracket{=}^{\mathclap{\eqref{eq:29}}}r\cos\beta+b-\frac{b\sin\beta\cos{(\theta-\Phi)}}{\sin{(\beta+\theta-\Phi)}}=r\cos\beta+b\frac{\sin{(\beta+\theta-\Phi)}-\sin\beta\cos{(\theta-\Phi)}}{\sin{(\beta+\theta-\Phi)}}\\
    &=r\cos\beta+b\frac{\cos\beta\sin{(\theta-\Phi)}}{\underbracket{\sin{(\beta+\theta-\Phi)}}_{\overset{\eqref{eqd0}}=d_0/r}},
\end{aligned}
\]
so
\begin{equation}\label{eq:33}
\frac{d\beta}{d\Phi}
{=}1+b\frac{\sin{(\theta-\Phi)}}{d_0}.
\end{equation}
Inserting \eqref{eq:32}, \eqref{eq:33} into \eqref{eq:30}, we get
\begin{equation}\label{eq:34}
\frac{dy_S}{d\Phi}=\frac{r}{2\cos^2{(\frac{\alpha-\beta}{2})}}\big[\sin{\beta}\frac{d_2-b\sin{(\theta+\Phi)}}{d_2}+\sin{\alpha}\frac{d_0+b\sin{(\theta-\Phi)}}{d_0}\big].
\end{equation}
If $\Phi=0$, then $d_0=d_2$, and \eqref{eq:34} becomes
\[
\frac{dy_S}{d\Phi}\bigr|_{\Phi=0}=\frac{r}{2\cos^2{(\frac{\alpha-\beta}{2})}}\Big[\sin\alpha+\sin\beta+\frac{b\sin\theta}{d_0}(\sin\alpha-\sin\beta)\Big]\overset{\pi/2\ge\alpha>\beta>0}>0,
\]
so 
$\displaystyle\big(\frac{dy_S}{d\Phi}-\frac{dy_P}{d\Phi}\big)\bigr|_{\Phi=0}=\frac{dy_S}{d\Phi}\bigr|_{\Phi=0}>0$, proving \Cref{lemma:fractionLemma2}.
\end{proof}


\section{Expansion in Case (a)}\label{sec:ProofTECase_a0_a1}
\subsection[Coordinates change, results from \texorpdfstring{\cite{cb}}{[CM06]} and expansion in the invariant cone by shears]{Coordinates change, results from \texorpdfstring{\cite{cb}}{[CM06]} and expansion in the invariant cone by shears}\label{SBSCoordChange}\hfill\strut\COMMENT{It would help to explain what happens where in this section}

This subsection 
gathers up some notions and results from \cite{cb}\COMMENT{Why mention \cite{hoalb}? All results are from \cite{cb}.\wentao{Yes, OK}} about successive collisions with a given circle arc that will be useful later.

Suppose $x=(\phi,\theta)\in M_r,dx=(d\phi,d\theta)\in T_x(M)$, In the arc-length parameter $s$ and angle with inward normal vector parameter $\varphi$, the coordinate change from $(\phi,\theta)$ to $(s,\varphi)$ is given by

\begin{equation}\label{eq:coordinateschange}
\begin{aligned}
\theta &= \frac{\pi}{2}-\varphi,\\
d\varphi &= -d\theta,\\
ds &= rd\phi.
\end{aligned}
\end{equation}
\cite[Equations (3.39), (3.40)]{cb} give the mirror equation for pre/post-collision wave front curvatures
$\mathcal{B}^{\pm}$ and the expansion of tangent vectors in the p-metric $\|dx\|_{p}=\cos{\varphi}|ds|=\sin{\theta}|ds|$:
\begin{equation}\label{eq:MirrorEquationExpansionFormula}
\begin{aligned}
-\frac{d\theta}{rd\phi}&\overbracket{=}^{\strut\mathclap{\eqref{eq:coordinateschange}}}\,\frac{d\varphi}{ds}=\mathcal{V}=\mathcal{B}^{-}\cos{\varphi}+\mathcal{K}=\mathcal{B}^{+}\cos{\varphi}-\mathcal{K},\\
\mathcal{B}^{+}&=\mathcal{B}^{-}+\frac{2\mathcal{K}}{\cos{\varphi}},\\
\frac{\|D_x\Fc(dx)\|_\p}{\|dx\|_\p}&=|1+\tau\mathcal{B^+}|,
\end{aligned}
\end{equation}
where $\mathcal{K}=-{1}/{r}$ is the curvature of the boundary\COMMENT{\red\small Search the source code for "'s". Almost all instances should be edited out. In the present case, "the billiard boundary's curvature" should become "the curvature of the billiard boundary". Again, this is not an isolated comment. If you find any occurrence of "'s" which you don't know how to fix, flag it so we can discuss it.} and $\tau$ is the distance between $p(x)$ and $p(\Fc(x))$.\COMMENT{\textcolor{red}{You sometimes forget to use math mode; a thing to fix.\wentao{OK}}}

The study of the stadium in \cite{cb} produces several insights about successive collisions with a given circle arc. We reproduce some of these.

\begin{remark}\label{proposition:ShearOneStepExpansion}
If $(\phi,\theta)=x\in M_r$, $\Fc(x)\in M_r$ are nonsingular, then 
$D_{x}\Fc=\begin{pmatrix}1 & 2\\ 0 & 1\end{pmatrix}$ in $\phi\theta$-coordinates, so the positive quadrant is preserved. Furthermore, if $dx$ is a tangent vector at $x$, then  $\tau=2r\cos{\varphi}$ gives 
\[\frac{\|D_x\Fc(dx)\|_\p}{\|dx\|_\p}=|1+\tau \Bcal^+|=|1+2r\cos{\varphi}\frac{\frac{-d\theta}{rd\phi}-\frac{1}{r}}{\cos{\varphi}}|=|-1-2\frac{d\theta}{d\phi}|.\]
Thus
if $\frac{d\theta}{d\phi}\in[c,+\infty]$ with $c\ge0$, then 
$\frac{\|D_x\Fc(dx)\|_\p}{\|dx\|_\p}\ge1+2c$.
\end{remark}
\begin{proposition}[{\cite[p.\ 223]{cb}}]\label{proposition:ImaginaryExpansion}
{
Consider an infinitesimal wave front that moves along the trajectory and experiences $k+1$ successive collisions with $\Gamma_r$. Denote by $\mathcal{B}^{-}_{0}$ and $\mathcal{B}^{-}_{i}$ the wave-front curvature before the very first collision and before the $i$-th collision, respectively, and by $\tau$ the distance between consecutive collision points. Then the total expansion (with respect to the p-metric) from the very first to the very last collision is
\begin{equation}\label{eqJfromCM}
    \mathcal{J}\dfn\prod_{i=0}^{k-1}|1+\tau(-\frac{4}{\tau}+\mathcal{B}_{i}^{-})|=\big|1-k\tau(-\frac{2}{\tau}+\mathcal{B}^{-}_0)\big|.
\end{equation}
}\end{proposition}
\begin{proof}
\cite[Excercise 8.29]{cb}.
\end{proof}
\begin{proposition}[{\cite[Excercise 8.30]{cb}}]\label{proposition:circleserieslength}There is a positive constant $c_r$ (depending on the length of $\Gamma_r$, hence on $r$, $\phistar$) such that in any sequence of $k+1\ge2$ successive collisions with $\Gamma_r$ the path length $k\tau$ between the very first and very last collision always exceeds $c_r$.
\end{proposition}
\begin{corollary}\label{corollary:shearlowerboundexpansion}
If $(\phi,\theta)=x\in {\Min_{r,k}}$, $k\ge2$, $(d\phi,d\theta)=dx \in T_xM$ such that $\frac{d\theta}{d\phi}\in[c,+\infty]$ with some $c>0$, and $\Fc^{k}(x)\in \Mrout$, $\Fc^{k-1}(x)\in \Fc^{-1}(\Mrout\smallsetminus\Mrin)$, then $\frac{\|D_{x}\Fc^{k-1}(dx)\|_\p}{\|dx\|_\p}\ge1+\frac{c_rc}{2}$ with $c_r$ from \cref{proposition:circleserieslength}.
\end{corollary}
\begin{proof}
By \cref{proposition:ImaginaryExpansion}, the mirror equation \eqref{eq:MirrorEquationExpansionFormula}, and $\tau=2r\cos{\varphi}=2r\sin{\theta}$,
\begin{align*}
\frac{\|D_x\Fc^{k-1}(dx)\|_\p}{\|dx\|_\p}&\overset{\eqref{eqJfromCM}}=|1-(k-1)\tau(-2/\tau+B^{-})|
\overset{\eqref{eq:MirrorEquationExpansionFormula}}=\Big|1-(k-1)\tau\Big(-\frac{2}{\tau}+\frac{\frac{-d\theta}{rd\phi}+\frac{1}{r}}{\cos{\varphi}}\Big)\Big|\\
\overset{\mathclap{\tau=2r\cos{\varphi}}}{\scalebox{6}[1]{$=$}}&\Big|1-(k-1)\tau\Big(-\frac{2}{2r\cos{\varphi}}+\frac{\frac{-d\theta}{rd\phi}+\frac{1}{r}}{\cos{\varphi}}\Big)\Big|=\Big|1+(k-1)\tau\frac{{d\theta}/{d\phi}}{r}\Big|.
\end{align*}
The reasoning in \cref{proposition:circleserieslength} also shows that there exists a constant $\overline{c}_r=\frac{1}{2}c_r$ such that for all $k\ge2$, the path length $(k-1)\tau$ between the first and next to last collisions always exceeds $\overline{c}_r$. Thus the conclusion follows.
\end{proof}
\subsection[A formula \texorpdfstring{(\cref{def:AformularForLowerBoundOfExpansion})}{} for expansion lower bound in Case (a)]{A formula \texorpdfstring{(\cref{def:AformularForLowerBoundOfExpansion})}{} for expansion lower bound in Case (a)}\label{subsec:AFELB}\strut\COMMENT{\small\red Say what the formula is. Not necessarily right here, but within at most five lines from here.\wentao{Having put \cref{def:AformularForLowerBoundOfExpansion} to subsection title.}}\hfill

From here on, we work towards the proofs of \cref{TEcase_a0,TEcase_a1}, and accordingly, subsequent statements are made in the context of these propositions.
\begin{notation}\label{def:AformularForLowerBoundOfExpansion}
For $x_0\in\Mrout$, $x_1\in\MRin$, $x_2\in\Mrin$, $n_1$ from \cref{def:MhatReturnOrbitSegment} and 
$\tau_0$, $\tau_1$, $d_0$, $d_1$, $d_2$ from \cref{def:fpclf}, considered as functions of $x_1\in\big\{(\Phi_1,\theta_1)\in\MRin\bigm|\sin\theta_1\ge2r/R\big\}$ (see \cref{fig:R1,fig:A_compactRegion,fig:heuristicproof}), let
\[\mathcal{I}\dfn-1+\frac{\tau_1}{d_0}\Big[\frac{2(\tau_0-d_0)}{d_1}-\frac{\tau_0+\tau_1-2d_0}{\tau_1}\Big].\] 
Note that if $d_1\ge2r$, then $n_1=0$, and if $n_1=0$, then $\tau_0+\tau_1>2d_i$, $i=0,2$. See \cref{fig:lengthfunctonFigure,fig:heuristicproof,remark:symmetrywiththetangentline}.
\end{notation}
\begin{lemma}[Orbit configuration and expansion estimate in Case (a) of \eqref{eqMainCases}]\label{lemma:ExpansionLowerBoundA0A1}
Consider case (a) of \eqref{eqMainCases}, i.e., suppose a nonsingular $x\in \hat{M}$  (\cref{def:SectionSetsDefs}) has a return orbit segment as in \cref{def:MhatReturnOrbitSegment}, which is $x,\Fc(x),\cdots,\Fc^{\sigma(x)}(x)=\hat{F}(x)\in \hat{M}$ with $x_{0}\in\Mrout$, $x_1=\Fc(x_0)=(\phi_1,\theta_1)\in \MRin$, $x_2\in\Mrin$ with $\sin\theta_1\ge 2r/R$. We now consider a  tangent vector at $x$ in the first or third quadrant in $\phi\theta$-coordinates,\COMMENT{Change all occurrences of "$\phi$, $\theta$-coordinates" to "$\phi\theta$-coordinates"!\wentao{done}}
\[dx=(d\phi,d\theta)\in\Big\{(d\phi,d\theta)\bigm|\frac{d\theta}{d\phi}\in\big[0,+\infty\big]\Big\}\]
The corresponding tangent vectors at $x_0,x_1,x_2$ are then defined as follows.
If $x\in\Mrin\cap\Mrout$, then $dx_0\dfn dx$ and $dx_2\dfn D_x\Fc^2(dx)$. Otherwise, $x\in\Fc^{-1}(\Mrout\setminus\Mrin)$, so $dx_0\dfn D_x\Fc(dx)$ and $dx_2\dfn D_{x}\Fc^3(dx)$. In either case, $dx_1\dfn D_{x_0}\Fc(dx_0)$. 

We then have the following.
\begin{itemize} 
    \item $\frac{\|dx_0\|_\p}{\|dx\|_\p}\ge 1$.\item $\frac{\|D_x\hat{F}(dx)\|_\p}{\|dx_2\|_\p}\ge 1$.
    \item If $\mathcal{I}<0$ (\cref{def:AformularForLowerBoundOfExpansion}), then $\frac{\|dx_2\|_\p}{\|dx_0\|_\p}\ge|\mathcal{I}|$.
\end{itemize}
\end{lemma}
\begin{proof}
We first note that the $dx_i$ for $i=0,2$ are in the first or third quadrant. For $i=0$ this is \cref{proposition:ShearOneStepExpansion}. For $i=2$, we invoke \cref{caseA}: $x_1=\Fc(x_0)$, $x_2=\Fc^2(x_0)$ and $D_{x_0}\Fc^2$ is a negative matrix in $\phi\theta$-coordinates. Hence $(d\phi_2,d\theta_2)=dx_2=D_{x_0}\Fc^2(dx_0)$ has $\frac{d\theta_2}{d\phi_2}>0$, that is $dx_2$ is also in the quadrant. 

If $x\in\Mrin\cap\Mrout$, then $x_0=x,dx_0\dfn dx$, hence  $\frac{\|dx_0\|_\p}{\|dx\|_\p}=1\ge 1$. Otherwise $x\in\Fc^{-1}(\Mrout\setminus\Mrin)$, so $x_0=\Fc(x)$, $D_{x}\Fc(dx)=dx_0=(d\phi_0,d\theta_0)$ and 
since $dx$ is in the first or third quadrant, \cref{proposition:ShearOneStepExpansion} gives $\frac{\|dx_0\|_\p}{\|dx\|_\p}\ge 1$---and likewise the second item.



To analyze $\frac{\|dx_2\|_\p}{\|dx_0\|_\p}$ we use the same notation as in \eqref{eq:MirrorEquationExpansionFormula} with subscripts. We denote by $\mathcal{B}^{\pm}_0,\mathcal{B}^{\pm}_1$ the infinitesimal wave front curvatures after/before collisions at $x_0$, $x_1$, respectively.
Then
\begin{equation}\label{eqB+d0}
\mathcal{B}^{+}_0\overset{\eqref{eq:MirrorEquationExpansionFormula}}=\frac{\frac{-d\theta_0}{rd\phi_0}-\frac{1}{r}}{\cos{\varphi_0}}\le-\frac{1}{r\cos{\varphi_0}}\overset{\eqref{eq:coordinateschange}}=-\frac{1}{r\sin{\theta_0}}\overset{\text{\tiny\cref{def:fpclf}}}{\scalebox{4}[1]{$=$}}-\frac1{d_0}<0.
\end{equation}
\begin{remark}

We first note that \eqref{eq:MirrorEquationExpansionFormula} gives the following. \begin{equation}\label{eqBBd}\mathcal{B}^{+}_1\overbracket{=}^{\text{\eqref{eq:MirrorEquationExpansionFormula}}}\mathcal{B}^{-}_1-\frac{2}{d_1}.\end{equation}
\cite[equation (3.31)]{cb} gives\begin{equation}\label{eqBBr}
\begin{aligned}
\frac{1}{\mathcal{B}_{1}^{-}}\overbracket{=}^{\mathclap{\substack{\text{\cite[equation (3.31)]{cb}}\\\mid}}}&\frac{1}{\mathcal{B}_{0}^{+}}+\tau_0
\end{aligned}
\end{equation}
\end{remark}
\eqref{eq:MirrorEquationExpansionFormula} gives the following:.\(\frac{\|dx_1\|_\p}{\|dx_0\|_\p}=|1+\tau_0\mathcal{B}^{+}_0|,\) then in turn, by re-indexing using $0$, $1$ and $2$ corresponding to $x_0$, $x_1$, $x_2$ from \cref{def:fpclf} respectively, we get the following.
\[\frac{\|dx_2\|_\p}{\|dx_1\|_\p}=|1+\tau_1\mathcal{B}^{+}_1|\overbracket{=}^{\mathclap{\mathcal{B}^{+}_1\overset{\text{\eqref{eq:MirrorEquationExpansionFormula}}}=\mathcal{B}^{-}_1-\frac{2}{d_1}}}\big|1+\tau_1(\mathcal{B}^{-}_1-\frac2{d_1})\big|\quad \overbracket{=}^{\mathclap{\text{\eqref{eqBBr}: }\frac{1}{\mathcal{B}_{1}^{-}}=\frac{1}{\mathcal{B}_{0}^{+}}+\tau_0}}\Big|1+\tau_1\big(\frac{1}{\frac{1}{\mathcal{B}^{+}_0}+\tau_0}-\frac{2}{d_1}\big)\Big|=\Big|1+\frac{\mathcal{B}^{+}_0\tau_1}{1+\mathcal{B}^{+}_0\tau_0}-\frac{2\tau_1}{d_1}\Big|.
\]
Therefore
\begin{equation}\label{eq:dx2overdx0}
\frac{\|dx_2\|_\p}{\|dx_0\|_\p}=\frac{\|dx_2\|_\p}{\|dx_1\|_\p}\frac{\|dx_1\|_\p}{\|dx_0\|_\p} =|1+\tau_0\mathcal{B}^{+}_0|\Big|1+\frac{\mathcal{B}^{+}_0\tau_1}{1+\mathcal{B}^{+}_0\tau_0}-\frac{2\tau_1}{d_1}\Big| =\Big|\underbracket{(\tau_0+\tau_1-\frac{2\tau_0\tau_1}{d_1})\mathcal{B}^{+}_0+1-\frac{2\tau_1}{d_1}}_{\nfd\mathcal{I}(\mathcal{B}^{+}_0)}\Big|
\end{equation}
In Case (a) of \eqref{eqMainCases} 
we have $d_1\ge2r$, so $\frac{2}{d_1}\le\frac{1}{r}\le\frac{1}{2r}+\frac{1}{2r}\le\frac{1}{\tau_0}+\frac{1}{\tau_1}$, hence $\tau_0+\tau_1-\frac{2\tau_0\tau_1}{d_1}\ge0$. Thus, 
\begin{equation}\label{eq:56}
\begin{split}
\mathcal{I}(\mathcal{B}^{+}_0) & \overbracket\le^{\eqref{eqB+d0}}(\tau_0+\tau_1-\frac{2\tau_0\tau_1}{d_1})(-\frac{1}{d_0})+1-\frac{2\tau_1}{d_1}\\
& =1-\frac{\tau_0+\tau_1}{d_0}+\frac{2\tau_0\tau_1}{d_0d_1}-\frac{2\tau_1}{d_1} = -1-\frac{\tau_0+\tau_1-2d_0}{d_0}+\frac{2\tau_1}{d_1}(\frac{\tau_0}{d_0}-1)\\
& =-1-\frac{\tau_0+\tau_1-2d_0}{d_0}+\frac{2\tau_1}{d_1}(\frac{\tau_0-d_0}{d_0}) = -1+\frac{\tau_1}{d_0}\Big[\frac{2(\tau_0-d_0)}{d_1}-\frac{\tau_0+\tau_1-2d_0}{\tau_1}\Big]=\mathcal{I}
\end{split}
\end{equation}
Hence \eqref{eq:56} combined with \eqref{eq:dx2overdx0} gives that if $\mathcal{I}<0$, then $\frac{\|dx_2\|_\p}{\|dx_0\|_\p}\ge|\mathcal{I}|$.
\end{proof}

\begin{proposition}\label{Prop:Ilowerbound}
For an orbit segment of $x\in\hat{M}$ from \cref{def:MhatReturnOrbitSegment} with $x_0\in\Mrout$, $\Fc(x_0)=x_1=(\Phi_1,\theta_1)\in\MRin,\,\sin{\theta_1}\ge\sqrt{4r/R}$, $x_2\in\Mrin$ in 
case (a0) in \eqref{eqMainCases}, we have $\mathcal{I}(x)<-1$ (\Cref{def:AformularForLowerBoundOfExpansion}), and hence $\frac{\|dx_2\|_\p}{\|dx_0\|_\p}>1$ by \cref{lemma:ExpansionLowerBoundA0A1}.
\end{proposition}
\begin{proof}
We invoke  \cref{prop:casea0largetheta,prop:casea0smallthetaR2,prop:casea0smallthetaR1} as follows (see also \cref{fig:R1_fraction_extreme}).

If $\sin\theta_1>1/2$, then $\mathcal{I}<-1.37$ (\cref{prop:casea0largetheta} below). If $\sin\theta_1\ge\sqrt{4r/R}$ and $x_1\in R_2$ in \cref{fig:R1}, then $\mathcal{I}<-1$ (\cref{prop:casea0smallthetaR2} below). If $1/2\ge\sin\theta_1\ge\sqrt{4r/R}$ and $x_1\in R_1$ in \cref{fig:R1}, then $\mathcal{I}<-1$ (\cref{prop:casea0smallthetaR1} below). Then, by symmetry on the entire region, $\mathcal{I}<-1$ for $x_1\in\big\{(\Phi_1,\theta_1)\in\Mrin\bigm|\sin\theta_1\ge\sqrt{4r/R}\big\}$. 
\end{proof}
By supplying the invoked facts, \cref{SBSProp53PartI,SBSProp53PartII} prove \cref{Prop:Ilowerbound}, dividing the proof according to large or small incidence angles.
\subsection{Proof of \texorpdfstring{\cref{Prop:Ilowerbound}}{Lower bound for I}, large angles}\label{SBSProp53PartI} We now consider the case when the returning orbit segment of $x$\COMMENT{\large I recommend to replace "$x$ return orbit segment" by "returning orbit segment of $x$" \bfseries\Large\color{red}throughout.} has a more transverse reflection on $\Gamma_R$, specifically, with $x_1=(\Phi_1,\theta_1)$, $\frac{\pi}{6}<\theta_1<\frac{5\pi}{6}$.\COMMENT{"$x_1=(\Phi_1,\theta_1), \frac{\pi}{6}<\theta_1<\frac{5\pi}{6}$" should be "$x_1=(\Phi_1,\theta_1)$, $\frac{\pi}{6}<\theta_1<\frac{5\pi}{6}$"---look at the tex file to see the difference; \Large\bfseries\color{red}this is a lesson that should be applied throughout.} The goal is to prove that $\mathcal{I}<-1.37$ (\cref{prop:casea0largetheta}). Less transverse reflections will be considered starting on page \pageref{SBSProp53PartII}.
\begin{lemma}\label{lemma:returnorbittheta1Large}
Assuming \eqref{eqJZHypCond} and $\phistar\in\big(0,\tan^{-1}(1/3)\big)$ consider $x\in \hat{M}$ and its returning orbit segment (\cref{def:MhatReturnOrbitSegment,def:fpclf,eqPtsFromeqMhatOrbSeg}) $x_0\in\Mrout$, $x_1=(\Phi_1,\theta_1)\in\MRin$, $x_2=(\phi_2,\theta_2)\in\Mrin$. If either $\frac{\pi}{2}<\phi_0<\frac{3\pi}{2}$ or $\frac{\pi}{2}<\phi_2<\frac{3\pi}{2}$, then $\theta_1\in \big(\frac{\pi}{6},\frac{5\pi}{6}\big)$, $\sin\theta_1>\sqrt{4r/R}$, $n_1=0$ in \eqref{eqPtsFromeqMhatOrbSeg}, hence $x_2=\Fc(x_1)$.
 

\end{lemma}
\begin{proof}
We will prove the first conclusion since the second conclusion proof will be just a symmetric argument of the first conclusion.

By symmetry we can assume $0<\theta_1\le\frac{\pi}{2}$. In \cref{fig:QonUpperHalfCircle}
, the corners of the lemon billiard have Cartesian coordinates $A=(-r\sin\phistar,-r\cos\phistar)$, $B=(r\sin\phistar,-r\cos\phistar)$. Let $C=(-r,0)$, $D=(r,0)$, both on $C_r$. 

For the return orbit segment 
of $x\in \hat{M}$ (with $x_0\in\Mrout$, $x_1=(\Phi_1,\theta_1)\in\MRin$, $x_2=(\phi_2,\theta_2)\in\Mrin$), write $P=p(\Fc^{-1}(x_2))$, $Q=p(x_2)$ (positions). Then $\Fc^{-1}(x_2)=(\Phi,\theta_1)\in\MRout$ for some $\Phi\in(-\Phistar,\Phistar)$.

On the other hand $x_2=(\phi_2,\theta_2)\in\Mrin,\phi_2\in(\frac{\pi}{2},\frac{3\pi}{2})$ implies $Q$ is in the upper half of $\Gamma_r$, i.e., $Q$ has a positive $y$-coordinate (\cref{fig:QonUpperHalfCircle}). 

Since $P$ is an interior point of the arc $\arc{AB}$, $Q$ is an interior point of the upper half of $C_r$ (\cref{fig:QonUpperHalfCircle}). The line segment $\overline{PQ}$ intersects the line segment $\overline{AB}$ at some $P_1$ in the interior of the line segment $\overline{AB}$, and $\overline{PQ}$ intersects $\overline{CD}$ at some $P_2$ in th interior of $\overline{CD}$. Note that the line $\overline{CD}$ goes through $O_r=(0,0)$ and is parallel to $\overline{AB}$. By some elementary geometry, $\measuredangle{BAD}=\frac{\pi}{4}-\frac{\phistar}{2}$ and $\theta_1+\Phi=\measuredangle{BP_1Q}$. 

Write $Q=(x_Q,y_Q),P=(x_P,y_P)$, $P_1=(x_{P_1},y_{P_1})$, $P_2=(x_{P_2},y_{P_2})$, $A=(x_A,y_A)$, $B=(x_B,y_B)$, $D=(x_D,y_D)$. Note that $y_Q-y_P>y_{P_2}-y_{P_1}=y_D-y_A>0$ and if $x_Q>x_P$ then $0<x_Q-x_P<x_D-x_A$. Hence if $0\ge x_Q-x_P$ then the line $PQ$ has a negative slope or infinite slope. Otherwise if $x_Q>x_P$ then the line $PQ$ has a positive slope larger than that of line $AD$ (\cref{fig:QonUpperHalfCircle}). In either case, 
$\theta_1+\Phi =\measuredangle{BP_1Q}>\measuredangle{BAD}=\frac{\pi}{4}-\frac{\phistar}{2}$, hence\COMMENT{\red See \url{https://www.overleaf.com/learn/latex/Aligning_equations_with_amsmath}}
\begin{equation}\label{eq:57}
\theta_1>\frac{\pi}{4}-\frac{\phistar}{2}-\Phi\overbracket{>}^{\mathclap{\strut\Phi<\Phistar}}\frac{\pi}{4}-\frac{\phistar}{2}-\Phistar=\frac{\pi}{4}-\frac{\phistar}{2}-\sin^{-1}{(\frac{r\sin\phistar}{R})}\overbracket{>}^{\mathclap{\bigstrut\phistar\in\big(0,\tan^{-1}(1/3)\big), R>1700r}}\frac{\pi}{6}.
\end{equation}

This was the case $\theta_1\in(0,\frac{\pi}{2}]$.
For $\theta_1\in[\frac{\pi}{2},\pi)$, by symmetry we get $\pi-\theta_1-\Phi=\measuredangle{AP_1Q}>\measuredangle{ABC}=\frac{\pi}{4}-\frac{\phistar}{2}$ and  $\theta_1<\frac{5\pi}{6}$ (similarly to \eqref{eq:57}).

Thus, $\sin\theta_1>0.5\overbracket{>}^{\mathclap{R>1700r}}\sqrt{4r/R}>\frac{r\sin\phistar}{R}$, so $n_1=0$  in \eqref{eqPtsFromeqMhatOrbSeg}
, that is, $x_1=\Fc^{-1}(x_2)$.
\end{proof}

\begin{figure}[h]
\begin{minipage}{.5\textwidth}
\begin{center}
\begin{tikzpicture}[xscale=0.6, yscale=0.6]
    \tkzDefPoint(0,0){Or};
    \tkzDefPoint(0,-4.8){Y};
    \pgfmathsetmacro{\rradius}{4.5};
    \pgfmathsetmacro{\bdist}{4.7};
    \tkzDefPoint(0,\bdist){OR};
    \pgfmathsetmacro{\phistardeg}{40};
    \pgfmathsetmacro{\XValueArc}{\rradius*sin(\phistardeg)};
    \pgfmathsetmacro{\YValueArc}{\rradius*cos(\phistardeg)};
    \tkzDefPoint(\XValueArc,-1.0*\YValueArc){B};
    \tkzDefPoint(-1.0*\XValueArc,-1.0*\YValueArc){A};
    \pgfmathsetmacro{\Rradius}{veclen(\XValueArc,\YValueArc+\bdist)};
    \tkzDrawArc[name path=Cr,very thick](Or,B)(A);
    \tkzDrawArc[name path=CR,very thick](OR,A)(B);
    \pgfmathsetmacro{\PXValue}{\Rradius*sin(10)};
    \pgfmathsetmacro{\PYValue}{\bdist-\Rradius*cos(10)};
    \tkzDefPoint(\PXValue,\PYValue){P};
    \pgfmathsetmacro{\QXValue}{\rradius*cos(10)};
    \pgfmathsetmacro{\QYValue}{\rradius*sin(10)};
    \tkzDefPoint(\PXValue,\PYValue){P};
    \tkzDefPoint(\QXValue,\QYValue){Q};
    \tkzDefPoint(-\rradius,0){C};
    \tkzDefPoint(\rradius,0){D};
    \path[name path=linePQ](P)--(Q);
    \begin{scope}[ultra thin,decoration={
    markings,
    mark=at position 0.55 with {\arrow{>}}}
    ] 
    \draw[blue,postaction={decorate}](P)--(Q);
    \end{scope}
    \draw[name path=lineCD,red,dashed,ultra thin](C)--(D);
    \draw[name path=lineAB,red,dashed,ultra thin](A)--(B);
    \draw[name path=lineORP,red,dashed,ultra thin](OR)--(P);
    \draw[name path=lineOrA,red,dashed,ultra thin](Or)--(A);
    \draw[name path=lineAB,red,dashed,ultra thin](B)--(A);
    \draw[name path=lineOROr,red,dashed,ultra thin](OR)--(Or);
    \draw[name path=lineAD,red,dashed,ultra thin](A)--(D);
    \path[name intersections={of=lineAB and linePQ,by={P1}}];
    \path[name intersections={of=lineCD and linePQ,by={P2}}];
    \tkzLabelPoint[above](OR){$O_R$};
    \tkzLabelPoint[above right](Q){$Q$};
    \tkzLabelPoint[below](P){$P$};
    \tkzLabelPoint[above](Or){$O_r$};
    \tkzLabelPoint[left](A){$A$};
    \tkzLabelPoint[right](B){$B$};
    \tkzLabelPoint[left](C){$C$};
    \tkzLabelPoint[right](D){$D$};
    \tkzLabelPoint[above left](P2){$P_2$};
    \tkzLabelPoint[above, yshift=1.1](P1){$P_1$};
    \tkzDrawPoints(A,B,C,D,Or,OR,P,Q,P1,P2);
    \tkzMarkAngle[arc=lll,ultra thin,size=0.6](C,Or,A);
    \tkzMarkAngle[arc=lll,ultra thin,size=0.4](B,P1,Q);
    \tkzMarkAngle[arc=lll,ultra thin,size=0.8](Or,OR,P);
    \tkzMarkAngle[arc=lll,ultra thin,size=0.8](B,A,D);
    \node at ($(OR) + (-83:1.6)$)  { $\Phi_1=\Phi$};
    \node at ($(Or) + (-160:1.4)$)  {\Large $\frac{\pi}{2}-\phistar$};
    \node at ($(A) + (10:1.8)$)  {$\frac{\pi}{4}-\frac{\phistar}{2}$};
    \node at ($(P1) + (18:1.4)$)  {\small $\theta_1+\Phi$};
\end{tikzpicture}
\end{center}
\caption{$p(x_2)$ is on the upper half circle. $O_r=(0,0)$, $O_R=(0,b)$,
\newline $A=(-r\sin\phistar,-r\cos\phistar)$,
\newline $B=(r\sin\phistar,-r\cos\phistar)$.
}\label{fig:QonUpperHalfCircle}
\end{minipage}\hfill
\begin{minipage}{.5\textwidth}
\begin{center}
\begin{tikzpicture}[xscale=0.6, yscale=0.6]
    \tkzDefPoint(0,0){Or};
    \tkzDefPoint(0,-4.8){Y};
    \pgfmathsetmacro{\rradius}{4.5};
    \pgfmathsetmacro{\bdist}{8.7};
    \tkzDefPoint(0,\bdist){OR};
    \pgfmathsetmacro{\phistardeg}{40};
    \pgfmathsetmacro{\XValueArc}{\rradius*sin(\phistardeg)};
    \pgfmathsetmacro{\YValueArc}{\rradius*cos(\phistardeg)};
    \tkzDefPoint(\XValueArc,-1.0*\YValueArc){B};
    \tkzDefPoint(-1.0*\XValueArc,-1.0*\YValueArc){A};
    \pgfmathsetmacro{\Rradius}{veclen(\XValueArc,\YValueArc+\bdist)};
    \clip(-\rradius-1.5,-\rradius-0.6) rectangle (\rradius+1.5,\bdist+0.4);
    \draw[name path=Cr,dashed,ultra thin] (Or) circle (\rradius);
    \tkzDrawArc[very thick](Or,B)(A);
    \tkzDrawArc[name path=CR,very thick](OR,A)(B);
    \pgfmathsetmacro{\PXValue}{\Rradius*sin(7)};
    \pgfmathsetmacro{\PYValue}{\bdist-\Rradius*cos(7)};
    \tkzDefPoint(\PXValue,\PYValue){P};
    \tkzDefPoint(0,\bdist-\Rradius){C};
    \draw[name path=PQ_Far,blue,ultra thin,dashed](P)--+(70:7)coordinate(PposX);
    \draw[name path=PT_Far,blue,ultra thin,dashed](P)--+(124:9);
    \draw[name path=PT0_Far,blue,ultra thin,dashed](P)--+(-56:0.7);
    \draw[name path=PQ0_Far,blue,ultra thin,dashed](P)--+(-110:1);
    \draw[red,dashed,->](P)--+(7:4)coordinate(PposX);
    \draw[red,dashed](P)--+(-173:6)coordinate(PnegX);
    \draw[name path=linePOR,red,dashed,->](P)--(OR);
    \draw[name path=linePOr,red,dashed](P)--(Or);
    \tkzDefPoint(0,-\rradius){Y};
    \draw[red,dashed](OR)--(Y);
    \tkzMarkAngle[arc=lll,ultra thin,size=1.7](Or,OR,P);
    \tkzMarkAngle[arc=lll,ultra thin,size=0.7](Y,Or,P);
    \tkzMarkAngle[arc=lll,ultra thin,size=0.8](OR,P,Or);
    
    \path[name intersections={of=Cr and PQ_Far,by={Q}}];
    \path[name intersections={of=Cr and PQ0_Far,by={Q0}}];
    \path[name intersections={of=Cr and PT0_Far,by={T0}}];
    \path[name intersections={of=Cr and PT_Far,by={T}}];
    \tkzMarkAngle[arc=lll,ultra thin,size=1.8](PposX,P,Q);
    \tkzMarkRightAngle(OR,P,PnegX);
    \node[] at ($(Or)+(-78:1)$)  {\Large $\alpha$};
    \node[] at ($(OR)+(-85:2.3)$)  {\Large $\Phi_1$};
    \node[] at ($(P)+(105:1.3)$)  {\Large {$\delta$}};
    \node[] at ($(P)+(35:2.3)$)  {\Large {$\theta_1$}};
    \tkzLabelPoint[left](A){$A$};
    \tkzLabelPoint[below](B){$B$};
    \tkzLabelPoint[above](C){$M$};
    \tkzLabelPoint[left](OR){$O_R$};
    \tkzLabelPoint[left](Or){$O_r$};
    \tkzLabelPoint[below left](P){$P$};
    \tkzLabelPoint[right](Q){$Q$};
    \tkzLabelPoint[left](T){$T$};
    \tkzLabelPoint[below](Q0){$Q_0$};
    \tkzLabelPoint[below](T0){$T_0$};
    \tkzLabelPoint[below](Y){$Y$};
    \tkzLabelPoint[below](PposX){$TP$};
    \tkzDrawPoints(A,Or,OR,P,B,Y,T,Q,T0,Q0,C);
    \begin{scope}[ultra thin,decoration={markings, mark=at position 0.55 with {\arrow{>}}}
    ] 
    \draw[blue,postaction={decorate}](P)--(Q);
    \draw[blue,postaction={decorate}](T)--(P);
    \end{scope}
\end{tikzpicture}
\caption{Normal coordinates at $P$}\label{fig:fractionEstimate2}
\end{center}
\end{minipage}
\end{figure}

Unlike \cref{def:StandardCoordinateTable}, we now define a new coordinate system for the billiard table with origin at the collision position $P\in\Gamma_R$.
\begin{remark}[Normal coordinates at $P$]\label{def:TPPORsystem}
With \cref{def:LemonTableConfiguration}, suppose that $M$ is the midpoint of $\Gamma_R$. Let $x_1=(\Phi_1,\theta_1)\in\MRin$ with $\theta_1\in[\sin^{-1}(2r/R),\pi-\sin^{-1}(2r/R)]$ in case (a) of \eqref{eqMainCases} and $P=p(x_1)\in\Gamma_R$ of the billiard boundary. As in \cref{fig:fractionEstimate2}, we will at times use coordinates centered at $P$ with axes tangent and normal to $C_R$, respectively. The length and angle variables 
    \begin{itemize}
        \item $\rho\dfn|O_rP|\overbracket=^{\mathclap{\text{cosine law}}}\sqrt{R^2+b^2-2bR\cos\Phi_1}$,
        \item $\alpha\dfn\tan^{-1}{\big(\frac{R\sin\Phi_1}{R\cos\Phi_1-b}\big)}$,
        \item $\delta\dfn\alpha-\Phi_1$
    \end{itemize}
in \cref{fig:fractionEstimate2} 
are smooth functions of $x_1$. Note that 
$P=(R\sin\Phi_1,b-R\cos\Phi_1)$, $O_r=(0,0)$, $O_R=(0,b)$, therefore, in \cref{fig:fractionEstimate2}, $\alpha$ is the angle between $\overrightarrow{O_rM}$ and $\overrightarrow{O_rP}$ given the counterclockwise orientation on the billiard table. Then $\alpha$, $\delta$ and $\Phi_1$ have the same sign, and $|\alpha|=\measuredangle{MO_rP}\ge\measuredangle{MO_RP}=|\Phi_1|$,  $|\delta|=\measuredangle{O_rPO_R}$.
\end{remark}



\begin{lemma}\label{lemma:newcoordinateslength}
    For the length functions $\tau_0$, $\tau_1$, $d_0$, $d_1$, $d_2$ from \cref{def:fpclf,def:AformularForLowerBoundOfExpansion} and $\rho$, $\alpha$, $\delta$ from \cref{def:TPPORsystem} as functions of $x_1=(\Phi_1,\theta_1)\in\Mrin$ with $\theta_1\in[\sin^{-1}(2r/R),\pi-\sin^{-1}(2r/R)]$, we have    \begin{equation}\label{eq:NewCoordinateLegnth}
        \begin{aligned}
        \tau_1&=\rho\sin{(\theta_1-\delta)}+\sqrt{r^2-\rho^2\cos^2{(\theta_1-\delta)}},\\
        d_2&=\sqrt{r^2-\rho^2\cos^2{(\theta_1-\delta)}},\\
        \tau_0&=\rho\sin{(\theta_1+\delta)}+\sqrt{r^2-\rho^2\cos^2{(\theta_1+\delta)}},\\
        d_0&=\sqrt{r^2-\rho^2\cos^2{(\theta_1+\delta)}},
        \end{aligned}
    \end{equation}
so 
$\tau_0,\tau_1,d_0,d_1,d_2$ are continuous functions of $x_1$ in $\mathfrak{D}^\circ\dfn R_1\cup R_2\cup R_3\cup R_4$ in \cref{fig:R1}.
\end{lemma}
\begin{proof}
    We first assume $\theta_1\in[\sin^{-1}(2r/R),\pi/2]$ and use normal coordinates at $P$ (\cref{def:TPPORsystem}). We have $P=(0,0)$, $O_r=(x_r,y_r)=(-\rho\sin\delta,\rho\cos\delta)$ and suppose $Q=p(x_2)=(x_Q,y_Q)$ with $x_Q>0$, $y_Q>0$. Then $Q_0=(x_{Q_0},y_{Q_0})$ with $x_{Q_0}<0$, $y_{Q_0}<0$ is another intersection point of the line $PQ$ with $C_r$.
    
$(x,y)=(x_Q,y_Q)$ and $(x,y)=(x_{Q_0},y_{Q_0})$ both satisfy
    \begin{align*}
         &\left.
        \begin{aligned}
            &(x-x_r)^2+(y-y_r)^2=r^2,\\
            &y=x\tan{\theta_1},\\
        \end{aligned}
        \right\} \text{ hence }  (\sec^2\theta_1)x^2-2(x_r+y_r\tan\theta_1)x+\underbracket{x^2_r+y^2_r}_{=\rho^2}-r^2=0,
    \end{align*}
and 
\(2x_r+2y_r\tan\theta_1=2(-\rho\sin\delta+\rho\cos\delta\tan\theta_1)=2\rho\sec\theta_1\sin{(\theta_1-\delta)}\). Thus, $x=x_Q$ and $x=x_{Q_0}$ both satisfy
    \begin{equation}\label{eq:60}
            x^2\sec^2\theta_1-2\rho\sin{(\theta_1-\delta)}x\sec\theta_1+\rho^2-r^2=0.
    \end{equation}
Since $\rho^2-r^2<0$, this quadratic equation $L^2-2\rho\sin{(\theta_1-\delta)}L+\rho^2-r^2=0$ for $L=x\sec\theta_1$ has one positive and one negative root. Since $\sec\theta_1>0$, $x_{Q_0}<0<x_Q$, they are \[\tau_1=|PQ|=x_Q\sec\theta_1=L_+=\rho\sin{(\theta_1-\delta)}+\sqrt{r^2-\rho^2\cos^2{(\theta_1-\delta)}}>0\] and $x_{Q_0}\sec\theta_1=L_-=\rho\sin{(\theta_1-\delta)}-\sqrt{r^2-\rho^2\cos^2{(\theta_1-\delta)}}<0$. Since 
$x_Q\sec\theta_1-x_{Q_0}\sec\theta_1=|PQ|+|PQ_0|=|QQ_0|=2d_2$ (see \cref{def:fpclf,fig:fractionEstimate2}), 
$d_2=\sqrt{r^2-\rho^2\cos^2{(\theta_1-\delta)}}$.

Similarly, in normal coordinates at $P$, if $(x_T,y_T)=T=p(x_0)=p(\Fc^{-1}(x_1))$ and $T_0$ is the other intersection of $PT$ with $C_r$, then $x_T<0<y_T$, $y_{T_0}<0<x_{T_0}$, and $(x,y)=(x_T,y_T)$,  $(x,y)=(x_{T_0},y_{T_0})$ both satisfy
    \begin{align*}
         &\left.
        \begin{aligned}
            &(x-x_r)^2+(y-y_r)^2=r^2,\\
            &y=-(\tan{\theta_1})x,\\
        \end{aligned}
        \right\} \text{ hence } (\sec^2\theta_1)x^2-2(x_r-y_r\tan\theta_1)x+\rho^2-r^2=0,    \end{align*}
and 
$2x_r-2y_r\tan\theta_1=2(-\rho\sin\delta-\rho\cos\delta\tan\theta_1)=-2\rho\sec\theta_1\sin{(\theta_1+\delta)}$. 
Thus, $x=x_T$ and $x=x_{T_0}$ both satisfy
    \begin{equation}\label{eq:61}
            (\sec^2\theta_1)x^2+2\rho\sin{(\theta_1+\delta)}(\sec\theta_1)x+\rho^2-r^2=0.
    \end{equation}
As before, this quadratic equation $L^2-2\rho\sin{(\theta_1+\delta)}L+\rho^2-r^2=0$ for $L=-x\sec\theta_1$ has roots of opposite sign:
\[\tau_0=|PQ|=-x_T\sec\theta_1=L_+=\rho\sin{(\theta_1+\delta)}+\sqrt{r^2-\rho^2\cos^2{(\theta_1+\delta)}}>0\] and  $-x_{T_0}\sec\theta_1=L_-=\rho\sin{(\theta_1+\delta)}-\sqrt{r^2-\rho^2\cos^2{(\theta_1+\delta)}}$. 
Then $-x_T\sec\theta_1+x_{T_0}\sec\theta_1=|PT|+|PT_0|=|TT_0|=2d_0$ (see \cref{def:fpclf,fig:fractionEstimate2}) gives 
$d_0=\sqrt{r^2-\rho^2\cos^2{(\theta_1+\delta)}}$.
    
This proves \eqref{eq:NewCoordinateLegnth} for $\theta_1\in[\sin^{-1}(2r/R),\pi/2]$.    
For $\theta_1\in[\pi/2,\pi-\sin^{-1}(2r/R)]$ we use symmetry. The value of $\tau_1$ for $\theta_1\in[\pi-\sin^{-1}(2r/R),\pi/2]$ equals the value of $\tau_0$ for $(\pi-\theta_1)\in[\sin^{-1}(2r/R),\pi/2]]$ and the value of $\tau_0$ for $\theta_1\in[\pi-\sin^{-1}(2r/R),\pi/2]$ equals the value of $\tau_1$ for $(\pi-\theta_1)\in[\sin^{-1}(2r/R),\pi/2]]$. Likewise, the value of $d_0$ for $\theta_1\in[\pi-\sin^{-1}(2r/R),\pi/2]$ equals $d_2$ for $(\pi-\theta_1)\in[\sin^{-1}(2r/R),\pi/2]$ and $d_2$ for $\theta_1\in[\pi-\sin^{-1}(2r/R),\pi/2]$ equls $d_0$ for $(\pi-\theta_1)\in[\sin^{-1}(2r/R),\pi/2]$.
Thus, in \eqref{eq:NewCoordinateLegnth} we can replace $\theta_1$ with $\pi-\theta_1$ and swap the index of $\tau_0$ and $\tau_1$ and of $d_0$ and $d_2$ to get the same equations as in \eqref{eq:NewCoordinateLegnth}.
\end{proof}
\begin{proposition}\label{prop:casea0largetheta}
    If \eqref{eqJZHypCond}, $\phistar\in\big(0,\tan^{-1}(1/3)\big)$ and $\theta_1\in\big(\frac{\pi}{6},\frac{5\pi}{6}\big)$ in \eqref{eqPtsFromeqMhatOrbSeg}
    , then  $\mathcal{I}<-1.37$ (see \cref{def:AformularForLowerBoundOfExpansion}).
\end{proposition}
\begin{proof}
 In \cref{fig:fractionEstimate2}, $P=p(x_1)$, $Q=p(x_2)$, $T=p(x_0)$ in $\Gamma_R$. By symmetry, we assume $\theta_1\in (\frac{\pi}{6},\frac{\pi}{2}]$. $\Phi_1\in(-\Phistar,\Phistar)$. Then $\measuredangle{P O_R Or}=|\Phi_1|$, $\measuredangle{P O_r Y}=|\alpha|$, $\measuredangle{O_R P O_r}=|\delta|$. 
 Using the coordinates and notation of \cref{def:TPPORsystem}, we have $\rho=|O_rP|$ and $\delta\in(-\phistar+\Phistar,+\phistar-\Phistar)\subset (-\phistar,+\phistar)$.

From \cref{def:AformularForLowerBoundOfExpansion} we note that 
\[
\mathcal{I} = -1+\frac{\tau_1}{d_0}\big[\frac{2(\tau_0-d_0)}{d_1}-\frac{\tau_0+\tau_1-2d_0}{\tau_1}\big]=-1-\frac{\tau_0+\tau_1-2d_0}{d_0}+\frac{\tau_1}{d_0}\frac{2(\tau_0-d_0)}{d_1}.
\]
To estimate the two fractions, we begin with \cref{lemma:newcoordinateslength} and some trigonometry:\COMMENT{\Large\bfseries\red Omit the denominator $d_0$ from \cref{{eq:casea0largetheta_trigonometrynegativedelta}}!!\newline\wentao{It will save two lines but add more sentence for conclusion.}\Large\color{red}\boris No, not at all, just leave \eqref{eq:casea0largetheta_trigonometrypositivedelta} exactly as it is.\Hrule Note and fix missing punctuation. Throughout.}
\begin{equation}\label{eq:casea0largetheta_trigonometrynegativedelta}
\begin{aligned}
\tau_0+\tau_1-2d_0& \overbracket{=}^{\mathclap{\eqref{eq:NewCoordinateLegnth}}}\rho\sin{(\theta_1-\delta)}+\rho\sin{(\theta_1+\delta)}-\sqrt{r^2-\rho^2\cos^2{(\theta_1+\delta)}}+ \sqrt{r^2-\rho^2\cos^2{(\theta_1-\delta)}}\\
& = 2\rho\sin\theta_1\cos\delta+\frac{r^2-\rho^2\cos^2{(\theta_1-\delta)}-(r^2-\rho^2\cos^2{(\theta_1+\delta)})}{\sqrt{r^2-\rho^2\cos^2{(\theta_1+\delta)}}+ \sqrt{r^2-\rho^2\cos^2{(\theta_1-\delta)}}}\\
& \overbracket{=}^{\mathllap{\cos^2x\frac12\cos2x+\frac12}}\rho\big[2\sin{\theta_1}\cos{\delta}+\frac{1}{2}\frac{\rho[\cos{(2\theta_1+2\delta)}-\cos{(2\theta_1-2\delta)}]}{\sqrt{r^2-\rho^2\cos^2{(\theta_1+\delta)}}+ \sqrt{r^2-\rho^2\cos^2{(\theta_1-\delta)}}}\big]\\
& \overbracket{=}^{\mathllap{2\sin x\sin y=\cos(x-y)-\cos(x+y)}}\rho\big[2\sin{\theta_1}\cos{\delta}-\frac{\rho\sin{(2\delta)}\sin{(2\theta_1)}}{\sqrt{r^2-\rho^2\cos^2{(\theta_1+\delta)}}+ \sqrt{r^2-\rho^2\cos^2{(\theta_1-\delta)}}}\big]\\
&\overbracket{=}^{\mathllap{\sin2x=2\sin x\cos x}}2\rho\sin{\theta_1}\cos{\delta}\big[1-\frac{2\rho\cos{\theta_1}\sin{\delta}}{\sqrt{r^2-\rho^2\cos^2{(\theta_1+\delta)}}+ \sqrt{r^2-\rho^2\cos^2{(\theta_1-\delta)}}}\big]
\end{aligned}
\end{equation}
Now, if $\delta\le0$, then $\displaystyle1-\frac{2\rho\cos{\theta_1}\sin{\delta}}{\sqrt{r^2-\rho^2\cos^2{(\theta_1+\delta)}}+ \sqrt{r^2-\rho^2\cos^2{(\theta_1-\delta)}}}\ge1$.
If $\delta>0$, then $0<\delta<\phistar<\frac{\pi}{6}$, $\frac{\phi}{6}<\theta_1<\frac{\pi}{2}$ and trigonometry give $\tan\delta<1/3$ and  $\cot\theta_1<\sqrt{3}$, so
\begin{equation}\label{eq:casea0largetheta_trigonometrypositivedelta}
\begin{multlined}
1-\frac{2\rho\cos{\theta_1}\sin{\delta}}{\sqrt{r^2-\rho^2\cos^2{(\theta_1+\delta)}}+ \sqrt{r^2-\rho^2\cos^2{(\theta_1-\delta)}}}
\overbracket{>}^{\mathclap{\rho<r}}  1-\frac{2\rho\cos{\theta_1}\sin{\delta}}{\sqrt{\rho^2-\rho^2\cos^2{(\theta_1+\delta)}}+ \sqrt{\rho^2-\rho^2\cos^2{(\theta_1-\delta)}}}\\
=  1-\frac{2\cos{\theta_1}\sin{\delta}}{\sin{(\theta_1+\delta)}+\sin{(\theta_1-\delta)}}=1-\cot{\theta_1}\tan{\delta}
>  1-\sqrt{3}\cdot\frac{1}{3}=\frac{\sqrt3-1}{\sqrt3}.
\end{multlined}
\end{equation}
$0<\phistar<\tan^{-1}{(\frac{1}{3})}$, $\theta_1\in(\frac{\pi}{6},\frac{\pi}{2}]$ together with $\rho>r\cos\phistar$ imply that $|\delta|<\phistar<\tan^{-1}{(\frac{1}{3})},\,\cos{\delta}>\cos{\phistar}>\cos\tan^{-1}\frac13=\sqrt8/3$,\COMMENT{Do we need 0.94? Or 0.4 above?\newline\wentao{no need now. we have 0.375 now in use. Do you agree no harm to mention 0.4,0.94?}\boris This is a VERY complicated paper already, so the writing should minimize any distraction and needless items. What is the absolute simplest way (for the READER) to present this here?\wentao{This is a very good question. We can make a remark suggesting that for $\phistar\in(0,\tan^{-1}(1/2)$ really our computation works! Now 0.94,0.4 are both removed.}} $\sin\theta_1>1/2$. Therefore, for either $\delta>0$ or $\delta\le 0$, we have
\begin{equation}\label{eq:casea0largethetalengthsfraction}
\frac{\tau_0+\tau_1-2d_0}{d_0}\overbracket{>}^{{\textstyle\bigstrut}\mathclap{(\ref{eq:casea0largetheta_trigonometrynegativedelta},\ref{eq:casea0largetheta_trigonometrypositivedelta})}}
\frac{\sqrt3-1}{\sqrt3}\frac{\overbracket{\rho}^{\mathclap{>r\cos\phistar}}}{\underbracket{d_0}_{\mathclap{\le r}}}\overbracket{2\sin{\theta_1}}^{>1}\overbracket{\cos{\delta}}^{>\cos\phistar}>\frac{\sqrt3-1}{\sqrt3}\cos^2\phistar>\frac89\frac{\sqrt3-1}{\sqrt3}>0.375.
\end{equation}
On the other hand, $\theta_1\in(\frac{\pi}{6},\frac{\pi}{2}]$ therefore $\sin{\theta_1}>\frac{1}{2}$, $-d_0\le\tau_0-d_0\le d_0$, $0<\tau_1\le 2r$, $R>1700r$ and $d_1=R\sin\theta_1>\frac12R>850r$ imply 
\begin{equation}\label{eq:casea0largethetalengthsfractionII}
    \begin{aligned}
        \frac{\overbracket{\tau_1}^{\le 2r}}{d_0}\frac{2(\tau_0-d_0)}{\underbracket{d_1}_{\mathclap{=R\sin\theta_1>\frac12R>850r}}}\le\frac{2}{850}\cdot2\big|\frac{\tau_0-d_0}{d_0}\big|<\frac{1}{425}\cdot2=\frac{4}{850}.
    \end{aligned}
\end{equation}
Hence, from \cref{def:AformularForLowerBoundOfExpansion}, \[
\mathcal{I} = -1+\frac{\tau_1}{d_0}\big[\frac{2(\tau_0-d_0)}{d_1}-\frac{\tau_0+\tau_1-2d_0}{\tau_1}\big]=-1-\frac{\tau_0+\tau_1-2d_0}{d_0}+\frac{\tau_1}{d_0}\frac{2(\tau_0-d_0)}{d_1}\overbracket{<}^{\mathrlap{\eqref{eq:casea0largethetalengthsfraction},\eqref{eq:casea0largethetalengthsfractionII}}\quad}
-1-0.375+\frac{4}{850}<-1.37.\qedhere\]
\end{proof}
For later use in \cref{theorem:UniformExpansionA0,corollary:halfquadrantcone}, we note:\begin{corollary}\label{Corollary:CollisionOnMrin0Mrin1}
  Suppose $\phistar\in\big(0,\tan^{-1}(1/3)\big)$ and \eqref{eqJZHypCond} holds. If the return orbit segment  \eqref{eqPtsFromeqMhatOrbSeg} of $x$ from \cref{def:MhatReturnOrbitSegment,def:fpclf} has $x_0\in \Mrout$, and $\Fc(x_0)=x_1=(\Phi_1,\theta_1)\in \MRin$, $\Fc(x_1)=x_2=(\phi_2,\theta_2)\in \textcolor{blue}{M^{\text{\upshape{in}}}_{r,0}}\cup \textcolor{blue}{M^{\text{\upshape{in}}}_{r,1}}$, then
  \begin{itemize}
      \item $\phi_2\in(\pi/2,3\pi/2)$, i.e. $p(x_2)$ is on the upper half of circle $C_r$ in \cref{fig:fractionEstimate2}.
      \item $\theta_1\in(\pi/6,5\pi/6)$.
      \item  $\mathcal{I}$ from \cref{def:AformularForLowerBoundOfExpansion} satisfies $\mathcal{I}<-1.37$.
  \end{itemize}
\end{corollary}
\begin{proof}
    In \cref{fig:MrN}, the leftmost point of the closure of \textcolor{blue}{$M^{\text{\upshape in}}_{r,1}$} has $\phi$-coordinate $\frac{2\pi}{3}-\phistar>\frac{\pi}{2}$, and by symmetry the rightmost point of the closure of \textcolor{blue}{$M^{\text{\upshape in}}_{r,1}$} has $\phi$-coordinate $\frac{4\pi}{3}+\phistar<\frac{3\pi}{2}$. The leftmost point of the closure of \textcolor{blue}{$M^{\text{\upshape in}}_{r,0}$} has $\phi$-coordinate $\pi-\phistar>\frac{\pi}{2}$ and the rightmost point of the closure of \textcolor{blue}{$M^{\text{\upshape in}}_{r,0}$} has $\phi$-coordinate $\pi+\phistar<\frac{3\pi}{2}$. Therefore,  $\theta_1\in (\frac{\pi}{6},\frac{5\pi}{6})$ by \cref{lemma:returnorbittheta1Large}. Therefore,\COMMENT{\small\bfseries\color{red}"Furthermore" is the wrong word; please look up what it means. "Therefore" is more likely the right word here. Please look it up as well. It would be good to search the file for both and make sure that each is always used correctly.\wentao{OK.}} $\mathcal{I}<-1.37$ by \cref{prop:casea0largetheta}.
\end{proof}
\begin{remark}
In the computation \eqref{eq:casea0largetheta_trigonometrypositivedelta} and \cref{lemma:returnorbittheta1Large}, we see that we can relax the condition in \eqref{eqJZHypCond} to be $\phistar<\pi/6$ and still obtain $\mathcal{I}<-1$. However, when $\phistar$ is close to $\pi/6$, in \eqref{eq:casea0largetheta_trigonometrypositivedelta}, $R$ may need to be arbitrarily large to yield $\mathcal{I}<-1$. We avoid computing $R_{\text{\upshape HF}}$ in such a complicated situation.






\end{remark}

\begin{figure}[h]
\begin{center}
\begin{sideways}
\begin{tikzpicture}[xscale=0.7, yscale=0.7]
    \tkzDefPoint(0,0){Or};
    \tkzDefPoint(0,-4.8){Y};
    \pgfmathsetmacro{\rradius}{4.5};
    \pgfmathsetmacro{\bdist}{8.7};
    \tkzDefPoint(0,\bdist){OR};
    \pgfmathsetmacro{\phistardeg}{40};
    \pgfmathsetmacro{\XValueArc}{\rradius*sin(\phistardeg)};
    \pgfmathsetmacro{\YValueArc}{\rradius*cos(\phistardeg)};
    \tkzDefPoint(\XValueArc,-1.0*\YValueArc){B};
    \tkzDefPoint(-1.0*\XValueArc,-1.0*\YValueArc){A};
    \pgfmathsetmacro{\Rradius}{veclen(\XValueArc,\YValueArc+\bdist)};
    \draw[name path=Cr,dashed,ultra thin] (Or) circle (\rradius);
    \tkzDrawArc[very thick](Or,B)(A);
    \tkzDrawArc[name path=CR, very thick](OR,A)(B);
    \pgfmathsetmacro{\PXValue}{\Rradius*sin(7)};
    \pgfmathsetmacro{\PYValue}{\bdist-\Rradius*cos(7)};
    \tkzDefPoint(\PXValue,\PYValue){P};
    \draw[name path=PQ_Far,blue,ultra thin](P)--+(70:7)coordinate(PposX);
    \draw[name path=PT_Far,blue,ultra thin](P)--+(124:9);
    \draw[name path=PT0_Far,blue,dashed,ultra thin](P)--+(-56:7);
    \draw[name path=PQ0_Far,blue,dashed,ultra thin](P)--+(-110:9);
    \draw[red,dashed,->](P)--+(7:4)coordinate(PposX);
    \draw[name path=PPnegX,red,dashed](P)--+(-173:6)coordinate(PnegX);
    \draw[name path=linePOR,red,dashed,->](P)--(OR);
    \draw[name path=linePOr,red,dashed](P)--(Or);
    \tkzDefPoint(0,-5.5){Y};
    \tkzDefPoint(0,\bdist-\Rradius){C};
    \draw[red,dashed](OR)--(Y);
    \tkzMarkAngle[arc=lll,ultra thin,size=1.7](Or,OR,P);
    \tkzMarkAngle[arc=lll,ultra thin,size=0.7](Y,Or,P);
    \tkzMarkAngle[arc=lll,ultra thin,size=0.8](OR,P,Or);
    \path[name intersections={of=Cr and PQ_Far,by={Q}}];
    \path[name intersections={of=Cr and PQ0_Far,by={Q0}}];
    \path[name intersections={of=Cr and PT0_Far,by={T0}}];
    \path[name intersections={of=Cr and PT_Far,by={T}}];
    \tkzMarkAngle[arc=lll,very thin,size=1.8](PposX,P,Q);
    \tkzMarkRightAngle(OR,P,PnegX);
    \node[rotate=-90] at ($(Or)+(-78:1)$)  {\Large $\alpha$};
    \node[rotate=-90] at ($(OR)+(-85:2.3)$)  {\Large \textcolor{blue}{$\Phi_1$}};
    \node[rotate=-90] at ($(P)+(105:1.3)$)  {\Large {$\delta$}};
    \node[rotate=-90] at ($(P)+(45:2.3)$)  {\Large {$\theta_1$}};
    
    \tkzLabelPoint[left,rotate=-90](A){$A$};
    \tkzLabelPoint[left,rotate=-90](B){$B$};
    \tkzLabelPoint[left,rotate=-90](OR){$O_R$};
    \tkzLabelPoint[below,rotate=-90](Or){$O_r$};
    \tkzLabelPoint[above,rotate=-90]([shift={(0,0.2)}]P){$P$};
    \tkzLabelPoint[above,rotate=-90](Q){$Q$};
    \tkzLabelPoint[left,rotate=-90](T){$T$};
    \tkzLabelPoint[below,rotate=-90](Q0){$Q_0$};
    \tkzLabelPoint[above,rotate=-90](T0){$T_0$};
    \tkzLabelPoint[above,rotate=-90](PposX){$TP$};
    \tkzDefLine[perpendicular=through Or](P,PnegX)\tkzGetPoint{M_far};
    \path[name path= lineOrMfar] (Or)--(M_far);
    \path[name path= PPnegX] (P)--(PnegX);
    \path[name path= PPposX] (P)--(PposX);
    \path[name intersections={of=lineOrMfar and PPnegX,by={M}}];
    \path[name intersections={of=Cr and PPnegX,by={A0}}];
    \path[name intersections={of=Cr and PPposX,by={B0}}];
    \coordinate (Or_mirror) at ($(Or)!2!(M)$);
    \tkzDrawArc[name path=mirrorCr, ultra thin, dashed](Or_mirror,A0)(B0);
    \tkzLabelPoint[below,rotate=-90](Or_mirror){$O'_r$};

\path[name intersections={of=mirrorCr and PQ0_Far,by={Q1}}];
    \path[name intersections={of=mirrorCr and PT0_Far,by={T1}}];
    \tkzLabelPoint[above,rotate=-90](Q1){$Q_1$};
    \tkzLabelPoint[above,rotate=-90](T1){$T_1$};
    \tkzLabelPoint[left,rotate=-90](C){$M$};    \tkzDrawPoints(A,Or,OR,P,B,T,Q,Or_mirror,A0,T0,Q0,T1,Q1,C);
    \begin{scope}[ultra thin,decoration={markings, mark=at position 0.55 with {\arrow{>}}}
    ] 
    \draw[blue,postaction={decorate}](P)--(Q);
    \draw[blue,postaction={decorate}](T)--(P);
    \end{scope}
    \tkzLabelPoint[above,rotate=-90]([shift={(0.15*\rradius,1.1*\rradius)}]Or){\Large $\Gamma_r$};
    \tkzLabelPoint[above,rotate=-90]([shift={(0.35*\rradius,-1.05*\rradius)}]Or_mirror){\Large $\Gamma_r'$};
\end{tikzpicture}
\end{sideways}
\captionsetup{type=figure}
\caption{Reflecting the billiard over the tangent line of $\Gamma_R$ at $P\dfn p(x_1)$. 
Here, $|TT_1|=\tau_0+\tau_1$, $|TT_0|=2d_0$, $|T_0T_1|=\tau_0+\tau_1-2d_0>0$, $|QQ_0|=2d_2$, $|Q_1Q_0|=\tau_0+\tau_1$, $|Q_0Q_1|=\tau_0+\tau_1-2d_2>0$. (See also \cite[Figure 7, equation (4.1)]{hoalb}.)}\label{fig:heuristicproof}
\end{center}
\end{figure}

\subsection{Proof of \texorpdfstring{\cref{Prop:Ilowerbound}}{Proposition of Lower bound for I}, small angles}\label{SBSProp53PartII}
We now consider the situation when the orbit segment \eqref{eqPtsFromeqMhatOrbSeg} of $x$ in \cref{def:MhatReturnOrbitSegment,def:fpclf} has a less transverse reflection on $\Gamma_R$, that is, with $x_1=(\Phi_1,\theta_1)\in R_1\cup R_2\cup R_3\cup R_4\subset\MRin\cap\MRout$ and $\sqrt{4r/R}\le\sin\theta_1
\le\frac{1}{2}$. We know that $\mathcal{I}<-1.37$ if $\sin\theta_1>1/2$ (\cref{prop:casea0largetheta}). 
Now we prove that $\mathcal{I}<-1$ if $x_1=(\Phi_1,\theta_1)\in R_1\cup R_2\cup R_3\cup R_4$ with $\sin\theta_1\ge\sqrt{4r/R}$. By symmetry, we can restrict our analysis to $\theta_1\in [\sin^{-1}(2r/R),\frac{\pi}{2}]$. i.e., the region $x_1\in R_1\cup R_2$  in \cref{fig:R1}. We deal with $x_1\in R_2$ in \cref{prop:casea0smallthetaR2} and with $x_1\in R_1$ in \cref{prop:casea0smallthetaR1}.

\begin{proposition}[$x_1\in R_2$]\label{prop:casea0smallthetaR2}
Suppose $\phistar\in\big(0,\tan^{-1}(1/3)\big)$,  \eqref{eqJZHypCond} holds, and the return orbit segment \eqref{eqPtsFromeqMhatOrbSeg} of $x$ from \cref{def:MhatReturnOrbitSegment,def:fpclf} has $x_0\in\Mrout$, $\Fc(x_0)=x_1=(\Phi_1,\theta_1)\in\MRin$, $x_2=\Fc(x_1)\in \Mrin$ with $\theta_1\in\big[\sin^{-1}{(2r/R)},\pi/2\big]$, and $\Phi_1\in(-\Phistar,0)$, that is, $x_1\in R_2$ in \cref{fig:R1} with $\sin\theta_1\ge2r/R$. Then $\mathcal{I}\le-1$  (\cref{def:AformularForLowerBoundOfExpansion}). In the special case where $\sin{\theta_1}\ge\sqrt{4r/R}$, we moreover obtain $\mathcal{I}<-1$.
\end{proposition}
\begin{proof}
For $x_1\in R_2$ in \cref{fig:R1}, thus, with $d_1=R\sin\theta_1\ge2r$ and $\tau_1>0$ in \cref{def:AformularForLowerBoundOfExpansion}, we can analyze $\mathcal{I}=-1+\frac{\tau_1}{d_0}\Big[\frac{2(\tau_0-d_0)}{d_1}-\frac{\tau_0+\tau_1-2d_0}{\tau_1}\Big]$ in the following three cases.
\newline
\textbf{Case i} $\tau_0\le d_0$, then since $\tau_0+\tau_1-2d_0>0$ (see \cref{remark:deffpclf}), we have $\frac{2(\tau_0-d_0)}{d_1}-\frac{\tau_0+\tau_1-2d_0}{\tau_1}<0$. Thus $\mathcal{I}<-1$ .
\newline
\textbf{Case ii} $\tau_0>d_0$ and $\tau_1\ge d_0$, then we have
\begin{equation}\label{eq:casea0smallthetaR2case2}
\begin{aligned}
& \frac{\tau_0+\tau_1-2d_0}{\tau_1}=1-\frac{2d_0-\tau_0}{\tau_1}\overbracket{\ge}^{\strut\mathclap{2d_0-\tau_0>0,\tau_1\ge d_0}} 1-\frac{2d_0-\tau_0}{d_0}=\frac{2(\tau_0-d_0)}{2d_0}\overbracket{\ge}^{\strut\mathclap{2d_0 \le 2r\le d_1,\tau_0>d_0}}\frac{2(\tau_0-d_0)}{d_1}
\end{aligned}
\end{equation}
Thus, \eqref{eq:casea0smallthetaR2case2} implies $\mathcal{I}\le-1$. In particular, if $\sin\theta_1\ge\sqrt{4r/R}>2r/R$, then $d_1=R\sin\theta_1>2r\ge2d_0$ and the last inequality in \eqref{eq:casea0smallthetaR2case2} becomes $\frac{2(\tau_0-d_0)}{2d_0}\overbracket{>}^{\strut\mathclap{2d_0 \le 2r< d_1,\tau_0>d_0}}\frac{2(\tau_0-d_0)}{d_1}$. Thus, if $\sin\theta_1\ge\sqrt{4r/R}$, then $\mathcal{I}<-1$.
\newline
\textbf{Case iii} $\tau_0>d_0>\tau_1$. We show that this cannot occur if $x_1\in R_2=\{(\phi_1,\theta_1)|\theta_1\in[\sin^{-1}{(2r/R)},\pi/2],\Phi_1\in(-\Phistar,0)\}$.

In \cref{fig:mirrorR2}, let $P=p(x_1)\in \Gamma_R$ be on the table boundary. $Q=p(x_2)=p(\Fc(x_1))$, $T=p(x_0)=p(\Fc^{-1}(x_1))$. In the coordinate system of \cref{def:StandardCoordinateTable}, $P$ has a negative $x$ coordinate. And since $O_R=(0,b)$, $O_r=(0,0)$, suppose $C_1$ and $C_2$ are the two intersections  of the line $O_RP$ with the circle $C_r$. We have $|C_1C_2|<2r$ and the arc $\arc{C_1TC_2}$ with a length less than half the perimeter of $C_r$. $Q$ and $T$ are on two sides of line $O_RP$. The arc $\arc{C_1QC_2}$ is having a length more than the half of the perimeter of $C_r$. 

Then the mirror symmetry of arc $\arc{C_1TC_2}$ with respect to the line $O_RP$ is contained inside the disk surrounded by $C_r$ . $T_1$ which is the mirror symmetry $T$ of the line $O_RP$ must be inside the disk surrounded by $C_r$ \big(dashed arc in \cref{fig:mirrorR2}  $\arc{C_1T_1C_2}$\big). 

Hence, in \cref{fig:mirrorR2}, $\tau_1=|PQ|\ge|PT_1|=|TP|=\tau_0$, 
so \textbf{Case iii}, $\tau_0>d_0>\tau_1$, cannot occur.
\end{proof}
\begin{figure}[h]
\begin{minipage}{.5\textwidth}
\begin{center}
\begin{tikzpicture}[xscale=0.5,yscale=0.5]
    \tkzDefPoint(0,0){Or};
    \tkzDefPoint(0,-4.8){Y};
    \pgfmathsetmacro{\rradius}{4.5};
    \pgfmathsetmacro{\bdist}{8.7};
    \tkzDefPoint(0,\bdist){OR};
    \pgfmathsetmacro{\phistardeg}{40};
    \pgfmathsetmacro{\XValueArc}{\rradius*sin(\phistardeg)};
    \pgfmathsetmacro{\YValueArc}{\rradius*cos(\phistardeg)};
    \tkzDefPoint(\XValueArc,-1.0*\YValueArc){B};
    \tkzDefPoint(-1.0*\XValueArc,-1.0*\YValueArc){A};
    \pgfmathsetmacro{\Rradius}{veclen(\XValueArc,\YValueArc+\bdist)};
    \draw[name path=Cr,ultra thin,dashed] (Or) circle (\rradius);
    \tkzDrawArc[name path=CR,very thick](OR,A)(B);
    \tkzDrawArc[very thick](Or,B)(A);
    \pgfmathsetmacro{\PXValue}{-\Rradius*sin(7)};
    \pgfmathsetmacro{\PYValue}{\bdist-\Rradius*cos(7)};
    \tkzDefPoint(\PXValue,\PYValue){P};
    \path[](OR) -- (P) -- ([turn]90:6)coordinate (PosX);
    \path[](OR) -- (P) -- ([turn]-90:3)coordinate (NegX);
    \path[](OR) -- (P) -- ([turn]0:1)coordinate (NegY);
    \path[](OR) -- (P) -- ([turn]155:9)coordinate (Q_far);
    \path[](OR) -- (P) -- ([turn]-155:9)coordinate (T_far);
    \draw[red,ultra thin,dashed,->](P)--(PosX);
    \draw[red,ultra thin,dashed](P)--(NegX);
    \draw[red,ultra thin,dashed](OR)--(NegY);
    \path[name path=Xpos](P)--(PosX);
    \path[name path=Xneg](P)--(NegX);
    \path[name path=Yneg](P)--(NegY);
    \path[name path=linePQ_far](P)--(Q_far);
    \path[name path=linePT_far](P)--(T_far);
    \path[name path=lineOR_P,ultra thin, dashed,->](P)--(OR);
    \path[name intersections={of=lineOR_P and Cr,by={D}}];
    \path[name intersections={of=linePQ_far and Cr,by={Q}}];
    \path[name intersections={of=linePT_far and Cr,by={T}}];
    \path[name intersections={of=Yneg and Cr,by={C}}];
    \begin{scope}[decoration={
    markings,
    mark=at position 0.5 with {\arrow{>}}}
    ] 
    \draw[blue,postaction={decorate}](T)--(P);
    \draw[blue,postaction={decorate}](P)--(Q);
    \end{scope}
    \coordinate (M) at ($(C)!0.5!(D)$);
    \coordinate (Or_mirror) at ($(Or)!2!(M)$);  
    \tkzDrawArc[name path= Crmirror,red,dashed](Or_mirror,C)(D);
    \tkzDefPointBy[projection=onto P--D](T)\tkzGetPoint{H};
    \coordinate (T1) at ($(T)!2!(H)$);
    \tkzDrawPoints(Or,OR,P,Q,T,T1,C,D);
    \tkzLabelPoint[above](OR){$O_R$};
    \tkzLabelPoint[above](Or){$O_r$};
    \tkzLabelPoint[above left](P){$P$};
    \tkzLabelPoint[right](Q){$Q$};
    \tkzLabelPoint[left](T){$T$};
    \tkzLabelPoint[left](T1){$T_1$};
    \tkzLabelPoint[below left](C){$C_1$};
    \tkzLabelPoint[above left](D){$C_2$};
\end{tikzpicture}
    \caption{Comparing $\tau_0,\tau_1$ for $x_1=(\Phi_1,\theta_1)\in R_2$ with $\Phi_1\in(-\Phistar,0)$, $\theta_1\in [\sin^{-1}{(2r/R)},\pi/2]$}\label{fig:mirrorR2}
\end{center}
\end{minipage}\hfill
\begin{minipage}{.5\textwidth}
\begin{center}
    \begin{tikzpicture}[xscale=0.55,yscale=.55]
        \pgfmathsetmacro{\PHISTAR}{1.5};
        \pgfmathsetmacro{\Dlt}{0.4};
        \pgfmathsetmacro{\Dltplus}{\Dlt+0.6}
        \pgfmathsetmacro{\WIDTH}{12};
        \tkzDefPoint(-\PHISTAR,0){LL};
        \tkzDefPoint(-\PHISTAR,\WIDTH){LU};
        \tkzDefPoint(\PHISTAR,0){RL};
        \tkzDefPoint(\PHISTAR,\WIDTH){RU};
        \tkzDefPoint(-\PHISTAR,\PHISTAR){X};
        \tkzDefPoint(\PHISTAR,\PHISTAR){IY};
        \tkzDefPoint(-\PHISTAR,\PHISTAR+\Dlt){LLD};
        \tkzDefPoint(\PHISTAR,\PHISTAR+\Dlt){RLD};
        \tkzDefPoint(0,\PHISTAR+\Dltplus){LLDp};
        \tkzDefPoint(\PHISTAR,\PHISTAR+\Dltplus){RLDp};
        \tkzDefPoint (0.3*\PHISTAR,\PHISTAR+\Dltplus){montonestart};
        \tkzDefPoint (0.3*\PHISTAR,0.333*\WIDTH){montoneend};
        \draw [thin] (LL) --(LU);
        \draw [thin] (LU) --(RU);
        \draw [thin] (RL) --(RU);
        \draw [thin] (LL) --(RL);
        \draw [blue,ultra thin,dashed] (LL) --(IY);
        \draw [red,ultra thin,dashed]  (X)--(RL);
        \draw [purple,ultra thin,dashed]  (LLD)--(RLD);
        \draw [green,ultra thin]  (LLDp)--(RLDp);
        \draw [green,ultra thin]  (montonestart)--(montoneend);
        \tkzLabelPoint[right](RU){$\pi/2$};
        \tkzLabelPoint[right](RLD){$\sin{\theta_1}=2r/R$};
        \tkzLabelPoint[right](RLDp){$\sin{\theta_1}=\sqrt{4r/R}$};
        \coordinate (RM) at ($(RL)!0.333!(RU)$);
        \coordinate (UpperM) at ($(LU)!0.5!(RU)$);
        \coordinate (LowerM) at ($(LL)!0.5!(RL)$);
        \coordinate (LM) at ($(LowerM)!0.333!(UpperM)$);
        \coordinate (x1) at ($(montonestart)!0.5!(montoneend)$);
        \tkzLabelPoint[right](RM){$\pi/6$};
        \draw [purple,ultra thin]  (LM)--(RM);
        \draw [dashed,ultra thin]  (LowerM)--(UpperM);
        \node[] at (0.5*\PHISTAR, 0.8*\PHISTAR+0.8*\Dlt+0.6*\WIDTH){\Large $R_1$};
        \node[] at (-0.5*\PHISTAR, 0.8*\PHISTAR+0.8*\Dlt+0.6*\WIDTH){\Large $R_2$};
        \tkzLabelPoint[below](LL){$-\Phistar$};
        \tkzLabelPoint[below](RL){$+\Phistar$};
        \tkzLabelPoint[below](LowerM){$0$};
        \tkzLabelPoint[left](LLDp){\tiny $x_M$};
        \tkzLabelPoint[below](montonestart){\tiny $x_L$};
        \tkzLabelPoint[above](montoneend){\tiny $x_U$};
        \tkzLabelPoint[right](x1){\tiny $x_1$};
        \tkzDrawPoints(LLDp,montonestart,montoneend,x1);
    \end{tikzpicture}
    \caption{For $x_1\in R_1$ with \newline $\frac{1}{2}\ge\sin{(\theta_1)}\ge\sqrt{4r/R}$, $\frac{2d_0-\tau_0}{\tau_1}$ \newline takes maximum value at $x_M$}\label{fig:R1_fraction_extreme}
\end{center}
\end{minipage}
\end{figure}

\begin{proposition}[$x_1\in R_1$]\label{prop:casea0smallthetaR1}
Suppose $\phistar\in\big(0,\tan^{-1}(1/3)\big)$, \eqref{eqJZHypCond} holds, and the return orbit segment \eqref{eqPtsFromeqMhatOrbSeg} of $x$ from \cref{def:MhatReturnOrbitSegment,def:fpclf} has $x_0\in\Mrout$, $\Fc(x_0)=x_1=(\Phi_1,\theta_1)\in\MRin$ and $x_2=\Fc(x_1)\in \Mrin$ with $\theta_1\in[\sin^{-1}{(\sqrt{4r/R})},\pi/6],\Phi_1\in[0,\Phistar)$, that is, $x_1\in R_1$ in \cref{fig:R1_fraction_extreme} with $1/2\ge\sin\theta_1\ge\sqrt{4r/R}$. Then $\mathcal{I}<-1$ (see \cref{def:AformularForLowerBoundOfExpansion}).
\end{proposition}


\begin{proof}
$x_1=(\Phi_1,\theta_1)\in R_1$ with $\theta_1\in[\sin^{-1}{(\sqrt{4r/R})},\frac{\pi}{2}]$ and $\Phi_1\in[0,\Phistar)$ is on a vertical line segment inside $R_1$ connected to its upper and lower end point $x_U=(\Phi_1,\pi/6)$ and $x_L=(\Phi_1,\sin^{-1}(\sqrt{4r/R}))$ as shown in \cref{fig:R1_fraction_extreme}. And $x_M=(0,\sin^{-1}(\sqrt{4r/R}))$ and $x_L$ are on the same green horizontal line in \cref{fig:R1_fraction_extreme}.

We consider $\tau_1(x_1)$, $\tau_0(x_1)$, $d_2(x_1)$, $d_0(x_1)$, $d_1(x_1)$ in \cref{def:fpclf} as functions of $x_1=(\Phi_1,\theta_1)$ on $R_1$. 

For $x_1$ being on the line segment $[x_L,x_U]$ (vertical green line in \cref{fig:R1_fraction_extreme}), in the context of \cref{monotone_p1,monotone_p2} since $\pi/6<\pi/2-\phistar$, \cref{monotone_p2,rmk:monotone_p2} imply that $(2d_0-\tau_0)\bigm|_{x_1=x_L}\ge (2d_0-\tau_0)(x_1)>0$ and \cref{monotone_p1} implies that $0<\tau_1\bigm|_{x_1=x_L}\le \tau_1(x_1)$. Therefore, $\frac{2d_0-\tau_0}{\tau_1}(x_1)\le \frac{2d_0-\tau_0}{\tau_1}\bigm|_{x_1=x_L}$.

On the other hand, for $x_1$ being on the line segment $[x_M,x_L]$ (horizontal green line in \cref{fig:R1_fraction_extreme}), \cref{monotone_p5} implies that $\frac{2d_0-\tau_0}{\tau_1}\bigm|_{x_1=x_L}\le \frac{2d_0-\tau_0}{\tau_1}\bigm|_{x_1=x_M}$.

Hence, for all $x_1\in R_1$ with $\sqrt{4r/R}\le\sin\theta_1\le1/2$, we have $\frac{2d_0-\tau_0}{\tau_1}(x_1)\le\frac{2d_0-\tau_0}{\tau_1}\bigm|_{x_1=x_M}$. That is $\frac{\tau_1+\tau_0-2d_0}{\tau_1}(x_1)=1-\frac{2d_0-\tau_0}{\tau_1}(x_1)\ge 1-\frac{2d_0-\tau_0}{\tau_1}\bigm|_{x_1=x_M}=\frac{\tau_1+\tau_0-2d_0}{\tau_1}\bigm|_{x_1=x_M}$. Then we will go to compute this lower bound $\frac{\tau_1+\tau_0-2d_0}{\tau_1}\bigm|_{x_1=x_M}$.

Note that if $(\Phi_M,\theta_M)=x_1=x_M=(0,\sin^{-1}(\sqrt{4r/R}))$, that is $p(x_1)=P=M$ the midpoint of the arc $\Gamma_R$ in \cref{fig:fractionEstimate2} with $\delta=0$, then in \eqref{eq:NewCoordinateLegnth} we have $\rho=|O_rP|=|O_rM|=r\cos\phistar+2R\sin^2{(\frac{\Phistar}{2})}>r\cos\phistar\overset{\phistar<\pi/3}>\frac{r}{2}$. We further get

\begin{equation}\label{eq:new73}
\big(\tau_0+\tau_1-2d_0\big)\bigm|_{x_1=x_M} \overset{\eqref{eq:NewCoordinateLegnth}\text{ with }\delta=0}{\scalebox{6}[1]{$=$}}2\rho\sin{\theta_M}.
\end{equation}
Therefore,
\begin{equation}\label{eq:new74}
    \begin{aligned}
        \big(\frac{\tau_1+\tau_0-2d_0}{\tau_1}\big)(x_1)&\ge\big(\frac{\tau_0+\tau_1-2d_0}{\tau_1}\big)\bigm|_{x_1=x_M}\overset{\eqref{eq:new73}}=\frac{2\rho\sin{\theta_M}}{\tau_1\bigm|_{x_1=x_M}}\overset{\rho>r/2}>r\frac{\sin{\theta_M}}{\tau_1\bigm|_{x_1=x_M}}\\
        &\overbracket{\ge}^{\mathllap{\tau_1\le2r}}\frac{\sin{\theta_M}}{2}\overset{\sin\theta_M=\sqrt{4r/R}}{\scalebox{8}[1]{$=$}}\frac{2r}{R\sin{\theta_M}}\overset{\sin\theta_1\ge\sin{\theta_M}}\ge\frac{2r}{d_1(x_1)}\overset{r\ge d_0}\ge\frac{2d_0}{d_1}(x_1)\overset{2d_0>\tau_0}>\frac{2(\tau_0-d_0)}{d_1}(x_1)
    \end{aligned}
\end{equation}
for all $x_1\in R_1$ with $\sqrt{4r/R}\le\sin\theta_1\le1/2$. Then $\mathcal{I}<-1$ from \cref{def:AformularForLowerBoundOfExpansion}.
\end{proof}

\subsection{Continuous extension to corners for the derivative of the second iterate of the billiard map}\label{subsec:continuousextension}
\begin{remark}[$D$ is negative on $\mathfrak D^\circ$]\label{rmk:extensionDefs}
\label{def:DF2matrix}
    In the orbit segment $x_0\in \Mrout, \Fc(x_0)=x_1=(\Phi_1,\theta_1)\in \MRin\cap\MRout$, $x_2=\Fc(x_1)\in\Mrin$ corresponding to cases (a) and (b) in \eqref{eqMainCases},  \cite[Equations (4.4), (4.5)]{hoalb} imply
\begin{equation}\label{eq:DF2matrix}
\begin{aligned}
D&=
\begin{bmatrix}
D_{11} & D_{12}\\
D_{21} & D_{22}
\end{bmatrix}
\dfn d_1 d_2D_{x_0}\Fc^2\\
&=
\begin{bmatrix}
2\tau_1(\tau_0-d_0)-d_1(\tau_0-d_0+\tau_1) & 2\tau_1\tau_0-d_1(\tau_0+\tau_1) \\
2(\tau_1-d_2)(\tau_0-d_0)-d_1(\tau_0-d_0+\tau_1-d_2) & 2(\tau_1-d_2)\tau_0-d_1(\tau_0+\tau_1-d_2)
\end{bmatrix}
\end{aligned}
\end{equation}with $d_0$, $d_1$, $\tau_0$, $\tau_1$, $\tau_2$ as in \cref{def:fpclf}.
If $x_1=(\Phi_1,\theta_1)\in\mathfrak{D}^\circ= R_1\cup R_2\cup R_3\cup R_4=\big\{(\Phi_1,\theta_1)\in\MRin\bigm|\sin\theta_1\ge2r/R\big\}\subset\MRin\cap \MRout$  (\cref{fig:R1}), then $D$ is a negative matrix by \cref{caseA}. Now in order to extend the negativity of matrix $D$ to closure$(\mathfrak{D}^\circ)$, we extend $\tau_0,\tau_1,d_0,d_1,d_2$ to $x_1=(\pm\Phistar,\theta_*)$\COMMENT{Is $\theta_*$ ideal notation? (Also: where is it defined?)\wentao{In \cref{def:extensionDefs}.}} on the boundary of $\mathfrak{D}^\circ$ with $p(x_1)=A$ or $B$ at a corner of the table.\COMMENT{\small\red This sounds like it has no purpose. Presumably, before this sentence you want to say that the aim is to extend negativity of $D$ to $\mathfrak D\dfn\overline{\mathfrak{D}}^\circ$. Is that the purpose?\wentao{Yes, some reason clause is now added.}}
\end{remark}

\begin{definition}[Continuous extensions]\label{def:extensionDefs}
The extensions of $\tau_0,\tau_1,d_0,d_1,d_2$ for $p(x_1)$ to the corners of the billiard table are given by the following (with $\theta_*\in [\sin^{-1}{(2r/R)},\pi-\sin^{-1}{(2r/R)}]$).
\par\smallskip\noindent\begin{minipage}{.495\textwidth}
For $p(x_1)=A$ in Fig.~\ref{fig:heuristicproof}, $x_1=(-\Phistar,\theta_*)\in\partial\MRin$  in Fig.~\ref{fig:R1}, 
\begin{equation}\label{eq:lengthextensionA}
\begin{aligned}
& \tau_0=r\sin{(\theta_*-\phistar+\Phistar)}+r|\sin{(\theta_*-\phistar+\Phistar)}|\text{,}\\
& d_0=r|\sin{(\theta_*-\phistar+\Phistar)}|\text{,}\\
& \tau_1=r\sin{(\theta_*+\phistar-\Phistar)}+r|\sin{(\theta_*+\phistar-\Phistar)}|\text{,}\\
& d_2=r|\sin{(\theta_*+\phistar-\Phistar)}|\text{,}\\
& d_1=R\sin{\theta_*}.
\end{aligned}
\end{equation}    
\end{minipage}
%
\hfill\vrule\hfill\begin{minipage}{.495\textwidth}
For $p(x_1)=B$ in Fig.~\ref{fig:heuristicproof}, $x_1=(\Phistar,\theta_*)\in\partial\MRin$ in Fig.~\ref{fig:R1}, 
\begin{equation}\label{eq:lengthextensionB}
\begin{aligned}
& \tau_0=r\sin{(\theta_*+\phistar-\Phistar)}+r|\sin{(\theta_*+\phistar-\Phistar)}|\text{,}\\
& d_0=r|\sin{(\theta_*+\phistar-\Phistar)}|\text{,}\\
& \tau_1=r\sin{(\theta_*-\phistar+\Phistar)}+r|\sin{(\theta_*-\phistar+\Phistar)}|\text{,}\\
& d_2=r|\sin{(\theta_*-\phistar+\Phistar)}|\text{,}\\
& d_1=R\sin{\theta_*}.
\end{aligned}
\end{equation}
\end{minipage}
\end{definition}



\begin{figure}[h]
\begin{center}
    \begin{tikzpicture}[xscale=.5, yscale=.5]
        \pgfmathsetmacro{\PHISTAR}{1.5};
        \pgfmathsetmacro{\Dlt}{0.4};
        \pgfmathsetmacro{\Dltplus}{\Dlt+0.6}
        \pgfmathsetmacro{\WIDTH}{12};
        \tkzDefPoint(-\PHISTAR,0){LL};
        \tkzDefPoint(-\PHISTAR,\WIDTH){LU};
        \tkzDefPoint(\PHISTAR,0){RL};
        \tkzDefPoint(\PHISTAR,\WIDTH){RU};
        \tkzDefPoint(-\PHISTAR,\PHISTAR){X};
        \tkzDefPoint(\PHISTAR,\WIDTH-\PHISTAR){Y};
        \tkzDefPoint(-\PHISTAR,\WIDTH-\PHISTAR){IX};
        \tkzDefPoint(\PHISTAR,\PHISTAR){IY};
        \tkzDefPoint(-\PHISTAR,\PHISTAR+\Dlt){LLD};
        \tkzDefPoint(\PHISTAR,\PHISTAR+\Dlt){RLD};
        \tkzDefPoint(-\PHISTAR,\PHISTAR+\Dlt){LLD};
        \tkzDefPoint(\PHISTAR,\PHISTAR+\Dlt){RLD};
        \tkzDefPoint(-\PHISTAR,\WIDTH-\PHISTAR-\Dlt){LUD};
        \tkzDefPoint(\PHISTAR,\WIDTH-\PHISTAR-\Dlt){RUD};
        \tkzDefPoint(-\PHISTAR,\PHISTAR+\Dltplus){LLDp};
        \tkzDefPoint(\PHISTAR,\PHISTAR+\Dltplus){RLDp};
        \tkzDefPoint(-\PHISTAR,\WIDTH-\PHISTAR-\Dltplus){LUDp};
        \tkzDefPoint(\PHISTAR,\WIDTH-\PHISTAR-\Dltplus){RUDp};
        \tkzDefPoint(-\PHISTAR,\PHISTAR+\Dltplus+0.8){LLDdp};
        \tkzDefPoint(\PHISTAR,\PHISTAR+\Dltplus+0.8){RLDdp};
        \tkzDefPoint(-\PHISTAR,\WIDTH-\PHISTAR-\Dltplus-0.8){LUDdp};
        \tkzDefPoint(\PHISTAR,\WIDTH-\PHISTAR-\Dltplus-0.8){RUDdp};
        \tkzDefPoint (0.3*\PHISTAR,\PHISTAR+\Dltplus){montonestart};
        \tkzDefPoint (0.3*\PHISTAR,0.333*\WIDTH){montoneend};
        \draw [thin] (LL) --(LU);
        \draw [thin] (LU) --(RU);
        \draw [thin] (RL) --(RU);
        \draw [thin] (LL) --(RL);
        \draw [blue,ultra thin,dashed] (LL) --(IY);
        \draw [red,ultra thin,dashed]  (X)--(RL);
        \draw [blue,ultra thin,dashed] (RU) --(IX);
        \draw [red,ultra thin,dashed]  (LU)--(Y);
        \draw [purple,thick]  (LLD)--(RLD);
        \draw [purple,thick]  (LUD)--(RUD);
        \draw [green,ultra thin]  (LLDp)--(RLDp);
        \draw [green,ultra thin]  (LUDp)--(RUDp);
        \draw [green,ultra thin]  (LLDdp)--(RLDdp);
        \draw [green,ultra thin]  (LUDdp)--(RUDdp);
        \tkzLabelPoint[right](RU){$\pi$};
        \tkzLabelPoint[right](RL){$0$};
        \tkzLabelPoint[right](RLD){$\sin{\theta_1}=2r/R$};
        \tkzLabelPoint[right](RUD){$\sin{\theta_1}=2r/R$};
        \tkzLabelPoint[right](RLDp){$\sin{\theta_1}=\sqrt{4r/R}$};
        \tkzLabelPoint[right](RUDp){$\sin{\theta_1}=\sqrt{4r/R}$};
        \tkzLabelPoint[right](RLDdp){$\theta_1=\phistar-\Phistar$};
        \tkzLabelPoint[right](RUDdp){$\theta_1=\pi-\phistar+\Phistar$};
        \coordinate (RM) at ($(RL)!0.333!(RU)$);
        \coordinate (UpperM) at ($(LU)!0.5!(RU)$);
        \coordinate (LowerM) at ($(LL)!0.5!(RL)$);
        \coordinate (LM) at ($(LowerM)!0.333!(UpperM)$);
        \draw [dashed,ultra thin]  (LowerM)--(UpperM);
        \tkzLabelPoint[below](LL){$-\Phistar$};
        \tkzLabelPoint[below](RL){$+\Phistar$};
        \tkzLabelPoint[below](LowerM){$0$};
        \fill [orange, opacity=15/30](LLDp) -- (RLDp) -- (RUDp) -- (LUDp) -- (LLDp) -- cycle;

        \fill [orange, opacity=15/30](LLDp) -- (RLDp) -- (RLD) -- (LLD) -- (LLDp) -- cycle;

        \fill [orange, opacity=15/30](LUDp) -- (RUDp) -- (RUD) -- (LUD) -- (LUDp) -- cycle;
        
        \node[] at (0,0.65*\WIDTH){\small  $(a0)$};
        \node[] at (0,0.35*\WIDTH){\small  $(a0)$};
        \node[] at ($(LLD)!0.5!(RLDp)$){\small  $(a1)$};
        \node[] at ($(RUD)!0.5!(LUDp)$){\small  $(a1)$};

    \end{tikzpicture}
    \caption{Compact region $\mathfrak{D}=\overline{R_1\cup R_2\cup R_3\cup R_4}=\big\{(\Phi_1,\theta_1)\bigm|-\Phistar\le\Phi_1\le\Phistar,\sin{\theta_1}\ge2r/R\big\}$ has $D_{11}<0$, $D_{21}<0$, $D_{12}<0$, $D_{22}<0$. \\ The orange colored region is for ($a0$) case:\, $\big\{(\Phi_1,\theta_1)\mid\sin{\theta_1}\ge\sqrt{4r/R}\big\}$.\\ The other orange colored region is for ($a1$) case:\, $\big\{(\Phi_1,\theta_1)\mid 2r/R\le\sin{\theta_1}<\sqrt{4r/R}\big\}$. }\label{fig:A_compactRegion}\label{fig:A1_region}
\end{center}
\end{figure}

\begin{lemma}\label{lemma:extendedDomainDF2}
Suppose $\phistar\in\big(0,\tan^{-1}(1/3)\big)$ and \eqref{eqJZHypCond} holds.
\begin{itemize}
    \item The extensions in \cref{def:extensionDefs} of $\tau_0$, $\tau_1$, $d_0$, $d_1$, $d_2$ from \cref{def:fpclf} are continuous on the extended domain $\mathfrak{D}\dfn\overline{\mathfrak{D}}^\circ=\,\overline{\!R_1\cup R_2\cup R_3 \cup R_4}=\big\{(\Phi_1,\theta_1)\bigm|-\Phistar\le\Phi_1\le\Phistar,\,\sin^{-1}{(2r/R)}\le\theta_1\le\pi-\sin^{-1}{(2r/R)}\big\}$\COMMENT{\large Correct?\wentao{yes.}} in \cref{fig:A_compactRegion}.\COMMENT{\small\red\bfseries\cref{fig:A_compactRegion} seems to represent $\mathfrak{D}$ incorrectly because it suggests that $\mathfrak{D}$ is the shaded region. The figure should be altered so $\mathfrak{D}$ is next to a vertical brace that spans the right vertical range.\Hrule Also, why is that figure so far from here?\wentao{edited.}}
    \item The matrix elements $D_{11}$, $D_{21}$, $D_{12}$, $D_{22}$ from 
    \eqref{eq:DF2matrix} are negative on $\mathfrak{D}$.\COMMENT{\red\bfseries\cref{fig:A_compactRegion} only mentions that two of them are negative.\wentao{Fixed}}
\end{itemize}
\end{lemma}
\begin{proof} We first check that \eqref{eq:lengthextensionA} and \eqref{eq:lengthextensionB} indeed provide continuous extensions of $\tau_0,\tau_1,d_0,d_1,d_2$.  \cref{lemma:newcoordinateslength} and \eqref{eq:NewCoordinateLegnth} show that 
    $\tau_0$, $\tau_1$, $d_0$, $d_2$ are continuous functions of $\rho,\alpha,\delta$ from \cref{def:TPPORsystem}. And $\rho,\alpha,\delta$ are also continuous functions of $x_1$ on $\Int{\mathfrak{D}}\subset\MRin\cap\MRout$.
    
    As $\Int{\mathfrak{D}}\ni x_1=(\Phi_1,\theta_1)\rightarrow(-\Phistar,\theta_*)\in\partial\MRin$ with $\theta_1\in [\sin^{-1}{(2r/R)},\pi-\sin^{-1}{(2r/R)}]$, by 
\cref{def:TPPORsystem} $\rho\rightarrow r$,  $\alpha\rightarrow -\phistar$, and $\delta\rightarrow -\phistar+\Phistar$. Then by \eqref{eq:NewCoordinateLegnth} and with $d_1$ from \cref{def:fpclf}, we get
\begin{equation}\label{eq:continuousExtension}
    \begin{aligned}
        \tau_1&=\rho\sin{(\theta_1-\delta)}+\sqrt{r^2-\rho^2\cos^2{(\theta_1-\delta)}}\longrightarrow r\sin{(\theta_*+\phistar-\Phistar)}+r|\sin{(\theta_*+\phistar-\Phistar)}|\text{,}\\
        d_2&=\sqrt{r^2-\rho^2\cos^2{(\theta_1-\delta)}}\longrightarrow r|\sin{(\theta_*+\phistar-\Phistar)}|\text{,}\\
        \tau_0&=\rho\sin{(\theta_1+\delta)}+\sqrt{r^2-\rho^2\cos^2{(\theta_1+\delta)}}\longrightarrow r\sin{(\theta_*-\phistar+\Phistar)}+r|\sin{(\theta_*-\phistar+\Phistar)}|\text{,}\\
        d_0&=\sqrt{r^2-\rho^2\cos^2{(\theta_1+\delta)}}\longrightarrow r|\sin{(\theta_*-\phistar+\Phistar)}|\text{,}\\
        d_1&=R\sin\theta_1\longrightarrow R\sin\theta_*
    \end{aligned}
\end{equation}
By the same reasoning, the limits of $\tau_0$, $\tau_1$, $d_0$, $d_0$, $d_1$ are given in \eqref{eq:lengthextensionB} as $\Int{\mathfrak{D}}\ni x_1=(\Phi_1,\theta_1)\rightarrow(\Phistar,\theta_*)\in\partial\MRin$. Therefore, with the extended definitions in \eqref{eq:lengthextensionA} and \eqref{eq:lengthextensionB} $\tau_0$, $\tau_1$, $d_0$, $d_1$, $d_2$ are continuous on $\mathfrak{D}$.

The remainder of the proof is dedicated to showing that $D$ is (entrywise) negative.
By \cref{caseA}, $D_{11}<0$, $D_{21}<0$ for $x_1\in\Big\{x_1=(\Phi_1,\theta_1)\bigm|-\Phistar<\Phi_1<\Phistar,\,\sin^{-1}{(2r/R)}\le\theta_1\le\pi-\sin^{-1}{(2r/R)}\Big\}$. It suffices to verify $D_{11}<0$, $D_{21}<0$ $D_{12}<0$, $D_{22}<0$ for $x_1=(\pm\Phistar,\theta_*)$ with $\sin\theta_*\ge2r/R$.\COMMENT{\Large Why?} We have the following cases.
\newline
\textbf{Case 1:} 
    $x_1=(-\Phistar,\theta_*)$, $p(x_1)=A$, $d_1=R\sin\theta_*\ge2r$ with $\tau_0,\tau_1,d_0,d_1,d_2$ as in \eqref{eq:lengthextensionA}.
\newline
\textbf{Subcase 1.1:}  $\theta_*\in\big(\phistar-\Phistar,\,\pi-\phistar+\Phistar\big)$.
By \eqref{eq:lengthextensionA} we have $\tau_0=2d_0=2r\sin{(\theta_*-\phistar+\Phistar)}>0$ and $\tau_1=2d_2=2r\sin{(\theta_*+\phistar-\Phistar)}>0$, hence
\begin{align*}\frac{2}{d_1}&\overset{d_1\ge2r}\le\frac{1}{r}\overset{0<\tau_0-d_0=d_0\le r}\le\frac{1}{\tau_0-d_0}\overset{\tau_1>0}<\frac{1}{\tau_0-d_0}+\frac{1}{\tau_1},\\
\frac{2}{d_1}&\overset{d_1\ge2r}\le\frac{1}{r}<\frac{2}{r}\overset{d_0\le r, d_2\le r}\le\frac{1}{d_0}+\frac{1}{d_2}\overset{\tau_0=2d_0,\tau_1=2d_2}{\scalebox{7}[1]{$=$}}\frac{1}{\tau_0-d_0}+\frac{1}{\tau_1-d_2},\\
\frac{2}{d_1}&\overset{d_1\ge2r}\le\frac{1}{r}\overset{0<\tau_1-d_2=d_2\le r}\le\frac{1}{\tau_1-d_2}\overset{\tau_0>0}<\frac{1}{\tau_0}+\frac{1}{\tau_1-d_2},\\
\frac{2}{d_1}&\overset{d_1\ge2r}\le\frac{1}{2r}+\frac{1}{2r}<\frac{1}{\tau_0}+\frac{1}{\tau_1},\text{ since }\big(\sin{(\theta_*-\phistar+\Phistar)},\sin{(\theta_*+\phistar-\Phistar)}\big)\neq(1,1).
\end{align*}
This and $d_0,d_1,d_2,\tau_0,\tau_1>0$ give
\begin{equation}\label{eq:DextensionCase1-1}
\begin{gathered}
\begin{multlined}[.85\textwidth]
D_{11}=2\tau_1(\tau_0-d_0)-d_1(\tau_0-d_0+\tau_1)\\
=(\tau_0-d_0)d_1\tau_1\big(\frac{2}{d_1}-\frac{1}{\tau_1}-\frac{1}{\tau_0-d_0}\big)
\overbracket{=}^{\mathclap{\tau_0=2d_0}}\underbracket{d_0d_1\tau_1}_{>0}\big(\overbracket{\frac{2}{d_1}-\frac{1}{\tau_1}-\frac{1}{\tau_0-d_0}}^{<0}\big)<0,
\end{multlined}\\
\begin{multlined}[.85\textwidth]
D_{21}=2(\tau_1-d_2)(\tau_0-d_0)-d_1(\tau_0-d_0+\tau_1-d_2)\\
=(\tau_1-d_2)(\tau_0-d_0)d_1(\frac{2}{d_1}-\frac{1}{\tau_1-d_2}-\frac{1}{\tau_0-d_0})
\overbracket{=}^{\mathclap{\tau_0=2d_0,\tau_1=2d_2}}\underbracket{d_2d_0d_1}_{>0}\big(\overbracket{\frac{2}{d_1}-\frac{1}{\tau_1-d_2}-\frac{1}{\tau_0-d_0}}^{<0}\big)<0,    
\end{multlined}\\
\begin{multlined}[.85\textwidth]
D_{22}=2(\tau_1-d_2)\tau_0-d_1(\tau_0+\tau_1-d_2)\\
=(\tau_1-d_2)d_1\tau_0(\frac{2}{d_1}-\frac{1}{\tau_0}-\frac{1}{\tau_1-d_2})=\overbracket{d_2d_1\tau_0}^{>0}\big(\overbracket{\frac{2}{d_1}-\frac{1}{\tau_0}-\frac{1}{\tau_1-d_2}}^{<0})<0,\\
\end{multlined}\\
\begin{multlined}[.85\textwidth]
D_{12}=2\tau_1\tau_0-d_1(\tau_0+\tau_1)=\overbracket{\tau_0\tau_1d_1}^{>0}\big(\overbracket{\frac{2}{d_1}-\frac{1}{\tau_0}-\frac{1}{\tau_1}}^{<0}\big)<0.
\end{multlined}
\end{gathered}
\end{equation}
\textbf{Subcase 1.2:} $\theta_*\in \big[\sin^{-1}(2r/R),\,\phistar-\Phistar\big]$. By \eqref{eq:lengthextensionA}, we have $\tau_0=0$,
\[
        d_0=r|\sin{(\theta_*-\phistar+\Phistar)}|=r\sin{(-\theta_*+\phistar-\Phistar)}\ge0
        \text{, and  }\tau_1=2d_2=2r\sin{(\theta_*+\phistar-\Phistar)}>0.
\]
Since $0\le -\theta_*+\phistar-\Phistar< \theta_*+\phistar-\Phistar\le2\phistar-2\Phistar<\pi/2$, we have \[d_2=r\sin{(\theta_*+\phistar-\Phistar)}>r\sin{(-\theta_*+\phistar-\Phistar)}=d_0\ge0\]
and
\begin{equation}\label{eq:DextensionCase1-2}
\begin{aligned}
& D_{11}=2\tau_1(\tau_0-d_0)-d_1(\tau_0-d_0+\tau_1)\overbracket=^{\mathclap{\tau_0=0,\ \tau_1=2d_2}}-2\tau_1d_0-d_1(-d_0+2d_2)\overbracket{<}^{\mathclap{d_2>d_0\ge0,\tau_1>0,d_1\ge2r}}0,\\
& D_{21}=2(\tau_1-d_2)(\tau_0-d_0)-d_1(\tau_0-d_0+\tau_1-d_2)\overbracket=^{\mathclap{\tau_0=0,\ \tau_1=2d_2}}-2d_2d_0-d_1(d_2-d_0)\overbracket{<}^{\mathclap{d_2>d_0\ge0,d_2>0,d_1\ge2r}}0.\\
& D_{22}=2\tau_0(\tau_1-d_2)-d_1(\tau_1-d_2+\tau_0)\overbracket=^{\mathclap{\tau_0=0,\ \tau_1=2d_2\qquad}}-d_1d_2\overbracket{<}^{\mathclap{\qquad d_2>d_0\ge0,d_1\ge2r}}0,\\
&D_{12}=2\tau_1\tau_0-d_1(\tau_0+\tau_1)\overbracket=^{\mathclap{\tau_0=0,\ \tau_1=2d_2\qquad}}-2d_1d_2\overbracket{<}^{\mathclap{\qquad d_2>d_0\ge0,d_1\ge2r}}0.
\end{aligned}
\end{equation}
\textbf{Subcase 1.3:} $\theta_*\in \big[\pi-\phistar+\Phistar,\,\pi-\sin^{-1}(2r/R)\big]$. By \eqref{eq:lengthextensionA}, we have 
\begin{align*}
    &\tau_1=0\text{ , }d_2=r|\sin{(\theta_*+\phistar-\Phistar-\pi)}| = r\sin{(\theta_*+\phistar-\Phistar-\pi)}\ge0
    \text{ and }\\&\tau_0=2d_0=2r|\sin{(\theta_*-\phistar+\Phistar)}|=2r\sin{(\pi-\theta_*+\phistar-\Phistar)}>0.
\end{align*}
Since $0\le\theta_*+\phistar-\Phistar-\pi<\pi-\theta_*+\phistar-\Phistar\le2(\phistar-\Phistar)<\pi/2$, we have \[d_0=r\sin{(\pi-\theta_*+\phistar-\Phistar)}>r\sin{(\theta_*+\phistar-\Phistar-\pi)}=d_2\ge0\] 
and
\begin{equation}\label{eq:DextensionCase1-3}
\begin{aligned}
&D_{11}=2\tau_1(\tau_0-d_0)-d_1(\tau_0-d_0+\tau_1)\overbracket=^{\mathclap{\tau_1=0,\ \tau_0=2d_0\qquad}}-d_1d_0\overbracket{<}^{\mathclap{\qquad d_0>0,d_1\ge2r}}0,\\
&D_{21}=2(\tau_1-d_2)(\tau_0-d_0)-d_1(\tau_0-d_0+\tau_1-d_2)\overbracket=^{\mathclap{\tau_1=0,\ \tau_0=2d_0}}-2d_2d_0-d_1(d_0-d_2)\overbracket{<}^{\mathclap{d_0>d_2\ge0,d_1\ge2r}}0.\\
&D_{22}=2(\tau_1-d_2)\tau_0-d_1(\tau_0+\tau_1-d_2)\overbracket=^{\mathclap{\tau_1=0,\ \tau_0=2d_0}}-4d_2d_0-d_1(2d_0-d_2)\overbracket{<}^{\mathclap{d_0>d_2\ge0,d_1\ge2r}}0\\
&D_{12}=2\tau_1\tau_0-d_1(\tau_0+\tau_1)\overbracket=^{\mathclap{\tau_1=0,\ \tau_0=2d_0\qquad}}-2d_0d_1\overbracket{<}^{\mathclap{\qquad d_0>d_2\ge0,d_1\ge2r}}0
\end{aligned}
\end{equation}
\textbf{Case 2} (symmetric to Case 1): $x_1=(\Phistar,\theta_*),p(x_1)=B,d_1=R\sin\theta_*\ge2r$, $\tau_0,\tau_1,d_0,d_1,d_2$ are given in \eqref{eq:lengthextensionB}.
\newline
\textbf{Subcase 2.1:} $\theta_*\in \big(\phistar-\Phistar,\,\pi-\phistar+\Phistar\big)$. By \eqref{eq:lengthextensionB}, we have $\tau_0=2d_0=2r\sin{(\theta_1-\phistar+\Phistar)}>0$ and
$\tau_1=2d_2=2r\sin{(\theta_1+\phistar-\Phistar)}>0$. Using the same reasoning as in \eqref{eq:DextensionCase1-1}, we see that $D_{11}<0$, $D_{21}<0$, $D_{12}<0$, $D_{22}<0$ in this subcase.
\newline
\textbf{Subcase 2.2:} $\theta_*\in [\sin^{-1}(2r/R),\,\phistar-\Phistar]$. By \eqref{eq:lengthextensionB}, we have \[\tau_1=0\text{ , }d_2=r|\sin{(\theta_*-\phistar+\Phistar)}|=r\sin{(-\theta_*+\phistar-\Phistar)}\text{ and }\tau_0=2d_0=2r\sin{(\theta_*+\phistar-\Phistar)}>0\]

Since $0\le -\theta_*+\phistar-\Phistar< \theta_*+\phistar-\Phistar\le2\phistar-2\Phistar<\pi/2$, we have \[d_0=r\sin{(\theta_*+\phistar-\Phistar)}>r\sin{(-\theta_*+\phistar-\Phistar)}\ge0\]

Using the same reasoning as in \eqref{eq:DextensionCase1-3}, we see $D_{11}<0$, $D_{21}<0$, $D_{12}<0$, $D_{22}<0$ in this subcase.
\newline
\textbf{Subcase 2.3:} $\theta_*\in \big[\pi-\phistar+\Phistar,\,\pi-\sin^{-1}(2r/R)\big]$. By \eqref{eq:lengthextensionB}, we have 
\begin{align*}
    &\tau_0=0\text{ , }d_0=r|\sin{(\theta_*+\phistar-\Phistar)}|=r\sin{(\theta_*+\phistar-\Phistar-\pi)}\ge0\\
    \text{ and }&\tau_1=2d_2=2r|\sin{(\theta_*-\phistar+\Phistar)}|=2r\sin{(\pi-\theta_*+\phistar-\Phistar)}>0
\end{align*}
Since $0\le\theta_*+\phistar-\Phistar-\pi<\pi-\theta_*+\phistar-\Phistar\le2(\phistar-\Phistar)<\pi/2$, we have \[d_2=r\sin{(\theta_*+\phistar-\Phistar-\pi)}>d_0=r\sin{(\pi-\theta_*+\phistar-\Phistar)}\ge0\] 
Using the same reasoning as in \eqref{eq:DextensionCase1-2}, we see $D_{11}<0$, $D_{21}<0$, $D_{12}<0$, $D_{22}<0$ in this subcase.

Hence, $D_{11}$, $D_{21}$, $D_{12}$, $D_{22}$ are negative on $\mathfrak{D}$.
\end{proof}



\begin{corollary}\label{corollary:cone-lowerbound}
If $\phistar\in\big(0,\tan^{-1}(1/3)\big)$ and \eqref{eqJZHypCond} holds, then for $D_{11}$, $D_{12}$, $D_{21}$, $D_{12}$ from \cref{def:DF2matrix,def:extensionDefs}, there exist  $\lambda_0(r,R,\phistar)>0$, $\lambda_1(r,R,\phistar)>0$ such that $\frac{D_{21}}{D_{11}}\ge \lambda_0(r,R,\phistar)$, $\frac{D_{22}}{D_{12}}\le \lambda_1(r,R,\phistar)$ for all $x_1$ in the compact region $\mathfrak{D}$ in \cref{lemma:extendedDomainDF2}.
\COMMENT{\large\bfseries\red\cref{fig:A_compactRegion} suggests that $\mathfrak D$ is the shaded region. It is not. I suggest a brace on the left to mark the right height.}
\end{corollary}
\begin{proof}
    By \cref{lemma:newcoordinateslength,lemma:extendedDomainDF2}, $\frac{D_{21}}{D_{11}}>0$ and is a continuous function of $x_1$ on the compact rectangle $\big\{x_1=(\Phi_1,\theta_1)\bigm|-\Phistar\le\Phi_1\le\Phistar,\sin^{-1}{(2r/R)}\le\theta_1\le(\pi-\sin^{-1}{(2r/R)})\big\}$ in \cref{fig:A_compactRegion}. Hence, on the compact rectangle $\frac{D_{21}}{D_{11}}$ attains its minimum  $\lambda_0(r,R,\phistar)>0$  determined by the billiard configuration. Similarly $\frac{D_{22}}{D_{12}}>0$ is a continuous function on $\mathfrak{D}$, so it attains a maximum $\lambda_1(r,R,\phistar)>0$ determined by the billiard configuration.
\end{proof}

\subsection[Narrow cone on \texorpdfstring{$\hat{M}$}{Mhat} and Proof of  \texorpdfstring{\cref{TEcase_a0}}{case a0 expansion}]{Invariance of half-quadrants and expansion in Case (a0)}\label{subsec:narrowcone}
\begin{proposition}[Strictly invariant half-quadrant]\label{corollary:halfquadrantcone}If $0<\phistar<\tan^{-1}{\big(1/3\big)}$ and \eqref{eqJZHypCond} holds, then the half-quadrant $\mathcal{HQ}_x(\mathrm{I,III})\dfn\big\{(d\phi,d\theta)\in T_x M\bigm|\frac{d\theta}{d\phi}\in[0,1]\big\}=C_x\subset T_x\hat{M}$ is strictly invariant under the differential $D\hat{F}$ of the return map (\cref{def:SectionSetsDefs}), i.e., $D_{x}\hat{F}(C_x)\subset\text{\big\{interior of }C_{\hat{F}(x)}\text{\big\}}$.
\end{proposition}
\begin{proof}
Suppose a nonsingular $x\in\hat{M}$ has a return orbit segment containing nonsingular $x_0\in\Mrout$, $x_1=(\Phi_1,\theta_1)\in\MRin$, $x_2\in\Mrin$. By \cref{def:MhatReturnOrbitSegment} we have 
\[
x_0\in\Mrout\xrightarrow{\Fc}x_1=(\Phi_1,\theta_1)\in\MRin\xrightarrow{\Fc}\cdots \xrightarrow{\Fc}\Fc^{n_1}(x_1)\in\MRout\xrightarrow{\Fc}x_2\in\Mrin.
\]
Write 
$\mathcal{Q}_x(\mathrm{I,III})=\big\{\frac{d\theta}{d\phi}\in[0,+\infty]\big\}$.
In cases (a1), (b) and (c) of \eqref{eqMainCases}, \eqref{eq:strictlyinvcone_casea1bc} in \cref{corollary:InvQuadrant} gives 
\[D_x\hat{F}(\mathcal{HQ}_x(\mathrm{I,III}))\subset D_x\hat{F}(\mathcal{Q}_x(\mathrm{I,III}))\overbracket{\subset}^{\text{\eqref{eq:strictlyinvcone_casea1bc}}}\big\{\text{interior of }\mathcal{HQ}_{\hat{F}(x)}(\mathrm{I,III})\big\}.\]

In case (a0) of \eqref{eqMainCases}, since $x\in\hat{M}=(\Mrin\cap\Mrout)\sqcup\Fc^{-1}(\Mrout\setminus\Mrin)$, either $x=x_0$ or $\Fc(x)=x_0$; in the latter case $x$, $x_0\in\Mrout$ are consecutive collisions on $\Gamma_r$, and $D_{x}\Fc=\begin{pmatrix}
1 & 2 \\
0 & 1 
\end{pmatrix}$ in $\phi\theta$-coordinates. So we have two possibilities:
\begin{equation}\label{eq:twoalternativefacts}
\begin{aligned}
    \text{If } x=x_0, &\text{ then } \mathcal{HQ}_x(\mathrm{I,III})=\mathcal{HQ}_{x_0}(\mathrm{I,III}).\\
    \text{If } x=\Fc^{-1}(x_0)\quad(x_0=\Fc(x)), &\text{ then } D_x\Fc(\mathcal{HQ}_x(\mathrm{I,III}))\subset\mathcal{HQ}_{x_0}(\mathrm{I,III}).
\end{aligned}
\end{equation}
Then also in (a0) there are two subcases based on nonsingular $x_2\in\Mrin$.
\newline
\textbf{Subcase (i):} $m(x_2)\le1$, i.e., $x_2\in \textcolor{blue}{M^{\text{\upshape{in}}}_{r,0}}\cup\textcolor{blue}{M^{\text{\upshape{in}}}_{r,1}}\subset\hat{M}$. In this subcase, $\hat{F}(x)=x_2\in\hat{M}$. From \cref{def:DF2matrix} and \eqref{eq:DF2matrix}, 
\begin{equation}\label{eq:narrowconereason1}
    \begin{aligned}
        D_{22}-D_{12}=d_2d_1-2\tau_0d_2=d_2(d_1-2\tau_0)&\underbracket{>}_{\mathclap{\tau_0<2r}}d_2(\overbracket{d_1}^{\mathclap{\substack{\text{By \cref{lemma:returnorbittheta1Large,Corollary:CollisionOnMrin0Mrin1}}\\\sin\theta_1>1/2 \text{ and }d_1=R\sin\theta_1>R/2>850r}}}-4r)>0. 
    \end{aligned}
\end{equation} Together with $D_{22}<0$, $D_{12}<0$, \eqref{eq:narrowconereason1} indicates that $\frac{D_{22}}{D_{12}}\in(0,1)$. And the billiard map (and its iteration $D_{x_0}\Fc^2$) maintains orientation (\cite[equation (2.29)]{cb}) which indicates that $D_{22}D_{11}-D_{12}D_{21}>0$. Since by \cref{caseA}, $D_{22}$, $D_{11}$, $D_{12}$, $D_{21}$ are all negative, we get the following.
\begin{equation}\label{eq:narrowconereason2}
    \frac{D_{22}}{D_{12}}-\frac{D_{21}}{D_{11}}=\frac{D_{22}D_{11}-D_{12}D_{21}}{D_{12}D_{11}}>0.
\end{equation}
Hence $0<\frac{D_{21}}{D_{11}}<\frac{D_{22}}{D_{12}}<1$.\newline For $dx_0=(d\phi_0,d\theta_0)\in\mathcal{Q}_{x_0}(\mathrm{I,III})$, that is, $\frac{d\theta_0}{d\phi_0}=k$ for some $k\in[0,\infty]$, $D_{x_0}\Fc^2(dx_0)=dx_2=(d\phi_2,d\theta_2)$ satisfies 
\[
\frac{D_{21}+kD_{22}}{D_{11}+kD_{12}}-\frac{D_{21}}{D_{11}}=\frac{\overbracket{k}^{\mathclap{\ge0}}(\overbracket{D_{11}D_{22}-D_{21}D_{12}}^{\mathrlap{>0}})}{\underbracket{D_{11}(D_{11}+kD_{12})}_{>0}}\ge0>\frac{(\overbracket{D_{11}D_{22}-D_{21}D_{12}}^{\mathrlap{>0}})}{\underbracket{D_{12}(D_{11}+kD_{12})}_{>0}}=
    \frac{D_{22}}{D_{12}}-\frac{D_{21}+kD_{22}}{D_{11}+kD_{12}}\rlap{, so}
\]
\begin{equation}\label{eq:narrowconereason3}
    0<\frac{D_{21}}{D_{11}}\le\frac{d\theta_2}{d\phi_2}=\frac{D_{21}+kD_{22}}{D_{11}+kD_{12}}<\frac{D_{22}}{D_{12}}<1,
\end{equation}
and $D_{x_0}\Fc^2(\mathcal{Q}_{x_0}(\mathrm{I,III}))\subset\big\{\text{interior of }\mathcal{HQ}_{x_2}(\mathrm{I,III})\big\}$. We get two conclusions:\newline If $x=x_0$, then 
\begin{align*}
D_x\hat{F}(\mathcal{HQ}_x(\mathrm{I,III}))&\overbracket{=}^{\mathclap{x=x_0,\Fc^2(x)=x_2=\hat{F}(x)}}D_x\Fc^2(\mathcal{HQ}_x(\mathrm{I,III}))\subset D_{x_0}\Fc^2(\mathcal{Q}_{x_0}(\mathrm{I,III}))\\
&\subset\big\{\text{interior of }\mathcal{HQ}_{x_2}(\mathrm{I,III})\big\}\overbracket=^{\mathclap{\hat{F}(x)=x_2}}\big\{\text{interior of }\mathcal{HQ}_{\hat{F}(x)}(\mathrm{I,III})\big\}.\\
\end{align*}\newline If $x=\Fc^{-1}(x_0)$, then
\begin{align*}
D_x\hat{F}(\mathcal{HQ}_x(\mathrm{I,III}))&\overbracket{=}^{\mathclap{\Fc(x)=x_0,\Fc^2(x_0)=x_2=\hat{F}(x)\quad}}D_x\Fc(\mathcal{HQ}_x(\mathrm{I,III}))^3\overbracket{=}^{\mathclap{\quad\text{chain rule and }x_0=\Fc(x)}}D_{x_0}\Fc^2(D_x\Fc(\mathcal{HQ}_x(\mathrm{I,III})))\overbracket{\subset}^{\mathclap{\text{\eqref{eq:twoalternativefacts}}}} D_{x_0}\Fc^2(\mathcal{HQ}_{x_0}(\mathrm{I,III}))\\
&\subset\big\{\text{interior of }\mathcal{HQ}_{x_2}(\mathrm{I,III})\big\}\overbracket=^{\mathclap{\hat{F}(x)=x_2}}\big\{\text{interior of }\mathcal{HQ}_{\hat{F}(x)}(\mathrm{I,III})\big\}.
\end{align*}

\textbf{Subcase (ii):} $m(x_2)\ge2$, i.e., $\Fc^{m(x_2)-1}(x_2)=\hat{F}(x)$. Here, $D_{x_2}\Fc^{m(x_2)-1}=\begin{pmatrix}
1 & 2(m(x_2)-1) \\
0 & 1 
\end{pmatrix}$ satisfies \[D_{x_2}\Fc^{m(x_2)-1}\big(\big\{\text{interior of }\mathcal{Q}_{x_2}(\mathrm{I,III})\big\}\big)\subset\big\{\text{interior of }\mathcal{HQ}_{\hat{F}(x)}(\mathrm{I,III})\big\}.\]
If $x=x_0$, then $x_2=\Fc^2(x)$, $\Fc^{m(x_2)+1}(x)=\hat{F}(x)$ and 
\begin{align*}
D_x\hat{F}(\mathcal{HQ}_x(\mathrm{I,III}))&\overbracket{=}^{\mathclap{\Fc^{m(x_2)+1}(x)=\hat{F}(x)}}D_x\Fc^{m(x_2)+1}(\mathcal{HQ}_x(\mathrm{I,III}))\overbracket{\subset}^{\mathclap{x=x_0,\Fc^2(x)=x_2}} D_{x_2}\Fc^{m(x_2)-1}(D_{x_0}\Fc^2(\mathcal{Q}_{x_0}(\mathrm{I,III})))\\
&\overbracket{\subset}^{\mathclap{D_{x_0}\Fc^2\text{ is a matrix with negative entries}}}D_{x_2}\Fc^{m(x_2)-1}\Big(\big\{\text{interior of }\mathcal{Q}_{x_2}(\mathrm{I,III})\big\}\Big)\subset\big\{\text{interior of }\mathcal{HQ}_{\hat{F}(x)}(\mathrm{I,III})\big\}.\\
\end{align*}
If $x=\Fc^{-1}(x_0)$, then $x_2=\Fc^3(x)$,  $\Fc^{m(x_2)+2}(x)=\hat{F}(x)$ and
\begin{align*}
D_x\hat{F}(\mathcal{HQ}_x(\mathrm{I,III}))&\overbracket{=}^{\mathclap{\Fc^{m(x_2)+2}=\hat{F}(x)}}D_x\Fc^{m(x_2)+2}(\mathcal{HQ}_x(\mathrm{I,III}))\overbracket{\subset}^{\mathclap{\Fc(x)=x_0,\Fc^2(x_0)=x_2\text{ and chain rule}}} D_{x_2}\Fc^{m(x_2)-1}(D_{x_0}\Fc^2(D_x\Fc(\mathcal{Q}_{x}(\mathrm{I,III}))))\overbracket{\subset}^{\mathclap{\text{\eqref{eq:twoalternativefacts}}}} D_{x_2}\Fc^{m(x_2)-1}(D_{x_0}\Fc^2(\mathcal{Q}_{x_0}(\mathrm{I,III})))\\
&\overbracket{\subset}^{\mathclap{D_{x_0}\Fc^2\text{ is a matrix with negative entries}}}D_{x_2}\Fc^{m(x_2)-1}\Big(\big\{\text{interior of }\mathcal{Q}_{x_2}(\mathrm{I,III})\big\}\Big)\subset\big\{\text{interior of }\mathcal{HQ}_{\hat{F}(x)}(\mathrm{I,III})\big\}.\\
\end{align*}In both \textbf{subcases (i) and (ii)}, $D_x\hat{F}(\mathcal{HQ}_x(\mathrm{I,III}))\subset\big\{\text{interior of }\mathcal{HQ}_{\hat{F}(x)}(\mathrm{I,III})\big\}$.
\end{proof}

\begin{theorem}[Case (a0) \cref{TEcase_a0}\eqref{itemTEcase_a0-1}, \eqref{itemTEcase_a0-2}]\label{theorem:UniformExpansionA0}
Suppose a nonsingular $x\in \hat{M}$ (\cref{def:SectionSetsDefs}) has a return orbit segment $x,\Fc(x),\cdots,\Fc^{\sigma(x)}(x)=\hat{F}(x)\in \hat{M}$ as in \cref{def:MhatReturnOrbitSegment} with $x_{0}\in\Mrout$, $x_1=\Fc(x_0)=(\phi_1,\theta_1)\in \MRin$, $x_2\in\Mrin$ and $\sin\theta_1\ge \sqrt{4r/R}$ (case (a0) of \eqref{eqMainCases}; \cref{TEcase_a0}).

Then there exists a $\lambda_c(r,R,\phistar)>0$ such that for all $x\in\hat{M}=(\Mrin\cap\Mrout)\sqcup\Fc^{-1}(\Mrout\setminus\Mrin)$  (\cref{def:SectionSetsDefs}) and any tangent vector $dx=(d\phi,d\theta)\in C_x=\mathcal{HQ}_{x}=\Big\{(d\phi,d\theta)\bigm|\frac{d\theta}{d\phi}\in\big[0,1\big]\Big\}$ at $x$ (the half-quadrant from \cref{corollary:halfquadrantcone}) we have 
\begin{itemize} 
    \item $D_x\hat{F}(dx)\in\big\{\text{interior of }\mathcal{HQ}_{\hat{F}(x)}\big\}$,
    \item $\frac{\|D_x\hat{F}(dx)\|}{\|dx\|_\p}>1+\lambda_c$.
\end{itemize}
\end{theorem} These two conclusions are exactly the statements of \cref{TEcase_a0}\eqref{itemTEcase_a0-1}, \eqref{itemTEcase_a0-2}.

\begin{proof}
By \cref{corollary:halfquadrantcone}, $D_x\hat{F}(dx)\in\big\{\text{interior of }\mathcal{HQ}_{\hat{F}(x)}\big\}$. And by \cref{lemma:ExpansionLowerBoundA0A1}, we have
\begin{itemize} 
    \item If $x\in\Mrin\cap\Mrout$, then $dx_0\dfn dx$ and $(d\phi_2,d\theta_2)=dx_2\dfn D_{x}\Fc^2(dx)$. Otherwise, $x\in\Fc^{-1}(\Mrout\setminus\Mrin)$, then $dx_0\dfn D_x\Fc(dx)$ and $(d\phi_2,d\theta_2)=dx_2\dfn D_x\Fc^3(dx)$. In either case, $dx_0$ is a tangent vector at $x_0$ and $dx_2$ is a tangent vector at $x_2$.
    \item $\frac{\|dx_0\|_\p}{\|dx\|_\p}\ge 1$ and $\frac{\|D_x\hat{F}(dx)\|_\p}{\|dx_2\|_\p}\ge 1$.
    \item If $\mathcal{I}<0$ (see \cref{def:AformularForLowerBoundOfExpansion}), then $\frac{\|dx_2\|_\p}{\|dx_0\|_\p}\ge|\mathcal{I}|$.
\end{itemize}

Then we have two subcases on $x_2$.
\newline
\textbf{Case i:} $m(x_2)\le1$, i.e., $x_2\in\textcolor{blue}{M^{\text{\upshape{in}}}_{r,0}}\sqcup\textcolor{blue}{M^{\text{\upshape{in}}}_{r,1}}$. \cref{Corollary:CollisionOnMrin0Mrin1} gives $\mathcal{I}<-1.37$. So, $\frac{\|D_x\hat{F}(dx)\|}{\|dx\|_\p}=\frac{\|dx_0\|_\p}{\|dx\|_\p}\frac{\|dx_2\|_\p}{\|dx_0\|_\p}\frac{\|D_x\hat{F}(dx)\|_\p}{\|dx_2\|_\p}>1.37$.
\newline
\textbf{Case ii:}  $m(x_2)\ge2$. \cref{Prop:Ilowerbound} gives $\mathcal{I}<-1$. Since the billiard map (and its iteration $D_{x_0}\Fc^2$) maintains orientation (\cite[equation (2.29)]{cb}), by \cref{corollary:cone-lowerbound}, $\frac{d\theta_2}{d\phi_2}\in [\frac{D_{21}}{D_{11}}(x_1),\frac{D_{22}}{D_{12}}(x_1)]\subset[\lambda_0(r,R,\phistar),\infty)$ for all $x_1=(\Phi_1,\theta_1)$ with $\sin\theta_1\ge\sqrt{4r/R}>2r/R$. Therefore, by \cref{corollary:shearlowerboundexpansion}, we have $\frac{\|D_x\hat{F}(dx)\|_\p}{\|dx_2\|_\p}\ge 1+\frac{c_r\lambda_0(r,R,\phistar)}{2}$. We also get $\frac{\|D_x\hat{F}(dx)\|_\p}{\|dx\|_\p}=\frac{\|dx_0\|_\p}{\|dx\|_\p}\frac{\|dx_2\|_\p}{\|dx_0\|_\p}\frac{\|D_x\hat{F}(dx)\|_\p}{\|dx_2\|_\p}>1+\frac{c_r\lambda_0(r,R,\phistar)}{2}$. 

Hence, we can choose $\lambda_c=\min{\big\{0.37,\frac{c_r\lambda_0(r,R,\phistar)}{2}\big\}}$, which depends on $r$, $R$, $\phistar$.
\end{proof}
\subsection[Proof of \texorpdfstring{\cref{TEcase_a1}}{case (a1) contraction control}]{Contraction control in Case (a1)}\label{subsec:ProofTECase_a1}

\begin{theorem}[Contraction control in Case (a1) of \eqref{eqMainCases} ]\label{theorem:contractioncontrolA1}
Suppose $x\in \hat{M}=(\Mrin\cap\Mrout)\sqcup\Fc^{-1}(\Mrout\setminus\Mrin)$ (\cref{def:SectionSetsDefs}) is not singular 
with $2r/R\le\sin\theta_1< \sqrt{4r/R}$ in \eqref{eqPtsFromeqMhatOrbSeg} (see \cref{fig:A1_region}), 
$dx=(d\phi,d\theta)\in C_x=\mathcal{HQ}_{x}=\Big\{(d\phi,d\theta)\bigm|\frac{d\theta}{d\phi}\in\big[0,1\big]\Big\}$ is a tangent vector at $x$ in the half-quadrant  (\cref{corollary:halfquadrantcone}). Then \cref{TEcase_a1}\eqref{itemTEcase_a1-1},\eqref{itemTEcase_a1-2} hold:
\begin{itemize} 
    \item $D_x\hat{F}(dx)\in\big\{\text{interior of }\mathcal{HQ}_{\hat{F}(x)}\big\}$,
    \item $\frac{\|D_x\hat{F}(dx)\|}{\|dx\|_\p}>0.26$.
\end{itemize}
\end{theorem} 
\begin{proof}
By \cref{corollary:halfquadrantcone}, $D_x\hat{F}(dx)\in\big\{\text{interior of }\mathcal{HQ}_{\hat{F}(x)}\big\}$.

By \cref{contraction_region}, we have $x_0\in\Nout$ and $x_2\in\Nin$. Then by \cref{def:SectionSetsDefs}, we have $x\in\Fc^{-1}(\Mrin\setminus\Mrout)$, $\Fc(x)=x_0$. Otherwise, $x=x_0\in\Mrin\cap\Mrout$ and $x_0\in\textcolor{red}{M^{\text{\upshape{out}}}_{r,0}}$, which contradicts \cref{Prop:Udisjoint} that $\Nout\cap\textcolor{red}{M^{\text{\upshape{out}}}_{r,0}}=\emptyset$. By \cref{lemma:ExpansionLowerBoundA0A1}, we have
\begin{itemize} 
    \item $x\in\Fc^{-1}(\Mrout\setminus\Mrin)$, $dx_0\dfn D_x\Fc(dx)$ and $(d\phi_2,d\theta_2)=dx_2\dfn D_x\Fc^3(dx)$. $dx_0$ is a tangent vector at $x_0$ and $dx_2$ is a tangent vector at $x_2$.
    \item $\frac{\|dx_0\|_\p}{\|dx\|_\p}\ge 1$ and $\frac{\|D_x\hat{F}(dx)\|_\p}{\|dx_2\|_\p}\ge 1$.
    \item If in \cref{def:AformularForLowerBoundOfExpansion} it has $\mathcal{I}<0$, then $\frac{\|dx_2\|_\p}{\|dx_0\|_\p}\ge|\mathcal{I}|$.
\end{itemize}
We assume $\theta_1\in(0,\pi/2]$, then $\theta_1\in\big[\sin^{-1}{(2r/R)},\sin^{-1}{\sqrt{4r/R}}\big)$. We have two subcases for $\Phi_1$.

\noindent\textbf{Subcase i}, $\Phi_1\in(-\Phistar,0)$, then $x_1\in R_2$ region of \cref{fig:R1}. By \cref{prop:casea0smallthetaR2}, we have $\mathcal{I}\le -1$.

\noindent\textbf{Subcase ii}, $\Phi_1\in[0,\Phistar)$, then $x_1\in R_1$ region of \cref{fig:R1}. Note that from \cref{def:AformularForLowerBoundOfExpansion}, we have $\tau_0+\tau_1>2d_0$, $\tau_1>0$, $d_0>0$ thus $\mathcal{I}<-1+\frac{\tau_1}{d_0}\cdot\frac{2(\tau_0-d_0)}{d_1}$.

For a fixed $\theta_1\in\Big[\sin^{-1}{(2r/R)},\sin^{-1}{\big(\sqrt{4r/R}\big)}\Big)$, as $\Phi_1$ varies from $0$ to $\Phistar$, $x_1=(\Phi_1,\theta_1)$ is variable depending on $\Phi_1\in[0,\Phistar)$, $\tau_1$, $d_0$, $d_2$ are functions of $x_1$ thus of $\Phi_1$. 

By \cref{monotone_p4}, $d_0$ as functions of $x_1$ is a monotone increasing function of $\Phi_1$. 

In \cref{subsec:ProofContractionRegion,rmk:CaseA1Bcoordinates}, we have $(r\sin{\phi_2},-r\cos\phi_2)=(x_Q,y_Q)=Q=p(x_2)$ by \eqref{eq:newXqXt}. By \eqref{eq:observation2(a1)(b)(c)region} in \cref{proposition:NinsizeInA1BCcases}\eqref{item01NinsizeInA1BCcases}, in case (a1) we further get the following.

\[0<\phi_2-\phistar\overbracket{<}^{\strut\mathclap{\text{\cref{proposition:NinsizeInA1BCcases}\eqref{item01NinsizeInA1BCcases},\eqref{eq:observation2(a1)(b)(c)region}}}}2.0001\sqrt{4r/R}\overbracket{<}^{\mathrlap{\text{\eqref{eqJZHypCond0}: }R>1700r}}0.0971.\]


Since $\phi_2\in(\phistar,2\pi-\phistar)$, $0<\phistar<\tan^{-1}{\big(1/3\big)}$, we have $\phistar<\phi_2<\phistar+0.0971<\pi/2$. Thus, $x_Q>0$, by \cref{monotone_p3}(1), have $\frac{d\tau_1}{d\Phi_1}=-\frac{bx_Q}{d_2}<0$, so $\tau_1$ is monotonically decreasing with respect to $\Phi_1$. 

By symmetry, $d_2\bigm|_{x_1=(0,\theta_1)}=d_0\bigm|_{x_1=(0,\theta_1)}$. Hence we have 
\begin{equation}\label{eq:a1Ibound1}
0<\frac{\tau_1}{d_0}\le \frac{\tau_1}{d_0}\bigm|_{x_1=(\Phi_1,\theta_1)}\overbracket{\le}^{\mathclap{\text{\cref{monotone_p3,monotone_p4}}}}\frac{\tau_1}{d_0}\bigm|_{x_1=(0,\theta_1)}=\frac{\tau_1}{d_2}\bigm|_{x_1=(0,\theta_1)}\overbracket{<}^{\mathclap{\text{\cref{remark:deffpclf}}}}2
\end{equation}

In \cref{fig:heuristicproof}, $P=p(x_1)$, let $M$ be the midpoint of arc segment $\arc{AB}$, we suppose $p(x_1)=P$ is on arc segment $MB$ with $\theta_1$ as the angle between the $\overrightarrow{PQ}$ and the tangent direction of $\Gamma_R$ at $P$. \COMMENT{This "moves" terminology is not helpful, but rampant. Over time, these occurrences should be rewritten. Low priority.\wentao{rewritten now.}}. Using normal coordinates at $P$ (\cref{def:TPPORsystem}), $\rho$, $\delta$ the length and angle variables of \cref{def:TPPORsystem,lemma:newcoordinateslength,fig:fractionEstimate2}, with $\Phi_1\in[0,\Phistar)$,\COMMENT{This sentence is {\Large incomprehensible}.\wentao{Much edited}} and with \textbf{Subcase ii}: $\Phi_1\in[0,\Phistar)$, $0\le\delta<\phistar$ we get the following.\COMMENT{The next display needs fixing\wentao{Much edited and issue addressed.}}
\begin{equation}\label{eq:a1Ibound2}
\begin{multlined}
     \frac{2(\tau_0-d_0)}{d_1}\overbracket=^{{\eqref{eq:NewCoordinateLegnth}}}\frac{2\rho\sin{(\theta_1+\delta)}}{d_1}\overbracket\le^{\mathrlap{\substack{0<\rho\le r, 2r\le d_1}}}\sin{(\theta_1+\delta)}<|\sin{\theta_1}|+|\sin{\delta}|<\sqrt{4r/R}+\sin{\phistar}\overbracket<^{\mathclap{R>1700r,0<\phistar<\tan^{-1}{(1/3)}}}0.37,\\
    \mathcal{I}=-1+\frac{\tau_1}{d_0}\frac{2(\tau_1-d_0)}{d_1}-\overbracket{\frac{\tau_0+\tau_1-2d_0}{d_0}}^{\mathllap{\text{\cref{remark:symmetrywiththetangentline}:}>0}}<-1+\overbracket{\frac{\tau_1}{d_0}}^{\strut\mathllap{\text{\eqref{eq:a1Ibound1}: }\in(0,2)}}\overbracket{\frac{2(\tau_1-d_0)}{d_1}}^{\mathrlap{<0.37}}<-1+2\times0.37=-0.26.
\end{multlined}
\end{equation}    

Therefore, $\mathcal{I}<-0.26$ in either \textbf{Subcase i} or \textbf{Subcase ii}. By \cref{lemma:ExpansionLowerBoundA0A1} we have $\frac{\|D_x\hat{F}(dx)\|_\p}{\|dx\|_\p}=\frac{\|D_x\hat{F}(dx)\|_\p}{\|dx_2\|_\p}\cdot\frac{\|dx_2\|_\p}{\|dx_0\|_\p}\cdot\frac{\|dx_0\|_\p}{\|dx\|_\p}>0.26$.\end{proof}

\section[Proof of \texorpdfstring{\cref{TEcase_b}}{case (b)}]{Contraction control in case (b)}\label{sec:ProofTECase_b}


In this section, we prove \cref{TEcase_b}\eqref{itemTEcase_b2}, 
 which covers Case (b) of \eqref{eqMainCases}. (This is \cref{subsec:ProofTECase_b};  \cref{TEcase_b}\eqref{itemTEcase_b1} is the conclusion of \cref{corollary:halfquadrantcone}.)
 Accordingly, all statements in this section are made in the context of \cref{TEcase_b}, which we now recall for the convenience of the reader. 

\hidesubsection{Standing assumptions for this section}
\begin{remark}\label{remark:standingassumptionInsection6}
We consider a nonsingular $x\in \hat{M}=(\Mrin\cap\Mrout)\sqcup\Fc^{-1}(\Mrout\setminus\Mrin)$ (\cref{def:SectionSetsDefs}) with return orbit segment $x,\Fc(x),\cdots,\Fc^{\sigma(x)}(x)=\hat{F}(x)\in \hat{M}$ as in \cref{def:MhatReturnOrbitSegment}, where $x_{0}\in\Mrout$, $x_1=\Fc(x_0)=(\phi_1,\theta_1)\in \MRin$, $x_2\in\Mrin$, $\sin\theta_1< 2r/R$, $n_1=0$ (i.e., Case (b) of \eqref{eqMainCases}),
and the tangent vector $dx=(d\phi,d\theta)\in\Big\{(d\phi,d\theta)\bigm|\frac{d\theta}{d\phi}\in\big[0,1\big]\Big\}$ at $x$ is in the half-quadrant cone in $\phi\theta$-coordinates. \end{remark}

\cref{lemma:ExpansionLowerBoundB} provides more information about this orbit segment.

\begin{theorem}[\cref{TEcase_b}\eqref{itemTEcase_b2}]\label{theorem:proofofTEcase_b}\label{subsec:ProofTECase_b}
    Under the standing assumptions for this section (\cref{remark:standingassumptionInsection6}), \[\dfrac{\|D_x\hat{F}(dx)\|_\p}{\|dx\|_\p}>0.05\].
\end{theorem}
\begin{proof}
From facts below, 
    $\displaystyle\frac{\|D_x\hat{F}(dx)\|_\p}{\|dx\|_\p}=\underbracket{\frac{\|D_x\hat{F}(dx)\|_\p}{\|dx_3\|_\p}}_{\mathclap{>1\text{ (\cref{lemma:ExpansionLowerBoundB}\eqref{itemExpansionLowerBoundB2})}}}\overbracket{\frac{\|dx_3\|_\p}{\|dx_2\|_\p}}^{\mathclap{\ge1 \text{ (\cref{proposition:dx3dx2expansion})}}}\underbracket{\frac{\|dx_2\|_\p}{\|dx_0\|_\p}}_{\mathclap{=|\mathcal{I}(k)|\text{ (\cref{lemma:ExpansionLowerBoundB}\eqref{itemExpansionLowerBoundB3})}}}\overbracket{\frac{\|dx_0\|_\p}{\|dx\|_\p}}^{\mathclap{\ge1\text{ (\cref{lemma:ExpansionLowerBoundB}\eqref{itemExpansionLowerBoundB1})}}}>|\mathcal{I}(k)|\underbracket>_{\mathrlap{\text{\cref{proposition:IktransionInCaseb}}}}0.05$.
\end{proof}
This provides the estimate in  \cref{TEcase_b}. The rest of this section  proves the ingredients just invoked. Readers wishing to skip the proofs of these facts can advance to page \pageref{eq:casebcase2minorterm4}.
\subsection{Lemmas for \texorpdfstring{\cref{TEcase_b}}{Case b}}\label{subsec:AFELBincaseb}
\begin{lemma}[Orbit configuration and contraction bound]\label{lemma:ExpansionLowerBoundB}
Under the standing assumptions for this section (\cref{remark:standingassumptionInsection6}), we have the following.
\begin{enumerate}
    \item\label{itemExpansionLowerBoundB1} $x\in\Fc^{-1}(\Mrout\setminus\Mrin)\subset M_r$, $x_0=\Fc(x)\in\Mrout$, $dx_0\dfn D_x\Fc(dx)$ and $dx_2\dfn D_{x}\Fc^3(dx)$, where $dx_0$ is a tangent vector at $x_0$ and $dx_2$ is a tangent vector at $x_2$. We have $\frac{\|dx_0\|_\p}{\|dx\|_\p}\ge1$.
    \item\label{itemExpansionLowerBoundB2} $m(x_2)\ge3$, $x_3\dfn\Fc(x_2)\in M_r$, $dx_3\dfn D_x\Fc^4(dx)$ and  $\Fc^{m(x_2)-2}(x_3)=\Fc^{m(x_2)-1}(x_2)\overset{\text{\cref{def:SectionSetsDefs}}}{\scalebox{7}[1]{$=$}}\hat{F}(x)\in\hat{M}$, i.e. $x_3=\Fc(x_2)$, $\cdots$, $\Fc^{m(x_2)-1}(x_2)$ are all in $M_r\setminus\Mrout$. We have $\frac{\|D_x\hat{F}(dx)\|_\p}{\|dx_3\|_\p}>1$.
    \item\label{itemExpansionLowerBoundB3} We use the same notation as in \eqref{eq:MirrorEquationExpansionFormula} with index subscript denoting by $\mathcal{B}^{\pm}_0,\mathcal{B}^{\pm}_1$ the infinitesimal wave front curvatures after/before collisions at $x_0$, $x_1$, respectively.
    We note that $\Bcal^+_0\in\big[-\frac{-4}{3d_0},\frac{-1}{d_0}\big]$.
    
    Then $\frac{\|dx_2\|_\p}{\|dx_0\|_\p}=\big|\mathcal{I}(k)\big|$ where $\mathcal{I}(k)\dfn(\tau_0+\tau_1-\frac{2\tau_0\tau_1}{d_1})(\frac{-k}{d_0})+1-\frac{2\tau_1}{d_1}$ and $k=1+\frac{1}{2+\frac{1}{{d\theta}/{d\phi}}}\in[1,\frac{4}{3}]$ (since $\frac{d\theta}{d\phi}\in[0,1]$), and $\tau_0$, $\tau_1$, $d_0$, $d_1$, $d_2$ are the length functions (\cref{def:fpclf}) of $x_1\in\big\{(\Phi_1,\theta_1)\in\MRin\cap\MRout\bigm|\sin\theta_1<2r/R\big\}=\mathcal{L}\sqcup\mathcal{U}$ in \cref{fig:A1contractionCase}. 
    Hence the expansion is determined by the cone vector $dx$ and length functions at $x_1$.
\end{enumerate}
\end{lemma}
\begin{figure}[ht]
\begin{center}
    \begin{tikzpicture}[xscale=0.6, yscale=0.6]
        \pgfmathsetmacro{\PHISTAR}{1.5};
        \pgfmathsetmacro{\Dlt}{0.6};
        \pgfmathsetmacro{\Dltplus}{\Dlt+0.6}
        \pgfmathsetmacro{\WIDTH}{12};
        \tkzDefPoint(-\PHISTAR,0){LL};
        \tkzDefPoint(-\PHISTAR,\WIDTH){LU};
        \tkzDefPoint(\PHISTAR,0){RL};
        \tkzDefPoint(\PHISTAR,\WIDTH){RU};
        \tkzDefPoint(-\PHISTAR,\PHISTAR){X};
        \tkzDefPoint(\PHISTAR,\WIDTH-\PHISTAR){Y};
        \tkzDefPoint(-\PHISTAR,\WIDTH-\PHISTAR){IX};
        \tkzDefPoint(\PHISTAR,\PHISTAR){IY};
        \tkzDefPoint(-\PHISTAR,\PHISTAR+\Dlt){LLD};
        \tkzDefPoint(\PHISTAR,\PHISTAR+\Dlt){RLD};
        \tkzDefPoint(-\PHISTAR,\PHISTAR+\Dlt){LLD};
        \tkzDefPoint(\PHISTAR,\PHISTAR+\Dlt){RLD};
        \tkzDefPoint(-\PHISTAR,\WIDTH-\PHISTAR-\Dlt){LUD};
        \tkzDefPoint(\PHISTAR,\WIDTH-\PHISTAR-\Dlt){RUD};
        \tkzDefPoint(-\PHISTAR,\PHISTAR+\Dltplus){LLDp};
        \tkzDefPoint(\PHISTAR,\PHISTAR+\Dltplus){RLDp};
        \tkzDefPoint(-\PHISTAR,\WIDTH-\PHISTAR-\Dltplus){LUDp};
        \tkzDefPoint(\PHISTAR,\WIDTH-\PHISTAR-\Dltplus){RUDp};
        \tkzDefPoint(-\PHISTAR,\PHISTAR+\Dltplus+0.8){LLDdp};
        \tkzDefPoint(\PHISTAR,\PHISTAR+\Dltplus+0.8){RLDdp};
        \tkzDefPoint(-\PHISTAR,\WIDTH-\PHISTAR-\Dltplus-0.8){LUDdp};
        \tkzDefPoint(\PHISTAR,\WIDTH-\PHISTAR-\Dltplus-0.8){RUDdp};
        \tkzDefPoint (0.3*\PHISTAR,\PHISTAR+\Dltplus){montonestart};
        \tkzDefPoint (0.3*\PHISTAR,0.333*\WIDTH){montoneend};
        \draw [thin] (LL) --(LU);
        \draw [thin] (LU) --(RU);
        \draw [thin] (RL) --(RU);
        \draw [thin] (LL) --(RL);
        \draw [blue,ultra thin,dashed] (LL) --(IY);
        \draw [red,ultra thin,dashed]  (X)--(RL);
        \draw [blue,ultra thin,dashed] (RU) --(IX);
        \draw [red,ultra thin,dashed]  (LU)--(Y);
        \draw [purple,thick]  (LLD)--(RLD);
        \draw [purple,thick]  (LUD)--(RUD);
        \tkzLabelPoint[right](RU){$\pi$};
        \tkzLabelPoint[right](RL){$0$};
        \tkzLabelPoint[right](RLD){$\sin{\theta_1}=2r/R$, $d_1=2r$};
        \tkzLabelPoint[right](RUD){$\sin{\theta_1}=2r/R$, $d_1=2r$};
        \tkzLabelPoint[right](IY){$\theta_1=\Phistar$, $d_1=r\sin{\phistar}$};
        \tkzLabelPoint[right](Y){$\theta_1=\pi-\Phistar$, $d_1=r\sin{\phistar}$};
        \coordinate (RM) at ($(RL)!0.5!(RU)$);
        \tkzLabelPoint[right](RM){$\frac{\pi}{2}$};
        \coordinate (UpperM) at ($(LU)!0.5!(RU)$);
        \coordinate (LowerM) at ($(LL)!0.5!(RL)$);
        \coordinate (LM) at ($(LU)!0.5!(LL)$);
        \draw [ultra thin,dashed]  (LM)--(RM);
        \tkzLabelPoint[below](LL){$-\Phistar$};
        \tkzLabelPoint[below](RL){$+\Phistar$};
        \tkzLabelPoint[below](LowerM){$0$};
        
        \tkzDefPoint(0,0.5*\PHISTAR){LMend};
        \tkzDefPoint(0,\WIDTH-0.5*\PHISTAR){HMend};
        \tkzDefPoint(\PHISTAR,0.5*\PHISTAR){thetamin};
        \tkzDefPoint(\PHISTAR,\WIDTH-0.5*\PHISTAR){thetamax};
        \draw[ultra thin,dashed] (LMend)--(thetamin);
        \draw[ultra thin,dashed] (HMend)--(thetamax);
        \tkzLabelPoint[right](thetamin){$\theta_1=\frac{\Phistar}{2}$, $d_1=R\sin{(\Phistar/2)}=\frac{R\sin{\Phistar}}{2\cos{(\Phistar/2)}}=\frac{r\sin{\phistar}}{2\cos{(\Phistar/2)}}>\frac{r\sin{\phistar}}{2}$};
        \tkzLabelPoint[right](thetamax){$\theta_1=\pi-\frac{\Phistar}{2}$, $d_1=R\sin{(\Phistar/2)}=\frac{R\sin{\Phistar}}{2\cos{(\Phistar/2)}}=\frac{r\sin{\phistar}}{2\cos{(\Phistar/2)}}>\frac{r\sin{\phistar}}{2}$};
        \fill [green, opacity=10/30](LLD) -- (RLD) -- (IY) -- (LMend) -- (X) -- cycle;
        \fill [green, opacity=10/30](LUD) -- (RUD) -- (Y) -- (HMend) -- (IX) -- cycle;       \tkzDrawPoints(X,Y,IX,IY,LLD,LUD,RLD,RUD,LMend,HMend,thetamin,thetamax);
        \node[] at ([shift={(0,0.5*\PHISTAR)}]LMend) {\small $\mathcal{L}$};
        \node[] at ([shift={(0,-0.5*\PHISTAR)}]HMend) {\small $\mathcal{U}$};
    \end{tikzpicture}
    \caption{$x_1\in\mathcal{U}\sqcup\mathcal{L}=\big\{(\Phi_1,\theta_1)\in\MRin\cap\MRout\bigm|\theta_1<\sin^{-1}{(2r/R)}\big\}$, the purview of \cref{TEcase_b}.}\label{fig:A1contractionCase}
\end{center}
\end{figure}
\begin{proof}\strut\COMMENT{\Large\color{red}Are there any ideas in this proof that one could explain?\wentao{Added one line summarizing tools used are from \cite{cb}}.}%
The proof mainly relies on the conclusions in \cref{SBSCoordChange} which are the tools from \cite{cb}.

\hidesubsection{Proof of \eqref{itemExpansionLowerBoundB1}} By \cref{contraction_region} in case (b) of \eqref{eqMainCases}, $x_0\in\Nout$ and \cref{Prop:Udisjoint} shows that $\Nout\cap\textcolor{red}{M^{\text{\upshape{out}}}_{r,0}}=\emptyset$. This means that it is impossible to have $x\in\Mrout\cap\Mrin=\textcolor{red}{M^{\text{\upshape{out}}}_{r,0}}$. Thus, $x\in\Fc^{-1}(\Mrout\setminus\Mrin)$, $x_0=\Fc(x)\in\Mrout$.
    
    Hence, $dx_0\dfn D_{x}\Fc(dx)$ is a tangent vector at $x_0$. Since $x$, $x_0$ are nonsingular on $M_r$ and $dx=(d\phi,d\theta)$ is a tangent vector in the quadrant cone, by \cref{proposition:ShearOneStepExpansion} we have $\frac{\|dx_0\|_\p}{\|dx\|_\p}\ge1$.
    
    We also have $x_1=\Fc(x_0)\in\MRin$, $dx_1\dfn D_{x_0}\Fc(dx_0)$. With $n_1=0$ in the orbit segment, we have $x_2=\Fc(x_1)\in\Mrin$, $dx_2\dfn D_{x_1}\Fc(dx_1)=D_{x_0}\Fc^2(dx_0)=D_{x}\Fc^3(dx)$, $dx_2$ is a tangent vector at $x_2$.

    \hidesubsection{Proof of \eqref{itemExpansionLowerBoundB2}} By \cref{contraction_region}, we have $x_2\in\Nin$ and by \cref{Prop:Udisjoint}, we have $\emptyset=\Nin\cap(\textcolor{blue}{M^{\text{\upshape{in}}}_{r,0}}\cup\textcolor{blue}{M^{\text{\upshape{in}}}_{r,1}}\cup\textcolor{blue}{M^{\text{\upshape{in}}}_{r,2}})$. Hence, $m(x_2)\ge3$. By definition of $m(x_2)$ and $\hat{M}$  (\cref{def:SectionSetsDefs}), $(\phi_3,\theta_3)=x_3\dfn\Fc(x_2),\Fc(x_3)$, $\cdots$, $\Fc^{m(x_2)-2}(x_3)=\Fc^{m(x_2)-1}(x_2)\in\hat{M}$ are nonsingular points on $M_r\setminus\Mrout$ since their $\Fc$-images are not on $M_R$. 
    
    Now $dx_3\dfn D_x\Fc^4(dx)=D_{x_2}\Fc(dx_2)$ is a tangent vector at $x_3$. By \cref{caseB}, $D_x\Fc^4$ is a negative matrix. Since $(d\phi_3,d\theta_3)=dx_3=D_x\Fc^4(dx)$ and $dx=(d\phi,d\theta)$ is in the quadrant with $\frac{d\theta}{d\phi}\ge0$, $(d\phi_3,d\theta_3)=dx_3$ is then in the interior of the quadrant, that is, $\frac{d\theta_3}{d\phi_3}>0$. Then \cref{proposition:ShearOneStepExpansion} implies $\frac{\|D_x\hat{F}(dx)\|_\p}{\|dx_3\|_\p}>1$.

    \hidesubsection{Proof of \eqref{itemExpansionLowerBoundB3}} Since $\frac{d\theta}{d\phi}\in[0,1]$, $(d\phi_0,d\theta_0)=dx_0=D_x\Fc(dx)$, $D_x\Fc=\begin{pmatrix}
        1 & 2\\
        0 & 1
    \end{pmatrix}$, that is $\begin{pmatrix} d\phi_0\\d\theta_0
    \end{pmatrix}=\begin{pmatrix}
        1 & 2\\
        0 & 1
    \end{pmatrix}\begin{pmatrix} d\phi\\d\theta
    \end{pmatrix}$, we obtain $\frac{d\theta_0}{d\phi_0}=\frac{d\theta}{d\phi+2d\theta}=\frac{1}{2+\frac{1}{\frac{d\theta}{d\phi}}}\in\big[0,\frac{1}{3}\big]$. Hence, by \eqref{eq:coordinateschange} and \eqref{eq:MirrorEquationExpansionFormula}, we have $\mathcal{B}^+_0\cos\varphi_0+\frac{1}{r}=\frac{d\varphi_0}{ds_0}=\frac{-1}{r}\frac{d\theta_0}{d\phi_0}\in[-\frac{1}{3r},0]$. Therefore, $\mathcal{B}^+_0=-\frac{1}{r\cos\varphi_0}(1+\frac{d\theta_0}{d\phi_0})$. For $k=1+\frac{1}{2+\frac{d\theta}{d\phi}}\in[1,\frac{4}{3}]$, $\mathcal{B}^+_0=\frac{-k}{r\cos\varphi_0}=\frac{-k}{r\sin\theta_0}\overbracket=^{\mathclap{\text{\cref{def:fpclf}}}}\frac{-k}{d_0}\in[-\frac{-4}{3d_0},\frac{-1}{d_0}]$. 
    
    Then 
    \begin{equation}\label{eq:Ik}
    \begin{aligned}
        \frac{\|dx_2\|_\p}{\|dx_0\|_\p}&\overset{\eqref{eq:dx2overdx0}}=|(\tau_0+\tau_1-\frac{2\tau_0\tau_1}{d_1})\Bcal^{+}_0+1-\frac{2\tau_1}{d_1}|=|(\tau_0+\tau_1-\frac{2\tau_0\tau_1}{d_1})(\frac{-k}{d_0})+1-\frac{2\tau_1}{d_1}|=|\mathcal{I}(k)|,
    \end{aligned}
    \end{equation}
     where $\mathcal{I}(k)=(\tau_0+\tau_1-\frac{2\tau_0\tau_1}{d_1})(\frac{-k}{d_0})+1-\frac{2\tau_1}{d_1}$.
\end{proof}




\begin{proposition}[Case (b): length restrictions]\label{proposition:lengthsrestrictionincaseB}
Under the standing assumptions for this section (\cref{remark:standingassumptionInsection6}), the length functions  $\tau_0$, $\tau_1$, $d_0$, $d_1$, $d_2$ of  $x_1\in\big\{(\Phi_1,\theta_1)\in\MRin\cap\MRout\bigm|\sin\theta_1<2r/R\big\}$ (\cref{def:fpclf,fig:A1contractionCase}) satisfy
\COMMENT{All references to the following equations need to be fixed now, in each case a number 1,2,3, or 4 must be added to the label according to which specific equation is being used.\wentao{Work In progress. Explaination is provided in the proof for each equation.}}
\begin{gather}\label{eq:lengthrestrictions1}
    \frac{1}{2}r\sin{\phistar}<\frac{r\sin{\phistar}}{2\cos{\frac{\Phistar}{2}}}=\frac{R\sin{\Phistar}}{2\cos{\frac{\Phistar}{2}}}=R\sin{\frac{\Phistar}{2}}\le d_1<2r,\\\label{eq:lengthrestrictions2}
    \tau_0<2d_1,\quad \tau_1<2d_1, \quad\text{and}\quad d_0+d_2<\tau_0+\tau_1<d_0+d_2+\frac{40.7r^2}{R\sin{\phistar}},\\\label{eq:lengthrestrictions3}
    0<\tau_0+\tau_1-2d_0<\frac{40.7r^2}{R\sin{\phistar}} \text{ and }0<\tau_0+\tau_1-2d_2<\frac{40.7r^2}{R\sin{\phistar}},\\\label{eq:lengthrestrictions4}
    |d_0-r\sin{\phistar}|<\frac{8.3r^2}{R\sin{\phistar}},\quad |d_2-r\sin{\phistar}|<\frac{8.3r^2}{R\sin{\phistar}}\quad\text{and}\quad|d_0-d_2|<\frac{16.6r^2}{R\sin{\phistar}}.
\end{gather}
\end{proposition}
\begin{proof}
    For $x_1=(\Phi_1,\theta_1)$ in case (b) of \eqref{eqMainCases} with $d_1<2r$, $n_1=0$, the smallest $\theta_1$ can be is $\frac{\Phistar}{2}$, which corresponds to the lowest point of the region $\mathcal{L}$ in \cref{fig:A1contractionCase}. The corresponding $d_1=R\sin(\Phistar/2)=\frac{R\sin\Phistar}{2\cos{(\Phistar/2)}}\overbracket{=}^{\mathllap{R\sin\Phistar=r\sin\phistar}}\frac{r\sin\phistar}{2\cos{(\Phistar/2)}}>\frac{r\sin\Phistar}{2}$, hence \eqref{eq:lengthrestrictions1}.

    In \eqref{eq:lengthrestrictions2},  $\tau_i<2d_1$, $i=0,1$ follows from \cref{remark:deffpclf}, and  
    $0<\tau_0+\tau_1-d_0-d_2<\frac{40.7r^2}{R\sin{\phistar}}$ is \cite[Equations (4.9)]{hoalb}.

    As to \eqref{eq:lengthrestrictions3}, $0<\tau_0+\tau_1-2d_i<\frac{40.7r^2}{R\sin{\phistar}}$ for $i=0,2$ is \cite[Equations (4.8)]{hoalb}.  
    
    Finally, $|d_i-r\sin{\phistar}|<\frac{8.3r^2}{R\sin{\phistar}}$, for $i=0,2$ is \cite[Equations (4.6)]{hoalb},
    and $|d_0-d_2|\le|d_0-r\sin{\phistar}|$+$|r\sin{\phistar}-d_2|<\frac{16.6r^2}{R\sin{\phistar}}$ by the triangle inequality.
\end{proof}
These simple estimates will be used throughout this section by referring to \cref{eq:lengthrestrictions1,eq:lengthrestrictions2,eq:lengthrestrictions3,eq:lengthrestrictions4}.
\subsection[\texorpdfstring{$x_2\rightarrow x_3$}{x2 to x3 expansion} expansion]{\texorpdfstring{$x_2\rightarrow x_3$}{x2 to x3} expansion}\label{subsec:x2Tox3Expansion} We now prove that $\frac{\|dx_3\|_\p}{\|dx_2\|_\p}\ge1$.
\begin{proposition}\label{Proposition:DF4ConeUpperBound}
\label{def:DF4matrix}
Under the standing assumptions for this section (\cref{remark:standingassumptionInsection6}), 
$G\dfn d_1d_2 D_x\Fc^4=
    \begin{pmatrix}
        G_{11} & G_{12}\\
        G_{21} & G_{22}
    \end{pmatrix}$ is a negative matrix with 
$0<\frac{G_{22}}{G_{12}}<1$. 
\end{proposition}
\begin{proof}
By \cref{lemma:ExpansionLowerBoundB}\eqref{itemExpansionLowerBoundB1}, $x=\Fc^{-1}(x_0)$ and by \cref{caseB}, $D_{x}\Fc^{4}=D_{x_{-1}}\Fc^{4}$ is a $2$ by $2$ matrix with negative entries. Hence $G$ is a negative matrix. And equations \cite[Equations (4.11)(4.13)(4.15)(4.17)]{hoalb} give the following.
\begin{equation}\label{eq:Gmatrix}\begin{aligned}
            G_{11}&=6(\tau_1-\frac{2}{3}d_2)(\tau_0-d_0)-3d_1(\tau_0+\tau_1-d_0-\frac{2}{3}d_2)\text{ \cite[Equation 4.13]{hoalb}},\\
            G_{21}&=2(\tau_1-d_2)(\tau_0-d_0)-d_1(\tau_0+\tau_1-d_0-d_2)\text{ \cite[Equation 4.11]{hoalb}},\\
            G_{12}&=18(\tau_1-\frac{2}{3}d_2)(\tau_0-\frac{2}{3}d_0)-9d_1(\tau_0+\tau_1-\frac{2}{3}d_0-\frac{2}{3}d_2)\text{ \cite[Equation 4.17]{hoalb}},\\
            G_{22}&=6(\tau_1-d_2)(\tau_0-\frac{2}{3}d_0)-3d_1(\tau_0+\tau_1-\frac{2}{3}d_0-d_2)\text{ \cite[Equation 4.15]{hoalb}}.
\end{aligned}\end{equation} With $G_{12}<0$, $G_{22}<0$, it suffices to prove that $G_{12}-G_{22}<0$.
\begin{equation}\label{eq:G12G22diff}
        \begin{aligned}
            G_{12}-G_{22}&\overset{\mathclap{\eqref{eq:Gmatrix}}}=(12\tau_1-6d_2)(\tau_0-\frac{2}{3}d_0)-d_1(6\tau_0+6\tau_1-4d_0-3d_2)\\
            &=12(\tau_1-\frac{1}{2}d_2)(\tau_0-\frac{2}{3}d_0)-6d_1(\tau_0-\frac{2}{3}d_0+\tau_1-\frac{1}{2}d_2)\\
            & =6d_1(\tau_0-\frac{2}{3}d_0)(\tau_1-\frac{1}{2}d_2)(\frac{2}{d_1}-\frac{1}{\tau_0-\frac{2}{3}d_0}-\frac{1}{\tau_1-\frac{1}{2}d_2}).
        \end{aligned}   \end{equation}
We divide this into two cases.

\noindent\textbf{Case 1}: If $(\tau_0-\frac{2}{3}d_0)(\tau_1-\frac{1}{2}d_2)\le0$, then
\[G_{12}-G_{22}=12(\tau_1-\frac{1}{2}d_2)(\tau_0-\frac{2}{3}d_0)-6d_1(\tau_0-\frac{2}{3}d_0+\tau_1-\frac{1}{2}d_2)\le-6d_1(\underbracket{\tau_0+\tau_1}_{\mathclap{>d_0+d_2\text{ by \eqref{eq:lengthrestrictions2}}}}-\frac{2}{3}d_0-\frac{1}{2}d_2)\overbracket{<}^{\mathllap{\eqref{eq:lengthrestrictions2}}}-6d_1(\frac{1}{3}d_0+\frac{1}{2}d_2)<0.\] 

\noindent\textbf{Case 2}: If $(\tau_0-\frac{2}{3}d_0)(\tau_1-\frac{1}{2}d_2)>0$, then 
$\tau_0+\tau_1\overset{\eqref{eq:lengthrestrictions2}}>d_0+d_2$ implies \(\tau_0-\frac{2}{3}d_0>0\) and \(\tau_1-\frac{1}{2}d_2>0\), so once we prove $\frac{2}{d_1}<\frac{1}{\tau_0-\frac{2}{3}d_0}+\frac{1}{\tau_1-\frac{1}{2}d_2}$, \eqref{eq:G12G22diff} yields \[G_{12}-G_{22}=6d_1(\tau_0-\frac{2}{3}d_0)(\tau_1-\frac{1}{2}d_2)(\frac{2}{d_1}-\frac{1}{\tau_0-\frac{2}{3}d_0}-\frac{1}{\tau_1-\frac{1}{2}d_2})<0.\]

It remains to prove $\frac{2}{d_1}<\frac{1}{\tau_0-\frac{2}{3}d_0}+\frac{1}{\tau_1-\frac{1}{2}d_2}$. We do so in five subcases,  using \cref{proposition:lengthsrestrictionincaseB} frequently.

\noindent\textbf{Subcase 2.1}: If $\tau_0-\frac{2}{3}d_0\ge\frac{10}{9}d_1$, $\tau_1-\frac{1}{2}d_2>0$, then 
$2d_1\overset{\eqref{eq:lengthrestrictions2}}{>}\tau_0\ge\frac{2}{3}d_0+\frac{10}{9}d_1\Rightarrow \frac{8}{9}d_1>\frac{2}{3}d_2$, so $d_1\ge \frac{3}{4}d_0$. Therefore, we have
\(\displaystyle d_0+d_2+\frac{40.7r^2}{R\sin\phistar}\overbracket{>}^{\mathclap{\eqref{eq:lengthrestrictions2}}}\tau_0+\tau_1\overbracket{\ge}^{\mathclap{\tau_0\ge\frac{2}{3}d_0+\frac{10}{9}d_1}}\frac{2}{3}d_0+\frac{10}{9}d_1+\tau_1\), which implies \(\displaystyle\tau_1<\frac{1}{3}d_0+d_2-\frac{10}{9}d_1+\frac{40.7r^2}{R\sin\phistar}\) and hence
\begin{align*}
    0\overbracket<^{\mathclap{\text{Subcase 2.1}}}\tau_1-\frac{1}{2}d_2&<\frac{1}{3}d_0+\frac{1}{2}d_2-\frac{10}{9}d_1+\frac{40.7r^2}{R\sin\phistar}\overbracket{\le}^{d_1\ge\frac{3}{4}d_0}\frac{1}{3}d_0+\frac{1}{2}d_2-\frac{10}{9}\cdot\frac{3}{4}d_0+\frac{40.7r^2}{R\sin\phistar}=\frac{1}{2}(d_2-d_0)+\frac{40.7r^2}{R\sin\phistar}\\
    &\overbracket{<}^{\mathclap{\eqref{eq:lengthrestrictions4}}}\frac{16.6r^2}{2R\sin\phistar}+\frac{40.7r^2}{R\sin\phistar}=\frac{49r^2}{R\sin\phistar}\overbracket{<}^{\mathclap{R>\frac{196r}{\sin^2{\phistar}}}}\frac{r\sin\phistar}{4}\overset{\eqref{eq:lengthrestrictions1}}<\frac{d_1}{2}.
\end{align*}
thus, $\frac{2}{d_1}<\frac{1}{\tau_1-\frac{1}{2}d_2}<\frac{1}{\tau_1-\frac{1}{2}d_2}+\frac{1}{\tau_0-\frac{2}{3}d_0}$.

\noindent\textbf{Subcase 2.2}: If $0<\tau_0-\frac{2}{3}d_0\le\frac{1}{2}d_1$, $\tau_1-\frac{1}{2}d_2>0$, then we directly get $\frac{2}{d_1}\le\frac{1}{\tau_0-\frac{2}{3}d_0}<\frac{1}{\tau_0-\frac{2}{3}d_0}+\frac{1}{\tau_1-\frac{1}{2}d_2}$.

\noindent\textbf{Subcase 2.3}: $d_1\le\tau_0-\frac{2}{3}d_0<\frac{10}{9}d_1$, $\tau_1-\frac{1}{2}d_2>0$.

From \eqref{eq:lengthrestrictions2} we have $2d_1\overbracket{>}^{\mathclap{\eqref{eq:lengthrestrictions2}}}\tau_0\overbracket{\ge}^{\mathrlap{\text{Subcase 2.3}}}\frac{2d_0}{3}+d_1$, so $d_1> \frac{2}{3}d_0$. \cref{proposition:lengthsrestrictionincaseB} also gives 
\begin{align*}
    0\overbracket<^{\mathclap{\text{Subcase 2.3}}}\tau_1-\frac{1}{2}d_2\overbracket{<}^{\mathclap{\eqref{eq:lengthrestrictions2}}}&d_0+\frac{1}{2}d_2-\tau_0+\frac{40.7r^2}{R\sin\phistar}\overbracket{\le}^{\mathrlap{\text{Subcase 2.3}}\quad} d_0+\frac{1}{2}d_2-(\frac{2}{3}d_0+d_1)+\frac{40.7r^2}{R\sin\phistar}=\frac{1}{3}d_0+\frac{1}{2}d_2-d_1+\frac{40.7r^2}{R\sin\phistar}
    \\
    \overbracket{<}^{\quad\mathllap{d_1>\frac{2}{3}d_0}}&\frac{1}{3}d_0+\frac{1}{2}d_2-\frac{2}{3}d_0+\frac{40.7r^2}{R\sin\phistar}=\frac{1}{2}d_2-\frac{1}{3}d_0+\frac{40.7r^2}{R\sin\phistar}
    \\ \overbracket<^{\quad\mathllap{\eqref{eq:lengthrestrictions4}}}&\frac{1}{2}(r\sin{\phistar}+\frac{8.3r^2}{R\sin{\phistar}})-\frac{1}{3}(r\sin{\phistar}-\frac{8.3r^2}{R\sin{\phistar}})+\frac{40.7r^2}{R\sin{\phistar}}=\frac{1}{6}r\sin{\phistar}+\frac{(\frac{5}{6}\times8.3+40.7)r^2}{R\sin{\phistar}}
    \\
    \overset{\quad\mathllap{R>\frac{166r}{\sin^2\phistar}>\frac{(\frac{5}{6}\times8.3+40.7)r}{(\frac{5}{11}-\frac{1}{6})\sin^2\phistar}}}<&\frac{10}{11}\times\frac{1}{2}r\sin\phistar 
    \\ \overbracket<^{\quad\mathllap{\eqref{eq:lengthrestrictions4}: d_1>\frac{1}{2}r\sin\phistar}}&\frac{10}{11}d_1.
\end{align*}
Thus, $\frac{1}{\tau_1-\frac{1}{2}d_2}>\frac{11}{10d_1}$. The definition of \textbf{Subcase 2.3} includes $\frac{1}{\tau_0-\frac{2}{3}d_0}>\frac{9}{10d_1}$, so $\frac{2}{d_1}=\frac{11}{10d_1}+\frac{9}{10d_1}<\frac{1}{\tau_1-\frac{1}{2}d_2}+\frac{1}{\tau_0-\frac{2}{3}d_0}$.

\noindent\textbf{Subcase 2.4}: $0.7d_1\le\tau_0-\frac{2}{3}d_0<d_1, \tau_1-\frac{1}{2}d_2>0$.

We have $2d_1\overbracket{>}^{\mathclap{\eqref{eq:lengthrestrictions2}}}\tau_0\ge\frac{2}{3}d_0+0.7d_1$, so $d_1>\frac{2}{3.9}d_0$. \cref{proposition:lengthsrestrictionincaseB} also gives\COMMENT{\red\bfseries This the first place where 1758 appears, and apparently the only one. We should refer to here from the discussion near \eqref{eqJZHypCond}. This is also the right time to comment on ``40.7,'' ``8.3,'' and ``39.''  How do they come about, what exact expressions do they replace, how close are they? Presumably $1758\approx39(40.7+\frac{8.3}{39}+\frac{8.3}{2})$.}
\begin{align*}
    0\overbracket{<}^{\mathclap{\text{Subcase 2.4}}}\tau_1-\frac{1}{2}d_2\overbracket{<}^{\mathclap{\eqref{eq:lengthrestrictions2}}}&d_0+\frac{1}{2}d_2-\tau_0+\frac{40.7r^2}{R\sin\phistar}\overbracket{\le}^{\mathclap{\substack{\text{Subcase 2.4}\\\tau_0\ge0.7d_1+\frac{2}{3}d_0\bigstrut}}} d_0+\frac{1}{2}d_2-(0.7d_1+\frac{2}{3}d_0)+\frac{40.7r^2}{R\sin\phistar}=\frac{1}{3}d_0+\frac{1}{2}d_2-0.7d_1+\frac{40.7r^2}{R\sin\phistar}\\
    \overbracket{<}^{\mathclap{d_1>\frac{2}{3.9}d_0}}&\frac{1}{2}d_2+(\frac{1}{3}-\frac{2\times0.7}{3.9})d_0+\frac{40.7r^2}{R\sin\phistar}=\frac{1}{2}d_2-\frac{1}{39}d_0+\frac{40.7r^2}{R\sin\phistar}\\ \overbracket{<}^{\mathclap{\eqref{eq:lengthrestrictions4}}}&\frac{1}{2}(r\sin{\phistar}+\frac{8.3r^2}{R\sin{\phistar}})-\frac{1}{39}(r\sin{\phistar}-\frac{8.3r^2}{R\sin{\phistar}})+\frac{40.7r^2}{R\sin{\phistar}}\\
    &=\big(\underbracket{-\frac{1}{39}+(40.7+\frac{8.3}{39}+\frac{8.3}{2})\frac{r}{R\sin^2\phistar}}_{<0\text{ when }R>\frac{1758r}{\sin^2{\phistar}}}\big)r\sin\phistar+\frac{1}{2}r\sin\phistar\overbracket{<}^{\mathclap{R>\frac{1758r}{\sin^2{\phistar}}}}\frac{1}{2}r\sin\phistar\overset{\eqref{eq:lengthrestrictions1}}<d_1.
\end{align*}
Together with $\frac{1}{\tau_0-\frac{2}{3}d_0}>\frac{1}{d_1}$ (from the definition of this subcase), this gives $\frac{2}{d_1}=\frac{1}{d_1}+\frac{1}{d_1}<\frac{1}{\tau_0-\frac{2}{3}d_0}+\frac{1}{\tau_1-\frac{1}{2}d_2}$.

\noindent\textbf{Subcase 2.5}: $\frac{1}{2}d_1<\tau_0-\frac{2}{3}d_0<0.7d_1$, $\tau_1-\frac{1}{2}d_2>0$.

$2d_1\overbracket{>}^{\mathllap{\eqref{eq:lengthrestrictions2}}}\tau_0\overbracket{>}^{\mathrlap{\text{Subcase 2.5}}}\frac{1}{2}d_1+\frac{2}{3}d_0$, so $d_1>\frac{4}{9}d_0$. \cref{proposition:lengthsrestrictionincaseB} also gives
\begin{align*}
    0\overbracket{<}^{\mathllap{\text{subcase 2.5}}}\tau_1-\frac{1}{2}d_2\overbracket{<}^{\eqref{eq:lengthrestrictions2}}&d_0+\frac{1}{2}d_2+\frac{40.7r^2}{R\sin\phistar}-\tau_0\overbracket{<}^{\mathclap{\tau_0>\frac{1}{2}d_1+\frac{2}{3}d_0}}d_0+\frac{1}{2}d_2+\frac{40.7r^2}{R\sin\phistar}-\frac{1}{2}d_1-\frac{2}{3}d_0\\
    \overbracket{<}^{\mathllap{d_1>\frac{4}{9}d_0}}&\frac{1}{3}d_0+\frac{1}{2}d_2-\frac{1}{2}\cdot\frac{4}{9}d_0+\frac{40.7r^2}{R\sin\phistar}=\frac{1}{9}d_0+\frac{1}{2}d_2+\frac{40.7r^2}{R\sin\phistar}\\ \overbracket{<}^{\mathllap{\eqref{eq:lengthrestrictions2}}}&\frac{1}{9}(r\sin\phistar+\frac{8.3r^2}{R\sin\phistar})+\frac{1}{2}(r\sin\phistar+\frac{8.3r^2}{R\sin\phistar})+\frac{40.7r^2}{R\sin\phistar}\\
    &=\frac{1}{9}r\sin{\phistar}+\frac{1}{2}r\sin{\phistar}+(\frac{1}{9}+\frac{1}{2})\frac{8.3r^2}{R\sin{\phistar}}+\frac{40.7r^2}{R\sin{\phistar}}\overbracket{<}^{\mathclap{R>\frac{174r}{\sin^2{\phistar}}}}\frac{7}{8}r\sin\phistar\overbracket{<}^{\mathllap{\eqref{eq:lengthrestrictions1}}}\frac{7}{4}d_1,
\end{align*}
where we used the fact that $R>\frac{174r}{\sin^2{\phistar}}\Rightarrow\big[(\frac{1}{9}+\frac{1}{2})\times8.3+40.7\big]\frac{r^2}{R\sin{\phistar}}<\big(\frac{7}{8}-\frac{1}{2}-\frac{1}{9}\big)r\sin{\phistar}$.

Together with $\frac{1}{\tau_0-\frac{2}{3}d_0}>\frac{10}{7d_1}$ (from the definition of the subcase), this gives $\frac{2}{d_1}=\frac{4}{7d_1}+\frac{10}{7d_1}<\frac{1}{\tau_0-\frac{2}{3}d_0}+\frac{1}{\tau_1-\frac{1}{2}d_2}$.
\end{proof}
\begin{proposition}\label{proposition:dx3dx2expansion}
Under the standing assumptions for this section (\cref{remark:standingassumptionInsection6}), we have $dx_2=(d\phi_2,d\theta_2)$ satisfies $\frac{d\theta_2}{d\phi_2}\in(-\infty,-1)\cup(0,+\infty)$ and  $\frac{\|dx_3\|_\p}{\|dx_2\|_\p}\ge1$.
\end{proposition}
\begin{proof}
By \cref{lemma:ExpansionLowerBoundB}\eqref{itemExpansionLowerBoundB1}, we have $x=\Fc^{-1}(x_0)=x_{-1}$. And by \cref{caseB}, $D_{x}\Fc^{4}=D_{x_{-1}}\Fc^{4}$ is a negative matrix. Then in \cref{lemma:ExpansionLowerBoundB}\eqref{itemExpansionLowerBoundB2}, for $dx=(d\phi,d\theta)$ in the cone $\big\{dx=(d\phi,d\theta)\bigm|\frac{d\theta}{d\phi}\in [0,1]\big\}$, we have $dx_3=D_{x}\Fc^4(dx)$ is a tangent vector at $x_3$ in the cone $\big\{dx_3=(d\phi_3,d\theta_3)\bigm|\frac{d\theta_3}{d\phi_3}\in(0,+\infty) \big\}$, $dx_2=D_{x}\Fc^{3}(dx)$ and $dx_3=D_{x_2}\Fc(dx_2)$ for $dx_2=(d\phi_2,d\theta_2)$ with $D_{x_2}\Fc=\begin{pmatrix}
1 & 2 \\
0 & 1 \end{pmatrix}$, that is, $\begin{pmatrix} d\phi_3\\d\theta_3
    \end{pmatrix}=\begin{pmatrix}
        1 & 2\\
        0 & 1
    \end{pmatrix}\begin{pmatrix} d\phi_2\\d\theta_2
    \end{pmatrix}$. $x_2$, $x_3$ are both nonsingular points in $M_r$ and we get a pair of equivalent ``continued-fraction'' equations:    \begin{equation}\label{eq:Twocontinuedfraction}
            \frac{d\theta_3}{d\phi_3}=\frac{d\theta_2}{d\phi_2+2d\theta_2} = \frac{1}{2}-\frac{1}{2}\frac{1}{1+2\frac{d\theta_2}{d\phi_2}},\qquad 
            \frac{d\theta_2}{d\phi_2}=\frac{d\theta_3}{d\phi_3-2d\theta_3} = -\frac{1}{2}+\frac{1}{2}\frac{1}{1-2\frac{d\theta_3}{d\phi_3}}.
    \end{equation}

For the negative matrix $G$ in \cref{def:DF4matrix} and with $\frac{d\theta}{d\phi}=k$ for some $k\in[0,1]$, 
\begin{align*}
    \begin{pmatrix} d\phi_3\\d\theta_3
    \end{pmatrix}=dx_3=D_x\Fc^4(dx)=\frac{1}{d_1d_2}\begin{pmatrix}
        G_{11} & G_{12}\\
        G_{21} & G_{22}
    \end{pmatrix}\begin{pmatrix} d\phi\\d\theta
    \end{pmatrix},
    \text{ so }\frac{d\theta_3}{d\phi_3}=\frac{G_{21}+kG_{22}}{G_{11}+kG_{12}}.
\end{align*}
Since the billiard map (and its iteration) is orientation-preserving \cite[Equation (2.29)]{cb}, the negative matrix $G$ in \cref{def:DF4matrix} has positive determinant: $G_{11}G_{22}-G_{12}G_{21}>0$. Hence $\frac{G_{22}}{G_{12}}-\frac{G_{21}}{G_{11}}=\frac{G_{11}G_{22}-G_{12}G_{21}}{G_{11}G_{12}}>0$ and\COMMENT{Format this to match the counterpart in Section 5.}
\begin{equation}\label{eq:casebdtheta3dphi3ratio}
\begin{aligned}
    \frac{d\theta_3}{d\phi_3}-\frac{G_{21}}{G_{11}}=\frac{G_{21}+kG_{22}}{G_{11}+kG_{12}}-\frac{G_{21}}{G_{11}}&=\frac{\overbracket{k}^{\mathclap{\ge0}}(\overbracket{G_{11}G_{22}-G_{21}G_{12}}^{\mathrlap{>0}})}{\underbracket{G_{11}(G_{11}+kG_{12})}_{>0}}\ge0,\\
    \frac{G_{22}}{G_{12}}-\frac{d\theta_3}{d\phi_3}=\frac{G_{22}}{G_{12}}-\frac{G_{21}+kG_{22}}{G_{11}+kG_{12}}&=\frac{(\overbracket{G_{11}G_{22}-G_{21}G_{12}}^{\mathrlap{>0}})}{\underbracket{G_{12}(G_{11}+kG_{12})}_{>0}}>0,\text{ hence}\\
    0<\frac{G_{21}}{G_{11}}\le\frac{d\theta_3}{d\phi_3}&=\frac{G_{21}+kG_{22}}{G_{11}+kG_{12}}<\frac{G_{22}}{G_{12}}\overbracket{<}^{\mathrlap{\text{\cref{def:DF4matrix}}}}1. 
\end{aligned}
\end{equation}
With $0<\frac{G_{21}}{G_{11}}\le\frac{d\theta_3}{d\phi_3}<\frac{G_{22}}{G_{12}}<1$, by \eqref{eq:Twocontinuedfraction} $\frac{d\theta_2}{d\phi_2}\in(-\infty,-1)\cup (0,+\infty]$ and we must analyze the following two cases:
\[
    \frac{d\theta_2}{d\phi_2}\in \big(0,+\infty\big]\overset{\eqref{eq:Twocontinuedfraction}}{\Longleftrightarrow} \frac{d\theta_3}{d\phi_3}\in\big(0,\frac{1}{2}\big]
    \quad\text{versus}\quad\frac{d\theta_2}{d\phi_2}\in \big[-\infty,-1\big) \overset{\eqref{eq:Twocontinuedfraction}}{\Longleftrightarrow} \frac{d\theta_3}{d\phi_3}\in\big[\frac{1}{2},1\big).
\]

In the first case, 
\cref{proposition:ShearOneStepExpansion} gives $\frac{\|dx_3\|_\p}{\|dx_2\|_\p}\ge1$.
%

In the second case, adapt the notations from \eqref{eq:MirrorEquationExpansionFormula}, denote by $\mathcal{B}_2^+$ the wave front curvature after collision at $x_2$, and note that the length of the free path between collisions $p(x_2)$ and $p(x_3)$ in the billiard table is $\tau_2=2r\sin\theta_2$. Then
\[1\overbracket{<}^{\mathclap{\frac{d\theta_2}{d\phi_2}<-1}}-1-2\frac{d\theta_2}{d\phi_2}\overbracket{=}^{\mathclap{\eqref{eq:coordinateschange}}}-1-2r\frac{-d\varphi_2}{ds_2}\overbracket{=}^{\mathclap{\eqref{eq:MirrorEquationExpansionFormula}}}-1+\mathcal{B}_2^{+}\overbracket{2r\cos\varphi_2}^{\mathclap{=2r\sin\theta_2=\tau_2}}-2\overbracket{r\mathcal{K}}^{=-1}=1+\tau_2\mathcal{B}^+_2\overbracket=^{\mathclap{\eqref{eq:MirrorEquationExpansionFormula}\ \ }}\frac{\|dx_3\|_\p}{\|dx_2\|_\p}.\qedhere\]
\end{proof}

\subsection[\texorpdfstring{$x_0\longrightarrow x_1\longrightarrow x_2$}{x0 to x1 to x2} Transition Contraction]{\texorpdfstring{$x_0\longrightarrow x_1\longrightarrow x_2$}{x0 to x1 to x2} contraction}\label{subsec:A1contractionCase}\hfill

\begin{definition}\label{def:casebhandgfunctions}
Under the standing assumptions for this section (\cref{remark:standingassumptionInsection6}), we define functions $g,h$ of $x_1\in\big\{(\Phi_1,\theta_1)\in\MRin\cap\MRout\bigm|\sin\theta_1<2r/R\big\}$  (\cref{fig:A1contractionCase}) by
\begin{equation}\label{eq:casebhandgfunctions}
        \begin{aligned}
            & h=\frac{\tau_0}{d_0}+\frac{\tau_1}{d_0}-2,\\
            & g=d_0-r\sin{\phistar}.
        \end{aligned}
    \end{equation}
\end{definition}
\begin{proposition}\label{proposition:casebhandgfunctions}
Under the standing assumptions for this section (\cref{remark:standingassumptionInsection6}), \cref{proposition:lengthsrestrictionincaseB} with \eqref{eqJZHypCond}: $ R>\frac{30000r}{\sin^2\phistar}$ implies that \COMMENT{???????\red Fix this statement!!! It makes no sense\wentao{Much rewritten,Looks better?}} the functions $g$, $h$ from \cref{def:casebhandgfunctions} and $d_0$ satisfy the following.
\begin{gather}\label{eq:casebghd0lengthrestrictions1}
    d_0-r\sin{\phistar}\overbracket{=}^{\mathclap{\eqref{eq:casebhandgfunctions}\ }}g\overbracket{\in}^{\mathclap{\ \eqref{eq:lengthrestrictions4}}}\big(-\frac{8.3 r^2}{R\sin\phistar},\frac{8.3r^2}{R\sin\phistar}\big)\overbracket{\subset}^{\mathclap{\quad R>\frac{30000r}{\sin^2\phistar}}} \big(-0.0003r\sin\phistar,0.0003r\sin\phistar\big),\\\label{eq:casebghd0lengthrestrictions2}
    0.9997r\sin\phistar\overbracket{<}^{\mathclap{\eqref{eq:casebghd0lengthrestrictions1}}}d_0\overbracket{<}^{\mathclap{\eqref{eq:casebghd0lengthrestrictions1}}}1.0003r\sin\phistar,\\\label{eq:casebghd0lengthrestrictions3}
    0\overbracket{<}^{\mathclap{\eqref{eq:lengthrestrictions3}\ }}h\overbracket{=}^{\mathclap{\strut\ \eqref{eq:casebhandgfunctions}}}\frac{\tau_0}{d_0}+\frac{\tau_1}{d_0}-2=\frac{\tau_0+\tau_1-2d_0}{d_0}\overbracket{<}^{\mathclap{\eqref{eq:lengthrestrictions3}\ }}\;\frac{\frac{40.7r^2}{R\sin\phistar}}{d_0}\overbracket{<}^{\mathclap{\ \ \eqref{eq:casebghd0lengthrestrictions2}}}\;\frac{\frac{40.7r^2}{R\sin\phistar}}{0.9997r\sin\phistar}<\frac{40.8r}{R\sin^2{\phistar}},\\\label{eq:casebghd0lengthrestrictions4}
    d_2-r\sin{\phistar}\overbracket{\in}^{\mathclap{\eqref{eq:lengthrestrictions4}}}\big(-\frac{8.3 r^2}{R\sin\phistar},\frac{8.3r^2}{R\sin\phistar}\big)\overbracket{\subset}^{\mathclap{\quad R>\frac{30000r}{\sin^2\phistar}}} \big(-0.0003r\sin\phistar,0.0003r\sin\phistar\big),\\\label{eq:casebghd0lengthrestrictions5}
    0.9997r\sin\phistar\overbracket{<}^{\mathclap{\eqref{eq:casebghd0lengthrestrictions4}\ }}d_2\overbracket{<}^{\mathclap{\ \eqref{eq:casebghd0lengthrestrictions4}}}1.0003r\sin\phistar.
    \end{gather}
\end{proposition}
\begin{proposition}\label{proposition:IktransionInCaseb}
Under the standing assumptions for this section (\cref{remark:standingassumptionInsection6}), $\mathcal{I}(k)=(\tau_0+\tau_1-\frac{2\tau_0\tau_1}{d_1})(\frac{-k}{d_0})+1-\frac{2\tau_1}{d_1}<-0.05$ for all $k\in[1,\frac{4}{3}]$ and all $x_1\in\big\{x_1=(\Phi_1,\theta_1)\in\MRin\cap\MRout\bigm|\sin\theta_1<2r/R\big\}$ (see \cref{fig:A1contractionCase}).
\end{proposition}
\begin{proof}
We first rewrite $\mathcal{I}(k)$ with $h$ from \cref{def:casebhandgfunctions}, that is,
\begin{equation}\label{eq:CaseBIkQuadraticPolynomial}
    \begin{aligned}
    \mathcal{I}(k)&=(\tau_0+\tau_1-\frac{2\tau_0\tau_1}{d_1})(\frac{-k}{d_0})+1-\frac{2\tau_1}{d_1}\overset{\eqref{eq:casebhandgfunctions}}=\big[(2+h)d_0-\frac{2\tau_0\tau_1}{d_1}\big]\big(\frac{-k}{d_0}\big)+1-\frac{2\tau_1}{d_1}\\
    &=-(2+h)k+\frac{2\tau_0\tau_1}{d_0d_1}k+1-\frac{2\tau_1}{d_1}=-2k+\frac{2\tau_0\tau_1}{d_0d_1}k+1-\frac{2\tau_1}{d_1}-hk\\
    &=-2k+1+\frac{2\tau_1}{d_1}(\frac{\tau_0}{d_0}k-1)-hk=-2k+1+\frac{2d_0}{d_1}\frac{\tau_1}{d_0}k(\frac{\tau_0}{d_0}-\frac{1}{k})-hk\\ \overbracket{=}^{\mathclap{\eqref{eq:casebhandgfunctions}: \frac{\tau_0}{d_0}=2+h-\frac{\tau_1}{d_0}}}& -2k+1+\frac{2d_0}{d_1}k\frac{\tau_1}{d_0}(2-\frac{1}{k}-\frac{\tau_1}{d_0})+\frac{2k\tau_1}{d_1}h-hk.
\end{aligned}
\end{equation}
Then, by symmetry, we assume that $x_1\in \mathcal{L}$ in \cref{fig:A1contractionCase} and consider two cases to analyze this: Case 1: $d_1\ge\frac{21}{32}r\sin{\phistar}$; Case 2: $d_1<\frac{21}{32}r\sin{\phistar}$.
\begin{figure}[h]
\begin{center}
    \begin{tikzpicture}
        \pgfmathsetmacro{\PHISTAR}{3};
        \pgfmathsetmacro{\Dlt}{0.9};
        \pgfmathsetmacro{\WIDTH}{5};
        \pgfmathsetmacro{\boudaryShift}{1.0};
        \pgfmathsetmacro{\LShift}{0.5};
        \tkzDefPoint(-\PHISTAR,0){LL};
        \tkzDefPoint(-\PHISTAR,\WIDTH){LU};
        \tkzDefPoint(\PHISTAR,0){RL};
        \tkzDefPoint(\PHISTAR,\WIDTH){RU};
        \tkzDefPoint(-\PHISTAR,\PHISTAR){X};
        \tkzDefPoint(\PHISTAR,\PHISTAR){IY};
        \tkzDefPoint(-\PHISTAR,\PHISTAR+\Dlt){LLD};
        \tkzDefPoint(\PHISTAR,\PHISTAR+\Dlt){RLD};
        \tkzDefPoint(0,0.5*\PHISTAR){LMD};
        \draw [ultra thick,red] (X) --(LLD);
        \draw [ultra thin,dashed] (LU) --(RU);
        \draw [ultra thin,dashed] (LU) --(LLD);
        \draw [ultra thin,dashed] (X) --(LL);
        \draw [ultra thin,dashed] (RU) --(RLD);
        \draw [ultra thin,dashed] (IY) --(RL);
        \draw [ultra thick,red] (IY) --(RLD);
        \draw [ultra thin,dashed] (LL) --(RL);
        \draw [red,ultra thick] (LMD) --(IY);
        \draw [red,ultra thick]  (X)--(LMD);
        \draw [red,ultra thick]  (LLD)--(RLD);
        \coordinate (LMDMark) at ([shift={(\PHISTAR,0)}]LMD);
        \coordinate (RBoundary) at ([shift={(-2*\boudaryShift,-\boudaryShift)}]IY);
        \coordinate (LBoundary) at ([shift={(2*\boudaryShift,-\boudaryShift)}]X);
         \coordinate (LlineBoundary) at ([shift={(2*\LShift,-\LShift)}]X);
         \coordinate (RlineBoundary) at ([shift={(-2*\LShift,-\LShift)}]IY);
        \coordinate (RlineBoundaryMark) at ([shift={(0,-\LShift)}]IY);
        \draw [ultra thick] (LBoundary) --(RBoundary);
        \draw [thick,blue] (LlineBoundary) --(RlineBoundary);
        \draw [ultra thin,blue,dashed] (RlineBoundary) --(RlineBoundaryMark);
        \coordinate (RBoundaryMark) at ([shift={(0,-\boudaryShift)}]IY);
        \draw [ultra thin,dashed]  (LMD)--(LMDMark);
        \draw [ultra thin,dashed] (RBoundary)--(RBoundaryMark);
        \fill [red, opacity=8/30](X) -- (LBoundary) -- (RBoundary) -- (IY) --(RLD)--(LLD)--cycle;
        \tkzDrawPoints(IY,RLD,LMD,RBoundaryMark,LMDMark,RlineBoundaryMark);
        \tkzLabelPoint[right](IY){$d_1=r\sin{\phistar}$};
        \tkzLabelPoint[right](RLD){$d_1=2r$};
        \tkzLabelPoint[right](RlineBoundaryMark){\textcolor{blue}{$d_1=\const$}};
        \tkzLabelPoint[right](RBoundaryMark){$d_1=\frac{21}{32}r\sin{\phistar}$};
        \tkzLabelPoint[right](LMDMark){$d_1=\frac{r\sin{\phistar}}{2\cos{\frac{\Phistar}{2}}}$};
        \tkzLabelPoint[below](LL){$-\Phistar$};
        \tkzLabelPoint[below](RL){$\Phistar$}
    \end{tikzpicture}
    \caption{
    \textbf{Case 1} Region $\mathcal{L}$ in \cref{fig:A1contractionCase} with $d_1\ge\frac{21}{32}r\sin{\phistar}$
    }\label{fig:b_contractionCase1}
\end{center}
\end{figure}
\hidesubsection{Case 1} $x_1\in\mathcal{L}$ and $2r>d_1\ge\frac{21}{32}r\sin\phistar$ (see \cref{fig:b_contractionCase1}).

Suppose $x_1$ is on the line $d_1=\const\ge\frac{21}{32}r\sin{\phistar}$  inside $\mathcal{L}$ (the blue line in \cref{fig:b_contractionCase1}).
Then 
$\frac{(a+b)^2}{4}-ab=\frac{(a-b)^2}{4}\ge0\Longrightarrow ab\le\frac{(a+b)^2}{4}$ implies 
\begin{equation}\label{arithmeticgeometricmean}
\overbracket{\frac{\tau_1}{d_0}}^{=a}\overbracket{(2-\frac{1}{k}-\frac{\tau_1}{d_0})}^{=b}\le\frac{1}{4}(\overbracket{2-\frac{1}{k}}^{=a+b})^2.
\end{equation}
Thus
, we get 
\begin{align*}
    \mathcal{I}(k)&\overbracket{=}^{\mathclap{\eqref{eq:CaseBIkQuadraticPolynomial}}}-2k+1+\frac{2d_0}{d_1}k\times\frac{\tau_1}{d_0}(2-\frac{1}{k}-\frac{\tau_1}{d_0})
    +\frac{2k\tau_1}{d_1}h-hk\\
    &\overbracket{\le}^{\mathclap{\substack{d_0>0,d_1>0,k>0\\\text{and }\eqref{arithmeticgeometricmean}}}}-2k+1+\frac{2d_0}{d_1}k\times\frac{1}{4}(2-\frac{1}{k})^2+\frac{2k\tau_1}{d_1}h-hk=-2k+1+\frac{kd_0}{2d_1}(2-\frac{1}{k})^2+\frac{2k\tau_1}{d_1}h-hk\\
    &\overbracket{=}^{\mathclap{\eqref{eq:casebhandgfunctions}}}-2k+1+\frac{k}{2d_1}(2-\frac{1}{k})^2(g+r\sin\phistar)+\frac{2k\tau_1}{d_1}h-hk\\
    &= -2k+1+\frac{kr\sin{\phistar}}{2d_1}(2-\frac{1}{k})^2+\frac{kg}{2d_1}(2-\frac{1}{k})^2+\frac{2k\tau_1}{d_1}h-hk\\
    &=\underbracket{\frac{r\sin{\phistar}}{2d_1}}_{\mathclap{\text{(Case 1): }\le\frac12\frac{32}{21}}}\frac{(2k-1)^2}{k}+1-2k+\frac{kg}{2d_1}(2-\frac{1}{k})^2+(\frac{2\tau_1}{d_1}-1)hk\\
    &\le\frac{32}{21}\frac{(2k-1)^2}{2k}+1-2k+\frac{kg}{2d_1}(2-\frac{1}{k})^2+(\frac{2\tau_1}{d_1}-1)hk=\frac{32}{21}(1-\frac{1}{2k})(2k-1)-2k+1+\frac{kg}{2d_1}(2-\frac{1}{k})^2+(\frac{2\tau_1}{d_1}-1)hk\\
    &=-\frac{43}{21}+\underbracket{\frac{22k}{21}+\frac{16}{21k}}_{\mathclap{\text{increases with $k\in[1,\frac{4}{3}]$}}}+\frac{kg}{2d_1}(2-\frac{1}{k})^2+(\frac{2\tau_1}{d_1}-1)hk
    \le\underbracket{-\frac{43}{21}+\frac{22}{21}\times\frac{4}{3}+\frac{16}{21}\times\frac{3}{4}}_{=-5/63<-0.079}+\underbracket{\frac{kg}{2d_1}(2-\frac{1}{k})^2}_{<0.0005}+\underbracket{(\frac{2\tau_1}{d_1}-1)hk}_{\mathclap{<0.00544}}<-0.07,
\end{align*}
since 
\[(\frac{2\tau_1}{d_1}-1)hk\overbracket{<}^{\mathclap{\substack{1\le k\le4/3,\text{\eqref{eq:casebghd0lengthrestrictions3}: }h>0\\\text{\eqref{eq:lengthrestrictions2}: } \tau_1<2d_1}}}(2\cdot2-1)\cdot\frac{4}{3}h=4h\overbracket{<}^{\mathclap{\text{\eqref{eq:casebghd0lengthrestrictions3}}\ }}\frac{4\times40.8r}{R\sin^2{\phistar}}=\frac{163.2r}{R\sin^2{\phistar}}\overbracket{<}^{\mathclap{\text{\eqref{eqJZHypCond}: }R>\frac{30000r}{\sin^2{\phistar}}}}0.00544,\]
and \[\frac{kg}{2d_1}(2-\frac{1}{k})^2=\frac{\overbracket{g}^{\mathclap{\text{\eqref{eq:casebghd0lengthrestrictions1}: }g<0.0003r\sin{\phistar}}}}{\underbracket{d_1}_{\mathllap{\text{(Case 1): } d_1\ge\frac{21}{32}r\sin\phistar}}}\times\underbracket{(2k+\frac{1}{2k}-2)}_{\mathclap{\text{\quad increases with }k\in[1,\frac{4}{3}]}}<\frac{0.0003r\sin{\phistar}}{\frac{21}{32}r\sin{\phistar}}\times(2\times\frac{4}{3}+\frac{1}{2\times\frac{4}{3}}-2)<0.0005.\]
\begin{figure}[h]
\begin{minipage}{.45\textwidth}
\begin{center}
    \begin{tikzpicture}[xscale=0.9,yscale=0.9]
        \pgfmathsetmacro{\PHISTAR}{3};
        \pgfmathsetmacro{\Dlt}{0.9};
        \pgfmathsetmacro{\WIDTH}{5};
        \pgfmathsetmacro{\boudaryShift}{1.0};
        \pgfmathsetmacro{\LShift}{1.25};
        \tkzDefPoint(-\PHISTAR,0){LL};
        \tkzDefPoint(-\PHISTAR,\WIDTH){LU};
        \tkzDefPoint(\PHISTAR,0){RL};
        \tkzDefPoint(\PHISTAR,\WIDTH){RU};
        \tkzDefPoint(-\PHISTAR,\PHISTAR){X};
        \tkzDefPoint(\PHISTAR,\PHISTAR){IY};
        \tkzDefPoint(-\PHISTAR,\PHISTAR+\Dlt){LLD};
        \tkzDefPoint(\PHISTAR,\PHISTAR+\Dlt){RLD};
        \tkzDefPoint(0,0.5*\PHISTAR){LMD};
        \draw [ultra thick,red] (X) --(LLD);
        \draw [ultra thin,dashed] (LU) --(RU);
        \draw [ultra thin,dashed] (LU) --(LLD);
        \draw [ultra thin,dashed] (X) --(LL);
        \draw [ultra thin,dashed] (RU) --(RLD);
        \draw [ultra thin,dashed] (IY) --(RL);
        \draw [ultra thick,red] (IY) --(RLD);
        \draw [ultra thin,dashed] (LL) --(RL);
        \draw [red,ultra thick] (LMD) --(IY);
        \draw [red,ultra thick]  (X)--(LMD);
        \draw [red,ultra thick]  (LLD)--(RLD);
        \coordinate (LMDMark) at ([shift={(\PHISTAR,0)}]LMD);
        \coordinate (RBoundary) at ([shift={(-2*\boudaryShift,-\boudaryShift)}]IY);
        \coordinate (LBoundary) at ([shift={(2*\boudaryShift,-\boudaryShift)}]X);
         \coordinate (LlineBoundary) at ([shift={(2*\LShift,-\LShift)}]X);
         \coordinate (RlineBoundary) at ([shift={(-2*\LShift,-\LShift)}]IY);
        \coordinate (RlineBoundaryMark) at ([shift={(0,-\LShift)}]X);
        \draw [ultra thick] (LBoundary) --(RBoundary);
        \draw [thick,blue] (LlineBoundary) --(RlineBoundary);
        \draw [ultra thin,blue,dashed] (LlineBoundary) --(RlineBoundaryMark);
        \coordinate (RBoundaryMark) at ([shift={(0,-\boudaryShift)}]IY);
        \draw [ultra thin,dashed]  (LMD)--(LMDMark);
        \draw [ultra thin,dashed] (RBoundary)--(RBoundaryMark);
        \fill [red, opacity=8/30] (LBoundary) -- (RBoundary) -- (LMD) --cycle;
        \tkzDrawPoints(IY,RLD,LMD,RBoundaryMark,LMDMark,RlineBoundaryMark,RlineBoundary,LlineBoundary);
        \tkzLabelPoint[below right](RlineBoundaryMark){Line $\mathfrak{l}$};
        \tkzLabelPoint[above left](RBoundary){\small $d_1=\frac{21}{32}r\sin{\phistar}$};
        \tkzLabelPoint[below left](LMDMark){\small $d_1=\frac{r\sin{\phistar}}{2\cos{\frac{\Phistar}{2}}},\theta_1=\frac{1}{2}\Phistar$};
        \tkzLabelPoint[below](LL){$-\Phistar$};
        \tkzLabelPoint[below](RL){$\Phistar$}
         \tkzLabelPoint[left](LlineBoundary){$\mathcal{W}$};
        \tkzLabelPoint[right](RlineBoundary){$\mathcal{E}$};
    \end{tikzpicture}
    \captionsetup{width=0.85\textwidth}\caption{
    \textbf{Case 2}: $x_1\in\mathfrak{l}\subset\mathcal{L}$ in \cref{fig:A1contractionCase} 
    with $d_1<\frac{21}{32}r\sin{\phistar}$.
    The blue line segment $\overline{\mathcal{W}\mathcal{E}}=\mathfrak{l}\dfn\Big\{(\Phi_1,\theta_1=\sin^{-1}(\frac{\const}{R}))\bigm|\Phi_1\in\big[\Phistar-2\theta_1,2\theta_1-\Phistar\big],\frac{r\sin\phistar}{2}<\frac{r\sin{\phistar}}{2\cos{\frac{\Phistar}{2}}}<d_1=\const<\frac{21r\sin{\phistar}}{32}\Big\}$ has endpoints $\mathcal{W}$, $\mathcal{E}$ (``west'' and ``east'').
    }
\label{fig:b_contractionCase2}
\end{center}
\end{minipage}\hfill\vrule
\begin{minipage}{.5\textwidth}
\begin{center}
\begin{tikzpicture}[xscale=0.65,yscale=.65]
    \tkzDefPoint(0,0){Or};
    \tkzDefPoint(0,-4.8){Y};
    \pgfmathsetmacro{\rradius}{4.5};
    \pgfmathsetmacro{\bdist}{6.3};
    \tkzDefPoint(0,\bdist){OR};
    \pgfmathsetmacro{\phistardeg}{40};
    \pgfmathsetmacro{\XValueArc}{\rradius*sin(\phistardeg)};
    \pgfmathsetmacro{\YValueArc}{\rradius*cos(\phistardeg)};
    \tkzDefPoint(\XValueArc,-1.0*\YValueArc){B};    \tkzDefPoint(-1.0*\XValueArc,-1.0*\YValueArc){A};
    \pgfmathsetmacro{\Rradius}{veclen(\XValueArc,\YValueArc+\bdist)};
    \tkzDrawArc[name path=Cr,very thick](Or,B)(A);
    \tkzDrawArc[name path=CR,very thick](OR,A)(B);
    \pgfmathsetmacro{\PXValue}{\Rradius*sin(7)};
    \pgfmathsetmacro{\PYValue}{\bdist-\Rradius*cos(7)};
    \tkzDefPoint(\PXValue,\PYValue){P};
    \tkzDefPoint(0,\bdist-\Rradius){C};
    \draw[blue,ultra thin] (P)--(B);
    \draw[red,ultra thin,dashed](OR)--(P);
    \draw[red,ultra thin,dashed](OR)--(A);
    \draw[red,ultra thin,dashed](OR)--(B);
    \tkzDefPointBy[projection=onto OR--P](A)\tkzGetPoint{M};
    \coordinate (Q_far) at ($(A)!2!(M)$);
    \path[name path=linePQ] (P)--(Q_far);
    \path[name intersections={of=Cr and linePQ,by={Q}}];
    \tkzDrawPoints(A,B,P,Or,OR,Q);
    \begin{scope}[ultra thin,decoration={
    markings,
    mark=at position 0.5 with {\arrow{>}}}
    ] 
    \draw[blue,postaction={decorate}] (A)--(P);
    \draw[blue,postaction={decorate}] (P)--(Q);
    \end{scope}
    \tkzLabelPoint[below,rotate=0](P){$p(\mathcal{E})=P_{\mathcal{E}}=(x_\mathcal{E},y_\mathcal{E})$};
    \tkzLabelPoint[left,rotate=0](A){$A$};
    \tkzLabelPoint[right,rotate=0](Q){$Q_\mathcal{E}$};
    \tkzLabelPoint[below,rotate=0](B){$B$};
    \tkzLabelPoint[left,rotate=0](Or){$O_r$};
    \tkzLabelPoint[left,rotate=0](OR){$O_R$};
    \tkzMarkAngle[arc=lll,ultra thin,size=1.2](A,OR,P);
    \tkzMarkAngle[arc=lll,ultra thin,size=2.5](P,OR,B);
    \node[rotate=0] at ($(OR)+(-95:2.2)$)  {\small $2\theta_1$};
    \node[rotate=0] at ($(OR)+(-80:4)$)  {\small $2\Phistar-2\theta_1$};
\end{tikzpicture}
\captionsetup{width=0.85\textwidth}\caption{If $x_1=\mathcal{E}$, then $p(x_0)$ is at corner $A$, $|P_{\mathcal{E}}Q_{\mathcal{E}}|=\tilde{\tau}_1$ and $\tilde{\tau}_0=|AP_{\mathcal{E}}|=2d_1=2\cdot\const$}\label{fig:Fraction_tau1_d0_min}
\end{center}
\end{minipage}
\end{figure}
\hidesubsection{Case 2} $x_1\in\mathcal{L}$ in \cref{fig:A1contractionCase} and $d_1<\frac{21}{32}r\sin\phistar$ (see \cref{fig:b_contractionCase2}). Note that in this case $\Phistar/2<\theta_1<\Phistar$

Suppose $x_1=(\Phi_1,\theta_1)$ is on the line segment \[\mathfrak{l}\dfn\Big\{(\Phi_1,\theta_1=\sin^{-1}{(\frac{\const}{R})})\bigm|\Phi_1\in\big(\Phistar-2\theta_1,2\theta_1-\Phistar\big),\frac{r\sin{\phistar}}{2\cos{\frac{\Phistar}{2}}}<d_1=\const<\frac{21r\sin{\phistar}}{32}\Big\}\]in region $\mathcal{L}$. $\mathfrak{l}$ has a left end point $\mathcal{W}=(2\theta_1-\Phistar,\theta_1)$ and a right end point $\mathcal{E}=(\Phistar-2\theta_1,\theta_1)$. Note that $\mathcal{W}$, $\mathcal{E}$ are singular points, so $x_1$ cannot be $\mathcal{W}$ or $\mathcal{E}$.

Under the standing assumptions for this section, in the coordinate system given in \cref{def:StandardCoordinateTable} (\cref{fig:CaseAMonotone,fig:Fraction_tau1_d0_min}), $x_0=(\phi_0,\theta_0)$, $x_2=(\phi_2,\theta_2)$, $P=p(x_1)=(x_P,y_P)$, $Q=p(x_2)=(x_Q,y_Q)=(r\sin\phi_2,-r\cos\phi_2)$ and $T_0=p(x_0)=(x_{T_0},y_{T_0})=(r\sin\phi_0,-r\cos\phi_0)$.

Since we consider $\tau_0$, $\tau_1$, $d_1$, $d_2$ and $d_0$ as functions of $x_1=(\Phi_1,\theta_1)$ on line $\mathfrak{l}$, by \cref{monotone_p3}\eqref{itemmonotone_p3-1}\eqref{itemmonotone_p3-3}, we know 
\begin{align*}
    \frac{d\tau_1}{d\Phi_1}&=-\frac{bx_Q}{d_2},\\
    \frac{d\tau_0}{d\Phi_1}&=-\frac{bx_{T_0}}{d_0}.
\end{align*}
\COMMENT{\textcolor{red}{I have fixed a number of like errors before: in a context like this, the comma must be inside the display, not after it. Likewise with periods, etc. Please fix any occurrence that you see.\newline\wentao{Work in progress\wentao{Now replaced by new proposition}}}}

By \cref{proposition:NinsizeInA1BCcases}, \eqref{eq:observation2(a1)(b)(c)region}, we have $0<\phi_2-\phistar<2.001\sqrt{4r/R}$ and by \cref{proposition:NoutsizeInA1BCcases}\eqref{item01NoutsizeInA1BCcases} and we have $|\phi_0-(2\pi-\phistar)|<2.001\sqrt{4r/R}$. Therefore, we further get the following.\COMMENT{Fix punctuation. \Large Throughout!\wentao{better now. Work in progress.}}
\begin{align*}
    \phistar<&\phi_2<\phistar+2.001\sqrt{4r/R}\overbracket{<}^{\mathllap{\text{\eqref{eqJZHypCond0}: }R>1700r}}\phistar+0.0971<\tan^{-1}{(1/3)}+0.0971<0.42<\pi/2,\\
    2\pi-\phistar>&\phi_0>2\pi-\phistar-2.001\sqrt{4r/R}\overbracket{>}^{\mathllap{\text{\eqref{eqJZHypCond0}: }R>1700r}}2\pi-\phistar-0.0971>2\pi-\tan^{-1}{(1/3)}-0.0971>5.86>3\cdot(\pi/2).
\end{align*} Therefore, with $x_Q=r\sin\phi_2>0$ and $x_{T_0}=r\sin\phi_0<0$, $\tau_1$ is a monotonically decreasing function of $\Phi_1\in(2\theta_1-\Phistar,\Phistar-2\theta_1)$ and $\tau_0$ is a monotonically increasing function of $\Phi_1\in(2\theta_1-\Phistar,\Phistar-2\theta_1)$. 

On the other hand, by \cref{monotone_p4} $d_0$ is a monotonically increasing function of $\Phi_1\in(2\theta_1-\Phistar,\Phistar-2\theta_1)$. Note that since $\mathfrak{l}\subset\MRin\cap\MRout$, the length functions in \cref{def:fpclf} are continuous on $\mathfrak{l}$. With $\theta_1$ fixed, we denote the length function limits at $\mathcal{E}$ by $\tilde{\tau}_0$, $\tilde{\tau}_1$, $\tilde{d}_0$ in the following.
\begin{equation}\label{eq:lengthfunc_Elimit}
\begin{aligned}
\tilde{\tau}_0&\dfn\lim_{\mathfrak{l}\ni x_1\to \mathcal{E}}\tau_0(x_1)=\lim_{\substack{x_1=(\Phi_1,\theta_1)\in\mathfrak{l}\\\Phi_1\to\Phistar-2\theta_1}}\tau_0,\\
    \tilde{\tau}_1&\dfn\lim_{\mathfrak{l}\ni x_1\to \mathcal{E}}\tau_1(x_1)=\lim_{\substack{x_1=(\Phi_1,\theta_1)\in\mathfrak{l}\\\Phi_1\to\Phistar-2\theta_1}}\tau_1,\\
    \tilde{d}_0&\dfn\lim_{\mathfrak{l}\ni x_1\to \mathcal{E}}d_0(x_1)=\lim_{\substack{x_1=(\Phi_1,\theta_1)\in\mathfrak{l}\\\Phi_1\to\Phistar-2\theta_1}}d_0.
\end{aligned}
\end{equation}
With $\theta_1$ fixed, $\frac{\tau_1}{d_0}$ is a decreasing function of $\Phi_1\in(2\theta_1-\Phistar,\Phistar-2\theta_1)$, therefore $\frac{\tau_1}{d_0}\ge\lim_{\mathfrak{l}\ni x_1\to \mathcal{E}}\frac{\tau_1}{d_0}(x_1)=\frac{\tilde{\tau}_1}{\tilde{d}_0}$. $\tilde{d}_0$ is the limit of $d_0(x_1)$ with $\mathfrak{l}\ni x_1\to\mathcal{E}$. Hence by \cref{proposition:casebhandgfunctions}\eqref{eq:casebghd0lengthrestrictions2}, $\tilde{d}_0$ and $\frac{\tilde{\tau}_1}{\tilde{d}_0}$ satisfy the following. \begin{equation}\label{eq:tau_d0_frac_min}
\begin{aligned}
    &0.9997r\sin\phistar\overbracket{\le}^{\text{\eqref{eq:casebghd0lengthrestrictions2}}}\tilde{d}_0\overbracket{\le}^{\text{\eqref{eq:casebghd0lengthrestrictions2}}}1.0003r\sin\phistar\\
    &\frac{\tau_1}{d_0}(x_1)\ge\frac{\tilde{\tau}_1}{\tilde{d}_0}\ge\frac{\tilde{\tau}_1}{1.0003r\sin\phistar}\\
\end{aligned}
\end{equation}

We observe that $\mathcal{E}$ in phase space represents the collision on $\Gamma_R$ with the billiard trajectory coming from corner $A$ (see \cref{fig:Fraction_tau1_d0_min}). In \cref{fig:Fraction_tau1_d0_min} with cartesian coordinates, 
let 
\[P_\mathcal{E}\dfn p(\mathcal{E})=\big(x_\mathcal{E},y_\mathcal{E}\big)=\big(R\sin{(\Phistar-2\theta_1)},b-R\cos{(\Phistar-2\theta_1)}\big).\]  We can observe that for any $x_1=(\Phi_1,\theta_1)\in\mathfrak{l}$, that is, $P=p(x_1)=\big(x_P,y_P\big)=\big(R\sin\Phi_1,b-R\cos\Phi_1\big)$ with $2\theta_1-\Phistar<\Phi_1<\Phistar-2\theta_1$. Therefore, $x_P<x_{\mathcal{E}}$, $y_P<y_{\mathcal{E}}$. Moreover with corner $B=(x_B,y_B)$ in \cref{fig:Fraction_tau1_d0_min}, since $(x_Q,y_Q)=(r\sin\phi_2,-r\cos\phi_2)$ with $0<\phi_2<\pi/2$, we have 
\begin{align*}
x_P<&x_{\mathcal{E}}<x_{B}=r\sin\phistar<x_Q\\
y_P<&y_{\mathcal{E}}<y_{B}=-r\cos\phistar<y_Q
\end{align*}
Hence $\forall x_1\in\mathfrak{l}$, $\tau_1(x_1)=|PQ|=\sqrt{(x_Q-x_P)^2+(y_Q-y_P)^2}>\sqrt{(x_B-x_\mathcal{E})^2+(y_B-y_\mathcal{E})^2}=|P_\mathcal{E}B|$. Since $\tau_1(x_1)$ is a continuous function for $x_1\in\mathfrak{l}$, as a limit \[\tilde{\tau}_1=\lim_{\mathfrak{l}\ni x_1\to \mathcal{E}}\tau_1(x_1)\ge|P_\mathcal{E}B|.\]

We now compute a lower bound of $|P_\mathcal{E}B|$.

In \cref{fig:Fraction_tau1_d0_min}, $\measuredangle{AO_RP_{\mathcal{E}}}=2\theta_1$, $\measuredangle{P_{\mathcal{E}}O_RB}=2\Phistar-2\theta_1$. Using trigonometry, we get the following.

\begin{equation}\label{eq:PBLowerBound}
\begin{aligned}
 |P_{\mathcal{E}}B|&=2R\sin{(\Phistar-\theta_1)}\overbracket{=}^{\mathclap{\sin{(\Phistar-\theta_1)}=\sin{\Phistar}\cos{\theta_1}-\cos{\Phistar}\sin{\theta_1}}}2R\sin{\Phistar}\cos{\theta_1}-2R\cos{\Phistar}\sin{\theta_1}\overbracket=^{\mathclap{r\sin\phistar=R\sin\Phistar}}2r\sin\phistar\cdot\cos\theta_1-2R\sin\theta_1\cos\Phistar\\
 \overbracket=^{\quad\mathllap{d_1=R\sin\Phistar}}&2r\sin\phistar\cos\theta_1-2d_1\cos\Phistar\overbracket>^{\mathclap{\Phistar/2<\theta_1<\Phistar}}2r\sin\phistar\cos\Phistar-2d_1\cos\Phistar=2r\sin\phistar\cos\Phistar(1-\frac{d_1}{r\sin\phistar})\\
 \overbracket>^{\quad\mathllap{d_1<\frac{21}{32}r\sin\phistar}}&2r\sin\phistar\cos\Phistar(1-\frac{21}{32})=\frac{11}{16}r\sin\phistar\sqrt{1-\sin^2\Phistar}\overbracket=^{r\sin\phistar=R\sin\Phistar}\frac{11}{16}r\sin\phistar\sqrt{1-\frac{r^2\sin^2\phistar}{R^2}}.
 \end{aligned}
 \end{equation}
\eqref{eq:tau_d0_frac_min} and \eqref{eq:PBLowerBound} imply several useful inequalities:
\begin{equation}\label{eq:tau1d0lowerbound5over8}
     \frac{\tau_1}{d_0}\overbracket\ge^{\mathclap{\strut\eqref{eq:tau_d0_frac_min}}}\frac{\tilde{\tau}_1}{\tilde{d}_0}\overbracket\ge^{\mathclap{\strut\eqref{eq:tau_d0_frac_min}}}\frac{\tilde{\tau}_1}{1.0003r\sin\phistar}\overbracket{>}^{\mathclap{\substack{\tilde{\tau}_1\ge |P_{\mathcal{E}}B|\\ \text{and \eqref{eq:PBLowerBound}} }}}\frac{11r\sin{\phistar}}{16\times 1.0003r\sin{\phistar}}\sqrt{1-\frac{r^2\sin^2\phistar}{R^2}}\overbracket>^{\mathrlap{\eqref{eqJZHypCond}:R>1700r}\quad}0.6872>\frac{5}{8}.
 \end{equation}
On the other hand, if $\mathfrak{l}\ni x_1\to \mathcal{E}$, then $p(x_0)\to A$, $p(x_1)\to P_{\mathcal{E}}$. Therefore,
\[\tilde{\tau}_0\overbracket{=}^{\mathclap{\eqref{eq:lengthfunc_Elimit}}}\lim_{\mathfrak{l}\ni x_1\to \mathcal{E}}\tau_1(x_1)\overbracket{=}^{\mathclap{\substack{\tau_1=\,\text{distance between }\\p(x_0),p(x_1)}}}\lim_{\mathfrak{l}\ni x_1\to \mathcal{E}}|p(x_0)p(x_1)|=|AP_{\mathcal{E}}|\overbracket{=}^{\mathclap{A,P_{\mathcal{E}}\text{ are on }\Gamma_R}}2d_1.\]

We extend the definitions of $g$, $h$ (\cref{def:casebhandgfunctions}) to $\tilde{h}$, $\tilde{g}$ on $\mathcal{E}$ by taking limits.\COMMENT{The first two equations in \eqref{eq:casebcase2ghlimitatE} are backwards; they should start with the thing being defined.}
\begin{equation}\label{eq:casebcase2ghlimitatE}
    \begin{aligned}
    \frac{\tilde{\tau}_0}{\tilde{d}_0}+\frac{\tilde{\tau}_1}{\tilde{d}_0}-2&\overbracket{=}^{\mathclap{\eqref{eq:lengthfunc_Elimit}}}\lim_{\mathfrak{l}\ni x_1\to \mathcal{E}}{\big(\frac{\tau_0}{d_0}+\frac{\tau_1}{d_0}-2\big)}\overbracket{=}^{\mathclap{\text{\cref{def:casebhandgfunctions}}}}\lim_{\mathfrak{l}\ni x_1\to \mathcal{E}}{h}\nfd\tilde{h},\\
    \tilde{d}_0-r\sin{\phistar}&\overbracket{=}^{\mathclap{\eqref{eq:lengthfunc_Elimit}}}\lim_{\mathfrak{l}\ni x_1\to \mathcal{E}}{\big(d_0-r\sin\phistar\big)}\overbracket{=}^{\mathclap{\text{\cref{def:casebhandgfunctions}}}}\lim_{\mathfrak{l}\ni x_1\to \mathcal{E}}{g}\nfd\tilde{g},\\
    \text{Therefore, }-\tilde{h}&=2-\frac{\tilde{\tau}_0}{\tilde{d}_0}-\frac{\tilde{\tau}_1}{\tilde{d}_0}=2-\frac{2d_1}{\tilde{d}_0}-\frac{\tilde{\tau}_1}{\tilde{d}_0},\\
    \frac{\tilde{\tau}_1}{\tilde{d}_0}&=2-\frac{2d_1}{\tilde{d}_0}+\tilde{h}.
\end{aligned}
\end{equation}
Since $g$, $h$ are continuous functions of $x_1$ on $\mathfrak{l}$, by \cref{proposition:casebhandgfunctions}\eqref{eq:casebghd0lengthrestrictions1}\eqref{eq:casebghd0lengthrestrictions3}, $\tilde{g}$, $\tilde{h}$ as limits of $x_1\to\mathcal{E}$ satisfy the following.
\begin{equation}\label{eq:tildeghfunctionsrestriction}
    \begin{aligned}
\tilde{g}\overbracket{\in}^{\mathclap{\eqref{eq:casebghd0lengthrestrictions1}}}\big[-\frac{8.3 r^2}{R\sin\phistar},\frac{8.3r^2}{R\sin\phistar}\big]\overbracket{\subset}^{\mathclap{\text{\eqref{eqJZHypCond}: }R>\frac{30000r}{\sin^2{\phistar}}}}\big[-0.0003r\sin\phistar,0.0003r\sin\phistar\big],\qquad
0\overbracket{\le}^{\eqref{eq:casebghd0lengthrestrictions3}}\tilde{h}\overbracket{\le}^{\eqref{eq:casebghd0lengthrestrictions3}}\frac{40.8r}{R\sin^2\phistar}.
\end{aligned}
\end{equation}
For $k\in[1,\frac{3}{4}]$, the quadratic polynomial $\mathfrak{Q}(t)\dfn t(2-\frac{1}{k}-t)$\COMMENT{It is NOT necessary to use different funny Q's for the various quadratic polynomials; it is clear that each is "disposable" and will only be used briefly---so there can be no confusion.} is a decreasing function of $t\ge\frac{5}{8}\ge1-\frac{1}{2k}$, so 
\begin{equation}\label{eq:IkUpperBoundCasebCase2}
\begin{aligned}
    \mathcal{I}(k)\overbracket{=}^{\mathclap{\eqref{eq:CaseBIkQuadraticPolynomial}}}&-2k+1+\frac{2d_0}{d_1}k\times\underbracket{\frac{\tau_1}{d_0}(2-\frac{1}{k}-\frac{\tau_1}{d_0})}_{=\mathfrak{Q}(\frac{\tau_1}{d_0})}+\frac{2k\tau_1}{d_1}h-hk\\ \overbracket<^{\mathclap{\eqref{eq:tau1d0lowerbound5over8}: \frac{\tau_1}{d_0}\ge\frac{\tilde{\tau}_1}{\tilde{d}_0}>\frac{5}{8}}}&-2k+1+\frac{2d_0}{d_1}k\times\frac{\tilde{\tau}_1}{\tilde{d}_0}(2-\frac{1}{k}-\frac{\tilde{\tau}_1}{\tilde{d}_0})+\frac{2k\tau_1}{d_1}h-hk\\
    \overbracket{=}^{\mathclap{\text{\eqref{eq:casebcase2ghlimitatE}: }\frac{\tilde{\tau}_1}{\tilde{d}_0}=2-\frac{2d_1}{\tilde{d}_0}+\tilde{h}}}&-2k+1+\frac{2d_0}{d_1}k\times(2-\frac{2d_1}{\tilde{d_0}}+\tilde{h})(-\tilde{h}+\frac{2d_1}{\tilde{d_0}}-\frac{1}{k})+\frac{2k\tau_1}{d_1}h-hk\\
    &=-2k+1+\underbracket{\frac{2d_0}{d_1}k(2-\frac{2d_1}{\tilde{d_0}})(\frac{2d_1}{\tilde{d_0}}-\frac{1}{k})}_{\mathrlap{=\frac{2\tilde{d_0}}{d_1}k(2-\frac{2d_1}{\tilde{d_0}})(\frac{2d_1}{\tilde{d_0}}-\frac{1}{k})+2\frac{d_0-\tilde{d_0}}{d_1}k(2-\frac{2d_1}{\tilde{d_0}})(\frac{2d_1}{\tilde{d_0}}-\frac{1}{k})}}+\tilde{h}\frac{2d_0}{d_1}k(\frac{2d_1}{\tilde{d_0}}-\frac{1}{k})-\frac{2d_0}{d_1}\tilde{h}k(2-\frac{2d_1}{\tilde{d_0}})-\frac{2d_0}{d_1}k(\tilde{h})^2
    +h(\frac{2k\tau_1}{d_1}-k)\\
    &=\underbracket{-2k+1+\frac{2\tilde{d_0}}{d_1}k(2-\frac{2d_1}{\tilde{d_0}})(\frac{2d_1}{\tilde{d_0}}-\frac{1}{k})}_{<-0.0639\text{ by \eqref{eq:Ikfirst3items} below}}\\&\underbracket{+\tilde{h}\frac{2d_0}{d_1}k(\frac{2d_1}{\tilde{d_0}}-\frac{1}{k})-\frac{2d_0}{d_1}\tilde{h}k(2-\frac{2d_1}{\tilde{d_0}})-\frac{2d_0}{d_1}k(\tilde{h})^2+h(\frac{2k\tau_1}{d_1}-k)
    +2\frac{d_0-\tilde{d_0}}{d_1}k(2-\frac{2d_1}{\tilde{d_0}})(\frac{2d_1}{\tilde{d_0}}-\frac{1}{k})}_{\leq0\text{ when }R\text{ is large}}
\end{aligned}
\end{equation}
To bound this, we begin with the sum $-2k+1+\frac{2\tilde{d_0}}{d_1}k(2-\frac{2d_1}{\tilde{d_0}})(\frac{2d_1}{\tilde{d_0}}-\frac{1}{k})$ of the first three terms on the right-hand side of \eqref{eq:IkUpperBoundCasebCase2}.
Let $\mathfrak{e}=\frac{d_1}{\tilde{d}_0}>0$, which is close to $\frac{d_1}{r\sin{\phistar}}=\frac{\const}{r\sin\phistar}\in(\frac{1}{2},\frac{21}{32})$ in this case. Then 
\begin{equation}\label{eq:Ikfirst3items}
\begin{aligned}
    -2k+1+\frac{2\tilde{d_0}}{d_1}k(2-\frac{2d_1}{\tilde{d_0}})(\frac{2d_1}{\tilde{d_0}}-\frac{1}{k})&=-2k+1+\frac{2k}{\mathfrak{e}}(2-2\mathfrak{e})(2\mathfrak{e}-\frac{1}{k})=-2k+1+2k\frac{1}{\mathfrak{e}}(-4\mathfrak{e}^2+\frac{2\mathfrak{e}}{k}+4\mathfrak{e}-\frac{2}{k})\\
    &=-2k+1-8k\mathfrak{e}+4+8k-\frac{4}{\mathfrak{e}}=6k+5-\frac{4}{\mathfrak{e}}-8k\mathfrak{e}\\
    \underbracket{\le}_{\mathllap{8ke+\frac{4}{e}\ge4\times2\sqrt{2ke\times\frac{1}{e}}=8\sqrt{2k}}}&6k+5-8\sqrt{2k}\underbracket{=}_{\mathclap{\text{Let }\mu=\sqrt{2k}\in\big[\sqrt{2},\sqrt{8/3}\big]}}3\mu^2-8\mu+5\\
    \le& 3\times\frac{8}{3}-8\times\sqrt{\frac{8}{3}}+5=13-8\times\sqrt{\frac{8}{3}}<-0.0639.
\end{aligned}
\end{equation}
We now show that the last five terms of the inequality \eqref{eq:IkUpperBoundCasebCase2} can be arbitrarily small when $R$ is large. 
We reorganize this sum of five terms as follows \eqref{eq:IkUpperBoundCasebCase2MinorTerms}.
\begin{equation}\label{eq:IkUpperBoundCasebCase2MinorTerms}
\begin{aligned}
&+\tilde{h}\frac{2d_0}{d_1}k(\frac{2d_1}{\tilde{d_0}}-\frac{1}{k})-\frac{2d_0}{d_1}\tilde{h}k(2-\frac{2d_1}{\tilde{d_0}})-\frac{2d_0}{d_1}k(\tilde{h})^2+h(\frac{2k\tau_1}{d_1}-k)+2\frac{d_0-\tilde{d_0}}{d_1}k(2-\frac{2d_1}{\tilde{d_0}})(\frac{2d_1}{\tilde{d_0}}-\frac{1}{k})\\
&=\tilde{h}k\frac{4d_0}{\tilde{d}_0}-\tilde{h}\frac{2d_0}{d_1}-\frac{4d_0}{d_1}\tilde{h}k+\frac{4d_0}{\tilde{d}_0}\tilde{h}k-\frac{2d_0}{d_1}k(\tilde{h})^2+h(\frac{2k\tau_1}{d_1}-k)+2\frac{d_0-\tilde{d_0}}{d_1}k(2-\frac{2d_1}{\tilde{d_0}})(\frac{2d_1}{\tilde{d_0}}-\frac{1}{k})\\
&=\underbracket{-\tilde{h}d_0\big[k(-\frac{8}{\tilde{d}_0}+\frac{4}{d_1})+\frac{2}{d_1}\big]}_{\mathclap{\eqref{eq:casebcase2minorterm1}:<\frac{63r}{R\sin^2\phistar}}}\underbracket{-\frac{2d_0}{d_1}k(\tilde{h})^2}_{\le0}+\underbracket{hk(\frac{2\tau_1}{d_1}-1)}_{\mathclap{\eqref{eq:casebcase2minorterm2}:<\frac{164r}{R\sin^2\phistar}}}+\underbracket{2\frac{d_0-\tilde{d_0}}{d_1}(2-\frac{2d_1}{\tilde{d_0}})k(\frac{2d_1}{\tilde{d_0}}-\frac{1}{k})}_{\mathclap{\text{\eqref{eq:casebcase2minorterm3},\eqref{eq:casebcase2minorterm4}: }<\frac{50r}{R\sin^2\phistar}}}
\end{aligned}
\end{equation}
For $k\in[1,\frac{4}{3}]$, we have the following upper bounds for these four terms: $\displaystyle-\frac{2d_0}{d_1}k(\tilde{h})^2\le0$,
\begin{equation}\label{eq:casebcase2minorterm1}
\begin{aligned}\underbracket{-\tilde{h}d_0}_{\mathllap{\eqref{eq:tildeghfunctionsrestriction}:\le0}}\Big[\overbracket{k}^{\mathrlap{\ge1}}(\frac{-8}{\underbracket{\tilde{d}_0}_{\qquad\mathllap{\substack{\eqref{eq:tau_d0_frac_min}:\,0.9997r\sin\phistar\le\tilde{d}_0\\\le1.0003r\sin{\phistar}}}}}+\underbracket{\frac{4}{d_1})+\frac{2}{d_1}}_{\mathrlap{\text{case (2) and \eqref{eq:lengthrestrictions1}: }\frac{21r\sin{\phistar}}{32}>d_1>0.5r\sin\phistar}}\Big]\overbracket{<}^{\mathrlap{1\le k\le4/3}}&-\tilde{h}d_0\Big(\frac{-8k}{\tilde{d}_0}+\frac{4k+2}{d_1}\Big)<-\tilde{h}d_0\big[-8\frac{4}{3\tilde{d}_0}+\frac{4\cdot1+2}{(21/32)r\sin\phistar}\big]\\
<&-\tilde{h}d_0\Big(-\frac{32}{3\cdot0.9997r\sin\phistar}+\frac{6}{(21/32)r\sin\phistar}\Big)\\
<&\frac{-1.521\tilde{h}d_0}{r\sin\phistar}\\
\overbracket{<}^{\mathllap{\substack{\eqref{eq:casebghd0lengthrestrictions2}:\,0.9997r\sin{\phistar}<d_0\\<1.0003r\sin{\phistar}}}}&-1.521\cdot0.9997\tilde{h}\overbracket{<}^{\mathclap{\eqref{eq:tildeghfunctionsrestriction}:\tilde{h}\le\frac{40.8r}{R\sin^2\phistar}}}\qquad\frac{63r}{R\sin^2\phistar},\end{aligned}
\end{equation}
and
\begin{equation}\label{eq:casebcase2minorterm2}
\overbracket{h}^{\mathllap{\eqref{eq:casebghd0lengthrestrictions3}:h>0}}\underbracket{k}_{\mathclap{k\ge1}}(\frac{2\tau_1}{d_1}-1)\overbracket<^{\mathclap{\substack{\text{\eqref{eq:lengthrestrictions2}}:\\\tau_1<2d_1\quad}}}3hk\underbracket{\le}_{\mathclap{h>0,k\le4/3}}4h\overbracket{<}^{\mathclap{{\textstyle\strut}\eqref{eq:casebghd0lengthrestrictions3}:\, h<\frac{40.8r}{R\sin^2{\phistar}}}}\qquad\frac{164r}{R\sin^2\phistar}.
\end{equation}
%
From case 2 condition $d_1<\frac{21}{32}r\sin\phistar$, \eqref{eq:lengthrestrictions2}: $0.5r\sin\phistar<d_1$, \eqref{eq:tildeghfunctionsrestriction}, \eqref{eq:tau_d0_frac_min}, \cref{proposition:casebhandgfunctions}\eqref{eq:casebghd0lengthrestrictions2}, \eqref{eq:casebcase2ghlimitatE} and \eqref{eq:tildeghfunctionsrestriction}, we obtain the following.
\begin{equation}\label{eq:casebcase2minorterm3}
\begin{aligned}
k(\frac{2d_1}{\tilde{d_0}}-\frac{1}{k})=\frac{2d_1k}{\tilde{d_0}}-1&\overbracket{\in}^{\mathclap{\substack{0.5r\sin\phistar<d_1<(21/32)r\sin\phistar\\ \eqref{eq:tau_d0_frac_min}:\,0.9997r\sin\phistar<\tilde{d}_0<1.0003r\sin\phistar\\1\le k\le4/3}}}(2\cdot\frac{0.5r\sin{\phistar}}{1.0003r\sin{\phistar}}-1,2\cdot\frac{(21/32)r\sin{\phistar}}{0.9997r\sin{\phistar}}\frac{4}{3}-1)\subset(-0.003,0.751),\\
2-\frac{2d_1}{\tilde{d_0}}&\overbracket{\in}^{\mathclap{\substack{0.5r\sin\phistar<d_1<(21/32)r\sin\phistar\\ \eqref{eq:tau_d0_frac_min}:\,0.9997r\sin\phistar<\tilde{d}_0<1.0003r\sin{\phistar}\\|}}}
(2-2\cdot\frac{(21/32)r\sin{\phistar}}{0.9997r\sin{\phistar}},2-2\times\frac{0.5\times r\sin{\phistar}}{1.0003r\sin{\phistar}})\subset(0.687,1.0003),\\
d_0-\tilde{d_0}&\overbracket{\in}^{\mathclap{\substack{\eqref{eq:casebghd0lengthrestrictions1}\eqref{eq:casebcase2ghlimitatE}:g-\tilde{g}=d_0-\tilde{d}_0\\\eqref{eq:tildeghfunctionsrestriction}:\tilde{g}\in[-\frac{8.3 r^2}{R\sin\phistar},\frac{8.3r^2}{R\sin\phistar}]\\\eqref{eq:casebghd0lengthrestrictions1}: g\in(-\frac{8.3 r^2}{R\sin\phistar},\frac{8.3r^2}{R\sin\phistar})\\|}}}(-\frac{16.6r^2}{R\sin{\phistar}}, \frac{16.6r^2}{R\sin{\phistar}}).
\end{aligned}
\end{equation} So, the conclusions of \eqref{eq:casebcase2minorterm3} consolidate into the following.
\begin{equation}\label{eq:casebcase2minorterm4}
\begin{aligned}
2\frac{d_0-\tilde{d_0}}{d_1}(2-\frac{2d_1}{\tilde{d_0}})k(\frac{2d_1}{\tilde{d_0}}-\frac{1}{k})&\le2\frac{|d_0-\tilde{d_0}|}{d_1}\bigm|(2-\frac{2d_1}{\tilde{d_0}})\bigm|\bigm|k(\frac{2d_1}{\tilde{d_0}}-\frac{1}{k})\bigm|<2\cdot\frac{\frac{16.6r^2}{R\sin{\phistar}}}{d_1}\cdot1.0003 \cdot0.751\\
&\overbracket{<}^{\mathllap{d_1>0.5r\sin\phistar}}2\cdot\frac{\frac{16.6r^2}{R\sin{\phistar}}}{0.5r\sin{\phistar}}\cdot1.0003\cdot0.751<\frac{50r}{R\sin^{2}{\phistar}}.
\end{aligned}
\end{equation}
This concludes \textbf{Case 2}, in which $d_1<\frac{21}{32}r\sin\phistar$, $k\in[1,\frac{4}{3}]$, $\mathcal{I}(k)$ in \eqref{eq:IkUpperBoundCasebCase2} satisfies the following. \[\mathcal{I}(k)\overbracket{<}^{\mathllap{\eqref{eq:Ikfirst3items},\eqref{eq:IkUpperBoundCasebCase2MinorTerms}}}-0.0639+\frac{63r}{R\sin^2\phistar}+\frac{164r}{R\sin^2\phistar}+\frac{50r}{R\sin^2\phistar}=-0.0639+\frac{277r}{R\sin^2\phistar}\overbracket<^{\mathclap{\eqref{eqJZHypCond}:R>\frac{30000r}{\sin^2\phistar}}}-0.05.\qedhere\]
\end{proof}
\section{Contraction control in case (c)
}\label{sec:ProofTECase_c}
This section is devoted to the proof of \cref{TEcase_c}, which involves the arduous work on the hardest case of those defined in \eqref{eqMainCases}. 

The first subsection defines functions which will be involved in lower bounds for the expansion in subsequent collisions, studies the return orbit segment involved in this case and basic formulas pertinent for the expansion computations, and begins estimates for the parameters involved. Later subsections then proceed according to the number $n_1$ of subsequent collisions on $\Gamma_R$  (\cref{def:MhatReturnOrbitSegment}) being at least 4, and equal to 3, 2, and 1, respectively. This becomes increasingly harder, and the last subsection alone takes up pages \pageref{SBSn1=1}--\pageref{sec:SACTWFELBWB}.

\begin{remark}[Standing assumptions]\label{REMCaseCAssumptions}
This section has as standing assumptions those of \cref{TEcase_c}: Fix $r,R$ satisfying the hyperbolicity condition \eqref{eqJZHypCond}, fix $\phistar\in(0,\tan^{-1}{(1/3)})$, $\hat F$ from \cref{def:SectionSetsDefs}, and the ``half-quadrant'' cone family $C_x\dfn\big\{(d\phi,d\theta)\bigm|\frac{d\theta}{d\phi}\in [0,1]\big\}$ from \cref{RMASUH}. We assume that $d_1<2r$ and that the orbit segment of an $x\in\hat M\aeq(\Mrin\cap\Mrout)\sqcup\Fc^{-1}(\Mrout\setminus\Mrin)$ from \eqref{def:Mhat} returning to $\hat M$ includes an orbit segment \eqref{eqPtsFromeqMhatOrbSeg} (page \pageref{eqPtsFromeqMhatOrbSeg}) with more than one collision with $\Gamma_R$.\COMMENT{Now that this is here, edit references to assumptions elsewhere}. Suppose $x$ has a return orbit segment $x$, $\Fc(x)$, $\cdots$, $\Fc^{\sigma(x)}(x)=\hat{F}(x)\in \hat{M}$ as in \cref{def:MhatReturnOrbitSegment} with $x_{0}\in\Mrout$, $x_1=\Fc(x_0)=(\phi_1,\theta_1)\in \MRin$, $x_2\in\Mrin$, $\sin\theta_1< 2r/R$, $n_1\ge1$ (i.e., case (c) of \eqref{eqMainCases}).
\end{remark}
We will prove:
\begin{enumerate}
    \item
    the cone family $C_x$ is strictly invariant under $D\hat{F}$, i.e., $D_x\hat{F}(C_x)\subset\text{\big\{interior of }C_{\hat{F}(x)}\text{\big\}}$.
    \item
    $\frac{\|D_{x}\hat{F}(dx)\|_\p}{\|dx\|_\p}>1-1743\big(r/R\big)-15450\big(r/R\big)^2>0.9$, for all $dx\in C_x$.
    \item
    when $n_1\ge4$, $\frac{\|D_x\hat{F}(dx)\|_\p}{\|dx\|_\p}>\frac{0.96}{1.01}n_1^2+5.655n_1-8.628>29.1$, for all $dx\in C_x$ (\cref{subsec:expansionLarge}).
\end{enumerate}
\subsection[A Formula for Expansion Lower Bound in \texorpdfstring{\cref{TEcase_c}}{case (c) transition}]{A candidate lower bound for expansion and parameter estimates}

\begin{notation}[Formulas for expansion in case (c) of \eqref{eqMainCases}]\label{def:AformularForLowerBoundOfExpansionForn1largerthan0}
Under our standing assumptions of this section (\cref{REMCaseCAssumptions}), for $x_0\in\Mrout$, $x_1\in\MRin$, $x_2\in\Mrin$, $n_1$ in \cref{REMCaseCAssumptions} and the free path and chord length functions $\tau_0$, $\tau_1$, $d_0$, $d_1$, $d_2$ from \cref{def:fpclf}. Denote by $\mathcal{B}^{\pm}_0$ the after/before wave front curvatures at $x_0$, as in \eqref{eq:MirrorEquationExpansionFormula} and let
\begin{equation}\label{eq:AformularForLowerBoundOfExpansionForn1largerthan0}
\begin{aligned}
\mathcal{II}_1(\mathcal{B}^+_0,n_1)&\dfn\big(1+\tau_0\mathcal{B}^+_0\big)+\big(\tau_1-\frac{2}{3}d_2\big)\Big[\big(2n_1+1-\frac{(2n_1+2)\tau_0}{d_1}\big)\mathcal{B}^+_0-\frac{2n_1+2}{d_1}\Big]\\
\mathcal{II}_2(\mathcal{B}^+_0,n_1)&\dfn 2n_1\big[1+(\tau_0-d_1)\mathcal{B}^+_0\big]
\end{aligned}
\end{equation}
\end{notation}

\begin{lemma}[Orbit configuration and expansion in \p-metric]\label{lemma:AformularForLowerBoundOfExpansionForn1largerthan0}
    With the standing assumptions of this section (\cref{REMCaseCAssumptions}), \COMMENT{\Large\color{red}Everything that precedes can be replaced by reference to \cref{REMCaseCAssumptions}. Check and implement.\wentao{much renovated statement}} we have the following.
    \begin{enumerate}
\item\label{itemAformularForLowerBoundOfExpansionForn1largerthan01} $x\in\Fc^{-1}(\Mrout\setminus\Mrin)\subset M_r$, $(\phi_0,\theta_0)=x_0=\Fc(x)\in\Mrout$, $x_1=\Fc(x_0)\in\MRin$. If we set
        $x_{1,R}\dfn\Fc(x_1)\in M_R$, $\cdots$, $ x_{n_1,R}\dfn\Fc^{n_1}(x_1)\in M_R$, then $\Fc(x_{n_1,R})=x_2\in\Mrin$, and we have $\frac{\|dx_0\|_\p}{\|dx\|_\p}\ge1$ and $\Bcal^+_0\in\big[-\frac{4}{3d_0},-\frac{1}{d_0}\big]$, where $d_0=r\sin\theta_0$ is from \cref{def:fpclf}.
        \item\label{itemAformularForLowerBoundOfExpansionForn1largerthan02} $dx_0\dfn D_x\Fc(dx)$ is a tangent vector at $x_0$, $dx_1\dfn D_{x_0}\Fc(dx_0)$ is a tangent vector at $x_1$, $dx_{n_1,R}\dfn
        D_{x_1}\Fc^{n_1}(dx_1)$ is a tangent vector at $x_{n_1,R}$, $dx_2\dfn D_{x_{n_1,R}}\Fc(dx_{n_1,R})$ is a tangent vector at $x_2$.
        \item\label{itemAformularForLowerBoundOfExpansionForn1largerthan03} $m(x_2)\ge3$ hence $x_3\dfn\Fc(x_2)\in M_r$. $dx_3\dfn D_{x_2}\Fc(dx_2)=D_x\Fc^{4+n_1}(dx)$ and \newline $\Fc^{m(x_2)-2}(x_3)=\Fc^{m(x_2)-1}(x_2)\overbracket{=}^{\mathllap{\text{\cref{def:SectionSetsDefs}}}}\hat{F}(x)\in\hat{M}$, i.e. $x_3=\Fc(x_2)$, $\cdots$, $\Fc^{m(x_2)-1}(x_2)$ are all in $M_r\setminus\Mrout$. We have $\frac{\|D_x\hat{F}(dx)\|_\p}{\|dx_3\|_\p}>1$.
        \item\label{itemAformularForLowerBoundOfExpansionForn1largerthan04} Computation of $\frac{\|dx_3\|_\p}{\|dx_0\|_\p}$. We use the notation from \eqref{eq:MirrorEquationExpansionFormula} with indices: denote by $\mathcal{B}^{\pm}_1$, $\mathcal{B}^{\pm}_{n_1,R}$, $\mathcal{B}^{\pm}_2$ and $\mathcal{B}^{\pm}_3$ the infinitesimal wave front curvatures after/before collisions at $x_1$, $x_{n_1,R}$, $x_2$ and $x_3$, respectively, $\tau_{1,\text{\upshape im}}\dfn-2n_1d_1$ and $\Bcal^+_{1,\text{\upshape im}}\dfn\Bcal^-_1-\frac{1}{d_1}$ (``\emph{im}'' means imaginary billiard trajectory in \cite[Figure 8.8, Equation (8.12)]{cb}).\COMMENT{Is there a figure in this article that we can refer to?} Then 
        \begin{equation}\label{eq:Expansionx3x0}\begin{aligned}E(n_1,\Bcal^+_0)\dfn\frac{\|dx_3\|_{\p}}{\|dx_0\|_{\p}}&=\bigm|1+\tau_0\Bcal^+_0\bigm|\bigm|1-2n_1d_1(\Bcal^-_1-\frac{1}{d_1})\bigm|\bigm|1+\tau_1\Bcal^{+}_{n_1,R}\bigm|\bigm|1+2d_2\Bcal^{+}_{2}\big|\\
        &=\bigm|1+\tau_0\Bcal^+_0\bigm|\bigm|1+\tau_{1,\text{\upshape im}}\Bcal^+_{1,\text{\upshape im}}\bigm|\bigm|1+\tau_1\Bcal^{+}_{n_1,R}\bigm|\bigm|1+2d_2\Bcal^{+}_{2}\big|\\&=3\bigm|\mathcal{II}_1(B^+_0,n_1)+\mathcal{II}_2(B^+_0,n_1)\big|,
        \end{aligned}\end{equation} where  $\mathcal{II}_1(\mathcal{B}^+_0,n_1)$, $\mathcal{II}_2(\mathcal{B}^+_0,n_1)$ and $\mathcal{B}^+_0$ are from \cref{def:AformularForLowerBoundOfExpansionForn1largerthan0}.
        \item\label{itemAformularForLowerBoundOfExpansionForn1largerthan05} 
         With $d_0$, $d_1$, $d_2$, $\tau_0$, $\tau_1$, $n_1$ as in \cref{def:fpclf,def:AformularForLowerBoundOfExpansionForn1largerthan0}, we have\begin{equation}\label{eq:B3minusEquation}
        \begin{multlined}[.85\textwidth]
            \mathcal{B}^{-}_3=\frac{-2}{d_2(1+2d_2\Bcal^+_2)}-\frac{1}{d_1(1+2d_2\Bcal_2^+)(1+\tau_1\Bcal^{+}_{n_1,R})}-\frac{1}{d_1(1+2d_2\Bcal^+_2)(1+\tau_1\Bcal^+_{n_1,R})(1+\tau_{1,\text{\upshape im}}\Bcal^+_{1,\text{\upshape im}})}\\
            +\frac{\Bcal^+_0}{(1+2d_2\Bcal^+_2)(1+\tau_1\Bcal^+_{n_1,R})(1+\tau_{1,\text{\upshape im}}\Bcal^+_{1,\text{\upshape im}})(1+\tau_0\Bcal_0^+)}.
        \end{multlined}\end{equation} 
    \end{enumerate}
\end{lemma}
\begin{proof}
\strut
   \hidesubsection{Proof of \eqref{itemAformularForLowerBoundOfExpansionForn1largerthan01}} For $x\in\hat{M}=\Fc^{-1}(\Mrout\setminus\Mrin)\sqcup(\Mrout\cap\Mrin)$ in \eqref{eqMhatOrbSeg} case (c), we first show that $x\notin \Mrout\cap\Mrin$. If we had  $x\in\Mrout\cap\Mrin$, then $x_0=x$ in \cref{def:MhatReturnOrbitSegment}. Since $x_0$ is nonsingular, we have $x_0=x\in\textcolor{red}{M^{\text{\upshape{out}}}_{r,0}}$. But by \cref{contraction_region} we have $x_0\in\Nout$ and \cref{Prop:Udisjoint} shows that $\Nout\cap\textcolor{red}{M^{\text{\upshape{out}}}_{r,0}}=\emptyset$. This means that it is impossible to have $x\in\Mrout\cap\Mrin$. Thus, $x\in\Fc^{-1}(\Mrout\setminus\Mrin)$, $x_0=\Fc(x)\in\Mrout$. 

    Hence, $dx_0\dfn D_{x}\Fc(dx)$ is a tangent vector at $x_0$. Since $x$, $x_0$ are nonsingular on $M_r$ and $dx=(d\phi,d\theta)$ is a tangent vector in the quadrant cone, by \cref{proposition:ShearOneStepExpansion} we have $\frac{\|dx_0\|_\p}{\|dx\|_\p}\ge1$. 
    
    Since $\frac{d\theta}{d\phi}\in[0,1]$, $(d\phi_0,d\theta_0)=dx_0=D_x\Fc(dx)$, $D_x\Fc=\begin{pmatrix}
        1 & 2\\
        0 & 1
    \end{pmatrix}$, that is $\begin{pmatrix} d\phi_0\\d\theta_0
    \end{pmatrix}=\begin{pmatrix}
        1 & 2\\
        0 & 1
    \end{pmatrix}\begin{pmatrix} d\phi\\d\theta
    \end{pmatrix}$, we have $\frac{d\theta_0}{d\phi_0}=\frac{d\theta}{d\phi+2d\theta}=\frac{1}{2+\frac{1}{\frac{d\theta}{d\phi}}}\in\big[0,\frac{1}{3}\big]$. Hence, by \eqref{eq:coordinateschange}, \eqref{eq:MirrorEquationExpansionFormula}, we have $\mathcal{B}^+_0\cos\varphi_0+\frac{1}{r}=\frac{d\varphi_0}{ds_0}=\frac{-1}{r}\frac{d\theta_0}{d\phi_0}\in[-\frac{1}{3r},0]$.
    
    Therefore, $\Bcal^+_0=-\frac{1}{r\cos\varphi_0}\big(1+\frac{d\theta_0}{d\phi_0}\big)\overbracket{\in}^{\mathclap{\substack{\frac{d\theta_0}{d\phi_0}\in[0,1/3]\text{, \eqref{eq:coordinateschange}: }\varphi_0=\pi/2-\theta_0 \text{ and \cref{def:fpclf}:}\\d_0=r\sin\theta_0=r\cos\varphi_0=}}}\big[-\frac{4}{3d_0},-\frac{1}{d_0}\big]$. This proves item \eqref{itemAformularForLowerBoundOfExpansionForn1largerthan01}.\COMMENT{\large\red$d_0=r\cos\varphi_0=r\sin\theta_0$ is inconsistent with \cref{def:fpclf}!\wentao{edited. Issue addressed.}}
    
    \hidesubsection{Proof of \eqref{itemAformularForLowerBoundOfExpansionForn1largerthan02}}
    We also have $x_1=\Fc(x_0)\in\MRin$, $dx_1\dfn D_{x_0}\Fc(dx_0)$. By \cref{def:MhatReturnOrbitSegment}, $ x_{n_1,R}\dfn\Fc^{n_1}(x_1)\in M_R$, $\Fc(x_{n_1,R})=x_2\in\Mrin$ is a sub-segment of \eqref{eqMhatOrbSeg} so that we can define $dx_{1,n_1}\dfn
    D_{x_1}\Fc^{n_1}(dx_1)$, $dx_2\dfn D_{x_{1,n_1}}\Fc(dx_{1,n_1})$, where $dx_{n_1,R}$ is a tangent vector at $x_{n_1,R}$, $dx_2$ is a tangent vector at $x_2$.
    This proves item \eqref{itemAformularForLowerBoundOfExpansionForn1largerthan02}.
    
    \hidesubsection{Proof of \eqref{itemAformularForLowerBoundOfExpansionForn1largerthan03}} By \cref{contraction_region}, we have $x_2\in\Nin$ and by \cref{Prop:Udisjoint}, we have $\emptyset=\Nin\cap(\textcolor{blue}{M^{\text{\upshape{in}}}_{r,0}}\cup\textcolor{blue}{M^{\text{\upshape{in}}}_{r,1}}\cup\textcolor{blue}{M^{\text{\upshape{in}}}_{r,2}})$. Hence, $x_2$ is nonsingular and $m(x_2)\ge3$. By the definition of $m(x_2)$ and $\hat{M}$ in \eqref{def:Mhat}, $(\phi_3,\theta_3)=x_3\dfn\Fc(x_2),\Fc(x_3)$, $\cdots$, $\Fc^{m(x_2)-2}(x_3)=\Fc^{m(x_2)-1}(x_2)\in\hat{M}$ are nonsingular points on $M_r\setminus\Mrout$ since their $\Fc$-images are not on $M_R$. Now $dx_3\dfn D_{x_2}\Fc(dx_2)=D_x\Fc^{4+n_1}(dx)$ and  $\Fc^{m(x_2)-2}(x_3)=\Fc^{m(x_2)-1}(x_2)\overset{\text{\cref{def:SectionSetsDefs}}}{\scalebox{7}[1]{$=$}}\hat{F}(x)\in\hat{M}$. 
    
    By \cref{caseB}, $D_x\Fc^{4+n_1}$ is a matrix with positive entries, since $(d\phi_3,d\theta_3)=dx_3=D_x\Fc^4(dx)$ and $dx=(d\phi,d\theta)$ is in the quadrant with $\frac{d\theta}{d\phi}\ge0$, $(d\phi_3,d\theta_3)=dx_3$ is then in the interior of the quadrant, that is, $\frac{d\theta_3}{d\phi_3}>0$. Since $D_x\hat{F}(dx)=D_{x_3}\Fc^{m(x_2)-2}(dx_3)$ and $x_3,\Fc(x_3),\cdots,\Fc^{m(x_2)-2}(x_3)$ are all in $M_r$, by \cref{proposition:ShearOneStepExpansion} we have $\frac{\|D_x\hat{F}(dx)\|_\p}{\|dx_3\|_\p}=\frac{\|D_{x_3}\Fc^{m(x_2)-2}(dx_3)\|_\p}{\|dx_3\|_\p}>1$. This proves item \eqref{itemAformularForLowerBoundOfExpansionForn1largerthan03}.
    
    \hidesubsection{Proof of \eqref{itemAformularForLowerBoundOfExpansionForn1largerthan04}}
    By applying \cite[exercises 8.28, 8.29]{cb} to $\Gamma_R$ and using \eqref{eq:MirrorEquationExpansionFormula}, \eqref{eqBBr}(\cite[equation (3.31)]{cb}), we get
    \begin{enumerate}[label={(\roman*)}]\label{Property:WavefrontsProofOf4}
        \item\label{item01ContinuedfractionResults}By \eqref{eq:MirrorEquationExpansionFormula}, $\frac{\|dx_2\|_\p}{\|dx_{n_1,R}\|_\p}=\big|1+\tau_1\Bcal^{+}_{n_1,R}\big|$, since $\tau_1$ is the distance between  $p(x_{n_1,R})=p(\Fcal^{-1}(x_2))$ and $p(x_2)$,
        \item\label{item02ContinuedfractionResults}By \eqref{eq:MirrorEquationExpansionFormula}, $\frac{\|dx_3\|_\p}{\|dx_2\|_\p}=\big|1+2d_2\Bcal^+_2\big|$, since $2d_2$ is the distance between $p(x_2)=p(\Fcal^{-1}(x_3))$ and $p(x_3)$,
        \item\label{item03ContinuedfractionResults}By \eqref{eqBBr}(\cite[equation (3.31)]{cb}), $\frac{1}{\Bcal^{+}_{n_1,R}}+\tau_1= 
        \frac{1}{\Bcal^{-}_2}$, thus $\Bcal^{-}_{2}=\frac{\Bcal^{+}_{n_1,R}}{1+\tau_1\Bcal^{+}_{n_1,R}}$,\COMMENT{\Large\red\bfseries Missing reference\wentao{fixed.}}
        \item\label{item04ContinuedfractionResults}By \eqref{eq:MirrorEquationExpansionFormula} with $\mathcal{K}=\frac{-1}{r}$, $\Bcal^{+}_2=
        \Bcal_2^{-}-\frac{2}{r\cos{\varphi_2}}=\Bcal^{-}_2-\frac{2}{d_2}$, since \eqref{eq:coordinateschange}: $\theta_2+\varphi_2=\pi/2$ and \cref{def:fpclf}: $d_2=r\sin\theta_2$,
        \item\label{item05ContinuedfractionResults}By \cite[exercise 8.28]{cb}, $\Bcal^{+}_{n_1,R}= 
        -\frac{1}{d_1}+\frac{1}{-2n_1d_1+\frac{1}{-\frac{1}{d_1}+\Bcal^{-}_1}}  = -\frac{1}{d_1}+\frac{\Bcal^{-}_1-\frac{1}{d_1}}{2n_1+1-2n_1d_1\Bcal^{-}_1}$,
        \item\label{item06ContinuedfractionResults}By \cite[exercise 8.29]{cb}, $\frac{\|dx_{n_1,R}\|_\p}{\|dx_1\|_\p}=
        \big|1-2n_1d_1(\Bcal^{-}_1+\frac{-1}{d_1})\big| =\big|1+2n_1-2n_1d_1\Bcal^{-}_1\big|$,
        \item\label{item07ContinuedfractionResults}By \eqref{eqBBr}(\cite[equation (3.31)]{cb}), $\frac{1}{\Bcal_{1}^{-}}=
        \frac{1}{\Bcal_{0}^{+}}+\tau_0$, thus $\Bcal^-_1=\frac{\Bcal^+_0}{1+\tau_0\Bcal^+_0}$,
        \item\label{item08ContinuedfractionResults}By \eqref{eq:MirrorEquationExpansionFormula}, $\frac{\|dx_1\|_\p}{\|dx_0\|_\p} 
        =\big|1+\tau_0\Bcal^{+}_0\big|$, since $\tau_0$ is the distance between $p(x_0)$ and $p(x_1)$.
        \end{enumerate}
        Therefore, \begin{align*}E(n_1,\Bcal^+_0)\dfn&\frac{\|dx_3\|_{\p}}{\|dx_0\|_{\p}}=\frac{\|dx_3\|_{\p}}{\|dx_2\|_{\p}}\frac{\|dx_2\|_{\p}}{\|dx_{n_1,R}\|_{\p}}\frac{\|dx_{n_1,R}\|_{\p}}{\|dx_1\|_{\p}}\frac{\|dx_1\|_{\p}}{\|dx_0\|_{\p}}\\
        \overbracket{=}^{\mathclap{\text{\ref{item01ContinuedfractionResults}\ref{item02ContinuedfractionResults}\ref{item04ContinuedfractionResults}\ref{item07ContinuedfractionResults}}}}&\bigm|1+2d_2\Bcal^+_2\bigm|\bigm|1+\tau_1\Bcal^{+}_{n_1,R}\bigm|\bigm|1-2n_1d_1(\Bcal^-_1-\frac{1}{d_1})\bigm|\bigm|1+\tau_0\Bcal^+_0\bigm|.
        \end{align*}        \begin{align*}\mathllap{\text{Next,\qquad}}
        \frac{\|dx_3\|_\p}{\|dx_{n_1,R}\|_\p}&=\frac{\|dx_2\|_\p}{\|dx_{n_1,R}\|_\p}\frac{\|dx_3\|_\p}{\|dx_2\|_\p}\overbracket=^{\text{(i)(ii)}}\big|1+\tau_1\Bcal^{+}_{n_1,R}\big|\big|1+2d_2\Bcal^+_2\big|\overbracket{=}^{\mathclap{\text{(iv)}}}\big|1+\tau_1\Bcal^{+}_{n_1,R}\big|\big|1+2d_2(\Bcal^-_2-\frac{2}{d_2})\big|\\
        &\overbracket=^{\text{(iii)}}\big|1+\tau_1\Bcal^{+}_{n_1,R}\big|\big|1+2d_2\big(\frac{\Bcal^{+}_{n_1,R}}{1+\tau_1\Bcal^{+}_{n_1,R}}-\frac{2}{d_2}\big)\big|=\big|1+\tau_1 \Bcal^{+}_{n_1,R}\big|\big|3-2d_2\frac{\Bcal_{n_1,R}^{+}}{1+\tau_1\Bcal^{+}_{n_1,R}}\big|\\
        &=\big|3+(3\tau_1-2d_2)\Bcal^{+}_{n_1,R}\big|\overset{\text{(v)}}{=}3\big|1+(\tau_1-\frac{2}{3}d_2)(\frac{-1}{d_1}+\frac{\Bcal^-_1-\frac{1}{d_1}}{2n_1+1-2n_1d_1\Bcal^{-}_1})\big|.
    \end{align*} 
    \begin{align*}\mathllap{\text{Therefore, }}
    \frac{\|dx_3\|_\p}{\|dx_1\|_\p}&=\frac{\|dx_3\|_\p}{\|dx_{n_1,R}\|_\p}\frac{\|dx_{n_1,R}\|_\p}{\|dx_1\|_\p}\overbracket=^{\text{(vi)}}3\big|1+(\tau_1-\frac{2d_2}{3})(\frac{-1}{d_1}+\frac{\Bcal^-_1-\frac{1}{d_1}}{2n_1+1-2n_1d_1\Bcal^{-}_1})\big|\big|1+2n_1-2n_1d_1\Bcal^{-}_1\big|\\
    &=3\big|1+2n_1-2n_1d_1\Bcal^{-}_1+(\tau_1-\frac{2}{3}d_2)(\Bcal^{-}_1-\frac{1}{d_1}-\frac{2n_1+1}{d_1}+2n_1\Bcal^{-}_1)\big|\\
    &=3\Big|1+2n_1-2n_1d_1\Bcal^{-}_1+\big(\tau_1-\frac{2}{3}d_2\big)\big[(2n_1+1)\Bcal^-_1-\frac{2n_1+2}{d_1}\big]\Big|\\
    &\overbracket{=}^{\mathclap{{\text{(vii)}}}}3\Big|1+2n_1-2n_1d_1\frac{\Bcal^+_0}{1+\tau_0\Bcal^+_0}+\big(\tau_1-\frac{2}{3}d_2\big)\big[(2n_1+1)\frac{\Bcal^+_0}{1+\tau_0\Bcal^+_0}-\frac{2n_1+2}{d_1}\big]\Big|\text{. }
    \end{align*}
\begin{equation}\label{eq:caseC0to3Expansion}
    \begin{aligned}\mathllap{\text{Hence \eqref{itemAformularForLowerBoundOfExpansionForn1largerthan04}:\ }}
    \frac{\|dx_3\|_\p}{\|dx_0\|_\p}&=\frac{\|dx_1\|_\p}{\|dx_0\|_\p}\frac{\|dx_3\|_\p}{\|dx_1\|_\p}\\
    &\overbracket{=}^{\mathclap{\text{(viii)}}}3\Big|1+\tau_0\Bcal^+_0\Big|\Big|1+2n_1-2n_1d_1\frac{\Bcal^+_0}{1+\tau_0\Bcal^+_0}+\big(\tau_1-\frac{2}{3}d_2\big)\big[(2n_1+1)\frac{\Bcal^+_0}{1+\tau_0\Bcal^+_0}-\frac{2n_1+2}{d_1}\big]\Big|\\
    & =3\Big|(1+\tau_0\Bcal^+_0)+2n_1(1+\tau_0\Bcal^+_0)-2n_1d_1\Bcal^+_0+\big(\tau_1-\frac{2}{3}d_2\big)\big[(2n_1+1)\Bcal^+_0-\frac{2n_1+2}{d_1}(1+\tau_0\Bcal^+_0)\big]\Big|\\
    & =3\Big|\underbracket{(1+\tau_0\Bcal^+_0)+(\tau_1-\frac{2}{3}d_2)\big[[(2n_1+1)-\frac{(2n_1+2)\tau_0}{d_1}]\Bcal^+_0-\frac{2n_1+2}{d_1}\big]}_{=\mathcal{II}_1(\Bcal^+_0,n_1)}+\underbracket{2n_1\big[1+(\tau_0-d_1)\Bcal^+_0\big]}_{=\mathcal{II}_2(\Bcal^+_0,n_1)}\Big|.
    \end{aligned}
\end{equation}

\hidesubsection{Proof of \eqref{itemAformularForLowerBoundOfExpansionForn1largerthan05}}\strut\COMMENT{\Large Note the overbracket fixes below; implement further, \color{red}such as, adding mathclap}\COMMENT{\color{red}Does ``\eqref{eqBBr}, \cite[equation (3.31)]{cb}'' below mean that both are needed? I suspect only \eqref{eqBBr} is. If correct, then delete ``\cite[equation (3.31)]{cb}''!}\COMMENT{Does ``since $2d_2$ is the distance between\dots'' explain the ``='' before or after? Both of those have added explanations on top anyway, so this seems confusing.}\COMMENT{\Large\red Missing reference}
\begin{align*}
\Bcal^{-}_3\overbracket{=}^{\Bigstrut\mathrlap{\text{\eqref{eqBBr} with } 2d_2 \text{ as the distance between } p(x_2)\text{ and }p(x_3)}}&\frac{\Bcal^+_2}{1+2d_2\Bcal^+_2},\\
\overbracket{=}^{\bigstrut\mathrlap{\text{\eqref{eq:MirrorEquationExpansionFormula} with } \mathcal{K}={-1}/{r},\ d_2=r\cos{\varphi_2}}}&\frac{-2/d_2}{1+2d_2\Bcal^+_2}+\frac{\Bcal^-_{2}}{1+2d_2\Bcal_2^+},\\
\overbracket{=}^{\bigstrut\mathclap{\text{\eqref{eqBBr}, \cite[equation (3.31)]{cb}}}}&\frac{-2/d_2}{1+2d_2\Bcal^+_2}+\frac{\overbracket{\Bcal^+_{n_1,R}}^{\mathclap{\text{\ref{item05ContinuedfractionResults} in proof of (4)}}}}{(1+2d_2\Bcal_2^+)(1+\tau_1\Bcal^+_{n_1,R})},\\
\overbracket{=}^{\mathclap{\text{\ref{item05ContinuedfractionResults} in proof of (4)}}}&\qquad\qquad\frac{-2/d_2}{1+2d_2\Bcal^+_2}+\frac{-\frac{1}{d_1}+\frac{\Bcal^{-}_1-\frac{1}{d_1}}{2n_1+1-2n_1d_1\Bcal^{-}_1}}{(1+2d_2\Bcal_2^+)(1+\tau_1\Bcal^+_{n_1,R})},\\
&=\frac{-2}{d_2(1+2d_2\Bcal^+_2)}-\frac{1}{d_1(1+2d_2\Bcal_2^+)(1+\tau_1\Bcal^+_{n_1,R})}-\frac{1}{d_1(1+2d_2\Bcal_2^+)(1+\tau_1\Bcal^+_{n_1,R})[1-2n_1d_1(\Bcal^-_1-\frac{1}{d_1})]}\\
&+\frac{\Bcal^-_1}{(1+2d_2\Bcal_2^+)(1+\tau_1\Bcal^+_{n_1,R})[1-2n_1d_1(\Bcal^-_1-\frac{1}{d_1})]}\\
\overbracket{=}^{\mathclap{\text{\eqref{eqBBr}, \cite[equation (3.31)]{cb}}}}&\qquad\qquad\frac{-2}{d_2(1+2d_2\Bcal^+_2)}-\frac{1}{d_1(1+2d_2\Bcal_2^+)(1+\tau_1\Bcal^+_{n_1,R})}-\frac{1}{d_1(1+2d_2\Bcal_2^+)(1+\tau_1\Bcal^+_{n_1,R})[1-2n_1d_1(\Bcal^-_1-\frac{1}{d_1})]}\\
&+\frac{\Bcal^+_0}{(1+2d_2\Bcal_2^+)(1+\tau_1\Bcal^+_{n_1,R})[1-2n_1d_1(\Bcal^-_1-\frac{1}{d_1})](1+\tau_0\Bcal^+_0)}.
\end{align*} This gives \eqref{eq:B3minusEquation}.
\end{proof}



\begin{figure}[h]
\begin{center}
    \begin{tikzpicture}[xscale=0.5,yscale=.5]
        \pgfmathsetmacro{\PHISTAR}{1.5};
        \pgfmathsetmacro{\Dlt}{0.6};
        \pgfmathsetmacro{\Dltplus}{\Dlt+0.6}
        \pgfmathsetmacro{\WIDTH}{8};
        \tkzDefPoint(-\PHISTAR,0){LL};
        \tkzDefPoint(-\PHISTAR,\WIDTH){LU};
        \tkzDefPoint(\PHISTAR,0){RL};
        \tkzDefPoint(\PHISTAR,\WIDTH){RU};
        \tkzDefPoint(-\PHISTAR,\PHISTAR){X};
        \tkzDefPoint(\PHISTAR,\WIDTH-\PHISTAR){Y};
        \tkzDefPoint(-\PHISTAR,\WIDTH-\PHISTAR){IX};
        \tkzDefPoint(\PHISTAR,\PHISTAR){IY};
        \tkzDefPoint(-\PHISTAR,\PHISTAR+\Dlt){LLD};
        \tkzDefPoint(\PHISTAR,\PHISTAR+\Dlt){RLD};
        \tkzDefPoint(-\PHISTAR,\PHISTAR+\Dlt){LLD};
        \tkzDefPoint(\PHISTAR,\PHISTAR+\Dlt){RLD};
        \tkzDefPoint(-\PHISTAR,\WIDTH-\PHISTAR-\Dlt){LUD};
        \tkzDefPoint(\PHISTAR,\WIDTH-\PHISTAR-\Dlt){RUD};
        \tkzDefPoint(-\PHISTAR,\PHISTAR+\Dltplus){LLDp};
        \tkzDefPoint(\PHISTAR,\PHISTAR+\Dltplus){RLDp};
        \tkzDefPoint(-\PHISTAR,\WIDTH-\PHISTAR-\Dltplus){LUDp};
        \tkzDefPoint(\PHISTAR,\WIDTH-\PHISTAR-\Dltplus){RUDp};
        \tkzDefPoint(-\PHISTAR,\PHISTAR+\Dltplus+0.8){LLDdp};
        \tkzDefPoint(\PHISTAR,\PHISTAR+\Dltplus+0.8){RLDdp};
        \tkzDefPoint(-\PHISTAR,\WIDTH-\PHISTAR-\Dltplus-0.8){LUDdp};
        \tkzDefPoint(\PHISTAR,\WIDTH-\PHISTAR-\Dltplus-0.8){RUDdp};
        \tkzDefPoint (0.3*\PHISTAR,\PHISTAR+\Dltplus){montonestart};
        \tkzDefPoint (0.3*\PHISTAR,0.333*\WIDTH){montoneend};
        \draw [thin] (LL) --(LU);
        \draw [thin] (LU) --(RU);
        \draw [thin] (RL) --(RU);
        \draw [thin] (LL) --(RL);
        \draw [blue,ultra thin,dashed] (LL) --(IY);
        \draw [red,ultra thin,dashed]  (X)--(RL);
        \draw [blue,ultra thin,dashed] (RU) --(IX);
        \draw [red,ultra thin,dashed]  (LU)--(Y);
        \tkzLabelPoint[right](RU){$\pi$};
        \tkzLabelPoint[right](RL){$0$};
        \tkzLabelPoint[right](IY){$\theta_1=\Phistar,d_1=r\sin{\phistar}$};
        \tkzLabelPoint[right](Y){$\theta_1=\pi-\Phistar,d_1=r\sin{\phistar}$};
        \coordinate (RM) at ($(RL)!0.5!(RU)$);
        \tkzLabelPoint[right](RM){${\pi}/{2}$};
        \coordinate (UpperM) at ($(LU)!0.5!(RU)$);
        \coordinate (LowerM) at ($(LL)!0.5!(RL)$);
        \coordinate (LM) at ($(LU)!0.5!(LL)$);
        \tkzLabelPoint[below](LL){$-\Phistar$};
        \tkzLabelPoint[below](RL){$+\Phistar$};
        \tkzLabelPoint[below](LowerM){$0$};
        \tkzDefPoint(0,0.5*\PHISTAR){LMend};
        \tkzDefPoint(0,\WIDTH-0.5*\PHISTAR){HMend};
        \fill [green, opacity=10/30] (LL) -- (LMend) -- (X) -- cycle;
        \fill [green, opacity=10/30] (RU) -- (HMend) -- (Y) -- cycle;
        \tkzDrawPoints(X,Y,IX,IY,LMend,HMend);
        \node[] at (0,0.5*\WIDTH){\Large $n_1=0$};
        \node[] at (-0.7*\PHISTAR,0.5*\PHISTAR){\small $n_1\ge1$};
        \node[] at (0.7*\PHISTAR,\WIDTH-0.5*\PHISTAR){\small $n_1\ge1$};
    \end{tikzpicture}
    \caption{Two components of \MRin\ where $x_1\in\MRin$ with $n_1+1\ge2$}\label{fig:x1_multi_collisions_caseC}
\end{center}
\end{figure}
\begin{proposition}[Case (c) 
parameter estimates
]\label{proposition:casecparametersestimate} With the standing assumptions of this section (\cref{REMCaseCAssumptions}),\COMMENT{\Large\red Edit if we state a standing assumption\wentao{OK, done}} for the length functions $\tau_0$, $\tau_1$, $d_0$, $d_1$, $d_2$ from \cref{def:AformularForLowerBoundOfExpansionForn1largerthan0}, under condition \eqref{eqJZHypCond}, we have the following parameter estimates.

If $n_1\ge1$, then
\begin{enumerate}[label={(\roman*)}]
        \item\label{itemcasecparametersestimate1} $\big|d_0-r\sin\phistar\big|\overbracket{<}^{\mathclap{\text{\cite[(3.20)]{hoalb}}}}\frac{17r^2\sin\phistar}{4R}$ and $\big|d_2-r\sin\phistar\big|\overbracket{<}^{\mathclap{\text{\cite[(3.20)]{hoalb}}}}\frac{17r^2\sin\phistar}{4R}$, (also see \cite[(5.11)]{hoalb}),
        \item\label{itemcasecparametersestimate2} $2r\sin{\phistar}\overbracket{<}^{\mathclap{\text{\cite[(3.16), (4.7)]{hoalb}}}}\tau_0+\tau_1+2n_1d_1\overbracket{<}^{\mathclap{\text{\cite[(3.16), (4.7)]{hoalb}}}}2r\sin{\phistar}+\frac{8.1r^2\sin{\phistar}}{R}$, (also see \cite[(5.12)]{hoalb}),
        \item\label{itemcasecparametersestimate3} $d_1\overbracket{>}^{\mathclap{\substack{\text{\upshape{(ii) and \cref{remark:deffpclf}: }}\\\tau_1<2d_1,\tau_0<2d_1}}}\frac{1}{n_1+2}r\sin{\phistar}$,
        \item\label{itemcasecparametersestimate4} $d_1<\frac{r\sin{\phistar}}{n_1}+\frac{8.1r}{2n_1R}r\sin{\phistar}$ by \ref{itemcasecparametersestimate2}: $2n_1d_1<\tau_0+\tau_1+2n_1d_1<2r\sin{\phistar}+\frac{8.1r^2\sin{\phistar}}{R}$,
        \item\label{itemcasecparametersestimate5} $\frac{r\sin\phistar}{d_0}\in(0.9975,1.0026)\text{ \upshape{and} } \frac{r\sin\phistar}{d_2}\in(0.9975,1.0026)$ by \ref{itemcasecparametersestimate1} and \eqref{eqJZHypCond}: $R>1700r$.
        \end{enumerate}
        
        If $n_1\ge2$, then
        \begin{enumerate}[label={(\roman*)}]
        \setcounter{enumi}{5}
        \item\label{itemcasecparametersestimate6} $\tau_0<\frac{2}{3}d_0+\frac{2}{3}d_1$ by \cite[Proposition 5.2]{hoalb},
        \item\label{itemcasecparametersestimate7} $\frac{2}{3}d_2+\frac{2}{3}d_1-\tau_1>\frac{1}{6}r\sin{\phistar}-\frac{16.6r}{3R}r\sin{\phistar}$ and $\frac{2}{3}d_0+\frac{2}{3}d_1-\tau_0>\frac{1}{6}r\sin{\phistar}-\frac{16.6r}{3R}r\sin{\phistar}$  by \cite[(5.15) and subsequent statement]{hoalb}.
    \end{enumerate} 
\end{proposition}
\begin{remark}

{\cref{proposition:casecparametersestimate}\ref{itemcasecparametersestimate1}} with $n_1\ge2$ is the claim of \cite[Inequalities (5.11)]{hoalb}. Its proof is based on the Mean-Value Theorem and \cite[Inequalities (3.20)]{hoalb}. This proof argument also works for $n_1=1$. 

{\cref{proposition:casecparametersestimate}\ref{itemcasecparametersestimate2}} with $n_1\ge2$ is the claim of \cite[(5.12)]{hoalb}. Its proof relies on \cite[Inequalities (3.16)(4.7)]{hoalb}. This proof argument works the same for $n_1=1$. 

In \cite{hoalb}, we also note that the statement following \cite[equation (5.19)]{hoalb} regards \cite[Inequalities (5.11)(5.12)]{hoalb} as established for $n_1=1$ to argue that $\tau_0\ge\frac{2}{3}d_0+\frac{1}{2}d_1$ and $\tau_1\ge\frac{1}{2}d_1+\frac{2}{3}d_2$ cannot simultaneously hold for $n_1=1$.

{\cref{proposition:casecparametersestimate}\ref{itemcasecparametersestimate5}} follows from \COMMENT{\large\red follows from?\wentao{OK.}} {\cref{proposition:casecparametersestimate}\ref{itemcasecparametersestimate1}}. Dividing by $r\sin\phistar$ on both sides of the two inequalities in \ref{itemcasecparametersestimate1} gives $\big|\frac{d_i}{r\sin\phistar}-1\bigm|<\frac{17r}{4R}$ for $i=0,2$. Then \eqref{eqJZHypCond}: $R>1700$ gives the conclusion of \ref{itemcasecparametersestimate5}.

{\cref{proposition:casecparametersestimate}\ref{itemcasecparametersestimate7}} is the claim of \cite[Inequalities (5.15)]{hoalb} and its following statement in \cite{hoalb}.
    
\end{remark}


\subsection[The expansion in Case (c) of \eqref{eqMainCases} with \texorpdfstring{$n_1\ge4$}{n1>3}]{The expansion in Case (c) of \texorpdfstring{\eqref{eqMainCases}}{orbit transition} with \texorpdfstring{$n_1\ge4$}{n1>3}}\label{subsec:expansionLarge}\hfill

\begin{proposition}\label{proposition:casecexpansionn1ge4}
    With the standing assumptions of this section (\cref{REMCaseCAssumptions}), if $n_1\ge4$, then $\mathcal{II}_1(\mathcal{B}^+_0,n_1)$, $\mathcal{II}_2(\mathcal{B}^+_0,n_1)$ in \cref{def:AformularForLowerBoundOfExpansionForn1largerthan0,lemma:AformularForLowerBoundOfExpansionForn1largerthan0} satisfy $\mathcal{II}_1(\Bcal^+_0,n_1)+\mathcal{II}_2(\Bcal^+_0,n_1)>\frac{0.32}{1.01}n_1^2+1.885n_1-2.876>9.73$. 
    
    Hence, $E(n_1,\Bcal^+_0)=\frac{\|dx_3\|_\p}{\|dx_0\|_\p}=3\big|\mathcal{II}_1(\mathcal{B}^+_0,n_1)+\mathcal{II}_2(\mathcal{B}^+_0,n_1)\big|>\frac{0.96}{1.01}n_1^2+5.655n_1-8.628>29.1$ in \cref{lemma:AformularForLowerBoundOfExpansionForn1largerthan0}(4).
\end{proposition}
\begin{proof}
        We first estimate $\mathcal{II}_2(\mathcal{B}^+_0,n_1)=2n_1\big[1+(\tau_0-d_1)\Bcal^+_0\big]$ with $n_1\ge1$.
        \begin{align*}
            \mathcal{II}_2(\mathcal{B}^+_0,n_1)&=2n_1\big[1+(\tau_0-d_1)\Bcal^+_0\big]\overbracket{>}^{\mathclap{\substack{\text{\cref{lemma:AformularForLowerBoundOfExpansionForn1largerthan0} (1): } \Bcal^+_0<0 \\\text{ \cref{remark:deffpclf}: } \tau_0-d_1<d_1}}}2n_1+2n_1d_1\Bcal_0^+\overbracket{>}^{\mathclap{\substack{\Bcal^+_0<0\text{ and}\\\text{\cref{proposition:casecparametersestimate} \ref{itemcasecparametersestimate4}}}}}2n_1+2n_1\big(\frac{r\sin\phistar}{n_1}+\frac{8,1r}{2n_1R}r\sin\phistar\big)\Bcal_0^+\\
            &=2n_1+\big(2+\frac{8.1r}{R}\big)r\sin\phistar\Bcal_0^+\overbracket{\ge}^{\mathclap{\Bcal_0^+\ge\frac{-4}{3d_0}}}2n_1-\big(2+\frac{8.1r}{R}\big)\frac{4r\sin\phistar}{3d_0}\\
            &\overbracket{>}^{\quad\mathllap{\text{\cref{proposition:casecparametersestimate}  }\ref{itemcasecparametersestimate5}}}2n_1-\frac{8\times1.0026}{3}-\frac{32.4r\times1.0026}{3R}\overbracket{>}^{\mathrlap{\eqref{eqJZHypCond}:R>1700r}}2n_1-2.693,
        \end{align*}
        that is,
        \begin{equation}\label{eq:II2B4estimate1}
        \begin{aligned}
            \text{if }n_1\ge1\text{, then }&\mathcal{II}_2(\mathcal{B}^+_0,n_1)>2n_1-2.693.
        \end{aligned}
        \end{equation}
        Next, we estimate $\mathcal{II}_1(\mathcal{B}^+_0,n_1)=\big(1+\tau_0\mathcal{B}^+_0\big)+\big(\tau_1-\frac{2}{3}d_2\big)\Big[\big(2n_1+1-\frac{(2n_1+2)\tau_0}{d_1}\big)\mathcal{B}^+_0-\frac{2n_1+2}{d_1}\Big]$, beginning with\COMMENT{\Large\bfseries\textcolor{red}{Some reasoning labels seem wrong. Check!!}\wentao{There was some mislabeling a while ago, but it was fixed. Now it seems correct labeling.}}
        \begin{align*}
        \tau_1-\frac{2}{3}d_2&\overbracket{<}^{\mathclap{\substack{\text{\cref{remark:deffpclf}: }\\\tau_1<2d_1}}}2d_1-\frac{2}{3}d_2\overbracket{<}^{\mathclap{\strut\text{\cref{proposition:casecparametersestimate} \ref{itemcasecparametersestimate1}}}}2d_1-\frac{2}{3}\big(r\sin\phistar-\frac{17r}{4R}r\sin\phistar\big)\\
        &\overbracket{<}^{\quad\mathllap{\text{\cref{proposition:casecparametersestimate}\ref{itemcasecparametersestimate4}}}}\big(\frac{2r\sin{\phistar}}{n_1}+\frac{8.1r\sin{\phistar}}{n_1}\frac{r}{R}\big)-\frac{2}{3}\big(r\sin\phistar-\frac{17r}{4R}r\sin\phistar\big)\\
        &\overbracket{\le}^{\quad\mathllap{n_1\ge4}}\frac{1}{2}r\sin{\phistar}+\frac{8.1r}{4R}r\sin{\phistar}-\frac{2}{3}r\sin{\phistar}+\frac{17r}{6R}r\sin{\phistar}\\
        &=\big(-\frac{1}{6}+\frac{58.3}{12}\cdot\frac{r}{R}\big)r\sin{\phistar}\overbracket{<}^{\mathclap{\eqref{eqJZHypCond}:R>1700r}}-0.16r\sin{\phistar}.
        \end{align*}
        Thus, 
\begin{equation}\label{eq:II1B4estimate1}
        \begin{aligned}
            \text{If }n_1\ge1\text{, then }\tau_1-\frac{2}{3}d_2<&\big(\frac{2r\sin\phistar}{n_1}+\frac{8.1r\sin\phistar}{n_1}\frac{r}{R}\big)-\frac{2}{3}(r\sin\phistar-\frac{17R}{4r}r\sin\phistar).\\
            \text{If }n_1\ge4\text{, then }\tau_1-\frac{2}{3}d_2<&-0.16r\sin\phistar.
        \end{aligned}
        \end{equation}
    We next use that $0<\frac{\tau_0}{d_1}<2$ (\cref{remark:deffpclf}), hence $(2n_1+1)-\frac{(2n_1+2)\tau_0}{d_1}> 2n_1+1-2(2n_1+2)=-2n_1-3$, and  $\Bcal^+_0<0$ (\cref{lemma:AformularForLowerBoundOfExpansionForn1largerthan0} (1)) to get
\begin{align*}
    \overbracket{\big[(2n_1+1)-\frac{(2n_1+2)\tau_0}{d_1}\big]}^{>-2n_1-3}\overbracket{\Bcal^+_0}^{<0}-\frac{2n_1+2}{d_1}<&-(2n_1+3)\Bcal^+_0-\frac{2n_1+2}{d_1}\overbracket\le^{\mathclap{\Bcal_0^+\ge{-4}/{3d_0}}}\frac{4}{3d_0}(2n_1+3)-\frac{2n_1+2}{d_1}\\
    \overbracket{<}^{\quad\mathllap{\text{\eqref{eqJZHypCond} and \cref{proposition:casecparametersestimate} \ref{itemcasecparametersestimate1} }: 0<r\sin\phistar-\frac{17r}{4R}r\sin\phistar<d_0}}&\frac{4(2n_1+3)}{3(r\sin{\phistar}-\frac{17r}{4R}r\sin{\phistar})}-\frac{2n_1+2}{d_1}\\
    \overbracket{<}^{\quad\mathllap{\text{\cref{proposition:casecparametersestimate}\ref{itemcasecparametersestimate1}: }0<d_1<\frac{r}{n_1}\sin{\phistar}+\frac{8.1r}{2n_1R}r\sin{\phistar}}}&\frac{4(2n_1+3)}{3(r\sin{\phistar}-\frac{17r}{4R}r\sin{\phistar})}-\frac{(2n_1+2)}{\frac{r}{n_1}\sin{\phistar}+\frac{8.1r}{2n_1R}r\sin{\phistar}}\\
    \overbracket{<}^{\quad\mathllap{\eqref{eqJZHypCond}\text{: R>1700r}}}&\big[\frac{4}{3\times0.99}(2n_1+3)-\frac{n_1(2n_1+2)}{1.01}\big]\frac{1}{r\sin{\phistar}}\\
    &=\big[-\frac{2}{1.01}n_1^2+(\frac{8}{3\times0.99}-\frac{2}{1.01})n_1+\frac{4}{0.99}\big]\frac{1}{r\sin\phistar}\\
    <&\underbracket{\big[-\frac{2}{1.01}n_1^2+0.7135n_1+\frac{4}{0.99}\big]}_{\text{decreases with } n_1\ge1}\frac{1}{r\sin\phistar}\overbracket<^{\mathclap{n_1\ge4}}\frac{-24}{r\sin\phistar}<0.
\end{align*}
Therefore:
\begin{equation}\label{eq:II1B4estimate2}
\begin{aligned}
\text{If }n_1\ge1\text{, then }&\big[(2n_1+1)-\frac{(2n_1+2)\tau_0}{d_1}\big]\Bcal^+_0-\frac{2n_1+2}{d_1}<\big[-\frac{2}{1.01}n_1^2+0.7135n_1+\frac{4}{0.99}\big]\frac{1}{r\sin\phistar}.\\
\text{If }n_1\ge4\text{, then }&\big[(2n_1+1)-\frac{(2n_1+2)\tau_0}{d_1}\big]\Bcal^+_0-\frac{2n_1+2}{d_1}<\frac{-24}{r\sin\phistar}<0.
\end{aligned}
\end{equation}
Furthermore, since $(n_1+1)\tau_0=\tau_0+n_1\tau_0\overbracket{<}^{\mathclap{\substack{\text{\cref{remark:deffpclf}:}\\\tau_0<2d_1}}}\tau_0+2n_1d_1+\tau_1\overbracket{<}^{\mathclap{\substack{\text{\cref{proposition:casecparametersestimate}\ref{itemcasecparametersestimate2}}}}}2r\sin\phistar+\frac{8.1r}{R}r\sin\phistar$, we have \[0<\tau_0<\frac{2r\sin\phistar}{n_1+1}+\frac{8.1r}{(n_1+1)R}r\sin\phistar=\big[\frac{2}{n_1+1}+\frac{8.1r}{(n_1+1)R}\big]r\sin\phistar,\] as well as
    \begin{align*}
    1+\underbracket{\tau_0}_{\mathclap{<[\frac{2}{n_1+1}+\frac{8.1r}{(n_1+1)R}]r\sin\phistar}}\overbracket{\Bcal_0^+}^{\mathclap{\text{By \cref{lemma:AformularForLowerBoundOfExpansionForn1largerthan0} \eqref{itemAformularForLowerBoundOfExpansionForn1largerthan01}, }\Bcal^+_0<0}}&>1+\big[\frac{2}{n_1+1}+\frac{8.1r}{(n_1+1)R}\big]r\sin\phistar\Bcal^+_0\\
    \overbracket{>}^{\quad\mathllap{\text{\cref{lemma:AformularForLowerBoundOfExpansionForn1largerthan0} (1), }\Bcal^+_0\in[\frac{-4}{3d_0},\frac{-1}{d_0}]}}&1+\big[\frac{2}{n_1+1}+\frac{8.1r}{(n_1+1)R}\big]\frac{-4r\sin\phistar}{3d_0}=1-\frac{4}{3}\big[\frac{2}{n_1+1}+\frac{8.1r}{(n_1+1)R}\big]\frac{r\sin\phistar}{d_0}\\
    \overbracket{>}^{\quad\mathllap{\text{\cref{proposition:casecparametersestimate} \ref{itemcasecparametersestimate5}}}}&1-\frac{4}{3}\big[\frac{2}{1+n_1}+\frac{8.1r}{(1+n_1)R}\big]\times1.0026\overbracket>^{\mathclap{\substack{\text{\eqref{eqJZHypCond}: R>1700r}\\\text{ and }n_1\ge4}}}1-\frac{4}{3}(\frac{2}{5}+\frac{8.1}{1700\cdot 5})\times1.0026>0.4640.
    \end{align*}
Thus:
    \begin{equation}\label{eq:II1B4estimate3}
    \begin{aligned}
        \text{If }n_1\ge1\text{, then }&1+\tau_0\Bcal^+_0>1+\big[\frac{2}{n_1+1}+\frac{8.1r}{(n_1+1)R}\big]\times\big(\frac{-4r\sin\phistar}{3d_0}\big).\\
        \text{If }n_1\ge4\text{, then }&1+\tau_0\Bcal^+_0>0.464.
    \end{aligned}
    \end{equation}
    Therefore, 
    if $n_1\ge4$, then \begin{align*}\mathcal{II}_1(\mathcal{B}^+_0,n_1)+\mathcal{II}_2(\mathcal{B}^+_0,n_1)&=\big(\overbracket{1+\tau_0\mathcal{B}^+_0}^{\mathclap{\text{\eqref{eq:II1B4estimate3}: }>0.464}}\big)+\big(\underbracket{\tau_1-\frac{2}{3}d_2}_{\mathclap{\text{\eqref{eq:II1B4estimate1}:}<-0.16r\sin\phistar<0}}\big)\Big[\overbracket{\big(2n_1+1-\frac{(2n_1+2)\tau_0}{d_1}\big)\mathcal{B}^+_0-\frac{2n_1+2}{d_1}}^{\mathclap{\text{\eqref{eq:II1B4estimate2}: }<\big(-\frac{2}{1.01}n_1^2+0.7135n_1+\frac{4}{0.99}\big)\frac{1}{r\sin\phistar}<\frac{-24}{r\sin\phistar}<0}}\Big]
    +\overbracket{\mathcal{II}_2(\mathcal{B}^+_0,n_1)}^{\mathclap{\text{\eqref{eq:II2B4estimate1}: }>2n_1-2.693}}\\
    &>0.464+(-0.16)\times(-\frac{2}{1.01}n_1^2+0.7135n_1+\frac{4}{0.99})+2n_1-2.693\\
    &>\frac{0.32}{1.01}n_1^2+1.885n_1-2.876.\qedhere
    \end{align*}
\end{proof}

\subsection[The expansion in Case (c) of \eqref{eqMainCases} with \texorpdfstring{$n_1=3$}{n1=3}]{The expansion in Case (c) of \texorpdfstring{\eqref{eqMainCases}}{orbit transition} with \texorpdfstring{$n_1=3$}{n1=3}}\hfill

By symmetry, assume that $x_1=(\Phi_1,\theta_1)$ with $n_1=3$ is in the lower component in \cref{fig:x1_multi_collisions_caseC} with $\theta_1\in (0,\frac{\pi}{2})$, i.e., in the shaded region in \cref{fig:x1_n1cases_caseC,fig:x1_n1equal3case_subcases} with $\theta_1\in\big(\frac{1}{5}\Phistar,\frac{1}{3}\Phistar\big)$.
We use a horizontal line $\theta_1=0.235\Phistar$ to separate this region into two parts: upper part $\theta_1\ge0.235\Phistar$, lower part $\theta_1<0.235\Phistar$, as shown in \cref{fig:x1_n1equal3case_subcases}.

\begin{figure}[h]
\begin{minipage}{.54\textwidth}
\begin{center}
    \begin{tikzpicture}[xscale=0.6, yscale=0.6]
        \pgfmathsetmacro{\PHISTAR}{6};
        \pgfmathsetmacro{\WIDTH}{8};
        \tkzDefPoint(-\PHISTAR,0){LL};
        \tkzDefPoint(-\PHISTAR,\WIDTH){LU};
        \tkzDefPoint(\PHISTAR,0){RL};
        \tkzDefPoint(\PHISTAR,\WIDTH){RU};
        \tkzDefPoint(-\PHISTAR,\PHISTAR){X};
        \tkzDefPoint(\PHISTAR,\WIDTH-\PHISTAR){Y};
        \tkzDefPoint(-\PHISTAR,\WIDTH-\PHISTAR){IX};
        \tkzDefPoint(\PHISTAR,\PHISTAR){IY};
        \tkzDefPoint(-\PHISTAR,0.5*\PHISTAR){S1};
        \tkzDefPoint(-\PHISTAR,\PHISTAR/3){S2};
        \tkzDefPoint(-\PHISTAR,\PHISTAR/4){S3};
        \tkzDefPoint(-\PHISTAR,\PHISTAR/5){L4};
        \draw [thin] (LL) --(LU);
        \draw [thin] (LU) --(RU);
        \draw [thin] (RL) --(RU);
        \draw [thin] (LL) --(RL);
        \draw [name path=lineLLIY, blue,ultra thin,dashed] (LL) --(IY);
        \draw [red,ultra thin,dashed]  (X)--(RL);
        \draw [red,ultra thin,dashed]  (S1)--(RL);
        \draw [name path=lineS2RL, red,ultra thin,dashed]  (S2)--(RL);
        \draw [name path=lineS3RL, red,ultra thin,dashed]  (S3)--(RL);
        \path[name intersections= 
        {of=lineLLIY and lineS3RL,by={LS}}];
        \path[name intersections=
        {of=lineLLIY and lineS2RL,by={HS}}];
        \tkzDrawPoints(S2,S3,HS,LS,L4);
        \draw [ultra thin,dashed] (L4) --(LS);
        \tkzLabelPoint[above left](S2){\small$\theta_1=\frac{1}{3}\Phistar$};
        \tkzLabelPoint[left](S3){\small$\theta_1=0.25\Phistar$};
        \tkzLabelPoint[below left](L4){\small $\theta_1=0.2\Phistar$};
        \node[] at (0,1.1*\PHISTAR){\Large $n_1=0$};
        \node[] at (-0.7*\PHISTAR,0.6*\PHISTAR){\Large $n_1=1$};
        \node[] at (-0.83*\PHISTAR,0.38*\PHISTAR){\Large $n_1=2$};
        \node[] at (-0.85*\PHISTAR,0.27*\PHISTAR){\small $n_1=3$};
        \node[] at (-0.90*\PHISTAR,0.13*\PHISTAR){ $n_1\ge4$};
        \fill [green, opacity=5/30](LS) -- (HS) -- (S2) -- (S3) -- cycle;
        \tkzLabelPoint[below](LL){$-\Phistar$};
        \tkzLabelPoint[below](RL){$+\Phistar$};
    \end{tikzpicture}
    \captionsetup{{width=.9\linewidth}}
    \caption{Quadrilateral of $n_1=3$, $\theta_1\in[0.2\Phistar,\frac{1}{3}\Phistar]$}\label{fig:x1_n1cases_caseC}
\end{center}
\end{minipage}\hfill
\begin{minipage}{0.45\textwidth}
\begin{center}
    \begin{tikzpicture}[xscale=2.2, yscale=2.2]
        \pgfmathsetmacro{\PHISTAR}{6};
        \pgfmathsetmacro{\WIDTH}{8};
        \tkzDefPoint(-\PHISTAR,0){LL};
        \tkzDefPoint(-\PHISTAR,\WIDTH){LU};
        \tkzDefPoint(\PHISTAR,0){RL};
        \tkzDefPoint(\PHISTAR,\WIDTH){RU};
        \tkzDefPoint(-\PHISTAR,\PHISTAR){X};
        \tkzDefPoint(\PHISTAR,\WIDTH-\PHISTAR){Y};
        \tkzDefPoint(-\PHISTAR,\WIDTH-\PHISTAR){IX};
        \tkzDefPoint(\PHISTAR,\PHISTAR){IY};
        \tkzDefPoint(-\PHISTAR,0.5*\PHISTAR){S1};
        \tkzDefPoint(-\PHISTAR,\PHISTAR/3){S2};
        \tkzDefPoint(-\PHISTAR,\PHISTAR/4){S3};
        \tkzDefPoint(-\PHISTAR,\PHISTAR/5){L4};
        \clip(-\PHISTAR-0.05,0.15*\PHISTAR) rectangle (-\PHISTAR+\PHISTAR*0.52, \PHISTAR*0.35);
        \draw [thin] (LL) --(LU);
        \draw [thin] (LU) --(RU);
        \draw [thin] (RL) --(RU);
        \draw [thin] (LL) --(RL);
        \draw [name path=lineLLIY, blue,ultra thin,dashed] (LL) --(IY);
        \draw [red,ultra thin,dashed]  (X)--(RL);
        \draw [red,ultra thin,dashed]  (S1)--(RL);
        \draw [name path=lineS2RL, red,ultra thin,dashed]  (S2)--(RL);
        \draw [name path=lineS3RL, red,ultra thin,dashed]  (S3)--(RL);
        \path[name intersections= 
        {of=lineLLIY and lineS3RL,by={LS}}];
        \path[name intersections=
        {of=lineLLIY and lineS2RL,by={HS}}];
        \tkzDefPoint(-\PHISTAR,0.235*\PHISTAR){L5};
        \tkzDefPoint(-\PHISTAR,0.222*\PHISTAR){L6};
        \path[name path=HorSeg] (L5)-- +(5,0);
        \path[name path=HorSeg1] (L6)-- +(5,0);
        \path[name intersections={of=lineS3RL and HorSeg,by={Seg1}}];
        \path[name intersections={of=lineLLIY and HorSeg,by={Seg2}}];
        \draw[very thin] (L5)--(Seg2);
        \tkzLabelPoint[above](Seg1){\small$\theta_1\ge0.235\Phistar$};
        \tkzLabelPoint[below](Seg1){\small$\theta_1<0.235\Phistar$};
        \path[name intersections={of=lineS3RL and HorSeg1,by={Seg3}}];
        \path[name intersections={of=lineLLIY and HorSeg1,by={Seg4}}];
        \draw[very thin,blue] (Seg3)--(Seg4);
        \tkzLabelPoint[below](Seg3){\small$\mathcal{W}$};
        \tkzLabelPoint[below](Seg4){\small$\mathcal{E}$};
         \tkzLabelSegment[below](Seg3,Seg4){\small{Line segment $\mathfrak{l}$}};  \tkzDrawPoints(S2,S3,HS,LS,L5,Seg1,Seg2,Seg3,Seg4);
        \fill [green, opacity=5/30](LS) -- (HS) -- (S2) -- (S3) -- cycle;
        \tkzLabelPoint[below](LL){$-\Phistar$};
        \tkzLabelPoint[below](RL){$+\Phistar$};
    \end{tikzpicture}
    \caption{Quadrilateral $n_1=3$ divided into two parts: upper part with $\theta_1\ge0.235\Phistar$ and lower part with $\theta_1<0.235\Phistar$. Line segment $\mathfrak{l}=\big\{(\Phi_1,\theta_1=\const)\big\}$ with endpoints $\mathcal{W}$, $\mathcal{E}$ in lower part.}\label{fig:x1_n1equal3case_subcases}
\end{center}
\end{minipage}\hfill
\end{figure}

\begin{proposition}\label{proposition:casecexpansionn1eq3}
    With the standing assumptions of this section (\cref{REMCaseCAssumptions}),if with $n_1=3$, then $\mathcal{II}_1(\Bcal^+_0,3)+\mathcal{II}_2(\Bcal^+_0,3)>4.28$ (see \cref{def:AformularForLowerBoundOfExpansionForn1largerthan0,lemma:AformularForLowerBoundOfExpansionForn1largerthan0}). Hence in (4) of \cref{lemma:AformularForLowerBoundOfExpansionForn1largerthan0}, $E(3,\Bcal^+_0)=\frac{\|dx_3\|_\p}{\|dx_0\|_\p}=3\big|\mathcal{II}_1(\mathcal{B}^+_0,3)+\mathcal{II}_2(\mathcal{B}^+_0,3)\big|>12.8$.
\end{proposition}
\begin{proof}
    We first estimate $\tau_1-\frac{2}{3}d_2$ for $x_1=(\Phi_1,\theta_1)$ with $n_1=3$. 
    \begin{equation}\label{eq:tau1minus2thirdsd2n1is3}
        \begin{aligned}
        &\left.\begin{aligned}
        &-\frac{2}{3}d_2\overbracket<^{\mathclap{\text{\cref{proposition:casecparametersestimate} (i)}}\strut}-\frac{2}{3}r\sin\phistar+\frac{2}{3}\cdot\frac{17r}{4R}r\sin\phistar\\
        &\tau_1\overbracket<^{\mathclap{\text{\cref{proposition:casecparametersestimate} \ref{itemcasecparametersestimate2}}}}2r\sin\phistar-6d_1-\tau_0+\frac{8.1r}{R}r\sin\phistar
        \end{aligned}
        \right\}\Rightarrow\tau_1-\frac{2}{3}d_2<\frac{4r\sin\phistar}{3}-6d_1-\tau_0+\frac{65.6r}{6R}r\sin\phistar
    \end{aligned}
    \end{equation}
    
    Then we analyze the two cases based on $\theta_1$ i.e. $x_1$ in two parts of \cref{fig:x1_n1equal3case_subcases}.
    
    \hidesubsection{\textbf{Case 1}} $0.235\Phistar\le\theta_1<\Phistar/3$ and $x_1$ is in the upper part of the $n_1=3$ quadrilateral  in \cref{fig:x1_n1equal3case_subcases}.
    \begin{align*}
        &\left.\begin{aligned}
        &0.235\Phistar\overbracket=^{\mathclap{\text{\cref{def:BasicNotations0}}}}0.235\sin^{-1}\big(\frac{r\sin\phistar}{R}\big)\overbracket<^{\mathclap{\text{\eqref{eqJZHypCond}: }R>1700r}\strut}0.2\\
        &\text{If }0<z<0.2\text{, then } \sinc(z)\dfn\frac{\sin z}z>\frac{1}{1.01}\\
        \end{aligned}
        \right\} \Rightarrow \frac{\sin{(0.235\Phistar)}}{0.235\Phistar}> \frac{1}{1.01} \\
        &\Rightarrow d_1=R\sin{\theta_1}\overbracket\ge^{\mathclap{0.235\Phistar\le\theta_1<\pi/2}}R\sin{(0.235\Phistar)}>\frac{0.235}{1.01}R\Phistar>\frac{0.235}{1.01}R\sin{\Phistar}=\frac{0.235}{1.01}r\sin{\phistar}\\
    \end{align*} 
    
    Hence in \textbf{case 1} we have
    \begin{equation}\label{eq:tau1minus2thirdsd2n1is3case1}
        \begin{aligned}
        \tau_1-\frac{2}{3}d_2&\overbracket{<}^{\mathclap{\eqref{eq:tau1minus2thirdsd2n1is3}}}\frac{4}{3}r\sin{\phistar}-6d_1-\tau_0+\frac{65.6r}{6R}r\sin{\phistar}\\
        &\overbracket{<}^{\mathllap{d_1>\frac{0.235}{1.01}r\sin\phistar}}\frac{4}{3}r\sin{\phistar}-\frac{6\times0.235}{1.01}r\sin{\phistar}+\frac{65.6r}{6R}r\sin{\phistar}\\
        &<(-0.0625+\frac{65.6r}{6R})r\sin{\phistar}\overbracket<^{\mathclap{\text{\eqref{eqJZHypCond}:} R>1700r}}-0.056r\sin{\phistar} 
        \end{aligned}
    \end{equation}
\hidesubsection{\textbf{Case 2}} $\theta_1\in \big(0.2\Phistar,0.235\Phistar\big)$ and $x_1$ is in the lower part of the $n_1=3$ quadrilateral in \cref{fig:x1_n1equal3case_subcases}.

In \cref{fig:x1_n1equal3case_subcases}, suppose that $x_1$ is on a line segment $\mathfrak{l}=\Big\{(
\Phi_1,\theta_1=\const
=\nu\Phistar)\bigm|\Phi_1\in\big[\Phi_{\mathcal{W}},\Phi_{\mathcal{E}}\big)\Big\}$ with end points $\mathcal{W}=(\Phi_{\mathcal{W}},\nu\Phistar)$ and $\mathcal{E}=(\Phi_{\mathcal{E}},\nu\Phistar)$ on the boundary of the lower part of the $n_1=3$ quadrilateral, $\Phi_{\mathcal{W}}\le\Phi_1<\Phi_{\mathcal{E}}<0$ and $\nu\in\big(0.2,0.235\big)$. Note that $\tau_0$ as a function of $x_1=(\Phi_1,\theta_1=\const)\in\mathfrak{l}$ and thus of $\Phi_1\in\big[\Phi_{\mathcal{W}},\Phi_{\mathcal{E}}\big)$.

In the coordinates of \cref{def:StandardCoordinateTable} for \cref{fig:n1_equal3_leftendpoint_case}, $(x_P,y_P)=P=p(x_1)=(R\sin{\Phi_1},b-R\cos{\Phi_1})$, $A=(x_A,y_A)=(-r\sin\phistar,-r\cos\phistar)$, $P_\mathcal{W}\dfn p(\mathcal{W})=(x_\mathcal{W},y_\mathcal{W})=(R\sin{(\Phi_{\mathcal{W}})},b-R\cos{(\Phi_{\mathcal{W}})})$. Since $-\pi/2<-\Phistar<\Phi_\mathcal{W}\le\Phi_1<0$ and by checking \cref{fig:n1_equal3_leftendpoint_case}, we have $x_A<x_\mathcal{W}\le x_P$ and $y_A>y_\mathcal{W}\ge y_P$.

Suppose that $(x_{T_0},y_{T_0})=T_0=p(x_0)=p(\Fcal^{-1}(x_1))$ with $x_0=(\phi_0,\theta_0)$. By \cref{contraction_region}, $x_0\in\Nin$. More specifically, we see $T_0$ is in some neighborhood of corner $A$ (see \cref{fig:n1_equal3_leftendpoint_case} for a belief). Furthermore, by \cref{proposition:NoutsizeInA1BCcases}\eqref{item02NoutsizeInA1BCcases} we have the following.
\[0<2\pi-\phistar-\phi_0<2.001\sqrt{4r/R}\overbracket{<}^{\mathllap{\text{\eqref{eqJZHypCond0}: }R>1700r}}0.0971\Longrightarrow 2\pi-\phistar>\phi_0>2\pi-\phistar-0.0971\overbracket{>}^{\mathclap{\phistar<\tan^{-1}(1/3)}}5.86>3\pi/2.\]

Hence, $(x_{T_0},y_{T_0})=(r\sin\phi_0,-r\cos\phi_0)$ satisfies $x_{T_0}<x_A$ and $y_{T_0}>y_A$, and
\begin{align*}
    &\left.\begin{aligned}
        &x_{T_0}<x_A<x_{\mathcal{W}}\le x_P\\
        &y_{T_0}>y_A>y_{\mathcal{W}}\ge y_P
        \end{aligned}
    \right\} \Rightarrow |AP_{\mathcal{W}}|=\sqrt{(x_A-x_\mathcal{W})^2+(y_A-y_\mathcal{W})^2}<\sqrt{(x_P-x_{T_0})^2+(y_{T_0}-y_P)^2}=|T_0P|=\tau_0
\end{align*}
So $\tau_0(x_1)>|AP_{\mathcal{W}}|$ for all $x_1\in\mathfrak{l}$. Since $p(\Fcal^4{\mathcal{W}})=B$, from \cref{fig:n1_equal3_leftendpoint_case} we get $|AP_\mathcal{W}|=2R\sin{\big(\frac{2\Phistar-8\theta_1}{2}\big)}$.\COMMENT{Do not forget punctuation.\wentao{OK}\Hrule\red\bfseries The purpose of the large overbracket is unclear. It's a jumble without indication of what is used how and where.\Hrule\large  Replace it by overbrackets on individual relations.\wentao{Done.}}

It is known that $\sinc{z}=\frac{\sin{z}}{z}>1/1.01$ for $z\in[0,0.2)$.

Since $\nu=\frac{\theta_1}{\Phistar}<0.235$, we get $0<\theta_1=\nu\overbracket{\Phistar}^{\strut\mathclap{r\sin\phistar=R\sin\Phistar\text{: }\Phistar=\sin^{-1}(r\sin\phistar/R)}}<0.235\sin^{-1}{(\frac{r\sin\phistar}{R})}\overbracket{<}^{\mathclap{\text{\eqref{eqJZHypCond}: } R>1700r}}0.2$ and further get the following.
\begin{equation}\label{eq:n1eq3d1bound}
d_1=R\sin\theta_1\overbracket{>}^{\mathllap{\sinc\theta_1>1/1.01}}\frac{R\theta_1}{1.01}\overbracket{=}^{\bigstrut\mathclap{\theta_1=\nu\Phistar}}\frac{R\nu\Phistar}{1.01}\overbracket{>}^{\bigstrut\mathclap{\Phistar>\sin\Phistar}}\frac{R\nu\sin\Phistar}{1.01}\overbracket{=}^{\bigstrut\mathclap{R\sin\Phistar=r\sin\phistar}}\frac{\nu}{1.01}r\sin\phistar
\end{equation}

Since $\nu=\frac{\theta_1}{\Phistar}>0.2$, we get $0<(1-4\nu)\Phistar<0.2\Phistar$, we get $\sinc\big[(1-4\nu)\Phistar\big]>\frac{1}{1.01}$ and we further get the following.
\begin{equation}\label{eq:n1eq3tau0bound}
\tau_0>|AP_\mathcal{W}|=2R\sin{\big(\frac{2\Phistar-8\theta_1}{2}\big)}\overbracket{=}^{\bigstrut\mathclap{\theta_1=\nu\Phistar}}2R\sin{\big[(1-4\nu)\Phistar\big]}\overbracket{>}^{\bigstrut\mathclap{\sinc\big[(1-4\nu)\Phistar\big]>\frac{1}{1.01}}}2R\frac{(1-4v)\Phistar}{1.01}>2R\frac{1-4v}{1.01}\sin{\Phistar}\overbracket{=}^{\bigstrut\mathclap{R\sin{\Phistar}=r\sin{\phistar}}}\frac{1-4v}{1.01}(2r\sin{\phistar}).
\end{equation}

Hence in \textbf{case 2} we have \begin{equation}\label{eq:tau1minus2thirdsd2n1is3case2}
        \begin{aligned}
        \tau_1-\frac{2}{3}d_2&\overbracket{<}^{\mathclap{\eqref{eq:tau1minus2thirdsd2n1is3}}}\frac{4}{3}r\sin{\phistar}-6d_1-\tau_0+\frac{65.6r}{6R}r\sin{\phistar}\\
        &\overbracket{<}^{\mathllap{\text{\eqref{eq:n1eq3tau0bound}: }\tau_0>\frac{(1-4\nu)2r\sin\phistar}{1.01}\text{, and \eqref{eq:n1eq3d1bound}: }d_1>\frac{\nu r\sin\phistar}{1.01}}}\frac{4}{3}r\sin{\phistar}-\frac{6\nu}{1.01}r\sin{\phistar}-\frac{1-4\nu}{1.01}2r\sin{\phistar}+\frac{65.6r}{R}r\sin{\phistar}\\
        &=(\frac{4}{3}-\frac{2}{1.01}+\frac{2\nu}{1.01})r\sin{\phistar}+\frac{65.6r}{6R}r\sin{\phistar}\\
        &\overbracket{<}^{\quad\mathllap{\nu<0.235\text{ and \eqref{eqJZHypCond}: }R>1700r}}-0.181r\sin\phistar+0.007r\sin\phistar<-0.17r\sin\phistar<-0.056r\sin\phistar.
        \end{aligned}
    \end{equation}

To summarize, both in \textbf{case 1} \eqref{eq:tau1minus2thirdsd2n1is3case1} and in \textbf{case 2} \eqref{eq:tau1minus2thirdsd2n1is3case2}, we get $\tau_1-\frac{2}{3}d_2<-0.056r\sin{\phistar}$. 
Furthermore,\COMMENT{\red\bfseries Fix this alignment jumble. Useless alignment, needlessly repeated arguments create confusion. Pay more attention to readability. Note that the last underbracket pretends that the first equation is not there. You made more work for yourself by not planning ahead. Instead there should be a reference to the top equation.\Hrule mathclap will help\wentao{much rewritten. but I feel we still need to align.}} 
\begin{multline*}
1+\tau_0\Bcal^+_0\underbracket{>}_{\mathclap{\text{\eqref{eq:II1B4estimate3}}}}1+\big[\frac{2}{n_1+1}+\frac{8.1r}{(n_1+1)R}\big]\times\big(\frac{-4r\sin\phistar}{3d_0}\big)\underbracket{=}_{\mathclap{n_1=3}}1+\big[0.5+\frac{8.1r}{4R}\big]\big(\frac{-4r\sin\phistar}{3d_0}\big)\underbracket{>}_{\mathrlap{\text{\cref{proposition:casecparametersestimate}\ref{itemcasecparametersestimate5}}}}1+\big(0.5+\frac{8.1r}{4R}\big)\frac{-4\times1.0026}{3}\underbracket{>}_{\mathclap{\text{\eqref{eqJZHypCond}: }R>1700r}}0.33,\\
\mathcal{II}_1(\mathcal{B}^+_0,n_1=3)=\overbracket{\big(1+\tau_0\mathcal{B}^+_0\big)}^{\mathllap{>0.33}}+\underbracket{\big(\tau_1-\frac{2}{3}d_2\big)}_{\mathllap{\text{\eqref{eq:tau1minus2thirdsd2n1is3case1}\eqref{eq:tau1minus2thirdsd2n1is3case2}: }<-0.056r\sin\phistar<0}}\underbracket{\Big[\big(2n_1+1-\frac{(2n_1+2)\tau_0}{d_1}\big)\mathcal{B}^+_0-\frac{2n_1+2}{d_1}\Big]\big|_{n_1=3}}_{\mathclap{\text{\eqref{eq:II1B4estimate2}: }<[-\frac{2}{1.01}n_1^2+0.7135n_1+\frac{4}{0.99}]\frac{1}{r\sin\phistar}\mid_{n_1=3}<\frac{-11.64}{r\sin\phistar}<0}}>0.33+0.056\times11.64>0.98,\\
\mathcal{II}_2(\mathcal{B}^+_0,n_1=3)\overbracket{>}^{\text{\eqref{eq:II2B4estimate1}}}(2n_1-2.693)\bigm|_{n_1=3}=3.307.
\end{multline*}
Therefore, $\mathcal{II}_1(\mathcal{B}^+_0,3)+\mathcal{II}_2(\mathcal{B}^+_0,3)>0.98+3.307=4.287$. This proves the claim.
\end{proof}
\begin{figure}[h]
\begin{center}
\begin{tikzpicture}[]
    \tkzDefPoint(0,0){Or};
    \tkzDefPoint(0,-4.8){Y};
    \pgfmathsetmacro{\rradius}{4.5};
    \pgfmathsetmacro{\phistardeg}{40};
    \pgfmathsetmacro{\XValueArc}{\rradius*sin(\phistardeg)};
    \pgfmathsetmacro{\YValueArc}{\rradius*cos(\phistardeg)};
    \pgfmathsetmacro{\Rradius}{
        \rradius*sin(\phistardeg)/sin(18)
    };
    \pgfmathsetmacro{\bdist}{(-\rradius*cos(\phistardeg))+cos(18)*\Rradius}
    \tkzDefPoint(0,\bdist){OR};
    \tkzDefPoint(\XValueArc,-1.0*\YValueArc){B};    \tkzDefPoint(-1.0*\XValueArc,-1.0*\YValueArc){A};
    \draw[name path=Cr,green,ultra thin] (Or) circle (\rradius);
    \tkzDrawArc[name path=CR](OR,A)(B);
    \pgfmathsetmacro{\PXValue}{-\Rradius*sin(14)};
    \pgfmathsetmacro{\PYValue}{\bdist-\Rradius*cos(14)};
    \tkzDefPoint(\PXValue,\PYValue){P};
    \pgfmathsetmacro{\PXValue}{-\Rradius*sin(6)};
    \pgfmathsetmacro{\PYValue}{\bdist-\Rradius*cos(6)};
    \tkzDefPoint(\PXValue,\PYValue){P1};
    \pgfmathsetmacro{\PXValue}{\Rradius*sin(2)};
    \pgfmathsetmacro{\PYValue}{\bdist-\Rradius*cos(2)};
    \tkzDefPoint(\PXValue,\PYValue){P2};
    \pgfmathsetmacro{\PXValue}{\Rradius*sin(10)};
    \pgfmathsetmacro{\PYValue}{\bdist-\Rradius*cos(10)};
    \tkzDefPoint(\PXValue,\PYValue){P3};
    \path[name path=Line_farT_P](P)-- +(162:3.5);
    \path[name intersections={of=Cr and Line_farT_P,by={T}}];
    \tkzDefPoint(0,\bdist-\Rradius){C};
    \draw[red,ultra thin,dashed](OR)--(P);
    \draw[red,ultra thin,dashed](OR)--(A);
    \draw[red,ultra thin,dashed](OR)--(B);
    \begin{scope}[ultra thin,decoration={markings,mark=at position 0.95 with {\arrow{>}}}] 
    \draw[blue,postaction={decorate}] (T)--(P);
    \draw[blue,postaction={decorate}] (P)--(P1);
    \draw[blue,postaction={decorate}] (P1)--(P2);
    \draw[blue,postaction={decorate}] (P2)--(P3);
    \draw[blue,postaction={decorate}] (P3)--(B);
    \end{scope}
    \tkzLabelPoint[above](P){$P_\mathcal{W}$};
    \tkzLabelPoint[below](A){$A$};
    \tkzLabelPoint[above left](T){$T_\mathcal{W}$};
    \tkzLabelPoint[above](P1){\small $p(\Fc(\mathcal{W}))$};
    \tkzLabelPoint[below](P2){\small $p(\Fc^2(\mathcal{W}))$};
    \tkzLabelPoint[above](P3){\small $p(\Fc^3(\mathcal{W}))$};
    \tkzLabelPoint[below right](B){\small $B=p(\Fc^4(\mathcal{W}))$};
    \tkzLabelPoint[below](Or){$O_r$};
    \tkzLabelPoint[above](OR){$O_R$};
    \tkzMarkAngle[arc=lll,ultra thin,size=4.5](A,OR,P);
    \tkzMarkAngle[arc=lll,ultra thin,size=2](P,OR,B);
    \tkzDrawPoints(A,B,P,Or,OR,T,P1,P2,P3);
    \node[] at ($(OR)+(-88:3.2)$)  {\Large $8\theta_1$};
    \node[] at ($(OR)+(-110:5.5)$)  { $2\Phistar-8\theta_1$};
\end{tikzpicture}
\caption{ For $\mathcal{W}=(\Phi_{\mathcal{W}},\theta_1)$ in \cref{fig:x1_n1equal3case_subcases},
$P_\mathcal{W}=p(\mathcal{W})$, $T_\mathcal{W}=p(\Fc^{-1}(\mathcal{W}))$, $p(\Fc^4(\mathcal{W}))=B$. $\measuredangle{AO_RP_\mathcal{W}}=2\Phistar-8\theta_1$, $\measuredangle{P_\mathcal{W}O_RB}=8\theta_1$, $|AP_\mathcal{W}|=2R\sin{(\frac{2\Phistar-8\theta_1}{2})}$}\label{fig:n1_equal3_leftendpoint_case}
\end{center}
\end{figure}

\subsection[Expansion in Case (c) of \eqref{eqMainCases} with \texorpdfstring{$n_1=2$}{n1=2}]{The expansion in Case (c) of \texorpdfstring{\eqref{eqMainCases}}{orbit transition} with \texorpdfstring{$n_1=2$}{n1=2}}\hfill
\begin{proposition}\label{proposition:casecexpansionn1eq2}
    With the standing assumptions of this section (\cref{REMCaseCAssumptions}),\COMMENT{\Large Not needed once we have a standing assumption, right?\wentao{OK, done.}} if $n_1=2$, then 
    $\mathcal{II}_1(\Bcal^+_0,2)+\mathcal{II}_2(\Bcal^+_0,2)>1.86$ (see \cref{def:AformularForLowerBoundOfExpansionForn1largerthan0,lemma:AformularForLowerBoundOfExpansionForn1largerthan0}). 
    
    Hence in \cref{lemma:AformularForLowerBoundOfExpansionForn1largerthan0}(4), $E(2,\Bcal^+_0)=\frac{\|dx_3\|_\p}{\|dx_0\|_\p}=3\big|\mathcal{II}_1(\mathcal{B}^+_0,2)+\mathcal{II}_2(\mathcal{B}^+_0,2)\big|>5.58$.
\end{proposition}
\begin{proof}
We consider two cases: $\tau_1-\frac{2}{3}d_2>0$ and $\tau_1-\frac{2}{3}d_2\le0$.

\hidesubsection[$\tau_1-\frac{2}{3}d_2>0$]{Case 1} We\COMMENT{\color{red}This use of \textbackslash hidesection is good for cases.\wentao{OK.}} first show that, in fact, $0<\tau_1-\frac{2}{3}d_2<0.06r\sin\phistar$:
\begin{align*}
        &\frac{2}{3}(r\sin{\phistar}-\frac{17r}{4R}r\sin{\phistar})+4d_1\overbracket{<}^{\mathclap{\strut\text{\cref{proposition:casecparametersestimate} \ref{itemcasecparametersestimate1}}}}\frac{2}{3}d_2+4d_1\overbracket{<}^{\mathclap{\tau_1>2d_2/3}}\tau_1+4d_1<\tau_1+4d_1+\tau_0\overbracket{<}^{\strut\mathclap{\text{\cref{proposition:casecparametersestimate} \ref{itemcasecparametersestimate2} with } n_1=2}}2r\sin\phistar+\frac{8.1r}{R}r\sin\phistar\\
        \Rightarrow&\frac{2}{3}r\sin{\phistar}-\frac{17r}{6R}r\sin{\phistar}+4d_1<2r\sin\phistar+\frac{8.1r}{R}r\sin\phistar\\
        \Rightarrow&4d_1<\frac{4}{3}r\sin{\phistar}+(\frac{17r}{6R}+\frac{8.1r}{R})r\sin{\phistar}\Rightarrow\frac{2}{3}d_1<\frac{2}{9}r\sin{\phistar}+(\frac{17r}{36R}+\frac{8.1r}{6R})r\sin{\phistar}\\
        \Rightarrow & 0<\tau_1-\frac{2}{3}d_2\overbracket{<}^{\strut\mathclap{\strut\text{\cref{proposition:casecparametersestimate} \ref{itemcasecparametersestimate7}}}}\frac{2}{3}d_1-\frac{1}{6}r\sin{\phistar}+\frac{16.6r}{3R}r\sin{\phistar}\overbracket{<}^{\bigstrut\mathclap{\frac{2}{3}d_1<\frac{2}{9}r\sin{\phistar}+(\frac{17r}{36R}+\frac{8.1r}{6R})r\sin{\phistar}}}\big(\frac{2}{9}-\frac{1}{6}\big)r\sin{\phistar}+\big(\frac{16.6r}{3R}+\frac{17r}{36R}+\frac{8.1r}{6R}\big)r\sin{\phistar}\overbracket{<}^{\mathclap{\text{\eqref{eqJZHypCond}: }R>1700r}}0.06r\sin\phistar,
    \end{align*}
as claimed. 
We now turn to $\mathcal{II}_1(\mathcal{B}^+_0,2)$ (see \eqref{eq:AformularForLowerBoundOfExpansionForn1largerthan0}). One term in it is
\begin{multline*}
        \big[(2n_1+1-\frac{(2n_1+2)\tau_0}{d_1})\Bcal^+_0-\frac{2n_1+2}{d_1}\big]\bigm|_{n_1=2}=5\Bcal^+_0-\Bcal^+_0\frac{6\tau_0}{d_1}-\frac{6}{d_1}\overbracket{>}^{\bigstrut\mathllap{\text{\cref{lemma:AformularForLowerBoundOfExpansionForn1largerthan0}(1): }\Bcal_0^+<0}}5\Bcal_0^+-\frac{6}{d_1}
        \overbracket{>}^{\bigstrut\mathclap{\substack{\text{\cref{lemma:AformularForLowerBoundOfExpansionForn1largerthan0}(1): }\Bcal_0^+\ge{-4}/{3d_0}\\\text{\cref{proposition:casecparametersestimate}\ref{itemcasecparametersestimate3} with }\\n_1=2\text{, }d_1>\frac{1}{4}r\sin\phistar}}}-\frac{5\times4}{3d_0}-\frac{6\times4}{r\sin\phistar}\\
        \underbracket{>}_{\strut\mathclap{\text{\cref{proposition:casecparametersestimate} \ref{itemcasecparametersestimate5}: }d_0>{r\sin\phistar}/{1.0026}}}\frac{-20\times1.0026}{3r\sin\phistar}-\frac{24}{r\sin\phistar}=\frac{-30.684}{r\sin\phistar}.
    \end{multline*}
So, \(\displaystyle\big(\tau_1-\frac{2}{3}d_2\big)\big[(2n_1+1-\frac{(2n_1+2)\tau_0}{d_1})\Bcal^+_0-\frac{2n_1+2}{d_1}\big]\bigm|_{n_1=2}\overbracket{>}^{\strut\mathclap{\tau_1-\frac{2}{3}d_2>0}}\big(\tau_1-\frac{2}{3}d_2\big)\big(\frac{-30.684}{r\sin\phistar}\big)\overbracket{>}^{\bigstrut\mathrlap{\tau_1-\frac{2}{3}d_2<0.06r\sin\phistar}}0.06\times(-30.684)>-1.842.\)
The\COMMENT{\red\bfseries Edit the next display the way I edited case 2 below\wentao{OK}} remaining term in $\mathcal{II}_1(\mathcal{B}^+_0,2)$  involves $\tau_0$, and we have the following.\COMMENT{Use overbracket and mathclap below.\wentao{Done.}} 
\[
\tau_0+2d_2\overbracket{<}^{\mathclap{\tau_1>\frac{2}{3}d_2}}\tau_0+3\tau_1=\tau_0+\tau_1+2\tau_1\overbracket{<}^{\mathllap{\text{\cref{remark:deffpclf}: }\tau_1<2d_1}}\tau_1+\tau_0+4d_1\overbracket{<}^{\strut\mathclap{\text{\cref{proposition:casecparametersestimate}\ref{itemcasecparametersestimate2} with }n_1=2}}2r\sin{\phistar}+\frac{8.1r}{R}r\sin{\phistar},
\]
hence
\[
\tau_0<2(r\sin\phistar-d_2)+\frac{8.1r}{R}r\sin{\phistar}\overbracket{<}^{\strut\mathclap{\text{\cref{proposition:casecparametersestimate}\ref{itemcasecparametersestimate1}}}}2\times\frac{17r}{4R}r\sin\phistar+\frac{8.1r}{R}r\sin\phistar=\frac{16.6r}{R}r\sin\phistar.
\]
So
\[
1+\tau_0\Bcal^+_0\overbracket{>}^{\mathclap{\Bcal^+_0<0}}1+\frac{16.6r}{R}\Bcal_0^+r\sin\phistar\overbracket{\ge}^{\mathclap{\Bcal^+_0\ge\frac{-4}{3d_0}}}1-\frac{4\times16.6r}{3R}\frac{r\sin{\phistar}}{d_0}\overbracket{>}^{\strut\mathclap{\text{\cref{proposition:casecparametersestimate}\ref{itemcasecparametersestimate5}}}}1-\frac{4\times16.6r}{3R}1.0026\overbracket{>}^{\mathclap{\eqref{eqJZHypCond}:R>1700r}}0.986
\]
and
\begin{multline*}
\mathcal{II}_2(B^+_0,2)\overbracket{=}^{\mathclap{\text{\cref{def:AformularForLowerBoundOfExpansionForn1largerthan0}}}}4\big[1+(\tau_0-d_1)\Bcal^+_0\big]=4(1+\tau_0\Bcal^+_0)-4d_1\Bcal^+_0>4\times0.986-4d_1\Bcal^+_0
        \\
        \overbracket{>}^{\strut\mathclap{\substack{\Bcal^+_0<0\text{ and \cref{proposition:casecparametersestimate}\ref{itemcasecparametersestimate3}:}\\n_1=2,\ d_1>(r\sin\phistar)/4}}}3.944-\Bcal^+_0r\sin\phistar
        \overbracket{\ge}^{\mathclap{-\Bcal^+_0\ge\frac{1}{d_0}}}3.944+\frac{r\sin\phistar}{d_0}\overbracket{>}^{\strut\mathrlap{\text{\cref{proposition:casecparametersestimate}\ref{itemcasecparametersestimate5}}}}3.944+0.9975=4.9415.
\end{multline*}
Thus, in case 1, we conclude
\begin{multline*}
\mathcal{II}_1(\mathcal{B}^+_0,2)+\mathcal{II}_2(\mathcal{B}^+_0,2)
=\overbracket{\mathcal{II}_2(\mathcal{B}^+_0,2)}^{>4.9415}+\overbracket{\big(1+\tau_0\mathcal{B}^+_0\big)}^{>0.986}+\overbracket{\big(\tau_1-\frac{2}{3}d_2\big)\Big[\big(2n_1+1-\frac{(2n_1+2)\tau_0}{d_1}\big)\mathcal{B}^+_0-\frac{2n_1+2}{d_1}\Big]\bigm|_{n_1=2}}^{>-1.842}
\\
>4.9415+0.986-1.842=4.0855>1.8698.
\end{multline*}
\hidesubsection[$\tau_1-\frac{2}{3}d_2\le0$]{Case 2}\strut\COMMENT{Use overbracket below.\wentao{Done.}}
\begin{align*}
    \big([(2n_1+1)-\frac{(2n_1+2)\tau_0}{d_1}]\Bcal^+_0-\frac{2n_1+2}{d_1}\big)&\overbracket=^{\strut\mathclap{n_1=2}}5\Bcal^+_0-\Bcal^+_0\frac{6\tau_0}{d_1}-\frac{6}{d_1}\overbracket{<}^{\mathclap{\substack{\text{\cref{remark:deffpclf}: }\tau_0<2d_1\\\text{and }\Bcal^+_0<0}}}-7\Bcal_0^+-\frac{6}{d_1}\overbracket{\le}^{\mathclap{-\Bcal^+_0\le\frac{4}{3d_0}}}-\frac{28}{3d_0}-\frac{6}{d_1}\\
    &\overbracket{<}^{\mathllap{\text{\cref{proposition:casecparametersestimate}\ref{itemcasecparametersestimate5}: }\frac{1}{d_0}<\frac{1.0026}{r\sin\phistar}\text{ and \ref{itemcasecparametersestimate3} with }n_1=2}}\frac{28\times1.0026}{3r\sin{\phistar}}-\frac{6}{(\frac{1}{2}+\frac{8.1r}{4R}r\sin{\phistar})}\overbracket{<}^{\strut\mathclap{\eqref{eqJZHypCond}:R>1700r}}\frac{-2.61}{r\sin\phistar}<0
\end{align*} Multiplying by $\tau_1-\frac{2}{3}d_2\le0$, we have 
\[\overbracket{\big(\tau_1-\frac{2}{3}d_2\big)}^{\le0}\overbracket{\big([(2n_1+1)-\frac{(2n_1+2)\tau_0}{d_1}]\Bcal^+_0-\frac{2n_1+2}{d_1}\big)}^{<0}\ge0.\] 
Next,\COMMENT{You are getting carried away with alignment acrobatics and it does not always help. I fixed most of the jumble, \red\large tidy up the corresponding part in case 1.\wentao{OK}} 
\(\displaystyle3\tau_0\overbracket{<}^{\mathclap{\text{\cref{def:fpclf}}:\ \tau_0<2d_1}}\tau_0+4d_1<\tau_0+4d_1+\tau_1\overbracket{<}^{\mathclap{\text{\cref{proposition:casecparametersestimate}\ref{itemcasecparametersestimate2}}}}2r\sin\phistar+\frac{8.1r}{R}r\sin\phistar\) implies \(\tau_0<\frac{2}{3}r\sin\phistar+\frac{2.7r}{R}r\sin\phistar\), so
\[
1+\tau_0\Bcal^+_0\overbracket{\ge}^{\mathclap{\Bcal^+_0\ge\frac{-4}{3d_0}}}1-\frac{4\tau_0}{3d_0}>1-\frac{4}{3}\frac{\frac{2}{3}r\sin{\phistar}+\frac{2.7r}{R}r\sin{\phistar}}{d_0}\overbracket{>}^{\bigstrut\mathclap{\text{\cref{proposition:casecparametersestimate}\ref{itemcasecparametersestimate5}}}}1-\frac{4}{3}(\frac{2}{3}+\frac{2.7r}{R})1.0026\overbracket{>}^{\mathclap{\text{\eqref{eqJZHypCond}}: R>1700r}}0.1066,
\]
and, again using \(\tau_0<\frac{2}{3}r\sin\phistar+\frac{2.7r}{R}r\sin\phistar\),
\begin{align*}\mathcal{II}_2(\Bcal^+_0,2)\overbracket{=}^{\mathclap{\text{\cref{def:AformularForLowerBoundOfExpansionForn1largerthan0}}}}&4\big(1+(\tau_0-d_1)\Bcal^+_0\big)\overbracket{>}^{\mathclap{\Bcal^+_0<0}}4\big(1+(\frac{2}{3}r\sin\phistar+\frac{2.7r}{R}r\sin\phistar-d_1)\Bcal^+_0\big)\\
\overbracket{>}^{\mathclap{\substack{\text{\cref{proposition:casecparametersestimate}\ref{itemcasecparametersestimate3}: }\\n_1=2, d_1>(r\sin\phistar)/4\\-\Bcal^+_0>0}}}&4\big(1+(\frac{2}{3}r\sin\phistar+\frac{2.7r}{R}r\sin\phistar-\frac{r\sin\phistar}{4})\Bcal^+_0\big)=4\Big[1+\big(\frac{5}{12}+\frac{2.7r}{R}\big)\Bcal^+_0r\sin\phistar\Big]\\
\overbracket{\ge}^{\mathclap{\Bcal^+_0\ge\frac{-4}{3d_0}}}&4\Big[1+\frac{-4}{3}\big(\frac{5}{12}+\frac{2.7r}{R}\big)\frac{r\sin\phistar}{d_0}\Big]\overbracket{>}^{\strut\mathclap{\text{\cref{proposition:casecparametersestimate}\ref{itemcasecparametersestimate5}}}}4\cdot\big[1-\frac{4\times1.0026}{3}\big(\frac{5}{12}+\frac{2.7r}{R}\big)\big]\overbracket{>}^{\hskip-2em\mathrlap{\text{\eqref{eqJZHypCond}: } R>1700r}}4\times0.4408=1.7632.\\
\end{align*}
Combining these, we get
\begin{multline*}
\mathcal{II}_1(\mathcal{B}^+_0,2)+\mathcal{II}_2(\mathcal{B}^+_0,2)=\overbracket{\mathcal{II}_2(\mathcal{B}^+_0,2)}^{>1.7632}+\overbracket{\big(1+\tau_0\mathcal{B}^+_0\big)}^{>0.1066}+\overbracket{\big(\tau_1-\frac{2}{3}d_2\big)\Big[\big(2n_1+1-\frac{(2n_1+2)\tau_0}{d_1}\big)\mathcal{B}^+_0-\frac{2n_1+2}{d_1}\Big]\bigm|_{n_1=2}}^{\ge0}
\\
>1.7632+0.1066+0=1.8698
\end{multline*}

Thus, $\mathcal{II}_1(\mathcal{B}^+_0,2)+\mathcal{II}_2(\mathcal{B}^+_0,2)>1.8698$ in both cases. This proves the claim.
\end{proof}

\subsection[The expansion in Case (c) of \eqref{eqMainCases} with \texorpdfstring{$n_1=1$}{n1=1}]{The expansion in Case (c) of \texorpdfstring{\eqref{eqMainCases}}{orbit transition} with \texorpdfstring{$n_1=1$}{n1=1}}\label{SBSn1=1}

\begin{definition}\label{def:II1kAndII2k}
    In \cref{def:AformularForLowerBoundOfExpansionForn1largerthan0},  with $\Bcal^+_0=\frac{-k}{d_0}$ for $k\in[1,\frac{4}{3}]$,\COMMENT{\Large Do we need any of the preceding to introduce this notation?\wentao{seems not, remove it?}} we let $\mathcal{II}_1(\frac{-k}{d_0},1)\nfd \mathcal{II}_1(k)$ and $\mathcal{II}_2(\frac{-k}{d_0},1)\nfd \mathcal{II}_1(k)$.\COMMENT{\Large What follows is not definition; take it out of the definition.\wentao{OK.}}
\end{definition}

From \eqref{eq:AformularForLowerBoundOfExpansionForn1largerthan0} get have the following formulas. \begin{equation}\label{eq:II1kAndII2k}
        \begin{aligned}
            \mathcal{II}_1(k)&=\mathcal{II}_1(\frac{-k}{d_0},1)=\Big[\big(1+\tau_0\Bcal^+_0\big)+\big(\tau_1-\frac{2}{3}d_2\big)\big[(2n_1+1-\frac{(2n_1+2)\tau_0}{d_1})B^+_0-\frac{2n_1+2}{d_1}\big]\Big]\bigm|_{n_1=1,\Bcal_0^+=\frac{-k}{d_0}}\\
            &=1+\frac{-k\tau_0}{d_0}+\big(\tau_1-\frac{2}{3}d_2\big)\big[(3-\frac{4\tau_0}{d_1})\Bcal^+_0-\frac{4}{d_1}\big]=1-\frac{k\tau_0}{d_0}+\big(\tau_1-\frac{2}{3}d_2\big)\big[(3-\frac{4\tau_0}{d_1})\frac{-k}{d_0}-\frac{4}{d_1}\big]\\
            &=1-k\frac{\tau_0}{d_0}+\big(\frac{\tau_1}{d_0}-\frac{2}{3}\frac{d_2}{d_0}\big)\big(\frac{3d_1}{d_0}-\frac{4\tau_0}{d_0}\big)\big(\frac{-kd_0}{d_1}\big)-\frac{4d_0}{d_1}\big(\frac{\tau_1}{d_0}-\frac{2}{3}\frac{d_2}{d_0}\big).\\
            \mathcal{II}_2(k)&=\mathcal{II}_2(\frac{-k}{d_0},1)=2n_1\big[1+(\tau_0-d_1)\Bcal^+_0\big]\bigm|_{n_1=1,\Bcal^+_0=\frac{-k}{d_0}}=2-\frac{2k\tau_0}{d_0}+\frac{2kd_1}{d_0}.
        \end{aligned}
    \end{equation}

\begin{definition}[Length fraction functions]\label{def:lengthfractionfunctionsn1eq1}
    With the standing assumptions of this section (\cref{REMCaseCAssumptions}) with $n_1=1$ and the length functions given in \cref{def:fpclf},\COMMENT{\Large Do we need any of the preceding to define these notations?? Omit if not.\wentao{OK.}} we define the length fraction functions to be \begin{enumerate}
        \item $h_1\dfn\frac{\tau_0}{d_0}+\frac{\tau_1}{d_0}+\frac{2d_1}{d_0}-2$\label{item01lengthfractionfunctionsn1eq1}
        \item $g_0\dfn d_0-r\sin{\phistar}$\label{item02lengthfractionfunctionsn1eq1}
        \item $g_1\dfn\frac{1}{d_0}-\frac{1}{r\sin{\phistar}}$\label{item03lengthfractionfunctionsn1eq1}
        \item $g_2\dfn d_2-r\sin{\phistar}$\label{item04lengthfractionfunctionsn1eq1}
        \item $g_3\dfn\frac{1}{d_2}-\frac{1}{r\sin{\phistar}}$.\label{item05lengthfractionfunctionsn1eq1}
    \end{enumerate}
 These are functions defined for $x_1$ in the quadrilateral with $n_1=1$ in \cref{fig:x1_n1cases_caseC}, also shown as \textcolor{blue}{$M^{\text{\upshape{in}}}_{R,1}$} in \cref{fig:x1_n1equal1case}.\COMMENT{\Large It sems that instead, they are defined whenever the denominators are not zero. We only USE them under some assumptions; those should be stated when needed.} 
\end{definition}
\begin{lemma}[Estimates for length fraction functions]\label{lemma:h1g1g0g1g2g3lengthfractionestimate}
    With the standing assumptions of this section (\cref{REMCaseCAssumptions}) and under condition \eqref{eqJZHypCond}, for the functions in \cref{def:lengthfractionfunctionsn1eq1} defined for $x_1$ in the $n_1=1$ region in \cref{fig:x1_n1cases_caseC} or equivalently the \textcolor{blue}{$M^{\text{\upshape{in}}}_{R,1}$} component in \cref{fig:x1_n1equal1case},\COMMENT{Edit to remove any standing assumptions.} we have \begin{enumerate}
        \item\label{item01h1g1g0g1g2g3lengthfractionestimate} $-8.53r/R<h_1<16.65r/R$,
        \item\label{item02h1g1g0g1g2g3lengthfractionestimate} $-\frac{17r}{4R}r\sin\phistar<g_0<\frac{17r}{4R}r\sin\phistar$,
        \item\label{item03h1g1g0g1g2g3lengthfractionestimate} $-\frac{17r}{4R}\frac{1}{0.9975r\sin\phistar}<g_1<\frac{17r}{4R}\frac{1}{0.9975r\sin\phistar}$,
        \item\label{item04h1g1g0g1g2g3lengthfractionestimate} $-\frac{17r}{4R}r\sin\phistar<g_2<\frac{17r}{4R}r\sin\phistar$,
        \item\label{item05h1g1g0g1g2g3lengthfractionestimate} $-\frac{17r}{4R}\frac{1}{0.9975r\sin\phistar}<g_3<\frac{17r}{4R}\frac{1}{0.9975r\sin\phistar}$
        \item\label{item06h1g1g0g1g2g3lengthfractionestimate} $-\frac{17r}{2R}r\sin\phistar<d_0-d_2=g_0-g_2<\frac{17r}{2R}r\sin\phistar$
    \end{enumerate}
\end{lemma}
\begin{proof} 
We provide item-by-item proofs.
\hidesubsection{Proof of (1)}$h_1\overbracket{=}^{\quad\mathllap{\text{\cref{def:lengthfractionfunctionsn1eq1}\eqref{item01lengthfractionfunctionsn1eq1}}}}\frac{\overbracket{\tau_0+\tau_1+2d_1}^{\mathrlap{\substack{\text{\cref{proposition:casecparametersestimate}\ref{itemcasecparametersestimate2} with }n_1=1: \\\in(2r\sin\phistar,\quad 2r\sin\phistar+(\frac{8.1r}{R})r\sin\phistar)}}\qquad\qquad\qquad\qquad}}{d_0}-2\in\big(\frac{2r\sin\phistar}{d_0}-2,\frac{2r\sin\phistar}{d_0}+\frac{8.1r}{R}\frac{r\sin\phistar}{d_0}-2\big),$
    
    Since 
    \begin{align*}
        \frac{2r\sin\phistar}{d_0}-2&=2\frac{\overbracket{r\sin\phistar-d_0}^{\mathclap{\substack{\text{\cref{proposition:casecparametersestimate}\ref{itemcasecparametersestimate1}: }\\ r\sin\phistar-d_0>-\frac{17r}{4R}r\sin\phistar}}}}{d_0}>-\frac{17r}{2R}\frac{r\sin\phistar}{d_0}\overbracket{>}^{\strut\mathclap{\text{\cref{proposition:casecparametersestimate}\ref{itemcasecparametersestimate5}}}}-\frac{17r}{2R}\frac{1}{0.9975}>-8.53r/R,\\
        \frac{2r\sin\phistar}{d_0}+\frac{8.1r}{R}\frac{r\sin\phistar}{d_0}-2&=2\frac{\overbracket{r\sin\phistar-d_0}^{\mathclap{\substack{\text{\cref{proposition:casecparametersestimate}\ref{itemcasecparametersestimate1}: }\\ r\sin\phistar-d_0<\frac{17r}{4R}r\sin\phistar}}}}{d_0}+\frac{8.1r}{R}\frac{r\sin\phistar}{d_0}<\frac{17r}{2R}\frac{r\sin\phistar}{d_0}+\frac{8.1r}{R}\frac{r\sin\phistar}{d_0}\\
        &\overbracket{<}^{\mathllap{\text{\cref{proposition:casecparametersestimate}\ref{itemcasecparametersestimate5}}}}\frac{17r}{2\times0.9975R}+\frac{8.1r}{0.9975R}<16.65r/R.
    \end{align*}
    Therefore, $-8.53r/R<h<16.65r/R$.
    \hidesubsection{Proof of (2) and (4)}\cref{proposition:casecparametersestimate}\ref{itemcasecparametersestimate1} implies that $|g_0|=|d_0-r\sin\phistar|<\frac{17r}{4R}r\sin\phistar$. Likewise for (4).
    \hidesubsection{Proof of (3) and (5)} $|g_1|\overbracket{=}^{\mathclap{\text{\cref{def:lengthfractionfunctionsn1eq1}\eqref{item01lengthfractionfunctionsn1eq1}}\quad}}\frac{|d_0-r\sin\phistar|}{d_0r\sin\phistar}\overbracket{<}^{\mathclap{\quad\text{\cref{proposition:casecparametersestimate}\ref{itemcasecparametersestimate1}}}}\frac{17r}{4R}\frac{r\sin\phistar}{d_0r\sin\phistar}=\frac{17r}{4R}\frac{1}{d_0}<\frac{17r}{4R}\frac{1}{0.9975r\sin\phistar}$. Likewise for (5).
    \hidesubsection{Proof of (6)} $g_0-g_2\overbracket{=}^{\mathclap{\text{\cref{def:lengthfractionfunctionsn1eq1}\eqref{item02lengthfractionfunctionsn1eq1}\eqref{item04lengthfractionfunctionsn1eq1}}}}d_0-d_2$ and by the triangle inequality 
    \[|d_0-d_2|\le|d_0-r\sin\phistar|+|r\sin\phistar-d_2|\overbracket{<}^{\mathclap{\text{\cref{proposition:casecparametersestimate}\ref{itemcasecparametersestimate1}}}}\frac{17r}{4R}r\sin\phistar+\frac{17r}{4R}r\sin\phistar=\frac{17r}{2R}r\sin\phistar.\qedhere\]

\end{proof}
\begin{proposition}\label{proposition:II1II2C1C2}
    With the standing assumptions of this section (\cref{REMCaseCAssumptions}), under condition \eqref{eqJZHypCond}, for the $\mathcal{II}_1(k)$, $\mathcal{II}_2(k)$ in \cref{def:II1kAndII2k} with $k\in[1,4/3]$ and $h_1$, $g_0$, $g_1$, $g_2$, $g_3$ from \cref{def:lengthfractionfunctionsn1eq1} defined for $x_1$ with $n_1=1$ in \cref{fig:x1_n1cases_caseC} or equivalently in the \textcolor{blue}{$M^{\text{\upshape{in}}}_{R,1}$} component in \cref{fig:x1_n1equal1case}, we have\COMMENT{Needs editing, but only after we decided about any standing assumptions.\wentao{I have edited and make the assumption in \cref{REMCaseCAssumptions}.}}
    \begin{equation}\label{eq:II1II2C1C2}
        \begin{aligned}
            \mathcal{II}_1(k)&=\frac{3+16k}{3}+\big(\frac{8-16k}{3}\big)\frac{r\sin{\phistar}}{d_1}+\frac{2kd_1}{r\sin{\phistar}}+\Big[\frac{(32k-12)r\sin{\phistar}}{3d_1}-10k\Big]\big(\frac{\tau_1}{d_0}\big)\\
            &\quad-\frac{4kr\sin{\phistar}}{d_1}\big(\frac{\tau_1}{d_0}\big)^2+\mathcal{C}_1(k)\\
            \mathcal{II}_2(k)&=2-4k+6k\frac{d_1}{r\sin{\phistar}}+2k\frac{\tau_1}{d_0}+\mathcal{C}_2(k)\\
            \mathcal{II}_1(k)+\mathcal{II}_2(k)&=\frac{9+4k}{3}+\big(\frac{8-16k}{3}\big)\frac{r\sin\phistar}{d_1}+\frac{8kd_1}{r\sin\phistar}+\Big[\frac{(32k-12)r\sin{\phistar}}{3d_1}-8k\Big]\big(\frac{\tau_1}{d_0}\big)\\
            &\quad-\frac{4kr\sin{\phistar}}{d_1}\big(\frac{\tau_1}{d_0}\big)^2+\mathcal{C}_1(k)+\mathcal{C}_2(k),
        \end{aligned}
    \end{equation}
    where \begin{equation}\label{eq:C1C2}
    \begin{aligned}
        \mathcal{C}_1(k)&=2kd_1g_1-\frac{2(g_2-g_0)}{3d_0}\frac{8kg_0}{d_1}+(\frac{\tau_1}{d_0}-\frac{2}{3})8k\frac{g_0}{d_1}-\frac{2(g_2-g_0)}{3d_0}(-11k+8k\frac{r\sin{\phistar}}{d_1})\\
        &\quad+(-4k\frac{\tau_1}{d_0})[-\frac{r\sin{\phistar}}{d_1}\frac{2(g_2-g_0)}{3d_0}+\frac{g_0}{d_1}(\frac{\tau_1}{d_0}-\frac{2}{3})-\frac{g_0}{d_1}\frac{2(g_2-g_0)}{3d_0}]-\frac{4g_0}{d_1}(\frac{\tau_1}{d_0}-\frac{2}{3})\\
        &\quad+\frac{8(g_2-g_0)}{3d_1}-kh_1+4h_1(\frac{\tau_1}{d_0}-\frac{2d_2}{3d_0})\frac{kd_0}{d_1},\\
        \mathcal{C}_2(k)&=6kd_1g_1-2kh_1.
    \end{aligned}
    \end{equation}
\end{proposition}
\begin{remark}\label{remark:II1II2Quadratic}
    For fixed $k\in[1,4/3]$, $\frac{d_1}{r\sin\phistar}=\const$, we define a quadratic function by\COMMENT{Use the same Q as previously.\wentao{OK.}} \begin{equation}\label{eq:II1II2Quadratic}
        \mathfrak{Q}(t)\dfn-\frac{4kr\sin\phistar}{d_1}t^2+\big[\frac{(32k-12)r\sin\phistar}{3d_1}-8k\big]t+\big(\frac{8-16k}{3}\big)\frac{r\sin\phistar}{d_1}+\frac{8kd_1}{r\sin{\phistar}}+\frac{4k+9}{3}.
    \end{equation}
\emph{\bfseries Alternate version; choose one:}    For fixed $k\in[1,4/3]$, define a quadratic function by \begin{equation}\label{eq:II1II2Quadratic2}
        \mathfrak{Q}(t)\dfn-12kDt^2+\big[(32k-12)D-8k\big]t+(8-16k)D+\frac{8k}{3D}+\frac{4k+9}{3},\quad\text{where}\quad D\dfn\frac{r\sin\phistar}{3d_1}.
    \end{equation}
    Then in \cref{proposition:II1II2C1C2}, \begin{equation}\label{eq:II1II2QuadraticAndC1C2}\mathcal{II}_1(k)+\mathcal{II}_2(k)=\mathfrak{Q}\big(\frac{\tau_1}{d_0}\big)+\mathcal{C}_1(k)+\mathcal{C}_2(k).\end{equation} 
    And we will see in \cref{proposition:C1C2lowerbound} that the lower bound of $\mathcal{C}_1(k)$ and $\mathcal{C}_2(k)$ can be arbitrarily close to $0$ as $R/r$ is large enough, and the same holds for the upper bound. 
\end{remark}
\begin{proof}
    We first compute $\mathcal{II}_1(k)$.\COMMENT{\small Lots of mathclap!\Hrule Always with overbracket or underbracket. Almost always without \textbackslash quad or 
    \textbackslash qquad} 
    \begin{align*}
        \mathcal{II}_1(k)&\overbracket{=}^{\mathclap{\text{\eqref{eq:II1kAndII2k}}}}1-k\frac{\tau_0}{d_0}+\big(\frac{\tau_1}{d_0}-\frac{2}{3}\frac{d_2}{d_0}\big)\big(\frac{3d_1}{d_0}-\frac{4\tau_0}{d_0}\big)\big(\frac{-kd_0}{d_1}\big)-\frac{4d_0}{d_1}\big(\frac{\tau_1}{d_0}-\frac{2}{3}\frac{d_2}{d_0}\big)\\
        &\overbracket{=}^{\mathclap{\substack{\text{\cref{def:lengthfractionfunctionsn1eq1} \eqref{item01lengthfractionfunctionsn1eq1}:}\\\frac{\tau_0}{d_0}=2+h_1-\frac{\tau_1}{d_0}-\frac{2d_1}{d_0}}}}1-k\big(2+h_1-\frac{\tau_1}{d_0}-\frac{2d_1}{d_0}\big)+\big(\frac{\tau_1}{d_0}-\frac{2}{3}\frac{d_2}{d_0}\big)\big[\frac{3d_1}{d_0}-4(2+h_1-\frac{\tau_1}{d_0}-\frac{2d_1}{d_0})\big]\big(\frac{-kd_0}{d_1}\big)-\frac{4d_0}{d_1}\big(\frac{\tau_1}{d_0}-\frac{2}{3}\frac{d_2}{d_0}\big)\\
        &=1-k\big(2+h_1-\frac{\tau_1}{d_0}-\frac{2d_1}{d_0}\big)+\big(\frac{\tau_1}{d_0}-\frac{2}{3}\frac{d_2}{d_0}\big)\big(\frac{3d_1}{d_0}-8-4h_1+\frac{4\tau_1}{d_0}+\frac{8d_1}{d_0}\big)\big(\frac{-kd_0}{d_1}\big)-\frac{4d_0}{d_1}\big(\frac{\tau_1}{d_0}-\frac{2}{3}\frac{d_2}{d_0}\big)\\
        &=1-k\big(2-\frac{\tau_1}{d_0}-\frac{2d_1}{d_0}\big)-kh_1+\big(\frac{\tau_1}{d_0}-\frac{2d_2}{3d_0}\big)\big(\frac{11d_1}{d_0}-8+\frac{4\tau_1}{d_0}\big)\big(\frac{-kd_0}{d_1}\big)\\
        &\quad+4h_1\big(\frac{\tau_1}{d_0}-\frac{2d_2}{3d_0}\big)\frac{kd_0}{d_1}-\frac{4d_0}{d_1}\big(\frac{\tau_1}{d_0}-\frac{2}{3}\frac{d_2}{d_0}\big),\\
    \end{align*}
    Using the functions from \cref{def:lengthfractionfunctionsn1eq1}, we split and recombine this:\COMMENT{Fix alignment. See next comment\wentao{Done.}}
    \begin{align*}
        \mathcal{II}_1(k)&=1-2k+k\frac{\tau_1}{d_0}+2k\frac{d_1}{d_0}+\big(\frac{\tau_1}{d_0}-\frac{2d_2}{3d_0}\big)\big(-11k+8k\frac{d_0}{d_1}\big)+\big(\frac{\tau_1}{d_0}-\frac{2d_2}{3d_0}\big)\big(\frac{-4k\tau_1}{d_0}\big)\big(\frac{d_0}{d_1}\big)-\frac{4d_0}{d_1}\big(\frac{\tau_1}{d_0}-\frac{2}{3}\frac{d_2}{d_0}\big)\\
        &\quad-kh_1+4h_1\big(\frac{\tau_1}{d_0}-\frac{2d_2}{3d_0}\big)
        \\
        &=1-2k+k\frac{\tau_1}{d_0}+2kd_1\big(\overbracket{\frac{1}{r\sin{\phistar}}+g_1}^{\mathclap{\text{\cref{def:lengthfractionfunctionsn1eq1}\eqref{item03lengthfractionfunctionsn1eq1}: }=\frac{1}{d_0}}}\big)+\big[\frac{\tau_1}{d_0}-\frac{2}{3}-\frac{2\overbracket{(g_2-g_0)}^{\mathclap{\substack{\text{\cref{def:lengthfractionfunctionsn1eq1}\eqref{item02lengthfractionfunctionsn1eq1}, \eqref{item04lengthfractionfunctionsn1eq1}: }\\ g_2-g_0=d_2-d_0}}}}{3d_0}\big]\big(-11k+\overbracket{8k\frac{r\sin{\phistar}}{d_1}+8k\frac{g_0}{d_1}}^{\mathclap{\text{\cref{def:lengthfractionfunctionsn1eq1}\eqref{item02lengthfractionfunctionsn1eq1}: }=8k\frac{d_0}{d_1}}}\big)
        \\
        &\quad+\big[\frac{\tau_1}{d_0}-\frac{2}{3}-\frac{2\overbracket{(g_2-g_0)}^{\mathclap{\text{\cref{def:lengthfractionfunctionsn1eq1}\eqref{item02lengthfractionfunctionsn1eq1}\eqref{item04lengthfractionfunctionsn1eq1}: } g_2-g_0=d_2-d_0\quad}}}{3d_0}\big]\big(-4k\frac{\tau_1}{d_0}\big)\big(\overbracket{\frac{r\sin{\phistar}}{d_1}+\frac{g_0}{d_1}}^{\mathclap{\text{\cref{def:lengthfractionfunctionsn1eq1}\eqref{item02lengthfractionfunctionsn1eq1}: }=\frac{d_0}{d_1}}}\big)-\big(\frac{4g_0}{d_1}+\frac{4r\sin{\phistar}}{d_1}\big)\big(\frac{\tau_1}{d_0}-\frac{2}{3}-\frac{2}{3}\frac{\overbracket{g_2-g_0}^{\mathclap{\text{\cref{def:lengthfractionfunctionsn1eq1}\eqref{item02lengthfractionfunctionsn1eq1}, \eqref{item04lengthfractionfunctionsn1eq1}: } g_2-g_0=d_2-d_0}}}{d_0}\big)
        \\
        &\quad-kh_1+4h_1\big(\frac{\tau_1}{d_0}-\frac{2d_2}{3d_0}\big)
        \\
        &\overbracket{=}^{\strut\mathclap{\substack{\text{Split  into terms}\\\text{with and without}\\ g_0,g_1,g_2,h_1}}}1-2k+k\frac{\tau_1}{d_0}+\frac{2kd_1}{r\sin{\phistar}}+\big(\frac{\tau_1}{d_0}-\frac{2}{3}\big)\big(-11k+8k\frac{r\sin{\phistar}}{d_1}\big)+\big(\frac{\tau_1}{d_0}-\frac{2}{3}\big)\big(-4k\frac{\tau_1}{d_0}\big)\frac{r\sin{\phistar}}{d_1}-\frac{4r\sin{\phistar}}{d_1}\big(\frac{\tau_1}{d_0}-\frac{2}{3}\big)\\
        &\quad+2kd_1g_1-\frac{2(g_2-g_0)}{3d_0}\big(-11k+8k\frac{r\sin\phistar}{d_1}+8k\frac{g_0}{d_1}\big)+\big(\frac{\tau_1}{d_0}-\frac{2}{3}\big)8k\frac{g_0}{d_1}+\frac{2(g_2-g_0)}{3d_0}\frac{4k\tau_1}{d_0}\big(\frac{r\sin\phistar}{d_1}+\frac{g_0}{d_1}\big)\\
        &\quad+\big(\frac{\tau_1}{d_0}-\frac{2}{3}\big)\big(-4k\frac{\tau_1}{d_0}\big)\frac{g_0}{d_1}-\frac{4g_0}{d_1}\big(\frac{\tau_1}{d_0}-\frac{2}{3}-\frac{2}{3}\frac{g_2-g_0}{d_0}\big)+\frac{4r\sin\phistar}{d_1}\frac{2(g_2-g_0)}{3d_0}\\
        &\quad-kh_1+4h_1\big(\frac{\tau_1}{d_0}-\frac{2d_2}{3d_0}\big)\\
        &=1-2k+\frac{2kd_1}{r\sin\phistar}+\frac{22k}{3}-\frac{16kr\sin\phistar}{3d_1}+\frac{8r\sin\phistar}{3d_1}+\big(k-11k+8k\frac{r\sin\phistar}{d_1}+\frac{8kr\sin\phistar}{3d_1}-\frac{4r\sin\phistar}{d_1}\big)\big(\frac{\tau_1}{d_0}\big)\\
        &\quad-\frac{4kr\sin\phistar}{d_1}\big(\frac{\tau_1}{d_0}\big)^2+\mathcal{C}_1(k)\\
        &=\frac{3+16k}{3}+\big(\frac{8-16k}{3}\big)\frac{r\sin{\phistar}}{d_1}+\frac{2kd_1}{r\sin{\phistar}}+\Big[\frac{(32k-12)r\sin{\phistar}}{3d_1}-10k\Big]\frac{\tau_1}{d_0}-\frac{4kr\sin{\phistar}}{d_1}\big(\frac{\tau_1}{d_0}\big)^2+\mathcal{C}_1(k),
    \end{align*}where (again, for ease of reference)
    \begin{multline*}
    \mathcal{C}_1(k)=2kd_1g_1-\frac{2(g_2-g_0)}{3d_0}\big(-11k+8k\frac{r\sin\phistar}{d_1}+8k\frac{g_0}{d_1}\big)+\big(\frac{\tau_1}{d_0}-\frac{2}{3}\big)8k\frac{g_0}{d_1}+\frac{2(g_2-g_0)}{3d_0}\frac{4k\tau_1}{d_0}\big(\frac{r\sin\phistar}{d_1}+\frac{g_0}{d_1}\big)
    \\
    +\big(\frac{\tau_1}{d_0}-\frac{2}{3}\big)\big(-4k\frac{\tau_1}{d_0}\big)\frac{g_0}{d_1}-\frac{4g_0}{d_1}\big(\frac{\tau_1}{d_0}-\frac{2}{3}-\frac{2}{3}\frac{g_2-g_0}{d_0}\big)+\frac{4r\sin\phistar}{d_1}\frac{2(g_2-g_0)}{3d_0}-kh_1+4h_1\big(\frac{\tau_1}{d_0}-\frac{2d_2}{3d_0}\big).
\end{multline*}

Similarly by introducing functions in \cref{def:lengthfractionfunctionsn1eq1}, we also have 
\begin{align*}
    \mathcal{II}_2(k)&\overbracket{=}^{\mathclap{\text{\eqref{eq:II1kAndII2k}}}}2-\frac{2k\tau_0}{d_0}+\frac{2kd_1}{d_0}\\
    &\overbracket{=}^{\mathclap{\substack{\text{\cref{def:lengthfractionfunctionsn1eq1} \eqref{item01lengthfractionfunctionsn1eq1}: }\\\frac{\tau_0}{d_0}=2+h_1-\frac{\tau_1}{d_0}-\frac{2d_1}{d_0}}}}2-2k\big(2-\frac{\tau_1}{d_0}-\frac{2d_1}{d_0}+h_1\big)+\frac{2kd_1}{d_0}=2-4k+\frac{6kd_1}{d_0}+\frac{2k\tau_1}{d_0}-2kh_1
    \\
    &\overbracket{=}^{\mathclap{\substack{\text{\cref{def:lengthfractionfunctionsn1eq1}\eqref{item03lengthfractionfunctionsn1eq1}}:\\\frac{1}{d_0}=\frac{1}{r\sin\phistar}+g_1}}}2-4k+6kd_1\big(\frac{1}{r\sin{\phistar}}+g_1\big)+\frac{2k\tau_1}{d_0}-2kh_1\\
    &=2-4k+6k\frac{d_1}{r\sin\phistar}+2k\frac{\tau_1}{d_0}+\underbracket{6kd_1g_1-2kh_1}_{\mathcal{C}_2(k)}=2-4k+6k\frac{d_1}{r\sin{\phistar}}+2k\frac{\tau_1}{d_0}+\mathcal{C}_2(k).\qedhere
\end{align*}
\end{proof}
\begin{proposition}\label{proposition:C1C2lowerbound}
    For $k\in[1,4/3]$, $\mathcal{C}_1(k)$ and $\mathcal{C}_2(k)$ from \eqref{eq:C1C2} 
    satisfy 
    \[\mathcal{C}_1(k)+\mathcal{C}_2(k)>-389\frac{r}{R}-483\big(\frac{r}{R}\big)^2.\] 
\end{proposition}
\begin{proof}
    Since  $\frac{1}{3}<\frac{d_1}{r\sin\phistar}<1$ by \cref{proposition:QUpperLowerRegionsLowerBoundsChatCtilde} \eqref{item01Qlowerboundn1eq1}, when $x_1$ is in the $n_1=1$ region of \cref{fig:x1_n1cases_caseC}, we have\COMMENT{\Large\color{red}Bad \TeX\  probably throughout the paper:  ``=\&'' This should {\bfseries ALWAYS} be ``\&=''---and then it should always be ``\&\textbackslash quad+'' to move the ``+'' to the right!\wentao{Checking.}}
    \begin{equation}\label{eq:C1C2reformed}
    \begin{aligned}
    \mathcal{C}_1(k)+\mathcal{C}_2(k)&\overbracket{=}^{\mathclap{\text{\eqref{eq:C1C2}}}}2kd_1g_1-\frac{2(g_2-g_0)}{3d_0}\big(-11k+8k\frac{r\sin\phistar}{d_1}+8k\frac{g_0}{d_1}\big)+\big(\frac{\tau_1}{d_0}-\frac{2}{3}\big)8k\frac{g_0}{d_1}+\frac{2(g_2-g_0)}{3d_0}\frac{4k\tau_1}{d_0}\big(\frac{r\sin\phistar}{d_1}+\frac{g_0}{d_1}\big)
    \\
        &\quad+\big(\frac{\tau_1}{d_0}-\frac{2}{3}\big)\big(-4k\frac{\tau_1}{d_0}\big)\frac{g_0}{d_1}-\frac{4g_0}{d_1}\big(\frac{\tau_1}{d_0}-\frac{2}{3}-\frac{2}{3}\frac{g_2-g_0}{d_0}\big)+\frac{4r\sin\phistar}{d_1}\frac{2(g_2-g_0)}{3d_0}-kh_1+4h_1\big(\frac{\tau_1}{d_0}-\frac{2d_2}{3d_0}\big)
        \\
        &\quad+6kd_1g_1-2kh_1
        \\
        &=8kd_1g_1+h_1\big(-3k+\frac{\tau_1}{d_0}-\frac{2d_2}{3d_0}\big)+\frac{2(g_2-g_0)}{3d_0}\big(\frac{4r\sin\phistar}{d_1}+\frac{4k\tau_1}{d_0}\frac{r\sin\phistar}{d_1}+\frac{4k\tau_1}{d_0}\frac{g_0}{d_1}\big)+\frac{2(g_2-g_0)}{3d_0}\frac{4g_0}{d_1}
        \\
        &\quad-\frac{2(g_2-g_0)}{3d_0}\big(-11k+8k\frac{r\sin\phistar}{d_1}+8k\frac{g_0}{d_1}\big)
        -\frac{4g_0}{d_1}\big(\frac{\tau_1}{d_0}-\frac{2}{3}\big)+\big(\frac{\tau_1}{d_0}-\frac{2}{3}\big)\big(-4k\frac{\tau_1}{d_0}\big)\frac{g_0}{d_1}+\big(\frac{\tau_1}{d_0}-\frac{2}{3}\big)8k\frac{g_0}{d_1}
        \\
        &=8kd_1g_1+h_1\big(-3k+\frac{\tau_1}{d_0}-\frac{2d_2}{3d_0}\big)+\big(\frac{\tau_1}{d_0}-\frac{2}{3}\big)\frac{g_0}{d_1}\big(8k-4-4k\frac{\tau_1}{d_0}\big)
        \\
        &\quad+\frac{2(g_2-g_0)}{3d_0}\big(\frac{4r\sin\phistar}{d_1}+\frac{4k\tau_1}{d_0}\frac{r\sin\phistar}{d_1}+\frac{4k\tau_1+4d_0}{d_0}\frac{g_0}{d_1}\big)-\frac{2(g_2-g_0)}{3d_0}\big(-11k+8k\frac{r\sin\phistar}{d_1}+8k\frac{g_0}{d_1}\big)\\
        &=8kd_1g_1+h_1\big(-3k+\frac{\tau_1}{d_0}-\frac{2d_2}{3d_0}\big)+\big(\frac{\tau_1}{d_0}-\frac{2}{3}\big)\frac{g_0}{d_1}\big(8k-4-4k\frac{\tau_1}{d_0}\big)\\
        &\quad+\frac{2(g_2-g_0)}{3d_0}\Big[11k-\frac{8kr\sin\phistar}{d_1}+\frac{4r\sin\phistar}{d_1}+\frac{4k\tau_1}{d_0}\frac{r\sin\phistar}{d_1}+\frac{4k\tau_1+(4-8k)d_0}{d_0}\frac{g_0}{d_1}\Big].
    \end{aligned}
    \end{equation}
Note that\COMMENT{\red Add lots of mathclap with overbracket or underbracket. With strut or bigstrut if needed.\wentao{checking.}} 
    \begin{equation}\label{eq:g0d1ratioestimate}
        -\frac{51r}{4R}\overbracket{<}^{\bigstrut\mathclap{d_1>(r\sin\phistar)/3}}-\frac{17r}{4R}\frac{r\sin\phistar}{d_1}\overbracket{<}^{\bigstrut\mathclap{\text{\cref{lemma:h1g1g0g1g2g3lengthfractionestimate}\eqref{item02h1g1g0g1g2g3lengthfractionestimate}}}}\frac{g_0}{d_1}\underbracket{<}_{\bigstrut\mathclap{\text{\cref{lemma:h1g1g0g1g2g3lengthfractionestimate}\eqref{item02h1g1g0g1g2g3lengthfractionestimate}}}}\frac{17r}{4R}\frac{r\sin\phistar}{d_1}\overbracket{<}^{\bigstrut\mathclap{d_1>(r\sin\phistar)/3}}\frac{51r}{4R},
    \end{equation} 
    and \begin{equation}\label{eq:g2minusg0overd0estimate}
        -\frac{17r}{3R}\cdot1.0026\overbracket{<}^{\strut\mathllap{\text{\cref{proposition:casecparametersestimate}\ref{itemcasecparametersestimate5}}}}-\frac{2}{3d_0}\frac{17r}{2R}r\sin\phistar\overbracket{<}^{\strut\mathllap{\text{\cref{lemma:h1g1g0g1g2g3lengthfractionestimate}\eqref{item06h1g1g0g1g2g3lengthfractionestimate}}}}\frac{2(g_2-g_0)}{3d_0}\overbracket{<}^{\strut\mathclap{\text{\cref{lemma:h1g1g0g1g2g3lengthfractionestimate}\eqref{item06h1g1g0g1g2g3lengthfractionestimate}}}}\frac{2}{3d_0}\frac{17r}{2R}r\sin\phistar\overbracket{<}^{\strut\mathclap{\text{\cref{proposition:casecparametersestimate}\ref{itemcasecparametersestimate5}}}}\frac{17r}{3R}\cdot1.0026.
    \end{equation}
Therefore,\COMMENT{Below are examples of how to get rid of unneeded ``alignments''\Large\red---apply the same approach elsewhere} $\displaystyle -\frac{20}{3}\overbracket{\le}^{\strut\mathclap{k\le4/3}}4-8k\overset{\tau_1>0}<\frac{4k\tau_1+(4-8k)d_0}{d_0}\overbracket{<}^{\strut\mathclap{\substack{\text{\cref{remark:deffpclf}}:\\\tau_1<2d_0}}}2\times4k+(4-8k)=4$. This and \eqref{eq:g0d1ratioestimate} imply
\begin{equation}\label{eq:4thtermg0d1ration}
    \begin{aligned}
    -\frac{85r}{R}=-\frac{20\times51r}{3\times4R}<\frac{4k\tau_1+(4-8k)d_0}{d_0}\frac{g_0}{d_1}<\frac{20\times51r}{3\times4R}=\frac{85r}{R}.
\end{aligned}
\end{equation}
This and $\displaystyle11k-\frac{8kr\sin\phistar}{d_1}+\frac{4r\sin\phistar}{d_1}+\frac{4k\tau_1}{d_0}\frac{r\sin\phistar}{d_1}=\underbracket{\overbracket{11k+\big(-8k+4+\frac{4k\tau_1}{d_0}\big)\frac{r\sin\phistar}{d_1}}^{\overbracket{\scriptstyle<}^{\mathclap{\text{\cref{remark:deffpclf}}:\ \tau_1<2d_0}}11k+\frac{4r\sin\phistar}{d_1}\overbracket{\scriptstyle<}^{\mathclap{d_1>\frac{r\sin\phistar}{3}}}11k+12}}_{\mathclap{\underbracket{\scriptstyle>}_{\mathclap{\tau_1>0}}11k+(4-8k)\frac{r\sin\phistar}{d_1}
    \underbracket{\scriptstyle>}_{\mathclap{4-8k<0\text{ and }d_1<r\sin\phistar}}11k+4-8k=3k+4}}
    $ give
\begin{equation}\label{eq:4thtermanalysis}
\begin{aligned}
    0\underbracket{<}_{\mathclap{\qquad\eqref{eqJZHypCond}:\ R>1700r,\ k\ge1}}3k+4-\frac{85r}{R}<\Big[11k-\frac{8kr\sin\phistar}{d_1}+\frac{4r\sin\phistar}{d_1}+&\frac{4k\tau_1}{d_0}\frac{r\sin\phistar}{d_1}+\frac{4k\tau_1+(4-8k)d_0}{d_0}\frac{g_0}{d_1}\Big]<11k+12+\frac{85r}{R}.\\
\end{aligned}
\end{equation}
Then, assuming \eqref{eqJZHypCond} (i.e., $R>1700r$), by \eqref{eq:g2minusg0overd0estimate}, \eqref{eq:4thtermg0d1ration}, we find\COMMENT{\Large Use multline below\wentao{Done.}}
\begin{equation}\label{eq:4thterminC1C2}
    \begin{multlined}
        -\frac{17\times1.0026r}{3R}\big(11k+12+\frac{85r}{R}\big)\\
        \overbracket{<}^{\mathllap{\eqref{eq:g2minusg0overd0estimate}, \eqref{eq:4thtermanalysis}}}\frac{2(g_2-g_0)}{3d_0}\Big[11k-\frac{8kr\sin\phistar}{d_1}+\frac{4r\sin\phistar}{d_1}+\frac{4k\tau_1}{d_0}\frac{r\sin\phistar}{d_1}+\frac{4k\tau_1+(4-8k)d_0}{d_0}\frac{g_0}{d_1}\Big]\\
        \overbracket{<}^{\mathllap{\eqref{eq:g2minusg0overd0estimate}, \eqref{eq:4thtermanalysis}}}\frac{17\times1.0026r}{3R}\big(11k+12+\frac{85r}{R}\big).
    \end{multlined}
\end{equation}Also,\COMMENT{Hard-to-read jumbles follow. Maybe OK...\Hrule More mathclap\Hrule\large Use multline for \eqref{eq:4thterminC1C2}, \eqref{eq:3rdterminC1C2}, \eqref{eq:2ndterminC1C2}, 
\eqref{eq:1stterminC1C2}\wentao{maybe the big curl bracket cannot be used in multline.}}
\begin{equation}\label{eq:3rdterminC1C2}
\begin{aligned}
&\left.
\begin{aligned}
&-4\overbracket{<}^{\mathclap{\text{\cref{remark:deffpclf}: }\tau_1<2d_0}}\big(8k-4-4k\frac{\tau_1}{d_0}\big)\overbracket{<}^{\mathclap{\tau_1>0}}8k-4\overbracket\le^{\strut\mathclap{k\le4/3}}\frac{20}{3}\\
&-\frac{2}{3}\overbracket{<}^{\mathclap{\strut\tau_1>0}}\frac{\tau_1}{d_0}-\frac{2}{3}\overbracket{<}^{\mathclap{\qquad\strut\text{\cref{remark:deffpclf}: }\tau_1<2d_0}}\frac{4}{3}
\end{aligned}
\right\}\Rightarrow-\frac{16}{3}<\big(\frac{\tau_1}{d_0}-\frac{2}{3}\big)\big(8k-4-4k\frac{\tau_1}{d_0}\big)<\frac{80}{9}
\overbracket{\implies}^{\mathclap{\text{\eqref{eq:g0d1ratioestimate}}}}\\
&-\frac{1020r}{9R}=-\frac{80\times51r}{9\times4R}<\big(\frac{\tau_1}{d_0}-\frac{2}{3}\big)\frac{g_0}{d_1}\big(8k-4-4k\frac{\tau_1}{d_0}\big)<\frac{80\times51r}{9\times4R}=\frac{1020r}{9R},
\end{aligned}
\end{equation}
and
\begin{equation}\label{eq:2ndterminC1C2}
    \begin{aligned}
    &\left.
    \begin{aligned}
    -3k-\frac{2d_2}{3d_0}&\overbracket{<}^{\strut\mathclap{\tau_1>0}}\big(-3k+\frac{\tau_1}{d_0}-\frac{2d_2}{3d_0}\big)\overbracket{<}^{\mathclap{\strut\text{\cref{remark:deffpclf}: }\tau_1<2d_0}}\big(-3k+2-\frac{2d_2}{3d_0}\big)\\
    \text{\cref{proposition:casecparametersestimate}\ref{itemcasecparametersestimate5}: }&0.9975<\frac{d_2}{r\sin\phistar}<1.0026,\quad\frac{1}{1.0026}<\frac{r\sin\phistar}{d_0}<\frac{1}{0.9975}
    \end{aligned}
    \right\}\\
    \Longrightarrow& -3k-\frac{2\times1.0026}{3\times0.9975}<-3k+\frac{\tau_1}{d_0}-\frac{2d_2}{3d_0}<-3k+2-\frac{2\times0.9975}{3\times1.0026}\overset{k\ge1}<0\\
    \overbracket{\Longrightarrow}^{\mathclap{\text{\cref{lemma:h1g1g0g1g2g3lengthfractionestimate}\eqref{item01h1g1g0g1g2g3lengthfractionestimate}}}}& \big(-3k-\frac{2\times1.0026}{3\times0.9975}\big)\big(16.65r/R\big)<h_1\big(-3k+\frac{\tau_1}{d_0}-\frac{2d_2}{3d_0}\big)<\big(-3k+2-\frac{2\times0.9975}{3\times1.0026}\big)\big(-8.53r/R\big),
\end{aligned}
\end{equation}
\begin{equation}\label{eq:1stterminC1C2}
    -8k\frac{17r}{4R}\frac{1.0026}{0.9975}\overbracket{<}^{\mathclap{\strut\text{\cref{proposition:casecparametersestimate}\ref{itemcasecparametersestimate5}}}}-8k\frac{17r}{4R}\frac{d_1}{0.9975r\sin\phistar}\overbracket{<}^{\mathclap{\strut\text{\cref{lemma:h1g1g0g1g2g3lengthfractionestimate}\eqref{item03h1g1g0g1g2g3lengthfractionestimate}\qquad}}}8kd_1g_1\overbracket<^{\strut\qquad\mathclap{\text{\cref{lemma:h1g1g0g1g2g3lengthfractionestimate}\eqref{item03h1g1g0g1g2g3lengthfractionestimate}}}}8k\frac{17r}{4R}\frac{d_1}{0.9975r\sin\phistar}\overbracket{<}^{\mathclap{\strut\text{\cref{proposition:casecparametersestimate}\ref{itemcasecparametersestimate5}}}}8k\frac{17r}{4R}\frac{1.0026}{0.9975}.
\end{equation}
In the end, from \eqref{eq:C1C2reformed}, by \eqref{eq:1stterminC1C2}, \eqref{eq:2ndterminC1C2}, \eqref{eq:3rdterminC1C2}, \eqref{eq:4thterminC1C2} we get
\begin{align*}
\mathcal{C}_1(k)+\mathcal{C}_2(k)&\overbracket{=}^{\mathclap{\eqref{eq:C1C2reformed}\quad}}\overbracket{8kd_1g_1}^{\eqref{eq:1stterminC1C2}}+\overbracket{h_1\big(-3k+\frac{\tau_1}{d_0}-\frac{2d_2}{3d_0}\big)}^{\eqref{eq:2ndterminC1C2}}+\overbracket{\big(\frac{\tau_1}{d_0}-\frac{2}{3}\big)\frac{g_0}{d_1}\big(8k-4-4k\frac{\tau_1}{d_0}\big)}^{\eqref{eq:3rdterminC1C2}}\\
&\quad+\underbracket{\frac{2(g_2-g_0)}{3d_0}\Big[11k-\frac{8kr\sin\phistar}{d_1}+\frac{4r\sin\phistar}{d_1}+\frac{4k\tau_1}{d_0}\frac{r\sin\phistar}{d_1}+\frac{4k\tau_1+(4-8k)d_0}{d_0}\frac{g_0}{d_1}\Big]}_{\eqref{eq:4thterminC1C2}}\\
&>-8k\frac{17r}{4R}\frac{1.0026}{0.9975}+\big(-3k-\frac{2\times1.0026}{3\times0.9975}\big)\big(\frac{16.65r}{R}\big)-\frac{1020r}{9R}-\frac{17\times1.0026r}{3R}\big(11k+12+\frac{85r}{R}\big)\\
&>\frac{r}{R}\big(-146.7k-192.7-483\frac{r}{R}\big)\overbracket{>}^{k\le4/3}-389\frac{r}{R}-483\big(\frac{r}{R}\big)^2\qedhere
\end{align*}
\end{proof}

\begin{figure}[h]
\begin{center}    
    \begin{tikzpicture}[xscale=1.0,yscale=1.0]
        \pgfmathsetmacro{\PHISTAR}{6};
        \pgfmathsetmacro{\WIDTH}{8};
        \tkzDefPoint(-\PHISTAR,0){LL};
        \tkzDefPoint(\PHISTAR,0){RL};
        \tkzDefPoint(-\PHISTAR,\PHISTAR){X};
        \tkzDefPoint(\PHISTAR,\PHISTAR){IY};
        \tkzDefPoint(-\PHISTAR,0.5*\PHISTAR){S1};
        \tkzDefPoint(\PHISTAR,0.5*\PHISTAR){RS1};
        \tkzDefPoint(-\PHISTAR,\PHISTAR/3){L4};
        \path[name path=lineLLIY] (LL)--(IY);
        \path[name path=lineS1RL] (S1)--(RL);
        \path[name path=lineRS1LL] (RS1)--(LL);
        \path[name path=lineXRL] (X)--(RL);
        \path[name intersections= 
        {of=lineLLIY and lineS1RL,by={LS}}];
        \path[name intersections=
        {of=lineLLIY and lineXRL,by={HS}}];
        \path[name intersections=
        {of=lineRS1LL and lineXRL,by={RLS}}];
        \draw[name path=IYRS1, dashed, thin] (IY)--(RS1);
        \draw[name path=RS1RLS, dashed, thin] (RS1)--(RLS);
        \draw[name path=RLSHS, dashed, thin] (RLS)--(HS);
        \draw[name path=HSIY, dashed, thin] (HS)--(IY);

        \tkzLabelPoint[left](L4){\small $\theta_1=\frac{1}{3}\Phistar$};
        \draw[name path=lineLSHS, very thick](LS)--(HS);
        \draw[name path=lineHSX, very thick](HS)--(X);
        \draw[name path=lineXS1, very thick](X)--(S1);
        \draw[name path=lineS1LS, very thick](S1)--(LS);
        \draw[ultra thin,dashed](L4)--(LS);
        \fill [orange, opacity=2/30](LS) -- (HS) -- (X) -- (S1) -- cycle;
        \tkzDefPoint(-\PHISTAR,0.47*\PHISTAR){L5};
        \path[name path=HorSeg] (L5)-- +(25,0);
        \path[name intersections={of=lineS1RL and HorSeg,by={Seg1}}];
        \path[name intersections={of=lineLSHS and HorSeg,by={Seg2}}];
        \draw[ultra thin,red,dashed](Seg1)--(Seg2);
        \draw[ultra thin,red,dashed](L5)--(Seg1);
        \tkzDefPoint(-\PHISTAR,0.4*\PHISTAR){L6};
        \path[name path=HorSeg2] (L6)-- +(15,0);
        \path[name intersections={of=lineS1LS and HorSeg2,by={Seg3}}];
        \path[name intersections={of=lineLSHS and HorSeg2,by={Seg4}}];
        \draw[ultra thin,blue](Seg3)--(Seg4);
        \path[name intersections={of=RLSHS and HorSeg2,by={Seg7}}];
        \path[name intersections={of=RS1RLS and HorSeg2,by={Seg8}}];
        \draw[blue, ultra thin] (Seg7)--(Seg8);
        \tkzLabelSegment[above](Seg7,Seg8){\small $\Fcal(\mathfrak{l}_L)$}; 
        \tkzLabelPoint[left, xshift=-0.7](Seg7){\small $\Fcal(\mathcal{W}_L)$};
        \tkzLabelPoint[right, xshift=2.7](Seg8){\small $\Fcal(\mathcal{E}_L)$};
        \tkzLabelPoint[above](L5){\small $\textcolor{red}{d_1=\frac{1}{2}r\sin{\phistar}}$};
        \tkzLabelPoint[above right](X){\small $d_1=r\sin{\phistar}$, $\theta_1=\Phistar$};
        \tkzDefPoint(-\PHISTAR,0.6*\PHISTAR){Seg5};
        \path[name path=HorSeg3] (Seg5)-- +(15,0);
        \path[name intersections={of=lineHSX and HorSeg3,by={Seg6}}];
        \draw[ultra thin,blue](Seg5)--(Seg6);
        \path[name intersections={of=HSIY and HorSeg3,by={Seg9}}];
        \path[name intersections={of=IYRS1 and HorSeg3,by={Seg10}}];
        \draw[blue, ultra thin] (Seg9)--(Seg10);
        \tkzLabelSegment[above](Seg9,Seg10){\small $\Fcal(\mathfrak{l}_U)$}; 
        \tkzLabelPoint[left, xshift=-0.7](Seg9){\small $\Fcal(\mathcal{W}_U)$};
        \tkzLabelPoint[right, xshift=2.7](Seg10){\small $\Fcal(\mathcal{E}_U)$};
        \tkzDrawPoints(X,LS,HS,S1,L4,Seg1,Seg2,L5,Seg3,Seg4,Seg5,Seg6,Seg7,Seg8,Seg9,Seg10,RL,LL,IY);
        \tkzLabelPoint[right, xshift=2.7](Seg4){\small $\mathcal{E}_L$};
        \tkzLabelPoint[right, xshift=2.7](Seg6){\small $\mathcal{E}_U$};
        \tkzLabelPoint[left, xshift=-0.7](Seg5){\small $\mathcal{W}_U$};
        \tkzLabelSegment[above](Seg5,Seg6){\small $\mathfrak{l}_U$}; 
        \tkzLabelPoint[left, xshift=-0.7](Seg3){\small $\mathcal{W}_L$};
        \tkzLabelSegment[above](Seg3,Seg4){\small{$\mathfrak{l}_L$}};
        \node[] at ($(X)+(-58:1)$)  { \textcolor{blue}{$M^{\text{\upshape{in}}}_{R,1}$}};
        \node[] at ($(IY)+(-122:1)$)  { \textcolor{red}{$M^{\text{\upshape{out}}}_{R,1}$}};
        \draw[dashed,ultra thin](LL)--(RL);
        \draw[dashed,ultra thin](RL)--(RS1);
        \draw[dashed,ultra thin](LL)--(S1);
        \draw[dashed,ultra thin](RLS)--(LL);
        \draw[dashed,ultra thin](RLS)--(RL);
        \draw[dashed,ultra thin](LS)--(LL);
        \draw[dashed,ultra thin](LS)--(RL);
        \tkzLabelPoint[below right](LL){\small $\Phi=-\Phistar, SW=(-\Phistar,0)$};
        \tkzLabelPoint[below](RL){\small $\Phi=+\Phistar$};
        \tkzLabelPoint[above](IY){\small $NE=(\Phistar,\Phistar)$};
        
    \end{tikzpicture}
    \caption{The \textcolor{blue}{$M^{\text{\upshape in}}_{R,1}$} component  from the $n_1=1$-region of \cref{fig:x1_n1cases_caseC}.}\label{fig:x1_n1equal1case}
    \end{center}
\end{figure}

\begin{proposition}\label{proposition:QUpperLowerRegionsLowerBoundsChatCtilde}
    With the standing assumptions of this section (\cref{REMCaseCAssumptions}), let $x_1=(\Phi_1,\theta_1)\in \textcolor{blue}{M^{\text{\upshape{in}}}_{R,1}}$  (\cref{fig:x1_n1equal1case}), i.e., $n_1=1$ (see \cref{fig:x1_n1cases_caseC}). Then this indeed falls into case (c) of \eqref{eqMainCases}, which is the hypothesis of \cref{TEcase_c}. Suppose $x_1$ is on the line $\mathfrak{l}=\Big\{(\Phi,\theta_1)\bigm|d_1=R\sin\theta_1=\const\Big\}\subset\textcolor{blue}{M^{\text{\upshape{in}}}_{R,1}}$. With $\mathfrak{Q}(t)$ from \eqref{eq:II1II2Quadratic},\COMMENT{Remember that we are considering a change in notation\wentao{Change is done now.}} and fixed $k\in[1,4/3]$, $\nu=\frac{d_1}{r\sin\phistar}$, we have the following.
    \begin{enumerate}
        \item\label{item01Qlowerboundn1eq1} $\nu\in(1/3,1)$.
        \item\label{item02Qlowerboundn1eq1}$\tau_0$, $\tau_1$, $d_0$, $d_1$ $d_2$ from \cref{def:fpclf} are functions of $x=(\Phi,\theta_1)\in\mathfrak{l}$. For constant $\theta_1$, $d_0$ is an increasing function of $\Phi$, and $\tau_1$ is a decreasing function of $\Phi$.
        \item\label{item03Qlowerboundn1eq1} If  $\nu>0.5$, that is, $x_1\in\mathfrak{l}=\mathfrak{l}_U=
        \Big\{x=(\Phi,\theta_1)\bigm|d_1=\nu r\sin\phistar>0.5r\sin\phistar\Big\}$ with endpoints $\mathcal{W}_U=(\Phi_{\mathcal{W}_U},\theta_1)$ (think "west"), $\mathcal{E}_U=(\Phi_{\mathcal{E}_U},\theta_1)$ (think "east") in \cref{fig:x1_n1equal1case} and $\Phi_{\mathcal{W}_U}<\Phi_{\mathcal{E}_U}$, then \[\mathfrak{Q}\big(\frac{\tau_1}{d_0}\bigm|_{x=x_1}\big)>\min{\Big\{\frac{11}{9},\frac{1}{3}+\mathcal{C}^{\mathcal{W}}_{U}(\nu,k)\Big\}},\] where, with $h_1,g_1$ from \cref{def:lengthfractionfunctionsn1eq1}\eqref{item01lengthfractionfunctionsn1eq1},\eqref{item03lengthfractionfunctionsn1eq1}, $\hat{h}_1(U)\dfn\lim_{\mathfrak{l}_U\ni x\to \mathcal{W}_U}h_1(x)$, $\hat{g}_1(U)\dfn\lim_{\mathfrak{l}_U\ni x\to \mathcal{W}_U}g_1(x)$,
        \begin{equation}\label{eq:CUWnuk}
        \begin{aligned}
        \mathcal{C}^{\mathcal{W}}_U(\nu,k)&=\Big[\frac{-16k-12}{3\nu}+8k\Big]\big(\hat{h}_1(U)-2d_1\hat{g}_1(U)\big)-\frac{4k}{\nu}\big(\hat{h}_1(U)-2d_1\hat{g}_1(U)\big)^2.
        \end{aligned}    
        \end{equation}
        
        \item\label{item04Qlowerboundn1eq1} If  $\nu\le 0.5$ that is $x_1\in\mathfrak{l}=\mathfrak{l}_L=
        \Big\{x=(\Phi,\theta_1)\bigm|d_1=\nu r\sin\phistar\le0.5r\sin\phistar\Big\}$ with endpoints $\mathcal{W}_L=(\Phi_{\mathcal{W}_L},\theta_1)$, $\mathcal{E}_L=(\Phi_{\mathcal{E}_L},\theta_1)$ in \cref{fig:x1_n1equal1case} and $\Phi_{\mathcal{W}_L}<\Phi_{\mathcal{E}_L}$, then \[\mathfrak{Q}\big(\frac{\tau_1}{d_0}\Bigm|_{x=x_1}\big)>\min{\Big\{\frac{11}{9}+\mathcal{C}^{\mathcal{E}}_{L}(\nu,k),\frac{1}{3}+\mathcal{C}^{\mathcal{W}}_{L}(\nu,k)\Big\}},
        \] 
        where
        \begin{equation}\label{eqL:CLWnukandCLEnukand}
        \begin{aligned}
            \mathcal{C}^{\mathcal{W}}_{L}(\nu,k)&=2d_1\hat{g}_1(L)\big(-24k+\frac{32k-12}{3\nu}\big)-\frac{16kd_1^2}{\nu}\hat{g}^2_1(L),
            \\
            \mathcal{C}^{\mathcal{E}}_{L}(\nu,k)&=\big(-8k+\frac{32k-12}{3\nu}\big)\big(\tilde{h}_1(L)-4d_1\tilde{g}_1(L)\big)-\frac{8k}{\nu}\big(\tilde{h}_1(L)-4d_1\tilde{g}_1(L)\big)\big(2-4\nu\big)-\frac{4k}{\nu}\big(\tilde{h}_1(L)-4d_1\tilde{g}_1(L)\big)^2,
        \end{aligned}   
        \end{equation}
        with (see \cref{def:lengthfractionfunctionsn1eq1}\eqref{item01lengthfractionfunctionsn1eq1},\eqref{item03lengthfractionfunctionsn1eq1})
        \begin{align*}
        \hat{h}_1(L)\dfn&\lim_{\mathfrak{l}_L\ni x\to \mathcal{W}_L}h_1(x),&
            \hat{g}_1(L)\dfn&\lim_{\mathfrak{l}_L\ni x\to \mathcal{W}_L}g_1(x),\\
            \tilde{h}_1(L)\dfn&\lim_{\mathfrak{l}_L\ni x\to \mathcal{E}_L}h_1(x),&
            \tilde{g}_1(L)\dfn&\lim_{\mathfrak{l}_L\ni x\to \mathcal{E}_L}g_1(x).
        \end{align*}
    \item\label{item05Qlowerboundn1eq1}  $\mathcal{C}_U^\mathcal{W}(\nu,k)$, $\mathcal{C}_L^\mathcal{E}(\nu,k)$, and $\mathcal{C}_L^\mathcal{W}(\nu,k)$  from \eqref{eq:CUWnuk},  \eqref{eqL:CLWnukandCLEnukand} satisfy \begin{align*}
        \mathcal{C}_U^\mathcal{W}(\nu,k)&>-192\big(\frac{r}{R}\big)-4667\big(\frac{r}{R}\big)^2,\\
        \mathcal{C}_L^\mathcal{W}(\nu,k)&>-50\big(\frac{r}{R}\big)-194\big(\frac{r}{R}\big)^2,\\
        \mathcal{C}_L^\mathcal{E}(\nu,k)&>-92\big(\frac{r}{R}\big)-7984\big(\frac{r}{R}\big)^2.
    \end{align*}
    \end{enumerate} 
\end{proposition}
\begin{proof}
These computations will run to page \pageref{proposition:casecexpansionn1eq1}.
        \hidesubsection{Proof of (1)}\cref{proposition:casecparametersestimate}\eqref{item03lengthfractionfunctionsn1eq1} with $n_1=1$ implies $d_1>\frac{1}{3}r\sin\phistar$. Also, since $p(x_1)=P$, $p(\Fc(x_1))=P_1$ are two consecutive collisions on $\Gamma_R$ with distance $2d_1$, $2d_1=|p(x_1)p(\Fc(x_1))|<|AB|=2r\sin\phistar$ (as shown in \cref{fig:n1_equal1_upper_leftendpoint_case,fig:n1_equal1_lower_leftendpoint_case,fig:n1_equal1_upper_leftendpoint_case2,fig:n1_equal1_lower_rightendpoint_case,fig:n1_equal1_upper_rightendpoint_case,fig:n1_equal1_upper_rightendpoint_case2}). Hence $\frac{1}{3}r\sin\phistar<d_1<r\sin\phistar$, $\nu=\frac{d_1}{r\sin\phistar}\in\big(\frac{1}{3},1\big)$.
        \hidesubsection{Proof of (2)}By \cref{def:fpclf}, $\tau_0=|p(x)p(\Fcal^{-1}(x))|$, $d_0$ are determined by the collision angle $\theta$ coordinate of $\Fcal^{-1}(x)\in\Mrout$, $\tau_1=|p(\Fcal(x))p(\Fcal^{2}(x))|$, $d_2$ are determined by the collision angle $\theta$ coordinate of $\Fcal^{2}(x)\in \Mrin$. Hence length functions $\tau_0$, $\tau_1$, $d_0$, $d_2$ defined in \cref{def:fpclf} are functions of $x\in\mathfrak{l}$.
        
        Since in \cref{fig:x1_n1equal1case} either $\mathfrak{l}=\mathfrak{l}_U\subset\MRin$ or $\mathfrak{l}=\mathfrak{l}_L\subset\MRin$ with some constant $\theta_1<\Phistar=\sin^{-1}(r\sin\phistar/R)\overset{R>r}<\sin^{-1}{(\sqrt{4r/R})}$, \cref{monotone_p4} implies that $d_0$ is a monotone increasing function of $\Phi$.

        Let $x\in\mathfrak{l}$, $\Fcal(x)\in\Fcal(\mathfrak{l})\subset\MRout$ hence $\Fcal^2(x)\in\Mrin$ then we can suppose $\Gamma_r\ni p(\Fcal^2(x))=(x_Q,y_Q)$ in the standard coordinate system defined by \cref{def:StandardCoordinateTable} for the billiard table.

        We note that $\Fcal\bigm|_\mathfrak{l}(x)=\Fcal(\Phi,\theta_1=\const)=(\Phi+2\theta_1,\theta_1=\const)$ is the diffeomorphism 
        between $\mathfrak{l}_U$ and $\Fcal(\mathfrak{l}_U)\subset\MRout$, between $\mathfrak{l}_L$ and $\Fcal(\mathfrak{l}_L)\subset\MRout$ in \cref{fig:x1_n1equal1case}. By \cref{monotone_p3}\eqref{itemmonotone_p3-1}, \eqref{eq1monotone_p3} and chain rule we get 
        \[
        \frac{d\tau_1}{d\Phi}\Bigm|_{x=(\Phi,\theta_1=\const)}\overbracket{=}^{\strut\mathclap{\text{chain rule}}}\frac{d\tau_1}{d\Phi\bigm|_{\Fcal(\mathfrak{l})}}\cdot\overbracket{\frac{d\Phi\bigm|_{\Fcal(\mathfrak{l})}}{d\Phi\bigm|_{\mathfrak{l}}}\Bigm|_{x}}^{\mathclap{\text{since }\Phi|_{\Fcal(\mathfrak{l})}=(\Phi+2\theta_1)|_{\mathfrak{l}},\ \frac{d(\Phi+2\const)}{d\Phi}=1}}=\frac{d\tau_1}{d\Phi\bigm|_{\Fcal(\mathfrak{l})}}\overbracket{=}^{\bigstrut\mathclap{\text{\cref{monotone_p3}\eqref{itemmonotone_p3-1}, \eqref{eq1monotone_p3}}}}\frac{-bx_Q}{d_2}
        \]
        We also recall that $I y_*=(\phistar,\phistar)$ in \cref{fig:MrN}, $\Fcal^2(x)\nfd(\phi_Q,\theta_Q)\in\Nin$ by \cref{contraction_region}. And more precisely in case (c) of \eqref{eqMainCases}, by \cref{proposition:NinsizeInA1BCcases}\eqref{item01NinsizeInA1BCcases} and \eqref{eq:observation2(a1)(b)(c)region} we get \COMMENT{?\wentao{edited using new proposition now.}}the following. 
        \[0<\phi_Q-\phistar\overbracket{<}^{\mathllap{\text{\cref{proposition:NinsizeInA1BCcases}\eqref{item01NinsizeInA1BCcases} and \eqref{eq:observation2(a1)(b)(c)region}}}}2.001\sqrt{4r/R}\overbracket{<}^{\mathrlap{\text{\eqref{eqJZHypCond}: }R>1700r}}0.0971.\]

        Therefore, we get $\phistar<\phi_Q<\phistar+0.0971<\tan^{-1}(1/3)+0.0971<\pi/2$. Hence in coordinate system \cref{def:StandardCoordinateTable}, $x_Q=r\sin{\phi_Q}>0$, $\frac{d\tau_1}{d\Phi}\Bigm|_{x=(\Phi,\theta_1=\const)}=-\frac{bx_Q}{d_2}<0$, that is, $\tau_1$ as a function of $x=(\Phi,\theta_1)\in\mathfrak{l}$ with constant $\theta_1$ is thus a decreasing function of $\Phi$.

        
\begin{figure}[h]
\begin{center}
\begin{tikzpicture}[xscale=1.4,yscale=1.4]
    \tkzDefPoint(0,0){Or};
    \tkzDefPoint(0,-4.8){Y};
    \clip(-5.7,-4.2) rectangle (5.7,-2.8);
    \pgfmathsetmacro{\rradius}{4.5};
    \pgfmathsetmacro{\phistardeg}{40};
    \pgfmathsetmacro{\XValueArc}{\rradius*sin(\phistardeg)};
    \pgfmathsetmacro{\YValueArc}{\rradius*cos(\phistardeg)};
    \pgfmathsetmacro{\Rradius}{
        \rradius*sin(\phistardeg)/sin(18)
    };
    \pgfmathsetmacro{\bdist}{(-\rradius*cos(\phistardeg))+cos(18)*\Rradius}
    \tkzDefPoint(0,\bdist){OR};
    \tkzDefPoint(\XValueArc,-1.0*\YValueArc){B};    \tkzDefPoint(-1.0*\XValueArc,-1.0*\YValueArc){A};
    \tkzDrawArc[name path=Cr,very thick](Or,B)(A);
    \tkzDrawArc[name path=CR,very thick](OR,A)(B);
    \pgfmathsetmacro{\PXValue}{\Rradius*sin(14)};
    \pgfmathsetmacro{\PYValue}{\bdist-\Rradius*cos(14)};
    \tkzDefPoint(\PXValue,\PYValue){P2};
    \tkzDefPointBy[projection=onto OR--P2](A)\tkzGetPoint{M};
    \coordinate (Q_far) at ($(A)!2!(M)$);
    \path[name path=Line_P2_Q_far](P2)-- (Q_far);
    \path[name intersections={of=Cr and Line_P2_Q_far,by={Q}}];
    \tkzDefPoint(0,\bdist-\Rradius){C};
    \begin{scope}[very thin,decoration={markings,mark=at position 0.5 with {\arrow{>}}}] 
    \draw[blue,postaction={decorate}] (P2)--(Q);
    \draw[blue,postaction={decorate}] (A)--(P2);
    \end{scope}
    \tkzLabelPoint[above right](A){\small $A=P=p(x=\mathcal{W}_{U})$};
    \tkzLabelPoint[below left](A){\small $\hat{\tau}_0(U)=|AP|=0$};
    \tkzLabelSegment[above right=0.2pt and 0.6pt](A,P2){$2\hat{d}_1(U)=2d_1$};
    \tkzLabelPoint[below, yshift=-3.5](P2){\small $P_1=p(\Fc(x=\mathcal{W}_{U}))$};
    \tkzLabelPoint[right](Q){\small $Q=p(\Fc^2(x=\mathcal{W}_{U}))$};
    \tkzLabelPoint[right](B){$B$};
    \tkzLabelPoint[below](Or){$O_r$};
    \tkzDrawPoints(A,B,Or,Q,P2);
    \tkzLabelSegment[above left=0.5pt and 3pt](P2,Q){$\hat{\tau}_1(U)=|P_1Q|$};
    
\end{tikzpicture}
\caption{With constant $\theta_1\ge0.5\Phistar$, $(\Phi_{\mathcal{W}_U},\theta_1)=x=\mathcal{W}_U$ is on the trajectory from corner $A$  }\label{fig:n1_equal1_upper_leftendpoint_case}

\end{center}
\end{figure}

\begin{figure}[h] 
\begin{center}
\begin{tikzpicture}[xscale=1.4,yscale=1.4]
    \tkzDefPoint(0,0){Or};
    \tkzDefPoint(0,-4.8){Y};
    \clip(-5.5,-4.3) rectangle (5.5,-2.3);
    \pgfmathsetmacro{\rradius}{4.5};
    \pgfmathsetmacro{\phistardeg}{40};
    \pgfmathsetmacro{\XValueArc}{\rradius*sin(\phistardeg)};
    \pgfmathsetmacro{\YValueArc}{\rradius*cos(\phistardeg)};
    \pgfmathsetmacro{\Rradius}{
        \rradius*sin(\phistardeg)/sin(18)
    };
    \pgfmathsetmacro{\bdist}{(-\rradius*cos(\phistardeg))+cos(18)*\Rradius}
    \tkzDefPoint(0,\bdist){OR};
    \tkzDefPoint(\XValueArc,-1.0*\YValueArc){B};    \tkzDefPoint(-1.0*\XValueArc,-1.0*\YValueArc){A};
    \tkzDrawArc[name path=Cr,very thick](Or,B)(A);
    \tkzDrawArc[name path=CR,very thick](OR,A)(B);
    \pgfmathsetmacro{\PXValue}{-\Rradius*sin(14)};
    \pgfmathsetmacro{\PYValue}{\bdist-\Rradius*cos(14)};
    \tkzDefPoint(\PXValue,\PYValue){P1};
    
    \tkzDefPointBy[projection=onto OR--P1](B)\tkzGetPoint{M};
    \coordinate (T_far) at ($(B)!2!(M)$);
    \path[name path=Line_P1_T_far](P1)-- (T_far);
    \path[name intersections={of=Cr and Line_P1_T_far,by={T}}];
    \tkzDefPoint(0,\bdist-\Rradius){C};
    \begin{scope}[ultra thin,decoration={markings,mark=at position 0.5 with {\arrow{>}}}] 
    \draw[blue,postaction={decorate}] (T)--(P1);
    \draw[blue,postaction={decorate}] (P1)--(B);
    \end{scope}
    \tkzLabelPoint[left](A){$A$};
    \tkzLabelPoint[below left](P1){\small $P=p(x=\mathcal{E}_U)$};
    \tkzLabelPoint[above right](T){\small $T=p(\Fc^{-1}(x=\mathcal{E}_U))$};
    \tkzLabelPoint[above left](B){$B=p(F(x=\mathcal{E}_U))$};
    \tkzLabelPoint[below](Or){$O_r$};
    \tkzLabelSegment[above left=0.5pt and 3pt](P1,B){$2d_1$};
    \tkzDrawPoints(A,B,Or,P1,T);
\end{tikzpicture}
\caption{With constant $\theta_1\ge0.5\Phistar$, $(\Phi_{\mathcal{E}_U},\theta_1)=x=\mathcal{E}_U$ is on a trajectory to corner $B$}\label{fig:n1_equal1_upper_rightendpoint_case}
\end{center}
\end{figure}
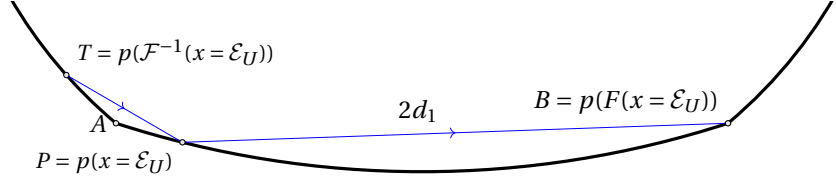

\begin{figure}[h]
\begin{center}
\begin{tikzpicture}[xscale=1.8,yscale=1.8]
    \tkzDefPoint(0,0){Or};
    \tkzDefPoint(0,-4.8){Y};
    \clip(-4.3,-4.2) rectangle (4.3,-2.9);
    \pgfmathsetmacro{\rradius}{4.5};
    \pgfmathsetmacro{\phistardeg}{40};
    \pgfmathsetmacro{\XValueArc}{\rradius*sin(\phistardeg)};
    \pgfmathsetmacro{\YValueArc}{\rradius*cos(\phistardeg)};
    \pgfmathsetmacro{\Rradius}{
        \rradius*sin(\phistardeg)/sin(18)
    };
    \pgfmathsetmacro{\bdist}{(-\rradius*cos(\phistardeg))+cos(18)*\Rradius}
    \tkzDefPoint(0,\bdist){OR};
    \tkzDefPoint(\XValueArc,-1.0*\YValueArc){B};    \tkzDefPoint(-1.0*\XValueArc,-1.0*\YValueArc){A};
    \tkzDrawArc[name path=Cr,very thick](Or,B)(A);
    \tkzDrawArc[name path=CR,very thick](OR,A)(B);
    \pgfmathsetmacro{\PXValue}{-\Rradius*sin(1)};
    \pgfmathsetmacro{\PYValue}{\bdist-\Rradius*cos(1)};
    \tkzDefPoint(\PXValue,\PYValue){P1};
    
    \tkzDefPointBy[projection=onto OR--P1](A)\tkzGetPoint{M};
    \coordinate (P2) at ($(A)!2!(M)$);
    
    \tkzDefPointBy[projection=onto OR--P2](P1)\tkzGetPoint{M2};
    \coordinate (Q_far) at ($(P1)!2!(M2)$);
    \path[name path=Line_P2_Q_far](P2)-- (Q_far);
    \path[name intersections={of=Cr and Line_P2_Q_far,by={Q}}];
    \tkzDefPoint(0,\bdist-\Rradius){C};
    \begin{scope}[ultra thin,decoration={markings,mark=at position 0.5 with {\arrow{>}}}] 
    \draw[blue,postaction={decorate}] (A)--(P1);
    \draw[blue,postaction={decorate}] (P1)--(P2);
    \draw[blue,postaction={decorate}] (P2)--(Q);
    \end{scope}
    \tkzLabelPoint[above right](A){\small $A=p\big(\Fc^{-1}(x=\mathcal{E}_U)\big)$};
    \tkzLabelPoint[below](P1){\small $P=p\big(x=\mathcal{E}_U\big)$};
    \tkzLabelPoint[below](P2){\small $p\big(\Fc(x=\mathcal{E}_U)\big)$};
    \tkzLabelPoint[above left](Q){\small $Q=p\big(\Fc^{2}(x=\mathcal{E}_U)\big)$};
    \tkzLabelPoint[right](B){$B$};
    \tkzLabelPoint[below](Or){$O_r$};
    \tkzLabelSegment[above left=0.5pt and 0.1pt](P1,P2){$2d_1$};
    \tkzDrawPoints(A,B,Or,P1,P2,Q);
    
\end{tikzpicture}
\caption{With constant $\theta_1<0.5\Phistar$ and $d_1>0.5r\sin\phistar$, $(\Phi_{\mathcal{E}_U},\theta_1)=x=\mathcal{E}_U$ is on a trajectory from corner $A$
} \label{fig:n1_equal1_upper_rightendpoint_case2}
\end{center}
\end{figure}
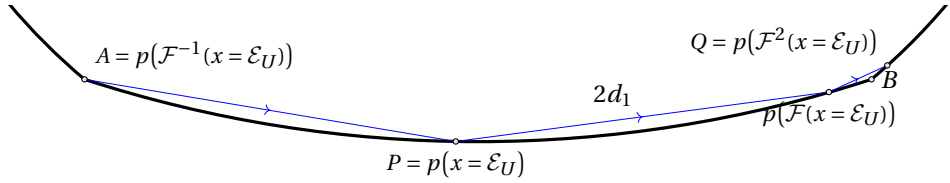

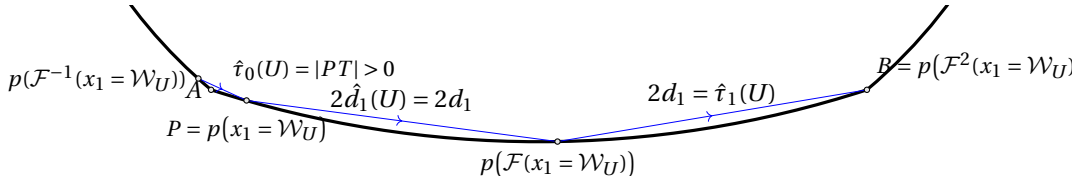
\begin{figure}[h]
\begin{center}
\begin{tikzpicture}[xscale=1.5,yscale=1.5]
    \tkzDefPoint(0,0){Or};
    \tkzDefPoint(0,-4.8){Y};
    \clip(-4.7,-4.2) rectangle (4.7,-2.7);
    \pgfmathsetmacro{\rradius}{4.5};
    \pgfmathsetmacro{\phistardeg}{40};
    \pgfmathsetmacro{\XValueArc}{\rradius*sin(\phistardeg)};
    \pgfmathsetmacro{\YValueArc}{\rradius*cos(\phistardeg)};
    \pgfmathsetmacro{\Rradius}{
        \rradius*sin(\phistardeg)/sin(18)
    };
    \pgfmathsetmacro{\bdist}{(-\rradius*cos(\phistardeg))+cos(18)*\Rradius}
    \tkzDefPoint(0,\bdist){OR};
    \tkzDefPoint(\XValueArc,-1.0*\YValueArc){B};    \tkzDefPoint(-1.0*\XValueArc,-1.0*\YValueArc){A};
    \tkzDrawArc[name path=Cr, very thick](Or,B)(A);
    \tkzDrawArc[name path=CR, very thick](OR,A)(B);
    \pgfmathsetmacro{\PXValue}{\Rradius*sin(1)};
    \pgfmathsetmacro{\PYValue}{\bdist-\Rradius*cos(1)};
    \tkzDefPoint(\PXValue,\PYValue){P2};
    
    \tkzDefPointBy[projection=onto OR--P2](B)\tkzGetPoint{M};
    \coordinate (P1) at ($(B)!2!(M)$);
    
    \tkzDefPointBy[projection=onto OR--P1](P2)\tkzGetPoint{M2};
    \coordinate (T_far) at ($(P2)!2!(M2)$);
    \path[name path=Line_P1_T_far](P1)-- (T_far);
    \path[name intersections={of=Cr and Line_P1_T_far,by={T}}];
    \tkzDefPoint(0,\bdist-\Rradius){C};
    \begin{scope}[ultra thin,decoration={markings,mark=at position 0.5 with {\arrow{>}}}] 
    \draw[blue,postaction={decorate}] (T)--(P1);
    \draw[blue,postaction={decorate}] (P1)--(P2);
    \draw[blue,postaction={decorate}] (P2)--(B);
    \end{scope}
    \tkzLabelPoint[left](A){$A$};
    \tkzLabelPoint[below =2.5pt](P1){\small $P=p\big(x_1=\mathcal{W}_U\big)$};
    \tkzLabelPoint[below](P2){\small $p\big(\Fc(x_1=\mathcal{W}_U)\big)$};
    \tkzLabelPoint[ left =0.8pt](T){\small $T=p(\Fc^{-1}(x_1=\mathcal{W}_U))$};
    \tkzLabelPoint[above right](B){\small $B=p\big(\Fc^2(x_1=\mathcal{W}_U)\big)$};
    \tkzLabelPoint[below](Or){$O_r$};
    \tkzLabelSegment[above =0.2pt](P1,P2){$2\hat{d}_1(U)=2d_1$};
    \tkzLabelSegment[above =0.5pt](P2,B){$2d_1=\hat{\tau}_1(U)$};
    \tkzLabelSegment[above right =0pt and 0.3pt](P1,T){\small $\hat{\tau}_0(U)=|PT|>0$};
    \tkzDrawPoints(A,B,P1,P2,T);
\end{tikzpicture}
\caption{
With constant $\theta_1<0.5\Phistar$ and $d_1>0.5r\sin\phistar$, $(\Phi_{\mathcal{W}_U},\theta_1)=x=\mathcal{W}_U$ is on a trajectory to corner $B$
}\label{fig:n1_equal1_upper_leftendpoint_case2}

\end{center}
\end{figure}

\hidesubsection{Proof of (3)} For $x_1\in \mathfrak{l}_U$ in \cref{fig:x1_n1equal1case} with $d_1>0.5r\sin\phistar$, we define\COMMENT{Why not ``$\hat{\tau}_i(U)\dfn\lim_{\mathfrak{l}_U\ni x\to \mathcal{W}_U}\tau_i\ge0$, $\hat{d}_1(U)\dfn\lim_{\mathfrak{l}_U\ni x\to \mathcal{W}_U}d_1$''? Alternatively, fit the pair of $\tau$ and the triple of $d$ side by side.}
\begin{equation}\label{eq:lengthlimitWU}
    \begin{aligned}
        \hat{\tau}_0(U)\dfn&\lim_{\mathfrak{l}_U\ni x\to \mathcal{W}_U}\tau_0\ge0
        \\
        \hat{\tau}_1(U)\dfn&\lim_{\mathfrak{l}_U\ni x\to \mathcal{W}_U}\tau_1\ge0
        \\
        \hat{d}_0(U)\dfn&\lim_{\mathfrak{l}_U\ni x\to \mathcal{W}_U}d_0>0
        \\
        \hat{d}_1(U)\dfn&\lim_{\mathfrak{l}_U\ni x\to \mathcal{W}_U}d_1=\const=\nu\cdot r\sin{\phistar}
        \\
        \hat{d}_2(U)\dfn&\lim_{\mathfrak{l}_U\ni x\to \mathcal{W}_U}d_2>0
    \end{aligned}
    \end{equation}
    
Note that these limits exist by the monotone properties proved in \cref{monotone_p3}\eqref{itemmonotone_p3-1}\eqref{itemmonotone_p3-3} for $\tau_0$, $\tau_1$ and in \cref{monotone_p4} for $d_0$ thus $d_2$ when considering functions of $x$ restricted to $\mathfrak{l}_U$ and $\Fc(\mathfrak{l}_U)$.\COMMENT{Because...\wentao{reason added.}} In \cref{fig:x1_n1equal1case}, $\mathfrak{l}_U$ and its image $\Fcal(\mathfrak{l}_U)$ are line segments with interiors containing no first-order singularity of $\Fcal$ (See \cref{def:BasicNotations0}). The line segment $\mathfrak{l}_U$ has endpoints $\mathcal{W}_U$, $\mathcal{E}_U$ possibly being on the boundary of $M_R$.  The line segment $\Fcal(\mathfrak{l}_U)$ has endpoints that are possible on the boundary $M_R$.

By \cref{contraction_region}$, \Fcal^{-1}(\mathfrak{l}_U)$ is a curve in $\Nout$ and $\Fcal^{2}(\mathfrak{l}_U)$ is a curve in $\Nin$ (see \cref{fig:MrN}) with endpoints possibly on the boundary of $M_r$. 
Then the collision angle and position limits for the endpoints can be defined by continuous extension from $\Fcal^{-1}(\mathfrak{l}_U)$, $\Fcal^{2}(\mathfrak{l}_U)$ on the interior of $M_r$ to the endpoints 
$\mathcal{W}_U$, $\mathcal{E}_U$; see \cref{fig:n1_equal1_upper_leftendpoint_case,fig:n1_equal1_upper_rightendpoint_case,fig:n1_equal1_upper_leftendpoint_case2,fig:n1_equal1_upper_rightendpoint_case2}.

Also, since by \eqref{eqJZHypCond} the $\delta_*$ in \cref{def:SectionSetsDefs} satisfies the following.\COMMENT{\Large\color{red}This sentence is broken, and so is the next line. What do you mean?\wentao{fixed. with $\delta_*$ updated }}
\[\delta_*=4.5\sqrt{r/R}\overbracket{<}^{\eqref{eqJZHypCond}}4.5\cdot\sqrt{\sin^2\phistar/30000}<4.5\sqrt{\phistar^2/30000}<0.025\phistar,\] which is bounded away from $\phistar$, \Nin,\Nout defined in \cref{def:SectionSetsDefs} as the $\delta_*$ union of neighborhoods of $(\phistar,\phistar)$, $(\phistar,\pi-\phistar)$, $(2\pi-\phistar,\phistar)$ and $(2\pi-\phistar,\pi-\phistar)$ must have $\theta$ coordinates  uniformly bounded away from $0$ and $\pi$ (see \cref{fig:MrN}).

Therefore, the $\theta$ coordinates of $\Fcal^{-1}(x)\bigm|_{x\in\mathfrak{l}_U}$, $\Fcal^{2}(x)\bigm|_{x\in\mathfrak{l}_U}$ are bounded away from the boundaries $\theta=0$ and $\theta=\pi$ of $M_r$. Thus, $d_0(x)\overset{\text{\cref{def:fpclf}}}=r\sin{\big(\theta \text{ coordinate of }\Fcal^{-1}(x)\big)}$ and $d_2(x)\overset{\text{\cref{def:fpclf}}}=r\sin{\big(\theta \text{ coordinate of }\Fcal^{2}(x)\big)}$ are bounded away from $0$ for all $x\in\mathfrak{l}_U$. Threfore, the limits $\hat{d}_0>0$, $\hat{d}_2>0$.
Define
\begin{equation}\label{eq:hhatUg1hatU}
    \begin{aligned}
    \hat{h}_1(U)\dfn&\lim_{\mathfrak{l}_U\ni x\to \mathcal{W}_U}h_1(x)\overbracket{=}^{\mathclap{\text{\cref{def:lengthfractionfunctionsn1eq1}\eqref{item01lengthfractionfunctionsn1eq1}}}}\lim_{\mathfrak{l}_U\ni x\to \mathcal{W}_U}\big(\frac{\tau_0}{d_0}+\frac{\tau_1}{d_0}+\frac{2d_1}{d_0}-2\big)\overbracket{=}^{\mathclap{\strut\eqref{eq:lengthlimitWU}\ }}\frac{\hat{\tau}_0(U)}{\hat{d}_0(U)}+\frac{\hat{\tau}_1(U)}{\hat{d}_0(U)}+\frac{2d_1}{\hat{d}_0(U)}-2
    \\
    \hat{g}_1(U)\dfn&\lim_{\mathfrak{l}_U\ni x\to \mathcal{W}_U}g_1(x)\overbracket{=}^{\mathclap{\text{\cref{def:lengthfractionfunctionsn1eq1}\eqref{item03lengthfractionfunctionsn1eq1}}}}\lim_{\mathfrak{l}_U\ni x\to \mathcal{W}_U}\big(\frac{1}{d_0}-\frac{1}{r\sin\phistar}\big)\overbracket=^{\mathclap{\eqref{eq:lengthlimitWU}}}\frac{1}{\hat{d}_0(U)}-\frac{1}{r\sin\phistar}
\end{aligned}
\end{equation}
On the other hand, by the monotone property \cref{proposition:QUpperLowerRegionsLowerBoundsChatCtilde}\eqref{item02Qlowerboundn1eq1}, we have 
\begin{align*}
    0<\frac{\tau_1}{d_0}\Bigm|_{x=x_1}\overbracket{<}^{\mathclap{\text{\cref{proposition:QUpperLowerRegionsLowerBoundsChatCtilde}\eqref{item02Qlowerboundn1eq1}}\strut\qquad\quad}}&\frac{\hat{\tau}_1(U)}{d_0\bigm|_{x=x_1}}\overbracket<^{\mathclap{\qquad\quad\strut\text{\cref{proposition:QUpperLowerRegionsLowerBoundsChatCtilde}\eqref{item02Qlowerboundn1eq1}}}}\frac{\hat{\tau}_1(U)}{\hat{d}_0(U)}\overset{\text{\eqref{eq:hhatUg1hatU}}}=2+\hat{h}_1(U)-\frac{\hat{\tau}_0(U)}{\hat{d}_0(U)}-\frac{2\hat{d}_1(U)}{\hat{d}_0(U)}\\
    \overbracket{\le}^{\mathclap{\text{\eqref{eq:lengthlimitWU}: }\hat{\tau}_0(U)\ge0, d_1=\nu r\sin\phistar}}&2+\hat{h}_1(U)-\frac{2\nu r\sin\phistar}{\hat{d}_0}\overset{\text{\eqref{eq:hhatUg1hatU}}}=2+\hat{h}_1(U)-2\nu r\sin\phistar\big(\hat{g}_1(U)+\frac{1}{r\sin\phistar}\big)\\
    &=2-2\nu+\hat{h}_1(U)-2\nu r\sin\phistar\hat{g}_1(U)\overbracket=^{\mathclap{d_1=\nu r\sin\phistar}}2-2\nu+\hat{h}_1(U)-2d_1\hat{g}_1(U)
\end{align*}

With fixed $k\in[1,4/3]$ and constant $\nu=\frac{d_1}{r\sin\phistar}$, the function $\mathfrak{Q}(t)$ from \eqref{eq:II1II2Quadratic} for $t\in (-\infty,+\infty)$ has only one extreme point, which is the global maximum point. Therefore, for $t\in[0,2-2\nu+\hat{h}_1(U)-2d_1\hat{g}_1(U)]$, $\mathfrak{Q}\big(t\big)$ cannot reach the minimum in the interior of $[0,2-2\nu+\hat{h}_1(U)-2d_1\hat{g}_1(U)]$. Hence, \begin{equation}\label{eq:UpperQlowerbound}
    \mathfrak{Q}\Big(\frac{\tau_1}{d_0}\bigm|_{x=x_1}\Big)>\min_{t\in\big[0,\,\,2-2\nu+\hat{h}_1(U)-2d_1\hat{g}_1(U)\big]}\mathfrak{Q}\big(t\big)=\min\Big\{\mathfrak{Q}\big(0\big),\,\,\mathfrak{Q}\big(2-2\nu+\hat{h}_1(U)-2d_1\hat{g}_1(U)\big)\Big\}.
\end{equation}
For $\frac{d_1}{r\sin\phistar}=\nu\in(0.5,1)$ we estimate \begin{equation}\label{eq:UpperEQlowerbound}
\begin{aligned}
\mathfrak{Q}\big(0\big)\overbracket{=}^{\mathllap{\strut\text{\eqref{eq:II1II2Quadratic}}}}&\frac{4k+9}{3}+\big(\frac{8-16k}{3}\big)\frac{r\sin{\phistar}}{d_1}+\frac{8kd_1}{r\sin{\phistar}}\\
    \overbracket{=}^{\mathclap{\frac{d_1}{r\sin\phistar}=\nu}}&3+\frac{4k}{3}-\big(\frac{16k}{3}-\frac{8}{3}\big)\frac{1}{\nu}+8k\nu\\
    \overbracket{>}^{\mathclap{\substack{ \text{For fixed }k\in[1,4/3],\ \frac{16k}{3}-\frac{8}{3}\ge0\text{, so }\\\nu\mapsto-\big(\frac{16k}{3}-\frac{8}{3}\big)\frac{1}{\nu}+8k\nu\text{ is}\\ 
    \text{increasing with }\nu\in[0.5,1)\strut}}}& 3+\frac{4k}{3}-\big(\frac{16k}{3}-\frac{8}{3}\big)\frac{1}{0.5}+8k\cdot0.5=\frac{25}{3}-\frac{16k}{3}\overbracket{\ge}^{\strut\mathclap{k\le4/3}\strut}\frac{25}{3}-\frac{64}{9}=\frac{11}{9}.\\
    \end{aligned}
    \end{equation} Then we further compute
    \begin{equation}\label{eq:UpperWQlowerbound}
    \begin{aligned}
    &\mathfrak{Q}\big(2-2\nu+\hat{h}_1(U)-2d_1\hat{g}_1(U)\big)\\
    \overbracket{=}^{\mathllap{\text{\eqref{eq:II1II2Quadratic}}}}&\frac{9+4k}{3}+\big(-\frac{16}{3}k+\frac{8}{3}\big)\frac{r\sin\phistar}{d_1}+\frac{8kd_1}{r\sin\phistar}-\frac{4kr\sin\phistar}{d_1}\big(2-2\nu+\hat{h}_1(U)-2d_1\hat{g}_1(U)\big)^2\\
    &+\big[\frac{(32k-12)r\sin\phistar}{3d_1}-8k\big]\big(2-2\nu+\hat{h}_1(U)-2d_1\hat{g}_1(U)\big)\\
    \overbracket{=}^{\mathllap{\nu=\frac{d_1}{r\sin\phistar}}}&\frac{9+4k}{3}+\big(\frac{-16k}{3}+\frac{8}{3}\big)\frac{1}{\nu}+8k\nu-\frac{4k}{\nu}(2-2\nu)^2+\big(-8k+\frac{32k-12}{3}\frac{1}{\nu}\big)\big(2-2\nu\big)\\
    &+\underbracket{\big[\frac{32k-12}{3\nu}-8k\big]\big(\hat{h}_1(U)-2d_1\hat{g}_1(U)\big)-\frac{8k}{\nu}\big(2-2\nu\big)\big(\hat{h}_1(U)-2d_1\hat{g}_1(U)\big)-\frac{4k}{\nu}\big(\hat{h}_1(U)-2d_1\hat{g}_1(U)\big)^2}_{=\mathcal{C}^{\mathcal{W}}_U(\nu,k)}\\    
    &=\frac{9+4k}{3}+32k-16k-\frac{64k-24}{3}+\big(\frac{-16k}{3}+\frac{8}{3}-16k+\frac{64k-24}{3}\big)\frac{1}{\nu}+\big(8k-16k+16k\big)\nu\\
    +&\underbracket{\big[\frac{-16k-12}{3\nu}+8k\big]\big(\hat{h}_1(U)-2d_1\hat{g}_1(U)\big)-\frac{4k}{\nu}\big(\hat{h}_1(U)-2d_1\hat{g}_1(U)\big)^2}_{=\mathcal{C}^{\mathcal{W}}_U(\nu,k)}\\
    &=-4k+11-\frac{16}{3}\frac{1}{\nu}+8k\nu+\mathcal{C}^{\mathcal{W}}_U(\nu,k)\underbracket{>}_{\mathrlap{\substack{\text{Given a } k\in[1,4/3],\\-\frac{16}{3}\frac{1}{\nu}+8k\nu \text{ increases}\\ \text{with }\nu\in[0.5,1]}}\qquad}-4k+11-\frac{16}{3\cdot0.5}+\cdot8k\cdot0.5+\mathcal{C}^{\mathcal{W}}_U(\nu,k)=\frac{1}{3}+\mathcal{C}^{\mathcal{W}}_U(\nu,k),
\end{aligned}\end{equation} where \begin{align*}
    \mathcal{C}^{\mathcal{W}}_U(\nu,k)=\Big[\frac{-16k-12}{3\nu}+8k\Big]\big(\hat{h}_1(U)-2d_1\hat{g}_1(U)\big)-\frac{4k}{\nu}\big(\hat{h}_1(U)-2d_1\hat{g}_1(U)\big)^2.
\end{align*}
Hence, by \eqref{eq:UpperQlowerbound}, \eqref{eq:UpperEQlowerbound} and \eqref{eq:UpperWQlowerbound} we have proved \eqref{item03Qlowerboundn1eq1}.
\hidesubsection{Proof of (4)} With $\nu\in(1/3,0.5]$, for $\mathfrak{l}_L$ in \cref{fig:x1_n1equal1case}, we define length functions limits on its endpoints $\mathcal{W}_L$ and $\mathcal{E}_L$ \newline as follows.

\smallskip\begin{minipage}{.47\textwidth}
    Limit at $\mathcal{W}_L$:
    \begin{equation}\label{eq:lengthlimitWL}
        \begin{aligned}
        \hat{\tau}_0(L)&\dfn\lim_{\mathfrak{l}_L\ni x\to \mathcal{W}_L}\tau_0\ge0\\
        \hat{\tau}_1(L)&\dfn\lim_{\mathfrak{l}_L\ni x\to \mathcal{W}_L}\tau_1\ge0\\
        \hat{d}_0(L)&\dfn\lim_{\mathfrak{l}_L\ni x\to \mathcal{W}_L}d_0>0\\
        \hat{d}_1(L)&\dfn\lim_{\mathfrak{l}_L\ni x\to \mathcal{W}_L}d_1=\const=\nu\cdot r\sin{\phistar}\\
        \hat{d}_2(L)&\dfn\lim_{\mathfrak{l}_L\ni x\to \mathcal{W}_L}d_2>0
    \end{aligned}
    \end{equation}
    \end{minipage}\vrule\quad
    \begin{minipage}{0.47\textwidth}
    Limit at $\mathcal{E}_L$:
    \begin{equation}\label{eq:lengthlimitEL}
        \begin{aligned}
        \tilde{\tau}_0(L)\dfn&\lim_{\mathfrak{l}_L\ni x\to \mathcal{E}_L}\tau_0\ge0\\
        \tilde{\tau}_1(L)\dfn&\lim_{\mathfrak{l}_L\ni x\to \mathcal{E}_L}\tau_1\ge0\\
        \tilde{d}_0(L)\dfn&\lim_{\mathfrak{l}_L\ni x\to \mathcal{E}_L}d_0>0\\
        \tilde{d}_1(L)\dfn&\lim_{\mathfrak{l}_L\ni x\to \mathcal{E}_L}d_1=\nu\cdot r\sin{\phistar}\\
        \tilde{d}_2(L)\dfn&\lim_{\mathfrak{l}_L\ni x\to \mathcal{E}_L}d_2>0
    \end{aligned}
    \end{equation}
    \end{minipage}
    
    These limits in \eqref{eq:lengthlimitWL}, \eqref{eq:lengthlimitEL} exist for the same reason as for the existence of limits in \eqref{eq:lengthlimitWU}. And the continuous extensions of these functions to the end-points $\mathcal{W}_U$, $\mathcal{E}_U$ can be seen in \cref{fig:n1_equal1_lower_leftendpoint_case,fig:n1_equal1_lower_rightendpoint_case}. 
    
    For the continuous extension to $\mathcal{W}_L$,  we observe that for the following. 
    
    $x\in\mathfrak{l}_L$, if $x\rightarrow\mathcal{W}_L$, then $p(\Fc^2(x))\rightarrow B$ which is the corner and $\Fcal(x)\rightarrow\Fcal(\mathcal{W}_L)$ which is an interior point of $M_R$ (see \cref{fig:x1_n1equal1case}), that is, $p(\Fcal(\mathcal{W}_L))$ is an interior point of $\Gamma_R$ (see \cref{fig:n1_equal1_lower_leftendpoint_case}).
    
    If $x\rightarrow\mathcal{E}_L$, then $p(\Fc^{-1}(x))\rightarrow A$ which is the corner and $x\rightarrow\mathcal{E}_L$ which is an interior point of $M_R$ (see \cref{fig:x1_n1equal1case}), that is, $p(\mathcal{E}_L)$ is an interior point of $\Gamma_R$ (see \cref{fig:n1_equal1_lower_rightendpoint_case}).
    Hence \begin{equation}\label{eq:tau0Ltau1Llimits}
        \begin{aligned}
        \hat{\tau}_1(L)&=\lim_{\mathfrak{l}_L\ni x\to \mathcal{W}_L}\overbracket{|p(\Fcal(x))p(\Fcal^2(x))|}^{=\tau_1\bigm|_x}=|p(\Fcal(\mathcal{W}_L))B|\overbracket{=}^{\qquad\qquad\mathllap{B\text{, } p(\Fcal(\mathcal{W}_L))\text{ are on }\Gamma_R}} 2\hat{d}_1(L)=2d_1,\\
        \tilde{\tau}_0(L)&=\lim_{\mathfrak{l}_L\ni x\to \mathcal{E}_L}\overbracket{|p(x)p(\Fcal^{-1}(x))|}^{{=\tau_0\bigm|_x}}=|p(\mathcal{E}_L)A|\overbracket{=}^{\qquad\qquad\mathllap{A\text{, } p(\mathcal{E}_L)\text{ are on }\Gamma_R}} 2\tilde{d}_1(L)=2d_1.
        \end{aligned}
    \end{equation}
    
    Similarly to \eqref{eq:hhatUg1hatU}, we define
    \begin{equation}\label{eq:htildehatLg1tildehatL}
    \begin{aligned}
    \hat{h}_1(L)\dfn&\lim_{\mathfrak{l}_U\ni x\to \mathcal{W}_L}h_1(x)\overbracket{=}^{\mathclap{\text{\cref{def:lengthfractionfunctionsn1eq1}\eqref{item01lengthfractionfunctionsn1eq1}}}}\lim_{\mathfrak{l}_U\ni x\to \mathcal{W}_L}\big(\frac{\tau_0}{d_0}+\frac{\tau_1}{d_0}+\frac{2d_1}{d_0}-2\big)\overset{\eqref{eq:lengthlimitWL}}=\frac{\hat{\tau}_0(L)}{\hat{d}_0(L)}+\frac{\hat{\tau}_1(L)}{\hat{d}_0(L)}+\frac{2d_1}{\hat{d}_0(L)}-2,\\
    \hat{g}_1(L)\dfn&\lim_{\mathfrak{l}_L\ni x\to \mathcal{W}_L}g_1(x)\overbracket{=}^{\mathclap{\text{\cref{def:lengthfractionfunctionsn1eq1}\eqref{item03lengthfractionfunctionsn1eq1}}}}\lim_{\mathfrak{l}_U\ni x\to \mathcal{W}_L}\big(\frac{1}{d_0}-\frac{1}{r\sin\phistar}\big)\overset{\eqref{eq:lengthlimitWL}}=\frac{1}{\hat{d}_0(L)}-\frac{1}{r\sin\phistar},\\
    \tilde{h}_1(L)\dfn&\lim_{\mathfrak{l}_L\ni x\to \mathcal{E}_L}h_1(x)\overbracket{=}^{\mathclap{\text{\cref{def:lengthfractionfunctionsn1eq1}\eqref{item01lengthfractionfunctionsn1eq1}}}}\lim_{\mathfrak{l}_U\ni x\to \mathcal{E}_L}\big(\frac{\tau_0}{d_0}+\frac{\tau_1}{d_0}+\frac{2d_1}{d_0}-2\big)\overset{\eqref{eq:lengthlimitEL}}=\frac{\tilde{\tau}_0(L)}{\tilde{d}_0(L)}+\frac{\tilde{\tau}_1(L)}{\tilde{d}_0(L)}+\frac{2d_1}{\tilde{d}_0(L)}-2,\\
    \tilde{g}_1(L)\dfn&\lim_{\mathfrak{l}_L\ni x\to \mathcal{E}_L}g_1(x)\overbracket{=}^{\mathclap{\text{\cref{def:lengthfractionfunctionsn1eq1}\eqref{item03lengthfractionfunctionsn1eq1}}}}\lim_{\mathfrak{l}_L\ni x\to \mathcal{E}_L}\big(\frac{1}{d_0}-\frac{1}{r\sin\phistar}\big)\overset{\eqref{eq:lengthlimitEL}}=\frac{1}{\tilde{d}_0(L)}-\frac{1}{r\sin\phistar}.
    \end{aligned}
\end{equation}

\begin{figure}[h]
\begin{center}
\begin{tikzpicture}[xscale=1.4,yscale=1.4]
    \tkzDefPoint(0,0){Or};
    \tkzDefPoint(0,-4.8){Y};
    \clip(-5.8,-4.23) rectangle (5.8,-2.85);
    \pgfmathsetmacro{\rradius}{4.5};
    \pgfmathsetmacro{\phistardeg}{40};
    \pgfmathsetmacro{\XValueArc}{\rradius*sin(\phistardeg)};
    \pgfmathsetmacro{\YValueArc}{\rradius*cos(\phistardeg)};
    \pgfmathsetmacro{\Rradius}{
        \rradius*sin(\phistardeg)/sin(18)
    };
    \pgfmathsetmacro{\bdist}{(-\rradius*cos(\phistardeg))+cos(18)*\Rradius}
    \tkzDefPoint(0,\bdist){OR};
    \tkzDefPoint(\XValueArc,-1.0*\YValueArc){B};    \tkzDefPoint(-1.0*\XValueArc,-1.0*\YValueArc){A};
    \tkzDrawArc[name path=Cr,very thick](Or,B)(A);
    \tkzDrawArc[name path=CR,very thick](OR,A)(B);
    \pgfmathsetmacro{\PXValue}{-\Rradius*sin(3)};
    \pgfmathsetmacro{\PYValue}{\bdist-\Rradius*cos(3)};
    \tkzDefPoint(\PXValue,\PYValue){P1};
    
    \tkzDefPointBy[projection=onto OR--P1](A)\tkzGetPoint{M};
    \coordinate (P2) at ($(A)!2!(M)$);
    
    \tkzDefPointBy[projection=onto OR--P2](P1)\tkzGetPoint{M2};
    \coordinate (Q_far) at ($(P1)!2!(M2)$);
    \path[name path=Line_P2_Q_far](P2)-- (Q_far);
    \path[name intersections={of=Cr and Line_P2_Q_far,by={Q}}];
    \tkzDefPoint(0,\bdist-\Rradius){C};
    \begin{scope}[ultra thin,decoration={markings,mark=at position 0.5 with {\arrow{>}}}] 
    \draw[blue,postaction={decorate}] (A)--(P1);
    \draw[blue,postaction={decorate}] (P1)--(P2);
    \draw[blue,postaction={decorate}] (P2)--(Q);
    \end{scope}
    \tkzLabelPoint[above right](A){\small $A=p\big(\Fc^{-1}(x=\mathcal{E}_L)\big)$};
    \tkzLabelPoint[below](P1){\small $P=p\big(x=\mathcal{E}_L\big)$};
    \tkzLabelPoint[below](P2){\small $P_1=p\big(\Fc(x=\mathcal{E}_L)\big)$};
    \tkzLabelPoint[right](Q){\small $Q=p\big(\Fc^{2}(x=\mathcal{E}_L)\big)$};
    \tkzLabelPoint[below](B){$B$};
    \tkzLabelPoint[below](Or){$O_r$};
    \tkzLabelSegment[above left=0.05pt and 0.03pt](P1,P2){\small $2d_1=|PP_1|$};
    \tkzLabelSegment[below left=0.5pt and 0.3pt](A,P1){\small $2d_1=|PA|=\tilde{\tau}_0(L)$};
     \tkzLabelSegment[above left=0.1pt and 0.2pt](P2,Q){\small $|P_1Q|=\tilde{\tau}_1(L)$};
    \tkzDrawPoints(A,B,Or,P1,P2,Q);
    
\end{tikzpicture}
\caption{With $d_1\le0.5r\sin\phistar$, $(\Phi_{\mathcal{E}_L},\theta_1)=x=\mathcal{E}_L$ is on a trajectory from corner $A$}\label{fig:n1_equal1_lower_rightendpoint_case}
\end{center}
\end{figure}
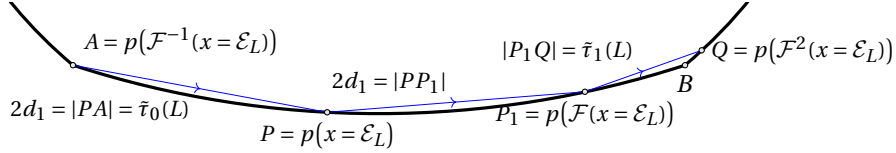

\begin{figure}[h]
\begin{center}
\begin{tikzpicture}[xscale=1.4,yscale=1.4]
    \tkzDefPoint(0,0){Or};
    \tkzDefPoint(0,-4.8){Y};
    \clip(-5.7,-4.23) rectangle (5.7,-2.85);
    \pgfmathsetmacro{\rradius}{4.5};
    \pgfmathsetmacro{\phistardeg}{40};
    \pgfmathsetmacro{\XValueArc}{\rradius*sin(\phistardeg)};
    \pgfmathsetmacro{\YValueArc}{\rradius*cos(\phistardeg)};
    \pgfmathsetmacro{\Rradius}{
        \rradius*sin(\phistardeg)/sin(18)
    };
    \pgfmathsetmacro{\bdist}{(-\rradius*cos(\phistardeg))+cos(18)*\Rradius}
    \tkzDefPoint(0,\bdist){OR};
    \tkzDefPoint(\XValueArc,-1.0*\YValueArc){B};    \tkzDefPoint(-1.0*\XValueArc,-1.0*\YValueArc){A};
    \tkzDrawArc[name path=Cr, very thick](Or,B)(A);
    \tkzDrawArc[name path=CR, very thick](OR,A)(B);
    \pgfmathsetmacro{\PXValue}{\Rradius*sin(3)};
    \pgfmathsetmacro{\PYValue}{\bdist-\Rradius*cos(3)};
    \tkzDefPoint(\PXValue,\PYValue){P2};
    
    \tkzDefPointBy[projection=onto OR--P2](B)\tkzGetPoint{M};
    \coordinate (P1) at ($(B)!2!(M)$);
    
    \tkzDefPointBy[projection=onto OR--P1](P2)\tkzGetPoint{M2};
    \coordinate (T_far) at ($(P2)!2!(M2)$);
    \path[name path=Line_P1_T_far](P1)-- (T_far);
    \path[name intersections={of=Cr and Line_P1_T_far,by={T}}];
    \tkzDefPoint(0,\bdist-\Rradius){C};
    \begin{scope}[ultra thin,decoration={markings,mark=at position 0.5 with {\arrow{>}}}] 
    \draw[blue,postaction={decorate}] (T)--(P1);
    \draw[blue,postaction={decorate}] (P1)--(P2);
    \draw[blue,postaction={decorate}] (P2)--(B);
    \end{scope}
    \tkzLabelPoint[below left](A){$A$};
    \tkzLabelPoint[below](P1){\small $P=p\big(x=\mathcal{W}_L\big)$};
    \tkzLabelPoint[below](P2){\small $P_1=p\big(\Fc(x=\mathcal{W}_L)\big)$};
    \tkzLabelPoint[left](T){\small $T=p(\Fc^{-1}(x=\mathcal{W}_L))$};
    \tkzLabelPoint[below right](B){\small $B=p\big(\Fc^2(x=\mathcal{W}_L)\big)$};
    \tkzLabelPoint[below](Or){$O_r$};
    \tkzLabelSegment[above =0.5pt](P1,P2){\small $2d_1=|PP_1|$};
    \tkzLabelSegment[above right= 0.5pt and 0.2pt](P1,T){\small $\hat{\tau}_0(L)=|TP|$};
    \tkzLabelSegment[above left= 0.5pt and 0.02pt](P2,B){\small $\hat{\tau}_1(L)=|P_1B|=2d_1$};
    \tkzDrawPoints(A,B,P1,P2,T);
\end{tikzpicture}
\caption{With $d_1\le0.5r\sin\phistar$, $(\Phi_{\mathcal{W}_L},\theta_1)=x=\mathcal{W}_L$ is on a trajectory to corner $B$}\label{fig:n1_equal1_lower_leftendpoint_case}


\end{center}
\end{figure}
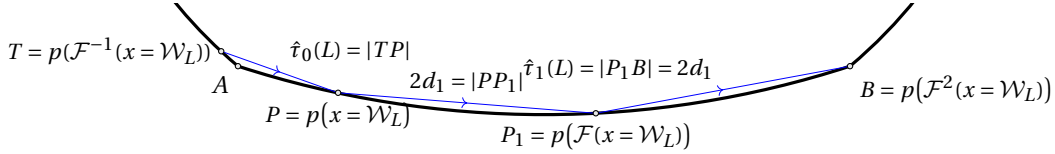

By the previous monotonicity property, item \eqref{item02Qlowerboundn1eq1}, we have
\[\frac{\tilde{\tau}_1(L)}{\tilde{d}_0(L)}\overbracket{<}^{\strut\mathllap{\substack{\text{item \eqref{item02Qlowerboundn1eq1}: }\\\tilde{d}_0(L)>d_0|_{x=x_1}\\|}}}\frac{\tilde{\tau}_1(L)}{d_0|_{x=x_1}}\overbracket{<}^{\strut\mathclap{\substack{\text{item \eqref{item02Qlowerboundn1eq1}: }\\\tilde{\tau}_1(L)<\tau_1|_{x=x_1}\\|}}}\quad\frac{\tau_1}{d_0}\Bigm|_{x=x_1}\overbracket{<}^{\strut\mathclap{\substack{\text{item \eqref{item02Qlowerboundn1eq1}: }\\\hat{\tau}_1(L)>\tau_1|_{x=x_1}\\|}}}\frac{\hat{\tau}_1(L)}{d_0\bigm|_{x=x_1}}\overbracket{<}^{\strut\mathrlap{\substack{\text{item \eqref{item02Qlowerboundn1eq1}: }\\\hat{d}_0(L)<d_0|_{x=x_1}\\|}}}\frac{\hat{\tau}_1(L)}{\hat{d}_0(L)}.\]
And for these two bounds we have 
\begin{align*}
    \frac{\hat{\tau}_1(L)}{\hat{d}_0(L)}\overbracket{=}^{\eqref{eq:tau0Ltau1Llimits}}&\frac{2d_1}{\hat{d}_0(L)}\overbracket{=}^{\eqref{eq:htildehatLg1tildehatL}}2\big(\frac{d_1}{r\sin{\phistar}}+d_1\hat{g}_1(L)\big)\overbracket{=}^{\mathclap{d_1=\nu r\sin\phistar}}2\nu+2d_1\hat{g}_1(L),\\
    \frac{\tilde{\tau}_1(L)}{\tilde{d}_0(L)}\overbracket{=}^{\eqref{eq:htildehatLg1tildehatL}}&2+\tilde{h}_1(L)-\frac{2d_1}{\tilde{d}_0(L)}-\frac{\tilde{\tau}_0(L)}{\tilde{d}_0(L)}\\
    \overbracket{=}^{\qquad\mathclap{\text{\eqref{eq:tau0Ltau1Llimits}: }\tilde{\tau}_0(L)=2d_1}}&2+\tilde{h}_1(L)-\frac{4d_1}{\tilde{d}_0(L)}\\
    \overbracket{=}^{\eqref{eq:htildehatLg1tildehatL}}&2+\tilde{h}_1(L)-4\big(\frac{d_1}{r\sin{\phistar}}+d_1\tilde{g}_1(L)\big)\overbracket{=}^{\mathclap{d_1=\nu r\sin\phistar}}2-4\nu+\tilde{h}_1(L)-4d_1\tilde{g}_1(L).\\
\end{align*}

With fixed $k\in[1,4/3]$ and constant $\nu=\frac{d_1}{r\sin\phistar}$, the function $\mathfrak{Q}(t)$ from \eqref{eq:II1II2Quadratic} for $t\in (-\infty,+\infty)$ has only one extreme point, which is the global maximum point. Therefore, for $t\in\big[2-4\nu+\tilde{h}_1(L)-4d_1\tilde{g}_1(L),\,\,2\nu+2d_1\hat{g}_1(L)\big]$, $\mathfrak{Q}\big(t\big)$ cannot reach the minimum in the interior of $\big[2-4\nu+\tilde{h}_1(L)-4d_1\tilde{g}_1(L),\,\,2\nu+2d_1\hat{g}_1(L)\big]$. Hence,\begin{equation}\label{eq:LowerQlowerbound}
    \mathfrak{Q}\Big(\frac{\tau_1}{d_0}\bigm|_{x=x_1}\Big)>\min_{t\in\big[2-4\nu+\tilde{h}_1(L)-4d_1\tilde{g}_1(L),\,\,2\nu+2d_1\hat{g}_1(L)\big]}\mathfrak{Q}\big(t\big)=\min{\Big\{\mathfrak{Q}\big(2-4\nu+\tilde{h}_1(L)-4d_1\tilde{g}_1(L)\big),\,\,\mathfrak{Q}\big(2\nu+2d_1\hat{g}_1(L)\big)\Big\}.}
\end{equation}
With $\frac{d_1}{r\sin\phistar}=\nu\in(1/3,0.5]$ we first compute $\mathfrak{Q}\big(2-4\nu+\tilde{h}_1(L)-4d_1\tilde{g}_1(L)\big)$.
\begin{equation}\label{eq:LowerEQlowerboundSimplified}
\begin{aligned}
&\mathclap{\qquad\mathfrak{Q}\big(2-4\nu+\tilde{h}_1(L)-4d_1\tilde{g}_1(L)\big)}\\
&\overbracket{=}^{\mathclap{\strut\text{\eqref{eq:II1II2Quadratic}}}}\frac{4k+9}{3}+\big(\frac{8-16k}{3}\big)\frac{r\sin{\phistar}}{d_1}+\frac{8kd_1}{r\sin{\phistar}}-\frac{4kr\sin\phistar}{d_1}\big(2-4\nu+\tilde{h}_1(L)-4d_1\tilde{g}_1(L)\big)^2\\
&\quad+\big[\frac{(32k-12)r\sin\phistar}{3d_1}-8k\big]\big(2-4\nu+\tilde{h}_1(L)-4d_1\tilde{g}_1(L)\big)\\
&=\frac{9+4k}{3}+\frac{8-16k}{3\nu}+8k\nu-\frac{4k(2-4\nu)^2}{\nu}-\frac{8k(2-4\nu)}{\nu}\big(\tilde{h}_1(L)-4d_1\tilde{g}_1(L)\big)\\
&\quad-\frac{4k\big(\tilde{h}_1(L)-4d_1\tilde{g}_1(L)\big)^2}{\nu}+\big(\frac{32k-12}{3\nu}-8k\big)\big(2-4\nu\big)+\big(\frac{32k-12}{3\nu}-8k\big)\big(\tilde{h}_1(L)-4d_1\tilde{g}_1(L)\big)\\
&=\frac{9+4k}{3}+64k-16k-\frac{128k-48}{3}+\big(\frac{8-16k}{3}-16k+\frac{64k-24}{3}\big)\frac{1}{\nu}+\big(8k-64k+32k\big)\nu\\
&\quad+\big(-8k+\frac{32k-12}{3\nu}\big)\big(\tilde{h}_1(L)-4d_1\tilde{g}_1(L)\big)-\frac{8k}{\nu}\big(\tilde{h}_1(L)-4d_1\tilde{g}_1(L)\big)\big(2-4\nu\big)-\frac{4k}{\nu}\big(\tilde{h}_1(L)-4d_1\tilde{g}_1(L)\big)^2\\
&=19+\frac{20k}{3}-\frac{16}{3\nu}-24k\nu+\underbracket{\big(24k+\frac{-16k-12}{3\nu}\big)\big(\tilde{h}_1(L)-4d_1\tilde{g}_1(L)\big)-\frac{4k}{\nu}\big(\tilde{h}_1(L)-4d_1\tilde{g}_1(L)\big)^2}_{\nfd\mathcal{C}^{\mathcal{E}}_{L}(\nu,k)}.
\end{aligned}
\end{equation}
For fixed $k\in[1,4/3]$, $\nu\mapsto19+\frac{20k}{3}-\frac{16}{3\nu}-24k\nu$ has only one extremum on $(0,+\infty)$, which is a maximum. Therefore, its minimum on $[1/3,1/2]$ 
is attained for $\nu=1/2$ or $\nu=1/3$, and
\begin{equation}\label{eq:LowerEQlowerbound}
\begin{aligned}
&\min_{\substack{k\in[1,4/3]\\\nu\in[1/3,1/2]}}{\Big\{\overbracket{19+\frac{20k}{3}-\frac{16}{3\nu}-24k\nu}^{\mathclap{\text{minimum attained }\text{at } \nu=1/2\text{ or }1/3}}\Big\}}\\
&=\min_{k\in[1,4/3]}{\Big\{\min{\big\{19+\frac{20k}{3}-\frac{16}{3\cdot0.5}-24k\cdot0.5,19+\frac{20k}{3}-\frac{16}{3\cdot(1/3)}-24k\cdot(1/3)\big\}}\Big\}}\\
&=\min_{k\in[1,4/3]}{\Big\{\min{\big\{\frac{25}{3}-\frac{16k}{3},3-\frac{4k}{3}\big\}}\Big\}}\overset{k\le4/3}=\frac{11}{9}\text{, hence}\\&
\mathfrak{Q}\big(2-4\nu+\tilde{h}_1(L)-4d_1\tilde{g}_1(L)\big)\overbracket{=}^{\mathllap{\eqref{eq:LowerEQlowerboundSimplified}}}19+\frac{20k}{3}-\frac{16}{3\nu}-24k\nu+\mathcal{C}^{\mathcal{E}}_{L}(\nu,k)\overbracket{\ge}^{\mathllap{k\in[1,4/3],\nu\in(1/3,1/2]}}\frac{11}{9}+\mathcal{C}^{\mathcal{E}}_{L}(\nu,k).
\end{aligned}
\end{equation}
We next compute 
\begin{equation}\label{eq:LowerWQlowerboundSimplified}
    \begin{aligned}
    \mathfrak{Q}\big(2\nu+2d_1\hat{g}_1(L)\big)&\overbracket{=}^{\mathclap{\strut\text{\eqref{eq:II1II2Quadratic}}}}\frac{4k+9}{3}+\big(\frac{8-16k}{3}\big)\frac{r\sin{\phistar}}{d_1}+\frac{8kd_1}{r\sin{\phistar}}-\frac{4kr\sin\phistar}{d_1}\big(2\nu+2d_1\hat{g}_1(L)\big)^2\\
        &\quad+\big[\frac{(32k-12)r\sin\phistar}{3d_1}-8k\big]\big(2\nu+2d_1\hat{g}_1(L)\big)\\
        &\overbracket{=}^{\mathllap{d_1=\nu r\sin\phistar}}\frac{9+4k}{3}+\big(\frac{8-16k}{3\nu}\big)+8k\nu+\big(\frac{32k-12}{3\nu}-8k\big)\big(2\nu+2d_1\hat{g}_1(L)\big)-\frac{4k}{\nu}\big(2\nu+2d_1\hat{g}_1(L)\big)^2\\
        &=\frac{9+4k}{3}+\frac{64k-24}{3}+\big(\frac{8-16k}{3\nu}\big)+(8k-16k-16k)\nu\\
        &\quad+\underbracket{\big(\frac{32k-12}{3\nu}-8k\big)\cdot2d_1\hat{g}_1(L)-\frac{4k\cdot8\nu d_1\hat{g}_1(L)}{\nu}-\frac{16kd_1^2(\hat{g}_1(L))^2}{\nu}}_{=\mathcal{C}^{\mathcal{W}}_{L}(\nu,k)}\\
        &=\frac{68k}{3}-5-\big(\frac{16k-8}{3\nu}\big)-24k\nu+\underbracket{\big(\frac{32k-12}{3\nu}-24k\big)\cdot2d_1\hat{g}_1(L)-\frac{16kd_1^2(\hat{g}_1(L))^2}{\nu}}_{\nfd\mathcal{C}^{\mathcal{W}}_{L}.(\nu,k)}.
    \end{aligned}
\end{equation}
Fix $k\in[1,4/3]$ and note that $\nu\mapsto\frac{68k}{3}-5-\big(\frac{16k-8}{3\nu}\big)-24k\nu$ has only one extreme point on $(0,+\infty)$, a maximum, so it attains its minimum on $[1/3,1/2]$ at 
$\nu=1/2$ or $\nu=1/3$. Thus,
\begin{equation}\label{eq:LowerWQlowerbound}
\begin{aligned}
&\min_{\substack{k\in[1,4/3]\\\nu\in[1/3,1/2]}}{\Big\{\overbracket{\frac{68k}{3}-5-\big(\frac{16k-8}{3\nu}\big)-24k\nu}^{{\mathclap{\text{minimum attained }\text{at } \nu=1/2\text{ or }1/3}}}\Big\}}\\
&=\min_{k\in[1,4/3]}{\Big\{\min{\big\{\frac{68k}{3}-5-\frac{16k-8}{3\cdot0.5}-24k\cdot0.5,\frac{68k}{3}-5-\big(\frac{16k-8}{3\cdot(1/3)}\big)-24k\cdot(1/3)\big\}}\Big\}}\\
&=\min_{k\in[1,4/3]}{\Big\{\min{\big\{\frac{1}{3},3-\frac{4k}{3}\big\}}\Big\}}\overbracket{=}^{\mathclap{k\le4/3}\strut}\frac{1}{3}.\\
\text{Hence }&\mathfrak{Q}\big(2\nu+2d_1\hat{g}_1(L)\big)\overbracket{=}^{\mathllap{\eqref{eq:LowerWQlowerboundSimplified}}}\frac{68k}{3}-5-\big(\frac{16k-8}{3\nu}\big)-24k\nu+\mathcal{C}^{\mathcal{W}}_{L}(\nu,k)\overbracket{\ge}^{\mathllap{k\in[1,4/3],\nu\in(1/3,1/2]}}\frac{1}{3}+\mathcal{C}^{\mathcal{W}}_{L}(\nu,k).
\end{aligned}    
\end{equation}

Now \eqref{eq:LowerQlowerbound}, \eqref{eq:LowerEQlowerbound} and \eqref{eq:LowerWQlowerbound} imply \eqref{item04Qlowerboundn1eq1}.
\hidesubsection{Proof of (5)} We first look at $\mathcal{C}^{\mathcal{W}}_{U}(\nu,k)$ for $k\in[1,4/3]$, $\nu\in(1/2,1)$.
We have \begin{equation}\label{eq:2ndproducttermCUWnuk}
    \begin{aligned}
    &\frac{-8k-12}{3}=\frac{-16k-12}{3\times0.5}+8k\overbracket{<}^{\mathclap{\substack{-16k-12<0,\\\nu>0.5}}}\frac{-16k-12}{3\nu}+8k\overbracket{<}^{\mathclap{\substack{-16k-12<0,\\\nu<1}}}\frac{-16k-12}{3\times1}+8k=\frac{8k-12}{3}\overbracket<^{\mathclap{k\le4/3}\strut}0. 
\end{aligned}
\end{equation}
On the other hand, by \eqref{eq:lengthlimitWU}, \eqref{eq:hhatUg1hatU} $\hat{h}_1(U)$, $\hat{g}_1(U)$ are the limits of $h_1(x)$, $g_1(x)$ which are continuous functions of $x\in \mathfrak{l}_U$. Hence:
\begin{equation}\label{eq:hat1Ug1Ubounds}
\begin{aligned}
&\text{\cref{lemma:h1g1g0g1g2g3lengthfractionestimate}\eqref{item01h1g1g0g1g2g3lengthfractionestimate}: } \forall x\in\mathfrak{l}_U\text{, }h_1(x)\in\big(-8.53r/R,\quad16.65r/R\big)\\
\Longrightarrow&\hat{h}_1(U)\in\big[-8.53r/R,\quad 16.65r/R\big]\\
&\text{and \cref{lemma:h1g1g0g1g2g3lengthfractionestimate}\eqref{item03h1g1g0g1g2g3lengthfractionestimate}: } \forall x\in\mathfrak{l}_U\text{, }g_1(x)\in\big(-\frac{17r}{4R}\frac{1}{0.9975r\sin\phistar},\quad \frac{17r}{4R}\frac{1}{0.9975r\sin\phistar}\big)\\
\Longrightarrow&\hat{g}_1(U)\in\big[-\frac{17r}{4R}\frac{1}{0.9975r\sin\phistar},\quad \frac{17r}{4R}\frac{1}{0.9975r\sin\phistar}\big],\\
\Longrightarrow&-2d_1\hat{g}_1(U)\in\big[-\frac{17r}{2R}\frac{d_1}{0.9975r\sin\phistar},\quad \frac{17r}{2R}\frac{d_1}{0.9975r\sin\phistar}\big]\overset{\nu=\frac{d_1}{r\sin\phistar}}=\big[-\frac{8.5\nu}{0.9925}\frac{r}{R},\quad +\frac{8.5\nu}{0.9975}\frac{r}{R}\big].\\
\end{aligned}\end{equation}Hence 
\begin{equation}\label{eq:1stproducttermCUWnuk}
\hat{h}_1(U)-2d_1\hat{g}_1(U)\in\Big[-\big(8.53+\frac{8.5\nu}{0.9975}\big)\frac{r}{R},\quad \big(16.65+\frac{8.5\nu}{0.9975}\big)\frac{r}{R}\Big].
\end{equation}
Therefore, 
\begin{equation}\label{eq:CWUnukbound}
\begin{aligned}
    \mathcal{C}^{\mathcal{W}}_{U}(\nu,k)&=\Big(\frac{-16k-12}{3\nu}+8k\Big)\big(\hat{h}_1(U)-2d_1\hat{g}_1(U)\big)-\frac{4k}{\nu}\big(\hat{h}_1(U)-2d_1\hat{g}_1(U)\big)^2\\
    \overbracket{>}^{\strut\mathclap{\text{\eqref{eq:1stproducttermCUWnuk},\eqref{eq:2ndproducttermCUWnuk}}}}&\big(\underbracket{\frac{-8k-12}{3}}_{<0}\big)\big(16.65+\frac{8.5\nu}{0.9975}\big)\frac{r}{R}-\frac{4k}{\nu}\Big[\big(16.65+\frac{8.5\nu}{0.9975}\big)\frac{r}{R}\Big]^2\\
    \overbracket{>}^{\mathclap{\nu<1}}&\big(\frac{-8k-12}{3}\big)\big(16.65+\frac{8.5}{0.9975}\big)\frac{r}{R}-4k\underbracket{\big(\frac{16.65}{\nu}+\frac{8.5}{0.9975}\big)\big(16.65+\frac{8.5\nu}{0.9975}\big)}_{<875 \text{ with }0.5\le\nu\le1}\big(\frac{r}{R}\big)^2\\
    &>-(68k+101)\frac{r}{R}-3500k\big(\frac{r}{R}\big)^2\overset{k\le4/3}>-192\frac{r}{R}-4667\big(\frac{r}{R}\big)^2.
\end{aligned}
\end{equation}
We next look at $\mathcal{C}^{\mathcal{W}}_L(\nu,k)$ for $k\in[1,4/3]$, $\nu\in(1/3,1/2]$. We have \begin{equation}\label{eq:1stproducttermCWLnuk}
\begin{aligned}
\frac{-8k-24}{3}=-24k+\frac{32k-12}{3\cdot0.5}\overbracket{\le}^{\strut\mathclap{32k-12>0,\nu\le1/2}}-24k+\frac{32k-12}{3\nu}\overbracket{<}^{\strut\mathclap{32k-12>0,\nu>1/3}}-24k+32k-12=8k-12\overbracket<^{\mathclap{\strut k\le4/3}}0.   
\end{aligned}    
\end{equation}
By \eqref{eq:lengthlimitWL}, \eqref{eq:htildehatLg1tildehatL} $\hat{g}_1(L)$ is the limit of $h_1(x)$,  which is a continuous function of $x\in \mathfrak{l}_L$. Hence we find
\begin{equation}\label{eq:2ndproducttermCWLnuk}
\begin{aligned}&\text{\cref{lemma:h1g1g0g1g2g3lengthfractionestimate}\eqref{item03h1g1g0g1g2g3lengthfractionestimate}: }\forall x\in\mathfrak{l}_L\text{, } g_1(x)\in \big(-\frac{17r}{4R}\frac{1}{0.9975r\sin\phistar},\ \frac{17r}{4R}\frac{1}{0.9975r\sin\phistar}\big)\\
\Longrightarrow&\hat{g}_1(L)\in\big[-\frac{17r}{4R}\frac{1}{0.9975r\sin\phistar},\ \frac{17r}{4R}\frac{1}{0.9975r\sin\phistar}\big]\\
\Longrightarrow&2d_1\hat{g}_1(L)\in\big[-\frac{17r}{2R}\frac{d_1}{0.9975r\sin\phistar},\ \frac{17r}{2R}\frac{d_1}{0.9975r\sin\phistar}\big] \overset{{\nu=\frac{d_1}{r\sin\phistar}}}{=}
\big[-\frac{8.5\nu}{0.9975}\frac{r}{R},\quad\frac{8.5\nu}{0.9975}\frac{r}{R}\big],
\end{aligned}
\end{equation}
Therefore \begin{equation}\label{eq:CWLnukbound}
    \begin{aligned}
        \mathcal{C}^{\mathcal{W}}_{L}(\nu,k)&=\underbracket{\big(\frac{32k-12}{3\nu}-24k\big)\cdot2d_1\hat{g}_1(L)}_{\eqref{eq:1stproducttermCWLnuk}: <0}-\frac{16kd_1^2(\hat{g}_1(L))^2}{\nu}\\
        \overbracket{>}^{\mathclap{\strut\text{\eqref{eq:1stproducttermCWLnuk},\eqref{eq:2ndproducttermCWLnuk}}}}&\frac{-8k-24}{3}\frac{8.5\nu}{0.9975}\frac{r}{R}-\frac{4k}{\nu}\big(\frac{8.5}{0.9975}\frac{r}{R}\big)^2\nu^2\\
        \overbracket{>}^{\mathllap{1/3<\nu\le0.5}}&\frac{\overbracket{-8k-24}^{<0}}{3}\frac{8.5\times0.5}{0.9975}\frac{r}{R}-4k\cdot0.5\cdot\big(\frac{8.5}{0.9975}\frac{r}{R}\big)^2\overbracket>^{\mathclap{k\le4/3}}-50\big(\frac{r}{R}\big)-194\big(\frac{r}{R}\big)^2.
    \end{aligned}
\end{equation}
We finally look at $\mathcal{C}^{\mathcal{E}}_L(\nu,k)$ for $k\in[1,4/3]$, $\nu\in(1/3,1/2]$. Here,
\begin{equation}\label{eq:2ndproducttermCLEnuk}
    -4\overbracket{\le}^{k\ge1}8k-12\overbracket{<}^{\nu>1/3}24k+\frac{\overbracket{-16k-12}^{<0}}{3\nu}\overbracket{\le}^{\nu\le0.5}24k+\frac{-16k-12}{3\cdot0.5}=\frac{40k-24}{3}\overbracket{\le}^{k\le4/3}\frac{88}{9}.
\end{equation}

On the other hand, by the same reason as for \eqref{eq:hat1Ug1Ubounds}, we have
\begin{align*}
    \tilde{h}_1(L)\in&\big[-8.53r/R,\  16.65r/R\big],\\
    -2d_1\tilde{g}_1(L)\in&\big[-\frac{8.5\nu}{0.9975}\frac{r}{R},\ \frac{8.5\nu}{0.9975}\frac{r}{R}\big].
\end{align*} Hence 
\begin{equation}\label{eq:1stproducttermCLEnuk}
\tilde{h}_1(L)-4d_1\tilde{g}_1(L)\in\Big[-\big(8.53+\frac{17\nu}{0.9975}\big)\frac{r}{R},\quad \big(16.65+\frac{17\nu}{0.9975}\big)\frac{r}{R}\Big].
\end{equation} Therefore\COMMENT{Does the next equation need to be numbered? I am pretty sure that it does not.}
\begin{equation}\label{eq:CELnukbound}
\begin{aligned}
    \mathcal{C}^{\mathcal{E}}_{L}(\nu,k)&=\big(24k+\frac{-16k-12}{3\nu}\big)\big(\tilde{h}_1(L)-4d_1\tilde{g}_1(L)\big)-\frac{4k}{\nu}\big(\tilde{h}_1(L)-4d_1\tilde{g}_1(L)\big)^2\\
    \overbracket{\ge}^{\mathclap{\text{\eqref{eq:1stproducttermCLEnuk},\eqref{eq:2ndproducttermCLEnuk}}}}&-\big(8.53+\frac{17\nu}{0.9975}\big)\frac{r}{R}\cdot\frac{88}{9}-\frac{4k}{\nu}\big(16.65+\frac{17\nu}{0.9975}\big)^2\big(\frac{r}{R}\big)^2\\
    \overbracket{\ge}^{\mathclap{\nu\le0.5}}&-\big(8.53+\frac{17\cdot0.5}{0.9975}\big)\frac{r}{R}\cdot\frac{88}{9}-4\overbracket{k}^{\mathclap{k\le4/3}}\underbracket{\big(16.65+\frac{17\nu}{0.9975}\big)\big(\frac{16.65}{\nu}+\frac{17}{0.9975}\big)}_{<1497\text{ with }\nu\in(1/3,1/2]}\big(\frac{r}{R}\big)^2\\
    &>-92\big(\frac{r}{R}\big)-7984\big(\frac{r}{R}\big)^2.
\end{aligned}
\end{equation}\eqref{eq:CWUnukbound},\eqref{eq:CELnukbound}, \eqref{eq:CWLnukbound} are the claim of \eqref{item05Qlowerboundn1eq1}\qedhere
\end{proof}
Finally we get the following conclusion.
\begin{proposition}\label{proposition:casecexpansionn1eq1}
        With the standing assumptions of this section (\cref{REMCaseCAssumptions}), if $n_1=1$, then $\mathcal{II}_1(\mathcal{B}^+_0,1)$, $\mathcal{II}_2(\mathcal{B}^+_0,1)$ in \cref{def:AformularForLowerBoundOfExpansionForn1largerthan0} and in \eqref{itemAformularForLowerBoundOfExpansionForn1largerthan04} of \cref{lemma:AformularForLowerBoundOfExpansionForn1largerthan0} satisfy \[E(1,\Bcal^+_0)=\frac{\|dx_3\|_\p}{\|dx_0\|_\p}=3\big|\mathcal{II}_1(\mathcal{B}^+_0,1)+\mathcal{II}_2(\mathcal{B}^+_0,1)\big|>1-1743\big(\frac{r}{R}\big)-15450\big(\frac{r}{R}\big)^2>0.9.\]
\end{proposition}
\begin{proof}
With $n_1=1$, $\frac{\|dx_3\|_\p}{\|dx_0\|_\p}\overbracket{=}^{\mathclap{\strut\text{\cref{lemma:AformularForLowerBoundOfExpansionForn1largerthan0}\eqref{itemAformularForLowerBoundOfExpansionForn1largerthan04}}}}3\big|\mathcal{II}_1(\mathcal{B}^+_0,1)+\mathcal{II}_2(\mathcal{B}^+_0,1)\big|\overbracket{=}^{\mathclap{\strut\text{\cref{def:II1kAndII2k}}}}3\big|\mathcal{II}_1(k)+\mathcal{II}_2(k)\big|\overbracket{=}^{\mathclap{\strut\text{\cref{remark:II1II2Quadratic},\eqref{eq:II1II2QuadraticAndC1C2}}}}3\big|\mathfrak{Q}\big(\frac{\tau_1}{d_0}\big)+\mathcal{C}_1(k)+\mathcal{C}_2(k)\big|$ with $\Bcal^+_0=\frac{-k}{d_0}$ for some $k\in[1,\frac{4}{3}]$.

With $\nu=\frac{d_1}{r\sin\phistar}$, by \cref{proposition:QUpperLowerRegionsLowerBoundsChatCtilde}\eqref{item03Qlowerboundn1eq1},\eqref{item04Qlowerboundn1eq1} and \eqref{item05Qlowerboundn1eq1}\begin{align*}
    \mathfrak{Q}\big(\frac{\tau_1}{d_0}\big)&\overbracket{>}^{\mathclap{\text{\eqref{item03Qlowerboundn1eq1},\eqref{item04Qlowerboundn1eq1}}}}\min{\Big\{\frac{11}{9},\frac{1}{3}+\mathcal{C}^{\mathcal{W}}_{U}(\nu,k),\frac{11}{9}+\mathcal{C}^{\mathcal{E}}_{L}(\nu,k),\frac{1}{3}+\mathcal{C}^{\mathcal{W}}_{L}(\nu,k)\Big\}}\\
    &\overbracket{>}^{\mathclap{\text{\eqref{item05Qlowerboundn1eq1}}}}\min{\Big\{\frac{11}{9},\frac{1}{3}-192\big(\frac{r}{R}\big)-4667\big(\frac{r}{R}\big)^2,\frac{11}{9}-92\big(\frac{r}{R}\big)-7984\big(\frac{r}{R}\big)^2,\frac{1}{3}-50\big(\frac{r}{R}\big)-194\big(\frac{r}{R}\big)^2\Big\}}\\
    &\overbracket{=}^{\mathclap{\strut\text{\eqref{eqJZHypCond}:}R>1700r}}\frac{1}{3}-192\big(\frac{r}{R}\big)-4667\big(\frac{r}{R}\big)^2.
\end{align*}
Hence \begin{align*}
    3\Big(\mathfrak{Q}\big(\frac{\tau_1}{d_0}\big)+\overbracket{\mathcal{C}_1(k)+\mathcal{C}_2(k)}^{\mathclap{\text{\cref{proposition:C1C2lowerbound}:}>-389(\frac{r}{R})-483(\frac{r}{R})^2}}\Big)&>3\Big(\frac{1}{3}-192\big(\frac{r}{R}\big)-4667\big(\frac{r}{R}\big)^2-389\big(\frac{r}{R}\big)-483\big(\frac{r}{R}\big)^2\Big)=3\Big(\frac{1}{3}-581\big(\frac{r}{R}\big)-5150\big(\frac{r}{R}\big)^2\Big)\\
    &=1-1743\big(\frac{r}{R}\big)-15450\big(\frac{r}{R}\big)^2\overbracket{>}^{\mathrlap{\eqref{eqJZHypCond}:R>30000r}}0.9,\\
    \frac{\|dx_3\|_\p}{\|dx_0\|_\p}&>1-1743\big(\frac{r}{R}\big)-15450\big(\frac{r}{R}\big)^2>0.9.\qedhere
\end{align*}
\end{proof}
\begin{theorem}[With the standing assumptions of this section (\cref{REMCaseCAssumptions})]\label{theorem:B3minusdx3dphi3ratiobound}
     The wave front curvature $\Bcal^-_3$ at $x_3$ in \cref{lemma:AformularForLowerBoundOfExpansionForn1largerthan0}\eqref{itemAformularForLowerBoundOfExpansionForn1largerthan05} and \eqref{eq:B3minusEquation} is uniformly bounded. i.e., $|\Bcal^-_3|<\const$ in case (c) of \eqref{eqMainCases} for some uniform constant.
     
     Furthermore, the tangent vector $dx_3$ in \cref{lemma:AformularForLowerBoundOfExpansionForn1largerthan0}\eqref{itemAformularForLowerBoundOfExpansionForn1largerthan03} $(d\phi_3,d\theta_3)=dx_3=D_{x_0}\Fc^{n_1+4}(dx_0)$ has the slope $\frac{d\theta_3}{d\phi_3}$ also uniformly bounded: there is a uniform upper bound $\lambda_2(r,R,\phistar)>0$ such that $0<\frac{d\theta_3}{d\phi_3}<\lambda_2(r,R,\phistar)$, where constant $\lambda_2(r,R,\phistar)$ is determined by lemon billiard configuration. 
     
     \COMMENT{\Large This reads like a proof rather than a Theorem. That makes is harder to follow. \red\bfseries Edit\wentao{Done}}
\end{theorem}
\begin{proof}
    It suffices to prove each term in \eqref{eq:B3minusEquation} is uniformly bounded. The first term:\COMMENT{\red Misalignment!}
\begin{equation}\label{eq:Bmins3firstterm}
    \begin{aligned}\bigm|\frac{-2}{d_2(1+2d_2\Bcal^+_2)}\bigm|&\overbracket{=}^{\mathclap{\eqref{eq:Expansionx3x0}}}\frac{2}{d_2}\frac{|1+\tau_0\Bcal^+_0||1+\tau_{1,\text{\upshape im}}\Bcal^+_{1,\text{\upshape im}}||1+\tau_1\Bcal^+_{n_1,R}|}{E(n_1,\Bcal^+_0)}\\
    &\overbracket{=}^{\mathllap{\text{\cref{lemma:AformularForLowerBoundOfExpansionForn1largerthan0}\eqref{itemAformularForLowerBoundOfExpansionForn1largerthan04}}}}\frac{|1+\tau_0\Bcal^+_0|}{E(n_1,\Bcal^+_0)}\bigm|1+(-2n_1d_1)(\Bcal^{-}_1-\frac{1}{d_1})\bigm||1+\tau_1\Bcal^+_{n_1,R}|\\
    &\overbracket{=}^{\mathllap{\text{Proof of \cref{lemma:AformularForLowerBoundOfExpansionForn1largerthan0}\eqref{itemAformularForLowerBoundOfExpansionForn1largerthan04}\ref{item04ContinuedfractionResults}}}}\frac{|1+\tau_0\Bcal^+_0|}{E(n_1,\Bcal^+_0)}\bigm|1+(-2n_1d_1)(\Bcal^{-}_1-\frac{1}{d_1})\bigm|\bigm|1+\tau_1[-\frac{1}{d_1}+\frac{\Bcal^{-}_1-\frac{1}{d_1}}{1-2n_1d_1(\Bcal^{-}_1-\frac{1}{d_1})}]\big|\\
    &=\frac{|1+\tau_0\Bcal^+_0|}{E(n_1,\Bcal^+_0)}\bigm|1-\frac{\tau_1}{d_1}+2n_1(\tau_1-d_1)(\Bcal^-_1-\frac{1}{d_1})\bigm|\\
    &=\frac{|1+\tau_0\Bcal^+_0|}{E(n_1,\Bcal^+_0)}\bigm|1+2n_1-(2n_1+1)\frac{\tau_1}{d_1}+2n_1(\tau_1-d_1)\Bcal^-_1\big|\\
    &\le\frac{|1+\tau_0\Bcal^+_0|}{E(n_1,\Bcal^+_0)}\bigm|1+2n_1-(2n_1+1)\frac{\tau_1}{d_1}\bigm|+\frac{|\Bcal^-_1||1+\tau_0\Bcal^+_0|}{E(n_1,\Bcal^+_0)}\bigm|2n_1(\tau_1-d_1)\big|\\
    &\overbracket{<}^{\mathllap{\substack{\text{\cref{remark:deffpclf}: }0<\tau_1<2d_1\text{ and }\\\text{\cite[equation (3.31)]{cb}: }\\\Bcal^-_1(1+\tau_0\Bcal^+_0)=\Bcal_0^+}}}\frac{|1+\tau_0\Bcal^+_0|(2n_1+1)}{E(n_1,\Bcal^+_0)}+\frac{2n_1d_1|\Bcal^+_0|}{E(n_1,\Bcal^+_0)}.
    \end{aligned}
    \end{equation}
Now we notice the following.
\begin{itemize}
    \item By \cref{lemma:AformularForLowerBoundOfExpansionForn1largerthan0}\eqref{itemAformularForLowerBoundOfExpansionForn1largerthan01}: $\Bcal^+_0\in[\frac{-4}{3d_0},\frac{-1}{d_0}]$ and \cref{proposition:casecparametersestimate}\ref{itemcasecparametersestimate1}: $|d_0-r\sin\phistar|<\frac{17R}{4R}r\sin\phistar$, thus $\Bcal^+_0$ is uniformly bounded.
    \item $1+\tau_0\Bcal^+_0\in(1-\frac{4\tau_0}{3d_0},1-\frac{\tau_0}{d_0})\subset(-\frac{5}{3},1)$ since $0<\tau_0<2d_0$. And $2n_1d_1<2R\Phistar$ the arc length of $\Gamma_R$.
    \item By \cref{proposition:casecexpansionn1ge4,proposition:casecexpansionn1eq3,proposition:casecexpansionn1eq2,proposition:casecexpansionn1eq1}, for every $n_1>0$, $E(n_1,\Bcal^+_0)>0.9$. And in particular by \cref{proposition:casecexpansionn1ge4}, $\frac{2n_1+1}{E(n_1,\Bcal^+_0)}\rightarrow0$ as $n_1\rightarrow+\infty$.
\end{itemize}
Thus, the terms on the right of \eqref{eq:Bmins3firstterm} are uniformly bounded. Hence so is $\big|\frac{-2}{d_2(1+2d_2\Bcal^+_2)}\big|$.

As to the second term in \eqref{eq:B3minusEquation}:
\begin{equation}\label{eq:Bmins3secondterm}
    \begin{aligned}
        \big|-\frac{1}{d_1(1+2d_2\Bcal_2^+)(1+\tau_1\Bcal^{+}_{n_1,R})}\bigm|\overbracket{=}^{\mathclap{\eqref{eq:Expansionx3x0}}}&\frac{1}{d_1}\frac{|1+\tau_0\Bcal^+_0||1+\tau_{1,\text{\upshape im}}\Bcal^+_{1,\text{\upshape im}}|}{E(n_1,\Bcal^+_0)}\\
        \overbracket{=}^{\mathllap{\text{\cref{lemma:AformularForLowerBoundOfExpansionForn1largerthan0}\eqref{itemAformularForLowerBoundOfExpansionForn1largerthan04}}}}&\frac{1}{d_1}\frac{|1+\tau_0\Bcal_0^+||1-2n_1d_1(\Bcal^-_1-\frac{1}{d_1})|}{E(n_1,\Bcal^+_0)}\\
        =&\frac{1}{d_1}\frac{|(1+2n_1)(1+\tau_0\Bcal^+_0)-2n_1d_1\Bcal^-_1(1+\tau_0\Bcal^+_0)|}{E(n_1,\Bcal_0^+)}\\
        \le&\frac{|1+\tau_0\Bcal^+_0|(2n_1+1)}{d_1E(n_1,\Bcal_0^+)}+\frac{2n_1|\Bcal^+_0|}{E(n_1,\Bcal^+_0)}\\
        \overbracket{<}^{\mathllap{\text{\cref{proposition:casecparametersestimate}\ref{itemcasecparametersestimate3}: } d_1>\frac{1}{n_1+2}r\sin\phistar}}&\frac{|1+\tau_0\Bcal^+_0|(2n_1+1)(n_1+2)}{r\sin\phistar E(n_1,\Bcal^+_0)}+\frac{2n_1|\Bcal^+_0|}{E(n_1,\Bcal^+_0)}.
    \end{aligned}
\end{equation}
As in the estimate of the first term, both $|\Bcal^+_0|$ and $|1+\tau_0\Bcal^+_0|$ are uniformly bounded. \cref{proposition:casecexpansionn1ge4,proposition:casecexpansionn1eq3,proposition:casecexpansionn1eq2,proposition:casecexpansionn1eq1} give $E(n_1,\Bcal^+_0)>0.9$ for each $n_1>0$. By \cref{proposition:casecexpansionn1ge4}, $\limsup_{n_1\rightarrow\infty}{\frac{(2n_1+1)(n_1+2)}{E(n_1,\Bcal^+_0)}}< C_{c}$ for some constant $C_{c}>0$, and $\lim_{n_1\rightarrow\infty}\frac{n_1}{E(n_1,\Bcal^+_0)}=0$. Therefore, the last two terms in \eqref{eq:Bmins3secondterm} are uniformly bounded, and hence so is $-\frac{1}{d_1(1+2d_2\Bcal_2^+)(1+\tau_1\Bcal^{+}_{n_1,R})}$.

The third term in \eqref{eq:B3minusEquation}: 
\begin{equation}\label{eq:B3minusthirdterm}
    \begin{aligned}
            \big|-\frac{1}{d_1(1+2d_2\Bcal^+_2)(1+\tau_1\Bcal^+_{n_1,R})(1+\tau_{1,\text{\upshape im}}\Bcal^+_{1,\text{\upshape im}})}\bigm|\overbracket{=}^{\mathclap{\eqref{eq:Expansionx3x0}}}&\frac{|1+\tau_0\Bcal_0^+|}{d_1E(n_1,\Bcal_0^+)}\\
            \overbracket{<}^{\mathllap{\text{\cref{proposition:casecparametersestimate}\ref{itemcasecparametersestimate3}: } d_1>\frac{1}{n_1+2}r\sin\phistar}}&\frac{|1+\tau_0\Bcal^+_0|(n_1+2)}{r\sin\phistar E(n_1,\Bcal_0^+)}.
    \end{aligned}
\end{equation}
As before, $|1+\tau_0\Bcal^+_0|$ is uniformly bounded. \cref{proposition:casecexpansionn1ge4,proposition:casecexpansionn1eq3,proposition:casecexpansionn1eq2,proposition:casecexpansionn1eq1} give $E(n_1,\Bcal^+_0)>0.9$ for each $n_1>0$. By \cref{proposition:casecexpansionn1ge4},  $\lim_{n_1\rightarrow\infty}\frac{n_1+2}{E(n_1,\Bcal^+_0)}=0$. Therefore, the right-hand side of \eqref{eq:B3minusthirdterm} is uniformly bounded, and hence so is the left-hand side $|-\frac{1}{d_1(1+2d_2\Bcal^+_2)(1+\tau_1\Bcal^+_{n_1,R})(1+\tau_{1,\text{\upshape im}}\Bcal^+_{1,\text{\upshape im}})}|$.

The fourth term in \eqref{eq:B3minusEquation} is $\big|\frac{\Bcal^+_0}{(1+2d_2\Bcal^+_2)(1+\tau_1\Bcal^+_{n_1,R})(1+\tau_{1,\text{\upshape im}}\Bcal^+_{1,\text{\upshape im}})(1+\tau_0\Bcal_0^+)}\bigm|=\frac{|\Bcal^+_0|}{E(n_1,\Bcal^+_0)}.$
The same proof as for the first term shows that $|\Bcal^+_0|$ is uniformly bounded by some constant. \cref{proposition:casecexpansionn1ge4,proposition:casecexpansionn1eq3,proposition:casecexpansionn1eq2,proposition:casecexpansionn1eq1} imply that $E(n_1,\Bcal^+_0)>0.9$ for each $n_1>0$.\COMMENT{\red It seems that we say ``\cref{proposition:casecexpansionn1ge4,proposition:casecexpansionn1eq3,proposition:casecexpansionn1eq2,proposition:casecexpansionn1eq1} imply that $E(n_1,\Bcal^+_0)>0.9$ for each $n_1>0$'' \emph{\Large\bfseries four} times. Once should suffice, right?\wentao{hmmm,kind of repeating for different reason.}\Hrule Maybe you can explain\dots} Therefore $|\frac{\Bcal^+_0}{(1+2d_2\Bcal^+_2)(1+\tau_1\Bcal^+_{n_1,R})(1+\tau_{1,\text{\upshape im}}\Bcal^+_{1,\text{\upshape im}})(1+\tau_0\Bcal_0^+)}|$ is uniformly bounded by some constant.

Now we get the fourth term in \eqref{eq:B3minusEquation} are uniformly bounded by some constant and we have already showed in \eqref{eq:B3minusthirdterm}\eqref{eq:Bmins3firstterm}\eqref{eq:Bmins3secondterm} that the first three terms in \eqref{eq:B3minusEquation} are uniformly bounded by some constant. Therefore $\Bcal^-_3$ is uniformly bounded by some positive constant. 

Furthermore, since \eqref{eq:MirrorEquationExpansionFormula}: $\frac{d\theta_3}{d\phi_3}=-r(\Bcal^-_3\cos\varphi_3-\frac{1}{r})$, for the $dx_3$ in \cref{lemma:AformularForLowerBoundOfExpansionForn1largerthan0}\eqref{itemAformularForLowerBoundOfExpansionForn1largerthan03}, $(d\phi_3,d\theta_3)=dx_3=D_{x_0}\Fc^{n_1+4}(dx_0)$ has $\frac{d\theta_3}{d\phi_3}$ also uniformly bounded. Since \cref{caseC} implies that $\frac{d\theta_3}{d\phi_3}\in(0,+\infty)$ (in the cone), there is a uniform upper bound $\lambda_2(r,R,\phistar)>0$ such that $0<\frac{d\theta_3}{d\phi_3}<\lambda_2(r,R,\phistar)$, where the constant $\lambda_2(r,R,\phistar)$ is determined by lemon billiard configuration.\qedhere
\end{proof}
This concludes the proof of \cref{TEcase_c}, the arduous work on the hardest case of those defined in \eqref{eqMainCases}. 

The subsections proceeded according to the number $n_1$ of subsequent collisions on $\Gamma_R$  (\cref{def:MhatReturnOrbitSegment}) being at least 4, and equal to 3, 2, and 1, respectively. This becomes increasingly harder, and this last subsection alone took up pages \pageref{SBSn1=1}--\pageref{sec:SACTWFELBWB}.

\NEWPAGE
\section[Feasible  \texorpdfstring{$\phistar$}{open angle} parameter for 1-petal billiard and uniform hyperbolicity of lemon billiard]{Feasible  \texorpdfstring{$\phistar$}{open angle} parameter for 1-petal billiard and uniform hyperbolicity of lemon billiard}\label{sec:SACTWFELBWB}
In this section, we determine values of the parameter $\phistar$ for which \cref{RMASUH} holds and determine the associated $R_{\text{\upshape HF}}(r,\phistar)$ in \eqref{eqJZHypCond}. We prove \cref{RMASUH} in \cref{RMASUHwithProof,Corollary:EuclideanExpansionMhat} and counterparts for other global sections in \cref{thm:MroutreturnUniformExpansion,,thm:UniformExpansionOnMhat1}) by comparing the lemon billiard expansion with the 1-petal billiard expansion.

\subsection[1-petal Billiard]{1-petal billiard}\label{SBS1petalbilliardorbits}\hfill

\begin{figure}[ht]
\strut\hfill
    \begin{tikzpicture}[xscale=0.65,yscale=.65]
    \tkzDefPoint(0,0){Or};
    \tkzDefPoint(0,-4.8){Y};
    \pgfmathsetmacro{\rradius}{4.5};
    \pgfmathsetmacro{\bdist}{6.3};
    \clip(-1.0*\rradius-0.2,-1.4*\rradius) rectangle (\rradius*1.0+0.2, \rradius*1.5);
    \tkzDefPoint(0,\bdist){OR};
    \pgfmathsetmacro{\phistardeg}{18.1};
    \pgfmathsetmacro{\XValueArc}{\rradius*sin(\phistardeg)};
    \pgfmathsetmacro{\YValueArc}{\rradius*cos(\phistardeg)};    
    \tkzDefPoint(\XValueArc,-1.0*\YValueArc){B};    \tkzDefPoint(-1.0*\XValueArc,-1.0*\YValueArc){A};
    \pgfmathsetmacro{\Rradius}{veclen(\XValueArc,\YValueArc+\bdist)};
    \pgfmathsetmacro{\XValueArcBOne}{\rradius*sin(3*\phistardeg)};
    \pgfmathsetmacro{\YValueArcBOne}{\rradius*cos(3*\phistardeg)};
    \pgfmathsetmacro{\XValueArcAOne}{-\rradius*sin(3*\phistardeg)};
    \pgfmathsetmacro{\YValueArcAOne}{\rradius*cos(3*\phistardeg)};
    
    \tkzDefPoint(\XValueArcBOne,-1.0*\YValueArcBOne){B1};
    \tkzDefPoint(-\XValueArcBOne,-1.0*\YValueArcBOne){A1};
    \tkzDrawArc[name path=Cr,very thick](Or,B)(A);
    \draw[very thick](A)--(B);
    \draw[very thin,blue,arrows={-Triangle[length=0.2cm]}](B) -- (B1);
    \draw[very thin,blue,arrows={-Triangle[length=0.2cm]}](B1) -- (B);
    \draw[very thin,blue,arrows={-Triangle[length=0.2cm]}](A) -- (A1);
    \draw[very thin,blue,arrows={-Triangle[length=0.2cm]}](A1) -- (A);
    \tkzLabelPoint[below,rotate=0](A){$A$};
    \tkzLabelPoint[left,rotate=0](A1){$A_1$};
    \tkzLabelPoint[below,rotate=0](B){$B$};
    \tkzLabelPoint[right,rotate=0](B1){$B_1$};
    \tkzLabelPoint[above,rotate=0](Or){$O_r$};
    \tkzDrawPoints(A,B,Or);
    \draw[dashed,ultra thin](Or)--+(0,-\rradius)coordinate(negYaxis);
    \draw[dashed,ultra thin](Or)--(A);
    \draw[dashed,ultra thin](Or)--(B);
    \tkzMarkAngle[arc=lll,ultra thin,size=1.1,mark=|](A,Or,Y);
    \tkzMarkAngle[arc=lll,ultra thin,size=1.2,mark=|](Y,Or,B);
    \node[rotate=0] at ($(Or)+(-80:1.8)$)  {\small $\phistar$};
    \node[rotate=0] at ($(Or)+(-0.6*\rradius,0.93*\rradius)$)  {\Large $\Gamma_r$};
    \node[rotate=0] at ($(Or)+(0,-1.18*\rradius)$)  {\Large $\PET(r,\phistar)$};
    \end{tikzpicture}\hfill
    \begin{tikzpicture}
        [xscale=0.65,yscale=.65]
    \tkzDefPoint(0,0){Or};
    \tkzDefPoint(0,-4.8){Y};
    \pgfmathsetmacro{\rradius}{4.5};
    \pgfmathsetmacro{\bdist}{6.3};
    \clip(-1.0*\rradius-0.2,-1.4*\rradius) rectangle (\rradius*1.0+0.2, \rradius*1.5);
    \tkzDefPoint(0,\bdist){OR};
    \pgfmathsetmacro{\phistardeg}{18.1};
    \pgfmathsetmacro{\XValueArc}{\rradius*sin(\phistardeg)};
    \pgfmathsetmacro{\YValueArc}{\rradius*cos(\phistardeg)};    
    \tkzDefPoint(\XValueArc,-1.0*\YValueArc){B};    \tkzDefPoint(-1.0*\XValueArc,-1.0*\YValueArc){A};
    \pgfmathsetmacro{\Rradius}{veclen(\XValueArc,\YValueArc+\bdist)};
    \pgfmathsetmacro{\XValueArcBOne}{\rradius*sin(3*\phistardeg)};
    \pgfmathsetmacro{\YValueArcBOne}{\rradius*cos(3*\phistardeg)};
    \pgfmathsetmacro{\XValueArcAOne}{-\rradius*sin(3*\phistardeg)};
    \pgfmathsetmacro{\YValueArcAOne}{\rradius*cos(3*\phistardeg)};
    
    \tkzDefPoint(\XValueArcBOne,-1.0*\YValueArcBOne){B1};
    \tkzDefPoint(-\XValueArcBOne,-1.0*\YValueArcBOne){A1};
    \tkzDrawArc[name path=Cr,very thick](Or,B)(A);
    \tkzDrawArc[name path=CR,very thick](OR,A)(B);
    \draw[very thin,blue,arrows={-Triangle[length=0.2cm]}](B) -- (B1);
    \draw[very thin,blue,arrows={-Triangle[length=0.2cm]}](B1) -- (B);
    \draw[very thin,blue,arrows={-Triangle[length=0.2cm]}](A) -- (A1);
    \draw[very thin,blue,arrows={-Triangle[length=0.2cm]}](A1) -- (A);
    \tkzLabelPoint[below,rotate=0](A){$A$};
    \tkzLabelPoint[left,rotate=0](A1){$A_1$};
    \tkzLabelPoint[below,rotate=0](B){$B$};
    \tkzLabelPoint[right,rotate=0](B1){$B_1$};
    \tkzLabelPoint[above,rotate=0](Or){$O_r$};
    \tkzLabelPoint[left,rotate=0](OR){$O_R$};
    \tkzDrawPoints(A,B,Or,OR);
    \draw[dashed,ultra thin](Or)--+(0,-\rradius)coordinate(negYaxis);
    \draw[dashed,ultra thin](Or)--(A);
    \draw[dashed,ultra thin](Or)--(B);
    \draw[dashed,ultra thin](OR)--(A);
    \draw[dashed,ultra thin](OR)--(B);
    \tkzMarkAngle[arc=lll,ultra thin,size=1.1,mark=|](A,Or,Y);
    \tkzMarkAngle[arc=lll,ultra thin,size=1.2,mark=|](Y,Or,B);
    \node[rotate=0] at ($(Or)+(-80:1.8)$)  {\small $\phistar$};
    \node[rotate=0] at ($(Or)+(0.6*\rradius,0.93*\rradius)$)  {\Large $\Gamma_r$};
    \node[rotate=0] at ($(Or)+(0,-1.18*\rradius)$)  {\Large $\LEM(r,R,\phistar)$};
    \end{tikzpicture}
    \hfill\strut
    \caption{A 1-petal billiard table $\PET(r,\phistar)$ and a lemon billiard table $\LEM(r,R,\phistar)$ with same parameters $r$, $\phistar$.}\label{fig:1-petalExtensionFromCornerAndLemon}
\end{figure}

In this subsection, we note that the exceptional $\phistar\in\big(0,\tan^{-1}{(1/3)}\big)$ in \cref{thm:MroutreturnUniformExpansion,RMASUH} are $\pi/n$ for $n\in\mathbb{N}$. This is related to two comparison billiards for the lemon billiard as follows. In the lemon billiard, the boundary arc that connects $A$ and $B$ in \cref{fig:1-petalExtensionFromCornerAndLemon} has curvature $1/R$. If the curvature is zero instead, then we obtain the 1-petal billiard on the left of \cref{fig:1-petalExtensionFromCornerAndLemon}, and if the curvature is $1/r$, then the billiard is perfectly circular. The exceptional $\phistar$s are exactly those for which $A,B$ are adjacent points of a periodic orbit in this circular billiard. 

Subsequent subsections make quantitative comparisons between the two billiards in \cref{fig:1-petalExtensionFromCornerAndLemon}, and we now explicitly adopt our standing notations for work on the 1-petal billiard.
\begin{notation}[Lemon billiard notations used for 1-petal billiard]\label{def:LemonAndPetTogether}\hfill
    \begin{enumerate}
        \item\label{item01LemonAndPetTogether} We denote by $\PET(r,\phistar)$ a 1-petal billiard table with corners $A$ and $B$, a flat boundary $\overline{AB}$ and a boundary $\Gamma_r$ that is an arc of circle $C_r$ centered on $O_r$, and $\measuredangle{AO_rB}=2\phistar<2\tan^{-1}\big(1/3\big)$ shown in \cref{fig:WojtkowskiBilliardTable,fig:1-petalExtensionFromCornerAndLemon}. We denote by $\bar F_{\phistar}$  the 1-petal billiard map of $\PET(r,\phistar)$.
        
        \item\label{item02LemonAndPetTogether} The 1-petal billiard collisions on $\Gamma_r$ are represented by $(\phi,\theta)\in M_r=[\phistar,2\pi-\phistar]\times[0,\pi]$, which is in \cref{fig:Mr} the same as the lemon billiard phase space $M_r$ in \cref{def:BasicNotations0}, with $\phi$ being the arc angle parameter of $C_r$ for the collision position and $\theta$ being the collision angle with the tangent direction of $\Gamma_r$. The collisions on $\overline{AB}$ are represented by $(s,\theta_1)\in M_{\text{\upshape{f}}}\dfn[-r\sin\phistar,r\sin\phistar]\times[0,\pi]$ with $s$ being the length parameter on $\overline{AB}$ and $\theta_1$ being the collision angle with tangent direction. The billiard map $\bar F_{\phistar}$ is a map from phase space $M_r\sqcup M_{\text{\upshape{f}}}\nfd\overline{M}$ to itself.
        
        \item\label{item03LemonAndPetTogether} For $\overline{x}\in M_{\text{\upshape{f}}}$, we define $p(\overline{x})$ as the collision position in $\PET(r,\phistar)$ in \cref{fig:1-petalExtensionFromCornerAndLemon}. 
        
        \item\label{item04LemonAndPetTogether} The $\Mrout$, $\Mrin$ in \cref{def:BasicNotations} defined for a $\LEM(r,R,\phistar)$ billiard map $\Fc$ have the similar meaning for the $\PET(r,\phistar)$ billiard map $\bar F_{\phistar}$ in \cref{fig:1-petalExtensionFromCornerAndLemon}. That is, they are the collisions that occur on the boundary $\Gamma_r$ with the next(previous) collision on another boundary. i.e. \begin{align*}
            \Int M_r\cap\Fc^{-1}(\Int M_R)\overbracket{=}^{\mathclap{\text{\cref{def:BasicNotations}}}}\Mrout&=\Int M_r\cap\bar F_{\phistar}^{-1}(\Int M_{\text{\upshape{f}}}),\\
            \Int M_r\cap\Fc(\Int M_R)\overbracket{=}^{\mathclap{\text{\cref{def:BasicNotations}}}}\Mrin&=\Int M_r\cap\bar F_{\phistar}(\Int M_{\text{\upshape{f}}}), 
        \end{align*}
        
        For $x\in \Int M_r$, $p(\Fcal(x))\in \Int\Gamma_r$ is equivalent to $p(\bar F_{\phistar}(x))\in\Int\Gamma_r$. Therefore for $x\in\Mrin$, \[
        \max{\big\{m\ge0\bigm|p(\bar F_{\phistar}^i(x))\in\operatorname{int}\Gamma_r\text{ for }0\le i\le m\big\}}=\max{\big\{m\ge0\bigm|p(\Fc^i(x))\in\operatorname{int}\Gamma_r\text{ for }0\le i\le m\big\}}\overbracket{=}^{\mathclap{\text{\cref{def:SectionSetsDefs}}}}m(x)
        \]
        
        \item\label{item05LemonAndPetTogether} For $\overline{x}_0=(\overline{\phi}_0,\overline{\theta}_0)\in\Mrout$ with $\overline{x}_1=\bar F_{\phistar}(\overline{x}_0)\in M_{\text{\upshape f}}$ and $(\overline{\phi}_2,\overline{\theta}_2)=\overline{x}_2=\bar F^2_{\phistar}(\overline{x}_0)\in\Mrin$ define $\overline{d}_0\dfn r\sin{\overline{\theta}_0}$, $\overline{d}_2\dfn r\sin{\overline{\theta}_2}$ and let $\overline{\tau}_0\dfn|p(\overline{x}_0)p(\overline{x}_1)|$, $\overline{\tau}_1\dfn|p(\overline{x}_1)p(\overline{x}_2)|$ be the Euclidean distance between the two consecutive collisions on $\PET(r,\phistar)$ billiard table in \cref{fig:1-petalExtensionFromCornerAndLemon}.
        
        \item\label{item06LemonAndPetTogether} With the same $\Min_{r,n}$, $\Mout_{r,n}$ and global section $\hat{M}= \big(\bigsqcup_{n\ge1}\Fc^{n-1}\Min_{r,n}\big)\sqcup\Min_{r,0}=\Fc^{-1}\big(\Mrout\smallsetminus{\Mrin}\big)\sqcup\big(\Mrout\cap\Mrin\big)$ in \cref{def:SectionSetsDefs} and \eqref{def:Mhat}, we define $\hat{\bar F}_{\phistar}$ the return map on $\hat{M}$ in $\PET(r,\phistar)$ table under $\bar F_{\phistar}$ dynamics. And we define $\tilde{\bar F}_{\phistar}$ the return map on $\Mrout$ in $\PET(r,\phistar)$ table under $\bar F_{\phistar}$ dynamics.
        
        
\item\label{item07LemonAndPetTogether} In \cref{fig:1-petalExtensionFromCornerAndLemon} both lemon billiard $\mathbb{L}(r,R,\phistar)$ and 1-petal billiard $\mathbb{P}(r,\phistar)$, if we consider $\overrightarrow{BA}$ as reflecting in $\Gamma_r$ at $A$  then $\overrightarrow{AA_1}$ is its next collision direction at $A_1\in\Gamma_r$. 
$\overrightarrow{BA}$ represents the boundary point $(2\pi-\phistar,\pi-\phistar)=Ix_*$ of $M_r$ in \cref{fig:PetalMroutMfMrinExtension1}, so 
$(2\pi-3\phistar,\pi-\phistar)\in M_r$ is the extended definition for the collision at $A_1$ as the next collision of $\overrightarrow{BA}$ in \cref{fig:1-petalExtensionFromCornerAndLemon}.\COMMENT{\large\bfseries Edit next paragraph carefully in the same way}

Symmetrically, if we consider of $\overrightarrow{AB}$ as a collision on the corner $B$ with direction $\overrightarrow{AB}$ which singular  boundary point of $M_r$, that is, $(\phistar,\phistar)=Iy_*$ in \cref{fig:Mr}, then $\overrightarrow{BB_1}$ is the next collision at $B_1\in\Gamma_r$. Therefore, in $M_r$ of \cref{fig:Mr}, $(3\phistar,\phistar)\in M_r$ in \cref{fig:Mr} is the collision at $B_1$ as the extended definition for the next collision of $\overrightarrow{AB}$ in \cref{fig:1-petalExtensionFromCornerAndLemon}.
        
    \end{enumerate}
\end{notation}

\begin{figure}[h]
\begin{minipage}{.6\textwidth}
    \begin{center}
    \begin{tikzpicture}[xscale=0.65,yscale=0.65]
        \pgfmathsetmacro{\PHISTAR}{0.8}
        \pgfmathsetmacro{\WIDTH}{7}
        \pgfmathsetmacro{\LENGTH}{7}
        \tkzDefPoint(-\LENGTH,0.5*\WIDTH){LLM};
        \tkzDefPoint(-\LENGTH,0){LLL};
        \tkzDefPoint(-\LENGTH+\PHISTAR,0){LL};
        \tkzDefPoint(-\LENGTH+\PHISTAR,\WIDTH){LU};
        \tkzDefPoint(-\LENGTH,\WIDTH){LLU};
        \tkzDefPoint(\LENGTH-\PHISTAR,0){RL};
        \tkzDefPoint(\LENGTH-\PHISTAR,\WIDTH){RU};
        \tkzDefPoint(\LENGTH,0){RRL};
        \tkzDefPoint(\LENGTH,\WIDTH){RRU};
        \tkzDefPoint(-\LENGTH+\PHISTAR,\PHISTAR){X};
        \tkzDefPoint(-\LENGTH+\PHISTAR,\WIDTH-\PHISTAR){IX};
        \tkzDefPoint(\LENGTH-\PHISTAR,\PHISTAR){IY};
        \tkzDefPoint(\LENGTH-\PHISTAR,\WIDTH-\PHISTAR){Y};

        \tkzDefPoint(-\LENGTH+3*\PHISTAR,2*\PHISTAR){C1};
        \tkzDefPoint(-\LENGTH+3*\PHISTAR,\WIDTH-2*\PHISTAR){C2};
        \tkzDefPoint(\LENGTH-3*\PHISTAR,\WIDTH-2*\PHISTAR){C3};
        \tkzDefPoint(\LENGTH-3*\PHISTAR,2*\PHISTAR){C4};
        
        \tkzLabelPoint[left](LLL){0};
        \tkzLabelPoint[left](LLM){$\frac{\pi}{2}$};
        \tkzLabelPoint[left](LLU){$\pi$};
        \tkzLabelPoint[below](LLL){0};
        \tkzLabelPoint[below](RRL){\ $2\pi$};
        \tkzLabelPoint[right](IY){$x_*=E_4$};
        \tkzLabelPoint[left](IX){$E_2=y_*$};
        \tkzLabelPoint[right](Y){\small $I x_*=E_3$};
        \tkzLabelPoint[left](X){\small $E_1=Iy_*$};
        
        \tkzLabelPoint[above left](C1){\small $C_1=(3\phistar,2\phistar)$};
        \tkzLabelPoint[below left](C2){\small $C_2=(3\phistar,\pi-2\phistar)$};
        \tkzLabelPoint[below](C3){\small $C_3=(2\pi-3\phistar,\pi-2\phistar)$};
        \tkzLabelPoint[above](C4){\small $C_4=(2\pi-3\phistar,2\phistar)$};
        \tkzLabelPoint[below right](LL){\small $D_1=(\phistar,0)$};
        \tkzLabelPoint[above right](LU){\small $D_2=(\phistar,\pi)$};
        \tkzLabelPoint[above left](RU){\small $D_3=(2\pi-\phistar,\pi)$};
        \tkzLabelPoint[below left](RL){\small $D_4=(2\pi-\phistar,0)$};

        \tkzLabelSegment[above right](C4,IY){\small $1$};
        \tkzLabelSegment[right](IY,RL){\small $2$};
        \tkzLabelSegment[below left](RL,C2){\small $3$};
        \tkzLabelSegment[below left](C2,IX){\small $4$};
        \tkzLabelSegment[left](IX,LU){\small $5$};
        \tkzLabelSegment[above right](LU,C4){\small $6$};

        \tkzLabelSegment[left](LL,X){\small $1$};
        \tkzLabelSegment[above left](X,C1){\small $2$};
        \tkzLabelSegment[above left](C1,RU){\small $3$};
        \tkzLabelSegment[right](RU,Y){\small $4$};
        \tkzLabelSegment[below right](Y,C3){\small $5$};
        \tkzLabelSegment[below right](C3,LL){\small $6$};


        \draw [thick,-triangle 45, color=purple](C4)--(IY);
        \draw [thick,-triangle 45, color=orange](IY)--(RL);
        \draw [thick,-triangle 45, color=green](RL)--(C2);
        \draw [thick,-triangle 45, color=red](C2)--(IX);
        \draw [thick,-triangle 45, color=brown](IX)--(LU);
        \draw [thick,-triangle 45, color=yellow](LU)--(C4);

        \draw [thick,-triangle 45, color=purple](LL)--(X);
        \draw [thick,-triangle 45, color=orange](X)--(C1);
        \draw [thick,-triangle 45, color=green](C1)--(RU);
        \draw [thick,-triangle 45, color=red](RU)--(Y);
        \draw [thick,-triangle 45, color=brown](Y)--(C3);
        \draw [thick,-triangle 45, color=yellow](C3)--(LL);
        
        \draw [ultra thin,dashed] (LL) --(LU);
        \draw [ultra thin,dashed] (LU) --(RU);
        \draw [ultra thin,dashed] (RL) --(RU);
        \draw [ultra thin,dashed] (LL) --(RL);
        \draw [ultra thin,dashed](LLL)--(LL);
        \draw [ultra thin,dashed](LLU)--(LU);
        \draw [ultra thin,dashed](LLU)--(LLL);
        \draw [ultra thin,dashed](RRL)--(RL);
        \draw [ultra thin,dashed](RRU)--(RU);
        \draw [ultra thin,dashed](RRU)--(RRL);
        \fill [blue, opacity=10/30](LL) -- (Y) -- (RU) -- (X) -- cycle;
        \fill [red, opacity=10/30](LU) -- (IY) -- (RL) -- (IX) -- cycle;
        
        \coordinate (little_x) at ($(RL)!0.65!(IY)$);
        \begin{scope}
            \clip (RL) -- (IY) -- (LU) -- cycle;
            \draw[fill=orange] circle[at=(little_x),radius=0.12];
        \end{scope}
        \coordinate (F1little_x) at ($(C1)!0.35!(X)$);
        \begin{scope}
            \clip (LL) -- (X) -- (RU) -- cycle;
            \draw[fill=orange] ellipse[rotate=27, at=(F1little_x),x radius=0.201, y radius =0.102];
        \end{scope}
        \tkzLabelPoint[left](little_x){\small $U_1$};
        \tkzLabelPoint[below ](F1little_x){\small $\bar F^2_{\phistar}(U_1)$};
        
        \coordinate (little_x_2) at ($(C4)!0.2!(LU)$);
        \begin{scope}
            \clip (IX) -- (LU) -- (IY) -- (RL) -- cycle;
            \draw[fill=yellow] ellipse[rotate=0, at=(little_x_2),radius=0.12];
        \end{scope}

        \coordinate (F1little_x_2) at ($(C3)!0.68!(LL)$);
        \begin{scope}
            \clip (LL) -- (X) -- (RU) -- (Y) -- cycle;
            \draw[fill=yellow] ellipse[rotate=30, at=(F1little_x_2),x radius=0.36, y radius =0.09];
        \end{scope}
        \tkzLabelPoint[above](little_x_2){\small $U_2$};
        \tkzLabelPoint[below](F1little_x_2){\small $\bar F^2_{\phistar}(U_2)$};

        \coordinate (little_x_3) at ($(C2)!0.9!(RL)$);
        \begin{scope}
            \clip (IX) -- (LU) -- (IY) -- (RL) -- cycle;
            \draw[fill=green] ellipse[rotate=0, at=(little_x_3),radius=0.12];
        \end{scope}

        \coordinate (F1little_x_3) at ($(C1)!0.13!(RU)$);
        \begin{scope}
            \clip (LL) -- (X) -- (RU) -- (Y) -- cycle;
            \draw[fill=green] ellipse[rotate=15, at=(F1little_x_3),x radius=0.23, y radius =0.06];
        \end{scope}
        \tkzLabelPoint[left](little_x_3){\small $U_3$};
        \tkzLabelPoint[above](F1little_x_3){\small $\bar F^2_{\phistar}(U_3)$};

        \tkzDrawPoints(X,IX,Y,IY,RL,RU,LL,LU,C1,C2,C3,C4);

    \end{tikzpicture}
    \end{center}
    In the figure, it is shown \caption{$\bar F_{\phistar}$ restricted to $\Mrout$ and with continuous extension to $\partial\Mrout$.\newline closure($\Mrin$) is the image of $\bar F_{\phistar}$ restricted to $M_{\text{\upshape f}}$ with continuous extension to $\partial M_{\text{\upshape f}}$.}\label{fig:PetalMroutMfMrinExtension1}    
\end{minipage}\hfill
\begin{minipage}{0.4\textwidth}
\begin{center}
    \begin{tikzpicture}[scale=.39]
        \pgfmathsetmacro{\PHISTAR}{2.5};
        \pgfmathsetmacro{\WIDTH}{15};
        \pgfmathsetmacro{\LENGTH}{7.5};
        \tkzDefPoint(-1.06*\LENGTH,0.5*\WIDTH){LLM};
        \tkzDefPoint(-\LENGTH,0){LLL};
        \tkzDefPoint(-\PHISTAR,0){LL};
        \tkzDefPoint(-\PHISTAR,\WIDTH){LU};
        \tkzDefPoint(-\LENGTH,\WIDTH){LLU};
        \tkzDefPoint(\PHISTAR,0){RL};
        \tkzDefPoint(\PHISTAR,\WIDTH){RU};
        \tkzDefPoint(\LENGTH,0){RRL};
        \tkzDefPoint(\LENGTH,\WIDTH){RRU};
        \tkzDefPoint(-\PHISTAR,\PHISTAR){X};
        \tkzDefPoint(-\PHISTAR,\WIDTH-\PHISTAR){IX};
        \tkzDefPoint(\PHISTAR,\WIDTH-\PHISTAR){Y};
        \tkzDefPoint(\PHISTAR,\PHISTAR){IY};
        \draw [thick] (LU) --(RU);
        \draw [thick] (LL) --(RL);

        \draw [thick,-triangle 45, color=orange](LL)--(X);
        \draw [thick,-triangle 45, color=green](X)--(IX);
        \draw [thick,-triangle 45, color=red](IX)--(LU);
        \draw [thick,-triangle 45, color=brown](RU)--(Y);
        \draw [thick,-triangle 45, color=yellow](Y)--(IY);
        \draw [thick,-triangle 45, color=purple](IY)--(RL);

        \tkzLabelPoint[left](X){$(-r\sin{\phistar},\phistar)$};
        \tkzLabelPoint[left](IX){$(-r\sin{\phistar},\pi-\phistar)$};
        
        \tkzLabelPoint[below](LL){\small $-r\sin{\phistar}$};
        \tkzLabelPoint[left](LL){$0$};
        \tkzLabelPoint[below left](RL){\small $+r\sin{\phistar}$};
        \tkzLabelPoint[above](LU){$\pi$};
        \tkzDrawPoints(LL,LU,X,Y,IX,IY,RU,RL)
        \tkzLabelPoint[above](LLM){\small closure$(\Mrin)\xleftarrow{\bar F_{\phistar,1}} M_{\text{\upshape f}}$};
        \tkzLabelPoint[below](LLM){\small closure$(\Mrout)\smallsetminus\{x_*,y_*\}\xrightarrow{\bar F_{\phistar,1}}M_{\text{\upshape f}}$};

        \coordinate (middle) at ($(RU)!0.5!(LL)$);
        \tkzLabelPoint[above](middle){$M_{\text{\upshape f}}$};
        
        \coordinate (F1little_x_1) at ($(LL)!0.3!(X)$);
        \coordinate (F1little_x_3) at ($(X)!0.12!(IX)$);
        \coordinate (F1little_x_2) at ($(Y)!0.7!(IY)$);

        \tkzLabelPoint[right](F1little_x_1){\small $\bar F_{\phistar}(U_1)$};
        \tkzLabelPoint[below left](F1little_x_2){\small $\bar F_{\phistar}(U_2)$};
        \tkzLabelPoint[right](F1little_x_3){\small $\bar F_{\phistar}(U_3)$};

        \begin{scope}
            \clip (LL) -- (LU) -- (RU) -- (RL) -- cycle;
            \draw[fill=orange] ellipse[rotate=27, at=(F1little_x_1),x radius=0.421, y radius =0.212];
            \draw[fill=green] ellipse[rotate=27, at=(F1little_x_3),x radius=0.501, y radius =0.272];
            \draw[fill=yellow] ellipse[rotate=50, at=(F1little_x_2),x radius=0.601, y radius =0.322];
        \end{scope}

        \tkzLabelSegment[left](LL,X){2};
        \tkzLabelSegment[left](X,IX){3};
        \tkzLabelSegment[left](IX,LU){4};
        \tkzLabelSegment[left](RU,Y){5};
        \tkzLabelSegment[left](Y,IY){6};
        \tkzLabelSegment[left](IY,RL){1};

    \end{tikzpicture}\caption{$\bar F_{\phistar}$ restricted to $M_{\text{\upshape f}}$ with continuous extension to $\partial M_{\text{\upshape f}}$.}\label{fig:PetalMroutMfMrinExtension2}
\end{center}
\end{minipage}
\end{figure}
\begin{remark}[1-petal billiard map continuous extension]\label{remark:petalMroutMfMinTransition}\hfill
    \begin{enumerate}
        \item\label{item01petalMroutMfMinTransition} \cref{fig:PetalMroutMfMrinExtension1} shows that $\partial \Mrout$ is mapped to $\partial M_{\text{\upshape f}}$ by the continuous extension of $\bar F_{\phistar}$'s restriction on $\Mrout$ with same color and number labeled straight-line segments being the image and preimage under this continuous extended map.
        \item\label{item02petalMroutMfMinTransition}\cref{fig:PetalMroutMfMrinExtension2} shows that $\partial M_{\text{\upshape f}}$ is mapped to $\partial \Mrin$ by the continuous extension of $\bar F_{\phistar}$'s restriction on $\partial M_{\text{\upshape f}}$ with same color and number labeled straight-line segments being the image and preimage under this continuous extended map.
        \item\label{item03petalMroutMfMinTransition} In \cref{fig:PetalMroutMfMrinExtension1}, line segment 1: $\overline{C_4E_4}\subset\partial\Mrout$ is $\bar F^2_{\phistar,1}$-diffeomorphic to line segment 1: $\overline{D_1E_1}\subset\partial\Mrin$.\newline Line segment 2: $\overline{E_4D_4}\subset\partial\Mrout$ is $\bar F^2_{\phistar,1}$-diffeomorphic to line segment 2: $\overline{E_1C_1}\subset\partial\Mrin$.\newline Line segment 3: $D_4C_2\subset\partial\Mrout$ is $\bar F^2_{\phistar,1}$-diffeomorphic to line segment 3: $\overline{C_1D_3}\subset\partial\Mrin$.\newline Line segment 4: $\overline{C_2E_2}\subset\partial\Mrout$ is $\bar F^2_{\phistar,1}$-diffeomorphic to line segment 4: $\overline{D_3E_3}\subset\partial\Mrin$.\newline Line segment 5: $\overline{E_2D_2}\subset\partial\Mrout$ is $\bar F^2_{\phistar,1}$-diffeomorphic to $\text{line segment } 5: \overline{E_3C_3}\subset\partial\Mrin$.\newline Line segment 6: $\overline{D_2C_4}\subset\partial\Mrout$ is $\bar F^2_{\phistar,1}$-diffeomorphic to line segment 6: $\overline{C_3D_1}\subset\partial\Mrin$.
    \end{enumerate}
\end{remark}

\begin{figure}[h!]
\begin{minipage}{.45\textwidth}
\begin{center}
    \begin{tikzpicture}[xscale=.4,yscale=.4]
        \pgfmathsetmacro{\PHISTAR}{0.8}
        \pgfmathsetmacro{\WIDTH}{7}
        \pgfmathsetmacro{\LENGTH}{7}
        \tkzDefPoint(-\LENGTH,0.5*\WIDTH){LLM};
        \tkzDefPoint(-\LENGTH,0){LLL};
        \tkzDefPoint(-\LENGTH+\PHISTAR,0){LL};
        \tkzDefPoint(-\LENGTH+\PHISTAR,\WIDTH){LU};
        \tkzDefPoint(-\LENGTH,\WIDTH){LLU};
        \tkzDefPoint(\LENGTH-\PHISTAR,0){RL};
        \tkzDefPoint(\LENGTH-\PHISTAR,\WIDTH){RU};
        \tkzDefPoint(\LENGTH,0){RRL};
        \tkzDefPoint(\LENGTH,\WIDTH){RRU};
        \tkzDefPoint(-\LENGTH+\PHISTAR,\PHISTAR){X};
        \tkzDefPoint(-\LENGTH+\PHISTAR,\WIDTH-\PHISTAR){IX};
        \tkzDefPoint(\LENGTH-\PHISTAR,\PHISTAR){IY};
        \tkzDefPoint(\LENGTH-\PHISTAR,\WIDTH-\PHISTAR){Y};
        \tkzLabelPoint[below](LL){$(\phistar,0)$};
        \tkzLabelPoint[below](RL){$(2\pi-\phistar,0)\qquad$};
        \tkzLabelPoint[above](LU){$(\phistar,\pi)$};
        \tkzLabelPoint[above](RU){$(2\pi-\phistar,\pi)$};
        \tkzLabelPoint[right](IY){$x_*$};
        \tkzLabelPoint[left](IX){$y_*$};
        \node[] at (0.15*\LENGTH,0.8*\WIDTH)  {\small $M^{1}_{r-out}\subset M_r\smallsetminus\Mrout$};
        \node[] at (-0.15*\LENGTH,0.2*\WIDTH)  {\small $M^{1}_{r-out}\subset M_r\smallsetminus\Mrout$};

        \draw [thick,green, -triangle 45] (IX) --(LL);
        \draw [thick,blue, -triangle 45] (LL) --(RL);
        \draw [thick,red, -triangle 45] (RL) --(IX);
        \draw [thin,dashed] (IX) --(LU);
        \draw [thin,dashed] (IY) --(RL);
        \draw [thick,Purple, -triangle 45] (RU) --(LU);
        \draw [thick,teal, -triangle 45] (LU) --(IY);
        \draw [thick,yellow, -triangle 45] (IY) --(RU);

        \fill [orange, opacity=3/30](LL) -- (RL) -- (IX) -- cycle;
        \fill [orange, opacity=3/30](LU) -- (IY) -- (RU) -- cycle;
    \end{tikzpicture}
    \end{center}
\end{minipage}\hfill
\begin{minipage}{.45\textwidth}
\begin{center}
    \begin{tikzpicture}[xscale=.4,yscale=.4]
        \pgfmathsetmacro{\PHISTAR}{0.8}
        \pgfmathsetmacro{\WIDTH}{7}
        \pgfmathsetmacro{\LENGTH}{7}
        \tkzDefPoint(-\LENGTH,0.5*\WIDTH){LLM};
        \tkzDefPoint(-\LENGTH,0){LLL};
        \tkzDefPoint(-\LENGTH+\PHISTAR,0){LL};
        \tkzDefPoint(-\LENGTH+\PHISTAR,\WIDTH){LU};
        \tkzDefPoint(-\LENGTH,\WIDTH){LLU};
        \tkzDefPoint(\LENGTH-\PHISTAR,0){RL};
        \tkzDefPoint(\LENGTH-\PHISTAR,\WIDTH){RU};
        \tkzDefPoint(\LENGTH,0){RRL};
        \tkzDefPoint(\LENGTH,\WIDTH){RRU};
        \tkzDefPoint(-\LENGTH+\PHISTAR,\PHISTAR){X};
        \tkzDefPoint(-\LENGTH+\PHISTAR,\WIDTH-\PHISTAR){IX};
        \tkzDefPoint(\LENGTH-\PHISTAR,\PHISTAR){IY};
        \tkzDefPoint(\LENGTH-\PHISTAR,\WIDTH-\PHISTAR){Y};
        \tkzLabelPoint[below](LL){$(\phistar,0)$};
        \tkzLabelPoint[below left](RL){$(2\pi-\phistar,0)$};
        \tkzLabelPoint[above](LLM){$\xrightarrow{\Fc=\begin{pmatrix}1 &2\\0 &1\end{pmatrix}}\qquad\quad$};
        \tkzLabelPoint[above](LU){$(\phistar,\pi)$};
        \tkzLabelPoint[above left](RU){$(2\pi-\phistar,\pi)$};
        \tkzLabelPoint[right](Y){$I x_*$};
        \tkzLabelPoint[left](X){$I y_*$};

        \draw [thick,green, -triangle 45] (Y) --(LL);
        \draw [thick,blue, -triangle 45] (LL) --(RL);
        \draw [thick,red, -triangle 45] (RL) --(Y);
        \draw [thin,dashed] (Y) --(RU);
        \draw [thin,dashed] (X) --(LL);
        \draw [thick,Purple, -triangle 45] (RU) --(LU);
        \draw [thick,teal, -triangle 45] (LU) --(X);
        \draw [thick,yellow, -triangle 45] (X) --(RU);
        \fill [orange, opacity=3/30](X) -- (LU) -- (RU) -- cycle;
        \fill [orange, opacity=3/30](LL) -- (RL) -- (Y) -- cycle;
        
    \end{tikzpicture}
    \end{center}
\end{minipage}\caption{$\bar F_{\phistar}$  with continuous extension to $\partial M_r\smallsetminus\Mrout$ by its definition restricted to $M_r\smallsetminus\Mrout$}\label{fig:MrminusMroutExtension}
\end{figure}

\begin{notation}[Continuous extension of $\bar F_{\phistar}$]\label{notation:PetContinuousExtension}\hfill
    \begin{enumerate}
    \item\label{item01PetContinuousExtension} (In \cref{fig:PetalMroutMfMrinExtension1,fig:PetalMroutMfMrinExtension2}) We define $\bar F_{\phistar,1}$ as the 1-petal $\PET(r,\phistar)$ billiard map $\bar F_{\phistar}$'s continuous extension to $\partial M_r$, $\partial M_{\text{\upshape f}}$ by the following continuous extended definitions:\newline For $y\in\partial\Mrout$ and $y\notin\big\{x_*,y_*\big\}$, $\bar F_{\phistar,1}(y)\dfn\lim\limits_{\Mrout\ni z\to y}\bar F_{\phistar}(z)$. Note that if $y\in\big\{x_*,y_*\big\}$, then this limit DNE. \newline For $y=(s,\theta)\in\partial M_{\text{\upshape f}}$, $\bar F_{\phistar,1}(y)\dfn\lim\limits_{\Int{M_{\text{\upshape f}}}\ni z\to y}\bar F_{\phistar}(z)$. Note that if $\theta\in\big\{0,\pi\big\}$, then this limit is $Ix_*$ or $Iy_*$.\newline For $y\in\Int(M_{\text{\upshape f}})\sqcup\Int(M_r)\smallsetminus\partial\Mrout$, $\bar F_{\phistar,1}(y)\dfn\bar F_{\phistar}(y)$.
    \item\label{item02PetContinuousExtension} (In \cref{fig:MrminusMroutExtension}) We define $\bar F_{\phistar,2}$ as the 1-petal $\PET(r,\phistar)$ billiard map $\bar F_{\phistar}$'s continuous extension to $\partial M_r$, $\partial M_{\text{\upshape f}}$ by the following continuous extended definitions:\newline For $y\in\partial\Mrout$, $\bar F_{\phistar,2}(y)\dfn\lim\limits_{(M_r\smallsetminus\Mrout)\ni z\to y}\bar F_{\phistar}(z)$. \newline For $y=(s,\theta)\in\partial M_{\text{\upshape f}}$, $\bar F_{\phistar,2}(y)\dfn\bar F_{\phistar,1}(y)$.\newline For $y\in\Int(M_{\text{\upshape f}})\sqcup\Int(M_r)\smallsetminus\partial\Mrout$, $\bar F_{\phistar,1}(y)\dfn\bar F_{\phistar,2}(y)=\bar F_{\phistar}(y)$.
    \end{enumerate}
\end{notation}
\begin{remark}[Multi-value caused by different continuous extension]\label{remark:1-petalContinuousExtension}\hfill
\begin{itemize}
    \item $\bar F_{\phistar,1}(y)$ and $\bar F_{\phistar,2}(y)$ are the same as defined except in $\partial\Mrout\smallsetminus\partial M_r$. Note that $F_{\phistar,2}$'s restriction on $M_r\smallsetminus\Mrout$ and its continuous extension on $\partial(M_r\smallsetminus\Mrout)$ are represented by $\begin{pmatrix}
        1,&2\\
        0,&1
    \end{pmatrix}$ (see \cref{fig:MrminusMroutExtension}).
    \item  $\overline{\tau}_0$, $\overline{d}_0$ in \cref{def:LemonAndPetTogether}\eqref{item05LemonAndPetTogether} as the length functions of $\overline{x}_0\in\Mrout$ have a continuous extended definition on closure$(\Mrout)\smallsetminus\{x_*,y_*\}$. 
    $\overline{\tau}_1$, $\overline{d}_2$ in \cref{def:LemonAndPetTogether}\eqref{item05LemonAndPetTogether} as the length functions of $\overline{x}_1\in M_{\text{\upshape f}}$ also have a continuous extended definition on the boundary of $M_{\text{\upshape f}}$.
\end{itemize}
\end{remark}

\begin{lemma}[Exponential expansion in a 1-petal billiard orbit segment]\label{lemma:PETExpExpansion}

Suppose $\phistar\in(0,\tan^{-1}(1/3)$, $\pi/\phistar\notin\mathbb{N}$ and $\tilde{\bar F}_{\phistar}^k(x)|^{n}_{k=0}\nfd \tilde{\overline{x}}_k|^n_{k\ge0}\in\Mrout$ is a finite orbit segment of $x\in\Mrout$ under the iteration of the return map $\tilde{\bar F}_{\phistar}$ (\cref{def:LemonAndPetTogether}\eqref{item06LemonAndPetTogether}) on $\Mrout$. If $\tilde{\bar F}_{\phistar}^k(x)|_{k=0,\cdots,n}$ never hits the corners of $\PET(r,\phistar)$,  then with $\overline{\tau}_{0,k}$, $\overline{\tau}_{1,k}$, $\overline{d}_{0,k}$ defined in \cref{def:LemonAndPetTogether}\eqref{item05LemonAndPetTogether} for each $\tilde{\overline{x}}_k$, $\prod_{k=1}^{n}\bigm|-1-\frac{\overline{\tau}_{0,k}+\overline{\tau}_{1,k}-2\overline{d}_{0,k}}{\overline{d}_{0,k}}\bigm|>\overline{C}\cdot\overline{E}^n$ with some constant $\overline{C}>0$ and $\overline{E}>1$. i.e. it has an exponential expansion growth.
\end{lemma}
\begin{proof}
    By the definitions of length functions $\overline{\tau}_0$, $\overline{\tau}_1$, $\overline{d}_0$ and points $\overline{x}_0$, $\overline{x}_1$, $\overline{x}_2$ in \cref{def:LemonAndPetTogether}\eqref{item05LemonAndPetTogether}, we have $\overline{\tau}_{0,k}+\overline{\tau}_{1,k}-2\overline{d}_{0,k}>0$ in \cref{fig:WojtkowskiBilliardTable}, which is the root cause of the defocusing mechanism in \cite[section 8.1,8.2]{cb}.  Suppose $d\overline{x}_0\in T_{\overline{x}_0}M_r$ in the first and third quadrant in $\phi\theta$-coordinates, then $d\overline{x}_2=D_{\overline{x}_0}\bar F^2_{\phistar}(d\overline{x}_0)\in T_{\overline{x}_2}M_r$ satisfies \begin{equation}\label{eq:lowerboundofPetalExpansion}
        \frac{\|d\overline{x}_2\|_{\p}}{\|d\overline{x}_0\|_{\p}}\ge\bigm|-1-\frac{\overline{\tau}_{0}+\overline{\tau}_{0}-2\overline{d}_0}{\overline{d}_0}\bigm|>1.
    \end{equation}

    Our proof will proceed to be based on the observations in \cref{fig:PetalMroutMfMrinExtension1}.
    
    We can see that if $\tilde{\overline{x}}_{k}\in\Mrout$ and is (in Euclidean distance) bounded away from the boundary segments $\overline{C_4E_4}$ and $\overline{C_2E_2}$ in \cref{fig:PetalMroutMfMrinExtension1} (segments 1 and 4 in \cref{fig:PetalMroutMfMrinExtension1}), then $\overline{\tau}_{0,k}+\overline{\tau}_{1,k}-2\overline{d}_{0,k}$ is bounded away from $0$. And the small enough neighborhood (contained in $\Mrout$) of $\overline{C_4E_4}$ ($\overline{C_2E_2}$) is diffeomorphically mapped by $\bar F^2_{\phistar}$ to the small enough neighborhood (contained in $\Mrin$) of $\overline{D_1E_1}$ ($\overline{D_3E_3}$) by \cref{remark:petalMroutMfMinTransition}\eqref{item03petalMroutMfMinTransition}. 
    
    Therefore, if $\bar{F}^2_{\phistar}(\tilde{\overline{x}}_k)\in\Mrin$ is bounded away from $\overline{D_1E_1}\cup\overline{D_3E_3}$, then $\overline{\tau}_{0}+\overline{\tau}_{0}-2\overline{d}_0$ is bounded away from $0$.

    Now on the contrary, we suppose $\bar{F}^2_{\phistar}(\tilde{\overline{x}}_k)=x\in\mathcal{N}\cap\Mrin$ with $\mathcal{N}$ being some neighborhood of $\overline{D_1E_1}$. In \cref{fig:PetalMroutMfMrinExtension1} note that $E_1=(\phistar,\phistar)$ is the endpoint for the line segment $\overline{D_1E_1}=\big\{(\phistar,\theta)\in M_r\bigm|\theta\in[0,\phistar]\big\}$. For $\forall x=(\phi,\theta)\in \overline{D_1E_1}$, let $\kappa(x)=\max{\big\{t\bigm|p(\bar F_{\phistar}^{i}(x))\in\Gamma_r,i=1,\cdots,t\big\}}$. Since $x,\cdots,\bar F_{\phistar}^{\kappa(x)}(x)$ are consecutive collisions on $\Gamma_r$, all the $\theta$ coordinates for $x,\cdots,\bar F_{\phistar}^{\kappa(x)}(x)$ are the same as for the $x$ and in the range $[0,\phistar]$.

    Therefore, if $x\in\mathcal{N}$ and $x$ is (in Euclidean distance) bounded away from $E_1=Iy_*=(\phistar,\phistar)$, then $\bar F_{\phistar}^{\kappa(x)}(x)$ is (in Euclidean distance) bounded away from the line segments $\overline{C_4E_4}$ and $\overline{E_2C_2}$ in \cref{fig:PetalMroutMfMrinExtension1}.

    Since $\pi/\phistar\notin\mathbb{N}$, $\bar F_{\phistar}^{k(Iy_*)}(Iy_*)\in\Mrout$ and is (in Euclidean distance) bounded away from the line segments $\overline{C_4E_4}$ and $\overline{C_2E_2}$. By continuity of $\bar F_{\phistar}$, there exists some $N(Iy_*)$: a neighborhood of $Iy_*$ such that $\bar F_{\phistar}^{\kappa(Iy_*)}(N(Iy_*))\subset\Mrout$ and $\bar F_{\phistar}^{\kappa(Iy_*)}(N(Iy_*))$ is (in Euclidean distance) bounded away from $(2\pi-\phistar,\phistar)$. Thus, $\bar F_{\phistar}^{\kappa(Iy_*)}(N(Iy_*)\cap \mathcal{N})$ is (in Euclidean distance) bounded away from $\overline{C_4C_4}$ (also $\overline{C_2E_2}$) for all .

    Hence, for all $x\in N(Iy_*)$ or $x\in\mathcal{N}\cap N(Iy_*)$, $\bar F_{\phistar}^{\kappa(x)}(x)$ is uniformly (in Euclidean distance) bounded away from line segments $\overline{C_4E_4}$ and $\overline{E_2C_2}$. 
    
    Therefore, there exists some uniform constant $\overline{\delta}>0$ so that for every two consecutive $\tilde{\overline{x}}_k, \tilde{\overline{x}}_{k+1}$, at least one of $\overline{\tau}_{0,k}+\overline{\tau}_{1,k}-2\overline{d}_{0,k}$, $\overline{\tau}_{0,k+1}+\overline{\tau}_{1,k+1}-2\overline{d}_{0,k+1}$ is larger than some constant $\overline{\delta}>0$. On the other hand for all $k$, $\overline{d}_{0,k}\le r$. We get the following.
    \[\big|-1-\frac{\overline{\tau}_{0,k}+\overline{\tau}_{1,k}-2\overline{d}_{0,k}}{\overline{d}_{0,k}}\bigm|\bigm|-1-\frac{\overline{\tau}_{0,k+1}+\overline{\tau}_{1,k+1}-2\overline{d}_{0,k+1}}{\overline{d}_{0,k+1}}\bigm|>\bigm|-1-\frac{\overline{\delta}}{r}\big|.\]
If we set $\overline{E}=(1+\frac{\overline{\delta}}{r})^{0.5}>1$ and $\overline{C}=(1+\frac{\overline{\delta}}{r})^{-1}$, then $\prod_{k=1}^{n}\bigm|-1-\frac{\overline{\tau}_{0,k}+\overline{\tau}_{1,k}-2\overline{d}_{0,k}}{\overline{d}_{0,k}}\bigm|>(1+\frac{\overline{\delta}}{r})^{\left\lfloor\frac{n}{2}\right\rfloor}>(1+\frac{\overline{\delta}}{r})^{\cdot{\frac{n}{2}-1}}=\overline{C}\cdot \overline{E}^n$.\qedhere
\end{proof}
\COMMENT{\large\red\bfseries Explain at start of this section what it is for\wentao{This section is rewritten/changed a lot. The subsection is absorbed to subsection 8.1: introducing 1-petal billiard its billiard map and continuous extension etc.}}

\begin{lemma}
[Approximation for consecutive collisions on $\Gamma_r$]\label{lemma:ApproximationForConsecutiveCollisions}\hfill

    For $\overline{x}=(\overline{\phi},\overline{\theta})\in 
    \Mrin$, if for the \cref{def:LemonAndPetTogether}\eqref{item04LemonAndPetTogether} defined $m(\overline{x})$, $\bar F^{m(\overline{x})}_{\phistar}\notin\partial M_r$ (not hitting corner), then we get the following conclusions.
    \begin{enumerate}
        \item There exists some $\overline{\delta}_0(\overline{x})>0$ such that if $x=(\phi,\theta)\in\Mrin$ satisfies
\begin{equation}\label{eq1:ApproximationForConsecutiveCollisions}
        \begin{aligned}
            |\phi-\overline{\phi}|<\overline{\delta}_0(\overline{x}),\\
            |\theta-\overline{\theta}|<\overline{\delta}_0(\overline{x}),
        \end{aligned}
    \end{equation} then the (\cref{def:LemonAndPetTogether}\eqref{item04LemonAndPetTogether} defined) $m(x)=m(\overline{x})$
    \item For $\forall\varepsilon>0$, there exists some $\delta_0(\overline{x},\varepsilon)\dfn\min\big\{\overline{\delta}_0(\overline{x}),\,\frac{\varepsilon}{2},\,\frac{\varepsilon}{2m(\overline{x})}\big\}$ such that if $x=(\phi,\theta)\in\Mrin$ satisfies the following:
        \begin{equation}\label{eq2:ApproximationForConsecutiveCollisions}
        \begin{aligned}
            |\phi-\overline{\phi}|<\delta_0(\overline{x},\varepsilon),\\
            |\theta-\overline{\theta}|<\delta_0(\overline{x},\varepsilon),
        \end{aligned}
    \end{equation} then the $\phi\theta$-coordinates of $\Fc^{m(x)}(\phi,\theta)=(\phi+2m(x)\theta,\quad\theta)\in\Mrout$ and $\bar F_{\phistar}^{m(\overline{x})}(\overline{\phi},\overline{\theta})=(\overline{\phi}+2m(\overline{x})\overline{\theta},\quad\overline{\theta})\in\Mrout$ are $\varepsilon$ close. i.e. $\big|\overline{\theta}-\theta\bigm|<\varepsilon$, $\big|\phi+2m(x)\theta-(\overline{\phi}+2m(\overline{x})\overline{\theta})\bigm|<\varepsilon$.
    \end{enumerate}
    
    If $0<\phistar<\tan^{-1}{(1/3)}$ and $\pi/\phistar\notin\mathbb{N}$, then $m_*\dfn\inf\big\{ k>0\bigm|\phistar+2k\phistar\in(-\phistar,+\phistar)(\text{ mod }2\pi)\big\}<\infty$, i.e., $\bar F_{\phistar}^{k}(Iy_*)\in M_{\textrm{f}}$ (see \cref{def:LemonAndPetTogether}\eqref{item02LemonAndPetTogether}, \cref{fig:Mr,fig:PetalMroutMfMrinExtension1}).
    \begin{enumerate}[start=3]
        \item There exists some $\overline{\delta}_0(Iy_*)>0$ such that if $x=(\phi,\theta)\in\Mrin$ satisfies        \begin{equation}\label{eq3:ApproximationForConsecutiveCollisions}
            \begin{aligned}
                |\phi-\overline{\phi}|<\overline{\delta}_0(Iy_*),\\
                |\theta-\overline{\theta}|      <\overline{\delta}_0(Iy_*),
            \end{aligned}
        \end{equation} 
        then the \cref{def:LemonAndPetTogether}\eqref{item04LemonAndPetTogether} defined $m(x)=m_*$.
        \item For $\forall\varepsilon>0$, there exists some $\delta_0(Iy_*,\varepsilon)\dfn\min\big\{\overline{\delta}_0(Iy_*),\,\frac{\varepsilon}{2},\,\frac{\varepsilon}{2m_*}\big\}$ such that if $x=(\phi,\theta)\in\Mrin$ satisfies the following:
        \begin{equation}\label{eq4:ApproximationForConsecutiveCollisions}
        \begin{aligned}
            |\phi-\overline{\phi}|<\delta_0(Iy_*,\varepsilon)\\
            |\theta-\overline{\theta}|<\delta_0(Iy_*,\varepsilon),
        \end{aligned}
    \end{equation} then the $\phi\theta$-coordinates of $\Fc^{m(x)}(\phi,\theta)=(\phi+2m(x)\theta,\,\theta)\in\Mrout$ and $\bar F_{\phistar}^{m_*}(\phistar,\phistar)=(\phistar+2m_*\phistar,\,\phistar)\in\Mrout$ are $\varepsilon$ close. i.e. $\big|\overline{\theta}-\theta\bigm|<\varepsilon$, $\big|\phi+2m(x)\theta-(\phistar+2m_*\phistar)\bigm|<\varepsilon$.
    \end{enumerate}
\end{lemma}
\begin{proof}
The first conclusion follows from the fact that the  $\overline{x}$ is in the interior of $\Min_{r,m(\overline{x})}$ which is a finite number of disconnected quadrilateral components (see \cref{fig:MrN}). It is easy to see that the second conclusion follows from the first conclusion.

The third conclusion follows from the fact that $Iy_*$ is on the boundary of $\Min_{r,m_*}$, which is a finite number of disconnected quadrilateral components. It is easy to see that the fourth conclusion follows from the third conclusion.
\end{proof}

\subsection[Trajectory closeness between lemon billiard and 1-petal billiard]{Trajectory closeness  between lemon billiard and 1-petal billiard}\strut\COMMENT{\red Subsection 8.3 needs explanations and guided tour what's going on, not just computations.\Hrule\large\bfseries Summarize this subsection here! A carefully written summary of the prose before the two main lemmas should serve, with explicit forward reference to both lemmas.\Hrule Should they be lemmas or propositions?\wentao{From usage purpose, it should be Lemmas.}}

\begin{notation}[\cref{fig:Trajectory_LemonCompareWithWojtkowski} collision position and angle]\label{def:flowcollsionpositionandangle}\label{Def9_2}\label{Def9_3} 


When considering  $\LEM(r,R,\phistar)$, we denote by $(\mathcal{X},\mathcal{Y},\mathcal{\Theta})$ the cartesian position and direction coordinates of $x\in M_R$.

When considering $\PET(r,\phistar)$, denote by 
$(\overline{\mathcal{X}},\overline{\mathcal{Y}},\overline{\mathcal{\Theta}})$ the position and direction of $\overline{x}\in M_{\text{\upshape{f}}}$, where $M_{\text{\upshape{f}}}$ is given in \cref{def:LemonAndPetTogether}\eqref{item02LemonAndPetTogether}.
\end{notation}

\begin{figure}[ht]
\begin{center}
\begin{tikzpicture}[xscale=0.8,yscale=0.8]
    \tkzDefPoint(0,0){Or};
    \tkzDefPoint(0,-4.8){Y};
    \pgfmathsetmacro{\rradius}{4.5};
    \clip(-2.4*\rradius,-1.1*\rradius) rectangle (1.3*\rradius, 0.03*\rradius);
    \pgfmathsetmacro{\phistardeg}{40};
    \pgfmathsetmacro{\XValueArc}{\rradius*sin(\phistardeg)};
    \pgfmathsetmacro{\YValueArc}{\rradius*cos(\phistardeg)};
    \pgfmathsetmacro{\Rradius}{
        \rradius*sin(\phistardeg)/sin(12)
    };
    \pgfmathsetmacro{\bdist}{(-\rradius*cos(\phistardeg))+cos(12)*\Rradius}
    \tkzDefPoint(0,\bdist){OR};
    \tkzDefPoint(\XValueArc,-1.0*\YValueArc){B};    \tkzDefPoint(-1.0*\XValueArc,-1.0*\YValueArc){A};
    \tkzDrawArc[name path=Cr, very thick](Or,B)(A);
    \tkzDrawArc[name path=CR](OR,A)(B);
    \pgfmathsetmacro{\PXValue}{-\Rradius*sin(3)};
    \pgfmathsetmacro{\PYValue}{\bdist-\Rradius*cos(3)};
    \tkzDefPoint(\PXValue,\PYValue){P};
    
    
    \path[name path=Line_P_T0_far](P)-- (176:6);
    \path[name path=Line_P_T_far](P)-- (172:6);
    \path[name intersections={of=Cr and Line_P_T0_far,by={T0}}];
    \path[name intersections={of=Cr and Line_P_T_far,by={T}}];
    \tkzDefPoint(0,\bdist-\Rradius){C};
    \path[name intersections={of=Cr and Line_P_T0_far,by={T0}}];
    \draw[name path=LineAB,thin] (A)--(B);
     \coordinate (PHorizon) at ([shift={(3,0)}]P);
    \draw[red,ultra thin,dashed] (P)--(PHorizon);
    \path[name intersections={of=Line_P_T0_far and LineAB,by={P_bar}}];
    \tkzDefPointBy[projection=onto Or--P](M)\tkzGetPoint{M2};
    
    \begin{scope}[ultra thin,decoration={markings,mark=at position 0.5 with {\arrow{>}}}] 
    \draw[red,postaction={decorate}] (T0)--(P_bar);
    \draw[blue,postaction={decorate}] (T)--(P);
    \end{scope}
    \tkzLabelPoint[left](A){$A$};
    \tkzLabelPoint[below](P){\small $P=p(\Fc(x))=(\mathcal{X}_1,\mathcal{Y}_1)$};
    \tkzLabelPoint[below left](T0){\small $T_0=p(\overline{x}_0)=(r\sin{\overline{\phi}_0},-r\cos{\overline{\phi}_0})$};
    \tkzLabelPoint[left](T){\small $T=p(x)=(r\sin{\phi},-r\cos{\phi})$};
    \tkzLabelPoint[below](Or){$O_r$};
    \tkzLabelPoint[below](B){$B$};
    \tkzLabelPoint[above right]([shift={(0.15,0.23)}]P_bar){\small $\overline{P}=p(\bar F_{\phistar}(\overline{x}_0))=(\overline{\mathcal{X}}_1,\overline{\mathcal{Y}}_1)$};
    \tkzLabelPoint[right]([shift={(0.08,0.15)}]P_bar){\small $\overline{\Theta}_1$};
    \tkzLabelPoint[right]([shift={(0.15,0.15)}]P){\small $\Theta_1$};
    \tkzDefPointBy[projection=onto B--A](T0)\tkzGetPoint{M0};
    \coordinate (T0Mirror) at ($(T0)!2!(M0)$);
    \coordinate (PbarArrow) at ($(T0Mirror)!1.16!(P_bar)$);
    \tkzDefPointBy[projection=onto P--OR](T)\tkzGetPoint{M};
    \coordinate (TMirror) at ($(T)!2!(M)$);
    \coordinate (PArrow) at ($(P)!0.18!(TMirror)$);
    \begin{scope}[ultra thin,decoration={markings,mark=at position 1.0 with {\arrow{>}}}] 
    \draw[red,postaction={decorate},ultra thin] (P_bar)--(PbarArrow);
    \draw[blue,postaction={decorate},ultra thin] (P)--(PArrow);
    \end{scope}
    \tkzMarkAngle[arc=lll,ultra thin,size=0.15](B,P_bar,PbarArrow);
    \tkzMarkAngle[arc=lll,ultra thin,size=0.18](PHorizon,P,PArrow);
    \tkzDrawPoints(A,B,Or,P,T0,T,P_bar);
\end{tikzpicture}
\caption{
$\overline{P}=(\overline{\mathcal{X}}_1,\overline{\mathcal{Y}}_1)=p(\bar F_{\phistar}(x_0))$ on $\PET(r,\phistar)$.
\quad $P=(\mathcal{X}_1,\mathcal{Y}_1)=p(\Fc(x))$ on $\LEM(r,R,\phistar)$.
} \label{fig:Trajectory_LemonCompareWithWojtkowski}
\end{center}
\end{figure}



We now come to our first trajectory approximation lemma, \cref{lemma:Trajectory_approximation_1}. \cref{lemma:Trajectory_approximation_1} supposes that two starting positions and directions on $\Gamma_r$ respectively on lemon billiard $\mathbb{L}(r,R,\phistar)$ and petal billiard $\mathbb{P}(r,\phistar)$ (see \cref{fig:Trajectory_LemonCompareWithWojtkowski}) are close enough and that $R$ is large enough and that their respective second collisions are on the arc boundary $\Gamma_R$ in case of the lemon billiard and on the flat boundary $\overline{AB}$ in case of the petal billiard. Then these respective second collisions are also close in collision position and angle.


\begin{lemma}[Shadowing 1: approximation for $\bar F_{\phistar}$ trajectory in \cref{fig:Trajectory_LemonCompareWithWojtkowski}]\label{lemma:Trajectory_approximation_1}
Consider the collision maps $\bar F_{\phistar}$ for $\PET(r,\phistar)$(\cref{def:LemonAndPetTogether}) with the continuous extensions $\bar F_{\phistar,1}$(\cref{notation:PetContinuousExtension}\eqref{item01PetContinuousExtension}) and $\Fc$ for $\LEM(r,R,\phistar)$ (\cref{def:LemonAndPetTogether}) and let $x=(\phi,\theta)\in\Mrout$, $\overline{x}_0=(\overline{\phi}_0,\overline{\theta}_0)\in$closure$(\Mrout)$ with $x_0\notin \{x_*,y_*\}$ (\cref{fig:PetalMroutMfMrinExtension1})
as well as $x_1=(\Phi_1,\theta_1)\nfd\Fc(x)\in M_R$. Write $p(x_1)=(\mathcal{X}_1,\mathcal{Y}_1)$, so $(\mathcal{X}_1,\mathcal{Y}_1,\Theta_1)$ represents the $\LEM(r,R,\phistar)$ collision position and reflection angle (\cref{def:flowcollsionpositionandangle}). 

Also write $\bar F_{\phistar,1}(\overline{x}_0)\nfd \overline{x}_1\in M_{\text{\upshape{f}}}$, $p(\overline{x}_1)=(\overline{\mathcal{X}}_1,\overline{\mathcal{Y}}_1)$, so $(\overline{\mathcal{X}}_1,\overline{\mathcal{Y}}_1,\overline{\Theta}_1)$ represents the $\PET(r,\phistar)$
collision position and reflection angle (\cref{def:flowcollsionpositionandangle}).

For $\forall \varepsilon,\varepsilon_1>0$, if $x$, $\overline{x}_0$ and $R$ satisfy \begin{equation}\label{eq:x0barx0ClosingCondition}
    \begin{aligned}
        |\theta-\overline{\theta}_0|<\mathtt{\delta}_1(\overline{\phi}_0,\overline{\theta}_0,\varepsilon,\varepsilon_1)\dfn\min&{\Big\{0.01\sin^2{(\overline{\phi}_0+\overline{\theta}_0)},-0.05\frac{\varepsilon}{r}\sin^3{(\overline{\phi}_0+\overline{\theta}_0)},0.4\varepsilon_1,\frac{\pi}{2}-\frac{\overline{\phi}_0+\overline{\theta}_0}{4},\frac{\overline{\phi}_0+\overline{\theta}_0}{4}\Big\}},\\
        |\phi-\overline{\phi}_0|<\mathtt{\delta}_2(\overline{\phi}_0,\overline{\theta}_0,\varepsilon,\varepsilon_1)=\min&{\Big\{0.01\sin^2{(\overline{\phi}_0+\overline{\theta}_0)},-0.05\frac{\varepsilon}{r}\sin^3{(\overline{\phi}_0+\overline{\theta}_0)},0.4\varepsilon_1,\frac{\pi}{2}-\frac{\overline{\phi}_0+\overline{\theta}_0}{4},\frac{\overline{\phi}_0+\overline{\theta}_0}{4}\Big\}},\\
        R>\mathtt{R}_1(\overline{\phi}_0,\overline{\theta}_0,\varepsilon,\varepsilon_1)\dfn\max&\Big\{r+\frac{23r^2}{-\varepsilon\sin^3{(\overline{\phi}_0+\overline{\theta}_0)}},r+\frac{500r}{\sin^2{(\overline{\phi}_0+\overline{\theta}_0)}},\frac{2r^2\sin^2{\phistar}}{\varepsilon},\frac{2.02r\sin{\phistar}}{0.4\varepsilon_1},\\
        &\qquad1700r
        ,\frac{4r}{\sin^2{(\frac{2\pi-\overline{\phi}_0-\overline{\theta}_0}{4})}},\frac{4r}{\sin^2{(\frac{\overline{\theta}_0+\overline{\phi}_0-\pi}{4})}}\Big\},
    \end{aligned}
\end{equation} then (see \eqref{eq:theta1theta1barApproximation})
\begin{equation}\label{eq:eq:x0barx0Closing}
        \begin{aligned}
            |\mathcal{X}_1-\overline{\mathcal{X}}_1|<\varepsilon,\quad
            |\mathcal{Y}_1-\overline{\mathcal{Y}}_1|<\varepsilon,\quad
            |\Theta_1-\overline{\Theta}_1|<\varepsilon_1.
        \end{aligned}
    \end{equation}
In particular, if\COMMENT{Pointless aligment in next equation. Do something better.\wentao{Now using multline saving/reducing one line cost.}}    \begin{equation}\label{eq:xorbitA0Condition}
        \begin{multlined}
            \big|\theta-\overline{\theta}_0\bigm|<\min\big\{\pi/2-(\overline{\phi}_0+\overline{\theta}_0)/4,(\overline{\phi}_0+\overline{\theta}_0)/4\big\},\quad
            \bigm|\phi-\overline{\phi}_0\bigm|<\min\big\{\pi/2-(\overline{\phi}_0+\overline{\theta}_0)/4,(\overline{\phi}_0+\overline{\theta}_0)/4\big\},\\
            R>\max\Big\{\frac{4r}{\sin^2{(\frac{2\pi-\overline{\phi}_0-\overline{\theta}_0}{4})}},\frac{4r}{\sin^2{(\frac{\overline{\theta}_0+\overline{\phi}_0-\pi}{4})}}\Big\},
        \end{multlined}
    \end{equation}
(this clearly follows from \eqref{eq:x0barx0ClosingCondition}), then $\sin{\theta_1}>\sqrt{4r/R}$. This means that if the $\LEM(r,R,\phistar)$ billiard $\hat{M}$ return orbit in \cref{def:MhatReturnOrbitSegment} with $x_0=x=(\phi,\theta)$ satisfies  \eqref{eq:x0barx0ClosingCondition}, then $x_1=\Fcal(x)=\Fcal(x_0)$ falls in the case (a0) of \eqref{eqMainCases}.
\end{lemma}
\begin{remark}\label{remark:sinPhi0plusTheta0isNegative}
        We notice that if $x=(\phi,\theta)\in\Mrout$ and $\overline{x}_0=(\overline{\phi}_0,\overline{\theta}_0)\in$ closure$(\Mrout)$ with $\overline{x}_0\notin\{x_*,y_*\}$, then by checking \cref{fig:PetalMroutMfMrinExtension1} we see $\pi<\phi+\theta<2\pi$ and $\pi<\overline{\phi}_0+\overline{\theta}_0<2\pi$.
\end{remark}
\begin{proof}[Proof of \cref{lemma:Trajectory_approximation_1}]
    In \cref{fig:Trajectory_LemonCompareWithWojtkowski} the coordinate system with $O_r=(0,0)$, $(0,b)$ and $\overline{AB}$ is on line $y=-r\cos\phistar$, we have the following coordinates of points.
    \begin{equation}\label{eq:xT0xTclose}
        \begin{aligned}
        (x_T,y_T)=T\dfn&p(x)=(r\sin\phi,-r\cos\phi)\in\Gamma_r\\ (x_{T_0},y_{T_0})=T_0\dfn&p(\overline{x}_0)=(r\sin{\overline{\phi}_0},-r\cos{\overline{\phi}_0})\in\Gamma_r\\
        P\dfn&p(\Fcal(x))=(\mathcal{X}_1,\mathcal{Y}_1)\in\Gamma_R\\
        \overline{P}\dfn&p(\bar F_{\phistar}(\overline{x}_0))=(\overline{\mathcal{X}}_1,\overline{\mathcal{Y}}_1)\in\overline{AB}\text{ with }\overline{\mathcal{Y}}_1=-r\cos\phistar.\\
        \end{aligned}
    \end{equation}

    On the other hand, by the geometric meaning of $\phi$ and $\theta$, given counterclockwise as the positive orientation for the billiard plane, the vector $\overrightarrow{TP}$ has an angle $\phi+\theta$ with respect to the positive axis $x$ and the vector $\overrightarrow{T_0\overline{P}}$ has an angle $\overline{\phi}_0+\overline{\theta}_0$ with respect to the positive axis $x$. Note that $2\pi-\overline{\Theta}_1$ is also the angle between $\overrightarrow{T_0\overline{P}}$ and the positive axis $x$ given the counterclockwise orientation. Therefore, 
    \begin{equation}\label{eq:X1barY1barTheta1bar}
        \begin{aligned}
            2\pi-\overline{\Theta}_1&=\overline{\phi}_0+\overline{\theta}_0,\\
            \frac{-r\cos\phi+r\cos\phistar}{r\sin\phi-\overline{\mathcal{X}}_1}\overset{\eqref{eq:xT0xTclose}}=\frac{y_T-\overline{\mathcal{Y}}_1}{x_T-\overline{\mathcal{X}}_1}&\overbracket{=}^{\mathrlap{\text{slope of }T\overline{P}}}\tan{(\overline{\phi}_0+\overline{\theta}_0)}=\frac{\sin{(\overline{\phi}_0+\overline{\theta}_0)}}{\cos{(\overline{\phi}_0+\overline{\theta}_0)}}\\
            \text{hence } \overline{P}=(\overline{\mathcal{X}}_1,\overbracket{\overline{\mathcal{Y}}_1}^{\mathllap{\eqref{eq:xT0xTclose}:=-r\cos\phistar}})&=(r\sin{\overline{\phi}_0}+\frac{r\cos{\overline{\phi}_0}-r\cos{\phistar}}{\sin{(\overline{\phi}_0+\overline{\theta}_0)}}\cos{(\overline{\phi}_0+\overline{\theta}_0)},-r\cos{\phistar}).
        \end{aligned}
    \end{equation} Let $t=|PT|$, since $P$ is on the circle $C_R$, we get the following equations for $(\mathcal{X}_1,\mathcal{Y}_1)$.
\begin{equation}\label{eq:X1Y1Polynomial}
\begin{aligned}
    &\left.\begin{aligned}\mathcal{X}_1&=r\sin{\phi}+t\cos{(\phi+\theta)}\\
        \mathcal{Y}_1&=-r\cos{\phi}+t\sin{(\phi+\theta)}\\
       R^2&=\mathcal{X}^2_1+(\mathcal{Y}_1-b)^2\\
    \end{aligned}\right\}\Longrightarrow R^2=\big(\overbracket{r\sin{\phi}+t\cos{(\phi+\theta)}}^{=\mathcal{X}_1}\big)^2+\big(\overbracket{-r\cos{\phi}+t\sin{(\phi+\theta)}-b}^{=\mathcal{Y}_1-b}\big)^2\\
    \Longrightarrow&R^2=r^2+t^2+b^2+\overbracket{2rt\sin\phi\cos{(\phi+\theta)}-2rt\cos\phi\sin{(\phi+\theta)}}^{-2rt\sin\theta}-2bt\sin{(\phi+\theta)}+2br\cos\phi\\
    \Longrightarrow&t^2-2rt\sin\theta-2bt\sin{(\phi+\theta)}+2br\cos\phi+r^2+b^2-R^2=0\\
    \overbracket{\Longrightarrow}^{\mathrlap{\substack{\text{cosine law: }\\R^2=r^2+b^2-2br\cos{\angle{O_RO_rB}}\\
    \measuredangle{O_RO_rB}=\pi-\phistar}}}&t^2-2\big(r\sin\theta+b\sin{(\phi+\theta)}\big)t+2br\cos\phi-2br\cos\phistar=0.
    \end{aligned} 
\end{equation}
    Since $\phi\in(\phistar,2\pi-\phistar)$, quadratic polynomial $t^2-2\big(r\sin\theta+b\sin{(\phi+\theta)}\big)t+2br\cos\phi-2br\cos\phistar$ has $|PT|$ as its nonnegative root $t_+$
    \begin{equation}\label{eq:PTposRoot}
    \begin{aligned}|PT|=t_+&=b\sin{(\phi+\theta)}+r\sin\theta+\sqrt{\big(b\sin{(\phi+\theta)}+r\sin\theta\big)^2+\big(2br\cos\phistar-2br\cos\phi\big)}\\
    &=\frac{2br\cos\phistar-2br\cos\phi}{-b\sin{(\phi+\theta)}-r\sin\theta+\sqrt{\big(b\sin{(\phi+\theta)}+r\sin\theta\big)^2+\big(2br\cos\phistar-2br\cos\phi\big)}}\\
    &=\frac{2r\cos\phistar-2r\cos\phi}{-\sin{(\phi+\theta)}-(r/b)\sin\theta+\sqrt{\big(\sin{(\phi+\theta)}+(r/b)\sin\theta\big)^2+\big(2(r/b)\cos\phistar-2(r/b)\cos\phi\big)}}.
    \end{aligned}
    \end{equation}
\newline Note that with $r$ fixed, as $R\rightarrow\infty$, $b\overbracket{>}^{\mathclap{\substack{b,r,R\text{ are edge lengths of }\\ \triangle O_rO_RB\text{ in \cref{fig:1-petalExtensionFromCornerAndLemon}}}}}R-r\rightarrow\infty$. So, with $(r/b)\rightarrow 0$, $\phi\rightarrow\overline{\phi}_0$, $\theta\rightarrow\overline{\theta}_0$, the right-hand side of the last equality in \eqref{eq:PTposRoot} has the limit $t_+\longrightarrow\frac{2r\cos\phistar-2r\cos{\overline{\phi}_0}}{-\sin{(\overline{\phi}_0+\overline{\theta}_0)}+|\sin{(\overline{\phi}_0+\overline{\theta}_0)}|}\overbracket{=}^{\mathclap{\substack{\pi<\overline{\phi}_0+\overline{\theta}_0<2\pi\\\sin{(\overline{\phi}_0+\overline{\theta}_0)}<0}}}\frac{r\cos\phistar-r\cos{\overline{\phi}_0}}{-\sin{(\overline{\phi}_0+\overline{\theta}_0)}}$. Therefore, if we write 
\begin{equation}\label{eq:xietaphitheta}
    \begin{aligned}
        \xi(\phi,\theta)&=(2r\cos{\phistar}-2r\cos{\phi})\cos{(\phi+\theta)},\\
        \eta(\phi,\theta)&=\sqrt{\sin^2{(\theta+\phi)}+(r/b)^2\sin^2{\theta}+2(r/b)\sin{\theta}\sin{(\theta+\phi)}+2(r/b)\cos{\phistar}-2(r/b)\cos{\phi}},\\
        &-\sin{(\theta+\phi)}-(r/b)\sin{\theta}\\
        \frac{\xi(\phi,\theta)}{\eta(\phi,\theta)}=t_+\cos{(\phi+\theta)}\overbracket{\longrightarrow}^{\mathclap{\substack{(r/b)\to0,\\\phi\to\overline{\phi}_0,\\\theta\to\overline{\theta}_0}}}&\frac{r\cos\phistar-r\cos{\overline{\phi}_0}}{-\sin{(\overline{\phi}_0+\overline{\theta}_0)}}\cos{(\overline{\phi}_0+\overline{\theta}_0)},
    \end{aligned}
\end{equation}
\newline then by \eqref{eq:PTposRoot}, \eqref{eq:xietaphitheta}, \eqref{eq:X1Y1Polynomial}, \eqref{eq:X1barY1barTheta1bar} we conclude
\begin{equation}\label{eq:X1X1bardiff}
\begin{multlined}\mathcal{X}_1\overbracket{=}^{\mathllap{\text{\eqref{eq:X1Y1Polynomial}}}}r\sin{\phi}+\overbracket{t_+\cos{(\phi+\theta)}}^{\text{\eqref{eq:xietaphitheta}: }=\frac{\xi(\phi,\theta)}{\eta(\phi,\theta)}}\overbracket{\longrightarrow}^{\mathclap{\substack{(r/b)\to0,\\\phi\to\overline{\phi}_0,\\\theta\to\overline{\theta}_0}}}r\sin\overline{\phi}_0+\frac{r\cos\phistar-r\cos{\overline{\phi}_0}}{-\sin{(\overline{\phi}_0+\overline{\theta}_0)}}\cos{(\overline{\phi}_0+\overline{\theta}_0)}\overbracket{=}^{\qquad\mathllap{\text{\eqref{eq:X1barY1barTheta1bar}}}}\overline{\mathcal{X}}_1,\\
    \mathcal{X}_1-\overline{\mathcal{X}}_1\overbracket{=}^{\mathllap{\text{\eqref{eq:X1barY1barTheta1bar}\eqref{eq:X1Y1Polynomial}}}}r\sin\phi-r\sin{\overline{\phi}_0}+\frac{\xi(\phi,\theta)}{\eta(\phi,\theta)}-\frac{r\cos\phistar-r\cos{\overline{\phi}_0}}{-\sin{(\overline{\phi}_0+\overline{\theta}_0)}}\cos{(\overline{\phi}_0+\overline{\theta}_0)}.\\
    \end{multlined}
\end{equation} For the $\xi(\phi,\theta)$ limit approximation, we observe that \begin{equation}\label{eq:xiapproximation}
    \begin{aligned}
        &\big|\xi(\phi,\theta)-(2r\cos{\phistar}-2r\cos{\overline{\phi}_0})\cos{(\overline{\phi}_0+\overline{\theta}_0)}\big|\\
        \overbracket{=}^{\mathclap{\text{\eqref{eq:xietaphitheta}}}}&\bigm|\big(2r\cos{\phistar}-2r\cos{\phi}\big)\cos{(\phi+\theta)}-\big(2r\cos{\phistar}-2r\cos{\overline{\phi}_0}\big)\cos{(\overline{\phi}_0+\overline{\theta}_0)}\big|\\
        &=2r\bigm|\cos{\phistar}\big(\cos{(\phi+\theta)}-\cos{(\overline{\phi}_0+\overline{\theta}_0)}\big)+\cos{\overline{\phi}_0}\cos{(\overline{\phi}_0+\overline{\theta}_0)}-\cos{\phi}\cos{(\phi+\theta)}\big|\\
        &=2r\bigm|\cos{\phistar}(\cos{(\phi+\theta)}-\cos{(\overline{\phi}_0+\overline{\theta}_0)})+\cos{\overline{\phi}_0}\cos{(\overline{\phi}_0+\overline{\theta}_0)}-\cos{\overline{\phi}_0}\cos{(\phi+\theta)}\\
        &+\cos{\overline{\phi}_0}\cos{(\phi+\theta)}-\cos{\phi}\cos{(\phi+\theta)}\big|\\
        \le&2r\cos{\phistar}\bigm|\overbracket{\cos{(\phi+\theta)}-\cos{(\overline{\phi}_0+\overline{\theta}_0)}}^{\mathclap{\substack{\text{Mean-value thm}\Rightarrow\text{there exists some }w_1\\  \cos{(\phi+\theta)}-\cos{(\overline{\phi}_0+\overline{\theta}_0)}=-(\phi+\theta-\overline{\phi}_0-\overline{\theta}_0)\sin{w_1}}}}\bigm|+2r\bigm|\cos{\overline{\phi}_0}\bigm|\bigm|\overbracket{\cos{(\phi+\theta)}-\cos{(\overline{\phi}_0+\overline{\theta}_0)}}^{\mathclap{\substack{\text{Mean-value thm}\Rightarrow \text{ there exists some }w_1\\  \cos{(\phi+\theta)}-\cos{(\overline{\phi}_0+\overline{\theta}_0)}=-(\phi+\theta-\overline{\phi}_0-\overline{\theta}_0)\sin{w_1}}}}\big|\\
&+2r\bigm|\cos{(\phi+\theta)}\bigm|\bigm|\overbracket{\cos{\overline{\phi}_0}-\cos\phi}^{{\mathclap{\substack{\text{Mean-value thm}\Rightarrow\text{ there exists some }w_2\\  \cos{\overline{\phi}_0}-\cos{\phi}=(\overline{\phi}_0-\phi)\sin{w_2}}}}}\big|\\
        \overbracket{\le}^{\mathclap{\text{Mean-value thm}\Rightarrow }}&2r\bigm|\phi+\theta-\overline{\phi}_0-\overline{\theta}_0\bigm|+2r\bigm|\phi+\theta-\overline{\phi}_0-\overline{\theta}_0\bigm|+2r\bigm|\overline{\phi}_0-\phi\big|\\
        \le&6r\bigm|\overline{\phi}_0-\phi\bigm|+4r\bigm|\overline{\theta}_0-\theta\big|.
    \end{aligned}
\end{equation} For the $\eta(\phi,\theta)$ limit approximation, we see that
\begin{equation}\label{eq:etaaproximation}
\begin{aligned}
    &\big|\eta(\phi,\theta)-(-2\sin{(\overline{\phi}_0+\overline{\theta}_0)})\big|\\\overbracket{=}^{\mathclap{\text{\eqref{eq:xietaphitheta}}}}&\bigm|\sqrt{\sin^2{(\theta+\phi)}+(r/b)^2\sin^2{\theta}+2(r/b)\sin{\theta}\sin{(\theta+\phi)}+2(r/b)\cos{\phistar}-2(r/b)\cos{\phi}}\\
    &-\sin{(\theta+\phi)}-(r/b)\sin{\theta}+2\sin{(\overline{\phi}_0+\overline{\theta}_0)}\big|\\
    \overbracket{\le}^{\mathclap{\substack{\text{Triangle inequality}}}}& \bigm|\sin{(\overline{\phi}_0+\overline{\theta}_0)}-\sin{(\theta+\phi)}\bigm|+\bigm|\frac{r\sin{\theta}}{b}\big|\\
    &+\Bigm|\sqrt{\sin^2{(\theta+\phi)}+\frac{r^2\sin^2{\theta}}{b^2}+\frac{2r\sin{\theta}\sin{(\theta+\phi)}}{b}+\frac{2r\cos{\phistar}-2r\cos{\phi}}{b}}+\sin{(\overline{\phi}_0+\overline{\theta}_0)}\Bigm|\\
    \overbracket{\le}^{\mathclap{\text{Mean-value thm}}}&|\theta+\phi-\overline{\theta}_0-\overline{\phi}_0|+\frac{r}{b}\\
    &+\Bigm|\frac{(r^2/b^2)\sin^2{\theta}+(2r/b)\sin{\theta}\sin{(\theta+\phi)}+(2r/b)(\cos{\phistar}-\cos{\phi})+\sin^2{(\theta+\phi)}-\sin^2{(\overline{\phi}_0+\overline{\theta}_0)}}{\underbracket{\sqrt{\sin^2{(\theta+\phi)}+\frac{r^2\sin^2{\theta}}{b^2}+\frac{2r\sin{\theta}\sin{(\theta+\phi)}}{b}+\frac{2r\cos{\phistar}-2r\cos{\phi}}{b}}-\sin{(\overline{\phi}_0+\overline{\theta}_0)}}_{\ge-\sin{(\overline{\phi}_0+\overline{\theta}_0)}>0\text{, by }\pi<\overline{\phi}_0+\overline{\theta}_0<2\pi}}\Big|\\
    \le&|\theta+\phi-\overline{\theta}_0-\overline{\phi}_0|+\frac{r}{b}+\frac{(r^2/b^2+6r/b)}{-\sin{(\overline{\phi}_0+\overline{\theta}_0)}}+\frac{|\sin^2{(\phi+\theta)}-\sin^2{(\overline{\phi}_0+\overline{\theta}_0)}|}{-\sin{(\overline{\phi}_0+\overline{\theta}_0)}}\\
    \overbracket{\le}^{\mathclap{\text{Mean-value thm}}}&|\theta+\phi-\overline{\theta}_0-\overline{\phi}_0|+\overbracket{\frac{r}{b}+\frac{(r^2/b^2+6r/b)}{-\sin{(\overline{\phi}_0+\overline{\theta}_0)}}}^{<(7.1r)/b}+\frac{2|\theta+\phi-\overline{\theta}_0-\overline{\phi}_0|}{-\sin{(\overline{\phi}_0+\overline{\theta}_0)}}\\
    \overbracket{<}^{\mathclap{\substack{b>R-r>1699r\\\text{and }\frac{1}{-\sin{(\overline{\phi}_0+\overline{\theta}_0)}}\ge1\\|}}}&\frac{7.1r}{b}\frac{1}{-\sin{(\overline{\phi}_0+\overline{\theta}_0)}}+\frac{3|\theta-\overline{\theta}_0|}{-\sin{(\overline{\phi}_0+\overline{\theta}_0)}}+\frac{3|\phi-\overline{\phi}_0|}{-\sin{(\overline{\phi}_0+\overline{\theta}_0)}}.\\
\end{aligned}
\end{equation}\newline Hence, if  
\begin{equation}\label{eq:thetaphiclose}
    \begin{aligned}
    \big|\theta-\overline{\theta}_0\big|<&0.01\sin^2{(\overline{\phi}_0+\overline{\theta}_0)},\\
    \big|\phi-\overline{\phi}_0\big|<&0.01\sin^2{(\overline{\phi}_0+\overline{\theta}_0)}, \\
    b>R-r>&\frac{500r}{\sin^2(\overline{\phi}_0+\overline{\theta}_0)}
    \text{, i.e., }\frac{r}{b}<\frac{1}{500}\sin^2{(\overline{\phi}_0+\overline{\theta}_0)},\\
    \end{aligned}
\end{equation}\newline then since we have the following estimate for $\eta(\phi,\theta)$,
\begin{align*}
    \big|\eta(\phi,\theta)-(-2\sin{(\overline{\phi}_0+\overline{\theta}_0)})\bigm|\overbracket{<}^{\mathclap{\eqref{eq:etaaproximation}}}&\frac{7.1r}{b}\frac{1}{-\sin{(\overline{\phi}_0+\overline{\theta}_0)}}+\frac{3|\theta-\overline{\theta}_0|}{-\sin{(\overline{\phi}_0+\overline{\theta}_0)}}+\frac{3|\phi-\overline{\phi}_0|}{-\sin{(\overline{\phi}_0+\overline{\theta}_0)}}\\
    \overbracket{<}^{\mathclap{\eqref{eq:thetaphiclose}}}&
    \frac{7.1}{500}\frac{\sin^2{(\overline{\phi}_0+\overline{\theta}_0)}}{-\sin{(\overline{\phi}_0+\overline{\theta}_0)}}+\frac{6\times0.01\sin^2{(\overline{\phi}_0+\overline{\theta}_0)}}{-\sin{(\overline{\phi}_0+\overline{\theta}_0)}}=-0.0742\sin{(\overline{\phi}_0+\overline{\theta}_0)},
\end{align*}we find a lower bound for $\eta(\phi,\theta)$: 
\begin{equation}\label{eq:etalowerbound}
    \begin{aligned}
    \eta(\phi,\theta)-(-2\sin{(\overline{\phi}_0+\overline{\theta}_0)})&>0.0742\sin{(\overline{\phi}_0+\overline{\theta}_0)}\\
    \eta(\phi,\theta)&>-1.9258\sin{(\overline{\phi}_0+\overline{\theta}_0)}>0.
    \end{aligned}
\end{equation}
For $\frac{\xi(\phi,\theta)}{\eta(\phi,\theta)}$, we see the following.\COMMENT{Misalignment below; fix. multline plus aligned might work\wentao{This jumble is too complicated..}}
\begin{equation}\label{eq:xietaratioApproximate}
    \begin{aligned}
        &\Big|\frac{\xi(\phi,\theta)}{\eta(\phi,\theta)}-\frac{2r(\cos\phistar-r\cos{\overline{\phi}_0})\cos{(\overline{\phi}_0+\overline{\theta}_0)}}{-2\sin{(\overline{\phi}_0+\overline{\theta}_0)}}\Bigm|=\Bigm|\frac{\xi(\phi,\theta)(-2\sin{(\overline{\phi}_0+\overline{\theta}_0)})-\eta(\phi,\theta)[2r(\cos{\phistar}-\cos{\overline{\phi}_0})\cos{(\overline{\phi}_0+\overline{\theta}_0)}]}{(-2\sin{(\overline{\phi}_0+\overline{\theta}_0)})\eta(\phi,\theta)}\Big|\\
        &=\Bigm|\frac{\xi(\phi,\theta)(-2\sin{(\overline{\phi}_0+\overline{\theta}_0)})-\xi(\phi,\theta)\eta(\phi,\theta)+\xi(\phi,\theta)\eta(\phi,\theta)-\eta(\phi,\theta)[2r(\cos{\phistar}-\cos{\overline{\phi}_0})\cos{(\overline{\phi}_0+\overline{\theta}_0)}]}
        {-2\sin{(\overline{\phi}_0+\overline{\theta}_0)}\eta(\phi,\theta)}\Big|\\
        &\le\frac{\overbracket{|\xi(\phi,\theta)|}^{\le4r}}{-2\sin{(\overline{\phi}_0+\overline{\theta}_0)}\underbracket{|\eta(\phi,\theta)|}_{\mathclap{\eqref{eq:etalowerbound}:>-1.9258\sin{(\overline{\phi}_0+\overline{\theta}_0)}>0}}}\underbracket{\Bigm|-2\sin{(\overline{\phi}_0+\overline{\theta}_0)}-\eta(\phi,\theta)\Bigm|}_{\text{LHS of \eqref{eq:etaaproximation}}}+\frac{\overbracket{\bigm|\xi(\phi,\theta)-(2r\cos{\phistar}-2r\cos{\overline{\phi}_0})\cos{(\overline{\phi}_0+\overline{\theta}_0)}\bigm|}^{\text{LHS of \eqref{eq:xiapproximation}}}}{-2\sin{(\overline{\phi}_0+\overline{\theta}_0)}}\\
        &\overbracket{<}^{\mathrlap{\substack{\text{\eqref{eq:xiapproximation}, \eqref{eq:etaaproximation}}\\ \text{\eqref{eq:etalowerbound}}}}}\frac{2r}{1.9258\sin^2{(\overline{\phi}_0+\overline{\theta}_0)}}\Big[\frac{7.1r}{b}\frac{1}{-\sin{(\overline{\phi}_0+\overline{\theta}_0)}}+\frac{3|\theta-\overline{\theta}_0|}{-\sin{(\overline{\phi}_0+\overline{\theta}_0)}}+\frac{3|\phi-\overline{\phi}_0|}{-\sin{(\overline{\phi}_0+\overline{\theta}_0)}}\Big]+\frac{6r\bigm|\overline{\phi}_0-\phi\bigm|+4r\bigm|\overline{\theta}_0-\theta\bigm|}{-2\sin{(\overline{\phi}_0+\overline{\theta}_0)}}.
    \end{aligned}
\end{equation}

We summarize: if
\begin{equation}\label{eq:X1closeCondition}
\begin{aligned}
    |\phi-\overline{\phi}_0|<&\min{\big\{\underbracket{0.01\sin^2{(\overline{\phi}_0+\overline{\theta}_0)}}_{\mathclap{\text{required in \eqref{eq:thetaphiclose}}}},\underbracket{-0.05\frac{\varepsilon}{r}\sin^3{(\overline{\phi}_0+\overline{\theta}_0)}}_{\text{used in \eqref{eq:X1X1BarClose}}}\big\}},\\ |\theta-\overline{\theta}_0|<&\min{\big\{\underbracket{0.01\sin^2{(\overline{\phi}_0+\overline{\theta}_0)}}_{\mathclap{\text{required in \eqref{eq:thetaphiclose}}}},\underbracket{-0.05\frac{\varepsilon}{r}\sin^3{(\overline{\phi}_0+\overline{\theta}_0)}}_{\text{used in \eqref{eq:X1X1BarClose}}}\big\}},\\
    R>&\max{\big\{\underbracket{1700r}_{\mathllap{\text{required by \eqref{eq:etaaproximation}}}},\underbracket{r+\frac{500r}{\sin^2(\overline{\phi}_0+\overline{\theta}_0)}}_{\mathclap{\text{required in \eqref{eq:thetaphiclose}}}},\underbracket{r+\frac{23r^2}{-\varepsilon\sin^3{(\overline{\phi}_0+\overline{\theta}_0)}}}_{\text{used in \eqref{eq:X1X1BarClose}}}\big\}},
\end{aligned} 
\end{equation} then
\begin{equation}\label{eq:X1X1BarClose}
    \begin{aligned}
        \big|\mathcal{X}_1-\overline{\mathcal{X}}_1\bigm|\overbracket{=}^{\mathclap{\eqref{eq:X1X1bardiff}}}&\bigm|r\sin\phi-r\sin{\overline{\phi}_0}+\frac{\xi(\phi,\theta)}{\eta(\phi,\theta)}-\frac{r\cos\phistar-r\cos{\overline{\phi}_0}}{-\sin{(\overline{\phi}_0+\overline{\theta}_0)}}\cos{(\overline{\phi}_0+\overline{\theta}_0)}\big|\\
        \le&r\underbracket{\bigm|\sin\phi-\sin{\overline{\phi}_0}\bigm|}_{\mathllap{\text{mean-value thm: }\le|\phi-\overline{\phi}_0|}}+\underbracket{\bigm|\frac{\xi(\phi,\theta)}{\eta(\phi,\theta)}-\frac{r\cos\phistar-r\cos{\overline{\phi}_0}}{-\sin{(\overline{\phi}_0+\overline{\theta}_0)}}\cos{(\overline{\phi}_0+\overline{\theta}_0)}\bigm|}_{\mathllap{\text{\eqref{eq:xietaratioApproximate}}}}\\
        \overbracket{<}^{{\mathllap{\text{\eqref{eq:xietaratioApproximate}}}}}&r|\phi-\overline{\phi}_0|+\frac{2r}{1.9258\sin^2{(\overline{\phi}_0+\overline{\theta}_0)}}\Big[\frac{7.1r}{b}\frac{1}{-\sin{(\overline{\phi}_0+\overline{\theta}_0)}}+\frac{3|\theta-\overline{\theta}_0|}{-\sin{(\overline{\phi}_0+\overline{\theta}_0)}}+\frac{3|\phi-\overline{\phi}_0|}{-\sin{(\overline{\phi}_0+\overline{\theta}_0)}}\Big]\\
        &+\frac{6r\bigm|\overline{\phi}_0-\phi\bigm|+4r\bigm|\overline{\theta}_0-\theta\bigm|}{-2\sin{(\overline{\phi}_0+\overline{\theta}_0)}}\\
        &=\big(\frac{6r}{-1.9258\sin^3{(\overline{\phi}_0+\overline{\theta}_0)}}+\frac{3r}{-\sin{(\overline{\phi}_0+\overline{\theta}_0)}}+r\big)\big|\phi-\overline{\phi}_0\big|\\
        &+\big(\frac{6r}{-1.9258\sin^3{(\overline{\phi}_0+\overline{\theta}_0)}}+\frac{2r}{-\sin{(\overline{\phi}_0+\overline{\theta}_0)}}+r\big)\big|\theta-\overline{\theta}_0\big|+\frac{14.2r^2}{-1.9258b\sin^3{(\overline{\phi}_0+\overline{\theta}_0)}}\\
        \overbracket{<}^{\mathclap{\frac{1}{-\sin{(\overline{\phi}_0+\overline{\theta}_0)}}\ge1}}&\frac{7.2r}{-\sin^3{(\overline{\phi}_0+\overline{\theta}_0)}}\big|\phi-\overline{\phi}_0\big|+\frac{6.2r}{-\sin^3{(\overline{\phi}_0+\overline{\theta}_0)}}\big|\theta-\overline{\theta}_0\big|+\frac{7.4r^2}{\underbracket{-b\sin^3{(\overline{\phi}_0+\overline{\theta}_0)}}_{\mathllap{\text{\eqref{eq:X1closeCondition}:} b>R-r>\frac{23r^2}{-\varepsilon\sin^3{(\overline{\phi}_0+\overline{\theta}_0)}}}}}\\
        \overbracket{<}^{\mathclap{\eqref{eq:X1closeCondition}}}&\frac{7.2r}{-\sin^3{(\overline{\phi}_0+\overline{\theta}_0)}}\frac{-0.05\varepsilon}{r}\sin^3{(\overline{\phi}_0+\overline{\theta}_0)}+\frac{6.2r}{-\sin^3{(\overline{\phi}_0+\overline{\theta}_0)}}\frac{-0.05\varepsilon}{r}\sin^3{(\overline{\phi}_0+\overline{\theta}_0)}+\frac{7.6\varepsilon}{23}<\varepsilon.
    \end{aligned}
\end{equation}
In the table 
$\LEM(r,R,\phistar)$, we denote by $M$ the midpoint of $\arc{AB}$ in \cref{fig:Trajectory_LemonCompareWithWojtkowski}. If \begin{equation}\label{eq:Y1closecondition}
    R>\frac{2r^2}{\varepsilon}\sin^2\phistar,
\end{equation} then\COMMENT{\red\bfseries Misalignment; use multline\wentao{Done.}}
\begin{equation}\label{eq:Y1Y1barApproximation}
    \begin{multlined}
        \big|\mathcal{Y}_1-\overline{\mathcal{Y}}_1\bigm|=\bigm|\mathcal{Y}_1-(-r\cos\phistar)\bigm|\le\text{distance}(M,|AB|)=2R\sin^2{(\frac{\Phistar}{2})}\\<2R\sin^2{\Phistar}\overbracket{=}^{\mathclap{r\sin\phistar=R\sin\Phistar}}2r\sin{\phistar}\sin{\Phistar}\overbracket{=}^{\qquad\mathllap{r\sin\phistar=R\sin\Phistar}}\frac{2r^2}{R}\sin^2{\phistar}\overbracket{<}^{\mathclap{\text{\eqref{eq:Y1closecondition}}}}\varepsilon.
    \end{multlined}
\end{equation} Since $\Fc(x)=(\Phi_1,\theta_1)$, $P=p(\Fc(x))\in\Gamma_R$, in \cref{fig:Trajectory_LemonCompareWithWojtkowski} the geometric meanings of $\phi$ and $\theta$ imply that $\overrightarrow{TP}$ has angle $2\pi-(\phi+\theta)\in(0,\pi)$ with positive axis $x$ given the clockwise orientation of the billiard plane. By the geometric meanings of $\Phi_1$ the tangent direction at $P$ of $\Gamma_R$ has angle $-\Phi_1$ with the positive axis $x$ given the clockwise orientation of the billiard plane. So, $2\pi+\Phi_1-(\phi+\theta)=\theta_1$ is the incident angle between $\overrightarrow{TP}$ and the tangent direction at $P$. By the geometric meaning of $\Theta_1$ in \cref{def:LemonAndPetTogether}, $\Theta_1-\Phi_1$ is the reflection angle $=\theta_1$. Hence\COMMENT{\small\red\bfseries Lots of wasted space below, redundancy, random alignment\wentao{much rewritten.}} \begin{equation}\label{eq:relectionAndincident}\
    \Phi_1\in\big(-\Phistar, +\Phistar\big),\quad
    \Theta_1-\Phi_1=\theta_1=2\pi+\Phi_1-(\phi+\theta),\quad
    \Theta_1=2\pi+2\Phi_1-(\phi+\theta),
\end{equation} So, if\COMMENT{Fix \eqref{eq:Theta1CloseCondition}. Like many equations it is a randomly arranged jumble and wastes space. Note and fix the problem with the comma---here and everywhere.\wentao{much rewritten.}} \begin{equation}\label{eq:Theta1CloseCondition}
        R>\frac{2.02r}{0.4\varepsilon_1},\,
        |\phi-\overline{\phi}_0|<0.3\varepsilon_1,\,
        |\theta-\overline{\theta}_0|<0.3\varepsilon_1,
\end{equation}
then since $\overbracket{2\pi-\overline{\Theta}_1}^{\mathclap{=\overline{\phi}_0+\overline{\theta}_0\ \eqref{eq:X1barY1barTheta1bar}\quad}}+\overbracket{\Theta_1-\Phi_1}^{\mathrlap{=2\pi+\Phi_1-(\phi+\theta)\text{ by }\ \eqref{eq:relectionAndincident}}}=2\pi+\Phi_1-(\phi+\theta)+\overline{\phi}_0+\overline{\theta}_0$, we get \COMMENT{\small\bfseries\color{red}Do these equations need to be stacked??? They fit on one line but now take up WAY too much space.\wentao{Much rewritten and it should be OK now.}}
\begin{equation}\label{eq:theta1theta1barApproximation}
    \begin{aligned}
        \big|\Theta_1-\overline{\Theta}_1\big|=&\big|2\Phi_1-(\phi+\theta)+\overline{\phi}_0+\overline{\theta}_0\big|\le 2\big|\Phi_1\big|+\big|\phi-\overline{\phi}_0\big|+\big|\theta-\overline{\theta}_0\big|<2|\Phistar|+\big|\phi-\overline{\phi}_0\big|+\big|\theta-\overline{\theta}_0\big|\\
        \overbracket{<}^{\mathclap{r\sin\phistar=R\sin\Phistar}}&2\sin^{-1}\big((r/R)\sin\phistar\big)+\big|\phi-\overline{\phi}_0\big|+\big|\theta-\overline{\theta}_0\big|\overbracket{<}^{\mathclap{\substack{R>1700r: 0<\sin^{-1}\big((r/R)\sin\phistar\big)<0.02\\ \text{If }0<z<0.2\text{, then } \sinc(z)>1/1.01}}}2\cdot1.01\sin{\big(\sin^{-1}{\big((r/R)\sin\phistar\big)}\big)}+\big|\phi-\overline{\phi}_0\big|+\big|\theta-\overline{\theta}_0\big|\\
        <&2.02(r/R)+\big|\phi-\overline{\phi}_0\big|+\big|\theta-\overline{\theta}_0\big|\overbracket{<}^{\mathclap{\eqref{eq:Theta1CloseCondition}}}2.02\frac{0.4\varepsilon_1}{2.02}+0.3\varepsilon_1+0.3\varepsilon_1=\varepsilon_1.
    \end{aligned}
\end{equation}

So far we have got the following.

\eqref{eq:X1closeCondition} implies \eqref{eq:X1X1BarClose}, which is the first inequality of \eqref{eq:eq:x0barx0Closing}. 

\eqref{eq:Y1closecondition} implies \eqref{eq:Y1Y1barApproximation}, which is the second inequality of \eqref{eq:eq:x0barx0Closing}.

\eqref{eq:Theta1CloseCondition} implies \eqref{eq:theta1theta1barApproximation}, which is the third inequality of \eqref{eq:eq:x0barx0Closing}. 

Since \eqref{eq:x0barx0ClosingCondition} implies \eqref{eq:X1closeCondition}\eqref{eq:Y1closecondition}\eqref{eq:Theta1CloseCondition}, we have therefore shown that \eqref{eq:x0barx0ClosingCondition} implies \eqref{eq:eq:x0barx0Closing}. 

We now show that \eqref{eq:xorbitA0Condition} implies $\sin\theta_1>\sqrt{4r/R}$. Note that if\COMMENT{\red\bfseries\large NOW EXPLAIN WHERE \eqref{eq:xorbitA0Condition} is used.\wentao{It's explained in the end of the proof.} Not good enough. I fixed it.}
\begin{equation}\label{eq:theta1a0condition1}
    \overline{\phi}_0-\phi>\overbracket{-(2\pi-\overline{\phi}_0-\overline{\theta}_0)/4}^{\text{\cref{remark:sinPhi0plusTheta0isNegative}: }<0},\quad 
    \overline{\theta}_0-\theta>\overbracket{-(2\pi-\overline{\phi}_0-\overline{\theta}_0)/4}^{\text{\cref{remark:sinPhi0plusTheta0isNegative}: }<0},\quad R>4r/\sin^2{(\overbracket{\frac{\pi}{2}-\frac{\overline{\phi}_0+\overline{\theta}_0}{4}}^{\mathclap{\text{\cref{remark:sinPhi0plusTheta0isNegative}: }\in(0,\pi/4)})}},
\end{equation} then
\begin{equation}\label{eq:theta1a0bound1}
    \begin{aligned}
        \theta_1-\sin^{-1}{(\sqrt{4r/R})}\overbracket{=}^{\mathclap{\text{\eqref{eq:relectionAndincident}}}}&2\pi+\Phi_1-(\phi+\theta)-\sin^{-1}{(\sqrt{4r/R})}
        =2\pi-(\overline{\phi}_0+\overline{\theta}_0)+(\overline{\phi}_0-\phi)+(\overline{\theta}_0-\theta)+\overbracket{\Phi_1}^{\mathclap{>-\Phistar=-\sin^{-1}{\big((r/R)\sin{\phistar}\big)}}}-\sin^{-1}{(\sqrt{4r/R})}\\
        >&2\pi-(\overline{\phi}_0+\overline{\theta}_0)+\underbracket{(\overline{\phi}_0-\phi)}_{\mathclap{\text{\eqref{eq:theta1a0bound1}: }>-(2\pi-\overline{\phi}_0-\overline{\theta}_0)/4}}+\overbracket{(\overline{\theta}_0-\theta)}^{\mathclap{\text{\eqref{eq:theta1a0bound1}: }>-(2\pi-\overline{\phi}_0-\overline{\theta}_0)/4}}\underbracket{-\sin^{-1}{\big((r/R)\sin{\phistar}\big)}}_{\mathclap{>-\sin^{-1}{\sqrt{4r/R}}}}-\sin^{-1}{(\sqrt{4r/R})}\\
        >&\pi-\frac{\overline{\phi}_0+\overline{\theta}_0}{2}-2\sin^{-1}(\sqrt{4r/R})\underbracket{>}_{\mathclap{\text{\eqref{eq:theta1a0bound1}: }R>4r/\sin^2{(\frac{\pi}{2}-\frac{\overline{\phi}_0+\overline{\theta}_0}{4})}}}0.
    \end{aligned}
\end{equation}

If \begin{equation}\label{eq:theta1a0condition2}
    \begin{aligned}
        -\overline{\phi}_0+\phi>-\frac{-\pi+\overline{\phi}_0+\overline{\theta}_0}{4},\quad
        -\overline{\theta}_0+\theta>-\frac{-\pi+\overline{\phi}_0+\overline{\theta}_0}{4},\quad
        R>\frac{4r}{\sin^2{(\frac{\overline{\phi}_0+\overline{\theta}_0}{4}-\frac{\pi}{4})}},
    \end{aligned}
\end{equation} then we get
\begin{equation}\label{eq:theta1a0bound2}
    \begin{aligned}
        \pi-\theta_1-\sin^{-1}{(\sqrt{\frac{4r}{R}})}\overbracket{=}^{\mathrlap{\text{\eqref{eq:relectionAndincident}: }\theta_1=2\pi+\Phi_1-(\phi+\theta)}}&-\pi+\Phi_1-(\phi+\theta)-\sin^{-1}{(\sqrt{\frac{4r}{R}})}=-\pi+(\overline{\phi}_0+\overline{\theta}_0)+(-\overline{\phi}_0+\phi)+(-\overline{\theta}_0+\theta)\overbracket{-\Phi_1}^{\mathllap{>-\Phistar=-\sin^{-1}{\big((r/R)\sin{\phistar}\big)}}}-\sin^{-1}{(\sqrt{\frac{4r}{R}})}\\
        >&-\pi+(\overline{\phi}_0+\overline{\theta}_0)+\underbracket{(-\overline{\phi}_0+\phi)}_{\mathclap{\text{\eqref{eq:theta1a0condition2}: }>-\frac{-\pi+\overline{\phi}_0+\overline{\theta}_0}{4}}}+\overbracket{(-\overline{\theta}_0+\theta)}^{\mathclap{\text{\eqref{eq:theta1a0condition2}: }>-\frac{-\pi+\overline{\phi}_0+\overline{\theta}_0}{4}}}\underbracket{-\sin^{-1}{\big((r/R)\sin{\phistar}\big)}}_{\mathclap{>-\sin^{-1}{\sqrt{4r/R}}}}-\sin^{-1}{(\sqrt{4r/R})}\\
        >&-\frac{\pi}{2}+\frac{\overline{\phi}_0+\overline{\theta}_0}{2}-2\sin^{-1}(\sqrt{4r/R})\underbracket{>}_{\mathclap{\text{\eqref{eq:theta1a0condition2}: }R>4r/\sin^2{(\frac{\overline{\phi}_0+\overline{\theta}_0}{4}-\frac{\pi}{4})}}}0.
    \end{aligned}
\end{equation}
Note that \eqref{eq:xorbitA0Condition} implies \eqref{eq:theta1a0condition1} and \eqref{eq:theta1a0condition2}. 

We also note that $\eqref{eq:theta1a0condition1} \implies \eqref{eq:theta1a0bound1}\implies \theta_1>\sin^{-1}{\sqrt{4r/R}}$ and $\eqref{eq:theta1a0condition2} \implies \eqref{eq:theta1a0bound2} \implies \theta_1<\pi-\sin^{-1}{\sqrt{4r/R}}$. Therefore, \eqref{eq:xorbitA0Condition} implies $\sin\theta_1>\sqrt{4r/R}$. This concludes the proof of \cref{lemma:Trajectory_approximation_1}.
\end{proof}

\begin{figure}[h]
\begin{center}
\begin{tikzpicture}[xscale=0.6,yscale=0.6]
    \tkzDefPoint(0,0){Or};
    \tkzDefPoint(0,-4.8){Y};
    \pgfmathsetmacro{\rradius}{4.5};
    \clip(-1.1*\rradius,-1.1*\rradius) rectangle (2.8*\rradius, 0.02*\rradius);
    \pgfmathsetmacro{\phistardeg}{40};
    \pgfmathsetmacro{\XValueArc}{\rradius*sin(\phistardeg)};
    \pgfmathsetmacro{\YValueArc}{\rradius*cos(\phistardeg)};
    \pgfmathsetmacro{\Rradius}{
        \rradius*sin(\phistardeg)/sin(18)
    };
    \pgfmathsetmacro{\bdist}{(-\rradius*cos(\phistardeg))+cos(18)*\Rradius}
    \tkzDefPoint(0,\bdist){OR};
    \tkzDefPoint(\XValueArc,-1.0*\YValueArc){B};    \tkzDefPoint(-1.0*\XValueArc,-1.0*\YValueArc){A};
    \tkzDrawArc[name path=Cr, thick](Or,B)(A);
    \tkzDrawArc[name path=CR, ultra thin](OR,A)(B);
    \pgfmathsetmacro{\PXValue}{-\Rradius*sin(3)};
    \pgfmathsetmacro{\PYValue}{\bdist-\Rradius*cos(3)};
    \tkzDefPoint(\PXValue,\PYValue){P};
    



    \path[name path=Line_P_Q0_far](P)-- (3:8);
    \path[name path=Line_P_Q_far](P)-- (1:8);
    \path[name intersections={of=Cr and Line_P_Q0_far,by={Q0}}];
    \path[name intersections={of=Cr and Line_P_Q_far,by={Q}}];
    \tkzDefPoint(0,\bdist-\Rradius){C};
    \path[name intersections={of=Cr and Line_P_Q0_far,by={Q0}}];
    \draw[name path=LineAB,thin] (A)--(B);
    \coordinate (PHorizon) at ([shift={(3,0)}]P);
    \draw[red,ultra thin,dashed] (P)--(PHorizon);
    
    \path[name intersections={of=Line_P_Q0_far and LineAB,by={P_bar}}];
    

    \tkzLabelPoint[below](A){$A$};
    \tkzLabelPoint[below](P){\small $p(x_1)=P$};
    \tkzLabelPoint[below](Or){$O_r$};
    \tkzLabelPoint[below](B){$B$};
    \tkzLabelPoint[above left]([shift={(-0.01,0.02)}]P_bar){\small $p(\overline{x}_1)=\overline{P}$};
    \tkzLabelPoint[right]([shift={(0.42,0.3)}]P_bar){\small $\overline{\Theta}_1$};
    \tkzLabelPoint[right]([shift={(0.2,0.13)}]P){\small $\Theta_1$};
    
    \tkzDefPointBy[projection=onto B--A](T0)\tkzGetPoint{M0};
    \coordinate (T0Mirror) at ($(T0)!2!(M0)$);
    \coordinate (PbarArrow) at ($(T0Mirror)!1.06!(P_bar)$);
    \tkzDefPointBy[projection=onto P--OR](T)\tkzGetPoint{M};
    \coordinate (TMirror) at ($(T)!2!(M)$);
    \coordinate (PArrow) at ($(P)!0.07!(TMirror)$);

    

    \tkzMarkAngle[arc=lll,ultra thin,size=0.55](B,P_bar,Q0);
    \tkzMarkAngle[arc=lll,ultra thin,size=0.22](PHorizon,P,Q);
    \tkzDefPoint(0,\bdist-\Rradius){C};
    \begin{scope}[ultra thin,decoration={markings,mark=at position 0.7 with {\arrow{>}}}] 
    \draw[red,postaction={decorate},ultra thin] (P_bar)--(Q0);
    \draw[blue,postaction={decorate},ultra thin] (P)--(Q);
    \end{scope}
    \tkzLabelPoint[above right](Q0){$\overline{Q}=p(\overline{x}_2)=(r\sin{\overline{\phi}_2},-r\cos{\overline{\phi}_2})$};
    \tkzLabelPoint[below right](Q){$Q=p(x_2)=(r\sin{\phi_2},-r\cos{\phi_2})$};
    \tkzDrawPoints(A,B,P,P_bar,Q0,Q,Or);
\end{tikzpicture}
\caption{
$\overline{P}=(\overline{\mathcal{X}}_1,\overline{\mathcal{Y}}_1)=p(\overline{x}_1)$ and $\overline{x}_2=\bar F_{\phistar}(\overline{x}_1),\overline{Q}=p(\overline{x}_2)$.\newline $P=(\mathcal{X}_1,\mathcal{Y}_1)=p(x_1)$, and $x_2=\Fc(x_1),Q=p(x_2)$
} \label{fig:Trajectory_LemonCompareWithWojtkowski2}
\end{center}
\end{figure}\hfill

We now are going to prove a second trajectory approximation lemma: \cref{lemma:Trajectory_approximation_2}. \cref{lemma:Trajectory_approximation_2} supposes in \cref{fig:Trajectory_LemonCompareWithWojtkowski2} one trajectory starting on $\Gamma_R$ of lemon billiard $\mathbb{L}(r,R,\phistar)$ and another trajectory starting on the flat boundary $\overline{AB}$ of the petal billiard $\mathbb{P}(r,\phistar)$. If these 2 trajectories start from close enough positions and with close enough directions, and $R$ is large enough, then these 2 trajectories respective second collisions on $\Gamma_r$ for the boundary shared by lemon billiard $\mathbb{L}(r,R,\phistar)$ and the petal billiard $\mathbb{P}(r,\phistar)$, also have close collision positions and angles (see \cref{fig:Trajectory_LemonCompareWithWojtkowski2}).\COMMENT{I edited the text before \cref{lemma:Trajectory_approximation_1}; edit this text likewise.}

\begin{lemma}[Shadowing 2, trajectory approximation in \cref{fig:Trajectory_LemonCompareWithWojtkowski2}]\label{lemma:Trajectory_approximation_2}\hfill

Using \cref{def:LemonAndPetTogether,def:flowcollsionpositionandangle} we suppose that $x_1\in \MRout\subset M_R$, $(\phi_2,\theta_2)=x_2=\Fcal(x_1)\in\Mrin$, $p(x_1)=(\mathcal{X}_1,\mathcal{Y}_1)$ and $(\mathcal{X}_1,\mathcal{Y}_1,\Theta_1)$ represents the position and direction of the $\LEM(r,R,\phistar)$ billiard flow collision on the boundary $\Gamma_R$. We suppose $\overline{P}=p(\overline{x}_1)=(\overline{\mathcal{X}}_1,\overline{\mathcal{Y}}_1)\in\overline{AB}$. i.e. $\big|\overline{\mathcal{X}}_1\big|\le r\sin\phistar$ and $(\overline{\phi}_2,\overline{\theta}_2)=\overline{x}_2\nfd\bar F_{\phistar,1}(\overline{x}_1)\in$ closure$(\Mrin)$ and $(\overline{\mathcal{X}}_1,\overline{\mathcal{Y}}_1,\overline{\Theta}_1)$ represents the position and direction of the $\PET(r,\phistar)$ billiard flow after collision on the boundary $\overline{AB}$  (in \cref{def:flowcollsionpositionandangle}) of $\overline{x}_1\in M_{\text{\upshape{f}}}$.



 Let $\overline{\constant}_+\dfn\max\big\{r\sin{\phistar}-\overline{\mathcal{X}}_1,r|\sin{(\phistar-\overline{\Theta}_1)}|\big\}$, $\overline{\constant}_-\dfn\max\big\{r\sin{\phistar}+\overline{\mathcal{X}}_1,r|\sin{(\phistar+\overline{\Theta}_1)}|\big\}$, $\overline{\constant}_0\dfn r^2\sin^2\phistar-\overline{\mathcal{X}}^2_1$ and 
 \begin{equation}\label{eq:APetalCollisionConstant}
     \overline{\constant}\dfn\max\Big\{\min\big\{\overline{\constant}^2_+,\overline{\constant}^2_-\big\}\cdot\chi_{\big\{|\overline{\mathcal{X}}_1|= r\sin\phistar\big\}},\overline{\constant}_0\Big\},
 \end{equation} where $\chi_{\big\{|\overline{\mathcal{X}}_1|= r\sin\phistar\big\}}$ is the indicator function. Note that $\overline{\constant}_-$, $\overline{\constant}_+$, $\overline{\constant}_0$, $\overline{\constant}$ depend on the $\overline{\mathcal{X}}_1$, $\overline{\Theta}_1$ the collision position and direction of $\overline{x}_1\in M_{\text{\upshape f}}$ on the boundary $\overline{AB}$ of $\PET(r,\phistar)$ and the parameter $\phistar$.

We have the following conclusion. For $\forall\varepsilon>0$ if \begin{equation}\label{eq:X1Y1Theta1CloseCondition}
    \begin{aligned}
        \big|\mathcal{X}_1-\overline{\mathcal{X}}_1\bigm|<\delta_3(\varepsilon,\overline{x}_1)\dfn&\left\{\begin{aligned}
        &\min\big\{\frac{0.3\sqrt{\overline{\constant}}}{6r+3\sqrt{\overline{\constant}}}\frac{1.4r}{3\pi}\varepsilon,\quad\frac{0.2\varepsilon}{3\pi}r,\quad\frac{0.015\sqrt{\overline{\constant}}r\sin{(2\phistar)}}{6r+3\sqrt{\overline{\constant}}},\quad 0.05r\sin{(2\phistar)}\big\}\text{\upshape, If }\overline{\constant}>0.\\
        &\min\big\{0.09(\frac{1.4}{3\pi})^2r\varepsilon^2,\quad\frac{0.2\varepsilon}{3\pi}r,\quad 0.000225r\sin^2{(2\phistar)}\big\}\text{\upshape , If }\overline{\constant}=0.
        \end{aligned}\right.\\
        \big|\mathcal{Y}_1-\overbracket{\overline{\mathcal{Y}}_1}^{\mathrlap{=-r\cos\phistar}}\bigm|<\delta_4(\varepsilon,\overline{x}_1)\dfn&\left\{\begin{aligned}
        &\min\big\{\frac{0.3\sqrt{\overline{\constant}}}{6r+3\sqrt{\overline{\constant}}}\frac{1.4r}{3\pi}\varepsilon,\quad\frac{0.2\varepsilon}{3\pi}r,\quad\frac{0.015\sqrt{\overline{\constant}}r\sin{(2\phistar)}}{6r+3\sqrt{\overline{\constant}}},\quad 0.05r\sin{(2\phistar)}\big\} \text{\upshape , If }\overline{\constant}>0.\\
            &\min\big\{0.045(\frac{1.4}{3\pi})^2r\varepsilon^2,\quad\frac{0.2\varepsilon}{3\pi}r,\quad \frac{0.000225}{2}r\sin^2{(2\phistar)}\big\}\text{\upshape , If }\overline{\constant}=0.\\
        \end{aligned}\right.\\
        \big|\Theta_1-\overline{\Theta}_1\bigm|<\delta_5(\varepsilon,\overline{x}_1)\dfn&\left\{\begin{aligned}
            &\min\big\{\frac{0.4\sqrt{\overline{\constant}}}{8r^2+2r\sqrt{\overline{\constant}}}\frac{1.4r}{3\pi}\varepsilon,\quad\frac{0.1\varepsilon}{4\pi+1},\quad\frac{0.02\sqrt{\overline{\constant}}r\sin{(2\phistar)}}{8r^2+2r\sqrt{\overline{\constant}}},\quad 0.05\sin{(2\phistar)}\big\} \text{\upshape , If }\overline{\constant}>0.\\
            &\min\big\{\frac{0.14\varepsilon}{3\pi},\quad\frac{0.1\varepsilon}{4\pi+1},\quad 0.005\sin{(2\phistar)}\big\}\text{\upshape , If }\overline{\constant}=0.\\
        \end{aligned}\right.
    \end{aligned}
\end{equation} then \begin{equation}\label{eq:phi2theta2close}
    \begin{aligned}
        |\overline{\phi}_2-\phi_2|<\varepsilon,\quad
        |\overline{\theta}_2-\theta_2|<\varepsilon.
    \end{aligned}
\end{equation}
\end{lemma}
\begin{proof}
We begin with a sufficient condition for $\big|t-\overline{t}\big|<\varepsilon$, \eqref{eq:ttbarvarepsilonCondition} (on page \pageref{eq:ttbarvarepsilonCondition}) below.

In \cref{fig:Trajectory_LemonCompareWithWojtkowski2} with coordinate system of \cref{def:StandardCoordinateTable}, $O_r=(0,0)$, $O_R=(0,b)$, $Q=p(x_2)=(r\sin{\phi_2},-r\cos{\phi_2})=(x_Q,y_Q)$, $\overline{Q}=(r\sin{\overline{\phi}_2},-r\cos{\overline{\phi}_2})=(x_{\overline{Q}},y_{\overline{Q}})$.
    
Write $\overline{t}=|\overline{P}\,\overline{Q}|$ and $t=|PQ|$. Note that the vector $\overrightarrow{\overline{P}\,\overline{Q}}$ has an angle $\overline{\Theta}_1$ with respect to the positive axis $x$, the vector $\overrightarrow{PQ}$ has an angle $\Theta_1$ with respect to the positive axis $x$ and $Q$, $\overline{Q}$ are on the circle $C_r$. We get the following.\COMMENT{Unhelpful alignment below}
\begin{equation}\label{eq:tau1lengthQuadraticEq}
    \begin{aligned}
        &\left.
        \begin{aligned}
            x_Q=\mathcal{X}_1+t\cos{\Theta_1}&=r\sin\phi_2\\
            y_Q=\mathcal{Y}_1+t\sin{\Theta_1}&=-r\cos\phi_2\\
        \end{aligned}\right\}\Longrightarrow(\overbracket{\mathcal{X}_1+t\cos{\Theta_1}}^{\mathllap{=r\sin{\phi_2}}})^2+(\overbracket{\mathcal{Y}_1+t\sin{\Theta_1}}^{=-r\cos{\phi_2}})^2=r^2\\
        &\Longrightarrow t^2+(2\mathcal{X}_1\cos{\Theta_1}+2\mathcal{Y}_1\sin{\Theta_1})t+\mathcal{X}_1^2+\mathcal{Y}_1^2-r^2=0.\\
    \end{aligned}
\end{equation} 
Since $(\mathcal{X}_1,\mathcal{Y}_1)=P$ is on the circle $C_R$ centered at $(0,b)$ and $\measuredangle{O_RO_rB}=\pi-\phistar$,\COMMENT{\Large\red\bfseries Break up alignment below. It is completely misaligned, and ``Therefore'' must be text rather than display math.\wentao{Much rewritten.}} by elementary geometry and Cosine law\COMMENT{Never capitalize Cosine law except at the start of a sentence. Check entire file for occurrences.} we ge the following. \begin{equation}\label{eq:tau1lengthQuadraticEq2}
    \begin{aligned}
        \left.
        \begin{aligned}
            \mathcal{X}^2_1+\mathcal{Y}^2_1-2b\mathcal{Y}_1+b^2=\mathcal{X}_1^2+(\mathcal{Y}_1-b)^2&=R^2\\
            b^2+r^2-2br\cos{(\angle{O_RO_rB})}=b^2+r^2+2br\cos\phistar&=R^2\\
        \end{aligned}
        \right\}\Rightarrow&\mathcal{X}^2_1+\mathcal{Y}^2_1-2b\mathcal{Y}_1+b^2=b^2+r^2+2br\cos\phistar\\
        \Rightarrow&\mathcal{X}^2_1+\mathcal{Y}^2_1-r^2=2b\mathcal{Y}_1+2br\cos\phistar
    \end{aligned}
\end{equation}

Therefore, the last quadratic equation in \eqref{eq:tau1lengthQuadraticEq} becomes
$t^2+(2\mathcal{X}_1\cos{\Theta_1}+2\mathcal{Y}_1\sin{\Theta_1})t+\underbracket{2b\mathcal{Y}_1+2br\cos\phistar}_{=\mathcal{X}_1^2+\mathcal{Y}_1^2-r^2}=0$.

Since $2b(\mathcal{Y}_1+r\cos\phistar)=\mathcal{X}_1^2+\mathcal{Y}_1^2-r^2<0$, the last quadratic equation in \eqref{eq:tau1lengthQuadraticEq2} has one positive root and one negative root. $|PQ|=t$ is its positive root.\COMMENT{\red\bfseries Also break this alignment up. No ``Where'' in display. Missing punctuation again. \eqref{eq:tau1barlengthQuadraticEq} also needs fixing.\wentao{It seems good now after editing?}} We get the following.
\begin{equation}\label{eq:tplus}
        |PQ|=t=-\mathcal{X}_1\cos{\Theta_1}-\mathcal{Y}_1\sin{\Theta_1}+\sqrt{\Delta},\quad
        \Delta=(\mathcal{X}_1\cos{\Theta_1}+\mathcal{Y}_1\sin{\Theta_1})^2-2b(\underbracket{\mathcal{Y}_1+r\cos{\phistar}}_{<0})>0,
\end{equation}
\begin{equation}\label{eq:tau1barlengthQuadraticEq}
    \begin{aligned}
        &\left.
        \begin{aligned}
            x_{\overline{Q}}=\overline{\mathcal{X}}_1+\overline{t}\cos{\overline{\Theta}_1}&=r\sin\overline{\phi}_2\\
            y_{\overline{Q}}=\overline{\mathcal{Y}}_1+\overline{t}\sin{\overline{\Theta}_1}&=-r\cos\overline{\phi}_2\\
            \overline{\mathcal{Y}}_1&=-r\cos{\phistar}\end{aligned}\right\}\Longrightarrow(\overbracket{\overline{\mathcal{X}}_1+\overline{t}\cos{\overline{\Theta}_1}}^{\mathllap{=r\sin{\overline{\phi}_2}}})^2+(\overbracket{-r\cos\phistar+\overline{t}\sin{\overline{\Theta}_1}}^{=-r\cos{\overline{\phi}_2}})^2=r^2\\
        \Longrightarrow& \overline{t}^2+(2\overline{\mathcal{X}}_1\cos{\overline{\Theta}_1}-2r\cos{\phistar}\sin{\overline{\Theta}_1})\overline{t}+\overline{\mathcal{X}}_1^2-r^2\sin^2{\phistar}=0.
    \end{aligned}
\end{equation}
Since $\overline{\mathcal{X}}_1^2-r^2\sin^2{\phistar}\le0$, the last quadratic equation in \eqref{eq:tau1barlengthQuadraticEq} has one nonnegative root: $\overline{t}=|\overline{P}\,\overline{Q}|$. Therefore, \begin{equation}\label{eq:tau1barlengthQuadraticEq2}
    \begin{aligned}
    |\overline{P}\,\overline{Q}|=\overline{t}=-\overline{\mathcal{X}}_1\cos{\overline{\Theta}_1}+r\cos\phistar\sin{\overline{\Theta}_1}+\sqrt{\overline{\Delta}},
    \end{aligned}
\end{equation}
where $\overline{\Delta}=(\overline{\mathcal{X}}_1\cos{\overline{\Theta}_1}-r\cos\phistar\sin{\overline{\Theta}_1})^2-(\overline{\mathcal{X}}^2_1-r^2\sin^2{\phistar})\ge0$, since $\overline{\mathcal{X}}_1\in[-r\sin{\phistar}, r\sin{\phistar}]$. 

To control the bound for $|t-\overline{t}|$ in \eqref{eq:DeltaCompareIneq7},\COMMENT{The reader has not yet seen \eqref{eq:DeltaCompareIneq7}, so edit} we will prove and use the following six inequalities: \cref{eq:DeltaCompareIneq1,eq:DeltaCompareIneq2,eq:DeltaCompareIneq3,eq:DeltaCompareIneq4,eq:DeltaCompareIneq5,eq:DeltaCompareIneq6_3}.\COMMENT{Those equations are missing punctuation}
\COMMENT{\large\color{red}Explain what is actually done and how and where\Hrule Almost everything below (the rest of the subsection) needs fixing.\wentao{More detailed info is provided.}}\COMMENT{\red\small\bfseries Virtually every one of these equations still has bad alignments.}
\begin{equation}\label{eq:DeltaCompareIneq1}
    \begin{aligned}
        |2b(\mathcal{Y}_1+r\cos{\phistar})-(\mathcal{X}^2_1-r^2\sin^2{\phistar})|=&|2b\mathcal{Y}_1-\mathcal{X}^2_1+2br\cos{\phistar}+r^2\sin^2{\phistar}|\\
        \overbracket{=}^{\mathllap{\eqref{eq:tau1lengthQuadraticEq2}:\mathcal{X}^2_1+\mathcal{Y}^2_1-2b\mathcal{Y}_1+b^2=R^2}}&|\mathcal{Y}^2_1+b^2-R^2+2br\cos{\phistar}+r^2\sin^2{\phistar}|\\
        \overbracket{=}^{\mathllap{b^2+r^2+2br\cos\phistar=R^2}}&|\mathcal{Y}^2_1-r^2+r^2\sin^2{\phistar}|=|\mathcal{Y}^2_1-r^2\cos^2{\phistar}|\\
        =&|\mathcal{Y}_1-r\cos{\phistar}||\mathcal{Y}_1+r\cos{\phistar}|\overbracket{<}^{\mathllap{|\mathcal{Y}_1-r\cos{\phistar}|<|\mathcal{Y}_1|+r<2r}}2r|\mathcal{Y}_1+r\cos{\phistar}|
    \end{aligned}
\end{equation}
\begin{equation}\label{eq:DeltaCompareIneq2}
    \begin{aligned}
        &|\mathcal{X}^2_1-r^2\sin^2{\phistar}-(\overline{\mathcal{X}}^2_1-r^2\sin^2{\phistar})|=|\mathcal{X}^2_1-\overline{\mathcal{X}}^2_1|=|\overline{\mathcal{X}}_1-\mathcal{X}_1|\overbracket{|\overline{\mathcal{X}}_1+\mathcal{X}_1|}^{\mathllap{|\overline{\mathcal{X}}_1+\mathcal{X}_1|<|\overline{\mathcal{X}}_1|+|\mathcal{X}_1|<2r}}<2r|\overline{\mathcal{X}}_1-\mathcal{X}_1|\\
    \end{aligned}
\end{equation}
\begin{equation}\label{eq:DeltaCompareIneq3}
    \begin{aligned}
        &\big|2b(\mathcal{Y}_1+r\cos{\phistar})-(\overline{\mathcal{X}}^2_1-r^2\sin^2{\phistar})\bigm|\\
        =&\bigm|2b(\mathcal{Y}_1+r\cos{\phistar})-(\mathcal{X}^2_1-r^2\sin^2{\phistar})+\mathcal{X}^2_1-r^2\sin^2{\phistar}-(\overline{\mathcal{X}}^2_1-r^2\sin^2{\phistar})\bigm|\\
        \le&\bigm|2b(\mathcal{Y}_1+r\cos{\phistar})-(\mathcal{X}^2_1-r^2\sin^2{\phistar})\bigm|+\bigm|\mathcal{X}^2_1-r^2\sin^2{\phistar}-(\overline{\mathcal{X}}^2_1-r^2\sin^2{\phistar})\bigm|\\
        \overbracket{<}^{\mathllap{\eqref{eq:DeltaCompareIneq1},\eqref{eq:DeltaCompareIneq2}}}&2r\bigm|\mathcal{Y}_1+r\cos{\phistar}\bigm|+2r\bigm|\overline{\mathcal{X}}_1-\mathcal{X}_1\bigm|\\
    \end{aligned}
\end{equation}
\begin{equation}\label{eq:DeltaCompareIneq4}
    \begin{aligned}
    &\bigm|\overline{\mathcal{X}}_1\cos{\overline{\Theta}_1}-r\cos{\phistar}\sin{\overline{\Theta}_1}-(\mathcal{X}_1\cos{\Theta_1}+\mathcal{Y}_1\sin{\Theta_1})\bigm|\\
    =&\big|\overline{\mathcal{X}}_1\cos{\overline{\Theta}_1}-\mathcal{X}_1\cos{\overline{\Theta}_1}+\mathcal{X}_1\cos{\overline{\Theta}_1}-\mathcal{X}_1\cos{\Theta_1}-r\cos{\phistar}\sin{\overline{\Theta}_1}-\mathcal{Y}_1\sin{\overline{\Theta}_1}+\mathcal{Y}_1\sin{\overline{\Theta}_1}-\mathcal{Y}_1\sin{\Theta_1}\bigm|\\
    \le& |\overline{\mathcal{X}}_1\cos{\overline{\Theta}_1}-\mathcal{X}_1\cos{\overline{\Theta}_1}|+\underbracket{|\mathcal{X}_1\cos{\overline{\Theta}_1}-\mathcal{X}_1\cos{\Theta_1}|}_{\qquad\qquad\mathllap{\text{By Mean-Value Theorem: }\le|\mathcal{X}_1||\overline{\Theta}_1-\Theta_1|< r|\overline{\Theta}_1-\Theta_1|}}+|-r\cos{\phistar}\sin{\overline{\Theta}_1}-\mathcal{Y}_1\sin{\overline{\Theta}_1}|+\underbracket{|\mathcal{Y}_1\sin{\overline{\Theta}_1}-\mathcal{Y}_1\sin{\Theta_1}|}_{\qquad\qquad\mathllap{\text{By Mean-Value Theorem: }\le|\mathcal{Y}_1||\overline{\Theta}_1-\Theta_1|<r|\overline{\Theta}_1-\Theta_1|}}\\
    \underbracket{<}_{\mathrlap{\substack{\text{By Mean-Value Theorem and }\\ |\mathcal{X}_1|<r, |\mathcal{Y}_1|<r}}}&\bigm|\overline{\mathcal{X}}_1-\mathcal{X}_1\bigm|+\bigm|r\cos{\phistar}+\mathcal{Y}_1\bigm|+2r\bigm|\Theta_1-\overline{\Theta}_1\bigm|
    \end{aligned}
\end{equation}

\begin{equation}\label{eq:DeltaCompareIneq5}
    \begin{aligned}
        &\big|(\overline{\mathcal{X}}_1\cos{\overline{\Theta}_1}-r\cos{\phistar}\sin{\overline{\Theta}_1})^2-(\mathcal{X}_1\cos{\Theta_1}+\mathcal{Y}_1\sin{\Theta_1})^2\big|\\
        =&\underbracket{\big|\overline{\mathcal{X}}_1\cos{\overline{\Theta}_1}-r\cos{\phistar}\sin{\overline{\Theta}_1}+\mathcal{X}_1\cos{\Theta_1}+\mathcal{Y}_1\sin{\Theta_1}\bigm|}_{<4r\text{ since } |\overline{\mathcal{X}}_1|<r,|\mathcal{X}_1|<r,|\mathcal{Y}_1|<r}\bigm|\overline{\mathcal{X}}_1\cos{\overline{\Theta}_1}-r\cos{\phistar}\sin{\overline{\Theta}_1}-(\mathcal{X}_1\cos{\Theta_1}+\mathcal{Y}_1\sin{\Theta_1})\big|\\<&4r\bigm|\overline{\mathcal{X}}_1\cos{\overline{\Theta}_1}-r\cos{\phistar}\sin{\overline{\Theta}_1}-(\mathcal{X}_1\cos{\Theta_1}+\mathcal{Y}_1\sin{\Theta_1})\big|\\
        \underbracket{<}_{\mathrlap{\eqref{eq:DeltaCompareIneq4}}}&4r\Big(\big|\overline{\mathcal{X}}_1-\mathcal{X}_1\bigm|+\bigm|r\cos{\phistar}+\mathcal{Y}_1\bigm|+2r\bigm|\Theta_1-\overline{\Theta}_1\bigm|\Big)
    \end{aligned}
\end{equation}

\begin{equation}\label{eq:deltadiff}
    \begin{aligned}
        \big|\Delta-\overline{\Delta}\bigm|\overbracket{=}^{\mathrlap{\text{\eqref{eq:tplus}, \eqref{eq:tau1barlengthQuadraticEq2}}}}&\bigm|(\mathcal{X}_1\cos{\Theta_1}+\mathcal{Y}_1\sin{\Theta_1})^2-2b(\mathcal{Y}_1+r\cos{\phistar})-\big[(\overline{\mathcal{X}}_1\cos{\overline{\Theta}_1}-r\cos\phistar\sin{\overline{\Theta}_1})^2-(\overline{\mathcal{X}}^2_1-r^2\sin^2{\phistar})\big]\big|\\
        \le&\bigm|(\mathcal{X}_1\cos{\Theta_1}+\mathcal{Y}_1\sin{\Theta_1})^2-(\overline{\mathcal{X}}_1\cos{\overline{\Theta}_1}-r\cos\phistar\sin{\overline{\Theta}_1})^2\bigm|+\bigm|-2b(\mathcal{Y}_1+r\cos{\phistar})+(\overline{\mathcal{X}}^2_1-r^2\sin^2{\phistar})\big|\\
        \overbracket{<}^{\mathllap{\eqref{eq:DeltaCompareIneq5},\eqref{eq:DeltaCompareIneq3}}}&4r\Big(\bigm|\overline{\mathcal{X}}_1-\mathcal{X}_1\bigm|+\bigm|r\cos{\phistar}+\mathcal{Y}_1\bigm|+2r\bigm|\Theta_1-\overline{\Theta}_1\bigm|\Big)+2r\bigm|\mathcal{Y}_1+r\cos{\phistar}\bigm|+2r\bigm|\overline{\mathcal{X}}_1-\mathcal{X}_1\big|\\
        =&2r\big(3\bigm|\overline{\mathcal{X}}_1-\mathcal{X}_1\bigm|+3\bigm|\mathcal{Y}_1+r\cos{\phistar}\bigm|+4r\bigm|\Theta_1-\overline{\Theta}_1\bigm|\big).
    \end{aligned}
\end{equation}

Since $\overline{\constant}_+=\max\big\{r\sin{\phistar}-\overline{\mathcal{X}}_1,r|\sin{(\phistar-\overline{\Theta}_1)}|\big\}$, $\overline{\constant}_-=\max\big\{r\sin{\phistar}+\overline{\mathcal{X}}_1,r|\sin{(\phistar+\overline{\Theta}_1)}|\big\}$ and $\overline{\constant}_0=r^2\sin^2\phistar-\overline{\mathcal{X}}^2_1$, $\overline{\constant}=\max\Big\{\min\big\{\overline{\constant}^2_+,\overline{\constant}^2_-\big\}\cdot\chi_{\big\{|\overline{\mathcal{X}}_1|= r\sin\phistar\big\}},\overline{\constant}_0\Big\}$, we get the following. \newline If $|\overline{\mathcal{X}}_1|<r\sin\phistar$, then \eqref{eq:tau1barlengthQuadraticEq2}: $\overline{\Delta}=(\overline{\mathcal{X}}_1\cos{\overline{\Theta}_1}-r\cos\phistar\sin{\overline{\Theta}_1})^2-(\overline{\mathcal{X}}^2_1-r^2\sin^2{\phistar})\ge r^2\sin^2{\phistar}-\overline{\mathcal{X}}^2_1>0$ and $\overline{\constant}>0$.\newline If $\overline{\mathcal{X}}_1=r\sin\phistar$ and $\phistar\ne\overline{\Theta}_1$, then \eqref{eq:tau1barlengthQuadraticEq2}: $\overline{\Delta}=(\overline{\mathcal{X}}_1\cos{\overline{\Theta}_1}-r\cos\phistar\sin{\overline{\Theta}_1})^2=r^2\sin^2{(\phistar-\overline{\Theta}_1)}>0$. (with $\overline{\Theta}_1\in[0,\pi]$), and $\overline{\constant}>0$.\newline If $\overline{\mathcal{X}}_1=-r\sin\phistar$ and $\phistar\ne\pi-\overline{\Theta}_1$, then \eqref{eq:tau1barlengthQuadraticEq2}: $\overline{\Delta}=(\overline{\mathcal{X}}_1\cos{\overline{\Theta}_1}-r\cos\phistar\sin{\overline{\Theta}_1})^2=r^2\sin^2{(\phistar+\overline{\Theta}_1)}>0$ (with $\overline{\Theta}_1\in[0,\pi]$), and $\overline{\constant}>0$\newline If $\overline{\mathcal{X}}_1=r\sin\phistar$ and $\phistar=\overline{\Theta}_1$, then \eqref{eq:tau1barlengthQuadraticEq2}: $\overline{\Delta}=(\overline{\mathcal{X}}_1\cos{\overline{\Theta}_1}-r\cos\phistar\sin{\overline{\Theta}_1})^2=r^2\sin^2{(\phistar-\overline{\Theta}_1)}=0$, and $\overline{\constant}=0$\newline If $\overline{\mathcal{X}}_1=-r\sin\phistar$ and $\phistar=\pi-\overline{\Theta}_1$, then \eqref{eq:tau1barlengthQuadraticEq2}: $\overline{\Delta}=(\overline{\mathcal{X}}_1\cos{\overline{\Theta}_1}-r\cos\phistar\sin{\overline{\Theta}_1})^2=r^2\sin^2{(\phistar+\overline{\Theta}_1)}=0$, and $\overline{\constant}=0$

Therefore, it is not hard to see that $\overline{\constant}_+\ge0$, $\overline{\constant}_-\ge0$, $\overline{\constant}_+\cdot\overline{\constant}_-\ge\overline{\constant}_0\ge0$ and if and only if $\overline{\constant}=0$, then $\overline{\Delta}=0$.

Hence, if $\overline{\constant}\ne 0$, then $\overline{\Delta}\ge\max\Big\{\min\big\{\overline{\constant}^2_+,\overline{\constant}^2_-\big\}\cdot\chi_{\{|\overline{\mathcal{X}}_1|= r\sin\phistar\}},\overline{\constant}_0\Big\}=\overline{\constant}>0$ is a positive constant, and we can compare/estimate the difference of $\sqrt{\Delta}$ and $\sqrt{\overline{\Delta}}$ in the following.
\begin{equation}\label{eq:DeltaCompareIneq6_1}
    \big|\sqrt{\Delta}-\sqrt{\overline{\Delta}}\bigm|=\frac{\bigm|\Delta-\overline{\Delta}\bigm|}{\sqrt{\Delta}+\sqrt{\overline{\Delta}}}<\frac{1}{\underbracket{\sqrt{\overline{\Delta}}}_{\mathrlap{\ge\sqrt{\overline{\constant}}}}}\bigm|\Delta-\overline{\Delta}\bigm|\le\frac{\bigm|\Delta-\overline{\Delta}\bigm|}{\sqrt{\overline{\constant}}}\underbracket{<}_{\eqref{eq:deltadiff}}\frac{2r}{\sqrt{\overline{\constant}}}\big(3\bigm|\overline{\mathcal{X}}_1-\mathcal{X}_1\bigm|+3\bigm|\mathcal{Y}_1+r\cos{\phistar}\bigm|+4r\bigm|\Theta_1-\overline{\Theta}_1\bigm|\big).
\end{equation}

On the other hand, if $\overline{\constant}=0$, then we have $\overline{\Delta}=0$ and the following.

\begin{equation}\label{eq:DeltaCompareIneq6_2}
    \begin{multlined}
        \big|\sqrt{\Delta}-\sqrt{\overline{\Delta}}\bigm|=\bigm|\sqrt{\Delta}\bigm|\overbracket{=}^{\eqref{eq:tplus}}\big[\underbracket{\bigm|\mathcal{X}_1\cos{\Theta_1}+\mathcal{Y}_1\sin{\Theta_1}\bigm|^2}_{=\alpha\ge0}\underbracket{-2b(\mathcal{Y}_1+r\cos{\phistar}}_{=\beta>0})\big]^{1/2}\\
        \overbracket{\le}^{\mathllap{\sqrt{\alpha+\beta}\le0.5(\sqrt{\alpha}+\sqrt{\beta})}}\frac{1}{2}(\bigm|\mathcal{X}_1\cos{\Theta_1}+\mathcal{Y}_1\sin{\Theta_1}\bigm|+\bigm|2b(\mathcal{Y}_1+r\cos{\phistar})\bigm|^{1/2}).
    \end{multlined}
\end{equation}

Notice that by $\overline\constant=0$, we have either $\overline{\mathcal{X}}_1=r\sin\phistar$ with $\phistar=\overline{\Theta}_1$ or $\overline{\mathcal{X}}_1=-r\sin\phistar$ with $\phistar=\pi-\overline{\Theta}_1$, then we use \eqref{eq:DeltaCompareIneq4} to get the following.
\begin{align*}
\big|\mathcal{X}_1\cos{\Theta_1}+\mathcal{Y}_1\sin{\Theta_1}\bigm|=&\bigm|\overbracket{\overline{\mathcal{X}}_1\cos{\overline{\Theta}_1}-r\cos{\phistar}\sin{\overline{\Theta}_1}}^{\mathclap{=r\sin{(\phistar-\overline{\Theta}_1)}=0\text{ with }\overline{\mathcal{X}}_1=r\sin\phistar,\,\phistar=\overline{\Theta}_1\text{ or }=r\sin{(\phistar+\overline{\Theta}_1)}=0\text{ with }\overline{\mathcal{X}}_1=-r\sin\phistar,\,\phistar=\pi-\overline{\Theta}_1}}-(\mathcal{X}_1\cos{\Theta_1}+\mathcal{Y}_1\sin{\Theta_1})\big|\\
\overbracket{<}^{\mathllap{\eqref{eq:DeltaCompareIneq4}}}&\bigm|\overline{\mathcal{X}}_1-\mathcal{X}_1\bigm|+\bigm|r\cos{\phistar}+\mathcal{Y}_1\bigm|+2r\bigm|\Theta_1-\overline{\Theta}_1\big|.
\end{align*}

And we also use \eqref{eq:DeltaCompareIneq1}\eqref{eq:DeltaCompareIneq2} to get the following.

\begin{multline*}
\big|2b(\mathcal{Y}_1+r\cos{\phistar})\bigm|= |2b(\mathcal{Y}_1+r\cos{\phistar})-(\mathcal{X}^2_1-r^2\sin^2{\phistar})|+|\mathcal{X}^2_1-r^2\sin^2{\phistar}|\overbracket{<}^{\mathllap{\eqref{eq:DeltaCompareIneq1}}}2r|\mathcal{Y}_1+r\cos{\phistar}|+|\mathcal{X}^2_1-r^2\sin^2{\phistar}|\\
\overbracket{<}^{\mathllap{\text{\eqref{eq:DeltaCompareIneq2} with }\overline{\mathcal{X}}^2_1=r^2\sin^2\phistar}}4r|\mathcal{Y}_1+r\cos{\phistar}|+2r\bigm|\overline{\mathcal{X}}_1-\mathcal{X}_1\big|.\\
\end{multline*}

Then \eqref{eq:DeltaCompareIneq6_2} together with \eqref{eq:DeltaCompareIneq6_1} we get
\begin{equation}\label{eq:DeltaCompareIneq6_3}
    \big|\sqrt{\Delta}-\sqrt{\overline{\Delta}}\bigm|<\left\{\begin{aligned}
    &\frac{2r}{\sqrt{\overline{\constant}}}\big(3\bigm|\overline{\mathcal{X}}_1-\mathcal{X}_1\bigm|+3\bigm|\mathcal{Y}_1+r\cos{\phistar}\bigm|+4r\bigm|\Theta_1-\overline{\Theta}_1\bigm|\big) \text{, if }\overline{\constant}>0.\\
    &0.5\bigm|\overline{\mathcal{X}}_1-\mathcal{X}_1\bigm|+0.5\bigm|r\cos{\phistar}+\mathcal{Y}_1\bigm|+r\bigm|\Theta_1-\overline{\Theta}_1\bigm|+(r|\mathcal{Y}_1+r\cos{\phistar}|+0.5r\bigm|\overline{\mathcal{X}}_1-\mathcal{X}_1\bigm|)^{1/2}\text{, if }\overline{\constant}=0.
    \end{aligned}
    \right.
\end{equation}

And in the following inequality \eqref{eq:DeltaCompareIneq7}, we estimate the bound for $|t-\overline{t}|$ by $\big|\overline{\mathcal{X}}_1-\mathcal{X}_1\big|$, $\big|\overline{\mathcal{Y}}_1-\mathcal{Y}_1\big|$ and $\big|\overline{\Theta}_1-\Theta_1\big|$.

\begin{equation}\label{eq:DeltaCompareIneq7}
    \begin{aligned}
        \big|t-\overline{t}\bigm|\overbracket{=}^{\mathllap{\text{\eqref{eq:tplus}, \eqref{eq:tau1barlengthQuadraticEq2}}}}&\bigm|-\mathcal{X}_1\cos{\Theta_1}-\mathcal{Y}_1\sin{\Theta_1}+\sqrt{\Delta}-(-\overline{\mathcal{X}}_1\cos{\overline{\Theta}_1}+r\cos\phistar\sin{\overline{\Theta}_1}+\sqrt{\overline{\Delta}})\big|\\
        \le&\underbracket{\bigm|-\mathcal{X}_1\cos{\Theta_1}-\mathcal{Y}_1\sin{\Theta_1}-(-\overline{\mathcal{X}}_1\cos{\overline{\Theta}_1}+r\cos\phistar\sin{\overline{\Theta}_1})\bigm|}_{\qquad\qquad\qquad\mathllap{\text{\eqref{eq:DeltaCompareIneq4}: }<|\overline{\mathcal{X}}_1-\mathcal{X}_1|+|r\cos{\phistar}+\mathcal{Y}_1|+2r|\Theta_1-\overline{\Theta}_1|}}+\underbracket{\bigm|\sqrt{\Delta}-\sqrt{\overline{\Delta}}\bigm|}_{\mathclap{\text{\eqref{eq:DeltaCompareIneq6_3} }}}\\
        \overbracket{<}^{\mathllap{\text{\eqref{eq:DeltaCompareIneq4},\eqref{eq:DeltaCompareIneq6_3}}}}&\left\{\begin{aligned}
            &\frac{6r+3\sqrt{\overline{\constant}}}{\sqrt{\overline{\constant}}}\bigm|\overline{\mathcal{X}}_1-\mathcal{X}_1\bigm|+\frac{6r+3\sqrt{\overline{\constant}}}{\sqrt{\overline{\constant}}}\bigm|\mathcal{Y}_1+r\cos{\phistar}\bigm|+\frac{8r^2+2r\sqrt{\overline{\constant}}}{\sqrt{\overline{\constant}}}\bigm|\Theta_1-\overline{\Theta}_1\bigm| \text{, if }\overline{\constant}>0.\\
            &1.5\bigm|\overline{\mathcal{X}}_1-\mathcal{X}_1\bigm|+1.5\bigm|r\cos{\phistar}+\mathcal{Y}_1\bigm|+3r\bigm|\Theta_1-\overline{\Theta}_1\bigm|+(r|\mathcal{Y}_1+r\cos{\phistar}|+0.5r\bigm|\overline{\mathcal{X}}_1-\mathcal{X}_1\bigm|)^{1/2}\text{, if }\overline{\constant}=0.
        \end{aligned}\right.
    \end{aligned}
\end{equation}

Therefore, we see that for $\forall \varepsilon>0$, if $\big|\mathcal{X}_1-\overline{\mathcal{X}}_1\big|$, $\big|\mathcal{Y}_1-\overline{\mathcal{Y}}_1\big|$, $\big|\Theta_1-\overline{\Theta}_1\big|$ satisfy
\begin{equation}\label{eq:ttbarvarepsilonCondition}
    \left\{\begin{aligned}
        &\big|\overline{\mathcal{X}}_1-\mathcal{X}_1\big|<\frac{0.3\sqrt{\overline{\constant}}}{6r+3\sqrt{\overline{\constant}}}\varepsilon,\quad \big|\mathcal{Y}_1-\overline{\mathcal{Y}}_1\big|<\frac{0.3\sqrt{\overline{\constant}}}{6r+3\sqrt{\overline{\constant}}}\varepsilon,\quad \bigm|\Theta_1-\overline{\Theta}_1\bigm|<\frac{0.4\sqrt{\overline{\constant}}}{8r^2+2r\sqrt{\overline{\constant}}}\varepsilon \text{, if }\overline{\constant}>0,\\
        &\big|\overline{\mathcal{X}}_1-\mathcal{X}_1\bigm|<\min\big\{\frac{0.09\varepsilon^2}{r},0.1\varepsilon\big\},\quad\bigm|\overline{\mathcal{Y}}_1-\mathcal{Y}_1\bigm|<\min\big\{\frac{0.045\varepsilon^2}{r},0.1\varepsilon\big\},\quad \bigm|\Theta_1-\overline{\Theta}_1\bigm|<0.1\frac{\varepsilon}{r}\text{, if }\overline{\constant}=0,
    \end{aligned}\right.
\end{equation} then $\big|t-\overline{t}\big|<\varepsilon$. This ends the first part of the proof of \cref{lemma:Trajectory_approximation_2}.


We will next see that $\big|\phi_2-\overline{\phi}_2\big|$ and $\big|\theta_2-\overline{\theta}_2\big|$ are controlled by $\big|\overline{\mathcal{X}}_1-\mathcal{X}_1\big|$, $\big|\overline{\mathcal{Y}}_1-\mathcal{Y}_1\big|$, $\big|\overline{\Theta}_1-\Theta_1\big|$ and $\big|t-\overline{t}\big|$ (see \eqref{eq:phi2phibar2diff} and \eqref{eq:theta2thetabar2diff} on pages \pageref{eq:phi2phibar2diff} and \pageref{eq:theta2thetabar2diff}, respectively). By \eqref{eq:DeltaCompareIneq7} the last of these four quantities is controlled by the other three, so $\big|\phi_2-\overline{\phi}_2\big|$, $\big|\theta_2-\overline{\theta}_2\big|$ is controlled by $\big|\overline{\mathcal{X}}_1-\mathcal{X}_1\big|$, $\big|\overline{\mathcal{Y}}_1-\mathcal{Y}_1\big|$ and $\big|\overline{\Theta}_1-\Theta_1\big|$ alone.

In order to prove \eqref{eq:phi2phibar2diff}, we will use inequality \eqref{eq:sindiffphi2phibar2} and \eqref{eq:phiphi2diffrange} with conditions \eqref{eq:X1Y1Theta1closeCond1} and \eqref{eq:X1Y1Theta1closeCond2} (it is clear that \eqref{eq:X1Y1Theta1closeCond2} implies \eqref{eq:X1Y1Theta1closeCond1}). 
Then inequality \eqref{eq:theta2thetabar2diff} follows from \COMMENT{\small\bfseries\red``follows'' or ``follows from''? Search through entire text for ``follows'' and make sure that in each instance it says exactly what it should, either ``follows'' or ``follows from.'' These mean different things.} inequality \eqref{eq:phi2phibar2diff}.

In order to prove \eqref{eq:phiphi2diffrange}, we also use inequality \eqref{eq:sindiffphi2phibar2}, conditions \eqref{eq:X1Y1Theta1closeCond2} and condition \eqref{eq:phi2phibar2diffrange}, where condition \eqref{eq:phi2phibar2diffrange} is guaranteed by inequalities \eqref{eq:cosdiffphi2phibar2_1}, \eqref{eq:cosdiffphi2phibar2_2}, \eqref{eq:cosdiffphi2phibar2_3} and conditions \eqref{eq:X1Y1Theta1closeCond1}, thus by conditions \eqref{eq:X1Y1Theta1closeCond2}.

Also note that in order to ensure that $\big|t-\overline{t}\big|<0.05r\sin{(2\phistar)}$\COMMENT{\large\red\bfseries This looks weird and distracting because $\bigm|$ is wrong here. Get rid of it here and everywhere else where it does not belong. Look up online what its purpose is. It may be suitable to use $\big|$ in some places for clarity, but $\bigm|$ creates spacing issues unsuitable for absolute values.\Hrule Here, there is no reason to make anything larger than normal.\wentao{fixed here. checking.}} in conditions \eqref{eq:X1Y1Theta1closeCond2} is satisfied, in addition to the conditions of  $\big|\overline{\mathcal{X}}_1-\mathcal{X}_1\big|$, $\big|\overline{\mathcal{Y}}_1-\mathcal{Y}_1\big|$, $\big|\overline{\Theta}_1-\Theta_1\big|$ in \eqref{eq:X1Y1Theta1closeCond2}, we impose additional conditions \eqref{eq:X1Y1Theta1closeCond3}.

Finally, given \eqref{eq:phi2phibar2diff} and \eqref{eq:theta2thetabar2diff}, for $\forall\varepsilon>0$ we impose the conditions \eqref{eq:ttbarvarepsilonCondition2} with additional conditions \eqref{eq:X1X1barY1Y1barTheta1Thetabar1Varepsilon} for $\big|\overline{\mathcal{X}}_1-\mathcal{X}_1\big|$, $\big|\overline{\mathcal{Y}}_1-\mathcal{Y}_1\big|$, $\big|\overline{\Theta}_1-\Theta_1\big|$ so that we get the desired $\varepsilon$-bound for $\big|\phi_2-\overline{\phi}_2\big|$ and $\big|\theta_2-\overline{\theta}_2\big|$ in \eqref{eq:maxphi2phibar2andtheta2theta2bardiff}, which is the conclusion \eqref{eq:phi2theta2close}.

Also, for $\forall\varepsilon>0$ we notice that the conditions for $|\overline{\mathcal{X}}_1-\mathcal{X}_1\bigm|$, $|\overline{\mathcal{Y}}_1-\mathcal{Y}_1\bigm|$, $|\overline{\Theta}_1-\Theta_1\bigm|$ in \eqref{eq:X1Y1Theta1CloseCondition} imply the various conditions in \eqref{eq:ttbarvarepsilonCondition2}, \eqref{eq:X1X1barY1Y1barTheta1Thetabar1Varepsilon}, \eqref{eq:X1Y1Theta1closeCond3}, \eqref{eq:X1Y1Theta1closeCond2} (thus \eqref{eq:X1Y1Theta1closeCond1} and \eqref{eq:phi2phibar2diffrange}). 

Therefore, we finish the proof for the claim that the conditions in \eqref{eq:X1Y1Theta1CloseCondition} imply the conclusion in \eqref{eq:phi2theta2close}.

We now first\COMMENT{``next first'' does not work} prove several trigonometric inequalities:  \eqref{eq:sindiffphi2phibar2}, \eqref{eq:cosdiffphi2phibar2_1}, \eqref{eq:cosdiffphi2phibar2_2} and \eqref{eq:cosdiffphi2phibar2_3}.

\begin{equation}\label{eq:sindiffphi2phibar2}
    \begin{aligned}
        r^2\bigm|\sin{(\phi_2-\overline{\phi}_2)}\bigm|&=r^2\bigm|\sin{(\overline{\phi}_2-\phi_2)}\bigm|=\bigm|-\underbracket{(r\sin{\phi_2})}_{\qquad\mathllap{\text{\eqref{eq:tau1lengthQuadraticEq}: }=\mathcal{X}_1+t\cos{\Theta_1}}}\overbracket{(r\cos{\overline{\phi}_2})}^{\mathllap{\text{\eqref{eq:tau1barlengthQuadraticEq}: }=r\cos{\phistar}-\overline{t}\sin{\overline{\Theta}_1}}}+\underbracket{(r\cos{\phi_2})}_{\qquad\mathllap{\text{\eqref{eq:tau1lengthQuadraticEq}: }=-\mathcal{Y}_1-t\sin{\Theta_1}}}\overbracket{(r\sin{\overline{\phi}_2})}^{\mathrlap{\text{\eqref{eq:tau1barlengthQuadraticEq}: }\overline{\mathcal{X}}_1+\overline{t}\cos{\overline{\Theta}_1}}\qquad}\bigm|\\
        &=\bigm|-
        (\mathcal{X}_1+t\cos{\Theta_1})(r\cos{\phistar}-\overline{t}\sin{\overline{\Theta}_1})-(\mathcal{Y}_1+t\sin{\Theta_1})(\overline{\mathcal{X}}_1+\overline{t}\cos{\overline{\Theta}_1})\bigm|\\
        &=\big|-\mathcal{X}_1r\cos{\phistar}-tr\cos{\phistar}\cos{\Theta_1}+\mathcal{X}_1\overline{t}\sin{\overline{\Theta}_1}+\cos{\Theta_1}\sin{\overline{\Theta}_1}t\overline{t}\\
        &-\overline{\mathcal{X}}_1\mathcal{Y}_1-\overline{\mathcal{X}}_1\sin{\Theta_1}t-\mathcal{Y}_1\cos{\overline{\Theta}_1}\overline{t}-\sin{\Theta_1}\cos{\overline{\Theta}_1}t\overline{t}\bigm|\\
        \le&\bigm|\underbracket{-\mathcal{X}_1r\cos{\phistar}-\overline{\mathcal{X}}_1\mathcal{Y}_1}_{\qquad\qquad\mathllap{=-\mathcal{X}_1r\cos{\phistar}-\mathcal{X}_1\mathcal{Y}_1+\mathcal{X}_1\mathcal{Y}_1-\overline{\mathcal{X}}_1\mathcal{Y}_1}}\bigm|+\bigm|\underbracket{\mathcal{X}_1\overline{t}\sin{\overline{\Theta}_1}-\overline{\mathcal{X}}_1\sin{\Theta_1}t}_{\mathrlap{=\mathcal{X}_1\overline{t}\sin{\overline{\Theta}_1}-\overline{\mathcal{X}}_1\overline{t}\sin{\overline{\Theta}_1}+\overline{\mathcal{X}}_1\overline{t}\sin{\overline{\Theta}_1}-\overline{\mathcal{X}}_1\sin{\Theta_1}t}}\bigm|\\
        &+\bigm|-tr\cos{\phistar}\cos{\Theta_1}-\mathcal{Y}_1\cos{\overline{\Theta}_1}\overline{t}\bigm|+\bigm|\cos{\Theta_1}\sin{\overline{\Theta}_1}-\sin{\Theta_1}\cos{\overline{\Theta}_1}\bigm|t\overline{t}\\
        \overbracket{\le}^{\mathllap{\text{apply triangle inequality}}}&\underbracket{\bigm|-\mathcal{X}_1r\cos{\phistar}-\mathcal{X}_1\mathcal{Y}_1\bigm|}_{< r|\mathcal{Y}_1+r\cos\phistar|}+\underbracket{\bigm|\mathcal{X}_1\mathcal{Y}_1-\overline{\mathcal{X}}_1\mathcal{Y}_1\bigm|}_{<r|\mathcal{X}_1-\overline{\mathcal{X}}_1|}+\underbracket{\bigm|\mathcal{X}_1\overline{t}\sin{\overline{\Theta}_1}-\overline{\mathcal{X}}_1\overline{t}\sin{\overline{\Theta}_1}\bigm|}_{\text{By }\overline{t}\le2r\text{, }<2r|\mathcal{X}_1-\overline{\mathcal{X}}_1|}\\
        &+\underbracket{\bigm|\overline{\mathcal{X}}_1\overline{t}\sin{\overline{\Theta}_1}-\overline{\mathcal{X}}_1\sin{\Theta_1}t\bigm|}_{\mathllap{\substack{=|\overline{\mathcal{X}}_1\overline{t}\sin{\overline{\Theta}_1}-\overline{\mathcal{X}}_1t\sin{\overline{\Theta}_1}+\overline{\mathcal{X}}_1t\sin{\overline{\Theta}_1}-\overline{\mathcal{X}}_1\sin{\Theta_1}t|\\ <r|t-\overline{t}|+2r^2|\Theta_1-\overline{\Theta}_1|}}}+\underbracket{\bigm|-tr\cos{\phistar}\cos{\Theta_1}+\overline{t}r\cos{\phistar}\cos{\Theta_1}\bigm|}_{\mathrlap{r|t-\overline{t}|}}\\
        &+\underbracket{\bigm|-tr\cos{\phistar}\cos{\Theta_1}-\mathcal{Y}_1\cos{\overline{\Theta}_1}\overline{t}\bigm|}_{\qquad\qquad\mathllap{\substack{=|-tr\cos{\phistar}\cos{\Theta_1}+tr\cos{\phistar}\cos{\overline{\Theta}_1}-tr\cos{\phistar}\cos{\overline{\Theta}_1}\\+\overline{t}r\cos{\phistar}\cos{\overline{\Theta}_1}-\overline{t}r\cos{\phistar}\cos{\overline{\Theta}_1}-\mathcal{Y}_1\cos{\overline{\Theta}_1}\overline{t}|}}}+\underbracket{t\overline{t}\bigm|\sin{(\Theta_1-\overline{\Theta}_1)}\bigm|}_{\mathrlap{\substack{\text{By Mean-Value Theorem and }t\le2r,\overline{t}\le2r: \\\le4r^2|\Theta_1-\overline{\Theta}_1|}}}\\
        <&3r\bigm|\mathcal{X}_1-\overline{\mathcal{X}}_1\bigm|+r\bigm|\mathcal{Y}_1+r\cos\phistar\bigm|+6r^2\bigm|\Theta_1-\overline{\Theta}_1\bigm|+2r\bigm|t-\overline{t}\big|\\
        &+\underbracket{\bigm|-tr\cos{\phistar}\cos{\Theta_1}+tr\cos{\phistar}\cos{\overline{\Theta}_1}\bigm|}_{\text{By} |t|\le2r \text{ and Mean-Value Theorem }<2r^2|\Theta_1-\overline{\Theta}_1|}+\underbracket{\bigm|-tr\cos{\phistar}\cos{\overline{\Theta}_1}+\overline{t}r\cos{\phistar}\cos{\overline{\Theta}_1}\bigm|}_{<r|t-\overline{t}|}\\
        &+\underbracket{\bigm|-\overline{t}r\cos{\phistar}\cos{\overline{\Theta}_1}-\mathcal{Y}_1\cos{\overline{\Theta}_1}\overline{t}\bigm|}_{<2r|\mathcal{Y}_1+r\cos\phistar|}\\
        <&3r\bigm|\mathcal{X}_1-\overline{\mathcal{X}}_1\bigm|+3r\bigm|\mathcal{Y}_1+r\cos\phistar\bigm|+8r^2\bigm|\Theta_1-\overline{\Theta}_1\bigm|+3r\bigm|t-\overline{t}\big|\\
    \end{aligned}
\end{equation}

\begin{equation}\label{eq:cosdiffphi2phibar2_1}
    \begin{aligned}
        r^2\cos{(\overline{\phi}_2-\phi_2)}&=\overbracket{(\mathcal{Y}_1+t\sin{\Theta_1})}^{\mathllap{\text{\eqref{eq:tau1lengthQuadraticEq}: }=-r\cos{\phi_2}}}\overbracket{(\overline{\mathcal{Y}}_1+\overline{t}\sin{\overline{\Theta}_1})}^{\text{\eqref{eq:tau1barlengthQuadraticEq}: }=-r\cos{\overline{\phi}_2}}+\overbracket{(\mathcal{X}_1+t\cos{\Theta_1})}^{\text{\eqref{eq:tau1lengthQuadraticEq}: }=r\sin{\phi_2}}\overbracket{(\overline{\mathcal{X}}_1+\overline{t}\cos{\overline{\Theta}_1})}^{\text{\eqref{eq:tau1barlengthQuadraticEq}: }=r\sin{\phi_2}}\\
        &=\overbracket{(\overline{\mathcal{Y}}_1+\overline{t}\sin{\overline{\Theta}_1})^2}^{=r^2\cos^2{\overline{\phi}_2}}+(\mathcal{Y}_1+t\sin{\Theta_1}-\overline{\mathcal{Y}}_1-\overline{t}\sin{\overline{\Theta}_1})(\overline{\mathcal{Y}}_1+\overline{t}\sin{\overline{\Theta}_1})\\
        &+\overbracket{(\overline{\mathcal{X}}_1+\overline{t}\cos{\overline{\Theta}_1})^2}^{=r^2\sin^2{\overline{\phi}_2}}+(\mathcal{X}_1+t\cos{\Theta_1}-\overline{\mathcal{X}}_1-\overline{t}\cos{\overline{\Theta}_1})(\overline{\mathcal{X}}_1+\overline{t}\cos{\overline{\Theta}_1})\\
        &=r^2-(-\mathcal{Y}_1-t\sin{\Theta_1}+\overline{\mathcal{Y}}_1+\overline{t}\sin{\overline{\Theta}_1})(\overline{\mathcal{Y}}_1+\overline{t}\sin{\overline{\Theta}_1})\\
        &-(-\mathcal{X}_1-t\cos{\Theta_1}+\overline{\mathcal{X}}_1+\overline{t}\cos{\overline{\Theta}_1})(\overline{\mathcal{X}}_1+\overline{t}\cos{\overline{\Theta}_1}),\\
    \end{aligned}
\end{equation}
\begin{equation}\label{eq:cosdiffphi2phibar2_2}
    \begin{aligned}
        (-\mathcal{Y}_1-t\sin{\Theta_1}+\overline{\mathcal{Y}}_1+\overline{t}\sin{\overline{\Theta}_1})(\overline{\mathcal{Y}}_1+\overline{t}\sin{\overline{\Theta}_1})\le&\bigm|-\mathcal{Y}_1-t\sin{\Theta_1}+\overline{\mathcal{Y}}_1+\overline{t}\sin{\overline{\Theta}_1}\bigm|\overbracket{\bigm|\overline{\mathcal{Y}}_1+\overline{t}\sin{\overline{\Theta}_1}\big|}^{\text{\eqref{eq:tau1barlengthQuadraticEq}: }=r|\cos{\overline{\phi}_2}|\le r}\\
        \le&r\bigm|\overline{\mathcal{Y}}_1-\mathcal{Y}_1-t\sin{\Theta_1}+\overline{t}\sin{\Theta_1}-\overline{t}\sin{\Theta_1}+\overline{t}\sin{\overline{\Theta}_1}\big|\\
        \overbracket{\le}^{\mathllap{\text{By Mean-Value Theorem and }t,\overline{t}\le2r}}&r\bigm|\overline{\mathcal{Y}}_1-\mathcal{Y}_1\bigm|+r\bigm|t-\overline{t}\bigm|+2r^2\bigm|\Theta_1-\overline{\Theta}_1\bigm|,
    \end{aligned}
\end{equation}
\begin{equation}\label{eq:cosdiffphi2phibar2_3}
    \begin{aligned}
        (-\mathcal{X}_1-t\cos{\Theta_1}+\overline{\mathcal{X}}_1+\overline{t}\cos{\overline{\Theta}_1})(\overline{\mathcal{X}}_1+\overline{t}\cos{\overline{\Theta}_1})\le&\bigm|-\mathcal{X}_1-t\cos{\Theta_1}+\overline{\mathcal{X}}_1+\overline{t}\cos{\overline{\Theta}_1}\bigm|\overbracket{\bigm|\overline{\mathcal{X}}_1+\overline{t}\cos{\overline{\Theta}_1}\bigm|}^{=r|\sin{\overline{\phi}_1}|\le r}\\
        \le&r\bigm|\overline{\mathcal{X}}_1-\mathcal{X}_1-t\cos{\Theta_1}+\overline{t}\cos{\Theta_1}-\overline{t}\cos{\Theta_1}+\overline{t}\cos{\overline{\Theta}_1}\bigm|\\
        \overbracket{\le}^{\mathllap{\text{By Mean-Value Theorem and }t,\overline{t}\le2r}}&r\bigm|\overline{\mathcal{X}}_1-\mathcal{X}_1\bigm|+r\bigm|\overline{t}-t\bigm|+2r^2\bigm|\Theta_1-\overline{\Theta}_1\bigm|.
    \end{aligned}
\end{equation}
By \eqref{eq:cosdiffphi2phibar2_1}, \eqref{eq:cosdiffphi2phibar2_2}, \eqref{eq:cosdiffphi2phibar2_3}, we see that if \begin{equation}\label{eq:X1Y1Theta1closeCond1}
\begin{aligned}
    \bigm|\overline{\mathcal{X}}_1-\mathcal{X}_1\bigm|&<0.1r,\\
    \bigm|\overline{\mathcal{Y}}_1-\mathcal{Y}_1\bigm|&<0.1r,\\
    \bigm|t-\overline{t}\bigm|&<0.1r,\\
    \bigm|\Theta_1-\overline{\Theta}_1\bigm|&<0.1,
\end{aligned}
\end{equation} then \begin{align*}
    r^2-r^2\cos{(\phi_2-\overline{\phi}_2)}\overbracket{=}^{\mathllap{\eqref{eq:cosdiffphi2phibar2_1}}}&(-\mathcal{Y}_1-t\sin{\Theta_1}+\overline{\mathcal{Y}}_1+\overline{t}\sin{\overline{\Theta}_1})(\overline{\mathcal{Y}}_1+\overline{t}\sin{\overline{\Theta}_1})\\
    &+(-\mathcal{X}_1-t\cos{\Theta_1}+\overline{\mathcal{X}}_1+\overline{t}\cos{\overline{\Theta}_1})(\overline{\mathcal{X}}_1+\overline{t}\cos{\overline{\Theta}_1})\\
    \overbracket{<}^{\mathllap{\eqref{eq:cosdiffphi2phibar2_2}, \eqref{eq:cosdiffphi2phibar2_3}}}&r\bigm|\overline{\mathcal{X}}_1-\mathcal{X}_1\bigm|+r\bigm|\overline{\mathcal{Y}}_1-\mathcal{Y}_1\bigm|+2r\bigm|\overline{t}-t\bigm|+4r^2|\Theta_1-\overline{\Theta}_1|\\
    \overbracket{<}^{\mathllap{\text{\eqref{eq:X1Y1Theta1closeCond1}}}}&r^2,
\end{align*} which implies $\cos{(\phi_2-\overline{\phi}_2)}>0$. Since $\phi_2,\overline{\phi}_2\in\big[\phistar,2\pi-\phistar\big]$, so 
\begin{equation}\label{eq:phi2phibar2diffrange}
\phi_2-\overline{\phi}_2\in[-2\pi+2\phistar,\,-\frac{3\pi}{2})\cup(-\frac{\pi}{2},\,\frac{\pi}{2})\cup(\frac{3\pi}{2},\,2\pi-2\phistar].\end{equation}\newline  With $\phistar<\pi/4$ and by \eqref{eq:sindiffphi2phibar2} if \begin{equation}\label{eq:X1Y1Theta1closeCond2}
    \begin{aligned}
        \bigm|\mathcal{X}_1-\overline{\mathcal{X}}_1\bigm|&<0.05r\sin{(2\phistar)},\\
        \bigm|\mathcal{Y}_1-\overbracket{\overline{\mathcal{Y}}_1}^{\mathrlap{=-r\cos\phistar}}\bigm|&<0.05r\sin{(2\phistar)},\\
        \bigm|\Theta_1-\overline{\Theta}_1\bigm|&<0.05\sin{(2\phistar)},\\
        \bigm|t-\overline{t}\bigm|&<0.05r\sin{(2\phistar)},
    \end{aligned}
\end{equation} then 
\begin{equation}\label{eq:phiphi2diffrange}
\begin{aligned}
    &r^2\bigm|\sin{(\phi_2-\overline{\phi}_2)}\bigm|\overbracket{<}^{\mathllap{\eqref{eq:sindiffphi2phibar2}}}3r\bigm|\mathcal{X}_1-\overline{\mathcal{X}}_1\bigm|+3r\bigm|\mathcal{Y}_1+r\cos\phistar\bigm|+8r^2\bigm|\Theta_1-\overline{\Theta}_1\bigm|+3r\bigm|t-\overline{t}\bigm|<r^2\sin(2\phistar)\\
    \Rightarrow& \bigm|\sin{(\phi_2-\overline{\phi}_2)}\bigm|<\sin{(2\phistar)}\overbracket{\Rightarrow}^{\mathclap{2\phistar<\pi/2}} \phi_2-\overline{\phi}_2\notin[-2\pi+2\phistar,\,-\frac{3\pi}{2})\cup(\frac{3\pi}{2},\,2\pi-2\phistar]\overbracket{\Rightarrow}^{\eqref{eq:phi2phibar2diffrange}} (\phi_2-\overline{\phi}_2)\in(-\pi/2,\pi/2).
\end{aligned}\end{equation} Therefore, since if $|z|<\pi/2$, then $|\sinc(z)|=\frac{|\sin{z}|}{|z|}>\frac{2}{\pi}$. we get the following estimate for $|\phi_2-\overline{\phi}_2|$. \begin{equation}\label{eq:phi2phibar2diff}
    \begin{aligned}
        \big|\phi_2-\overline{\phi}_2\bigm|&<\frac{\pi}{2}\bigm|\sin{(\phi_2-\overline{\phi}_2)}\bigm|\overbracket{<}^{\eqref{eq:sindiffphi2phibar2}}\frac{\pi}{2}\big(\frac{3}{r}\bigm|\mathcal{X}_1-\overline{\mathcal{X}}_1\bigm|+\frac{3}{r}\bigm|\mathcal{Y}_1-\overline{\mathcal{Y}}_1\bigm|+8\bigm|\Theta_1-\overline{\Theta}_1\bigm|+\frac{3}{r}\bigm|t-\overline{t}\bigm|\big).\\
    \end{aligned} 
\end{equation} 
To satisfy $\big|t-\overline{t}\big|<0.05r\sin{(2\phistar)}$ in \eqref{eq:X1Y1Theta1closeCond2},  
by \eqref{eq:DeltaCompareIneq7} and \eqref{eq:ttbarvarepsilonCondition} it suffices to make $\big|\mathcal{X}_1-\overline{\mathcal{X}}_1\big|$, $\big|\mathcal{Y}_1-\overline{\mathcal{Y}}_1\big|$ and $\big|\Theta_1-\overline{\Theta}_1\big|$ satisfy
\begin{equation}\label{eq:X1Y1Theta1closeCond3}
    \left\{\begin{aligned}
        &\big|\overline{\mathcal{X}}_1-\mathcal{X}_1\big|<\frac{0.015\sqrt{\overline{\constant}}r\sin{(2\phistar)}}{6r+3\sqrt{\overline{\constant}}},\quad \big|\mathcal{Y}_1-\overline{\mathcal{Y}}_1\big|<\frac{0.015\sqrt{\overline{\constant}}r\sin{(2\phistar)}}{6r+3\sqrt{\overline{\constant}}},\quad \big|\Theta_1-\overline{\Theta}_1\big|<\frac{0.02\sqrt{\overline{\constant}}r\sin{(2\phistar)}}{8r^2+2r\sqrt{\overline{\constant}}} \text{, if }\overline{\constant}>0,\\
        &\big|\overline{\mathcal{X}}_1-\mathcal{X}_1\big|<0.000225r\sin^2{(2\phistar)},\quad\big|\overline{\mathcal{Y}}_1-\mathcal{Y}_1\big|<\frac{0.000225}{2}r\sin^2{(2\phistar)},\quad \big|\Theta_1-\overline{\Theta}_1\big|<0.005\sin{(2\phistar)}\text{, if }\overline{\constant}=0,\\
    \end{aligned}\right.
\end{equation}

In \cref{fig:Trajectory_LemonCompareWithWojtkowski2} the tangent direction at $\overline{P}$ of $\Gamma_r$ has angle $\overline{\phi}_2$ and $\overrightarrow{\overline{P}\,\overline{Q}}$ has angle $\overline{\Theta}_1$, the collision angle $\overline{\theta}_2=\overline{\phi}_2-\overline{\Theta}_1$. The tangent direction at $P$ of $\Gamma_r$ has angle $\phi_2$ and $\overrightarrow{PQ}$ has angle $\Theta_1$, the collision angle $\theta_2=\phi_2-\Theta_1$. we get the following estimate for $|\theta_2-\overline{\theta}_2|$. \begin{equation}\label{eq:theta2thetabar2diff}
    \begin{aligned}
    \big|\theta_2-\overline{\theta}_2\big|&=\bigm|\phi_2-\Theta_1-(\overline{\phi}_2-\overline{\Theta}_1)\bigm|=\bigm|\phi_2-\overline{\phi}_2+\overline{\Theta}_1-\Theta_1\big|\\
    &\le\bigm|\phi_2-\overline{\phi}_2\bigm|+\bigm|\overline{\Theta}_1-\Theta_1\big|\\
    &\overbracket{\le}^{\mathllap{\eqref{eq:phi2phibar2diff}}}\frac{3\pi}{2r}\bigm|\mathcal{X}_1-\overline{\mathcal{X}}_1\bigm|+\frac{3\pi}{2r}\bigm|\mathcal{Y}_1-\overline{\mathcal{Y}}_1\bigm|+(4\pi+1)\bigm|\Theta_1-\overline{\Theta}_1\bigm|+\frac{3\pi}{2r}\bigm|t-\overline{t}\big|.
    \end{aligned}
    \end{equation}
We conclude by explaining how \eqref{eq:phi2phibar2diff} and \eqref{eq:theta2thetabar2diff} imply \cref{lemma:Trajectory_approximation_2}.
For $ \varepsilon>0$, by \eqref{eq:ttbarvarepsilonCondition} if $\big|\mathcal{X}_1-\overline{\mathcal{X}}_1\big|$, $\big|\mathcal{Y}_1-\overline{\mathcal{Y}}_1\big|$ and $\big|\Theta_1-\overline{\Theta}_1\big|$ satisfy 
\begin{equation}\label{eq:ttbarvarepsilonCondition2}
    \left\{\begin{aligned}
        &\big|\overline{\mathcal{X}}_1-\mathcal{X}_1\bigm|<\frac{0.3\sqrt{\overline{\constant}}}{6r+3\sqrt{\overline{\constant}}}\frac{1.4r}{3\pi}\varepsilon,\quad \bigm|\mathcal{Y}_1-\overline{\mathcal{Y}}_1\bigm|<\frac{0.3\sqrt{\overline{\constant}}}{6r+3\sqrt{\overline{\constant}}}\frac{1.4r}{3\pi}\varepsilon,\quad \bigm|\Theta_1-\overline{\Theta}_1\bigm|<\frac{0.4\sqrt{\overline{\constant}}}{8r^2+2r\sqrt{\overline{\constant}}}\frac{1.4r}{3\pi}\varepsilon \text{, if }\overline{\constant}>0,\\
        &\big|\overline{\mathcal{X}}_1-\mathcal{X}_1\bigm|<0.09(\frac{1.4}{3\pi})^2r\varepsilon^2,\quad\bigm|\overline{\mathcal{Y}}_1-\mathcal{Y}_1\bigm|<0.045(\frac{1.4}{3\pi})^2r\varepsilon^2,\quad \bigm|\Theta_1-\overline{\Theta}_1\bigm|<\frac{0.14\varepsilon}{3\pi}\text{, if }\overline{\constant}=0,\\
    \end{aligned}\right.
\end{equation} then $|t-\overline{t}|<\frac{1.4r}{3\pi}\varepsilon$.

Furthermore, if $\big|\mathcal{X}_1-\overline{\mathcal{X}}_1\big|$, $\big|\mathcal{Y}_1-\overline{\mathcal{Y}}_1\big|$ and $\big|\Theta_1-\overline{\Theta}_1\big|$ also  satisfy the following 
\begin{equation}\label{eq:X1X1barY1Y1barTheta1Thetabar1Varepsilon}
\begin{aligned}
    \big|\mathcal{X}_1-\overline{\mathcal{X}}_1\bigm|&<\frac{2r}{3\pi}\cdot0.1\varepsilon,\\
    \big|\mathcal{Y}_1-\overline{\mathcal{Y}}_1\bigm|&<\frac{2r}{3\pi}\cdot0.1\varepsilon,\\
    \big|\Theta_1-\overline{\Theta}_1\bigm|&<\frac{1}{4\pi+1}\cdot0.1\varepsilon,
\end{aligned}
\end{equation} then 
\begin{equation}\label{eq:maxphi2phibar2andtheta2theta2bardiff}
    \max{\big\{\bigm|\phi_2-\overline{\phi}_2\bigm|,\bigm|\theta_2-\overline{\theta}_2\bigm|\big\}}\overbracket{\le}^{\mathllap{\eqref{eq:theta2thetabar2diff},\eqref{eq:phi2phibar2diff}}}\underbracket{\frac{3\pi}{2r}\bigm|\mathcal{X}_1-\overline{\mathcal{X}}_1\bigm|}_{\mathclap{\text{\eqref{eq:X1X1barY1Y1barTheta1Thetabar1Varepsilon}: }<0.1\varepsilon}}+\underbracket{\frac{3\pi}{2r}\bigm|\mathcal{Y}_1-\overline{\mathcal{Y}}_1\bigm|}_{\mathclap{\text{\eqref{eq:X1X1barY1Y1barTheta1Thetabar1Varepsilon}: }<0.1\varepsilon}}+\underbracket{(4\pi+1)\bigm|\Theta_1-\overline{\Theta}_1\bigm|}_{\mathclap{\text{\eqref{eq:X1X1barY1Y1barTheta1Thetabar1Varepsilon}: }<0.1\varepsilon}}+\underbracket{\frac{3\pi}{2r}\bigm|t-\overline{t}\bigm|}_{\mathclap{\eqref{eq:ttbarvarepsilonCondition2}: <\frac{3\pi}{2r}\frac{1.4r}{3\pi}\varepsilon=0.7\varepsilon}}<\varepsilon.
\end{equation}
This is \eqref{eq:phi2theta2close} and hence concludes the proof of \cref{lemma:Trajectory_approximation_2}.
\end{proof}
\begin{remark}[Approximation procedure]\label{remark:shadowingConclusions}
For $\forall \varepsilon>0$, suppose $\bar y_1,\cdots, \bar y_n$ is a finite orbit segment of the 1-petal billiard $\mathbb{P}(r,\phistar)$, which are all smooth collisions (not hitting any corner). By \cref{lemma:ApproximationForConsecutiveCollisions,lemma:Trajectory_approximation_1,lemma:Trajectory_approximation_2} $y_{i+1}$ and $\bar y_{i+1}$ being $\varepsilon$-close can be guaranteed by setting $y_{i}$ and $\bar y_{i}$ to be $\delta(\epsilon)$ close enough and setting $R$ large enough. Recursively from $i=n,\cdots,1$, we can find some $\bar\delta$ and $R_*>0$ so that $y_1$, $\bar y_1$ being close enough and $R>R_*$ ensures that all $y_{i}$ and $\bar y_{i}$, $i=1,\cdots,n$ are $\varepsilon$ close enough. Both $\overline{\delta}$ $R_*$ depend on the orbit segment $\bar y_1,\cdots, \bar y_n$.
\end{remark}
\subsection[Closeness of expansion between lemon billiard and 1-petal billiard]{Closeness of expansion between lemon billiard and 1-petal billiard}
    For each $x_0=(\phi_0,\theta_0)\in \Mrout$ and each $\overline{x}_0=(\overline{\phi}_0,\overline{\theta}_0)\in\Mrout$, as long as $R$ is large enough and $x_0$ and $\overline{x}_0$ are close enough, for lemon billiard $\mathbb{L}(r,R,\phistar)$ with its billiard map $\Fc$ and 1-petal billiard $\mathbb{P}(r,\phistar)$ with its billiard map $\bar\Fc_{\phistar}$, we have shown by \cref{lemma:Trajectory_approximation_1,lemma:Trajectory_approximation_2} that with $\Fc^{2}(x_0)\nfd x_2\in\Mrin$, $\bar F^2_{\phistar}(\overline{x}_0)\nfd\overline{x}_2\in\Mrin$. $x_2$ and $\overline{x}_2$ are close.

 Suppose $dx_0\in T_{x_0}M_r$, $d\overline{x}_0\in T_{\overline{x}_0} M_r$, $dx_2=D_{x_0}\Fc^2(dx_0)$ and $d\overline{x}_2=D_{\overline{x}_0}\bar F^2_{\phistar}(d\overline{x}_0)$, $\frac{\|dx_2\|_{\p}}{\|dx_0\|_{\p}}$ the p-metric expansion for the lemon billiard phase space tangent vector $dx_0$ in the first and third quadrants at $x_0$ also has $\mathcal{I}$ from \cref{def:AformularForLowerBoundOfExpansion} the lower bound by \cref{lemma:ExpansionLowerBoundA0A1}.

We now show that as long as $x_0$, $\overline{x}_0$ are close enough, $\mathcal{I}$ is close to $\big|-1-\frac{\overline{\tau}_{0}+\overline{\tau}_{0}-2\overline{d}_0}{\overline{d}_0}\big|$, which is the lower bound of $\frac{\|d\overline{x}_2\|_{\p}}{\|d\overline{x}_0\|_{\p}}$ in \eqref{eq:lowerboundofPetalExpansion}. 



\begin{lemma}[Expansion of $\PET(r,\phistar)$ and $\LEM(r,R,\phistar)$ trajectories]\label{lemma:expansionbetweenPETandLEM} With \cref{def:LemonAndPetTogether}\eqref{item05LemonAndPetTogether}, $\overline{x}_0=(\overline{\phi}_0$, $\overline{\theta}_0)\in\Mrout$ with $\overline{x}_1=\bar F_{\phistar}(\overline{x}_0)\in M_{\text{\upshape f}}$, $\overline{d}_0= r\sin{\overline{\theta}_0}$, $\overline{\tau}_0=|p(\overline{x}_0)p(\overline{x}_1)|$, $\overline{d}_2= r\sin{\overline{\theta}_2}$, $\overline{\tau}_1=|p(\overline{x}_1)p(\overline{x}_2)|$, in the $\mathbb{P}(r,\phistar)$ billiard table in \cref{fig:1-petalExtensionFromCornerAndLemon} with coordinates in \cref{def:StandardCoordinateTable} let $p(\overline{x}_1)=(\overline{\mathcal{X}}_1,\overline{\mathcal{Y}}_1)$. We suppose that $p(\overline{x}_1)=(\overline{\mathcal{X}}_1,\overline{\mathcal{Y}}_1)$ and $(\overline{\mathcal{X}}_1,\overline{\mathcal{Y}}_1,\overline{\Theta}_1)$ (in \cref{def:flowcollsionpositionandangle}) be the position and direction of the $\PET(r,\phistar)$ billiard flow after collision at $p(\overline{x}_1)$.

Also, with \cref{def:fpclf}: $x=(\phi,\theta)\in\Mrout$, $(\Phi_1,\theta_1)=x_1=\Fcal(x)\in M_R$ and $d_0= r\sin{\theta}$, $d_1=R\sin{\theta_1}$, $\tau_0=|p(x_1)p(x)|$,then in the $\mathbb{L}(r,R,\phistar)$ billiard table in \cref{fig:1-petalExtensionFromCornerAndLemon} with coordinates in \cref{def:StandardCoordinateTable} let $p(x_1)=(\mathcal{X}_1,\mathcal{Y}_1)$. We suppose $(\mathcal{X}_1,\mathcal{Y}_1,\Theta_1)$ to be the position and direction of the $\LEM(r,R,\phistar)$ (in \cref{def:flowcollsionpositionandangle}) billiard flow after collision at $p(x_1)$. Let $\overline{\constant}$ be the constant given in \eqref{eq:APetalCollisionConstant} depending on $\overline{x}_1$. Then We have the following.

$\overline{x}_0\in\Mrout$ implies $p(\overline{x}_1)\in\Int\overline{AB}$ in \cref{fig:1-petalExtensionFromCornerAndLemon}, that is, $|\overline{\mathcal{X}}_1|<r\sin\phistar$. For $\forall\varepsilon_2>0$ if \begin{equation}\label{eq:ExpansionClosingCondition_i}
        \begin{aligned}
            \big|\theta-\overline{\theta}_0\bigm|<&\delta_6(\overline{\phi}_0,\overline{\theta}_0,\varepsilon_2,\overline{x}_1)\dfn\min\Bigm\{0.5\sin{\overline{\theta}_0},\,\frac{1}{32}\varepsilon_2\sin^2{\overline{\theta}_0},\,\frac{-\varepsilon_2\sin^3{(\overline{\phi}_0+\overline{\theta}_0)}\sin{\overline{\theta}_0}}{64},\\
            &\mathtt{\delta}_1\Big(\overline{\phi}_0,\,\overline{\theta}_0,\,\frac{0.3\sqrt{\overline{\constant}}}{6r+3\sqrt{\overline{\constant}}}\cdot0.25\varepsilon_2r\sin{\overline{\theta}_0},\,\frac{0.4\sqrt{\overline{\constant}}}{8r^2+2r\sqrt{\overline{\constant}}}\cdot0.25\varepsilon_2r\sin{\overline{\theta}_0}\Big)\Bigm\},\\
            \big|\phi-\overline{\phi}_0\bigm|<&\delta_7(\overline{\phi}_0,\overline{\theta}_0,\varepsilon_2,\overline{x}_1)\dfn\min\Bigm\{\frac{-\varepsilon_2\sin^3{(\overline{\phi}_0+\overline{\theta}_0)}\sin{\overline{\theta}_0}}{64},\,\\
            &\,\mathtt{\delta}_2\Big(\overline{\phi}_0,\,\overline{\theta}_0,\,\frac{0.3\sqrt{\overline{\constant}}}{6r+3\sqrt{\overline{\constant}}}\cdot0.25\varepsilon_2r\sin{\overline{\theta}_0},\,\frac{0.4\sqrt{\overline{\constant}}}{8r^2+2r\sqrt{\overline{\constant}}}\cdot0.25\varepsilon_2r\sin{\overline{\theta}_0}\Big)\Bigm\},\\
            R>&\mathtt{R}_2(\overline{\phi}_0,\overline{\theta}_0,\varepsilon_2,\overline{x}_1)\dfn\max\Bigm\{\,\frac{1}{64\varepsilon^2_2}r,\,r+\frac{64r}{-\varepsilon_2\sin^3{(\overline{\phi}_0+\overline{\theta}_0)}\sin{\overline{\theta}_0}},\\
            &\mathtt{R}_1\Big(\overline{\phi}_0,\,\overline{\theta}_0,\,\frac{0.3\sqrt{\overline{\constant}}}{6r+3\sqrt{\overline{\constant}}}\cdot0.25\varepsilon_2r\sin{\overline{\theta}_0},\,\frac{0.4\sqrt{\overline{\constant}}}{8r^2+2r\sqrt{\overline{\constant}}}\cdot0.25\varepsilon_2r\sin{\overline{\theta}_0}\Big)\Bigm\},
        \end{aligned}
    \end{equation} where $\delta_1$, $\delta_2$ and $\mathtt{R}_1$ are defined in \eqref{eq:x0barx0ClosingCondition}, then $R$, $(\overline{\phi}_0,\overline{\theta}_0)$ and $(\phi,\theta)$ satisfy condition \eqref{eq:xorbitA0Condition}, i.e. $\sin{\theta_1}>\sqrt{4r/R}$. \newline Hence the $\LEM(r,R,\phistar)$ billiard return orbit on $\hat{M}$ in \cref{def:MhatReturnOrbitSegment} with $(\phi_0,\theta_0)=x_0=x=(\phi,\theta)$, $x_1=\Fcal(x)=\Fcal(x_0)$ is in the case (a0) in \eqref{eqMainCases} and $\mathcal{I}$ from \cref{def:AformularForLowerBoundOfExpansion} as function of $x_1$ thus also of $x$ satisfies 
    \begin{equation}
    \begin{aligned}
        \Big|\mathcal{I}-\big(-1-\frac{\overline{\tau}_0+\overline{\tau}_1-2\overline{d}_0}{\overline{d}_0}\big)\Big|<\varepsilon_2
    \end{aligned}
    \end{equation}
\end{lemma}
\begin{remark}\label{remarke:positiveCollisionConstant}
    Since $|\mathcal{X}_1|<r\sin\phistar$, \eqref{eq:APetalCollisionConstant} implies the constant $\overline{\constant}>0$.
\end{remark}
\begin{proof}[Proof of \cref{lemma:expansionbetweenPETandLEM}]
    Conditions in \eqref{eq:ExpansionClosingCondition_i} imply $\overline{\phi}_0$, $\overline{\theta}_0$, $R$ satisfy condition \eqref{eq:xorbitA0Condition} thus $d_1=R\sin{\theta_1}>\sqrt{4rR}$.
        \begin{equation}\label{eq:LemonexpansionIcompareWithPET}
        \begin{aligned}
            \Big|\mathcal{I}-\big(-1-\frac{\overline{\tau}_0+\overline{\tau}_1-2\overline{d}_0}{\overline{d}_0}\big)\Big|&=\Big|\overbracket{-1+\frac{\tau_1}{d_0}\Big[\frac{2(\tau_0-d_0)}{d_1}-\frac{\tau_0+\tau_1-2d_0}{\tau_1}\Big]}^{{\text{\cref{def:AformularForLowerBoundOfExpansion}}}\nfd\mathcal{I}}-\big(-1-\frac{\overline{\tau}_0+\overline{\tau}_1-2\overline{d}_0}{\overline{d}_0}\big)\Big|\\
            &=\Big|\frac{\overline{\tau}_0+\overline{\tau}_1}{\overline{d}_0}-\frac{\tau_0+\tau_1}{d_0}+\frac{\tau_1}{d_0}\frac{2(\tau_0-d_0)}{d_1}\Big|\\
            &\le\bigm|\frac{\overline{\tau}_0+\overline{\tau}_1}{\overline{d}_0}-\frac{\tau_0+\tau_1}{d_0}\bigm|+\underbracket{\bigm|\frac{\tau_1}{d_1}\frac{2(\tau_0-d_0)}{d_0}\bigm|}_{\substack{<2\sqrt{r/R}\\\text{since }\tau_1\le2r,d_1=R\sin\theta_1>\sqrt{4rR},\\\text{\cref{remark:deffpclf}:}-d_0<\tau_0-d_0<d_0}}\\
            &<\frac{1}{d_0\overline{d}_0}\bigm|\overline{d}_0(\tau_0+\tau_1)-d_0(\overline{\tau}_0+\overline{\tau}_1)\bigm|+2\sqrt{r/R}\\
            &=\frac{1}{d_0\overline{d}_0}\bigm|\overline{d}_0(\tau_0+\tau_1)-d_0(\tau_0+\tau_1)+d_0(\tau_0+\tau_1)-d_0(\overline{\tau}_0+\overline{\tau}_1)\bigm|+2\sqrt{r/R}\\
            &\le\frac{1}{d_0\overline{d}_0}\bigm|\tau_0+\tau_1\bigm|\bigm|d_0-\overline{d}_0\bigm|+\frac{1}{\overline{d}_0}\bigm|\tau_0-\overline{\tau}_0\bigm|+\frac{1}{\overline{d}_0}\bigm|\tau_1-\overline{\tau}_1\bigm|+2\sqrt{r/R}
        \end{aligned}                    
        \end{equation} Note that the same as in \cref{remark:sinPhi0plusTheta0isNegative,lemma:Trajectory_approximation_1}, $(\phi_0,\theta_0)$ satisfies $\phi_0+\theta_0\in(\pi,2\pi)$.
    
    For the first term on the right-hand side of the last inequality in \eqref{eq:LemonexpansionIcompareWithPET}, if $x=(\phi,\theta)$ satisfies the following \begin{equation}\label{eq:FirstTermInLemonexpansionIcompareWithPET}
        \begin{aligned}
            \big|\theta-\overline{\theta}_0\bigm|<&\min\Bigm\{0.5\sin{\overline{\theta}_0},\,\frac{1}{32}\varepsilon_2\sin^2{\overline{\theta}_0}\Bigm\},
        \end{aligned}
    \end{equation} then
    \begin{equation}\label{eq:1sttermexpansioncomparison}
    \begin{aligned}
        \big|d_0-\overline{d}_0\bigm|&=\bigm|r\sin{\theta}-r\sin{\overline{\theta}_0}\bigm|\overbracket{\le}^{\mathclap{\text{Mean-Value Theorem}}}r\bigm|\theta-\overline{\theta}_0\bigm|\overbracket{<}^{\mathclap{\eqref{eq:FirstTermInLemonexpansionIcompareWithPET}}}0.5r\sin{\overline{\theta}_0}\Longrightarrow r\sin{\theta}>0.5r\sin{\overline{\theta}_0}\\
        \Longrightarrow\frac{\bigm|\tau_0+\tau_1\bigm|\bigm|d_0-\overline{d}_0\bigm|}{d_0\overline{d}_0}&=\frac{\tau_0+\tau_1}{d_0\overline{d}_0}\bigm|d_0-\overline{d}_0\bigm|\le\frac{4r}{\underbracket{d_0}_{\mathllap{>0.5r\sin{\overline{\theta}_0}}}\underbracket{\overline{d}_0}_{\mathrlap{=r\sin{\overline{\theta}_0}}}}\bigm|d_0-\overline{d}_0\bigm|<\frac{8r}{r^2\sin^2{\overline{\theta}_0}}\bigm|d_0-\overline{d}_0\bigm|\\
        \overbracket{<}^{\mathclap{\text{Mean-Value Theorem}\qquad\qquad}}&\frac{8r^2}{r^2\sin^2{\overline{\theta}_0}}\bigm|\theta-\overline{\theta}_0\bigm|\overbracket{<}^{\mathclap{\eqref{eq:FirstTermInLemonexpansionIcompareWithPET}}}\frac{\varepsilon_2}{4}\\
    \end{aligned}
    \end{equation}

    For the second term on the right-hand side of the last inequality in \eqref{eq:LemonexpansionIcompareWithPET}, if 
\begin{equation}\label{eq:secondtermtau0taubar0closecondition}
        \begin{aligned}
            \big|\theta-\overline{\theta}_0\bigm|<&\min\big\{0.01\sin^2{(\overline{\phi}_0+\overline{\theta}_0)},\,\frac{-\varepsilon_2\sin^3{(\overline{\phi}_0+\overline{\theta}_0)}\sin{\overline{\theta}_0}}{64}\big\}\\
            \bigm|\phi-\overline{\phi}_0\bigm|<&\min\big\{0.01\sin^2{(\overline{\phi}_0+\overline{\theta}_0)},\,\frac{-\varepsilon_2\sin^3{(\overline{\phi}_0+\overline{\theta}_0)}\sin{\overline{\theta}_0}}{64}\big\}\\
            R&>\max\big\{r+\frac{500r}{\sin^2(\overline{\phi}_0+\overline{\theta}_0)},1700r,r+\frac{64r}{-\varepsilon_2\sin^3{(\overline{\phi}_0+\overline{\theta}_0)}\sin{\overline{\theta}_0}}\big\},\\
    \end{aligned}
    \end{equation} then \eqref{eq:thetaphiclose} and thus \eqref{eq:etalowerbound}, \eqref{eq:etaaproximation} hold.  Notice that in \cref{fig:Trajectory_LemonCompareWithWojtkowski}, $\overrightarrow{T_0\bar P}$ has directional angle $\overline{\phi}_0+\overline{\theta}_0\in(0,2\pi)$ by \cref{remark:sinPhi0plusTheta0isNegative} and the distance from $T_0$ to line $\overline{AB}$ is $-r\cos{\overline{\phi}_0}+r\cos\phistar$. Therefore, by elementary geometry and \eqref{eq:PTposRoot} we get the following. \begin{align*}
    \tau_0&=t_+\overbracket{=}^{\mathclap{\eqref{eq:PTposRoot}}}\frac{2r\cos\phistar-2r\cos\phi}{-\sin{(\phi+\theta)}-(r/b)\sin\theta+\sqrt{\big(\sin{(\phi+\theta)}+(r/b)\sin\theta\big)^2+\big(2(r/b)\cos\phistar-2(r/b)\cos\phi\big)}}\nfd\frac{\upsilon(\phi,\theta)}{\eta(\phi,\theta)},\\
        \overline{\tau}_0&\overbracket{=}^{\mathclap{\text{elementary geometry}}}\quad\qquad\frac{r\cos\phistar-r\cos{\overline{\phi}_0}}{-\sin{(\overline{\phi}_0+\overline{\theta}_0)}}=\frac{2r\cos\phistar-2r\cos{\overline{\phi}_0}}{-2\sin{(\overline{\phi}_0+\overline{\theta}_0)}},
    \end{align*}\newline where $\upsilon(\phi,\theta)=2r\cos\phistar-2r\cos\phi$, and $\eta(\phi,\theta)$ is defined from \eqref{eq:xietaphitheta}\[\eta(\phi,\theta)=-\sin{(\phi+\theta)}-(r/b)\sin\theta+\sqrt{\big(\sin{(\phi+\theta)}+(r/b)\sin\theta\big)^2+\big(2(r/b)\cos\phistar-2(r/b)\cos\phi\big)}.\] And by the Mean-Value Theorem we have
    \begin{equation}\label{eq:mvtforvarsigma}
        \begin{aligned}
            &\big|\upsilon(\phi,\theta)-(2r\cos{\phistar}-2r\cos{\overline{\phi}_0})\bigm|=2r\bigm|\cos\phi-\cos{\overline{\phi}_0}\bigm|\overbracket\le^{\mathclap{\text{Mean-Value Theorem}}}2r\bigm|\phi-\overline{\phi}_0\bigm|.
        \end{aligned}
    \end{equation} And \begin{equation}\label{eq:2ndtermexpansioncomparison}
    \begin{aligned}
        \big|\tau_0-\overline{\tau}_0\bigm|&=\Bigm|\frac{\upsilon(\phi,\theta)}{\eta(\phi,\theta)}-\frac{2r(\cos\phistar-\cos{\overline{\phi}_0})}{-2\sin{(\overline{\phi}_0+\overline{\theta}_0)}}\Big|\\
        &=\Bigm|\frac{\upsilon(\phi,\theta)\big[-2\sin{(\overline{\phi}_0+\overline{\theta}_0)}\big]-\eta(\phi,\theta\big[2r(\cos{\phistar}-\cos{\overline{\phi}_0})\big]}{(-2\sin{(\overline{\phi}_0+\overline{\theta}_0)})\eta(\phi,\theta)}\Big|\\
        &=\Big|\frac{\upsilon(\phi,\theta)\big[-2\sin{(\overline{\phi}_0+\overline{\theta}_0)}\big]-\upsilon(\phi,\theta)\eta(\phi,\theta)+\upsilon(\phi,\theta)\eta(\phi,\theta)-\eta(\phi,\theta)[2r(\cos{\phistar}-\cos{\overline{\phi}_0})]}
        {-2\sin{(\overline{\phi}_0+\overline{\theta}_0)}\eta(\phi,\theta)}\Big|\\
        \le&\frac{\overbracket{|\upsilon(\phi,\theta)|}^{\le4r}}{-2\sin{(\overline{\phi}_0+\overline{\theta}_0)}\underbracket{|\eta(\phi,\theta)|}_{\mathclap{\eqref{eq:etalowerbound}:\ >-1.9258\sin{(\overline{\phi}_0+\overline{\theta}_0)}>0}}}\underbracket{\Bigm|-2\sin{(\overline{\phi}_0+\overline{\theta}_0)}-\eta(\phi,\theta)\Bigm|}_{\text{LHS of \eqref{eq:etaaproximation}}}+\frac{\overbracket{\bigm|\upsilon(\phi,\theta)-(2r\cos{\phistar}-2r\cos{\overline{\phi}_0})\bigm|}^{\text{LHS of \eqref{eq:mvtforvarsigma}}}}{-2\sin{(\overline{\phi}_0+\overline{\theta}_0)}}\\
        \overbracket{<}^{\mathclap{\text{\eqref{eq:mvtforvarsigma}, \eqref{eq:etaaproximation} and \eqref{eq:etalowerbound}}}}&\frac{2r}{1.9258\sin^2{(\overline{\phi}_0+\overline{\theta}_0)}}\Big[\frac{7.1r}{-b\sin{(\overline{\phi}_0+\overline{\theta}_0)}}+\frac{3|\theta-\overline{\theta}_0|}{-\sin{(\overline{\phi}_0+\overline{\theta}_0)}}+\frac{3|\phi-\overline{\phi}_0|}{-\sin{(\overline{\phi}_0+\overline{\theta}_0)}}\Big]+\frac{2r\bigm|\overline{\phi}_0-\phi\bigm|}{-2\sin{(\overline{\phi}_0+\overline{\theta}_0)}}\\
        \underbracket{<}_{\mathllap{\frac{1}{-\sin{(\overline{\phi}_0+\overline{\theta}_0)}}\ge1}}&\underbracket{\frac{5.2r}{-\sin^3{(\overline{\phi}_0+\overline{\theta}_0)}}\big|\phi-\overline{\phi}_0\big|}_{\text{\eqref{eq:secondtermtau0taubar0closecondition}:}<\frac{5.2\varepsilon_2}{64}r\sin{\overline{\theta}_0}}+\underbracket{\frac{3.2r}{-\sin^3{(\overline{\phi}_0+\overline{\theta}_0)}}\big|\theta-\overline{\theta}_0\big|}_{\text{\eqref{eq:secondtermtau0taubar0closecondition}:}<\frac{3.2\varepsilon_2}{64}r\sin{\overline{\theta}_0}}+\frac{7.4r^2}{\underbracket{-b\sin^3{(\overline{\phi}_0+\overline{\theta}_0)}}_{\mathclap{\quad\text{\eqref{eq:secondtermtau0taubar0closecondition}: } b>R-r>\frac{64r}{-\varepsilon_2\sin^3{(\overline{\phi}_0+\overline{\theta}_0)}\sin{\overline{\theta}_0}}}}}\\\overbracket{<}^{\mathclap{\text{\eqref{eq:secondtermtau0taubar0closecondition}:}}}&\big(5.2/64+3.2/64+7.4/64\big)r\sin{\overline{\theta}_0}<\frac{\varepsilon_2}{4}\overline{d}_0.\\
    \end{aligned}
    \end{equation}\newline Therefore, $\frac{|\tau_0-\overline{\tau}_0|}{\overline{d}_0}<\frac{\varepsilon_2}{4}$.
    
    For the third term on the right-hand side of the last inequality in \eqref{eq:LemonexpansionIcompareWithPET}, if 
    \begin{equation}\label{eq:thirdtermclosingcondition}
        \begin{aligned}
            |\theta-\overline{\theta}_0|<&\mathtt{\delta}_1\Big(\overline{\phi}_0,\,\overline{\theta}_0,\,\frac{0.3\sqrt{\overline{\constant}}}{6r+3\sqrt{\overline{\constant}}}\cdot0.25\varepsilon_2r\sin{\overline{\theta}_0},\,\frac{0.4\sqrt{\overline{\constant}}}{8r^2+2r\sqrt{\overline{\constant}}}\cdot0.25\varepsilon_2r\sin{\overline{\theta}_0}\Big),\\
            |\phi-\overline{\phi}_0|<&\mathtt{\delta}_2\Big(\overline{\phi}_0,\,\overline{\theta}_0,\,\frac{0.3\sqrt{\overline{\constant}}}{6r+3\sqrt{\overline{\constant}}}\cdot0.25\varepsilon_2r\sin{\overline{\theta}_0},\,\frac{0.4\sqrt{\overline{\constant}}}{8r^2+2r\sqrt{\overline{\constant}}}\cdot0.25\varepsilon_2r\sin{\overline{\theta}_0}\Big),\\
            R&>\mathtt{R}_1\Big(\overline{\phi}_0,\,\overline{\theta}_0,\,\frac{0.3\sqrt{\overline{\constant}}}{6r+3\sqrt{\overline{\constant}}}\cdot0.25\varepsilon_2r\sin{\overline{\theta}_0},\,\frac{0.4\sqrt{\overline{\constant}}}{8r^2+2r\sqrt{\overline{\constant}}}\cdot0.25\varepsilon_2r\sin{\overline{\theta}_0}\Big),
        \end{aligned}
    \end{equation} with $\mathtt{\delta}_1$, $\mathtt{\delta}_2$, $\mathtt{R}_1$ given in \eqref{eq:x0barx0ClosingCondition}, then \cref{lemma:Trajectory_approximation_1} conclude \eqref{eq:eq:x0barx0Closing} that the conditions in \eqref{eq:thirdtermclosingcondition} with  \cref{remarke:positiveCollisionConstant}: $\overline{\constant}>0$ imply \begin{equation}\label{eq:X1Y1Theta1closingforTau1}\begin{aligned}
        \big|\mathcal{X}_1-\overline{\mathcal{X}}_1\bigm|<&\frac{0.3\sqrt{\overline{\constant}}}{6r+3\sqrt{\overline{\constant}}}\cdot0.25\varepsilon_2r\sin{\overline{\theta}_0},\\
        \big|\mathcal{Y}_1-\overline{\mathcal{Y}}_1\bigm|<&\frac{0.3\sqrt{\overline{\constant}}}{6r+3\sqrt{\overline{\constant}}}\cdot0.25\varepsilon_2r\sin{\overline{\theta}_0},\\
        \big|\Theta_1-\overline{\Theta}_1\bigm|<&\frac{0.4\sqrt{\overline{\constant}}}{8r^2+2r\sqrt{\overline{\constant}}}\cdot0.25\varepsilon_2r\sin{\overline{\theta}_0}.\\
    \end{aligned}\end{equation}\newline And
\eqref{eq:ttbarvarepsilonCondition} concludes that the conditions of \eqref{eq:X1Y1Theta1closingforTau1} imply \begin{equation}\label{eq:3rdtermexpansioncomparison}
    \begin{aligned}
        |\tau_1-\overline{\tau}_1|=\overbracket{|t-\overline{t}|}^{\text{in \eqref{eq:ttbarvarepsilonCondition}}}<&\frac{\varepsilon_2}{4}r\sin{\overline{\theta}_0}=\frac{\varepsilon_2}{4}\overline{d}_0,\\
        \frac{|\tau_1-\overline{\tau}_1|}{\overline{d}_0}<\frac{\varepsilon_2}{4}.
    \end{aligned}
\end{equation}
Therefore, the conditions in \eqref{eq:thirdtermclosingcondition} imply that $\frac{|\tau_1-\overline{\tau}_1|}{\overline{d}_0}<\frac{\varepsilon_2}{4}$.

It is clear that if \begin{equation}\label{eq:fourthtermclosingcondition}
    R>\frac{1}{64\varepsilon^2_2}r, 
\end{equation} then the fourth term on the right-hand side of the last inequality in \eqref{eq:LemonexpansionIcompareWithPET}: $2\sqrt{r/R}<\frac{\varepsilon_2}{4}$. 

The conditions in \eqref{eq:ExpansionClosingCondition_i} imply all conditions in \eqref{eq:FirstTermInLemonexpansionIcompareWithPET}, \eqref{eq:secondtermtau0taubar0closecondition}, \eqref{eq:thirdtermclosingcondition}, \eqref{eq:fourthtermclosingcondition} to establish the following inequality.
\begin{align*}
            \Big|\mathcal{I}-\big(-1-\frac{\overline{\tau}_0+\overline{\tau}_1-2\overline{d}_0}{\overline{d}_0}\big)\Big|\overset{\eqref{eq:LemonexpansionIcompareWithPET}}{<}&\underbracket{\frac{1}{d_0\overline{d}_0}\bigm|\tau_0+\tau_1\bigm|\bigm|d_0-\overline{d}_0\bigm|}_{\eqref{eq:1sttermexpansioncomparison}:<0.25\varepsilon_2}+\underbracket{\frac{1}{\overline{d}_0}\bigm|\tau_0-\overline{\tau}_0\bigm|}_{\eqref{eq:2ndtermexpansioncomparison}:<0.25\varepsilon_2}+\underbracket{\frac{1}{\overline{d}_0}\bigm|\tau_1-\overline{\tau}_1\bigm|}_{\eqref{eq:3rdtermexpansioncomparison}:<0.25\varepsilon_2}+\underbracket{2\sqrt{r/R}}_{\eqref{eq:fourthtermclosingcondition}:<0.25\varepsilon_2}<\varepsilon_2.\qedhere\\
\end{align*}

\end{proof}

\subsection[Uniform expansion on the section (proof of the main theorems)]{Uniform expansion on the section (proof of the main theorems)}\label{subsec:UniformExpansionOnSection}

Now we start to derive $R_{\text{\upshape HF}}$ in \cref{def:AcomputationForRHF} for which we obtain uniform hyperbolicity. It depends on how closely we need to approximate an $\mathbb{L}(r,R,\phistar)$ trajectory by a $\mathbb{P}(r,\phistar)$ trajectory, and it is certainly not optimal.

\begin{definition}[$\overline{l}$ the \emph{sufficient number of returns} and $\phistar$ \emph{feasible parameter} for $\PET(r,\phistar)$]\label{def:ApettrajectoryExpansion}\label{def:ExceptionPhistar}
For a 1-petal billiard $\PET(r,\phistar)$ (\cref{fig:1-petalExtensionFromCornerAndLemon,def:LemonAndPetTogether}\eqref{item01LemonAndPetTogether}) with $\phistar\in(0,\tan^{-1}(1/3))$ and $\pi/\phistar\notin\mathbb{N}$, we call $\overline{l}$ the \emph{sufficient number of returns} for $\PET(r,\phistar)$ to gain $25$ times expansion and we call $\phistar$ a \emph{feasible parameter} for $\PET(r,\phistar)$
if there is a finite orbit segment $y_{B_1},\cdots\bar F_{\phistar}^{n}(y_{B_1})$ starting at  $y_0=y_{B_1}\dfn(3\phistar,\phistar)$ (
the collision $\overrightarrow{BB_1}$ in \cref{fig:1-petalExtensionFromCornerAndLemon}) such that
\begin{enumerate}
    \item\label{item01ApettrajectoryExpansion} the orbit segment never hits corners and has $\overline{l}$ points $\overline{y}_{n_1},\overline{y}_{n_2}\cdots,\overline{y}_{n_{\overline{l}}}$ in $\Mrout$ with $\overline{y}_{n_{k+1}}=\tilde{\bar F}_{\phistar}(\overline{y}_{n_{k}})=\Fc^{n_{k+1}-n_k}(\overline{y}_{n_{k}})$ with $n_0=0$, $n_{k+1}>n_k$ for $k=0,1,\cdots,\overline{l}-1$, where $\tilde{\bar F}_{\phistar}$ is the return map to $\Mrout$ from \cref{def:LemonAndPetTogether}\eqref{item06LemonAndPetTogether};
    
    \item\label{item02ApettrajectoryExpansion} 
    for 
    $\overline{\tau}_{0,k}\dfn\overline{\tau}_0(\overline{y}_{n_k})$, $\overline{\tau}_{1,k}\dfn\overline{\tau}_0(\overline{y}_{n_k})$, $\overline{d}_{0,k}\dfn\overline{d}_0(\overline{y}_{n_k})$,
    $\overline{d}_{2,k}\dfn\overline{d}_2(\overline{y}_{n_k})$ from \cref{def:LemonAndPetTogether}\eqref{item05LemonAndPetTogether} we have the expansion estimate
    \begin{equation}\label{eq:enoughExpansion}
            \prod_{k=1}^{\overline{l}}\bigm|-1-\frac{\overline{\tau}_{0,k}+\overline{\tau}_{1,k}-2\overline{d}_{0,k}}{\overline{d}_{0,k}}\bigm|>25.
    \end{equation} 
    \end{enumerate}
\end{definition}
\begin{remark}[feasible parameter for $\PET(r,\phistar)$]\label{suffcientparameterForPET} With $\pi/\phistar\notin\mathbb{N}$, \cref{lemma:PETExpExpansion} ensures that for a given 1-petal billiard $\PET(r,\phistar)$, $\prod_{k=1}^{l}\bigm|-1-\frac{\overline{\tau}_{0,k}+\overline{\tau}_{1,k}-2\overline{d}_{0,k}}{\overline{d}_{0,k}}\big|$ grows exponentially  with the length of any orbit segment in \cref{def:ApettrajectoryExpansion}\eqref{item01ApettrajectoryExpansion} that avoids corners. Therefore, it is reasonable to conjecture that almost every $\phistar\in(0,\tan^{-1}(1/3))$ is a \emph{feasible parameter} for $\PET(r,\phistar)$ to have $25$ times expansion.
\end{remark}

\begin{remark}[Examples of feasible parameter for $\PET(r,\phistar)$] We have conducted several numerical experiments with an algorithm to check if a parameter $\phistar$ is a \emph{feasible parameter} for $\PET(r,\phistar)$.\newline$\phistar=\pi/13.8$ is a \emph{feasible parameter} with sufficient number of returns $\overline{l}=3$, with the expansion expression in the right-hand side of \eqref{eq:enoughExpansion} being $\approx 58.708988$.\newline  $\phistar=\pi/14.3$ is a \emph{feasible parameter} with sufficient number of returns $\overline{l}=3$, with the expansion expression in the right-hand side of \eqref{eq:enoughExpansion} being $\approx 26.949882$.\newline $\phistar=\pi/15.5$ is a \emph{feasible parameter} with sufficient number of returns $\overline{l}=3$, with the expansion expression in the right-hand side of \eqref{eq:enoughExpansion} being $\approx 44.100198$.\newline $\phistar=\pi/16.935$ is a \emph{feasible parameter} with sufficient number of returns $\overline{l}=2$, with the expansion expression in the right-hand side of \eqref{eq:enoughExpansion} being $\approx 25.615915$.\newline $\phistar=\pi/28.962$ is a \emph{feasible parameter} with sufficient number of returns $\overline{l}=2$, with the expansion expression in the right-hand side of \eqref{eq:enoughExpansion} being $\approx 27.953443$.\newline $\phistar=\pi/36.9$ is a \emph{feasible parameter} with sufficient number of returns $\overline{l}=2$, with the expansion expression in the right-hand side of \eqref{eq:enoughExpansion} being $\approx 26.712283$.\newline Being \emph{feasible} is an open condition which means if a parameter $\phistar$ is \emph{feasible}, then so are the parameters in some neighborhood of $\phistar$.
\end{remark}

\begin{remark}[Generalized feasible parameter for $\PET(r,\phistar)$]\label{def:generalizedsufficientparameter} 
    
    We can generalize the orbit segment in \cref{def:ApettrajectoryExpansion}\eqref{item01ApettrajectoryExpansion} to allow containing points hitting corners by using the continuous extension of \cref{notation:PetContinuousExtension}\eqref{item01PetContinuousExtension} and \eqref{item02PetContinuousExtension} to define the next collision of those  hitting corners. By \cref{remark:1-petalContinuousExtension}, $\overline{\tau}_0$, $\overline{\tau}_1$, $\overline{d}_0$ and $\overline{d}_2$ all have continuous extended definitions for points on the boundary of phase space and boundary of $\Mrout$, which represent hitting corner collisions. If such a generalized to corner hitting finite orbit segment satisfies condition \eqref{eq:enoughExpansion}, we also call such $\phistar$ to be the \emph{generalized feasible parameter} for $\PET(r,\phistar)$. We conjecture that  except $\pi/\phistar\in\mathbb{N}$ every $\phistar\in(0,\tan^{-1}(1/3))$ is a generalized feasible parameter for $\PET(r,\phistar)$.
\end{remark}

\begin{remark}[Computation of $R_{\text{\upshape HF}}$]\label{def:AcomputationForRHF} Consider $x_{B_1}=\overline{y}_{n_0=0}$, $\cdots$, $\overline{y}_{n_{\overline{l}}}$ from \cref{def:ApettrajectoryExpansion}. When $1\le k\le\overline{l}$, let $\bar F_{\phistar}(\overline{y}_{n_k})=\overline{x}_k$ and $(\overline{\mathcal{X}}_k,\overline{\mathcal{Y}}_k)=p(\overline{x}_k)$ in \cref{fig:Trajectory_LemonCompareWithWojtkowski2} with coordinate system given in \cref{def:StandardCoordinateTable}, suppose $(\overline{\phi}_k,\overline{\theta}_k)=\overline{y}_{n_k}$ and with $\overline{\tau}_0$, $\overline{\tau}_1$, $\overline{d}_0$ given in \cref{def:LemonAndPetTogether}\eqref{item05LemonAndPetTogether}, let \begin{equation}\label{eq:ExpansionkInPetalBilliards}
    \overline{E}_k\dfn-1-\frac{\overline{\tau}_0+\overline{\tau}_1-2\overline{d}_0}{\overline{d}_0}.
\end{equation}

We initially set $R_{\text{\upshape HF}}=\frac{30000r}{\sin^2\phistar}$. We have the following finitely many steps to iteratively update $R_{\text{\upshape HF}}$ and compute $\varepsilon_k$ by iterating $k$ decreasingly from $\overline{l}$ to $1$.
\begin{equation}\label{eq:RHFcomputation}
\begin{aligned}
        \text{For }&k=\overline{l},\cdots,1\\
        &\text{set }\varepsilon_{k-1,1}=\left\{\begin{aligned}&\min\Big\{\delta_6(\overline{\phi}_k,\overline{\theta}_k,\big|1-0.9^{\frac{1}{\overline{l}}}\bigm|\bigm|\overline{E}_k\bigm|,\overline{x}_k),\delta_7(\overline{\phi}_k,\overline{\theta}_k,\bigm|1-0.9^{\frac{1}{\overline{l}}}\bigm|\bigm|\overline{E}_k\bigm|,r,\overline{x}_k),\varepsilon_k\Big\}\text{ if }k\ne\overline{l},\\&\min\Big\{\delta_6(\overline{\phi}_k,\overline{\theta}_k,\bigm|1-0.9^{\frac{1}{\overline{l}}}\bigm|\bigm|\overline{E}_k\bigm|,\overline{x}_k),\delta_7(\overline{\phi}_k,\overline{\theta}_k,\bigm|1-0.9^{\frac{1}{\overline{l}}}\bigm|\bigm|\overline{E}_k\bigm|,r,\overline{x}_k)\Big\}\text{ if }k=\overline{l},
        \end{aligned}\right.\\
        &\text{set }R_{\text{\upshape HF}}=\max\big\{R_{\text{\upshape HF}},\,r+\mathtt{R}_2(\overline{\phi}_k,\overline{\theta}_k,\bigm|1-0.9^{\frac{1}{\overline{l}}}\bigm|\bigm|\overline{E}_k\bigm|,\overline{\mathcal{X}}_k)\big\},\\
        &\text{set } \varepsilon_{k-1,2}=\left\{\begin{aligned}&\min\big\{\delta_3(\varepsilon_{k},\overline{x}_k),\delta_4(\varepsilon_{k},\overline{x}_k)\big\},\text{ if }k\ne\overline{l},\\&\varepsilon_{k-1,1},\text{ if }k=\overline{l},\end{aligned}\right.\\
        &\text{set }\varepsilon_{k-1,3}=\left\{\begin{aligned}&\delta_5(\varepsilon_{k},\overline{x}_k),\text{ if }k\ne\overline{l},\\&\varepsilon_{k-1,1},\text{ if }k=\overline{l},\end{aligned}\right.\\
        &\text{set }\varepsilon_{k-1,4}=\left\{\begin{aligned}&\min\Big\{\delta_1(\overline{\phi}_k,\overline{\theta}_k,\varepsilon_{k-1,2},\varepsilon_{k-1,3}),\delta_2(\overline{\phi}_k,\overline{\theta}_k,\varepsilon_{k-1,2},\varepsilon_{k-1,3})\Big\},\text{ if }k\ne\overline{l},\\&\varepsilon_{k-1,1},\text{ if }k=\overline{l},\end{aligned}\right.\\
        &\text{set }R_{\text{\upshape HF}}=\max\big\{R_{\text{\upshape HF}},r+\mathtt{R}_1(\overline{\phi}_k,\overline{\theta}_k,\varepsilon_{k-1,2},\varepsilon_{k-1,3})\big\},\\
        &\text{set }\varepsilon_{k-1}=\left\{\begin{aligned}&\delta_0\Big(\bar F^2_{\phistar}(\overline{y}_{n_{k-1}}),\min\big\{\varepsilon_{k-1,1},\varepsilon_{k-1,4}\big\}\Big)\text{, if }k>1,\\
        &\delta_0\Big(Iy_*,\min\big\{\varepsilon_{k-1,1},\varepsilon_{k-1,4}\big\}\Big) \text{, if }k=1,
        \end{aligned}\right.
    \end{aligned}
\end{equation}
where $\delta_6$, $\delta_7$, $\mathtt{R}_2$ are defined in \cref{lemma:expansionbetweenPETandLEM} \eqref{eq:ExpansionClosingCondition_i}, $\delta_3$, $\delta_4$, $\delta_5$ are defined in \cref{lemma:Trajectory_approximation_2} \eqref{eq:X1Y1Theta1CloseCondition}, $\delta_1$, $\delta_2$, $\mathtt{R}_1$ are defined in \cref{lemma:Trajectory_approximation_1} \eqref{eq:x0barx0ClosingCondition}. $\delta_0$ is defined in \cref{lemma:ApproximationForConsecutiveCollisions} \eqref{eq2:ApproximationForConsecutiveCollisions} \eqref{eq4:ApproximationForConsecutiveCollisions}. The procedure \eqref{eq:RHFcomputation} terminates with an updated $R_{\text{\upshape HF}}$ and a precision $\varepsilon_0$. 
We then finally update $R_{\text{\upshape HF}}=\max\big\{R_{\text{\upshape HF}}, \underbracket{\frac{20.25r}{\varepsilon^2_0}}_{\mathclap{R>20.25r/\varepsilon^2_0\text{ guarantees each point in }\Nout \text{ is }\varepsilon_0-\text{close to either } (\phistar,\phistar)\text{ or }(2\pi-\phistar,\pi-\phistar)} }\big\}$. (see \cref{def:SectionSetsDefs} defined bad regions $\Nout$ and $\delta_*$)
\end{remark}

\begin{remark}\label{remark:GaranteeLargeTheta}
    By making $R_{\text{\upshape HF}}(r,\phistar)$ as large as in \cref{def:AcomputationForRHF} and by \cref{lemma:Trajectory_approximation_1}\eqref{eq:xorbitA0Condition} we see that we have ensured that if $R>R_{\text{\upshape HF}}(r,\phistar)$ for each $x\in\Nin$ with $\Nin$ as the "Bad" region given in \cref{def:SectionSetsDefs} and as the neighborhoods of $(\phistar,\phistar)$ and $(2\pi-2\phistar,\pi-\phistar)$, there exists some $j(x)\ge0$ such that $p(\Fc^i(x))\in\Gamma_r$, $i=0,\cdots,j(x)$ and $p(\Fc^{j(x)}(x))\in\Gamma_R$ with $\Fc^{j(x)}(x)=(\Phi,\theta)\in M_R$ having $\theta>\sin^{-1}{(\sqrt{r/R})}$. Therefore, orbits starting from $\Nin$ ( $\Nout$ ) cannot have 2 consecutive returns to $\Mrin$ ( $\Mrout$ ) being returns to $\Nin$ ( $\Nout$ ).
\end{remark}
\begin{proposition}[$R_{\text{\upshape{HF}}}$ ensures expansion for orbits starting from $\Fc^{-1}(\Nout)$]\label{proposition:ExpansionWins}\hfill

    Suppose that the lemon billiard $\mathbb{L}(r,R,\phistar)$ has $\phistar\in(0,\tan^{-1}{(1/3)})$ and $\phistar$ is the feasible parameter in \cref{def:ApettrajectoryExpansion} for 1-petal billiard $\PET(r,\phistar)$. With $R_{\text{\upshape HF}}(r,\phistar)>0$ computed in \cref{def:AcomputationForRHF}, suppose $R$ satisfies \eqref{eqJZHypCond}: $R>R_{\text{\upshape HF}}(r,\phistar)$.
     
    For $x\in\hat{M}$ with its return orbit segment being in cases (a1)(b)(c) of \eqref{eqMainCases}, $x_0=\Fc(x)$, $x_1$, $x_2$, $dx_2$ and $dx_3=D_{x_2}\Fc(dx_2)$ are respectively defined in \cref{lemma:ExpansionLowerBoundA0A1,lemma:ExpansionLowerBoundB,lemma:AformularForLowerBoundOfExpansionForn1largerthan0}. We suppose $dx\in C_x$ given in \cref{RMASUH}. For the return map $\hat{F}$ on $\hat{M}$ given in \cref{def:SectionSetsDefs}, we have the following.
    \[\frac{\|D_x\hat{F}^{\overline{l}}(dx)\|_{\p}}{\|dx_3\|_{\p}}>25\times0.9=22.5.\]
\end{proposition}
\begin{proof}
     \cref{def:ApettrajectoryExpansion,def:AcomputationForRHF,lemma:Trajectory_approximation_1,lemma:Trajectory_approximation_2,lemma:expansionbetweenPETandLEM}
     ensure that for $\hat{F}$ iterates, after the cases (a1)(b)(c) orbit segment (in \eqref{eqMainCases}) there will have at least $\overline{l}$ times case (a0) return orbit segment (in \eqref{eqMainCases}) following. Suppose that in the $j$-th following return orbit segment, the $x_1$ in \cref{def:MhatReturnOrbitSegment} is $x_{1,j}=(\Phi_{1,j},\theta_{1,j})$, for $j=1,\cdots,\overline{l}$. Being in case (a0) means $\sin\theta_{1,j}>\sqrt{4r/R}$ with $\mathcal{I}_j$ given in \cref{def:AformularForLowerBoundOfExpansion} so that by \cref{lemma:ExpansionLowerBoundA0A1} $|\mathcal{I}_j|$ is the lower bound of expansion between the $(j-1)$-th and $j-$th return.

     With $\overline{E}_j$ given in \eqref{eq:ExpansionkInPetalBilliards} of  \cref{def:AcomputationForRHF}, \cref{lemma:expansionbetweenPETandLEM} ensures that the $j$-th following orbit segment of $\hat{M}$ return has $\mathcal{I}_j$  satisfying $\big|\mathcal{I}_j-\overline{E}_j\bigm|<\bigm|1-0.9^{\frac{1}{\overline{l}}}\bigm|\bigm|\overline{E}_j\big|$ so that $\big|\mathcal{I}_j\bigm|\ge\bigm|\overline{E}_j\bigm|-\bigm|\overline{E}_j-\mathcal{I}_j\bigm|>(1-|1-0.9^{\frac{1}{\overline{l}}}|)\bigm|\overline{E}_j\bigm|=0.9^{\frac{1}{\overline{l}}}\bigm|\overline{E}_j\big|$, for $j=1,\cdots \overline{l}$. Therefore, \[\frac{\|D_x\hat{F}^{\overline{l}}(dx)\|_{\p}}{\|dx_3\|_{\p}}\overbracket{>}^{\text{\cref{lemma:ExpansionLowerBoundA0A1}}}\prod_{j=1}^{\overline{l}}\bigm|\mathcal{I}_j\bigm|>\prod_{j=1}^{\overline{l}}\bigm|(0.9)^{\frac{1}{\overline{l}}}E_j\bigm|\overbracket{>}^{\mathrlap{\text{\cref{def:ApettrajectoryExpansion}\eqref{item02ApettrajectoryExpansion}}}}25\times0.9=22.5\qedhere\]
\end{proof}
\begin{definition}[Uniformly hyperbolic lemon billiard]\label{def:UniformHyperbolicLemonBilliard}
    Now we define our \emph{uniform hyperbolic lemon billiard} to be the $\mathbb{L}(r,R,\phistar)$ with $\phistar\in(0,\tan^{-1}{(1/3)})$, $\phistar$ being the \emph{feasible parameter} in \cref{def:ApettrajectoryExpansion}, and with $R_{\text{\upshape HF}}(r,\phistar)>0$ given in \cref{def:AcomputationForRHF}.
    
\end{definition}
\begin{theorem}[\cref{RMASUH}, uniform expansion of $\hat F$]\label{RMASUHwithProof}\label{thm:UniformExpansionOnMhat}
For the lemon billiards $\LEM(r,R,\phistar)$ satisfying \cref{def:UniformHyperbolicLemonBilliard} and for the return map $\hat{F}$ to $\hat{M}$ (\cref{def:Mhatdefinition}), we have the following. 
\begin{enumerate}
    \item\label{item01RMASUHwithProof} The cone family $C_x\dfn\big\{(d\phi,d\theta)\bigm|\frac{d\theta}{d\phi}\in [0,1]\big\}$ is strictly invariant under the return map $\hat{F}$ from \cref{def:SectionSetsDefs}.
    \item\label{item02RMASUHwithProof}If $x\in\hat{M}$ is nonsingular, $dx\in C_x$ and $n\in\mathbb{N}$, then $\dfrac{\|D_x\hat{F}^n(dx)\|_\p}{\|dx\|_\p}>c\Lambda^n$, where $c=(0.05)/\overline{\Lambda}$, $\Lambda=(\overline{\Lambda})^{\frac{1}{1+\overline{l}}}$.
\end{enumerate}
\end{theorem}
\begin{proof}
    The first conclusion that $C_x$ is strictly invariant is proved by \cref{corollary:InvQuadrant}, \eqref{eq:strictlyinvcone_casea1bc}, and \cref{corollary:halfquadrantcone}. 
    
    For the second conclusion, it remains to show the uniform exponential expansion for nonsingular $x\in\hat{M}$ with its return orbit segment cases \eqref{eqMainCases}. 
    
    In case (a0) of  \eqref{eqMainCases}, $\frac{\|D_x\hat{F}(dx)\|_{\p}}{\|dx\|_{\p}}>1+\lambda_c$ with constant $\lambda_c>0$ from \cref{TEcase_a0}.

    In cases (a1)(b)(c) of \eqref{eqMainCases}, for the $x_2$ and $dx_2$ defined in \cref{lemma:ExpansionLowerBoundA0A1,lemma:ExpansionLowerBoundB,lemma:AformularForLowerBoundOfExpansionForn1largerthan0}, with $dx_3=D_{x_2}\Fc(dx_2)$, \cref{lemma:ExpansionLowerBoundA0A1,lemma:ExpansionLowerBoundB,lemma:AformularForLowerBoundOfExpansionForn1largerthan0} give $\frac{\|dx_3\|_{\p}}{\|dx\|_{\p}}>0.05$.  Then with our $R_{\text{\upshape HF}}$ in \cref{def:AcomputationForRHF}\eqref{eq:RHFcomputation} and the $\mathbb{L}(r,\phistar)$ trajectory orbit section segment $\big\{y_{n_k}\big\}_{k=0}^{\overline{l}}$ in \cref{def:ApettrajectoryExpansion}, \cref{proposition:ExpansionWins} gives $\frac{\|D_x\hat{F}^{\overline{l}}(dx)\|_{\p}}{\|dx_3\|_{\p}}>25\times0.9=22.5$. Therefore, in cases (a1)(b)(c), $\frac{\|D_x\hat{F}^{\overline{l}+1}(dx)\|_{\p}}{\|dx\|_{\p}}>0.05\times22.5=1.125$.

    Let $\overline{\Lambda}=\min\{1+\lambda_c,1.125\}$, where $\lambda_c$ is from \cref{TEcase_a0}. \cref{def:ApettrajectoryExpansion,def:AcomputationForRHF,lemma:Trajectory_approximation_1}\eqref{eq:xorbitA0Condition} ensure that among every consecutive $\overline{l}+1$ times $\hat{M}$ return orbit segment defined in \eqref{eqMhatOrbSeg} there can exist at most 1 time return orbit segment not in case (a0) of \eqref{eqMainCases}. Therefore, We can conclude that for all nonsingular $x$ and $\forall n\ge1$, \[\frac{\|D_x\hat{F}^{n}(dx)\|_{\p}}{\|dx\|_{\p}}>(\overline{\Lambda})^{\left\lfloor\frac{n}{\overline{l}+1}\right\rfloor}(0.05).\qedhere\] 
\end{proof}
\begin{remark}
    \cref{RMASUHwithProof} is \cref{RMASUH} and we expect that the feasible parameter $\phistar$ given in \cref{def:ApettrajectoryExpansion} has full measure in $(0,\tan^{-1}(1/3))$ (see \cref{suffcientparameterForPET}). We can also generalize the computation in \cref{def:AcomputationForRHF} to compute $R_{\text{\upshape HF}}(r,\phistar)>0$ for the generalized feasible parameter $\phistar$ given in \cref{def:generalizedsufficientparameter}.
\end{remark}

\begin{corollary}[Corollary of \cref{RMASUHwithProof},\cref{RMASUH}, uniform expansion of $\hat F$ for Euclidean metric]\label{corollary:sinpmetricexpansion}\label{corollary:sinpmetricexpansionMhat}
There exist constant $\hat{c}>0$ and $\Lambda>1$ such that $\dfrac{\sin\theta\|D_x\hat{F}^n(dx)\|_\p}{\|dx\|_\p}>\hat{c}\Lambda^n$ for all nonsingular $x=(\phi,\theta)\in\hat{M}$, $dx\in C_x$ and $n\in\mathbb{N}$.
\end{corollary}
\begin{proof}
    Let $x=\hat{x}_0=(\phi,\theta)\in\hat{M}$, $(d\phi,d\theta)=dx\in C_{x}$, $\hat{x}_n=\hat{F}^n(x)$, $d\hat{x}_n=D_x\hat{F}^n(dx)$, $n\ge0$.
    \begin{equation}\label{eq:pexpansion_sin}
\displaystyle\sin\theta\frac{\|d\hat{x}_n\|_\p}{\|dx\|_\p}=\frac{\|d\hat{x}_n\|_\p}{\|d\hat{x}_1\|_\p}\frac{\|d\hat{x}_1\|_\p\sin{\theta}}{\|dx\|_\p}
\end{equation}
For the orbit segment of $x$, defined in \cref{def:MhatReturnOrbitSegment} and \eqref{eqMhatOrbSeg} with points $x_0=(\phi_0,\theta_0)\in\Mrout$, $x_1=(\Phi_1,\theta_1)\in\MRin$ and $x_2\in\Mrin$ in \eqref{eqPtsFromeqMhatOrbSeg}. Note that since for either $x\in\Mrin\cap\Mrout$ or $\Fc^{-1}(\Mrout\smallsetminus\Mrin)$, we have $\theta=\theta_0$ since either $x=x_0$ or $x$, $x_0$ are consecutive collisions on $\Gamma_r$. 

We analyze the two cases classified by $\theta_1$.

Case (i), $\sin{\theta_1}<\sqrt{4r/R}$. This encompasses cases (a1)(b)(c) of \eqref{eqMainCases}. \cref{TEcase_a1,TEcase_b,TEcase_c} imply that $\frac{\|d\hat{x}_1\|_\p}{\|dx\|_\p}>0.05$.

By \cref{contraction_region},  we get $x_0\in\Nout$ and the following.
\[\mid\theta_0-(2\pi-\phistar)\mid<\delta_*=4.5\sqrt{4r/R}\overbracket{<}^{\eqref{eqJZHypCond}}4.5\cdot\sqrt{\sin^2\phistar/30000}<4.5\sqrt{\phistar^2/30000}<0.025\phistar.\]

or \[\mid\theta_0-\phistar\mid<\delta_*=4.5\sqrt{4r/R}\overbracket{<}^{\eqref{eqJZHypCond}}4.5\cdot\sqrt{\sin^2\phistar/30000}<4.5\sqrt{\phistar^2/30000}<0.025\phistar.\]

Since $0<\phistar<\tan^{-1}(1/3)$, we see $\theta_0$ is bounded away from $0$ and $\pi$. Thus, $\sin\theta_0$ is bounded away from $0$. That is, there exists some constant $\constant_0>0$ for all $x_0\in\Nout$, $\sin\theta_0>\constant_0$. Therefore, in this case, there exists $C_0>0$ such that $\frac{\|d\hat{x}_1\|_\p\sin{\theta}}{\|dx\|_\p}=\frac{\|d\hat{x}_1\|_\p\sin{\theta_0}}{\|dx\|_\p}>C_0$. 


Case (ii), $\sin{\theta_1}\ge\sqrt{4r/R}$. This is case (a0) of \eqref{eqMainCases}. Therefore, \cref{TEcase_a0} gives $\frac{\|d\hat{x}_1\|_\p\sin{\theta_0}}{\|dx\|_\p}>\sin\theta_0$. 

On the other hand, by \cref{def:AformularForLowerBoundOfExpansion} and computation in \cref{subsec:AFELB},
\begin{equation}\label{eq:casea0onestepsinexpansion}
\begin{aligned}
\frac{\|d\hat{x}_1\|_\p\sin{\theta_0}}{\|dx\|_\p}&\underbracket{>}_{\mathclap{\substack{\text{\cref{lemma:ExpansionLowerBoundA0A1}}\\\text{\cref{def:AformularForLowerBoundOfExpansion}}}}}\Big|-1+\frac{\tau_1}{d_1}\frac{2(\tau_0-d_0)}{d_0}-\frac{\tau_0+\tau_1-2d_0}{d_0}\Big|\sin\theta_0\\
    & =\Big|\underbracket{-\sin\theta_0}_{\mathllap{\rightarrow0}}+\underbracket{\frac{\tau_1}{d_1}\frac{2(\tau_0-d_0)}{d_0}\sin\theta_0}_{\mathllap{\text{\cref{remark:deffpclf}: }\rightarrow0}}-\underbracket{\frac{\tau_0+\tau_1-2d_0}{d_0}\sin\theta_0}_{\mathclap{\substack{\quad=\frac{\tau_0+\tau_1-2d_0}{r}\text{ since }d_0=r\sin{\theta_0}\\
    \text{Since }d_0\rightarrow0,\rho\rightarrow r\text{, \eqref{eq:NewCoordinateLegnth}: }\theta_1\rightarrow\phistar-\Phistar\text{ or }\pi-\phistar+\Phistar\\\rho\rightarrow r\text{ in \eqref{eq:continuousExtension}: }\tau_0+\tau_1\rightarrow 2r\sin{(2\phistar-2\Phistar)}}}}\Big|\xrightarrow{\sin\theta_0\rightarrow0}2\sin{(2\phistar-2\Phistar)}
\end{aligned}
\end{equation}
If $\sin{\theta_0}\rightarrow0$, then $d_0=r\sin\theta_0\rightarrow0$ and $0\le\tau_0<2d_0\rightarrow 0$. This implies the collision position $p(x_1)$ will be approximating the corners, that is, $p(x_1)\rightarrow A \cup B$.

Thus in equations \eqref{eq:lengthextensionA}, \eqref{eq:lengthextensionB}, $\theta_0\longrightarrow \{0\} \cup \{\pi\}$ will yield $\rho\rightarrow r$, $d_0\rightarrow 0$, $\theta_1\rightarrow\{\phistar-\Phistar\}\cup\{ \pi-\phistar+\Phistar\}$,  $\tau_0\rightarrow0$ and $\tau_1\rightarrow 2r\sin{(2\phistar-2\Phistar)}$.

Since $\frac{\tau_1}{d_1}\in[0,2],\frac{2(\tau_0-d_0)}{d_0}\in[-2,2]$, the right-hand side of the last equality in \eqref{eq:casea0onestepsinexpansion} $\rightarrow 2\sin{(2\phistar-2\Phistar)}$ and is bounded away from 0 as $\sin\theta_0\rightarrow0$. Hence, there exists $C_1>0$ such that $\frac{\|d\hat{x}_1\|_\p\sin{\theta}}{\|dx\|_\p}>C_1$ in case (ii). 

in both case (i) and case (ii), there exists $C=\min\big\{C_0,C_1\big\}>0$ such that $\frac{\|d\hat{x}_1\|_\p\sin{\theta}}{\|dx\|_\p}=\frac{\|d\hat{x}_1\|_\p\sin{\theta_0}}{\|dx\|_\p}>C$. Thus \eqref{eq:pexpansion_sin} gives\COMMENT{Use something like the overbracket construction at the end of this equation in many other places below, such as, \cref{eq:pmetricConvertEuclideanmetric}, end of proof of \cref{thm:MroutreturnUniformExpansion}.\wentao{edited now.}} 
\[\sin\theta\frac{\|d\hat{x}_n\|_\p}{\|dx\|_\p}=\frac{\|d\hat{x}_n\|_\p}{\|d\hat{x}_1\|_\p}\frac{\|d\hat{x}_1\|_\p\sin{\theta}}{\|dx\|_\p}\overbracket{>}^{\mathrlap{\text{\cref{RMASUHwithProof}}}}C\cdot c\Lambda^{n-1}=\overbracket{cC\Lambda^{-1}}^{\mathclap{=\hat{c}}}\Lambda^{n}\qedhere\]
\end{proof}

\begin{corollary}[of \cref{RMASUHwithProof}, \cref{RMASUH} uniform expansion in Euclidean Metric] \label{Corollary:EuclideanExpansionMhat}\label{Col:EuclideanExpansionMhat} Under the return map $\hat{F}$ on $\hat{M}$ the vectors in the unstable cone $\mathcal{HQ}_x(\mathrm{I,III})=C_x$ defined in \cref{RMASUHwithProof} expand uniformly with respect to the Euclidean metric $\|\,\|$, i.e., for nonsingular $x\in\hat{M}$, $dx=(d\phi,d\theta)\in C_x$ and all $n\ge0$, there exist $c_0>0,\Lambda>1$ such that $\dfrac{\|D_x\hat{F}^n(dx)\|}{\|dx\|}>c_0\Lambda^n$.
\end{corollary}
\begin{proof}
    Let $x=(\phi,\theta)\in\hat{M}$, $(d\phi,d\theta)=(d\hat{\phi}_0,d\hat{\theta}_0)=d\hat{x}\in C_{x}$, $\hat{x}_n=\hat{F}^n(\hat{x}_0)$, $(d\hat{\phi}_n,d\hat{\theta}_n)=d\hat{x}_n=D_x\hat{F}^n(dx),\mathcal{V}_n=\frac{d\hat{\theta}_n}{d\hat{\phi}_n}\in[0,1]$, $n\ge0$.
Reasoning as in \cite[equation (8.21)]{cb}, it gives
\begin{equation}\label{eq:pmetricConvertEuclideanmetric}
\displaystyle\frac{\|d\hat{x}_n\|}{\|dx\|}=\frac{\|d\hat{x}_n\|_\p}{\|dx\|_\p}\frac{\sin{\theta}}{\sin{\hat{\theta}_n}}\frac{\sqrt{1+\mathcal{V}^2_n}}{\sqrt{1+\mathcal{V}^2_0}}\ge\sin{\theta}\frac{\|d\hat{x}_n\|_\p}{\|dx\|_\p}\frac{1}{\sqrt{2}}\underbracket{>}_{\text{\cref{corollary:sinpmetricexpansion}}}\frac{\hat{c}}{\sqrt{2}}\Lambda^n\nfd \underbracket{c_0}_{\mathrlap{=\hat{c}/\sqrt{2}}}\Lambda^n.\qedhere
\end{equation}
\end{proof}

\begin{theorem}[Uniform expansion on section $\hat{M}_1$]\label{thm:UniformExpansionOnMhat1} For the lemon billiards $\LEM(r,R,\phistar)$ satisfying \cref{def:UniformHyperbolicLemonBilliard} and for the return map $\hat{F}_1$ on $\hat{M}_1\dfn(\Mrin\cap\Mrout)\sqcup\Fc(\Mrin\smallsetminus\Mrout)$ we have the following.


\begin{enumerate}[label={(\roman*)}]
    \item\label{item01:UniformExpansionOnotherSections}
    For nonsingular $x\in\hat{M}_1$, if $dx\in \hat{C}^{u}_{1,x}\dfn\begin{cases}
        \big\{(d\phi,d\theta)\in T_x\hat{M}_1\bigm|\frac{d\theta}{d\phi}\in[0,1]\big\}\text{, if } x\in \Mrin\cap\Mrout,\\
            \big\{(d\phi,d\theta)\in T_x\hat{M}_1\bigm|\frac{d\theta}{d\phi}\in[0,\hat{\lambda}]\big\} \text{ elsewhere,}
    \end{cases}$\newline where $\hat{\lambda}=\max\{1,\lambda_2\}$ with $\lambda_2$ from \cref{theorem:B3minusdx3dphi3ratiobound}, then $D_x\hat{F}_1(dx)\in\big\{\text{interior of }\hat{C}^{u}_{1,\hat{F}_1(x)}\big\}$.
    \item\label{item02:UniformExpansionOnotherSections} For nonsingular $x\in\hat{M}_1$, $dx\in \hat{C}^u_{1,x}$, there exist some $c_1>0$ and some $\Lambda_1>1$ such that $\frac{\|D_x\hat{F}_1^n(dx)\|_\p}{\|dx\|_\p}>c_1\Lambda_1^n$, $\forall n\ge1$.
    \end{enumerate}
\end{theorem}
\begin{proof}
    We use positive quadrant/half quadrant notation from \cref{corollary:halfquadrantcone}: $\mathcal{Q}_x(\mathrm{I,III})=\big\{(d\phi,d\theta)\in T_x M\bigm|\frac{d\theta}{d\phi}\in[0,+\infty]\big\}$, $\mathcal{HQ}_x(\mathrm{I,III})=\big\{(d\phi,d\theta)\in T_x M\bigm|\frac{d\theta}{d\phi}\in[0,1]\big\}$.
\newline
    \textbf{Proof of (i)} is similar to \cref{corollary:InvQuadrant}. For nonsingular $x\in\hat{M}_1$, $\hat{F}_1$ return orbit segment (similar to \eqref{eqMhatOrbSeg}) is 
    \begin{equation}\label{eqhatF1returnorbit}\underbracket{x}_{\in\hat{M}_1},\underbracket{\Fc(x),\cdots}_{\notin\hat{M}_1},\underbracket{\Fc^{\sigma_1(x)}(x)=\hat{F}_1(x)}_{\in\hat{M}_1}\text{ where }\sigma_1(x)\dfn\inf\big\{k>0\mid \Fc^{k}(x)\in\hat{M}_1\big\}.\end{equation} 
    Note that it is not hard see that $\sigma_1(x)\ge2$. For $\forall dx\in\mathcal{Q}_x(\mathrm{I,III})$, 
    if $x\in\Min_{r,0}\sqcup\overbracket{\Fc(\Min_{r,1})}^{\mathclap{=\Mout_{r,1}}}$, then $x_0\dfn x$, $dx_0\dfn dx$. And if $x\in\Fc(\Min_{r,i})$ for some $i\ge2$, then $\Mrout\ni x_0\dfn\Fc^{m(\Fc^{-1}(x))-1}(x)$ and $dx_0\dfn D_{x}\Fc^{m(\Fc^{-1}(x))-1}(dx)\in \mathcal{Q}_{x_0}(\mathrm{I,III})$ since $D_{x}\Fc^{m(\Fc^{-1}(x))-1}=\begin{pmatrix}
1 & 2(m(\Fc^{-1}(x))-1) \\
0 & 1 
\end{pmatrix}$.

Therefore, the same as \eqref{eqPtsFromeqMhatOrbSeg}, the orbit segment \eqref{eqhatF1returnorbit} contains $x_0$, $x_1=\Fc(x_0)\in\MRin$, $x_2=\Fc^{n_1+1}(x_1)\in\Mrin$ with $n_1\ge0$. We define $dx_1=D_{x_0}\Fc(dx_0)$, $dx_2=D_{x_1}\Fc^{n_1+1}(dx_1)$. There are length functions $\tau_0,\tau_1,d_0,d_1,d_2$ defined in \cref{def:fpclf}, as well as three cases of $d_1$, $n_1$ by \eqref{eqMainCases}.

    For case (a): $d_1\ge2r$, \cref{caseA}: $D_{x_0}\Fc^2$ is a negative matrix. Since $dx_0\in\mathcal{Q}_{x_0}(\mathrm{I,III})$, 
    \begin{equation}\label{eq:Fhat1dx2caseA}
    dx_2=D_{x_0}\Fc^2(dx_0)\in\big\{\text{ interior of }\mathcal{Q}_{x_2}\big\}
    \end{equation} 
    
    
    If $x_2\in\Mrin\cap\Mrout$, then from \cref{def:DF2matrix} and \eqref{eq:DF2matrix} $d_1d_2D_{x_0}\Fc^2=\begin{bmatrix}
    D_{11} & D_{12}\\
    D_{21} & D_{22}
    \end{bmatrix}$ satisfies $0<\frac{D_{21}}{D_{11}}<\frac{D_{22}}{D_{12}}<1$, thus by \eqref{eq:narrowconereason1}, \eqref{eq:narrowconereason2} and \eqref{eq:narrowconereason3} so that We get the following.
    \begin{equation}\label{eq:Mhat1DF2}
        D_{x}\hat{F}_1(dx)=D_{x_0}\Fc^2(dx_0)\in\big\{\text{interior of }\mathcal{HQ}_{x_2}(\mathrm{I,III})\big\}\subset\big\{\text{interior of }\hat{C}^u_{1,\hat{F}_1(x)}\big\}.
    \end{equation}\newline
    
    Otherwise, $x_2\in\Mrin\smallsetminus\Mrout$, then by \eqref{eq:Fhat1dx2caseA} and $D_{x_2}\Fc=\begin{bmatrix}
    1 & 2\\
    0 & 1
    \end{bmatrix}$, we get the following. 
    \begin{equation}
    D_x\hat{F}_1(dx)=D_{x_2}\Fc(\underbracket{dx_2}_{\mathclap{\eqref{eq:Fhat1dx2caseA}: \in \mathcal{Q}_{x_2}(\mathrm{I,III})}})\in\big\{\text{interior of }\mathcal{HQ}_{\hat{F}_1(x)}(\mathrm{I,III})\big\}\subset\big\{\text{interior of }\hat{C}^u_{1,\hat{F}_1(x)}\big\}.\end{equation}

    For cases (b) and (c), per \cref{contraction_region}, $x_0\in\Nout$, $x_2\in\Nin$. Since \cref{Prop:Udisjoint}, $\Nout\cap(\Mout_{r,0}\sqcup\Mout_{r,1}\sqcup\Mout_{r,2})=\emptyset$ and $\Nin\cap(\Min_{r,0}\sqcup\Min_{r,1}\sqcup\Min_{r,2})=\emptyset$. Therefore we get $x\in\Fc(\Mrin\smallsetminus\Mrout)$, $x\neq x_0$, $x_2\notin\Mrin\cap\Mrout$,$\hat{F}_1(x)=\Fc(x_2)$ and $x_0=\Fc^{m(\Fc^{-1}(x))-1}(x)$ with $m(\Fc^{-1}(x))\ge3$. Since there exists an $x_{-1}\dfn\Fc^{-1}(x_0)$ s.t. $x_{-1}=\Fc^{m(\Fc^{-1}(x))-2}(x)$ and $dx_{-1}=D_{x}\Fc^{m(\Fc^{-1}(x))-2}(\underbracket{dx}_{\mathclap{\in\hat{C}^u_{1,x}\subset\mathcal{Q}_x}})$ with $D_{x}\Fc^{m(\Fc^{-1}(x))-2}=\begin{pmatrix}
1 & 2(m(\Fc^{-1}(x))-2) \\
0 & 1 
\end{pmatrix}$, we get $dx_{-1}\in\mathcal{HQ}_{x_{-1}}(\mathrm{I,III})$.

In case (b), $G=D_{x_{-1}}\Fc^4$ is a negative matrix by \cref{caseB}. \cref{Proposition:DF4ConeUpperBound} and \eqref{eq:casebdtheta3dphi3ratio} show for $x_3=\Fc^4(x_{-1})=\hat{F}_1(x)$ that $(d\phi_3,d\theta_3)\nfd dx_3=D_x\hat{F}_1(dx)=D_{x_{-1}}\Fc^4(D_{x}\Fc^{m(\Fc^{-1}(x))-2}(dx))=D_{x_{-1}}\Fc^4(dx_{-1})$ satisfies $0<\frac{d\theta_3}{d\phi_3}<1$, so $dx_3$ is in the interior of $\hat{C}^u_{1,\hat{F}_1(x)}$.

In case (c), $D_{x_{-1}}\Fc^{n_1+4}$ is a positive matrix by \cref{caseC}. \cref{theorem:B3minusdx3dphi3ratiobound} shows for $x_3=\Fc^4(x_{-1})=\hat{F}_1(x)$ that \[D_x\hat{F}_1(dx)=D_{x_{-1}}\Fc^{4+n_1}(D_{x}\Fc^{m(\Fc^{-1}(x))-2}(dx))=D_{x_{-1}}\Fc^{4+n_1}(dx_{-1})=dx_3\nfd(d\phi_3,d\theta_3)\] and $0<\frac{d\theta_3}{d\phi_3}<\lambda_2$, so $D_x\hat{F}_1(dx)\in\Int(\hat{C}^u_{1,\hat{F}_1(x)})$.
\newline
\textbf{Proof of (ii)}: The argument is the same as in \cref{thm:MroutreturnUniformExpansion}. The return orbit segment \eqref{eqhatF1returnorbit} contains $x_0\in\Mrout$, $x_1=\Fc(x_0)\in\MRin$, $x_2=\Fc^{n_1+1}(x_1)\in\Mrin$. And if $x_0\in\Nin$, $\hat{F}_1(x)=\Fc(x_2)\nfd x_3$ and $x_{-1}=\Fc^{-1}(x_0)$ are also contained in the return orbit segment \eqref{eqhatF1returnorbit}.

The conclusions of \cref{TEcase_a0,TEcase_a1,TEcase_b,TEcase_c} also hold in the present context:

In case (a0) and (a1) $d_1\ge2r$, $\frac{\|D_x\hat{F}_1(dx)\|_{\p}}{\|dx\|_{\p}}\ge\frac{\|dx_2\|_{\p}}{\|dx_0\|_{\p}}>0.26$.

In case (b) $d_1<2r$ and $n_1=0$, $\frac{\|D_x\hat{F}_1(dx)\|_{\p}}{\|dx\|_{\p}}\ge\frac{\|dx_3\|_{\p}}{\|dx_{-1}\|_{\p}}>0.05$.

In case (c) $d_1<2r$ and $n_1\ge1$, $\frac{\|D_x\hat{F}_1(dx)\|_{\p}}{\|dx\|_{\p}}\ge\frac{\|dx_3\|_{\p}}{\|dx_{-1}\|_{\p}}>0.9$.

\cref{caseA,caseB,caseC,TEcase_a0,TEcase_a1,TEcase_b,TEcase_c,corollary:halfquadrantcone,Corollary:CollisionOnMrin0Mrin1} imply the following:
\begin{itemize}
    \item For nongsingular $x\in\hat{M}_1$, we denote by $k_1(x)\dfn\inf\big\{k>0\bigm|\Fc^{k_1(x)}(x)\in\hat{M}\big\}$ and $d\hat{x}\dfn D_x\Fc^{k_1(x)}(dx)$. For $\forall dx\in C^u_{1,x}$, $d\hat{x}_1=D_x\Fc^{k(x)}(dx)\in C_x$ and $\frac{\|d\hat{x}_1\|_{\p}}{\|dx\|_{\p}}>0.05$,
    \item For nongsingular $x\in\hat{M}$, we denote by $k(x)\dfn\inf\big\{k>0\bigm|\Fc^{k(x)}(x)\in\hat{M}_1\big\}$ and $d\hat{x}_1\dfn D_x\Fc^{k(x)}(dx)$. For $\forall dx\in C_{x}$, $d\hat{x}_1=D_x\Fc^{k(x)}(dx)\in \hat{C}^u_{1,x}$ and $\frac{\|d\hat{x}_1\|_{\p}}{\|dx\|_{\p}}>0.05$,
\end{itemize}
where $C_x$ is the cone defined in \cref{RMASUHwithProof}.
    
Since we can find that every 2 consecutive return orbit segments \eqref{eqhatF1returnorbit} contain at least one element in $\hat{M}$. The $\p$-metric uniform exponential expansion on $\hat{M}$ by \cref{RMASUHwithProof} implies the $\p$-metric uniform exponential expansion on $\hat{M}_1$.
\end{proof}

\begin{corollary}[Corollary of \cref{thm:UniformExpansionOnMhat1}]\label{corollary:sinpmetricexpansionMhat1}
There exist constants $\hat{c}_1>0$ and $\Lambda_1>1$ such that $\dfrac{\sin\theta\|D_x\hat{F}_1^n(dx)\|_\p}{\|dx\|_\p}>\hat{c}_1\Lambda_1^n$ for all nonsingular $x=(\phi,\theta)\in\hat{M}_1$, $dx\in \hat{C}^u_{1,x}$ and $n\in\mathbb{N}$.
\end{corollary}
\begin{proof}
The proof is the same as \cref{corollary:sinpmetricexpansion}.
\end{proof}

\begin{corollary}[
Uniform expansion in Euclidean Metric]\label{Col:EuclideanExpansionMhat1}
Under the return map $\hat{F}_1$ on $\hat{M}_1$ the vectors in the unstable cone $\hat{C}^u_{1,x}$ from \cref{thm:UniformExpansionOnMhat1} expand uniformly with respect to the Euclidean metric $\|\,\|$, i.e., for nonsingular $x\in\hat{M}_1$, $dx=(d\phi,d\theta)\in \hat{C}^u_{1,x}$ and all $n\ge0$, there exist $c_2>0,\Lambda_1>1$ such that $\dfrac{\|D_x\hat{F}_1^n(dx)\|}{\|dx\|}>c_2\Lambda_1^n$.
\end{corollary}
\begin{proof}
    The proof is the same as \cref{Corollary:EuclideanExpansionMhat}.
\end{proof}

\begin{corollary}[Expansion between $\hat{M}$ and $\hat{M}_1$ in $\p$-metric and Euclidean metric]\label{corollary:ExpansionBetweenMhatminusMhatplus}\hfill 

\begin{enumerate}
    \item\label{item01ExpansionBetweenMhatminusMhatplus} There exists a constant $c_3>0$ such that if $(\phi,\theta)=x\in\hat{M}$ is nonsingular, $dx\in C_x$ as defined in \cref{RMASUHwithProof}, $j(x)=\inf\big\{i>0\bigm|\Fc^i(x)\in\hat{M}_1\big\}$, then $\frac{\sin\theta\|D_x\Fc^{j(x)}(dx)\|_{\p}}{\|dx\|_{\p}}>c_3$ and $\frac{\|D_x\Fc^{j(x)}(dx)\|}{\|dx\|}>c_3$.
    \item\label{item02ExpansionBetweenMhatminusMhatplus} There exists a constant $c_4>0$ such that if  $(\phi,\theta)=x\in\hat{M}_1$ is nonsingular, $dx\in \hat{C}^u_{1,x}$ as defined in \cref{thm:UniformExpansionOnMhat1},  $j_1(x)=\inf\big\{i>0\bigm|\Fc^i(x)\in\hat{M}\big\}$, then $\frac{\|D_x\Fc^{j_1(x)}(dx)\|_{\p}}{\|dx\|_{\p}}>c_4$ and $\frac{\|D_x\Fc^{j_1(x)}(dx)\|}{\|dx\|}>c_4$.
\end{enumerate}
\end{corollary}
\begin{proof}
    For the first conclusion, note that the trajectory segment $x\in\hat{M}_1,\Fc(x),\cdots,\Fc^{j_1(x)}(x)\in\hat{M}_1$ contains $x_0\in\Mrout$, $x_1=\Fc(x_0)\in \MRin$ and $x_2=\Fc^{n_1+1}(x_1)\in\Mrin$ with $n_1=\inf\big\{k\ge0:p(\Fc^i(x_1)),i=0,\cdots,k\in\Gamma_R\big\}$ similar to that in \cref{def:MhatReturnOrbitSegment}. Thus, the same analysis for \cref{TEcase_a0,TEcase_a1,TEcase_b,TEcase_c,corollary:sinpmetricexpansion,Corollary:EuclideanExpansionMhat} yields the first conclusion, 

    For the second conclusion, if the trajectory segment $x\in\hat{M}_1,\Fc(x),\cdots,\Fc^{j(x)}(x)\in\hat{M}$ are all collisions on $\Gamma_r$, then since $\hat{C}^u_{1,x}$ is a subset of positive quadrant $\mathcal{Q}_x(\mathrm{I,III})$ and $D_x\Fc^{j(x)}=\begin{pmatrix}
1 & 2j(x) \\
0 & 1 
\end{pmatrix}$ in $\phi\theta$-coordinates, $\frac{\|D_x\Fc^{j_1(x)}(dx)\|_{\p}}{\|dx\|_{\p}}>1$ and $\frac{\|D_x\Fc^{j_1(x)}(dx)\|}{\|dx\|}>1$. Otherwise, $x\in\hat{M}$, $\Fc(x)$, $\cdots$, $\Fc^{j(x)}(x)\in\hat{M}$ contains $x_0\in\Mrout$, $x_1=\Fc(x_0)\in \MRin$ and $x_2=\Fc^{n_1+1}(x_1)\in\Mrin$ with $n_1=\max\big\{k\ge0:p(\Fc^i(x_1))\in\Gamma_R\text{ for } i=0,\cdots,k\big\}$, then the same analysis for \cref{TEcase_a0,TEcase_a1,TEcase_b,TEcase_c,RMASUHwithProof,corollary:sinpmetricexpansion,Corollary:EuclideanExpansionMhat} yields the second conclusion. 
\end{proof}

\begin{theorem}[Uniform expansion on \Mrout]\label{thm:MroutreturnUniformExpansion}
    For the lemon billiards $\LEM(r,R,\phistar)$ satisfying \cref{def:UniformHyperbolicLemonBilliard} and for the return map $\hat{F}_1$ on, the return map $\tilde{F}$ on $\tilde{M}=\Mrout$ satisfies the following.
    \begin{enumerate}[label={(\roman*)}]
    \item\label{item01:MroutreturnUniformExpansion} For nonsingular $x\in\Mrout$, if $dx\in \tilde{C}^{u}_x\dfn\begin{cases}
            \big\{(d\phi,d\theta)\in T_x\Mrout\bigm|\frac{d\theta}{d\phi}\in[0,1]\big\}\text{, if } x\in \Mrout\smallsetminus\Nout,\\
            \big\{(d\phi,d\theta)\in T_x\Mrout\bigm|\frac{d\theta}{d\phi}\in[0,\frac{1}{3}]\big\} \text{, if } x\in\Nout, 
        \end{cases}$ \newline then $D_x\tilde{F}(dx)\in\big\{\text{interior of }\tilde{C}^{u}_{\tilde{F}(x)}\big\}$.
    \item\label{item02:MroutreturnUniformExpansion} For nonsingular $x\in\Mrout$ and $dx\in \tilde{C}^{u}_x$, there exists $\tilde{c}>0$ such that $\frac{\|D_x\tilde{F}^n(dx)\|_\p}{\|dx\|_\p}>\tilde{c}\Lambda^n$ for all $n\ge1$, where $\Lambda$ is as in \cref{thm:UniformExpansionOnMhat}.
    \end{enumerate}
\end{theorem}
\begin{proof}
\textbf{\ref{item01:MroutreturnUniformExpansion}}: 
For any nonsingular $y\in\Mrout$, the $y$ return to $\Mrout$ orbit segment is defined as follows. \begin{equation}\label{eq:yMRoutreturnorbitsegment}
\underbracket{y=x_0}_{\mathclap{\in\Mrout}},\underbracket{\Fc(y)=x_1,\cdots}_{\notin\Mrout},\underbracket{\Fc^{\tilde{\sigma}(y)}(y)=\tilde{F}(y)\nfd\tilde{y}_1}_{\in \Mrout},\text{ where }\tilde{\sigma}(y)\dfn\inf\{k>0\mid \Fc^{k}(y)\in\Mrout\}.
\end{equation} and similar as \eqref{eqPtsFromeqMhatOrbSeg}, \eqref{eq:yMRoutreturnorbitsegment} also contains $x_1=\Fc(y)=\Fc(x_0)\in\MRin$ and $x_2=\Fc^{n_1+1}(x_1)\in\Mrin$. In \eqref{eq:yMRoutreturnorbitsegment} there are also the same (a0)(a1)(b)(c) cases in \eqref{eqMainCases} based on $x_1=(\Phi_1,\theta_1)$ and $n_1$. \cref{corollary:InvQuadrant,corollary:halfquadrantcone}, \eqref{eq:strictlyinvcone_casea1bc} and the same argument as in \cref{thm:UniformExpansionOnMhat} imply the strict invariance of $\tilde{C}^{u}_x$.

\textbf{\ref{item02:MroutreturnUniformExpansion}}: 
In \eqref{eq:yMRoutreturnorbitsegment}, there is a $\hat{y}_1\in\hat{M}$ with 
\begin{equation}\label{eq:hatydefinitionInOrbit}
\begin{aligned}
    \hat{y}_1&=\left\{
    \begin{aligned}
    &\tilde{y}_1,\text{ if }\tilde{y}_1\in\Mrin\cap\Mrout\\ 
    &\Fc^{-1}(\tilde{y}_1).\text{ if }\tilde{y}_1\in\Mrout\smallsetminus\Mrin.
    \end{aligned}\right.
\end{aligned}
\end{equation}
Suppose $T_{\hat{y}_1}M\ni\hat{y}_1=\Fc^{t_1}(y)$ with $t_1\in\mathbb{N}$. Then for $dy\in \tilde{C}^u_{y}$, let $\hat{y}_1\dfn D_y\Fc^{t_1}(dy)$ and $d\tilde{y}_1\dfn D_y\tilde{F}(dy)$.

For all cases (a0)(a1)(b)(c) of \eqref{eq:yMRoutreturnorbitsegment}, the analysis as in \cref{lemma:ExpansionLowerBoundA0A1,theorem:UniformExpansionA0,theorem:proofofTEcase_b,theorem:contractioncontrolA1,lemma:ExpansionLowerBoundB,lemma:AformularForLowerBoundOfExpansionForn1largerthan0,proposition:casecexpansionn1eq1,proposition:casecexpansionn1eq2,proposition:casecexpansionn1eq3,proposition:casecexpansionn1ge4} implies that $\frac{\|d\hat{y}_1\|_{\p}}{\|dy\|_{\p}}>0.05$. 

The $n$ times $\tilde{F}$ return orbit segment for $y$ is \[\underbracket{y\dfn\tilde{y}_0}_{\in\Mrout},\cdots,\underbracket{\tilde{y}_1\dfn\tilde{F}(y)}_{\in\Mrout},\cdots,\underbracket{\tilde{y}_n\dfn\tilde{F}^n(y)}_{\in\Mrout}\] and it contains $\underbracket{\hat{y}_1}_{\in\hat{M}},\cdots,\underbracket{\hat{y}_n}_{\in\hat{M}}$ with each $\hat{y}_k$ defined similarly as in \eqref{eq:hatydefinitionInOrbit}. \begin{align*}
    \hat{y}_k&=\left\{
    \begin{aligned}
    &\tilde{y}_k,\text{ if }\tilde{y}_k\in\Mrin\cap\Mrout\\ 
    &\Fc^{-1}(\tilde{y}_k).\text{ if }\tilde{y}_k\in\Mrout\smallsetminus\Mrin.
    \end{aligned}\right.
\end{align*}

With $d\tilde{y}_k=D_y\tilde{F}^{k}(dy)$, $T_{\hat{y}_k}M\ni d\hat{y}_k=D_{\tilde{y}_{k-1}}\Fc^{t_k}(d\tilde{y}_{k-1})$, $t_k>0$, $k=1,\cdots,n$. 

By \cref{corollary:halfquadrantcone}, $d\hat{y}_n\in C_{\hat{y}_n}=\mathcal{HQ}_{\hat{y}_n}(\mathrm{I,III})$. And since either $d\hat{y}_n=d\tilde{y}_n$ or $D_{\hat{y}_n}\Fc(d\hat{y}_n)=d\tilde{y}_n$ with $D_{\hat{y}_n}\Fc=\begin{pmatrix}
1 & 2 \\
0 & 1 
\end{pmatrix}$, $\|d\tilde{y}_n\|_{\p}\ge \|d\hat{y}_n\|_{\p}$ by \cref{proposition:ShearOneStepExpansion}. We get 
\[\frac{\|D_y\tilde{F}^n(dy)\|_\p}{\|dy\|_{\p}}=\overbracket{\frac{\|d\tilde{y}_n\|_{\p}}{\|d\hat{y}_n\|_{\p}}}^{\ge1}\cdot\underbracket{\frac{\|d\hat{y}_n\|_{\p}}{\|d\hat{y}_1\|_{\p}}}_{\mathclap{\text{\cref{thm:UniformExpansionOnMhat}: }>c\Lambda^{n-1}}}\cdot\overbracket{\frac{\|d\hat{y}_1\|_{\p}}{\|dy\|_{\p}}}^{>0.05}>1\cdot c\Lambda^{n-1}\cdot0.05=\underbracket{\tilde{c}}_{\mathrlap{0.05c/\Lambda}}{\Lambda}^n.\qedhere
\]
\end{proof}
\begin{corollary}\label{corollary:MroutreturnUniformExpansion}In the context of \cref{thm:MroutreturnUniformExpansion} the following hold. Under the return map $\tilde{F}$ on $\Mrout$ the vectors in the unstable cone $\tilde{C}_x^u$ defined in \cref{thm:MroutreturnUniformExpansion} expand uniformly with respect to the Euclidean metric $\|\,\|$, i.e., for nonsingular $(\phi,\theta)=x\in\MRout$, $dx=(d\phi,d\theta)\in \tilde{C}^u_x$ and all $n\ge0$, there exist $\tilde{c}_0,\tilde{c}_1>0,\tilde{\Lambda}_0>1$ such that
\[\dfrac{\sin\theta\|D_x\tilde{F}^n(dx)\|_{\p}}{\|dx\|_{\p}}>\tilde{c}_0\tilde{\Lambda}_0^n,\qquad
\dfrac{\|D_x\tilde{F}^n(dx)\|}{\|dx\|}>\tilde{c}_1\tilde{\Lambda}_0^n.
\]
\end{corollary}
\begin{proof}
    The proof is the same as \cref{corollary:sinpmetricexpansionMhat,Col:EuclideanExpansionMhat} using \eqref{eq:casea0onestepsinexpansion}.
\end{proof}

\begin{lemma}[Symmetries between $\hat{M}_1$ and $\hat{M}$ from \cref{RMASUHwithProof} (see \cref{fig:hatUplusUminus})]\label{lemma:symmetriesbetweenMhatMhat1} Both symmetries $I$ and $J$ from \eqref{eq:InversionSymmetryDef} and  \eqref{eq:FoldSymmetryDef} conjugate $\hat{F}$ on $\hat{M}$ and $\hat{F}^{-1}$ on $\hat{M}_1$. i.e. $I\circ \hat{F}_1=\hat{F}^{-1}\circ I$, $J\circ \hat{F}_1=\hat{F}^{-1}\circ J$.
\end{lemma}
\begin{proof}
    Note that both $I$, $J$ are isometries between $\hat{M}$ and $\hat{M}_1$ since $J\Mrin=\Mrout=I\Mrin$ and $J\Mout_{r,i}=\Min_{r,i}=I\Mout_{r,i}$, $\forall i\ge0$, therefore $I(\hat{M})=\hat{M}_1=J(\hat{M})$ (also see \cref{fig:MrN}). 
    
    $\forall$ nonsingular point $x\in\hat{M}_1$, let $y=I(x)\in\hat{M}$ and then $x=I(y)$. Since $I$, $J$ are conjugates between $\Fc$ and $\Fc^{-1}$ and $I\circ I=Id=J\circ J$. 
    
    Therefore, $\forall k>0$, $\Fc^k(x)=\Fc^k(I(y))\overbracket{=}^{\mathclap{\eqref{eq:InversionSymmetryDef}}}I\Fc^{-k}(y)$, and thus \begin{equation}\label{eq:F1hatFhatConjugate}
        \Fc^{k}(x)\in\hat{M}_1\text{ if only if }\Fc^{-k}(y)\in I^{-1}(\hat{M}_1)=\hat{M}. 
    \end{equation}
    Note that $\sigma_1(x)\ge2$ and by \eqref{eqhatF1returnorbit}: $\sigma_1(x)=\inf\{k>0\mid\Fc^k(x)\in\hat{M}_1\}$. This means that for $k=1,\cdots,\sigma_1(x)-1$, $\Fc^{k}(x)\notin\hat{M}_1$ and $\Fc^{\sigma_1(x)}(x)\in\hat{M}_1$.
    
    Therefore, by \eqref{eq:F1hatFhatConjugate} for $k=1,\cdots,\sigma_1(x)-1$, $\Fc^{-k}(y)\notin\hat{M}$ and $\Fc^{-\sigma_1(x)}(y)\in\hat{M}$. With $y\in\hat{M}$ and by definition \eqref{eq:F1hatFhatConjugate}, we see $\sigma_1(x)=\sigma(\hat{F}^{-1}(y))$. We further get the following. \begin{align*}I\hat{F}_1(x)&\overbracket{=}^{\mathclap{\eqref{eq:F1hatFhatConjugate}}}I\circ \Fc^{\sigma_1(x)}(x)\overbracket{=}^{\mathclap{\eqref{eq:InversionSymmetryDef}}}\Fc^{-\sigma_1(x)}\circ I(x)\overbracket{=}^{\mathclap{I(x)=y}}\Fc^{-\sigma_1(x)}(y)\overbracket{=}^{\mathclap{\sigma_1(x)=\sigma(\hat{F}^{-1}(y))}}\hat{F}^{-1}(y)\underbracket{=}_{\mathclap{y=I(x)}}\hat{F}^{-1}I(x)\\
    \end{align*}
    
    Let $z=J(x)$, by the exactly same preceding reasoning and using \eqref{eq:FoldSymmetryDef}, we get $J\hat{F}_1(x)=\hat{F}^{-1}J(x)$.
\end{proof}

Using symmetries $I$ and $J$, we conclude that the stable cone family on $\hat{M}$ is strictly invariant under $\hat{F}^{-1}$.
\bibliographystyle{alpha}
\bibliography{lemonExpansion}  
\end{document}